\magnification 1100 \input eplain.tex \input miniltx.tex \input pictex \input amssym \input color.sty \input graphicx
\def \nocirc {}   \font \bbfive = bbm5 \font \bbseven = bbm7 \font \bbten =
bbm10 \font \eightbf = cmbx8 \font \eighti = cmmi8 \skewchar \eighti = '177 \font \eightit = cmti8 \font \eightrm = cmr8
\font \eightsl = cmsl8 \font \eightsy = cmsy8 \skewchar \eightsy = '60 \font \eighttt = cmtt8 \hyphenchar \eighttt = -1
  \font \sixi = cmmi6 \skewchar \sixi = '177 \font \sixrm = cmr6 \font \sixsy =
cmsy6 \skewchar \sixsy = '60 \font \tensc = cmcsc10   \font \titlefont =
cmbx12 \scriptfont \bffam = \bbseven \scriptscriptfont \bffam = \bbfive \textfont \bffam = \bbten \newskip \ttglue \def
\eightpoint {\def \rm {\fam 0 \eightrm }\textfont 0= \eightrm \scriptfont 0 = \sixrm \scriptscriptfont 0 = \fiverm
\textfont 1 = \eighti \scriptfont 1 = \sixi \scriptscriptfont 1 = \fivei \textfont 2 = \eightsy \scriptfont 2 = \sixsy
\scriptscriptfont 2 = \fivesy \textfont 3 = \tenex \scriptfont 3 = \tenex \scriptscriptfont 3 = \tenex \def \it {\fam
\itfam \eightit }\textfont \itfam = \eightit \def \sl {\fam \slfam \eightsl }\textfont \slfam = \eightsl \def \bf {\fam
\bffam \eightbf }\textfont \bffam = \bbseven \scriptfont \bffam = \bbfive \scriptscriptfont \bffam = \bbfive \def \tt
{\fam \ttfam \eighttt }\textfont \ttfam = \eighttt \tt \ttglue = .5em plus.25em minus.15em \normalbaselineskip = 9pt
\def \MF {{\manual opqr}\-{\manual stuq}}\let \sc = \sixrm \let \big = \eightbig \setbox \strutbox = \hbox {\vrule
height7pt depth2pt width0pt}\normalbaselines \rm } \def \setfont #1{\font \auxfont = #1 \auxfont } \font \bigrm = cmr12
scaled 1200 \font \amsbig = msam10 scaled 1200 \font \euler = eufm10 \font \eulersmall = eufm10 scaled 833 \font \sixtt
= cmtt6 \font \headlnfont = cmbcsc10 scaled 1200 \font \chapfont = cmbcsc10 scaled 1200 \font \chapnumfont = cmbcsc10
scaled 900 \font \rs = rsfs10 scaled 950 \font \rsfootnote = rsfs10 scaled 800 \font \rssmall = rsfs10 scaled 600 \font
\smallcal = cmsy9 \font \rsten =rsfs10 \font \baita = cmbx10 scaled 3000 \hyphenation {u-ni-que equip-ped associ-ati-ve
alge-bra-ic ame-na-bi-li-ty ho-mo-mor-phism non-empty no-where pa-ra-me-tri-zed cros-sed} \def \blank {}  \def \newpage {\vfill \eject } \def \ucase #1{\edef \auxvar {\uppercase {#1}}\auxvar } \def
\lcase #1{\edef \auxvar {\lowercase {#1}}\auxvar } \def \noindentNext {\parindent 0pt \everypar {\parindent 20pt
\everypar {}}} \def \goOdd {} \def \syschapter #1#2{\endgroup \nopagenumbers \vfill \eject \goOdd \null \vskip 3truecm
\global \advance \secno by 1 \stno = 0 \begingroup \baselineskip 18pt \raggedright \pretolerance 10000 \hbadness =5000
\noindent \nocirc \Htarget {chapter.\number \secno }{\chapfont \number \secno .\enspace \ucase {\spaceskip =0.5em #1}}
\par \endgroup \begingroup \write 2 {\string \ChapterEntry {#1} {\number \pageno }} \Headlines {{\chapnumfont \number
\secno .}\ #2}{\GlobalRunninghead } \vskip 1truecm \noindentNext } \def \chapter #1 \par {\syschapter {#1}{#1}} \long
\def \state #1 #2\par {\begingroup \def \InsideBlock {} \medbreak \goodbreak \noindent \advseqnumbering \bf \current
.\enspace #1.\enspace \sl #2\par \endgroup \medbreak } \def \definition #1\par {\state Definition \rm #1\par } \def
\nostate {\medbreak \goodbreak \noindent \advseqnumbering {\bf \current .\enspace }} \def \nrem {\bigskip \bigskip
\noindent {\it Notes and remarks.\enspace }} \long \def \Proof #1\endProof {\begingroup \def \InsideBlock {} \medbreak
\noindent {\it Proof.\enspace }#1 \ifmmode \eqno \square $$ \else \hfill $\square $ \looseness = -1 \fi \medbreak
\endgroup } \def \omitDoubleDollar {\ifmmode \else This should have been Math Mode \end \fi } \def \$#1{#1 $$$$ #1} \def
\explica #1#2{\mathrel {\buildrel \hbox {\sixrm (#2)} \over #1}} \def \pilar #1{\vrule height #1 width 0pt} \def \stake
#1{\vrule depth #1 width 0pt} \newcount \fnctr \fnctr = 0 \def \fn #1{\global \advance \fnctr by 1 \edef \footnumb
{$^{\number \fnctr }$}\footnote {\footnumb }{\eightpoint #1\par \vskip -10pt}} \def \ifundef #1{\expandafter \ifx
\csname #1\endcsname \relax } \def \empty {} \def \ifempty #1{\edef \aux {#1}\ifx \aux \empty } \def \ControlSequence
#1{\ifundef {#1}\raise 3pt \hbox {\sixtt \textbackslash #1}\message {*.*.* Label [#1] is undefined *.*.*}\else \csname
#1\endcsname \fi } \newcount \secno \secno = 0 \newcount \stno \stno = 0 \newcount \eqcntr \eqcntr = 0 \def \track
#1#2#3{\ifundef {#1}\else \hbox {\sixtt [#2 \string #3 ] }\fi } \def \advseqnumbering {\global \advance \stno by 1
\global \eqcntr =0 \relax } \def \current {\number \secno \ifnum \number \stno = 0 \else .\number \stno \fi } \def
\laberr #1{\message {*.*.* RELABEL CHECKED FALSE for #1 *.*.*} RELABEL CHECKED FALSE FOR #1, EXITING.  \end } \def
\deflabel #1#2{\ifundef {#1}\global \expandafter \edef \csname #1\endcsname {#2}\Htarget {link.#1}{}\else \edef
\deflabelaux {\expandafter \csname #1\endcsname }\edef \deflabelbux {#2}\ifx \deflabelaux \deflabelbux \else \laberr
{#1=(\deflabelaux )=(\deflabelbux )} \fi \fi \track {showlabel}{*}{#1}} \def \equationmark #1 {\ifundef {InsideBlock}
\advseqnumbering \eqno {(\current )} \deflabel {#1}{\current } \else \global \advance \eqcntr by 1 \edef \subeqmarkaux
{\current .\number \eqcntr } \eqno {(\subeqmarkaux )} \deflabel {#1}{\subeqmarkaux } \fi } \def \label #1 {\deflabel
{#1}{\current }} \def \aToL #1{\ifcase #1 O\or I\or D\or T\or Q\or C\or G\or T\or B\or N\fi } \def \ntoP
#1#2#3#4#5;{\aToL #1\aToL #2\aToL #3\aToL #4#5} \def \CodedRefToNumber #1#2{\def \aux {#1\expandafter \ntoP
#2;}\ControlSequence {\aux }} \def \getparam #1/#2/#3/#4/#5;{\global \edef \paramOne {#1}\global \edef \paramTwo
{#2}\global \edef \paramThree {#3}\global \edef \paramFour {#4}} \def \sysref #1#2#3{\ifempty {#3}\def \extra {}\edef
\au {#1}\edef \yr {#2}\else \def \extra {, #1}\edef \au {#2}\edef \yr {#3}\fi [{\rm \Hyperref {paper}{\au /\yr
}{\CodedRefToNumber {\au }{\yr }}\extra }]} \def \ref #1{\getparam #1////;\sysref {\paramOne }{\paramTwo }{\paramThree
}} \def \syscite #1#2{\Hyperref {link}{#1}{\ControlSequence {#1}\ifempty {#2}\else .#2\fi }} \def \cite #1{\getparam
#1////;(\syscite {\paramOne }{\paramTwo })} \def \_#1#2{\explica {#1}{\getparam #2////;\syscite {\paramOne }{\paramTwo
}}} \def \=#1{\_={#1}} \newcount \zitemno \zitemno = 0  \def \Item #1{\smallskip \item {{\rm
#1}}} \def \izitem {\global \zitemno = 0} \def \zitemplus {\global \advance \zitemno by 1 \relax } \def \rzitem
{\romannumeral \zitemno } \def \rzitemplus {\zitemplus \rzitem } \def \zitem {\Item {{\rm (\rzitemplus )}}} \def \Zitem
{\izitem \zitem } \newcount \nitemno \nitemno = 0  \def \nitem {\global \advance \nitemno by
1 \Item {{\rm (\number \nitemno )}}} \newcount \aitemno \aitemno = -1 \def \boxlet #1{\hbox to 6.5pt{\hfill #1\hfill }}
\def \iaitem {\aitemno = -1} \def \aitemconv {\ifcase \aitemno a\or b\or c\or d\or e\or f\or g\or h\or i\or j\or k\or
l\or m\or n\or o\or p\or q\or r\or s\or t\or u\or v\or w\or x\or y\or z\else zzz\fi } \def \aitem {\global \advance
\aitemno by 1\Item {(\boxlet \aitemconv )}} \def \zitemBismark #1 {\deflabel {#1Loc}{\rzitem }\deflabel {#1}{\current
.\rzitem }} \def \aitemBismark #1 {\deflabel {#1Loc}{\aitemconv }\deflabel {#1}{\current .\aitemconv }}  \hsize = 12.9truecm \vsize = 20.6truecm \voffset = -0.5truecm \output {\hoffset
=1.8truecm \voffset =1.3truecm \shipout \vbox {\makeheadline \pagebody \makefootline } \advancepageno \ifnum
\outputpenalty >-20000\else \dosupereject \fi } \def \partpage #1#2{ \begingroup \newpage \goOdd \headline {\hfill }
\nopagenumbers \null \vfill \centerline {\baita PART #1} \bigskip \centerline {\baita ---} \bigskip \centerline {\baita
#2} \vfill \null \newpage \endgroup \write 2 {\vskip 1cm \centerline {\titlefont #1 -- #2} \vskip 0.5cm}} \def
\Headlines #1#2{\nopagenumbers \global \headline {\ifodd \pageno \headlnfont \lcase {#1} \hfill \folio \else \headlnfont
\folio \hfill \lcase {#2}\fi }} \def\write#1#2{} \newcount \indexno \indexno 0 \def \ChapterEntry #1#2{\medskip
\noindent \advance \indexno by 1 \line {\hbox to 0.5cm{\hfill \number \indexno .}\enspace \Hyperref {chapter}{\number
\indexno }{#1} \dotfill \enspace #2}} \def \OtherIndexEntry #1#2{\medskip \noindent \advance \indexno by 1 \line {\hbox
to 0.5cm{\hfill }\enspace \Hyperref {extra}{#1}{#1} \dotfill \enspace #2}} \def \subjex #1#2{\Htarget {subj.#2}{\it
#1}\write 1 {#2/\number \pageno }} \def \subj #1{\subjex {#1}{#1}} \def \*{\otimes } \def \ds {\displaystyle } \def \and
{\hbox {,\quad and \quad }} \def \calcat #1{\,{\vrule height8pt depth4pt}_{\,#1}} \def \imply {\kern 7pt \Rightarrow
\kern 7pt} \def \iff {\ \Leftrightarrow \ } \def \umlaut #1{{\accent "7F #1}} \def \breve #1{{\accent "15 #1}} \def
\"#1{{\it #1}\/} \def \inv {^{-1}} \def \for #1{,\quad \forall \,#1} \def \bfdef #1{\global \expandafter \edef \csname
#1\endcsname {{\bf #1}}} \bfdef N \bfdef Z \bfdef C \bfdef R \def \rsbox #1{{\mathchoice {\hbox {\rs #1}} {\hbox {\rs
#1}} {\hbox {\rssmall #1}} {\hbox {\rssmall #1}}}} \def \eubox #1{{\mathchoice {\hbox {\euler #1}} {\hbox {\euler #1}}
{\hbox {\eulersmall #1}} {\hbox {\eulersmall #1}}}} \def \trepa #1#2{\buildrel #1 \over {#2}} \def \>#1{\setbox 0\hbox
{$#1$}\dimen 1=\dp 0 \advance \dimen 1 by 2pt \trepa {#1 \stake {\dimen 1 }} {\ifdim \wd 0<8pt\rightarrow \else
\longrightarrow \fi }} \long \def \pix #1{\noindent \hfill \beginpicture #1 \endpicture \hfill \null } \def \fix
{\noindent {\amsbig I} \kern 4pt } \def \resp #1{(resp.~#1)} \def \med #1{\mathop {\textstyle #1}\limits } \def
\medoplus {\med \bigoplus }  \def \medprod {\med \prod } \def \medvee {\med \bigvee } \def
\somalim #1{\mathop {\textstyle \sum }\limits _{#1}} \def \somanol #1{\mathop {\textstyle \sum }\nolimits _{#1}} \def
\soma #1{\mathchoice {\somalim {#1}} {\somanol {#1}} {\somanol {#1}} {\somanol {#1}}} \def \usomalim #1#2{\mathop
{\textstyle \sum }\limits _{#1}^{#2}} \def \usomanol #1#2{\mathop {\textstyle \sum }\nolimits _{#1}^{#2}} \def \somain
{\usoma {i=1}{n}} \def \modulogrande {\vrule height 10pt depth 6pt width 0.6pt} \def \normagrande {\modulogrande \kern
1.3pt\modulogrande } \def \normsum #1{\kern 3pt\normagrande #1 \kern 2pt\normagrande \kern 3pt} \def \lrm #1{_{\hbox
{\sixrm #1}}} \def \vertequal {\vrule height 6pt depth 0pt width 0.5pt \kern 1.5pt \vrule height 6pt depth 0pt width
0.5pt} \def \essIdeal {essential ideal} \def \nonDegAlg {non-degenerate} \def \nonDegRep {non-degenerate} \def
\nonDegHomo {non-degenerate} \def \diagram #1{\advseqnumbering \centerline {\eightrm \current .\enspace Diagram.\enspace
#1} \bigskip }  \def \D #1{\SymbolDom _{#1}} \def \Di #1{\D {#1\inv }} \def \SymbolAction {\theta }
\def \Th {\SymbolAction } \def \th #1{\Th _{#1}} \def \thi #1{\th {#1\inv }} \def \Graph {\hbox {Graph}} \def \Orb
{\hbox {Orb}} \def \F {{\bf F}} \def \IX {{\cal I}(X)} \def \bool #1{[#1]}  \def \Part
{\hbox {\rs P}} \def \IXt {{\rm pHomeo}(X)} \def \OG {\{0,1\}^G} \def \OuG {\Omega _1} \def \ideal {\mathrel
\trianglelefteq } \def \vezes {{\,\cdot \,}} \def \usoma #1#2{\mathchoice {\usomalim {#1}{#2}} {\usomanol {#1}{#2}}
{\usomanol {#1}{#2}} {\usomanol {#1}{#2}}} \def \K {{\bf K}} \def \mstar {{\eightrm (*-)}\kern 1pt} \def \tth #1{\Theta
_{#1}} \def \tthi #1{\tth {#1\inv }} \def \Mult {{\cal M}} \def \span {{\rm span}} \def \rt #1{{\ifundef {rtimes} \times
\else \rtimes \fi }} \def \Rt #1{{\ifundef {rtimes} \times \else \rtimes \fi _{#1}}} \def \Pr {u} \def \pr #1{\Pr _{#1}}
\def \pri #1{\pr {{#1}\inv }} \def \e #1{e _{#1}} \def \ei #1{e _{#1\inv }} \def \Kpar {\K \underpar (G)} \def \Apar
{A\underpar (G)} \def \ponto {{\cdot }} \def \clspan #1{[#1]} \def \IA {{\rm pAut}(A)} \def \art #1{{\ifundef {rtimes}
\times \else \rtimes \fi }\lrm {alg}} \def \maxnorm #1{\Vert #1\Vert \lrm {max}} \def \Lin {\hbox {\rs L}\kern 1pt} \def
\Gen {\hbox {\rs G}} \def \mx {{\vee }} \def \underpar {\lrm {par}} \def \Rel {{\cal R}} \def \SSat {{\cal R}\lrm {sat}}
\def \ORSat {\Omega _{\SSat }} \def \CstarparRel #1#2{C^*\underpar (#1,#2)} \def \ee #1{\varepsilon _{#1}} \def \G {G}
\def \OR {\Omega _\Rel } \def \FRel {{\cal F}_\Rel } \def \RME {\RiefMor equi\-va\-} \def \convrg #1{\trepa {#1\to
\infty } \longrightarrow } \def \Clspan #1{\big [#1\big ]} \def \AdjOp {\rsbox {L}\,} \def \Kp {{\cal K}} \def \subsub
#1{\kern -1pt_{_{#1}}} \def \ip #1#2#3{\langle #2,#3\rangle \subsub {#1}} \def \bip #1#2#3{\big \langle #2,#3\big
\rangle \subsub {#1}} \def \Bun {\rsbox {B}} \def \hj {\hat \jmath } \def \CCB {C_c(\Bun )} \def \CB {C^*(\Bun )} \def
\lred {\lrm {red}} \def \redrt #1{{\ifundef {rtimes} \times \else \rtimes \fi \lred }} \def \CrB {\Cr {\Bun }} \def
\tensor #1{\mathchoice {\mathop {\otimes }\limits #1} {\otimes #1}{\otimes #1}{\otimes #1}} \def \tmin {\tensor {\lmin
}} \def \tmax {\tensor {\lmax }} \def \lreg {\lambda ^{\kern -1pt\scriptscriptstyle G}} \def \Cr #1{\CstarRed (#1)} \def
\Ker {{\rm Ker}} \def \weird {\hbox {\rsten S\kern 1pt}} \def \tm {\tmax } \def \restr #1{\,{\vrule height 6pt depth
5pt}_{\,#1}} \def \arw {\to } \def \stateend {} \def \stateRef #1#2{\state #1 \label #2} \def \aproxNet {Cesaro net}
\def \sFb {sub-bun} \def \subFellBun {Fell \sFb } \def \segura #1{\mathop {\textstyle #1}} \def \Abun {{\rsbox {A}\,}}
\def \Jbun {{\mathchoice {\hbox {\smallcal J}\,} {\hbox {\smallcal J}\,} {\hbox {\rssmall J}} {\hbox {\rssmall J}}}}
\def \preFellBun {pre-Fell-bun} \def \PreFellBun {Pre-Fell-bun} \def \Cun {\rsbox {C}} \def \lmin {\lrm {min}} \def
\lmax {\lrm {max}} \def \BigCr #1{\CstarRed \Big (#1\Big )} \def \barsum #1{\overline {\pilar {9pt}#1}} \def \E
#1#2{e_{#1,#2}} \def \uE #1#2{1\otimes \E {#1}{#2}} \def \smsh {\Bun \kern 1pt\raise 1pt \hbox {$\scriptstyle \sharp
$}\kern 1pt G} \def \smshz {\Bun \kern 1pt\raise -2pt \hbox {$\buildrel \sharp \over {\scriptscriptstyle 0}$}\kern 1pt
G} \def \smshid {\Bun \kern 1pt{\mathchoice {\scriptstyle \flat }{\scriptstyle \flat }{\scriptscriptstyle \flat
}{\scriptscriptstyle \flat }}\kern .3pt G} \def \rTh {\Delta } \def \rEta {\Gamma } \def \iinv {^{^{-1}}\kern -2pt} \def
\rth #1{\rTh _{#1}} \def \rthi #1{\rth {#1\inv }} \def \RiefMor {Morita-Rieffel-} \def \barspan {\overline {\pilar
{6pt}\span }} \def \setmenos {{\setminus }} \def \ON {\Omega _\NicaRel } \def \NicaAlg {N(P)} \def \NicaRel {{\cal N}}
\def \FNRel {{\cal F}_\NicaRel } \def \WH {{\cal W}(P)} \def \s #1{s_{#1}} \def \p #1{p_{#1}} \def \TE {{\cal T}_E} \def
\TTE {\widetilde {\cal T}_E} \def \FinPath {E^\ast } \def \InfPath {E^\infty } \def \PathExponent {\sharp } \def \Path
{E^\PathExponent } \def \OE {\Omega _E} \def \CstarRed {C^*\lred } \def \x {\times } \def \FullPath {\widetilde
E^\PathExponent } \def \loc #1#2{{\rm loc}_{#2}(#1)} \def \locstd {\loc \omega g}  \def \CstarPath {E^\flat
} \def \CoefAlg {C_0(\Path )}

                    \def \GlobalRunninghead {Partial Dynamical Systems and Fell Bundles} \nopagenumbers

	\expandafter\edef\csname hypers@fe\endcsname{\catcode
												 `\noexpand @=\the\catcode`\@}%
	\catcode`\@=11
	%
	%
	\ifx\hyperd@ne\hyper@ndefined
	 \global\let\hyperd@ne=\relax
	\else
	 \errhelp{hyperbasics.tex needs to be included only once outside
			  of any {...} or \begingroup...\endgroup. You have tried to
			  include it more than once. If the previous include was indeed
			  outside any groupings, continue and all will be well.}%
	 \errmessage{Input this file only once!}%
	 \endinput
	\fi
	%
	%
	\def\hyperv@rsion{8}%
	%
	%
	\newread\hyperf@le
	\def\hyperf@lename{\jobname.hrf}%
	\immediate\openin\hyperf@le\hyperf@lename\relax
	\ifeof\hyperf@le\relax
	 \immediate\closein\hyperf@le\relax
	\else
	 \immediate\closein\hyperf@le\relax
	 \input \hyperf@lename
	\fi
	%
	%
	\newwrite\hyperf@le
	\immediate\openout\hyperf@le\hyperf@lename
	%
	%
	\newtoks\hypert@ks
	%
	%
	\edef\hypert@mp{\catcode`\noexpand\#=\the\catcode`\#}%
	\catcode`\#=12
	\def\hyperh@sh{#}%
	\hypert@mp
	\let\hypert@mp=\relax
	\let\hyper@nd=\relax
	\def\hyperstr@pquote"#1"#2\hyper@nd{\ifx\hyper@ndefined#2\hyper@ndefined#1\else
										\ifx\hyper@ndefined#1\hyper@ndefined
										\hyperstr@pquote#2"\hyper@nd\else
										#1\hyperstr@pquote"#2"\hyper@nd\fi\fi}%
	\def\hyperstr@pblank" #1 #2\hyper@nd"{\ifx\hyper@ndefined#2\hyper@ndefined#1\else
										\ifx\hyper@ndefined#1\hyper@ndefined
										\hyperstr@pblank"#2 \hyper@nd"\else
										#1\hyperstr@pblank" #2 \hyper@nd"\fi\fi}
	\long\def\hyper@nchor#1#2{\edef\hyperm@cro{html:<A #1>}%
							  \special\expandafter{\hyperm@cro}%
							  {#2}}%
	\def\hyper@atm@ning#1->#2\hyper@nd{#2}
	\def\hyperlink#1{\edef\hypert@mp{#1}%
				   \edef\hypert@mp{\expandafter\hyper@atm@ning\meaning\hypert@mp
								   \hyper@nd}%
				   \edef\hypert@mp"{ \expandafter\hyperstr@pquote\expandafter"%
								   \hypert@mp"\hyper@nd}%
				   \edef\hypert@mp{\expandafter\hyperstr@pblank\expandafter%
								   "\hypert@mp" \hyper@nd"}%
				   \hyper@nchor{href=\expandafter"\hypert@mp"}}%
	\def\hypertarget#1{\edef\hypert@mp{#1}%
				   \edef\hypert@mp{\expandafter\hyper@atm@ning\meaning\hypert@mp
								   \hyper@nd}%
				   \edef\hypert@mp"{ \expandafter\hyperstr@pquote\expandafter"%
								   \hypert@mp"\hyper@nd}%
				   \edef\hypert@mp{\expandafter\hyperstr@pblank\expandafter%
								   "\hypert@mp" \hyper@nd"}%
				   \hyper@nchor{name=\expandafter"\hypert@mp"}}%
	\def\hyperref{\afterassignment\hyperr@f\let\hyperp@ram}
	\def\hyperr@f{\ifx\hyperp@ram{\iffalse}\fi
				   \expandafter\expandafter\expandafter\hyperr@@
				   \expandafter{%
				  \else
				   \iffalse}\fi
				   \ifx\hyperp@ram\hyper@ndefined
					 \message{Undefined reference}%
					 \def\hyperp@r@m{{}{undefined}{}}%
				   \else
					 \edef\hyperp@r@m{\hyperp@ram}%
				   \fi
				   \expandafter\expandafter\expandafter\hyperr@@
				   \expandafter\hyperp@r@m
				  \fi}%
	\def\hyperr@@#1#2#3{\ifx\hyper@ndefined#1\hyper@ndefined
						\hypert@ks\expandafter{\hyperh@sh#2.#3}%
						\else
						 \ifx\hyper@ndefined#2#3\hyper@ndefined
						  \hypert@ks{#1}%
						 \else
						  \def\hypert@mp{#1}%
						  \hypert@ks\expandafter\expandafter\expandafter
						  {\expandafter\hypert@mp\hyperh@sh#2.#3}%
						 \fi
						\fi
						\expandafter\hyperlink\expandafter{\the\hypert@ks}}%
	\def\hyperdef#1#2#3{{\global\escapechar=`\\\relax
						 \edef\hypert@mp{\hyperstr@pquote"#2.#3"\hyper@nd}%
						 \expandafter\ifx\csname hyperd@\meaning\hypert@mp
						 \endcsname
						 \relax
						 \expandafter\gdef\csname hyperd@\meaning\hypert@mp
						 \endcsname{}%
						 \gdef#1{{}{\hyperstr@pquote"#2"\hyper@nd}%
								   {\hyperstr@pquote"#3"\hyper@nd}}%
						 \immediate\write\hyperf@le{\def\noexpand#1{#1}}%
						 \xdef\hypert@mp{\global\let\noexpand\hypert@mp=\relax
										 \noexpand\hypertarget{\hypert@mp}}%
						 \global\hypert@ks={\hypert@mp}%
						 \else
						 \message\expandafter{'\hypert@mp' duplicate}%
						 \global\let\hypert@mp=\relax
						 \global\hypert@ks={\hyperdef{#1}{#2}{#3@}}%
						 \fi}\the\hypert@ks}%

	\def\hyper@nique#1#2#3#4{\global\escapechar=`\\\relax
						 \edef\hypert@mp{\hyperstr@pquote"#2.#3"\hyper@nd}%
						 \expandafter\ifx\csname hyperd@\meaning\hypert@mp
						 \endcsname
						 \relax
						 \gdef#1{{}{\hyperstr@pquote"#2"\hyper@nd}%
								   {\hyperstr@pquote"#3"\hyper@nd}}%
						 \global\let\hypert@mp=\relax
						 #4%
						 \else
						 \global\let\hypert@mp=\relax
						 \hyper@nique{#1}{#2}{#3@}{#4}%
						 \fi
						 }%

	\let\hyper@@@@=\relax
	\def\hyper@@{\let\hyper@@@=\relax}%
	\hyper@@
	\def\hyper@{\relax\let\hyper@@@\noexpand\hyper@\noexpand}%
	\def\hyperpr@ref{\hyper@@\hyperref}
	\def\hyperpr@def{\hyper@@\hyperdef}

	\let\href\hyperlink
	
	%
	%
	\hypers@fe

\definecolor {linkColor}{rgb}{0.22,0.28,0.84} \def \Htarget #1#2{\hypertarget {#1}{#2}} \def \Hyperref #1#2#3{\hyperref
{}{#1}{#2}{\color {linkColor}#3}} \def \Href #1#2{\href {#1}{\color {linkColor}#2}} \begingroup \catcode `\@ =0 \catcode
`\\=11 @global@def@textbackslash{\} @endgroup

\begingroup
\obeylines
\leftskip = 0pt plus 2fil \rightskip = 0pt plus 1fil
\null \nocirc

							  \vskip 2.8cm
					  \setfont {cmr10 scaled 2500}
							 \parskip 16pt
					   Partial Dynamical Systems
							Fell Bundles and
							  Applications

\vfill \eject \blank \null \nocirc

							  \vskip 2.8cm
					   Partial Dynamical Systems
							Fell Bundles and
							  Applications

							   \vskip 1cm
					  \setfont {cmr12 scaled 1200}
							  \parskip 6pt
								Ruy Exel
				 Universidade Federal de Santa Catarina
         Partially supported by \Href {http://www.cnpq.br}{CNPq}
			\Href {mailto:exel@mtm.ufsc.br}{exel@mtm.ufsc.br}

								 \vfill
							  --- 2014 ---

\eject \endgroup \nocirc \centerline {INDEX} \vskip 3truecm \ChapterEntry {Introduction} {5} \vskip 1cm \hbox to\hsize
{\hss \titlefont I -- PARTIAL ACTIONS\hss } \vskip 0.5cm \ChapterEntry {Partial actions} {9} \ChapterEntry {Restriction
and globalization} {14} \ChapterEntry {Inverse semigroups} {18} \ChapterEntry {Topological partial dynamical systems}
{23} \ChapterEntry {Algebraic partial dynamical systems} {29} \ChapterEntry {Multipliers} {39} \ChapterEntry {Crossed
products} {42} \ChapterEntry {Partial group representations} {50} \ChapterEntry {Partial group algebras} {60}
\ChapterEntry {C*-algebraic partial dynamical systems} {67} \ChapterEntry {Partial isometries} {72} \ChapterEntry
{Covariant representations of C*-algebraic dynamical systems} {92} \ChapterEntry {Partial representations subject to
relations} {100} \ChapterEntry {Hilbert modules and Morita-Rieffel-equi\-va\-lence} {112} \vskip 1cm \hbox to\hsize
{\hss \titlefont II -- FELL BUNDLES\hss } \vskip 0.5cm \ChapterEntry {Fell bundles} {123} \ChapterEntry {Reduced
cross-sectional algebras} {137} \ChapterEntry {Fell's absorption principle} {147} \ChapterEntry {Graded C*-algebras}
{153} \ChapterEntry {Amenability for Fell bundles} {159} \ChapterEntry {Functoriality for Fell bundles} {169}
\ChapterEntry {Functoriality for partial actions} {188} \ChapterEntry {Ideals in graded algebras} {195} \ChapterEntry
{Pre-Fell-bundles} {202} \ChapterEntry {Tensor products of Fell bundles} {207} \ChapterEntry {Smash product} {219}
\ChapterEntry {Stable Fell bundles as partial crossed products} {225} \ChapterEntry {Globalization in the C*-context}
{234} \ChapterEntry {Topologically free partial actions} {247} \vskip 1cm \hbox to\hsize {\hss \titlefont III --
APPLICATIONS\hss } \vskip 0.5cm \ChapterEntry {Dilating partial representations} {257} \ChapterEntry {Semigroups of
isometries} {261} \ChapterEntry {Quasi-lattice ordered groups} {270} \ChapterEntry {C*-algebras generated by
semi\-groups of isometries} {282} \ChapterEntry {Wiener-Hopf C*-algebras} {287} \ChapterEntry {The Toeplitz C*-algebra
of a graph} {300} \ChapterEntry {Path spaces} {310} \ChapterEntry {Graph C*-algebras} {325} \OtherIndexEntry
{References} {342} \OtherIndexEntry {Subject Index} {347}

\newpage

\blank \ifodd \pageno \else \blank \fi

\begingroup

\chapter Introduction

The concept of a partial dynamical system has been an essential part of Mathematics since at least the late 1800's when,
thanks to the work of Cauchy, Lindel\umlaut of, Lipschitz and Picard, we know that given a Lipschitz vector field $X$ on
an open subset $U\subseteq {\bf R}^n$, the initial value problem $$ f'(t) = X\big (f(t)\big ), \quad f(0) = x_0, $$
admits a unique solution for every $x_0$ in $U$, defined on some open interval about zero.

Assuming that we extend the above solution $f$ to the maximal possible interval, and if we write $\phi _t(x_0)$ for
$f(t)$, then each $\phi _t$ is a diffeomorphism between open subsets of $U$.  Moreover, if $x$ is in the domain of $\phi
_t$, and if $\phi _t(x)$ is in the domain of $\phi _s$, it is easy to see that $x$ lies in the domain of $\phi _{s+t}$,
and that $$ \phi _{s+t}(x) = \phi _s\big (\phi _t(x)\big ).  $$ This is to say that, defining the composition $\phi
_s\circ \phi _t$ on the largest possible domain where it makes sense, one has that $$ \phi _s\circ \phi _t \subseteq
\phi _{s+t}, $$ meaning that $\phi _{s+t}$ extends $\phi _s\circ \phi _t$.  This extension property is the central piece
in defining the notion of a \"{partial dynamical system}, the main object of study in the present book.

Since dynamical systems permeate virtually all of mankind's most important scientific advances, a wide variety of
methods have been used in their study.  Here we shall adopt an algebraic point of view to study partial dynamical
systems, occasionally veering towards a functional analytic perspective.

According to this approach we will extend our reach in order to encompass partial actions on several categories, notably
sets, topological spaces, algebras and C*-algebras.

One of our main goals is to study graded C*-algebras from the point of view of partial actions.  The fundamental
connection between these concepts is established via the notion of \"{crossed product} (known to algebraists as
\"{skew-group algebra}) in the sense that, given a partial action of a group $G$ on an algebra $A$, the crossed product
of $A$ by $G$ is a graded algebra.  While there are many graded algebras which cannot be built out of a partial action,
as above, the number of those who can is surprisingly large.  Firstly, there is a vast quantity of graded algebras
which, when looked at with the appropriate bias, simply happen to be a partial crossed product.  Secondly, and more
importantly, we will see that any given graded algebra satisfying quite general hypotheses is necessarily a partial
crossed product!

Once a graded algebra is described as a partial crossed product, we offer various tools to study it, but we also
dedicate a large part of our attention to the study of graded algebras per se, mainly through a very clever device
introduced by J.~M.~G.~Fell under the name of \"{C*-algebraic bundles}, but which is now more commonly known as \"{Fell
bundles}.  A Fell bundle may be seen essentially as a graded algebra which has been disassembled in such a way that we
are left only with the scattered resulting parts.

Our study of Fell bundles consists of two essentially disjoint disciplines.  On the one hand we study its internal
structure and, on the other, we discuss the various ways in which a Fell bundle may be re-assembled to form a
C*-algebra.  The main structural result we present is that every separable Fell bundle with stable unit fiber algebra
must necessarily arise as the semi-direct product bundle for a partial action of the base group on its unit fiber
algebra.  The study of reassembly, on the other hand, is done via the notions of cross-sectional algebras and
amenability.

A number of applications are presented to the study of C*-algebras, notably C*-algebras generated by semigroups of
isometries, and the now standard class of graph C*-algebras.

Although not discussed here, the reader may find several other situations where well known C*-algebras may be described
as partial crossed products.  Among these we mention: \smallskip \item {$\bullet $} Bunce-Deddens algebras \ref
{Exel/1994b}, \item {$\bullet $} AF-algebras \ref {Exel/1995}, \item {$\bullet $} the Bost-Connes algebra \ref
{BostConnes/1995}, \item {$\bullet $} Exel-Laca algebras \ref {ExelLaca/1999}, \item {$\bullet $} C*-algebras associated
to right-angled Artin groups \ref {CrispLaca/2007}, \item {$\bullet $} Hecke algebras \ref {Exel/2008a}, \item {$\bullet
$} algebras associated with integral domains \ref {BoavaExel/2010}, and \item {$\bullet $} algebras associated to
dynamical systems of type $(m,n)$ \ref {AraExelKatsura/2011}.

Besides, there are numerous other developments involving partial actions whose absence in this book should be
acknowledged.  First and foremost we should mention that we have chosen to restrict ourselves to discrete groups
(i.e. groups without any topology), completely avoiding partial actions of topological groups, even though the latter is
a well studied theory.  See, for example, \ref {Exel/1997a} and \ref {Abadie/2003}.

Although the computation of the K-theory groups of partial crossed product algebras is one of the main focus of the
first two papers on the subject, namely \ref {Exel/1994a} and \ref {McClanahan/1995}, there is no mention here of these
developments.

Twisted partial actions and projective partial representations are also absent, even though they form an important part
of the theory.  See, e.g.~\null \ref {Exel/1997a}, \ref {Sieben/1998}, \ref {DokuchaevExelSimon/2008}, \ref
{BagioLazzarinPaques/2010}, \ref {DokuchaevExelSimon/2010}, \ref {PaquesSantAna/2010}, \ref {DokuchaevNovikov/2010},
\ref {BussExel/2011}, \ref {BussExel/2012} and \ref {DokuchaevNovikov/2012}.

Even though we briefly discuss the relationship between partial actions and inverse semigroups, many interesting
developments have been left out, such as \ref {Sieben/1997}, \ref {Exel/1998}, \ref {Sieben/1998}, \ref {Exel/2008b},
\ref {ExelVieira/2010}, \ref {BussExel/2011}, \ref {BussExel/2012} and \ref {BussExelMeyer/2012}.  The absence of any
mention of the close relationship between partial actions and groupoids \ref {Abadie/2004}, must also be pointed out.

The study of KMS states for gauge actions on partial crossed products studied in \ref {ExelLaca/2003} is also missing
here.

The author gratefully acknowledges financial support from CNPq (Conselho Nacional de Desenvolvimento Cient\'{\i }fico e
Tecnol\'ogico -- Brazil).

\partpage {I}{PARTIAL ACTIONS}

\chapter Partial actions

The notion of a group action applies to virtually every category in Mathematics, the most basic being the category of
sets.  Correspondingly we shall start our development by focusing on partial actions of groups on sets.

\medskip \fix Throughout this chapter we will therefore fix a group $G$, with unit denoted 1, and a set $X$.

\definition \label DefPA A \subj {partial action} of\/ $G$ on $X$ is a pair $$ \Th = \big (\{\D g\}_{g\in G}, \{\th
g\}_{g\in G}\big ) $$ consisting of a collection $\{\D g\}_{g\in G}$ of subsets of $X$, and a collection $\{\th
g\}_{g\in G}$ of maps, $$ \th g: \Di g \to \D g, $$ such that \izitem \zitem \zitemBismark PAIdentity $\D 1= X$, and
$\th 1$ is the identity map, \zitem \zitemBismark PAExtend $\th g\circ \th h \subseteq \th {gh}$, for all $g$ and $h$ in
$G$.  \medskip \noindent By a \subj {partial dynamical system} we shall mean a quadruple $$ \big (X,\ G,\ \{\D g\}_{g\in
G},\ \{\th g\}_{g\in G}\big ) $$ where $X$ is a set, $G$ is a group, and $\big (\{\D g\}_{g\in G}, \{\th g\}_{g\in
G}\big )$ is a partial action of $G$ on $X$.  In case every $\D g=X$, we will say that $\Th $ is a \subj {global
action}, or that we have a \subj {global dynamical system}.

\bigskip We should observe that the composition ``$\th g\circ \th h$'' appearing in \cite {PAExtend} is not defined in
the traditional way since the image of $\th h$ is not necessarily contained in the domain of $\th g$.  Instead this
composition is meant to refer to the map whose domain is the set of all elements $x$ in $X$ for which the expression $$
\th g\big (\th h(x)\big ) $$ makes sense.  For this, $x$ must be in $\Di h$ (the domain of $\th h$), and $\th h(x)$ must
be in $\Di g$ (the domain of $\th g$).  In other words, the domain of $\th g\circ \th h$ is the set $$ \{x\in \Di h: \th
h(x)\in \Di g\} = \th h\inv (\Di g) = \th h\inv (\D h \cap \Di g).  \equationmark DomainCompos $$ For each such $x$, we
of course define $(\th g\circ \th h)(x)=\th g\big (\th h(x)\big )$.

\pix {\setcoordinatesystem units <0.0061truecm, 0.0061truecm> point at 3000 0 \setplotarea x from -1000 to 1000, y from
-500 to 500

\put {$\Di h$} <-15pt,-5pt> at -800 0 \ellipticalarc axes ratio 4:3 360 degrees from -1000 0 center at -800 0

\put {$\th h\inv (\D h \cap \Di g)$} at -850 300 \setdots <1.5pt> \arrow <0.15cm> [0.25,0.75] from -850 250 to -750 100
\setsolid

\setdashes <1.5pt> \ellipticalarc axes ratio 4:3 105 degrees from -800 150 center at -600 100 \setsolid

\put {$\D h$} <-15pt,-5pt> at -100 000 \ellipticalarc axes ratio 4:3 360 degrees from -300 0 center at -100 0

\put {$\Di g$} <15pt,5pt> at 100 100 \ellipticalarc axes ratio 4:3 360 degrees from 300 100 center at 100 100

\put {$\D g$} <15pt,5pt> at 800 100 \ellipticalarc axes ratio 4:3 360 degrees from 1000 100 center at 800 100

\setdashes <1.5pt> \ellipticalarc axes ratio 4:3 105 degrees from 800 -50 center at 600 0 \setsolid

\put {$ \D h \cap \Di g$} <0pt,0pt> at -150 300 \setdots <1.5pt> \arrow <0.15cm> [0.25,0.75] from -150 250 to -50 100
\setsolid

\setquadratic \put {$\th h$} at -450 250 \plot -650 150 -450 200 -250 150 / \arrow <0.15cm> [0.25,0.75] from -250 150 to
-230 140

\put {$\th g$} at 450 350 \plot 650 250 450 300 250 250 / \arrow <0.15cm> [0.25,0.75] from 650 250 to 660 245

\put {$\th g\circ \th h$} at 0 -350 \plot -670 30 0 -290 700 40 / \arrow <0.15cm> [0.25,0.75] from 700 40 to 710 50 }
\diagram {Composing partially defined functions.}

Still referring to \cite {PAExtend}, the symbol ``$\subseteq $'' appearing there is meant to express the fact that the
function on the right-hand-side is an extension\fn {If one defines a function as a \"{set of ordered pairs}, in the
usual technical way, then the symbol ``$\subseteq $'' should indeed be interpreted simply as set inclusion.}  of the one
on the left-hand-side.  In other words \cite {PAExtend} requires that $\th {gh}$ be an extension of $\th g\circ \th h$.

It is possible to rephrase the definition of partial actions just mentioning the collection $\{\th g\}_{g\in G}$ of
partially defined\fn {By a partially defined map on $X$ we mean any map between two subsets of $X$.}  maps on $X$,
without emphasizing the collection of sets $\{\D g\}_{g\in G}$.  In this case we could denote a posteriori the domain of
$\thi g$ by $\D g$, and an axiom should be added to require that the range of $\th g$ be contained in $\D g$.

\state Proposition \label Inverses Given a partial action $\Th $ of\/ $G$ on $X$, as above, one has that each $\th g$ is
a bijection from $\Di g$ onto $\D g$ and, moreover, $ \thi g = \th g\inv .  $

\Proof By \cite {PAExtend}, we have that $\thi g \circ \th g$ is a restriction of $\th 1$, which is the identity map by
\cite {PAIdentity}.  Thus $\thi g \circ \th g$ is the identity on its domain, which is clearly $\Di g$.  Similarly $\th
g\circ \thi g $ is the identity on $\D g$. This concludes the proof.  \endProof

The following provides an equivalent definition of partial actions.

\state Proposition \label SecondPADef Let $\{\D g\}_{g\in G}$ be a collection of subsets of $X$, and let $\{\th
g\}_{g\in G}$ be a collection of maps $$ \th g: \Di g \to \D g.  $$ Then $\big (\{\D g\}_{g\in G}, \{\th g\}_{g\in
G}\big )$ is a partial action of\/ $G$ on $X$ if and only if, in addition to \cite {PAIdentity}, for all $g$ and $h$ in
$G$, one has that: \izitem \zitem \zitemBismark PAIntersecContains $\th g\big (\Di g \cap \D h\big ) \subseteq \D {gh}$,
\zitem \zitemBismark PACompos $\th g\big (\th h(x)\big ) = \th {gh}(x)$, for all $x\in \Di h \cap \Di {(gh)}$.

\Proof Before we begin, we should notice that \cite {PAIntersecContainsLoc} justifies \cite {PAComposLoc} in the
following sense: for all $x$ as in \cite {PAComposLoc}, the fact that $x$ is in $\Di h$ tells us that $\th h(x)$ is well
defined, while $x$ being in $\Di {(gh)}$ implies that $\th {gh}(x)$ is also well defined.  In addition notice that $$
\th h(x) \in \th h\big (\Di h \cap \Di {(gh)}\big ) \_\subseteq {PAIntersecContainsLoc} \Di g, $$ whence $\th h(x)$
indeed lies in the domain of $\th g$.

Assuming we have a partial action, we have already seen that the domain of $\th g\circ \th h$ is precisely $\th h\inv
(\D h \cap \Di g)$.  Since this map is extended by $\th {gh}$, whose domain is $\Di {(gh)}$, we deduce that $$ \th h\inv
(\D h \cap \Di g) \subseteq \Di {(gh)}.  $$ With the change of variables $h:= g\inv $, and $g:= h\inv $, and using \cite
{Inverses}, we obtain \cite {PAIntersecContainsLoc}.

Given $x \in \Di h \cap \Di {(gh)}$, notice that by \cite {PAIntersecContainsLoc}, we have $\th h(x) \in \D h \cap \Di
g$, whence $x$ lies in the domain of $\th g\circ \th h$, and then \cite {PAComposLoc} follows from \cite {PAExtend}.

Conversely, assuming \cite {PAIdentity}, \cite {PAIntersecContainsLoc} and \cite {PAComposLoc}, let us first prove that
$$ \thi g = \th g\inv \for g\in G.  \equationmark AnInverseFact $$ In fact, with $h=g\inv $, point \cite {PAComposLoc}
states that $$ \th g\big (\thi g(x)\big )=x \for x\in \D g, $$ which says that $\thi g$ is a right inverse for $\th g$.
Replacing $g$ by $g\inv $, we see that $\thi g$ is also a left inverse for $\th g$, thus proving \cite {AnInverseFact}.

Let us now prove \cite {PAExtend}.  For this, let $x$ be in the domain of $\th g\circ \th h$, which, as already seen,
means that $$ x\in \th h\inv (\D h \cap \Di g) \={AnInverseFact} \thi h(\D h \cap \Di g) \_\subseteq
{PAIntersecContainsLoc} \D {h\inv g\inv }, $$ thus showing that the domain of $\th g\circ \th h$ is contained in the
domain of $\th {gh}$.  Since $x$ is evidently also in the domain of $\th h$, we conclude that $$ x \in \Di h \cap \D
{h\inv g\inv }, $$ and then \cite {PAComposLoc} implies that $ \th g\big (\th h(x)\big )= \th {gh}(x), $ which means
that $ \th {gh} $ indeed extends $ \th g\circ \th h, $ as desired.  \endProof

Condition \cite {PAIntersecContains} may be slightly improved, as follows:

\state Proposition \label NewPAIntersect Let $ \Th = \big (\{\D g\}_{g\in G}, \{\th g\}_{g\in G}\big ) $ be a partial
action of\/ $G$ on $X$. Then $$ \th g\big (\Di g \cap \D h\big ) = \D g \cap \D {gh} \for g,h\in G.  $$

\Proof Combining \cite {PAIntersecContains} with the fact that the range of $\th g$ is contained in $\D g$, we have that
$$ \th g\big (\Di g \cap \D h\big ) \subseteq \D g \cap \D {gh}.  \equationmark OneInclusion $$ Applying $\thi g$ to
both sides of the above inclusion, and using \cite {Inverses}, we then deduce that $$ \Di g \cap \D h \subseteq \thi
g(\D {gh} \cap \D g).  $$ Replacing $g$ by $g\inv $, and $h$ by $gh$, we then deduce the reverse inclusion in \cite
{OneInclusion}, hence equality.  \endProof

The natural notion of equivalence for partial actions is part of our next:

\definition \label DefineEquivar Let $G$ be a group, and suppose that, for each $i=1,2$, we are given a partial action $
\Th ^i = \big (\{\D g^i\}_{g\in G}, \{\th g^i\}_{g\in G}\big ) $ of $G$ on a set $X^i$.  A map $$ \phi :X^1\to X^2 $$
will be said to be \subjex {$G$-equivariant}{$G$-equivariant map} when, for all $g$ in $G$, one has that \izitem \zitem
$\phi (\D g^1) \subseteq \D g^2$, and \zitem $\phi \big (\th g^1(x)\big ) = \th g^2\big (\phi (x)\big )$, for all $x$ in
$\Di g^1$.  \medskip \noindent If moreover $\phi $ is bijective and $\phi \inv $ is also $G$-equivariant, we will say
that $\phi $ is an \subj {equivalence of partial actions}.  If such an equivalence exists, we will say that $\Th ^1$ is
equivalent to $\Th ^2$.

Observe that a $G$-equivariant bijective map $\phi $ from $X^1$ to $X^2$ may fail to be an equivalence since property
\cite {DefineEquivar/i} may fail for $\phi \inv $.

On the other hand, it is easy to see that a necessary and sufficient condition for a $G$-equivariant bijective map to be
an equivalence is that equality hold in \cite {DefineEquivar/i}.

\definition \label DefGraph Let $ \Th = \big (\{\D g\}_{g\in G}, \{\th g\}_{g\in G}\big ) $ be a partial action of $G$
on $X$. The \subjex {graph}{graph of a partial action} of $\Th $ is defined to be the set $$ \Graph (\Th ) = \{(y,g,x)
\in X\times G\times X : x\in \Di g,\ \th g(x) = y\}.  $$

The graph of a partial action evidently encodes all of the information contained in the partial action itself. It is
therefore a very useful device in the study of partial actions and we will re-encounter it often in the sequel.
Meanwhile let us remark that it has a natural groupoid \ref {Abadie/2004} structure, where two elements $(w,h,z)$ and
$(y,g,x)$ in $\Graph (\Th )$ may be multiplied if and only if $z=y$, in which case we put $$ (w,h,z)(y,g,x) = (w,hg,x),
$$ while $$ (y,g,x)\inv = (x,g\inv ,y).  $$ We leave it for the reader to check that $\Graph (\Th )$ is indeed a
groupoid with these operations.

\medskip The usual notion of invariance for group actions has a counterpart in partial actions, as follows.

\definition \label DefineInvariant Given a partial action $ \Th = \big (\{\D g\}_{g\in G}, \{\th g\}_{g\in G}\big ) $ of
$G$ on $X$, we will say that a given subset $Y\subseteq X$ is \subjex {invariant}{invariant set under a partial action}
under $\Th $, if $$ \th g(Y\cap \Di g)\subseteq Y \for g\in G.  $$

As in the global case we have:

\state Proposition \label RestrToInvar Given a partial action $ \Th = \big (\{\D g\}_{g\in G}, \{\th g\}_{g\in G}\big )
$ of\/ $G$ on $X$, and an invariant subset $X'\subseteq X$, let $$ \D g' = X'\cap \D g \for g\in G, $$ and let $\th g'$
be the restriction of $\th g$ to $\Di g'$.  Then \izitem \zitem $ \Th ' = \big (\{\D g'\}_{g\in G}, \{\th g'\}_{g\in
G}\big ) $ is a partial action of\/ $G$ on $X'$, and \zitem the inclusion $X' \hookrightarrow X$ is a $G$-equivariant
map.

\Proof By hypothesis it is clear that each $\th g'$ maps $\Di g'$ into $\D g'$.  Observing that \cite {PAIdentity} is
evident for $\Th '$, it suffices to check conditions \cite {SecondPADef/i--ii}.  For all $g$ and $h$ in $G$, we have
that $$ \th g'(\Di g'\cap \D h') = \th g(X'\cap \Di g\cap \D h) \$\subseteq \ \th g(X'\cap \Di g) \ \cap \ \th g(\Di
g\cap \D h)\ \subseteq \ X' \cap \D {gh} = \D {gh}', $$ thus verifying \cite {PAIntersecContains}.  We leave the easy
proof of \cite {PACompos} for the reader.  Point (ii) is also easy to check.  \endProof

\nrem Partial actions of the group ${\bf Z}$ on C*-algebras were first considered by the author in \ref {Exel/1994a}.
Soon afterward McClanahan \ref {McClanahan/1995} generalized this notion for arbitrary groups.  Partial actions of
groups on arbitrary sets were introduced in \ref {Exel/1998}.

\chapter Restriction and globalization

One of the easiest ways to produce nontrivial examples of partial actions is via the process of \"{restriction} to not
necessarily invariant subsets, which we would now like to describe.

Suppose we are given a global action $\eta $ of a group $G$ on a set $Y$.  Suppose further that $X$ is a given subset of
$Y$ and we wish to restrict $\eta $ to an action of $G$ on $X$.  Evidently this requires that $X$ be an $\eta
$-invariant subset.

\pix {\setcoordinatesystem units <0.005truecm, 0.005truecm> point at 0 0 \setplotarea x from -1000 to 1000, y from -600
to 600

\put {$Y$} at 0 380 \ellipticalarc axes ratio 4:3 360 degrees from -600 0 center at 0 0

\put {$X$} at -350 -10 \ellipticalarc axes ratio 4:3 360 degrees from -450 0 center at -150 0

\put {$\eta _g(X)$} at 300 -10 \ellipticalarc axes ratio 4:3 360 degrees from 450 0 center at 150 0

\put {$\D g$} at 0 10 }

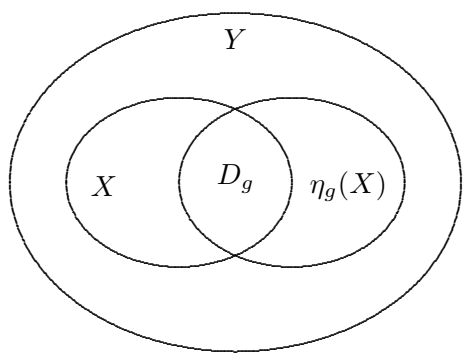
\diagram {Restricting a global action to a partial action.}

However, even when $X$ is not $\eta $-invariant, we may restrict $\eta $ to a \"{partial action} of $G$ on $X$, by
setting $$ \D g= \eta _g(X)\cap X, $$ and, after observing that $\eta _g(\Di g) = \D g$, we may let $$ \th g: \Di g \to
\D g $$ be the restriction of $\eta _g$ to $\Di g$, for each $g$ in $G$.  The easy verification that this indeed leads
to a partial action is left to the reader.

\definition \label DefineRestriction The partial action $\Th $ of $G$ on $X$ defined above is called the \subjex
{restriction}{restriction of a global action} of the global action $\eta $ to $X$.

It is readily seen that, under the above conditions, the inclusion of $X$ into $Y$ is a $G$-equivariant map.

\medskip

Let us now discuss the relationship between the graph of a global action and that of its restriction to a subset.

\state Proposition \label GraphRestr Let $\eta $ be a global action of a group $G$ on a set $Y$ and let $\Th $ be its
restriction to a subset $X\subseteq Y$.  Then $$ \Graph (\Th ) = \Graph (\eta )\cap (X\times G\times X).  $$

\Proof Left to the reader.  \endProof

Taking a different point of view, one may start with a partial action $\Th $ of $G$ on $X$ and wonder about the
existence of a global action $\eta $ on some set $Y$ containing $X$, in such a way that $\Th $ is the restriction of
$\eta $.  Assuming that $\eta $ does exist, the \subj {orbit} of $X$ in $Y$, namely $$ \Orb (X) := \bigcup _{g\in G}
\eta _g(X) $$ will evidently be $\eta $-invariant and the restriction of $\eta $ to $\Orb (X)$ will be another global
action which, further restricted to $X$, will again produce the original partial action $\Th $.

\definition \label DefineGlobalization Let $\eta $ be a global action of $G$ on a set $Y$, and let $\Th $ be the partial
action obtained by restricting $\eta $ to a subset $X\subseteq Y$.  If the orbit of $X$ coincides with $Y$, we will say
that $\eta $ is a \subj {globalization}, or an \subj {enveloping action} for $\Th $.

\state Theorem \label GlobaSets Every partial action admits a globalization, which is unique in the following sense: if
$\Th $ is a partial action of $G$ on $X$, and we are given globalizations $\eta ^i$ acting on sets $Y^i$, for $i=1,2$,
then there exists an equivalence $$ \phi :Y^1\to Y^2, $$ such that $\phi $ coincides with the identity on $X$.

\Proof Given a partial action $\Th = \big (\{\D g\}_{g\in G}, \{\th g\}_{g\in G}\big )$ of $G$ on a set $X$, define an
equivalence relation on $G\times X$ by $$ (g,x) \sim (h,y) \iff x\in \D {g\inv h}, \hbox { and } \th {h\inv g}(x)=y.  $$

That this is indeed a reflexive and symmetric relation is apparent, while the transitivity property may be verified as
follows: assuming that $$ (g,x) \sim (h,y) \sim (k,z), $$ we have that $y\in \D {h\inv k}$, and $\th {k\inv h}(y)=z$.
We then have that $$ x = \th {g\inv h}(y) \in \th {g\inv h}(\D {h\inv g} \cap \D {h\inv k}) = \D {g\inv h} \cap \D
{g\inv k}, $$ whence by \cite {PACompos} $$ \th {k\inv g}(x) = \th {k\inv h}\big (\th {h\inv g}(x)\big ) = \th {k\inv
h}(y) = z, $$ showing that $(g,x) \sim (k,z)$.

Let $\tilde X$ be the quotient of $G\times X$ by this equivalence relation, and let us denote by $[g,x]$ the equivalence
class of each $(g,x)\in G\times X$.  The map $$ \iota : x\in X \mapsto [1,x] \in \tilde X \equationmark MapXToGloba $$
is clearly injective and hence we may identify $X$ with its image $X'$ in $\tilde X$ under $\iota $.  One then verifies
that the expression $$ \tau _g([h,x]) = [gh,x] \for g,h\in G \for x\in X, $$ gives a well defined global action $\tau $
of $G$ on $\tilde X$.

We claim that $$ \tau _g(X') \cap X' = \iota (\D g).  $$ In fact, an element of the form $[1,x]$ lies in $\tau _g(X')$
if and only if $ [1,x] = [g,y], $ for some $y$ in $X$, but this is so iff $x\in \D g$, and $\thi g(x)=y$, proving the
claim.

In order to verify that $\Th $ indeed corresponds to the restriction of $\tau $ to $X$, let $x\in \Di g$, and observe
that $$ \big (1,\th g(x)\big ) \sim (g,x), $$ whence $$ \iota (\th g(x)) = [1,\th g(x)] = [g,x] = \tau _g([1,x]) = \tau
_g(\iota (x)\big ).  $$ This is to say that $\iota \circ \th g = \tau _g\circ \iota $, which means that $\th g = \tau
_g$, on $\Di g$, up to the above identification of $X$ with $X'$.

Observing that $\tau _g([1,x])=[g,x]$, we immediately see that the orbit of $X$ coincides with $\tilde X$, so $\tau $ is
indeed a globalization of $\Th $.

If $\eta $ is another globalization of $\Th $, acting on the set $Y\supseteq X$, define a map $$ \phi :(g,x)\in G\times
X \mapsto \eta _g(x)\in Y, $$ and observe that $$ \phi (g,x) = \phi (h,y) \iff \eta _g(x) = \eta _h(y) \iff x = \eta
_{g\inv h}(y), $$ Under these conditions we have that $ x \in \eta _{g\inv h}(X)\cap X = \D {g\inv h}, $ and $ y = \eta
_{h\inv g}(x) = \th {h\inv g}(x), $ meaning that $(g,x)\sim (h,y)$.

Consequently the map $\phi $, above, factors through the quotient $\tilde X$, providing an injective map $$ \tilde \phi
:\tilde X \to Y.  $$ Being a globalization, $Y$ coincides with the orbit of $X$ under $\eta $, from where one deduces
that $\tilde \phi $ is surjective.

We leave it for the reader to prove that $\tilde \phi $ is $G$-equivariant and coincides with the identity on the
respective copies of $X$ in $\tilde X$ and $Y$.  This proves that $\tilde X$ and $Y$ are equivalent globalizations.
\endProof

The process of restricting a global action $\eta $ to a non-invariant subset, studied above, could be generalized to the
case in which $\eta $, itself, is a partial action.  More precisely, suppose we are given a partial action $$ \eta =
\big (\{Y_g\}_{g\in G}, \{\eta _g\}_{g\in G}\big ) $$ of a group $G$ on a set $Y$, and let $X\subseteq Y$ be a not
necessarily invariant subset.  Defining $$ \D g = Y_g\cap X\,\cap \,\eta _g\big (Y_{g\inv }\cap X\big ) \for g\in G, $$
one may prove that $\eta _g(\Di g)=\D g$.  If we then let $\th g$ be the restriction of $\eta _g$ to $\Di g$, one has
that $$ \Th = \big (\{\D g\}_{g\in G}, \{\th g\}_{g\in G}\big ) $$ is a partial action of $G$ on $X$, with respect to
which the inclusion $X\hookrightarrow Y$ is a $G$-equivariant map.

Contrary to the case of restricting a global action, we have found no interesting applications of the above idea.  Thus,
until this situation changes, we will not bother to discuss it any further.

\nrem The notion of enveloping actions and the proof of Theorem \cite {GlobaSets} is due to Abadie and it first appeared
in his PhD thesis \ref {Abadie/1999}, which led to the article \ref {Abadie/2003}.

\chapter Inverse semigroups

In the study of partial actions, the notion of inverse semigroups plays a predominant role.  So let us briefly introduce
this concept, referring the reader to \ref {Lawson/1998} for proofs and further details.

\definition An \subj {inverse semigroup} is a nonempty set $S$ equipped with a binary associative operation (i.e., $S$
is a semigroup) such that, for every $s$ in $S$, there exists a unique element $s^*$ in $S$, such that $$ ss^*s = s \and
s^*ss^* = s^*.  $$

Given an inverse semigroup $S$, one may prove that the collection of idempotent elements in $S$, namely $$ E(S) = \{s\in
S: s^2=s\}, $$ is a commutative sub-semigroup.  Under the partial order relation defined in $E(S)$ by $$ e\leq f \iff
ef=e \for e,f\in E(S), $$ $E(S)$ becomes a \subj {semilattice}, meaning that for every $e$ and $f$ in $E(S)$, there
exists a largest element which is smaller than $e$ and $f$, namely $ef$.  One therefore often refers to $E(S)$ as the
\subj {idempotent semilattice} of $S$.

There is also a natural partial order relation defined on $S$ itself by $$ s\leq t \iff ts^*s=s \iff ss^*t = s \for
s,t\in S.  \equationmark DefinePOinISG $$ This is compatible with the multiplication operation in the sense that $$
s\leq t,\ s'\leq t' \imply ss'\leq tt'.  \equationmark CompatibleOrder $$

One of the main examples of inverse semigroups is as follows: given a set $X$, two subsets $C,D\subseteq X$, and a
bijective map $$ f:C\to D, $$ we will say that $f$ is a \subj {partial symmetry} of $X$.  The set $$ \IX = \{f: f \hbox
{ is a partial symmetry of } X\} \equationmark DefineParSym $$ may be turned into a semigroup by equipping it with the
operation of composition, where, as before, the composition of two partially defined maps is defined on the largest
domain where it makes sense.

It may be easily proven that $\IX $ is an inverse semigroup, where, for every $f$ in $\IX $, one has that $f^*$ is the
inverse of $f$.

The idempotent elements of $\IX $ are just the identity maps defined on subsets of $X$.  The order among idempotents
happens to be the same as the order of inclusion of their domains.

More generally, the natural order among general elements of $\IX $ is the order given by ``extension'', meaning that,
for $f$ and $g$ in $\IX $, one has that $f\leq g$, if and only if $g$ is an extension of $f$, in symbols $$ f\leq g \iff
f\subseteq g \for f,g\in \IX .  $$

The classical Wagner-Preston Theorem \ref {Wagner/1952}, \ref {Preston/1954} asserts that any inverse semigroup is
isomorphic to a *-invariant sub-semigroup of $\IX $, for some $X$.  This may be considered as the version for inverse
semigroups of the well known Cayley Theorem for groups.

Given a partial action $\Th $ of a group $G$ on a set $X$, notice that each $\th g$ is an element of $\IX $.

\state Proposition \label PAasPR Let $G$ be a group, $X$ be a set, and $$ \Th :G\to \IX $$ be a map.  Then $\Th $ is a
partial action of $G$ on $X$ if and only if, for every $g$ and $h$ in $G$, one has that \izitem \zitem $\th 1$ is the
identity map of $X$, \zitem $\thi g = (\th g)^*$, \zitem $\th g \th h \thi h = \th {gh}\thi h$, \zitem $\thi g\th g \th
h = \thi g\th {gh}$.

\Proof Assuming that $\Th $ is a partial action, one immediately checks (i) and (ii).  As for (iii), observe that,
following \cite {DomainCompos}, the domain of $\th {gh}\thi h$ is $$ \th h(\Di h \cap \D {(gh)\inv }) \={NewPAIntersect}
\D h \cap \Di g, $$ which is clearly also the domain of $\th g \th h \thi h$.  For any $x$ in this common domain we have
that $$ \thi h(x) \in \thi h(\D h \cap \Di g) = \Di h \cap \D {h\inv g\inv }, $$ so (iii) follows immediately from \cite
{PACompos}.  Finally, (iv) follows from (iii) and (ii), in view of the fact that the star operation reverses
multiplication.

Conversely, suppose that $\Th $ satisfies (i--iv) and denote the range of $\th g$ by $\D g$.  Since $\thi g$ is the
inverse of $\th g$ by (i), we have that the range of $\thi g$ is the domain of $\th g$, meaning that $\th g$ is a map $$
\th g: \Di g \to \D g.  $$

Since \cite {PAIdentity} is granted, let us verify \cite {PAExtend}.  For this, let us suppose we are given $g$ and $h$
in $G$, as well as an element $x$ in the domain of $\th g\th h$.  Then the element $y:= \th h(x)$ clearly lies in the
domain of $\th g \th h \thi h$, and hence (iii) implies that it is also in the domain of $\th {gh} \thi h$, which is to
say that $x$ is in the domain of $\th {gh}$.  Moreover $$ \th g\big (\th h(x)\big ) = \th g\big (\th h(\thi h(y))\big )
= \th {gh}\big (\thi h(y)\big ) = \th {gh}(x), $$ proving that $\th g\th h \subseteq \th {gh}$.  \endProof

The language of inverse semigroups is especially well suited to the introduction of our next example.  In order to
describe it, let $X$ be a set and let $\{f_\lambda \}_{\lambda \in \Lambda }$ be any collection of partial symmetries of
$X$.

Letting $\F = \F (\Lambda )$ be the free group on the index set $\Lambda $, our plan is to construct a partial action of
$\F $ on $X$ in the form of a map $$ \Th : \F \to \IX .  $$

As a first step, let us define $$ \th \lambda := f_\lambda \and \thi \lambda := f_\lambda \inv \for \lambda \in \Lambda
.  $$ Given any element $g\in \F $, write $$ g = x_1x_2\ldots x_n, $$ in \subj {reduced form}, meaning that each $
x_j\in \Lambda \cup \Lambda \inv , $ and $x_{j+1}\neq x_j\inv $.  It is well known that a unique such decomposition of
$g$ exists.  We then put $$ \th g := \th {x_1}\th {x_2}\ldots \th {x_n}, \equationmark FreeGpAction $$ with the
convention that, if $g=1$, then its reduced form is ``empty'', and $\th 1$ is the identity function on $X$.

\state Proposition The map $\Th $ defined in \cite {FreeGpAction}, above, is a partial action of\/ $\F $ on $X$.

\Proof We leave it for the reader to prove \cite {PAIdentity}.  In order to verify \cite {PAExtend}, pick $g$ and $h$ in
$\F $, and write $$ g = x_nx_{n-1}\ldots x_2x_1 \and h = y_1y_2\ldots y_{m-1}y_m, $$ in reduced form, as above (for
reasons which will soon become clear, we have chosen to reverse the indices in the reduced form of $g$).

Let $p$ be the number of cancellations occurring when performing the multiplication $gh$, meaning that $$ x_i = y_i\inv
\for i=1,\ldots ,p, $$ and that $p$ is maximal with this property.  Define $g'$, $h'$, and $k$ in $\F $, as follows $$ g
= \overbrace {x_n\ldots x_{p+1}}^{g'}\ \overbrace {x_p\ldots x_1}^k \and h = \overbrace {y_1\ldots y_p}^{k\inv }\
\overbrace {y_{p+1}\ldots y_m}^{h'}, $$ so that $$ gh = g'h' = x_n\ldots x_{p+1}\ y_{p+1}\ldots y_m $$ is in reduced
form.  Denoting the identity map on $X$ by $id_X$, we then clearly have that $$ \th k\th k\inv \subseteq id_X, $$ so $$
\th g \th h = \th {g'}\th k \th k\inv \th h' \_\subseteq {CompatibleOrder} \th {g'} \, id_X \, \th {h'} = \th {g'}\th
{h'} = \th {g'h'} = \th {gh}, \equationmark MotivateSemiSat $$ concluding the proof.  \endProof

Observe that in \cite {MotivateSemiSat} we were allowed to use that $$ \th {g'}\th {h'} = \th {g'h'} $$ because the
juxtaposition of the reduced forms of $g'$ and $h'$ turned out to be the reduced form of $g'h'$.  Equivalently, $|g'h'|
= |g'|+|h'|$, where $|\cdot |$ refers to \subj {word length}.

\definition A \subj {length function} on a group $G$ is a function $ \ell : G \to \R _+ $ such that \Zitem $\ell (1)=0$,
and \zitem $\ell (gh)\leq \ell (g)+ \ell (h)$, for all $g$ and $h$ in $G$.

Evidently $|\cdot |$ is a length function for the free group.

\definition \label DefineSemiSatPAct Let $G$ be a group equipped with a length function $\ell $.  A partial action $\Th
$ of $G$ is said to be \subjex {semi-saturated}{semi-saturated partial action} (with respect to the given length
function $\ell $) if $$ \ell (gh)=\ell (g)+\ell (h) \quad \Longrightarrow \quad \th g\th h = \th {gh} \for g,h\in G.  $$

In the free group, observe that the condition that $|gh| = |g|+|h|$ means that the juxtaposition of the reduced forms of
$g$ and $h$ is precisely the reduced form of $gh$.  So the partial action $\Th $ defined in \cite {FreeGpAction} is
easily seen to be semi-saturated.

Summarizing we have:

\state Proposition \label ActionsOfFreeGroups Let $X$ be a set and let $\{f_\lambda \}_{\lambda \in \Lambda }$ be any
collection of partial symmetries of $X$.  Then there exists a unique semi-saturated partial action $\Th $ of\/ $\F
(\Lambda )$ on $X$ such that $$ \th \lambda = f_\lambda \for \lambda \in \Lambda .  $$

\Proof The existence was proved above and we leave the proof of uniqueness as an easy exercise. \endProof

\nrem Semi-saturated partial actions are related to Quigg and Raeburn's notion of \"{multiplicativity}, as defined in
\ref {Definition 5.1/QuiggRaeburn/1997}.  The relationship between partial actions and inverse semigroups is further
discussed in \ref {Sieben/1997}, \ref {Exel/1998}, \ref {Exel/2008b}, \ref {ExelVieira/2010}, \ref {BussExel/2011}, \ref
{BussExel/2012} and \ref {BussExelMeyer/2012}.

\chapter Topological partial dynamical systems

We shall now be concerned with the topological aspects of partial actions, so, \fix throughout this chapter, we will fix
a group $G$ and topological space $X$.

\definition \label DefTopoPa A \subj {topological partial action} of the group $G$ on the topological space $X$ is a
partial action $\Th = \big (\{\D g\}_{g\in G}, \{\th g\}_{g\in G}\big )$ on the underlying set $X$, such that each $\D
g$ is \"{open} in $X$, and each $\th g$ is a \"{homeomorphism}.  By a \subj {topological partial dynamical system} we
shall mean a partial dynamical system $$ \big (X,\ G,\ \{\D g\}_{g\in G},\ \{\th g\}_{g\in G}\big ) $$ where $X$ is a
topological space and $\big (\{\D g\}_{g\in G}, \{\th g\}_{g\in G}\big )$ is a topological partial action of $G$ on $X$.

When it is understood that we are working in the category of topological spaces and there is no chance for confusion we
will drop the adjective \"{topological} and simply say \"{partial action} or \"{partial dynamical system}.

The notion of topological partial actions may also take into account topological groups (see \ref {Exel/1997a} and \ref
{Abadie/2003}) but, for the sake of simplicity, we will only be concerned with discrete groups in this book.

Recall from \cite {DefineParSym} that $\IX $ denotes the inverse semigroup formed by all partial symmetries of $X$.

\definition We will say that a partial symmetry $f\in \IX $ is a \subj {partial homeomorphism} of $X$, if the domain and
range of $f$ are open subsets of $X$, and $f$ is a homeomorphism from its domain to its range.  We will denote by $\IXt
$ the collection of all partial homeomorphisms of $X$.  It is evident that $\IXt $ is an inverse sub-semigroup of $\IX
$.

As an immediate consequence of \cite {PAasPR} we have:

\state Proposition \label PAasPRTop Let $G$ be a group, $X$ be a topological space, and $$ \Th :G\to \IXt $$ be a map.
Then $\Th $ is a topological partial action of $G$ on $X$ if and only if conditions (i--iv) of \cite {PAasPR} are
fulfilled.

We will say that two topological partial actions of $G$ are \subjex {topologically equivalent}{topologically equivalent
partial actions} if they are equivalent in the sense of Definition \cite {DefineEquivar} and, in addition, the map $\phi
$ mentioned there is a homeomorphism.

The process of restriction may be applied in the context of topological spaces as follows: suppose we are given a global
topological action $\eta $ of a group $G$ on a space $Y$, as well as an \"{open subset} $X\subseteq Y$.

Let $\Th $ be the restriction of $\eta $ to $X$, as defined in \cite {DefineRestriction} and recall that $$ \D g= \eta
_g(X)\cap X.  $$ Since $\eta _g$ is a homeomorphism we have that $\eta _g(X)$ is open in $Y$, so $\D g$ is an open
subset of $X$.  Also, since the restricted maps $\th g$ are obviously homeomorphisms, we conclude that $\Th $ is a
topological partial action of $G$ on $X$.

\definition \label DefineTopoGlobalization Let $\eta $ be a topological global action of $G$ on a space $Y$, and let
$\Th $ be the topological partial action obtained by restricting $\eta $ to an open subset $X\subseteq Y$.  If the orbit
of $X$ coincides with $Y$, we will say that $\eta $ is a \subj {topological globalization} of $\Th $.

\state Proposition \label TopoGlob Every topological partial action admits a topological globalization, unique up to
topological equivalence.

\Proof Given a topological partial action $\Th $ of the group $G$ on the space $X$, let $\tau $ and $\tilde X$ be as in
the proof of \cite {GlobaSets}.

Viewing $G$ as a discrete space, consider the product topology on $G\times X$ and let us equip $\tilde X$ with the
quotient topology.  We then claim that $\tau $ is a topological globalization of $\Th $.  In other words we must prove
that each $\tau _g$ is a homeomorphism, the map $\iota $ mentioned in \cite {MapXToGloba} is a homeomorphism onto its
image, and $\iota (X)$ is open.

The first assertion follows easily from the fact that, for $g\in G$, the map $$ (h,x)\in G\times X \mapsto (gh,x)\in
G\times X $$ is a homeomorphism respecting the equivalence relation ``$\sim $''.

Observe that $\iota (x) = \pi (1,x)$, for every $x$ in $X$, where $$ \pi :G\times X \to G\times X/{\sim } $$ is the
quotient map, so $\iota $ is clearly continuous.  In order to prove that $\iota $ is a homeomorphism onto its image it
is therefore enough to show that it is an open map.  With this purpose in mind let $U\subseteq X$ be an open set, and
let us prove that $\iota (U)$ is open in $\tilde X$ (in fact, it would be enough to prove that $\iota (U)$ is open in
$\iota (X)$).

By definition of the quotient topology, $\iota (U)$ is open in $\tilde X$ if and only if $\pi \inv \big (\iota (U)\big
)$ is open in $G\times X$.  Moreover observe that $$ (g,x) \in \pi \inv \big (\iota (U)\big ) $$ if and only if $ (g,x)
\sim (1,y), $ for some $y\in U$, but this is clearly equivalent to saying that $x\in \Di g$, and $\th g(x)\in U$, or
that $x\in \th g\inv (\D g\cap U)$.  So $$ \pi \inv \big (\iota (U)\big ) = \bigcup _{g\in G} \ \{g\}\times \th g\inv
(\D g\cap U), $$ which is then seen to be open in $G\times X$, whence $\iota (U)$ is open in $\tilde X$, as claimed.
This proves that $\iota $ is a homeomorphism onto its image, and it also implies that $\iota (X)$ is an open subset of
$\tilde X$.  Consequently $\tau $ is a topological globalization of $\Th $, taking care of the existence part of the
statement.

In order to prove uniqueness, let us be given another topological globalization of $\Th $, say $\eta $, acting on the
space $Y$.  In particular $\eta $ is a set-theoretical globalization of $\Th $, so by the uniqueness part of \cite
{GlobaSets}, there exists an equivalence $$ \phi : Y \to \tilde X, $$ coinciding with the identity map on $X$ (or
rather, its canonical copies within $\tilde X$ and $Y$).  In order to conclude the proof it is now enough to prove that
$\phi $ is a homeomorphism.

Given $y\in Y$, use the fact that the orbit of $X$ coincides with $Y$ to find some $g$ in $G$ such that $y\in \eta
_g(X)$.  For any other $y'\in \eta _g(X)$, we then have that $x':=\eta _g\inv (y')\in X$, so $$ \phi (y') = \phi \big
(\eta _g(x')\big ) = \tau _g \big (\phi (x')\big ) = \tau _g (x') = \tau _g \big (\eta _g\inv (y')\big ).  $$ Since the
map $$ y'\in \eta _g(X) \mapsto \tau _g \big (\eta _g\inv (y')\big )\in \tau _g(X) $$ is clearly continuous and
coincides with $\phi $ on $\eta _g(X)$, which is an open neighborhood of $y$, we see that $\phi $ is continuous at $y$.
Similarly one may prove that the inverse of $\phi $ is continuous, and hence that $\phi $ is a homeomorphism, as
claimed.  \endProof

The reader should be warned that there is a catch in the above result: the quotient topology defined on $\tilde X$ may
not be Hausdorff, even if $X$ is Hausdorff.  In fact it is easy to characterize the occurrence of this problem, as seen
in the next:

\state Proposition \label TopoGlobIffClosed A topological partial action of a group $G$ on a Hausdorff space $X$ admits
a Hausdorff globalization if and only if its graph is closed in $X\times G\times X$, where $G$ is given the discrete
topology.

\Proof Assuming that $\eta $ is a globalization of $\Th $ on a Hausdorff space $Y$, observe that $$ \Graph (\eta ) =
\{(y,g,x)\in Y\times G\times Y: y = \eta _g(x)\}.  $$ This should be contrasted with Definition \cite {DefGraph}, as
here there is no need to require $x$ to lie in the domain of $\eta _g$, a globally defined map!  For this reason we see
that $\Graph (\eta )$ is closed in $Y\times G\times Y$ and it then follows from \cite {GraphRestr} that $\Graph (\Th )$
is closed in $X\times G\times X$.

Conversely, assuming that $\Graph (\Th )$ is closed, let $\eta $ be any topological globalization of $\Th $, acting on
the space $Y$ (recall from \cite {TopoGlob} that there is only one such).

We will prove that $Y$ is Hausdorff.  For this, let us be given two distinct points $y_1$ and $y_2$ in $Y$, and write $
y_i = \eta _{g_i}(x_i), $ with $g_i\in G$, and $x_i\in X$.  Noticing that $$ \eta _{g_1}(x_1) = y_1 \neq y_2 = \eta
_{g_2}(x_2), $$ one sees that $\eta _{g_2\inv g_1}(x_1) \neq x_2$, whence $(x_2,g_2\inv g_1,x_1)$ is not in the graph of
$\Th $.  From the hypothesis it follows that there are open subsets $U_1,U_2\subseteq X$, with $x_i\in U_i$, such that
$$ U_2\times \{g_2\inv g_1\}\times U_1 \ \cap \ \Graph (\Th ) = \emptyset .  $$

Since $X$ is open in $Y$, we have that each $U_i$ is open in $Y$ and hence so is $\eta _{g_i}(U_i)$.  Since $x_i\in
U_i$, we have that $$ y_i = \eta _{g_i}(x_i) \in \eta _{g_i}(U_i).  $$

In order to verify Hausdorff's axiom, it now suffices to prove that $\eta _{g_1}(U_1)$ and $\eta _{g_2}(U_2)$ are
disjoint.  Arguing by contradiction, suppose that $$ y\in \eta _{g_1}(U_1) \cap \eta _{g_2}(U_2), $$ so we may write $y
= \eta _{g_1}(z_1) = \eta _{g_2}(z_2)$, with $z_i\in U_i$.  It follows that $\eta _{g_2\inv g_1}(z_1) = z_2$, but since
both $z_1$ and $z_2$ lie in $X$, we indeed have that $$ \th {g_2\inv g_1}(z_1) = z_2, $$ whence $$ (z_2,g_2\inv g_1,z_1)
\in \Graph (\Th ) \ \cap \ U_2\times \{g_2\inv g_1\}\times U_1 = \emptyset , $$ a contradiction.  This proves that $\eta
_{g_1}(U_1)$ and $\eta _{g_2}(U_2)$ are disjoint and hence that $Y$ is Hausdorff.  \endProof

We will see that in many interesting examples, the $\D g$ are closed besides being open.  In this case we have:

\state Proposition \label ClopenGlobalizes A topological partial dynamical system $$ \big (X,\ G,\ \{\D g\}_{g\in G},\
\{\th g\}_{g\in G}\big ), $$ such that each $\D g$ is closed, always admits a Hausdorff globalization.

\Proof If $\{(y_i,g_i,x_i)\}_i$ is a net in $\Graph (\Th )$, converging to $(y,g,x)$, then $g_i$ is eventually constant,
since $G$ has the discrete topology.  So we may assume, without loss of generality that $g_i=g$, for all $i$.  Then
every $x_i$ lie in $\Di g$, and hence $x$ is in $\Di g$, by hypothesis.  Consequently $$ y = \lim _iy_i = \lim _i\th
g(x_i) = \th g(x), $$ so $(y,g,x)$ lies in $\Graph (\Th )$, which is therefore proven to be closed.  The conclusion then
follows from \cite {TopoGlobIffClosed}.  \endProof

Let us now give an important example of a topological partial dynamical system based on Bernoulli's action.  For this,
let $G$ be a group and consider the compact topological space $$ \OG $$ relative to the product topology.

It is well known that $\OG $ is naturally equivalent to $\Part (G)$, the power set of $G$, in such a way that, when a
given $ \omega \in \OG $ is seen as a subset of $G$, the Boolean value of the expression ``$g\in \omega $'' is given by
the coordinate $\omega _g$.  In symbols $$ \omega _g=\bool {g\in \omega }, \equationmark VectorsAndSets $$ where we are
using brackets to denote Boolean value.  In what follows we will often identify $\OG $ and $\Part (G)$ based on this
correspondence.

For each $g\in G$, and each $\omega \in \Part (G)$, let us indicate by $$ g\omega = \{gh: h\in \omega \}, $$ as usual,
and let us consider the mapping $$ \eta _g : \omega \in \Part (G) \mapsto g\omega \in \Part (G).  \equationmark
BernoulliAction $$ It is easy to see that each $\eta _g$ is a homeomorphism and that $\eta $ is a topological global
action of $G$ on $\OG $.

\definition \label DefineGlobalBernouli The action $\eta $ of $G$ on $\OG $, defined above, will be called the \subjex
{global Bernoulli action}{Bernoulli global action} of $G$.

The next concept to be defined is based on the observation that the set $$ \OuG = \big \{\omega \in \OG : 1\in \omega
\big \} \equationmark DefineOuG $$ is a compact open subset of $\OG $.

\definition \label PartialBernoulli The \subjex {partial Bernoulli action}{Bernoulli partial action} of a group $G$ is
the topological partial action $\beta $ of $G$ on $\OuG $ obtained by restricting the global Bernoulli action to $\OuG
$, according to \cite {DefineRestriction}.

Denoting by $$ \D g = \OuG \cap \eta _g(\OuG ), \equationmark FirstDescriptBernouliDoms $$ as demanded by \cite
{DefineRestriction}, observe that an element $\omega $ in $\eta _g(\OuG )$ is characterized by the fact that $g\in
\omega $.  Thus $$ \D g = \{\omega \in \OuG : g\in \omega \} = \big \{\omega \in \OG : 1,g\in \omega \big \}.
\equationmark DescriptBernouliDoms $$

As already noted, besides being open, $\OuG $ is a compact set, hence $\D g$ is compact for every $g$ in $G$.

The partial Bernoulli action will prove to be of utmost importance in our study of partial representations subject to
relations.

\nrem Propositions \cite {TopoGlob} and \cite {TopoGlobIffClosed} were proved by Abadie in his PhD thesis \ref
{Abadie/1999}.

\chapter Algebraic partial dynamical systems

\label AgebraicSection

With this chapter we start our study of partial actions of groups on algebraic structures.  Among these, we would
eventually like to include rings, algebras, *-algebras and C*-algebras.  We thus begin with a device meant to allow us
to treat rings and algebras in the same footing.

\medbreak \fix We will assume, throughout, that $\K $ is a unital commutative ring.  Whenever convenient we will assume
that $\K $ is equipped with a \subj {conjugation}, that is, an involutive automorphism $$ r\in \K \ \mapsto \ \bar r\in
\K , \equationmark KConjuga $$ which will be fixed form now on.  In the absence of more interesting conjugations one
could take the identity map by default.  When $\K $ is the field of complex numbers we will always choose the standard
complex conjugation.

\definition By a \subj {$\K $-algebra} we will mean a ring $A$ equipped with the structure of a left $\K $-module, such
that, for all $\lambda \in \K $, and all $a,b\in A$, one has \izitem \zitem $1a=a$, \zitem $(\lambda a)b = a(\lambda b)
= \lambda (ab)$.

When $\K $ is a field, this is nothing but the usual concept of an algebra over a field.  On the other hand, any ring
$A$ may be seen as a $\Z $-algebra, as long as we equip it with the obvious left $\Z $-module structure.  With this
device we may therefore treat rings and algebras in the same footing.

Ideals in $\K $-algebras will always be assumed to be $\K $-sub-modules (that is, closed under multiplication by
elements of $\K $) and homomorphisms between $\K $-algebras will always be assumed to be $\K $-linear.

\definition \label DefineStarAlgebra By a \subj {*-algebra} over $\K $ we mean a $\K $-algebra equipped with an \subj
{involution}, namely a map $$ a\in A \mapsto a^*\in A, $$ such that, for all $a,b\in A$, and all $\lambda \in \K $, one
has \izitem \zitem $(a^*)^*=a$, \zitem $(a+\lambda b)^* = a^* + \bar \lambda b^*$, \zitem $(ab)^* = b^*a^*$.

Ideals in *-algebras will always be assumed to be \subj {self-adjoint} (that is, closed under *) and homomorphisms
between *-algebras will always be assumed to be \subjex {*-homomorphisms}{*-homomorphism}, (that is, preserving the
involution).

Among the most important examples of *-algebras, we mention the group ring $\K (G)$, for any group $G$, equipped with
the involution $$ \Big (\soma {g\in G}\lambda _g\delta _g\Big )^* = \soma {g\in G}\bar \lambda _g\delta _{g\inv }.  $$

\definition An \subjex {algebraic}{algebraic partial action} \resp {\"{*-algebraic}} partial action of a group $G$ on a
$\K $-algebra \resp {*-algebra} $A$ is a partial action $$ \Th = \big (\{\D g\}_{g\in G}, \{\th g\}_{g\in G}\big ) $$ on
the underlying set $A$, such that each $\D g$ is a \"{two-sided ideal} \resp {\"{self-adjoint two-sided ideal}} in $A$,
and each $\th g$ is an \"{isomorphism} \resp {\"{*-isomor\-phism}}.  Taking the point of view of \cite {DefPA} or \cite
{DefTopoPa}, one likewise defines the notions of \subjex {algebraic}{algebraic partial dynamical system} and
\"{*-algebraic partial dynamical systems}.

When the category one is working with is understood, be it $\K $-algebras or *-algebras, and when there is no chance for
confusion, we will drop the adjectives \"{algebraic} or \"{*-algebraic} and simply say \"{partial action} or \"{partial
dynamical system}.

We will say that two algebraic \resp {*-algebraic} partial actions of $G$ are \subjex {algebraically
equivalent}{algebraically equivalent actions} \resp {\"{*-algebraically equivalent}} if they are equivalent in the sense
of Definition \cite {DefineEquivar} and, in addition, the map $\phi $ mentioned there is a homomorphism \resp
{*-homomorphism}.

The process of restriction may also be applied in the algebraic context: given a \hbox {\mstar algebraic} global action
$\eta $ of $G$ on a \mstar algebra B, and given a two-sided ideal $A\ideal B$ (supposed to be self-adjoint in the
*-algebra case), let $\Th $ be the restriction of $\eta $ to $A$, as defined in \cite {DefineRestriction}.  Since $$ \D
g= \eta _g(A)\cap A, $$ we have that $\D g$ is a (self-adjoint) two-sided ideal of $A$ and, since the restricted maps
$\th g$ are obviously \mstar homomorphisms, we conclude that $\Th $ is a \mstar algebraic partial action of $G$ on $A$.

Reversing this process we may likewise define the notion of \mstar algebraic globalization, except that we need to adapt
the notion of orbit employed in \cite {DefineGlobalization} to the algebraic context.  Before we spell out the
appropriate definition of algebraic globalization let us discuss a minor technical point.

\state Proposition \label AlgebraOrbit Let $\eta $ be an algebraic global action of the group $G$ on the $\K $-algebra
$B$, and let $A$ be a two-sided ideal of $B$.  Then the smallest $\eta $-invariant subalgebra of $B$ containing $A$ is
$$ \soma {g\in G}\eta _g(A).  $$

\Proof Since each $\eta _g(A)$ is a two-sided ideal in $B$, the sum of all such is also a two-sided ideal in $B$, and
hence also a subalgebra, evidently $\eta $-invariant.  That it is the smallest $\eta $-invariant subalgebra is also
clear.  \endProof

If, in addition to the hypotheses of the result above, $\eta $ is a *-algebraic global action and $A$ is a self-adjoint
ideal, then it is easily verified that $ \soma {g\in G}\eta _g(A) $ is the smallest $\eta $-invariant *-subalgebra of
$B$.

The notion of algebraic globalization is thus adapted from \cite {DefineGlobalization} as follows:

\definition \label DefineAlgGloba Let $\eta $ be a \mstar algebraic global action of $G$ on a \mstar algebra $B$, and
let $A$ be a (self-adjoint) two-sided ideal of $B$.  Also let $\Th $ be the partial action obtained by restricting $\eta
$ to $A$.  If $$ B = \soma {g\in G}\eta _g(A), $$ we will say that $\eta $ is a \subjex {\mstar algebraic
globalization}{algebraic globalization} of $\Th $.

\fix In order to avoid cluttering the exposition we will now restrict ourselves to the category of $\K $-algebras, even
though most of our results from now on are valid for *-algebras as well.

\bigskip Unlike the case of partial actions on sets or topological spaces, algebraic globalizations do not always
exists, and when they do, uniqueness may fail.  Postponing the existence question, let us present an example of a
partial action admitting many non-equivalent algebraic globalizations.

\state Example \label NotUniqueGloba \rm This is essentially an example of actions on $\K $-modules masqueraded as $\K
$-algebras.  Given a left $\K $-module $M$, we may view it as a $\K $-algebra by introducing the trivial multiplication
operation, namely, $$ ab = 0 \for a,b\in M.  \equationmark GrotesqueMult $$ Notice that every $\K $-sub-module of $M$ is
an ideal and that any $\K $-linear map between two such $\K $-modules is automatically a homomorphism.

This said, let us assume that $\K =\R $, the field of real numbers, and let us choose a real number $$ \alpha \in \R
\setminus {\bf Q}.  $$

Consider the action $\eta ^\alpha $ of\/ $\Z $ on $\R ^2$ defined in terms of the group homomorphism $$ \eta ^\alpha :
n\in \Z \mapsto \left [\matrix {\cos (2n\pi \alpha ) & -\sin (2n\pi \alpha ) \cr \pilar {12pt} \sin (2n\pi \alpha ) &
\hfill \cos (2n\pi \alpha )}\right ] \in \hbox {GL}_2(\R ).  $$ In other words, each $\eta ^\alpha _n$ acts by rotating
the plane by an angle of $2n\pi \alpha $.

If, as above, we make $\R ^2$ into an algebra by introducing the trivial multiplication, we may view $\eta ^\alpha $ as
a global algebraic action of $\Z $ on $\R ^2$.

Let $A$ be the one-dimensional subspace of $\R ^2$ spanned by the vector $(1,0)$, hence a two-sided ideal, and let $\Th
$ be the restriction of $\eta ^\alpha $ to $A$.  Since $\eta ^\alpha _n$ has no real eigenvectors for nonzero $n$, we
have that $$ \D n := \eta ^\alpha _n(A)\cap A = \{0\}, $$ and consequently $\th n$ is the zero map, except of course for
$\th 0$, which is the identity map on $\D 0 = A$.  It is also easy to see that $$ \R ^2 = \soma {n\in \Z } \eta ^\alpha
_n(A), $$ so $\eta ^\alpha $ is an algebraic globalization of $\Th $.

Notice that the complex eigenvalues for $\eta ^\alpha _1$ are ${\rm e}^{\pm 2\pi i\alpha }$ so, given two distinct
irrational numbers $\alpha _1$ and $\alpha _2$ in $(0,1/2)$, the corresponding global actions $\eta ^{\alpha _1}$ and
$\eta ^{\alpha _2}$ are not equivalent.  Therefore $\Th $ admits an uncountable number of non-equivalent algebraic
globalizations.

\bigskip Many aspects of partial actions become greatly simplified when the ideals involved are unital.

\fix In order to discuss these aspects, let us fix, for the time being, a partial action $\Th = \big (\{\D g\}_{g\in G},
\{\th g\}_{g\in G}\big )$ of a group $G$ on an algebra $A$, such that $\D g$ is unital, with unit $$ 1_g\in \D g, $$ for
all $g$ in $G$.

We will now prove a few technical lemmas which will later be useful to study the globalization question under the
present hypotheses.

\state Lemma \label IdealUnits For every $g,h\in G$, one has that $$ \th g(1_{g\inv }1_h) = 1_g1_{gh}.  $$

\Proof This follows immediately from the fact that $1_{g\inv }1_h$ is the unit for $\Di g\cap \D h$, while $1_g1_{gh}$
is the unit for $\D g\cap \D {gh}$, and that $\th g$ is an isomorphisms between these ideals by \cite {NewPAIntersect}.
\endProof

Observe that $1_g$ is a central idempotent in $A$, as is the case for the unit of any two-sided ideal.  Moreover, since
$\D g = 1_gA$, the correspondence $$ \tth g: a\in A\mapsto \th g(1_{g\inv }a)\in A, \equationmark ThetaTilde $$ gives a
well defined endomorphism of $A$.

The range of $\tth g$ is clearly $\D g$, so we have that $$ \tth g(a) = \tth g(a)1_g \for a\in A.  \equationmark
ThetaRangeDg $$

\state Lemma \label FormularThetaTilde For every $g,h\in G$, and every $a\in A$, one has that $$ \tth g\big (\tth
h(a)\big ) = 1_g\tth {gh}(a).  $$

\Proof We have $$ \tth g\big (\tth h(a)\big ) = \th g\big (1_{g\inv }\th h(1_{h\inv }a)\big ) = \th g\big (1_{g\inv
}1_h\th h(1_{h\inv }a)\big ) \={IdealUnits} $$$$= \th g\big (\th h(1_{h\inv g\inv }1_{h\inv })\,\th h(1_{h\inv }a)\big )
= \th g\big (\th h(1_{h\inv g\inv }1_{h\inv }a)\big ) \={PACompos} $$$$= \th {gh}(1_{h\inv g\inv }1_{h\inv }a) = \th
{gh}(1_{h\inv g\inv }1_{h\inv }) \vezes \th {gh}(1_{h\inv g\inv }a) \={IdealUnits} $$$$= 1_{gh}1_g \tth {gh}(a) = 1_g
\tth {gh}(a).  \omitDoubleDollar \endProof

\state Lemma \label UnityInSum Suppose that $A$ is an algebra which is a (not necessarily direct) sum of a finite number
of unital ideals. Then $A$ is a unital algebra.

\Proof By induction on the number of summands it is enough to consider two ideals, say $A = J +K$.  Letting $1_J$ and
$1_K$ denote the units of $J$ and $K$, respectively, one proves that $1_J$ and $1_K$ are central idempotents in $A$ and
that $$ 1_J + 1_K - 1_J \vezes 1_K $$ is the unit for $A.$ \endProof

We are now ready to tackle our first globalization result for algebraic partial actions.  So far we will only treat the
unital case, leaving the highly involved non-unital case for later.

\state Theorem \label GlobaUnital Let $\Th = \big (\{\D g\}_{g\in G}, \{\th g\}_{g\in G}\big )$ be a partial action of a
group $G$ on a \underbar {unital} algebra $A$.  Then a necessary and sufficient condition for $\Th $ to admit a
globalization is that $\D g$ be unital for every $g$ in $G$.  Moreover, in this case the globalization is unique up to
algebraic equivalence.

\def \1#1{1\kern -1pt_{_{#1}}}

\Proof Let $\eta $ be a globalization of $\Th $, acting on the algebra $B$.  Denote by $\1A$ the unit of $A$ and observe
that, for each $g$ in $G$, one has that $$ \1A\eta _g(\1A) \in A\cap \eta _g(A) = \D g.  $$ One then easily checks that
$\1A\eta _g(\1A)$ is the unit for $\D g$, which is therefore proven to be a unital ideal.

Conversely, supposing that each $\D g$ is unital, let us build a globalization for $\Th $.  As a first step we consider
the algebra $ A^G, $ formed by all functions $$ f: G \to A, $$ with pointwise operations.  Consider the global action
$\eta $ of $G$ on $A^G$ defined by $$ \eta _g(f)|_h = f(g\inv h) \for g,h\in G \for f\in A^G, $$ and let $$ \iota :A\to
A^G $$ be the homomorphism defined by $$ \iota (a)|_g = \thi g(1_ga) = \tthi g(a) \for a\in A \for g\in G.  $$

Since $\iota (a)|_1 = a$, we have that $\iota $ is injective, so we may identify $A$ with its image $$ A':=\iota
(A)\subseteq A^G.  $$

We have therefore embedded $A$ into the larger algebra $A^G$, where a global action of $G$ is defined but,
unfortunately, $A'$ is not necessarily an ideal in $A^G$, so we need to do some more work if we are to obtain the
desired globalization of $\Th $.  The trick is to find an $\eta $-invariant subalgebra of $A^G$, containing $A'$ as an
ideal.  In doing so the crucial technical ingredient is to check that for every $g,h\in G$, one has that $$ \eta
_g(A')\vezes \eta _h(A') \subseteq \eta _g(A')\cap \eta _h(A').  \equationmark WhyIsIdeal $$ In order to prove this let
$a,b\in A$, and notice that, for all $k\in G$, one has $$ \eta _g\big (\iota (a)\big )|_k\vezes \eta _h\big (\iota
(b)\big )|_k = \iota (a)|_{g\inv k}\vezes \iota (b)|_{h\inv k} \$= \tth {k\inv g}(a)\vezes \tth {k\inv h}(b) = (\star ).
$$ This can be developed in two different ways.  On the one hand $$ (\star ) \={ThetaRangeDg} \tth {k\inv g}(a) 1_{k\inv
h} \tth {k\inv h}(b) \={FormularThetaTilde} \tth {k\inv h}\big (\tth {h\inv g}(a)\big ) \vezes \tth {k\inv h}(b) \$=
\tth {k\inv h}\big (\tth {h\inv g}(a)b\big ) = \iota \big (\tth {h\inv g}(a)b\big )|_{h\inv k} = \eta _h\Big (\iota \big
(\tth {h\inv g}(a)b\big )\Big )|_k.  $$ This shows that $$ \eta _g\big (\iota (a)\big )\vezes \eta _h\big (\iota (b)\big
) = \eta _h\Big (\iota \big (\tth {h\inv g}(a)b\big )\Big ) \in \eta _h(A').  $$ On the other hand we may write $$
(\star ) \={ThetaRangeDg} \tth {k\inv g}(a) 1_{k\inv g} \tth {k\inv h}(b) \={FormularThetaTilde} \tth {k\inv g}(a)
\vezes \tth {k\inv g}\big (\tth {g\inv h}(b)\big ) \$= \tth {k\inv g}\big (a\tth {g\inv h}(b)\big ) = \iota \big (a\tth
{g\inv h}(b)\big )|_{g\inv k} = \eta _g\Big (\iota \big (a\tth {g\inv h}(b)\big )\Big )|_k, $$ whence $$ \eta _g\big
(\iota (a)\big )\vezes \eta _h\big (\iota (b)\big ) = \eta _g\Big (\iota \big (a\tth {g\inv h}(b)\big )\Big ) \in \eta
_g(A'), $$ thus proving \cite {WhyIsIdeal}.

Let $B$ be the subset of $A^G$ defined by $$ B = \soma {g\in G} \eta _g(A').  $$

Using \cite {WhyIsIdeal} one immediately sees that $B$ is a subalgebra of $A^G$.  Moreover, by plugging $h=1$ in \cite
{WhyIsIdeal} one concludes that $A'$ is a left ideal of $B$, while with $g=1$, \cite {WhyIsIdeal} implies that $A'$ is a
right ideal of $B$, so $A'$ is a two-sided ideal.

It is clear that $B$ is invariant relative to the global action $\eta $, so we may restrict $\eta $ to $B$, obtaining
another global action which, by abuse of language, we still denote by $\eta $.  The very definition of $B$ then implies
that it is the smallest $\eta $-invariant subalgebra of $A^G$ containing $A'$, as seen in \cite {AlgebraOrbit}.

In order to complete the existence part of the proof it is now enough to check that the restriction of $\eta $ to $A'$
is equivalent to $\Th $.  As a first step let us verify that $$ \eta _g(A')\cap A' = \iota (\D g).  \equationmark
Intersesao $$

Given $x\in \eta _g(A')\cap A'$, write $x = \iota (a) = \eta _g(\iota (b))$, for $a,b\in A$, and notice that $$ a =
\iota (a)|_1 = \eta _g\big (\iota (b)\big )|_1 = \iota (b)|_{g\inv } = \tth g(b) \in \D g, $$ so $x\in \iota (\D g)$.
Conversely, if $x = \iota (a)$, with $a\in \D g$, let $b = \thi g(a)$, and notice that, for all $h\in G$, $$ \eta _g\big
(\iota (b)\big )|_h = \iota (b)|_{g\inv h} = \tth {h\inv g}(b) \$= \tth {h\inv g}\big (\tthi g(a)\big )
\={FormularThetaTilde} \tthi h(a) 1_{h\inv g}.  \equationmark IntersesaoHelp $$ Observe, however that $$ \tthi h(a) =
\thi h(1_ha) \in \thi h(\D h\cap \D g) \={NewPAIntersect} \Di h\cap \D {h\inv g}, $$ so we deduce from \cite
{IntersesaoHelp} that $$ \eta _g\big (\iota (b)\big )|_h = \tthi h(a) = \iota (a)|_h, $$ whence $$ x = \iota (a) = \eta
_g\big (\iota (b)\big ) \in \eta _g(A')\cap A', $$ proving \cite {Intersesao}.

We must now prove that the restriction of each $\eta _g$ to $\iota (\Di g)$ corresponds to $\th g$.  For this, let $a\in
\Di g$ and observe that, for all $k$ in $G$ we have $$ \eta _g\big (\iota (a)\big )|_k = \iota (a)|_{g\inv k} = \th
{k\inv g}(1_{g\inv k}a) \={NewPAIntersect} \thi k\big (\th g(1_{g\inv k}a)\big ) = (\dagger ).  $$ Notice that the
argument of $\thi k$, above, may be rewritten as $$ \th g(1_{g\inv k}a) = \th g(1_{g\inv k}1_{g\inv }a) = \th g(1_{g\inv
k}1_{g\inv })\th g(a) \={IdealUnits} $$$$ = 1_k1_g\th g(a) = 1_k\th g(a), $$ so $$ (\dagger ) = \thi k\big (1_k\th
g(a)\big ) = \tthi k\big (\th g(a)\big ) = \iota (\th g(a))|_k.  $$

This shows that $ \eta _g\big (\iota (a)\big ) = \iota \big (\th g(a)\big ), $ which is to say that restriction of $\eta
_g$ to $\iota (\D g)$ corresponds to $\th g$, as desired.

This takes care of existence, so let us now focus on proving uniqueness.  Therefore we suppose we are given two
globalizations of $\Th $, say $\eta $ and $\eta '$, acting respectively on the algebras $B$ and $B'$.  As a first step
we will prove that given $a_1,a_2,\ldots ,a_n\in A$, and $g_1,g_2,\ldots ,g_n\in G$, one has that $$ \usoma {i=1}n \eta
_{g_i}(a_i) = 0 \imply \usoma {i=1}n \eta '_{g_i}(a_i) = 0.  \equationmark ImplicationForUniqueness $$

We begin by claiming that, for all $g\in G$, and all $a$ and $b$ in $A$, one has that $$ b\eta _g(a) = b\tth g(a).
\equationmark MinClaim $$ In fact observe that $$ b\eta _g(a) = \eta _g\big (\eta _{g\inv }(b) a\big ), $$ so, since $A$
is an ideal in $B$, we see that $$ b\eta _g(a) \in A\cap \eta _g(A) = \D g.  $$ Therefore, recalling that $\eta _g$
extends $\th g$, we have that $$ b\eta _g(a) = b\eta _g(a)1_g = b\eta _g(a)\eta _g(1_{g\inv }) = b\eta _g(a1_{g\inv }) =
b\th {g}(a1_{g\inv }) = b\tth g(a), $$ proving \cite {MinClaim}.  If we now apply this for $h\inv g$, where $h$ is
another element of $G$, and then compute $\eta _h$ on both sides of the resulting expression, we obtain $$ \eta
_h(b)\eta _g(a) = \eta _h\big (b\tth {h\inv g}(a)\big ).  \equationmark MinClaimDois $$ Therefore, under the hypothesis
of \cite {ImplicationForUniqueness}, we have for all $b\in A$, and $h\in G$, that $$ 0 = \eta _h(b)\vezes \usoma {i=1}n
\eta _{g_i}(a_i) \={MinClaimDois} \usoma {i=1}n \eta _h\big (b\tth {h\inv g_i}(a_i)\big ) = \eta _h\big (\usoma {i=1}n
b\tth {h\inv g_i}(a_i)\big ).  $$ As a consequence we deduce that $$ \usoma {i=1}n b\tth {h\inv g_i}(a_i)=0.  $$

Let us denote by $Z= \usoma {i=1}n \eta '_{g_i}(a_i)$, so that proving \cite {ImplicationForUniqueness} amounts to
proving that $Z=0$.  Employing \cite {MinClaimDois}, this time for $\eta '$, it then follows that $$ 0= \eta '_h\big
(\usoma {i=1}n b\tth {h\inv g_i}(a_i)\big ) = \eta '_h(b)\vezes \usoma {i=1}n \eta '_{g_i}(a_i)= \eta '_h(b)Z.  $$ Since
$B'= \soma {h\in G}\eta '_h(A)$, this shows that $Z$ is in the annihilator ideal of $B'$.  Notice that $Z$ lies in $$
\usoma {i=1}n \eta '_{g_i}(A), $$ which is a finite sum of unital ideals of $B$, and hence itself a unital ideal by
\cite {UnityInSum}.  From this we easily see that $Z=0$, hence concluding the proof of \cite {ImplicationForUniqueness}.

As a most important consequence, it follows that there exists a well defined $\K $-linear map $\phi :B\to B'$, such that
$$ \phi \big (\usoma {i=1}n \eta _{g_i}(a_i)\big ) = \usoma {i=1}n \eta '_{g_i}(a_i), $$ for all $a_1,a_2,\ldots ,a_n\in
A$, and $g_1,g_2,\ldots ,g_n\in G$.  By reversing the implication in \cite {ImplicationForUniqueness}, we see that $\phi
$ admits an inverse, and so it is a bijection from $B$ onto $B'$.

We leave it for the reader to prove the easy facts that $\phi $ is $G$-equivariant and restricts to the identity over
$A$.

Finally, to see that $\phi $ is a homomorphism, we must prove that $ \phi (xy) = \phi (x)\phi (y), $ for all $x,y\in B$,
but it is clearly enough to consider $x=\eta _g(a)$, and $y=\eta _h(b)$, with $g,h\in G$, and $a,b\in A$.  We have $$
\phi (xy) = \phi \big (\eta _g(a)\eta _h(b)\big ) \={MinClaimDois} \phi \big (\eta _g(a\tth {g\inv h}(b)\big ) \$= \eta
'_g\big (a\tth {g\inv h}(b)\big ) \={MinClaimDois} \eta '_g(a)\eta '_h(b) = \phi (x)\phi (y).  \omitDoubleDollar
\endProof

With this result it is easy to produce examples of algebraic partial actions possessing no algebraic globalizations:

\state Example \label NoGlobalInAlg \rm Take any unital algebra $A$ containing a non-unital two-sided ideal $J$, and
consider the partial action $\Th $ of the group $\Z _2 = \{0,1\}$ on $A$, defined by setting $\th 0=id_A$, and $\th
1=id_J$, so that $\D 0 = A$, and $\D 1= J$.  By \cite {GlobaUnital}, the absence of a unit in $\D 1$ precludes the
existence of a globalization for $\Th $.

\nrem Partial actions on C*-algebras were studied long before partial actions on rings and algebras.  The reason for
this is the fact that the first proofs of the associativity property for C*-algebraic partial crossed products (defined
later in this book) uses the existence of approximate identities, which are absent in the purely algebraic case.

Having realized this difficulty early on, I have discussed it with many colleagues in the Algebra community but it took
a while before anyone was convinced of the relevance of this question.  I eventually managed to attract the attention of
Misha Dokuchaev, and together we were able to understand exactly what was going on \ref {DokuchaevExel/2005}.  Since
then, Misha has been able to convince a host of people to study the algebraic aspects of partial actions.  Unfortunately
the growing body of knowledge in that direction is poorly represented in this book and it might soon deserve a book of
its own.

Theorem \cite {GlobaUnital} was proved in \ref {DokuchaevExel/2005}.  It is an early attempt to replicate Abadie's
previous work on enveloping actions in the realm of C*-algebras \ref {Abadie/2003}, which we will discuss later in this
book.  Further attempts were made in \ref {DokuchaevExelSimon/2010} and \ref {AbadieDokuchaevExelSimon/2008}.

Although not covered by this work, the notion of \"{twisted} partial actions of groups on $\K $-algebras has also been
considered \ref {Exel/1997a}, \ref {DokuchaevExelSimon/2008}, \ref {DokuchaevExelSimon/2010}.

\chapter Multipliers

The notion of multiplier algebras is a well established concept among C*-algebraists, but its purely algebraic version
has not appeared too often in the specialized literature.  Since multipliers are an important tool in the treatment of
algebraic partial actions we have included this brief chapter to subsidize our future work.  See \ref
{DokuchaevExel/2005} for more details.

\medskip \fix Throughout this chapter we shall let $A$ be a fixed $\K $-algebra, sometimes assumed to be a *-algebra.

\definition \label DefineMultiplier Let $L$ and $R$ be $\K $-linear maps from $A$ to itself.  We shall say that the pair
$(L,R)$ is a \subjex {multiplier}{multiplier pair} of $A$ if, for every $a$ and $b$ in $A$, one has that \izitem \zitem
$L(ab) = L(a)b$, \zitem $R(ab) = aR(b)$, \zitem $R(a)b = a L(b)$.

If $A$ is a unital algebra, and if $(L,R)$ is a multiplier of $A$, notice that, $$ R(1) = R(1)1 \={DefineMultiplier/iii}
1L(1) = L(1).  $$ Moreover, letting $m:= R(1) = L(1)$, we have for every $a$ in $A$ that $$ L(a) = L(1a) = L(1)a = ma,
$$ while $$ R(a) = R(a1) = aR(1) = am, $$ so we see that $(L,R)=(L_m,R_m)$, where $L_m$ is the operator of left
multiplication by $m$, and $R_m$ is the operator of right multiplication by $m$.

\definition \label Multipalgebra The \subj {multiplier algebra} of $A$ is the set $\Mult (A)$ consisting of all
multipliers $(L,R)$ of $A$.  Given $(L,R)$ and $(L',R')$ in $\Mult (A)$, and $\lambda \in \K $, we define $$ \lambda
(L,R) = (\lambda L, \lambda R), $$ $$ (L,R) + (L',R') = (L+L', R+R'), $$ $$ (L,R) (L', R') = (L \circ L', R' \circ R).
$$

We leave it for the reader to prove that $\Mult (A)$ is an associative $\K $-algebra with the above operations.  It is
moreover a unital algebra with unit $(id_A, id_A)$.

If $A$ is a *-algebra, and $T:A\to A$ is a linear map, define a mapping $T^*:A\to A$ by the formula $$ T^*(a) = \big
(T(a^*)\big )^*.  $$ It is then easy to see that $T^*$ is also linear.

\state Proposition If $A$ is a *-algebra then so is $\Mult (A)$, when equipped with the involution $$ (L,R)^*:=
(R^*,L^*) \for (L,R)\in \Mult (A).  $$

\Proof Left for the reader. \endProof

Assuming that $A$ is an ideal in some other algebra $B$, let $m$ be a fixed element of $B$ and consider the maps $$ L_m:
a\in A \mapsto ma \in A \and R_m: a\in A \mapsto am \in A.  $$

It is an easy exercise to show that the pair $(L_m,R_m)$ is a multiplier of $A$ and that the correspondence $$ \mu _B:
m\in B \mapsto (L_m,R_m)\in \Mult (A), \equationmark ExternalMultiplier $$ is a homomorphism.

\definition We shall say that $A$ is an \subj {\essIdeal } in $B$ if the map $\mu _B$ mentioned above is injective.

In general the kernel of $\mu _B$ is the intersection of the left and right annihilators of $A$ in $B$.  Therefore $A$
is an {\essIdeal } in $B$ if and only if, for every nonzero element $a$ in $A$, there exists some $b$ in $B$ such that
either $ab\neq 0$ or $ba\neq 0$.

\definition We shall say that $A$ is a \subjex {\nonDegAlg }{\nonDegAlg { }algebra} algebra if there is no nonzero
element $a$ in $A$ such that $ab=ba=0$, for all $b$ in $A$.  Equivalently, if $A$ is an {\essIdeal } in itself.

\state Proposition \izitem \zitem $\mu _A$ is an isomorphism onto $\Mult (A)$ if and only if $A$ is unital.  \zitem The
range of $\mu _A$, henceforth denoted by $A'$, is an ideal of $\Mult (A)$.  \zitem If $A$ is {\nonDegAlg } then $A'$ is
isomorphic to $A$.

\Proof Left for the reader. \endProof

Thus, when $A$ is an {\essIdeal } in some algebra $B$, a situation which is only possible when $A$ is {\nonDegAlg }, we
have that $B$ is isomorphic to a subalgebra of $\Mult (A)$, namely the range of $\mu _B$, containing $A'$.  In other
words, up to isomorphism, all algebras $B$ containing $A$ as an {\essIdeal } are to be found among the subalgebras of
$\Mult (A)$ containing $A'$.

\medskip Given two multipliers $(L^1,R^1)$ and $(L^2,R^2)$ in $\Mult (A)$ we shall be concerned with the validity of the
formula $$ R^2L^1 \explica ={?} L^1R^2.  \equationmark Extraassoc $$ We will see that this is the crucial property
governing the associativity of partial crossed products studied below.

To see the connection between associativity and property \cite {Extraassoc} above, notice that if $A$ is an ideal in
some algebra $B$, and if for all $i=1,2$ we let $(L^i,R^i)=\mu _B(m_i)$, for some $m_i\in B$, then for every $a$ in $A$
we have that $$ R^2\big (L^1(a)\big ) = (m_1a)m_2 ,\quad {\rm while} \quad L^1\big (R^2(a)\big ) = m_1(am_2), $$ so that
\cite {Extraassoc} holds as a consequence of the associativity property of $B$.

However \cite {Extraassoc} may just as well fail: take $A$ to be any $\K $-module equipped with the trivial
multiplication operation described in \cite {GrotesqueMult}.  In this case, observe that any pair $(L,R)$ of linear
operators on $A$ would constitute a multiplier and one should clearly not expect \cite {Extraassoc} to hold in such a
generality!

\state Proposition \label Condforextraassoc If $A$ is either {\nonDegAlg } or idempotent then \cite {Extraassoc} holds
for any pair of multipliers $(L^1,R^1)$ and $(L^2,R^2)$.

\Proof Given $a,b\in A$, we have that $$ R^2\big (L^1 (a)\big )b = L^1 (a)L^2(b) = L^1 \big (aL^2(b)\big ) = L^1 \big
(R^2(a)b\big ) = L^1 \big (R^2(a)\big )b.  $$ This shows that $R^2\big (L^1 (a)\big )-L^1 \big (R^2(a)\big )$ lies in
the left annihilator of $A$.  With a similar argument one shows that this also lies in the right annihilator of $A$.  So
this proves \cite {Extraassoc} under the assumption that $A$ is {\nonDegAlg }.

Next suppose we are given $a_1,a_2\in A$.  Letting $a=a_1a_2$, notice that $$ R^2\big (L^1 (a)\big ) = R^2\big (L^1
(a_1a_2)\big ) = R^2(L^1 (a_1)a_2) = L^1 (a_1)R^2(a_2) \$= L^1 \big (a_1R^2(a_2)\big ) = L^1 \big (R^2(a_1a_2)\big ) =
L^1 \big (R^2(a)\big ).  $$ Assuming that $A$ is idempotent, we have that every element of $A$ is a sum of terms of the
form $a_1a_2$, whence the conclusion.  \endProof

\chapter Crossed products

\label CrossProdChap

The classical notion of crossed product\fn {Algebraists sometimes reserve the term \"{crossed product} for situations
when, besides a group action, a cocycle is also involved.  When only a group action is present, as above, the algebraic
literature usually favors the expression \"{skew-group algebra}.}  is a useful tool to construct interesting examples of
algebras.  Its basic ingredient is a global action $\eta $ of a group $G$ on an algebra $A$.  One then defines the
crossed product algebra $A\rt \eta G$, sometimes also denoted $A\Rt \eta G$ when we want to make the action explicit, to
consist of all finite formal linear combinations $$ \soma {g\in G}a_g\delta _g, \equationmark GlobalPolynomial $$ (that
is $a_g=0$, except for finitely many $g$'s), where the $\delta _g$ are seen as place markers\fn {Technically speaking,
$A\rt \eta G$ is the set of all finitely supported functions from $G$ to $A$.  Denoting by $a_g\delta _g$ the function
supported on the singleton $\{g\}$, and whose value at $g$ is the element $a_g$, the term in \cite {GlobalPolynomial}
corresponds simply to the function sending $g$ to $a_g$.  Note however that the expression ``$\delta _g$'' has no
meaning in itself according to the present convention.}.

A multiplication operation on $A\rt \eta G$ is then defined by $$ (a\delta _g)(b\delta _h) = a\eta _g(b)\delta _{gh},
\equationmark ProdInCrossProd $$ for all $a$ and $b$ in $A$, and all $g$ and $h$ in $G$.

\def \sixbox #1{\hbox {\sixrm #1}} \def \twolines #1#2{{\buildrel \sixbox {#1} \over {\sixbox {#2}}}}

The intuitive idea behind the above definition of product is to think of the $\delta _g$'s as invertible elements
implementing the given action.  A somewhat imprecise but highly enlightening calculation motivating this definition is
as follows $$ (a\delta _g)(b\delta _h) = \underbrace {a\delta _gb\delta _h}_\twolines {eliminate}{parentheses} = a\delta
_gb\underbrace {\delta _g\inv \delta _g}_{\sixbox { insert 1}} \delta _h = a \kern -5pt\underbrace {\delta _gb\delta
_g\inv }_\twolines {view as}{conjugation} \kern 4pt \underbrace {\delta _g\delta _h}_\twolines {apply the}{group law} =
a\eta _g(b)\delta _{gh}.  $$

It is our intention to develop a similar construction for partial actions.

\medskip \noindent \fix So let us fix, for the remainder of this chapter, a group $G$, a $\K $-algebra $A$, and an
algebraic partial action $\Th = \big (\{\D g\}_{g\in G}, \{\th g\}_{g\in G}\big )$ of $G$ on $A$.

Evidently the main difficulty in generalizing the above construction of the crossed product is that, should $\eta _g$ be
a partially defined map, the term $\eta _g(b)$ appearing in \cite {ProdInCrossProd} is only defined for $b$ in the
domain of $\eta _g$.

\definition \label DefCrossProd The \subj {crossed product}\fn {In order to emphasize that the action $\Th $ is partial,
we will sometimes use the expression \subj {partial crossed product} to refer to this construction.}  of the algebra $A$
by the group $G$ under the partial action $\Th = \big (\{\D g\}_{g\in G}, \{\th g\}_{g\in G}\big )$ is the algebra $A\rt
\Th G$, sometimes also denoted $A\Rt \Th G$ when we want to make $\Th $ explicit, consisting of all finite formal linear
combinations $$ \soma {g\in G}a_g\delta _g, \equationmark Polynomial $$ with $a_g\in \D g$, for all $g$ in $G$.
Addition and scalar multiplication are defined in the obvious way, while multiplication is determined by $$ (a \delta
_g) (b \delta _h) = \th g \big (\thi g(a)b\big ) \delta _{gh}, \equationmark ProdInParCrossProd $$ for all $g$ and $h$
in $G$, and for all $a\in \D g$, and $b\in \D h$.

Here is another imprecise but likewise enlightening calculation, this time motivating the definition of the
multiplication in \cite {ProdInParCrossProd}: $$ (a\delta _g)(b\delta _h) = \underbrace {a\delta _gb\delta _h}_\twolines
{eliminate}{parentheses} = \underbrace {\delta _g\delta _g\inv }_{\sixbox { insert 1}} a\delta _gb \underbrace {\delta
_g\inv \delta _g}_{\sixbox {insert 1}}\delta _h = \delta _g \kern -4pt \underbrace {\delta _g\inv a\delta _g}_\twolines
{view as}{conjugation} \kern -4pt b\delta _g\inv \kern -7pt \underbrace {\delta _g\delta _h}_\twolines {apply the}{group
law}\$= \underbrace {\delta _g \thi g(a)b\delta _g\inv }_\twolines {view again as}{conjugation} \kern 0pt \delta _{gh} =
\th g \big (\thi g(a)b\big ) \delta _{gh}.  \vrule depth 30pt width 0pt $$

Although slightly longer than \cite {ProdInCrossProd}, the definition of the multiplication in \cite
{ProdInParCrossProd} has the advantage of avoiding the temptation of applying a function to elements outside its domain.
In fact, on the one hand notice that since $a\in \D g$, the reference to $\thi g(a)$ in \cite {ProdInParCrossProd} is
perfectly legal.  On the other hand, since $\thi g(a)$ belongs to the \"{ideal} $\Di g$, we see that $\thi g(a)b$ is in
$\Di g$ as well, so the reference to $\th g \big (\thi g(a)b\big )$ is legal too.

As a final syntactic check, recall that in \cite {Polynomial}, each $a_g$ must lie in $\D g$, so we need to ensure that
$\th g \big (\thi g(a)b\big )$ lies in $\D {gh}$, if we are to allow it to stand besides $\delta _{gh}$.  But this is
also granted because $$ \thi g(a)b \in \Di g \cap \D h, $$ so $$ \th g \big (\thi g(a)b\big ) \in \th g(\Di g \cap \D h)
\={NewPAIntersect} \D g\cap \D {gh}.  $$

After extending the multiplication defined above to a bi-linear map on $A\rt \Th G$, we must worry about the
associativity property, but unfortunately it does not always hold.

In order to identify the origin of this problem let us begin with the trivial remark that $A\rt \Th G$ is associative if
and only if $$ (a \delta _g \ b\delta _h)\ c \delta _k = a \delta _g \ ( b \delta _h \ c \delta _k), \equationmark
EqAssoc $$ for all $g,h,k\in G$, and all $a\in \D g$, $b\in \D h$, and $c\in \D k$.  Focusing on the left-hand-side
above we have $$ (a \delta _g \ b\delta _h)\ c \delta _k = \th g\big (\thi g(a)b\big )\delta _{gh} \ c \delta _k \$= \th
{gh}\Big (\th {h\inv g\inv }\Big [\th g\big (\thi g(a)b\big )\Big ]c\Big )\delta _{ghk} = (\star _{_1}).  $$ Observe
that the term within square brackets above satisfies $$ \th g\big (\thi g(a)b\big ) \in \th g(\Di g\cap \D h) = \D g\cap
\D {gh}, $$ which is precisely the set on which \cite {PACompos} yields $\th {h\inv g\inv } = \thi h\thi g$.  Therefore
we may cancel out the composition ``$\thi g \th g$'', so that $$ (\star _{_1})= \th {gh}\Big (\ \underline {\vrule width
0pt depth 8pt\th {h\inv }\big (\thi g(a)b}\big ) \ c\Big )\delta _{ghk} = (\star _{_2}).  $$ Observe also that the
underlined term above satisfies $$ \th {h\inv }\big (\thi g(a)b\big ) \in \thi h(\Di g\cap \D h) = \D {h\inv g\inv }\cap
\Di h, $$ which is where $\th { gh} = \th g\th h$, again according to \cite {PACompos}.  We thus finally obtain that $$
(\star _{_2}) = \th g\Big [\th h\Big (\th {h\inv }\big (\thi g(a)b\big )c\Big )\Big ]\delta _{ghk}.  $$

On the other hand, the right-hand-side of \cite {EqAssoc} equals $$ a \delta _g \ ( b \delta _h \ c \delta _k) = a
\delta _g \Big ( \th h\big (\thi h(b)c\big ) \delta _{hk}\Big ) = \th g\Big [\thi g(a) \th h\big (\thi h(b)c\big ) \Big
]\delta _{ghk}.  $$ This said, we see that equation \cite {EqAssoc} is equivalent to $$ \th g\Big [\thi g(a) \th h\big
(\thi h(b)c\big ) \Big ] = \th g\Big [\th h\Big (\th {h\inv }\big (\thi g(a)b\big )c\Big )\Big ], $$ which is clearly
the same as $$ \thi g(a) \th h\big (\thi h(b)c\big ) = \th h\Big (\th {h\inv }\big (\thi g(a)b\big )c\Big ).  $$

Since the only occurrence of $g$ and $a$, above, is in the term $\thi g(a)$, me may substitute $a$ for $\thi g(a)$
without altering the logical content of this expression.  We have therefore proven:

\state Lemma \label NessSufAssoc A necessary and sufficient condition for $A\rt \Th G$ to be associative is that, for
all $h\in G$, $b\in \D h$, and $a,c\in A$, one has $$ a\th h\big (\thi h(b)c\big ) = \th h\Big (\th {h\inv }\big (ab\big
)c\Big ).  \equationmark CleanEqAssoc $$

The above condition might still look a bit messy at first, but it may be given a clean interpretation in terms of
multipliers.  With this goal in mind, consider the multiplier of $\D h$ given by $(L^1,R^1) = \mu _A(a)$, as defined in
\cite {ExternalMultiplier}.

There is another relevant multiplier of $\D h$ given as follows: initially consider the multiplier $$ \mu _A(c) =
(L_c,R_c) \in \Mult (\Di h).  $$ Since $\th h$ is an isomorphism from $\Di h$ to $\D h$, we may transfer the above to a
multiplier $(L^2,R^2)\in \Mult (\D h)$, by setting $$ L^2 = \th h L_c \thi h \and R^2 = \th h R_c \thi h.  $$ Given
$b\in \D h$, notice that $$ L^1\big (R^2(b)\big ) = a\th h \big (\thi h(b)c\big ) \and R^2\big (L^1(b)\big ) = \th h
\big (\thi h(ab)c\big ), \equationmark MessyMultipliers $$ so that \cite {CleanEqAssoc} is precisely expressing the
commutativity of $L^1$ and $R^2$, as discussed in \cite {Extraassoc}.

We thus arrive at the main result of this chapter:

\state Theorem \label TheoremOnAssoc Given an algebraic partial action $\Th = \big (\{\D g\}_{g\in G}, \{\th g\}_{g\in
G}\big )$ of the group $G$ on the algebra $A$, a sufficient condition for $A\rt \Th G$ to be an associative algebra is
that each $\D g$ be either {\nonDegAlg } or idempotent.

\Proof Follows immediately from \cite {MessyMultipliers}, \cite {NessSufAssoc}, and \cite {Condforextraassoc}.
\endProof

We refer the reader to \ref {Proposition 3.6/DokuchaevExel/2005} for an example of a non-associative partial crossed
product.

In virtually all of the examples of interest to us, sufficient conditions for associativity will be present, so we will
not have to deal with non-associative algebras.  However, the next few general results in this chapter do not need
associativity, so we will temporarily be working with possibly non-associative algebras.

\medskip The following elementary fact should be noticed:

\state Proposition \label EmbeddIntoCrossProd Given an algebraic partial action $\Th $ of the group $G$ on the algebra
$A$, the correspondence $$ a\in A \mapsto a\delta _1\in A\rt \Th G $$ is an injective homomorphism.

We will therefore often identify $A$ with its image in $A\rt \Th G$ under the above map.

Let us now briefly study crossed products by *-algebraic partial actions.

\state Proposition \label StarInCrossProd Given a *-algebraic partial action $$ \Th = \big (\{\D g\}_{g\in G}, \{\th
g\}_{g\in G}\big ) $$ of the group $G$ on the *-algebra $A$, let $$ *: A\rt \Th G \to A\rt \Th G, $$ be the unique
conjugate-linear map such that $$ (a\delta _g)^* = \thi g(a^*)\delta _{g\inv } \for g\in G \for a\in \D g.  $$ Then this
operation satisfies \cite {DefineStarAlgebra/i--iii}, so $A\rt \Th G$ is a (possibly non-associa\-tive) *-algebra.

\Proof We leave the easy proofs of \cite {DefineStarAlgebra/i--ii} to the reader and concentrate on \cite
{DefineStarAlgebra/iii}.  We must then prove that $$ (a\delta _g\ b\delta _h)^* = (b\delta _h)^* (a\delta _g)^* \for
g,h\in G\for a\in \D g \for b\in \D h.  \equationmark StarAntiMult $$ Observe that the left-hand-side above equals $$
(a\delta _g\ b\delta _h)^* = \Big (\th g \big (\thi g(a)b\big ) \delta _{gh}\Big )^* \$= \th {h\inv g\inv }\Big (\th g
\big (b^*\thi g(a^*)\big ) \Big ) \delta _{h\inv g\inv } = (\star ).  $$ Notice that the term between the outermost pair
of parenthesis above satisfies $$ \th g \big (b^*\thi g(a^*)\big ) \in \th g (\D h \cap \Di g) \in \D {gh}\cap \D g, $$
where $\th {h\inv g\inv }$ coincides with $\thi h\thi g$, by \cite {PACompos}.  So $$ (\star ) = \thi h \big (b^*\thi
g(a^*)\big )\delta _{h\inv g\inv }.  $$

Speaking of the right-hand-side of \cite {StarAntiMult}, we have $$ (b\delta _h)^* (a\delta _g)^* = \big (\thi
h(b^*)\delta _{h\inv }\big )\ \big (\thi g(a^*)\delta _{g\inv }\big ) \$= \thi h\big (b^*\thi g(a^*)\big ) \delta
_{h\inv g\inv }, $$ which coincides with $(\star )$ hence proving \cite {StarAntiMult}.  This concludes the proof.
\endProof

Under the conditions of the above Proposition, it is clear that the canonical embedding of $A$ into $A\rt \Th G$
described in \cite {EmbeddIntoCrossProd} is a *-homomorphism.

\medskip

One of the main applications of partial actions is to the theory of graded algebras.  Let us therefore formally
introduce this important concept.

\definition \label AlgGradedALg Let $B$ be a \mstar algebra and let $G$ be a group.  By a \subj {$G$-grading} of $B$ we
shall mean an independent collection $ \{B_g\}_{g\in G} $ of $\K $-sub-modules of $B$, such that $B=\bigoplus _{g\in G}
B_g$, and $$ B_g B_h \subseteq B_{gh} \for g,h\in G.  $$ In the *-algebra case we also require that $$ B_g^* \subseteq
B_{g\inv } \for g\in G.  $$ Given such a $G$-grading, we say that $B$ is a \subj {$G$-graded algebra}, and each $B_g$ is
called a \subj {grading space}, or a \subj {homogeneous space}.

The next result states that a partial crossed product algebra has a canonical grading.  Its proof is left as an easy
exercise to the reader.

\state Proposition \label PAisGraded Given a \mstar algebraic partial action $$ \Th = \big (\{\D g\}_{g\in G}, \{\th
g\}_{g\in G}\big ) $$ of\/ $G$ on the \mstar algebra $A$, let $B_g$ be the subset of $A\rt \Th G$ given by $$ B_g = \D
g\delta _g.  $$ Then $\{B_g\}_{g\in G}$ is a $G$-grading of $A\rt \Th G$.

As already claimed, partial actions will often be used to describe graded algebras.  Given a graded algebra $B=\bigoplus
_{g\in G} B_g$, the plan is to look for sufficient conditions to ensure the existence of a partial action of $G$ on\fn
{Recall that, in the grading of $A\rt \Th G$ provided by \cite {PAisGraded}, the grading subspace corresponding to the
unit group element is naturally isomorphic to $A$.}  $B_1$ such that $B$ is isomorphic to $B_1\rt \Th G$.

As an example let $M_n(\K )$ denote the algebra of all $n\times n$ matrices with coefficients in $\K $, and consider the
$\Z $-grading of $M_n(\K )$ given by $$ B_k = \span \{e_{i,j}: i-j=k\} \for k\in \Z , $$ where the $e_{i,j}$ are the
standard matrix units.  $$ \left [\matrix { 0 & 0 & \bigstar & 0 & 0 & 0 \cr 0 & 0 & 0 & \bigstar & 0 & 0 \cr 0 & 0 & 0
& 0 & \bigstar & 0 \cr 0 & 0 & 0 & 0 & 0 & \bigstar \cr 0 & 0 & 0 & 0 & 0 & 0 \cr 0 & 0 & 0 & 0 & 0 & 0 \cr } \right ]
$$ \hfill {\eightrm An illustration of $\scriptstyle B_{-2}$ in $\scriptstyle M_6(\K )$.} \hfill \null

\goodbreak \bigskip

One may immediately rule out the possibility that this grading comes from a global action since $B_k=\{0\}$, for all
$|k|\geq n$.  However we will soon see that it may be effectively described as the crossed product algebra for a partial
action of $\Z $ on $\K ^n$.

\medskip

Let us now discuss a few elementary functorial properties of crossed products.

\state Proposition \label Functoriality Let $G$ be a group and suppose we are given \mstar alge\-braic partial dynamical
systems $$ \Th ^i = \big (A^i,\ G,\ \{\D g^i\}_{g\in G},\ \{\th g^i\}_{g\in G}\big ), $$ for $i=1,2$.  Suppose, in
addition, that $\varphi :A^1\to A^2$ is a $G$-equivariant \mstar homomorphism.  Then there is a graded\fn {A map $\psi $
between two graded algebras $B=\bigoplus _{g\in G} B_g$ and $C=\bigoplus _{g\in G} C_g$ is said to be \subjex
{graded}{graded map} if $\psi (B_g)\subseteq C_g$, for every $g$ in $G$.}  \mstar homomorphism $$ \psi :A^1\rt {\Th ^1}G
\to A^2\rt {\Th ^2}G, $$ such that $$ \psi (a\delta _g) = \varphi (a)\delta _g \for a\in \D g^1, $$ (we believe the
context here is enough to determine the appropriate interpretations of the expressions ``$a\delta _g$'' and ``$\varphi
(a)\delta _g$'', above, so as to dispense with heavier notation such as ``$a\delta _g^1$'' and ``$\varphi (a)\delta
_g^2$'').

\Proof Left to the reader.  \endProof

Injectivity and surjectivity of the induced map is discussed in the next Proposition.  The proof is routine so we again
leave it for the reader.

\state Proposition \label InjSurjCross Under the conditions of \cite {Functoriality} one has: \izitem \zitem If $\varphi
$ is injective, then so is $\psi $.  \zitem If $\varphi (\D g^1) = \D g^2$, for every $g$ in $G$, then $\psi $ is
surjective.

\bigskip There are some recurring calculations with elements of $A\rt \Th G$ which are useful to know in advance.  So,
before concluding this chapter, let us collect some of these for later reference.

\state Proposition \label LotsAFormulas In what follows, every time we write $a\delta _g$, it is assumed that $g$ is an
arbitrary element of $G$ and $a$ is an arbitrary element of $\D g$, satisfying explicitly stated conditions, if any.  If
an expression involves a ``*'', we assume we are speaking of a *-algebraic partial dynamical system. \goodbreak
\noindent Wherever $1_g$ is mentioned, we tacitly assume that $\D g$ is unital with unit $1_ g$.  In this case recall
that $\tth g$ was defined in \cite {ThetaTilde}.  \def \lign #1 ={\hbox to 100pt{\hfill $#1$\hfill } = \ \ } \medskip
\iaitem \aitem $\lign (a\delta _1) (b\delta _h) = ab \delta _h$, \aitem $\lign (a\delta _g) (b\delta _h) = a\th g(b)
\delta _{gh}$, \quad provided $b\in \Di g\cap \D h$, \aitemBismark LotsPAreceGloba \aitem $\lign \big (\th g(a)\delta
_g\big ) (b\delta _h) = \th g(ab) \delta _{gh}$, \quad provided $a\in \Di g$, \aitem $\lign (u\delta _g) (a\delta _1)
(v\delta _{g\inv }) = u\th g(av) \delta _1$, \aitem $\lign (a\delta _g)^* (b\delta _h) = \thi g(a^*b)\delta _{g\inv h}$,
\aitemBismark LotsLefInnProd \aitem $\lign (a\delta _g) (b\delta _g)^* = ab^*\delta _1$, \aitemBismark LotsRigInnProd
\aitem $\lign (1_g\delta _1) (a\delta _g) = a\delta _g$, \aitem $\lign (a\delta _g) (1_{g\inv }\delta _1) = a\delta _g$,
\aitem $\lign (1_g\delta _g) (a\delta _1) (1_{g\inv }\delta _{g\inv }) = \th g(a) \delta _1$, \qquad provided $a\in \Di
g$, \aitem $\lign (a\delta _g) (b\delta _h) = a\th g(1_{g\inv }b) \delta _{gh} \ = \ a\tth g(b) \delta _{gh}$.

\Proof Left for the reader. \endProof

\nrem Theorem \cite {TheoremOnAssoc} was first proved in \ref {DokuchaevExel/2005}.  It will be used again to give a
proof of the associativity of C*-algebraic partial crossed products which, unlike the original proofs, does not use
approximate identities, relying only on the fact that C*-algebras are {\nonDegAlg }.

Similar results on the associativity of \"{twisted} partial crossed products can be found in \ref {Proposition
2.4/Exel/1997a} and \ref {Theorem 2.4/DokuchaevExelSimon/2008}.

\chapter Partial group representations

Given a partial action $\Th = \big (\{\D g\}_{g\in G}, \{\th g\}_{g\in G}\big )$ of a group $G$ on an algebra $A$, we
have seen that $A\rt \Th G$ consists of elements of the form $$ \soma {g\in G}a_g\delta _g, $$ with $a_g\in \D g$.  As
already mentioned the $\delta _g$ appearing above are just place markers and thus have no meaning by themselves.
However, in case $\Th $ is a global action and $A$ is a unital algebra, the elements $$ v_g:= 1\delta _g \for g\in G, $$
make perfectly good sense and are indeed crucial in the study of skew products.  In that case it is easy to see that
each $v_g$ is an invertible element in $A\rt \Th G$ and, in addition, we have that $$ v_{gh}=v_gv_h \for g,h\in G, $$ so
$v$ may be seen as a group representation of $G$ in $A\rt \Th G$.  Identifying $A$ with its canonical copy in $A\rt \Th
G$, we also have that $$ v_g a v_g\inv = \th g(a) \for g\in G\for a\in A, $$ which is to say that $\th g$ coincides with
the inner automorphism $Ad_{v_g}$ on $A$.

In the case of a partial action the above definition of $v_g$ does not make sense, unless we assume that each $\D g$ is
a unital ideal in which case we may redefine $$ \pr g := 1_g\delta _g \for g\in G, $$ where $1_g$ is the unit of $\D g$.
With the appropriate modifications we will see that the $\pr g$ play as important a role relative to partial action as
the $v_g$ does in the global case.

\definition \label DefPR Given a group $G$ and a unital algebra $B$, we shall say that a map $$ \Pr : G \to B $$ is a
\subj {partial representation} of $G$ in $B$ if, for all $g$ and $h$ in $G$, one has that \izitem \zitem $\pr 1 = 1$,
\zitemBismark PrOnOne \zitem $\pr g\pr h\pri h = \pr {gh}\pri h$, \zitemBismark ObserveRight \zitem $\pri g\pr g\pr h =
\pri g\pr {gh}$. \zitemBismark ObserveLeft \medskip \noindent If $B$ is a *-algebra and $\Pr $ moreover satisfies \zitem
$\pri g = (\pr g)^*$, \zitemBismark StarPrepCond \medskip \noindent we will say that $\Pr $ is a \subj {*-partial
representation}.

The above definition of partial group representation may be generalized to maps from $G$ to any unital semigroup (also
called monoid), with or without an involution.  Specializing this to the multiplicative semigroup of an algebra $B$, we
evidently recover the above definition.

\nostate \label LeaveOutOneAxiom Given a unital *-algebra $B$, and a map $\Pr : G\to B$ satisfying \cite {StarPrepCond},
notice that \cite {ObserveRight} is equivalent to \cite {ObserveLeft}.  Thus, for such a map to be proven a *-partial
representation, one may choose to check only one axiom among \cite {ObserveRight} and \cite {ObserveLeft}.

\medskip For comparison purposes, let us recall a well known concept:

\definition Given a group $G$ and a unital algebra $B$, we shall say that a map $ v: G \to B $ is a \subj {group
representation} of $G$ in $B$ if, for all $g$ and $h$ in $G$, one has that \izitem \zitem $v_1 = 1$, and \zitem $v_gv_h
= v_{gh}$.  \medskip \noindent If $B$ is a *-algebra and $v$ moreover satisfies \zitem $v_{g\inv } = (v_g)^*$, \medskip
\noindent we will say that $v$ is a \subj {unitary representation}.

Needless to say, a (unitary) group representation is a \mstar partial group representation, but we shall encounter many
examples of partial representations which are not group representations.

\medskip A rule of thumb to memorize condition \cite {ObserveRight}, above, is to think of it as the usual axiom for a
group representation, namely ``$\pr g\pr h = \pr {gh}$'', but which is being ``observed'' by the term $\pri h$ on the
right.  The observer apparently does not take any part in the computation and, in case it is invertible, one could
effectively cancel it out and be left with the traditional ``$\pr g\pr h = \pr {gh}$''.  A similar comment could of
course be made with respect to the left-hand observer $\pri g$ in \cite {ObserveLeft}.

Returning to the discussion at the beginning of the present chapter, we now present an important source of examples of
partial representations.

\state Proposition \label PRinCrossProd Let $\Th = \big (\{\D g\}_{g\in G}, \{\th g\}_{g\in G}\big )$ be a \mstar
partial action of a group $G$ on a \mstar algebra $A$, such that each $\D g$ is unital, with unit denoted $1_g$.  Then
the map $$ \Pr : g\in G \mapsto 1_g\delta _g \in A\rt \Th G, $$ is a \mstar partial representation.

\Proof Condition \cite {PrOnOne} is evident.  As for \cite {ObserveLeft}, let $g,h\in G$, and observe that, using \cite
{LotsAFormulas} we have $$ \pr g\pr h = (1_g \delta _g) (1_h \delta _h) = 1_g \th g(1_{g\inv } 1_h) \delta _{gh} = 1_g
1_{gh} \delta _{gh}, $$ so $$ \pri g \pr g\pr h = (1_{g\inv }\delta _{g\inv })(1_g 1_{gh} \delta _{gh}) = 1_{g\inv }\thi
g(1_g 1_{gh}) \delta _{h} \$= 1_{g\inv }1_h \delta _{h} = (\star ).  $$ On the other hand, $$ \pri g \pr {gh} =
(1_{g\inv } \delta _{g\inv }) (1_{gh} \delta _{gh}) = 1_{g\inv } \thi g(1_g1_{gh}) \delta _h = 1_{g\inv }1_h \delta
_{h}, $$ which coincides with $(\star )$, hence proving \cite {ObserveLeft}.  The proof of \cite {ObserveRight} follows
along similar lines.

In the *-algebra case, observe that $1_g^*$ is also a unit for $\D g$ and, since an algebra has at most one unit, we
must have $1_g^*=1_g$.  Therefore $$ (\pr g)^* = (1_g\delta _g)^* = \thi g(1_g^*)\delta _{g\inv } = \thi g(1_g)\delta
_{g\inv } = 1_{g\inv }\delta _{g\inv } = \pri g.  $$ proving \cite {StarPrepCond}.  \endProof

Let us now discuss another source of examples of partial representations, but let us first observe that, given an
idempotent element $p$ in an algebra $B$, then $pBp$ is a unital subalgebra of $B$, with unit $p$.  If $B$ is moreover a
*-algebra and $p$ is self-adjoint (that is, $p^*=p$), then $pBp$ is clearly also a *-algebra.

\state Proposition \label CompressPrep Let $B$ be a unital \mstar algebra, and $v$ be a (unitary) group representation
of\/ $G$ in $B$.  Suppose that $p$ is a (self-adjoint) idempotent element of $B$ such that $v_gpv_{g\inv }$ commutes
with $p$, for every $g$ in $G$.  Then the formula $$ \pr g = pv_gp \for g\in G $$ defines a \mstar partial
representation of $G$ in the unital \mstar algebra $pBp$.

\Proof With respect to \cite {PrOnOne}, we clearly have that $\pr 1=p$, which is the unit of $pBp$.  Given $g$ and $h$
in $G$ we have $$ \def \linha #1{\,\underline {#1\stake {5pt}}\,} \pri g\pr g \pr h = p \linha {v_{g\inv } pv_g} \linha
{p} v_hp = p \linha {p} \linha {v_{g\inv } pv_g} v_hp \$= pv_{g\inv } pv_{gh}p = \pri g\pr {gh}, $$ proving \cite
{ObserveLeft}, while \cite {ObserveRight} may be proved similarly.  In the *-algebra case we have $$ \pri g = pv_{g\inv
}p = pv_g^*p = (pv_gp)^* = (\pr g)^*, $$ proving \cite {StarPrepCond}.  \endProof

After developing the necessary tools we will come back to this construction proving that the above process may sometimes
be reversed, so that certain partial representations may be \"{dilated} to unitary group representations.

\medskip

Notice that a partial action may itself be viewed as a *-partial representation of $G$ in the inverse semigroup $\IX $,
according to \cite {PAasPR}.  In what follows we will show that, conversely, partial representations have a
characterization similar to the definition of partial actions given in \cite {DefPA}.

\state Proposition \label characPRViaIneq Let $G$ be a group, $S$ be a unital inverse semigroup, and $\Pr :G\to S$ be a
map.  Then $\Pr $ is a *-partial representation\fn {Here we are considering a generalized notion of partial
representations, taking values in a semigroup with involution, rather than a *-algebra, as already commented after \cite
{DefPR}.  In any case, when we say that $\Pr $ is a *-partial representation, all we mean is that \cite {DefPR/i--iv}
are satisfied.}  if and only if, for all $g,h\in G$, one has that \izitem \zitem $\pr 1 = 1$, \zitem $\pri g = (\pr
g)^*$, \zitem $\pr g\pr h\leq \pr {gh}$.

\Proof Assume conditions (i--iii) hold.  Given $g$ and $h$ in $G$, we interpret (iii) from the point of view of \cite
{DefinePOinISG} obtaining $$ \pr g \pr h = \pr {gh}\pri h \pri g \pr g \pr h.  \equationmark InterprerIneq $$ Focusing
on \cite {ObserveRight}, we right multiply the above equation by $\pri h$ to obtain $$ \pr g \pr h \pri h = \pr {gh}\pri
h \ \pri g \pr g \ \pr h \pri h \$= \pr {gh}\pri h \ \pr h \pri h \ \pri g \pr g = \pr {gh}\pri h \ \pri g \pr g =
\ldots \equationmark InterruptCalc $$ which is almost what we want, except for the extraneous term $\pri g \pr g$ on the
right-hand-side.  Letting $e$ denote the initial projection of $\pr {gh}\pri h$, namely $$ e = \pr h \pr {h\inv g\inv
}\pr {gh}\pri h , $$ we claim that $$ e \leq \pri g \pr g.  \equationmark CompareIniProjs $$

Should the claim be verified, we could continue from \cite {InterruptCalc} as follows: $$ \ldots = \pr {gh}\pri h \ e\
\pri g \pr g = \pr {gh}\pri h \ e = \pr {gh}\pri h, $$ thus proving the desired axiom \cite {ObserveRight}.
Unfortunately the proof we found for \cite {CompareIniProjs} is a bit cumbersome.  It starts by applying \cite
{InterprerIneq} with $gh$ and $h\inv $ playing the roles of $g$ and $h$, respectively, producing $$ \pr {gh} \pr {h\inv
} = \pr g \pr h \pr {h\inv g\inv } \pr {gh} \pri h.  \equationmark InterprerIneqTwo $$

With an eye in \cite {CompareIniProjs} we then compute $$ e \pri g \pr g = \pr h \pr {h\inv g\inv } \pr {gh}\pri h \pri
g \pr g = \ldots $$ Replacing $\pri h$ by $\pri h\pr h\pri h$, and using the fact that idempotents in $S$ commute, the
above equals $$ \ldots = \pr h \pr {h\inv g\inv } \pr {gh}\pri h \pri g \pr g \pr h \pri h = \ldots $$ Repeating this
procedure, but now replacing $\pr {gh}$ by $\pr {gh} \pr {h\inv g\inv }\pr {gh}$, we get $$\def \quad {\kern 4pt}\def
\ccr {\pilar {12pt}\cr } \matrix { \ldots & = & \pr h \pr {h\inv g\inv } \underbrace {\pr {gh} \pri h \pri g \pr g \pr
h}_{(\InterprerIneq )} \pr {h\inv g\inv }\pr {gh}\pri h = \ccr & = & \pr h \pr {h\inv g\inv } \hfill \underbrace {\pr g
\pr h \kern 25pt \pr {h\inv g\inv }\pr {gh}\pri h}_{(\InterprerIneqTwo )} = \ccr & = & \pr h \pr {h\inv g\inv } \hfill
\pr {gh} \pri h = \hfill \ccr & = & e. \hfill \ccr } $$

This proves \cite {CompareIniProjs} and, as already mentioned, \cite {ObserveRight} follows.  Point \cite {ObserveLeft}
is now an easy consequence of \cite {ObserveRight} and (ii).

In order to prove the converse, given that $\Pr $ is a *-partial representation, we have, for all $g$ and $h$ in $G$,
that $$ \pr {gh}\pri h \pri g \pr g \pr h \={ObserveRight} \pr g\ \pr h \pri h \ \pri g \pr g \ \pr h \$= \pr g\ \pri g
\pr g \ \pr h \pri h \ \pr h = \pr g \pr h, $$ proving (iii).  Points (i) and (ii) are evident from the hypothesis.
\endProof

Suggested by \cite {DefineSemiSatPAct} we have the following important concept.

\definition \label DefineSemiSatPRept Let $G$ be a group equipped with a length function $\ell $.  A partial
representation $\Pr $ of $G$ in a unital algebra $B$ is said to be \subjex {semi-saturated}{semi-saturated partial
representation} (with respect to the given length function $\ell $) if $$ \ell (gh)=\ell (g)+\ell (h) \quad
\Longrightarrow \quad \pr g\pr h = \pr {gh} \for g,h\in G.  $$

Among the many interesting algebraic properties of partial representations we note the following:

\state Proposition Let $\Pr $ be a given \mstar partial representation of a group $G$ in a unital \mstar algebra $B$.
Denoting by $$ \e g:= \pr g \pri g, \equationmark DefineEg $$ one has for every $g,h\in G$, that \izitem \zitem $\pr g
\pri g \pr g = \pr g$, \zitemBismark PrgIsPiso \zitem $\e g$ is a (self-adjoint) idempotent, \zitem $\pr g\e h= \e {gh}
\pr g$, \zitemBismark PrepComutRel \zitem $\e g\e h= \e h\e g$, \zitemBismark MainTheEgCommute \zitem $\pr g\pr h= \e
g\pr {gh}$. \zitemBismark MultWithProj

\Proof The first point follows easily from \cite {PrOnOne} and \cite {ObserveRight}.  Point (ii) then follows
immediately from \cite {PrgIsPisoLoc} and the observation, in the *-algebra case, that $$ (\e g)^* = (\pr g \pri g)^* =
(\pri g)^* (\pr g)^* = \pr g \pri g = \e g.  $$ As for (iii) we have $$ \pr g\e h = \pr g\pr h \pri h \={ObserveRight}
\pr {gh} \pri h \={PrgIsPisoLoc} $$$$ = \pr {gh} \pri {(gh)} \pr {gh} \pri h \={ObserveLeft} \pr {gh} \pri {(gh)} \pr g
= \e {gh} \pr g.  $$ To prove \cite {MainTheEgCommuteLoc} we compute $$ \e g \e h = \pr g \pri g \e h
\={PrepComutRelLoc} \pr g \e {g\inv h} \pri g \={PrepComutRelLoc} \e {gg\inv h} \pr g \pri g = \e h \e g.  $$ We finally
have $$ \pr g\pr h \={PrgIsPisoLoc} \pr g \pri g \pr g \pr h \={ObserveLeft} \e g\pr {gh}.  \omitDoubleDollar \endProof

Even though a certain disregard for the group law is the defining feature of partial representations, sometimes the
group law is duly respected:

\state Proposition \label RegardDisregard Let $\Pr $ be a partial representation of the group $G$ in a unital $\K
$-algebra $B$.  Also let $h$ be a fixed element of $G$.  Then the following are equivalent: \izitem \zitem $\pr h$ is
left-invertible, \zitem $\pri h \pr h = 1$, \zitem $\pr g \pr h = \pr {gh}$, for all $g$ in $G$, \zitem $\pri h$ is
right-invertible, \zitem $\pri h \pr g = \pr {h\inv g}$, for all $g$ in $G$.

\Proof \def \implic #1#2{\noindent (#1) $\Rightarrow $ (#2).}\null \implic {i}{ii} Let $v$ be a left-inverse for $\pr
h$, meaning that $v\pr h = 1$.  Then $$ 1 = v\pr h \={PrgIsPiso} v \pr h \pri h \pr h = \pri h \pr h.  $$

\medskip \implic {ii}{iii} For every $g$ in $G$, we have that $$ \pr g \pr h = \pr g \pr h \pri h \pr h \={ObserveRight}
\pr {gh} \pri h \pr h = \pr {gh}.  $$

\medskip \implic {iii}{i} Plugging $g=h\inv $ in (iii), we have $$ \pr g \pr h = \pr {gh} = \pr 1 = 1, $$ so $\pr g$ is
a left-inverse for $\pr h$.

\medskip \implic {ii}{v} For every $g$ in $G$, we have that $$ \pri h \pr g = \pri h \pr h \pri h \pr g \={ObserveLeft}
\pri h \pr h \pr {h\inv g} = \pr {h\inv g}.  $$

\medskip \implic {v}{iv} Taking $g=h$ in (v), we have $$ \pri h \pr g = \pr {h\inv g} = \pr 1 = 1, $$ so $\pr g$ is a
right-inverse for $\pri h$.

\medskip \implic {iv}{ii} Let $v$ be a right-inverse for $\pri h$, meaning that $\pri hv = 1$.  Then $$ 1 = \pri h v =
\pri h \pr h \pri h v= \pri h \pr h.  \omitDoubleDollar \endProof

One of the most important roles played by partial representations is as part of covariant representations of algebraic
partial dynamical systems, a notion to be introduced next.

\definition \label DefineCovarRep By a \subj {covariant representation} of a \mstar algebraic partial action $$ \Th =
\big (\{\D g\}_{g\in G},\ \{\th g\}_{g\in G}\big ) $$ in a unital \mstar algebra $B$, we shall mean a pair $(\pi ,\Pr
)$, where $\pi :A \to B$ is a \mstar homomorphism, $\Pr $ is a \mstar partial representation of $G$ in $B$, and, $$ \pr
g \pi (a) \pri g = \pi \big (\th g(a)\big ) \for g\in G\for a\in \Di g.  $$

Under the assumptions of \cite {PRinCrossProd} we have seen that $\Pr $ is a partial representation of $G$ in $A\rt \Th
G$.  If, in addition, we denote by $\pi $ the canonical inclusion of $A$ in $A\rt \Th G$, one may easily prove that the
pair $(\pi ,\Pr )$ is a covariant representation of $\Th $ in $A\rt \Th G$.  In doing so, the conclusions of \cite
{LotsAFormulas} will make this task a lot easier.

Let us now present a few easy consequences of the above definition.

\state Proposition \label CovarTopRepInAlg Let $$ \Th = \big (\{\D g\}_{g\in G},\ \{\th g\}_{g\in G}\big ) $$ be a
\mstar algebraic partial action of a group $G$ on a \mstar algebra $A$.  Also let $(\pi ,\Pr )$ be a covariant
representation of\/ $\Th $ in a unital \mstar algebra $B$, and denote by $\e g = \pr g \pri g$, as usual.  Then, for
every $g$ in $G$, one has that: \izitem \zitem $\pi (a) = \e g \pi (a) = \pi (a) \e g$, \ for every $a\in \D g$,
\zitemBismark EgAsUnits \zitem $\pr g \pi (a) = \pi \big (\th g(a)\big ) \pr g$, \ for every $a\in \Di g$, \zitemBismark
OtherCovFormula \zitem the linear mapping $\pi \times \Pr : A\rt \Th G \to B$ determined by $$ (\pi \times \Pr )(a\delta
_g) = \pi (a)\pr g \for g\in G\for a\in \D g, $$ is a \mstar homomorphism. \zitemBismark FromCovarToRep

\Proof Given $a$ in $\D g$, we have $$ \e g \pi (a) \e g = \pr g \pri g \pi (a) \pr g \pri g = \pr g \pi \big (\thi
g(a)\big ) \pri g \$= \pi \big (\th g(\thi g(a))\big ) = \pi (a), $$ from where one easily deduces \cite {EgAsUnitsLoc}.
Now, if $a\in \Di g$, then $$ \pr g\pi (a) \={EgAsUnitsLoc} \pr g\pi (a) \ei g = \pr g\pi (a) \pri g \pr g = \pi \big
(\th g(a)\big ) \pr g, $$ proving \cite {OtherCovFormulaLoc}.  As for \cite {FromCovarToRepLoc}, we have for all $a\in
\D g$ and $b\in \D h$, that $$ (\pi \times \Pr )(a\delta _g)\ (\pi \times \Pr )(b\delta _h) = \pi (a)\pr g \pi (b)\pr h
\={OtherCovFormulaLoc} \pr g \pi \big (\thi g(a)\big ) \pi (b)\pr h \$= \pr g \pi \big (\thi g(a)b\big )\pr h = \pi \Big
(\th g\big (\thi g(a)b\big )\Big )\pr g\pr h \={MultWithProj} $$$$ = \pi \Big (\th g\big (\thi g(a)b\big )\Big )\e g\pr
{gh} \={EgAsUnitsLoc} (\pi \times \Pr )\Big (\th g\big (\thi g(a)b\big )\delta _{gh}\Big ) \$= (\pi \times \Pr )(a\delta
_g\ b\delta _h), $$ so we see that $\pi \times \Pr $ is multiplicative.  In the *-algebraic case, we have $$ \big ((\pi
\times \Pr )(a\delta _g)\big )^* = (\pi (a)\pr g)^* = \pr g^* \pi (a)^* \={EgAsUnits} \pri g \pi (a^*)\pr g\pri g \$=
\pi \big (\thi g(a^*)\big )\pri g = (\pi \times \Pr )\big (\thi g(a^*)\delta _{g\inv }\big ) = (\pi \times \Pr )\big
((a\delta _g)^*\big ), $$ proving that $\pi \times \Pr $ is a *-homomorphism.  \endProof

Therefore we see that covariant representations lead to \"{representations of} (simply meaning homomorphisms defined on)
the partial crossed product algebra.  One could then ask to what extent does every homomorphism defined on $A\rt \Th G$
arise from a covariant representation. We will later see that this is always the case in the category of C*-algebras but
for now we simply note that, if every $\D g$ is unital, then a homomorphism $$ \phi : A\rt \Th G \to B $$ into a unital
algebra $B$ gives rise to a covariant representation $(\pi ,\Pr )$ in a trivial way: set $\pi (a) = \phi (a\delta _1)$,
and define $ \pr g = \phi (1_g\delta _g).  $ One would then have for all $a$ in $\D g$, that $$ (\pi \times \Pr
)(a\delta _g) = \pi (a)\pr g = \phi (a\delta _1)\phi (1_g\delta _g) = \phi (a\delta _1 \ 1_g\delta _g) \={LotsAFormulas}
\phi (a\delta _g), $$ hence proving that $\phi =(\pi \times \Pr )$.

\medskip In case one is attempting to provide a concrete model for a given partial crossed product algebra, perhaps
trying to prove it to be isomorphic to some well known algebra, it is useful to know when are homomorphisms defined on
the crossed product injective.  The following result may therefore come in handy.

\state Proposition \label GaugeUniqueness Let $B$ be a $G$-graded \mstar algebra and let $(\pi ,\Pr )$ be a covariant
representation of a \mstar algebraic partial dynamical system $$ \big (A,\ G,\ \{\D g\}_{g\in G},\ \{\th g\}_{g\in
G}\big ), $$ in $B$, such that $\pi (A)\subseteq B_1$, and $\pr g\in B_g$, for all $g\in G$.  Then $\pi \times \Pr $ is
a graded homomorphisms (meaning that it takes the $G$-grading subspace $\D g\delta _g$ into $B_g$), and moreover the
following are equivalent: \izitem \zitem $\pi $ is one-to-one, \zitem $\pi \times \Pr $ is one-to-one.

\Proof For all $a\in \D g$, notice that $$ (\pi \times \Pr )(a\delta _g) = \pi (a_g)\pr g\in B_1B_g \subseteq B_g,
\equationmark IsoIsGraded $$ so $\pi \times \Pr $ is a graded homomorphisms, as desired.

Since the restriction of $\pi \times \Pr $ to the canonical copy of $A$ in $A\rt \Th G$ coincides with $\pi $, it is
obvious that (ii) implies (i).

Conversely, assuming that (i) holds, let $ a = \soma {g\in G}a_g\delta _g \in A\rt \Th G, $ be in the kernel of $\pi
\times \Pr $.  We then have $$ 0 = (\pi \times \Pr )(a)=\soma {g\in G}\pi (a_g)\pr g.  \equationmark WriteInKernel $$

Since each $\pi (a_g)\pr g\in B_g$, by \cite {IsoIsGraded}, and since the $B_g$ are independent, we deduce that $\pi
(a_g)\pr g=0$, for all $g$ in $G$.  Recalling that the $\pi (a_g)\in \D g$, we then have $$ \pi (a_g) \={EgAsUnits} \pi
(a_g) \e g = \pi (a_g) \pr g \pri g = 0.  $$ The assumed injectivity of $\pi $ then implies that all $a_g=0$, so $a=0$,
proving that $\pi \times \Pr $ is injective.  \endProof

\nrem A first rudimentary notion of partial representations was introduced by McClanahan in \ref {McClanahan/1995}, as
part of the ingredients of covariant representations for partial dynamical systems.  A refinement of this idea was
subsequently presented by Quigg and Raeburn in \ref {QuiggRaeburn/1997}, where the expression \"{partial representation}
was coined.  The current set of axioms, as described in \cite {DefPR}, was introduced in \ref {Exel/1998}, under the
presence of an involution, and in \ref {DokuchaevExelPiccione/1999}, otherwise.  The characterization of partial
representations given in \cite {characPRViaIneq} is closely related to the original definition given by Quigg and
Raeburn.

\chapter Partial group algebras

Given a group $G$, one of the main reasons why one studies the classical group algebra $\K (G)$ is that its
representation theory is equivalent to that of $G$.  Having expanded the notion of group representations to include
partial ones, we will now introduce the partial group algebra of $G$ whose role relative to partial representations of
$G$ will be shown to parallel that of $\K (G)$.

\medskip \fix To begin, let us temporarily fix a partial representation $\Pr $ of a given group $G$ in a unital algebra
$B$.  Recall from \cite {MainTheEgCommute} that the $\e g$ (defined to be $\pr g\pri g$) form a commutative set of
idempotents, so the subalgebra $A$ of $B$ generated by all of the $\e g$'s is a commutative algebra.  Denote by $\D g$
the ideal of $A$ generated by each $\e g$, that is $$ \D g = A\e g.  $$ Since $\e g$ is idempotent, we see that $\D g$
is a unital ideal, with $\e g$ playing the role of the unit.  For each $g$ in $G$, let us also consider the map $$ \th
g: \Di g \to \D g, $$ defined by $\th g(a) =\pr g a \pri g$, for all $a$ in $\Di g$.

\state Proposition \label SystemFromPrep Under the above conditions one has that $$ \big (A,\ G,\ \{\D g\}_{g\in G},\
\{\th g\}_{g\in G}\big ) $$ is an algebraic partial dynamical system.

\Proof For $a,b\in \Di g$, notice that $$ \th g(a) \th g(b) = \pr g a \pri g \pr g b \pri g = \pr g a \ei g b \pri g =
\pr g a b \pri g = \th g(ab), $$ so $\Th $ is a homomorphism.  We next need to check condition \cite {PAIdentity}, which
is evident, as well as condition \cite {PAExtend}.  So let us be given $g$ and $h$ in $G$, and let $a$ be in the domain
of $\th g\circ \th h$, which we have seen in \cite {DomainCompos} to be $\th h\inv (\D h \cap \Di g)$.  This means that
$a = a\ei h $, and also that $$ \pr h a \pri h = \th h(a) = \th h(a) \e h \ei g = \pr h a \pri h \e h \ei g \$= \pr h a
\pri h\ei g \={PrepComutRel} \pr h a \e {h\inv g\inv }\pri h = \pr h a \ei {(gh)}\pri h.  $$ It follows that $$ a = \ei
h a \ei h = \pri h \pr h a \pri h \pr h \$= \pri h \pr h a \ei {(gh)}\pri h \pr h = a \ei {(gh)}\ei h \in \Di {(gh)}
\cap \Di h, $$ so we see that $a$ lies in the domain of $\th {gh}$.  Moreover, $$ \th g\big (\th h(a)\big ) = \pr g \pr
h a \pri h \pri g = \pr g \pr h \ei h a \ei h\pri h \pri g = (\star ), $$ but $$ \pr g \pr h \ei h = \pr g \pr h \pri h
\pr h \={ObserveRight} \pr {gh} \pri h \pr h = \pr {gh} \ei h, $$ and similarly $\ei h\pri h \pri g = \ei h\pri {(gh)}$.
Therefore $$ (\star ) = \pr {gh} \ei h a \ei h \pri {(gh)} = \pr {gh} a \pri {(gh)} = \th {gh}(a), $$ proving that $\th
g\circ \th h \subseteq \th {gh}$, as desired.  \endProof

Denoting by $\iota $ the inclusion of $A$ in $B$, it is then evident that $(\iota ,\Pr )$ is a covariant representation
of the above partial dynamical system in $B$, so we may use \cite {FromCovarToRep} to conclude that $$ \iota \times \Pr
: A\rt \Th G \to B \equationmark MapToBzero $$ is a homomorphism.  The range of $\iota \times \Pr $ is clearly the
subalgebra $B_0$ of $B$ generated by the range of $\Pr $, so $B_0$ is a quotient of the partial crossed product.  In
general we cannot assert that $B_0$ is isomorphic to the partial crossed product, but under certain conditions this may
be guaranteed.

\state Theorem \label GradedAlgsAsCrossProd Let $G$ be a group and let $B=\bigoplus _{g\in G} B_g$ be a unital
$G$-graded algebra which is generated by the range of a partial representation $\Pr $ of\/ $G$, such that $ \pr g\in
B_g, $ for all $g$ in $G$.  Then $B_1$ is commutative and there exists a partial action of\/ $G$ on $B_1$ such that $$ B
\simeq B_1\rt \Th G, $$ as graded algebras.

\Proof Let $A$ be the subalgebra of $B$ generated by the $\e g$.  For every $g$ in $G$ one has that $$ \e g = \pr g \pri
g \in B_g B_{g\inv } \subseteq B_1, $$ so $A\subseteq B_1$.  Considering the algebraic partial dynamical system arising
from $\Pr $ as in \cite {SystemFromPrep}, we have that the homomorphism $\iota \times \Pr $ described in \cite
{MapToBzero} is one-to-one by \cite {GaugeUniqueness}.  As already seen, the range of $\iota \times \Pr $ is the algebra
generated by the range of $\Pr $, which is assumed to be $B$, whence $\iota \times \Pr $ is also onto, thus proving that
$$ A\rt \Th G \simeq B.  $$

Being a surjective graded homomorphism it is then clear that $\iota \times \Pr $ satisfies $$ (\iota \times \Pr )(\D
g\delta _g) = B_g \for g\in G, $$ and, in particular $B_1 = (\iota \times \Pr )(A\delta _1)$.  So $B_1$ is a commutative
algebra isomorphic to $A$.  We may then transfer the partial action on $A$ over to $B_1$ via this isomorphism, so that
$$ B_1\rt \Th G \simeq A\rt \Th G \simeq B.  \omitDoubleDollar \endProof

As an example let us return to the grading $\{B_n\}_{n\in \Z }$ of $M_n(\K )$ discussed at the end of chapter \cite
{CrossProdChap}.  Consider the element $$ v = \left [\matrix { 0 & 0 & 0 & \cdots & 0 & 0 \cr 1 & 0 & 0 & \cdots & 0 & 0
\cr 0 & 1 & 0 & \cdots & 0 & 0 \cr \vrule height 20pt depth 14pt width 0pt \vdots & \vdots & \ddots & & \vdots & \vdots
\cr 0 & 0 & 0 & \cdots & 1 & 0 }\right ] \in B_1 \subseteq M_n(\K ), $$ and let $\Pr :\Z \to M_n(\K )$ be given by $$
\pr n = \left \{\matrix { v^n, & \hbox { if } n\geq 0,\cr \pilar {20pt}(v^*)^{|n|}, & \hbox { if } n<0, }\right .  $$
where $v^*$ refers to the transpose of $v$.  It is relatively easy to check that $\Pr $ is a partial representation of
$\Z $ in $M_n(\K )$ satisfying the conditions of \cite {GradedAlgsAsCrossProd} and, since $B_0$ is isomorphic to $\K
^n$, we conclude that $$ M_n(\K ) \simeq \K ^n\rt \Th \Z .  $$

The above partial action of $\Z $ on $\K ^n$ may be shown to be precisely the semi-saturated partial action of $\Z $
given by \cite {ActionsOfFreeGroups} in terms of the partial automorphism $f$ of $\K ^n$ defined by $$ f(\lambda _1,
\lambda _2, \ldots , \lambda _{n-1}, 0) = f(0, \lambda _1, \lambda _2, \ldots , \lambda _{n-1}), $$ whose domain is, of
course, the set of vectors in $\K ^n$ with zero in the last coordinate.

\definition \label DefineKPar Let $G$ be a group.  The \subj {partial group algebra} of $G$, denoted $\Kpar $, is the
universal unital $\K $-algebra generated by a set of symbols $ \{[g] : g\in G\}, $ subject to the relations $$ [1] = 1,
\qquad [g][h][h\inv ] = [gh][h\inv ] \and [g\inv ][g][h] = [g\inv ][gh].  $$ for all $g$ and $h$ in $G$.

The universal property of $\Kpar $ may be phrased as follows:

\state Proposition \label UnivPRep The correspondence $$ g\in G\mapsto [g]\in \Kpar $$ is a partial representation,
which we will call the \subj {universal partial representation}.  In addition, for any partial representation $\Pr $ of
$G$ in a unital $\K $-algebra $B$, there exists a unique homomorphism $$ \phi :\Kpar \to B, $$ such that $\pr g = \phi
([g])$, for all $g\in G$.

In case $\K $ is equipped with a conjugation, as in \cite {KConjuga}, we have the following:

\state Proposition For every group $G$, one has that $\Kpar $ is a *-algebra under a unique involution satisfying $$
[g]^* = [g\inv ] \for g\in G, $$ whence the universal partial representation is a *-representation.

\Proof Let $B= \Kpar ^*$, meaning the conjugate-opposite algebra, in which the multiplication operation is replaced by
$$ x\star y = yx \for x,y\in \Kpar , $$ and the scalar multiplication is replaced by $$ \lambda \cdot x = \bar \lambda x
\for \lambda \in \K \for x\in \Kpar .  $$

One may then easily prove that the mapping $$ \Pr : g\in G \mapsto [g\inv ]\in B $$ is a partial representation, so the
universal property \cite {UnivPRep} implies the existence of a homomorphism $ \sigma : \Kpar \to B $ such that $\sigma
([g]) = [g\inv ]$, for all $g$ in $G$.  The required involution is then defined by $$ a^*:= \sigma (a) \for a\in \Kpar .
\omitDoubleDollar \endProof

It is our next immediate goal to analyze the partial dynamical system arising from the universal partial representation
of $G$, as in \cite {SystemFromPrep}.  Our method will be to first construct an abstract partial dynamical system and
later prove it to be the one we are looking for.  In particular we will be able to describe $\Kpar $ as a partial
crossed product algebra.

The tip we will follow is that the idempotents $\e g = [g][g\inv ]$ in $\Kpar $ are not likely to satisfy any algebraic
relation other than the fact that they commute and $\e 1=1$.

\definition \label DefineApar Let $\Apar $ be the universal unital $\K $-algebra generated by a set of symbols $$ {\cal
E} := \{\varepsilon _g: g\in G\}, $$ subject to the relations stating that the $\varepsilon _g$ are commuting
idempotents, and that $\varepsilon _1=1$.

Notice that the only ingredient of the group structure of $G$ which is relevant for the above definition is the singling
out of the unit element.  This definition would therefore make sense for any set with a distinguished element.

Observing that $\Apar $ is a commutative algebra, consider, for each $g\in G$, the ideal of $\Apar $ generated by
$\varepsilon _g$, namely $$ \D g = \Apar \varepsilon _g.  $$

Noticing that $$ {\cal E}_g := \{\varepsilon _{gh}\varepsilon _g: h\in G\} $$ is a set of commuting idempotents in $\D
g$, and that for $h=1$, the corresponding idempotent there is the unit of $\D g$, we may invoke the universal property
of $\Apar $ to conclude that there exists a homomorphism $ \tth g:\Apar \to \D g, $ such that $$ \tth g(\varepsilon _h)
= \varepsilon _{gh}\varepsilon _g \for h\in G.  $$

\state Proposition \label DynamicalSysOnAPap For each $g$ in $G$, denote by $\th g$ the restriction of\/ $\tth g$ to
$\Di g$.  Then $ \big (\{\D g\}_{g\in G},\ \{\th g\}_{g\in G}\big ) $ is a partial action of\/ $G$ on $\Apar $.

\Proof Initially observe that, for all $h\in G$, one has $$ \th g(\varepsilon _h\varepsilon _{g\inv }) = \tth
g(\varepsilon _h) \tth g(\varepsilon _{g\inv }) = \varepsilon _{gh}\varepsilon _g\ \varepsilon _1\varepsilon _g =
\varepsilon _{gh}\varepsilon _g, $$ which the reader might want to compare with \cite {IdealUnits}.

To see that each $\th g$ is an isomorphism, let $h\in G$, and notice that, using the calculation above, we have $$ \thi
g \big (\th g(\varepsilon _h\varepsilon _{g\inv })\big ) = \thi g \big (\varepsilon _{gh}\varepsilon _g) = \varepsilon
_h\varepsilon _{g\inv }.  $$ Since $\Di g$ is generated, as an algebra, by the set $\{\varepsilon _h\varepsilon _{g\inv
}: h\in G\}$, we conclude that $\thi g \th g$ is the identity on $\Di g$, whence $\thi g$ is the inverse of $\th g$,
showing that $\th g$ is invertible.

It is evident that $\th 1$ is the identity of $\Apar $, so the statement will be proved once we check \cite {PAExtend}.
For this, let $g,h\in G$, and recall from \cite {DomainCompos} that the domain of $\th g\th h$ is given by $$ \th h\inv
(\D h \cap \Di g) = \thi h(\Apar \varepsilon _h\varepsilon _{g\inv }) \$= \Apar \varepsilon _{h\inv }\varepsilon _{h\inv
g\inv } = \Di h \cap \D {h\inv g\inv }, $$ which we then see is a subset of the domain of $\th {gh}$, as needed.  For
any $k\in G$, we have $$ \th g\big (\th h(\varepsilon _k\varepsilon _{h\inv }\varepsilon _{h\inv g\inv })\big ) = \th
g(\varepsilon _{hk}\varepsilon _h\varepsilon _{g\inv })\big ) \$= \varepsilon _{ghk}\varepsilon _{gh}\varepsilon _g =
\th {gh}(\varepsilon _k\varepsilon _{h\inv }\varepsilon _{h\inv g\inv }).  $$

Since the elements considered above generate $\Di h \cap \D {h\inv g\inv }$, as an algebra, we see that $\th g\th h$
coincides with $\th {gh}$ on the domain of the latter, so $\th g\th h\subseteq \th {gh}$, concluding the proof.
\endProof

We may then form the crossed product $\Apar \rt \Th G$.

\state Theorem \label ParGpAlgAsCP For every group $G$, there exists an isomorphism $$ \Phi : \Kpar \to \Apar \rt \Th G,
$$ such that $\Phi ([g]) = \varepsilon _g\delta _g$, for all $g\in G$.

\Proof Being under the hypothesis of \cite {PRinCrossProd}, we have that the map $$ g\in G \mapsto \varepsilon _g\delta
_g\in \Apar \rt \Th G \ $$ is a partial representation of $G$.  By the universal property of $\Kpar $ we then obtain a
homomorphism $$ \Phi : \Kpar \to \Apar \rt \Th G, $$ such that $\phi ([g]) = \varepsilon _g\delta _g$, for all $g\in G$,
and it now suffices to prove that $\Phi $ is bijective.

In order to produce an inverse of $\Phi $, we will build a covariant representation $(\pi ,\Pr )$ of our partial
dynamical system in $\Kpar $.  For $\pi $ we take the homomorphism from $\Apar $ to $\Kpar $ obtained via the universal
property of the former, such that $$ \pi (\varepsilon _g) = [g][g\inv ] \for g\in G, $$ while, for $\Pr $, we take the
universal partial representation, namely $$ \pr g = [g] \for g\in G.  $$

To check that $(\pi ,\Pr )$ is in fact a covariant representation, let $g,h\in G$, and notice that $$ \pr g \pi
(\varepsilon _h\varepsilon _{g\inv }) \pri g = [g][h][h\inv ][g\inv ] \={PrepComutRel} [gh][(gh)\inv ][g][g\inv ] \$=
\pi (\varepsilon _{gh}\varepsilon _g) = \pi (\th g\big (\varepsilon _h\varepsilon _{g\inv })\big ), $$ as desired.  By
\cite {FromCovarToRep} we obtain a homomorphism $\pi \times \Pr $ from $\Apar \rt \Th G$ to $\Kpar $, such that $$ (\pi
\times \Pr )(a\delta _g) = \pi (a)[g] \for g\in G\for a\in \D g.  $$ On the one hand we have, for all $g,h\in G$, that
$$ \Phi \big ((\pi \times \Pr )(\varepsilon _h\varepsilon _g\delta _g)\big ) = \Phi \big ([h][h\inv ][g]\big ) =
(\varepsilon _h\delta _h)(\varepsilon _{h\inv }\delta _{h\inv })(\varepsilon _g\delta _g) = \varepsilon _h\varepsilon
_g\delta _g, $$ from where we conclude that $\Phi \circ (\pi \times \Pr )$ is the identity.  On the other hand $$ (\pi
\times \Pr )\big (\Phi ([g])\big ) = (\pi \times \Pr )(\varepsilon _g\delta _g) = [g][g\inv ][g] = [g], $$ so $(\pi
\times \Pr )\circ \Phi $ is also the identity, whence $\Phi $ is an isomorphism.  \endProof

One may now easily verify that the partial dynamical system provided by \cite {SystemFromPrep} in terms of the universal
partial representation is equivalent to the one given in \cite {DynamicalSysOnAPap}.

\medskip

Before concluding this chapter, let us mention, without proofs, another interesting feature of the partial group
algebra.  Recall that if $G$ is a finite abelian group and $\K $ is the field of complex numbers then the classical
group algebra $\K (G)$ is isomorphic to $\K ^{|G|}$.  In particular, the only feature of $G$ retained by its complex
group algebra is the number of elements in $G$.  When it comes to partial group algebras the situation is completely
different.

\state Theorem {\rm \ref {Corollary 4.5/DokuchaevExelPiccione/1999}} Let $G$ and $H$ be two finite abelian groups and
let $\K $ be an integral domain whose characteristic does not divide $|G|$. If the partial group algebras $\Kpar $ and\/
$\K \underpar (H)$ are isomorphic, then $G$ and $H$ are isomorphic groups.

\nrem The concept of partial group algebra was introduced in \ref {Definition 2.4/DokuchaevExelPiccione/1999}.  It is
the purely algebraic version of the corresponding notion of \"{partial group C*-algebra}, previously introduced in \ref
{Definition 6.4/Exel/1998}.

\chapter C*-algebraic partial dynamical systems

In this chapter we will adapt the construction of the partial crossed product to the category of C*-algebras.  We begin
with a quick review of basic concepts.

\definition \label DefineCstarAlg A \subj {C*-algebra} is a *-algebra $A$ over the field of complex numbers, equipped
with a norm $\Vert \ponto \Vert $, with respect to which it is a Banach space, and such that for all $a$ and $b$ in $A,$
one has that \izitem \zitem $\Vert ab\Vert \leq \Vert a\Vert \Vert b\Vert $, \zitem $\Vert a^*\Vert = \Vert a\Vert $,
\zitem $\Vert a^*a\Vert = \Vert a\Vert ^2$.

There are many references for the basic theory of C*-algebras where the interested reader will find the basic results,
such as \ref {Pedersen/1979}, \ref {Arveson/1981}, \ref {Murphy/1990} and \ref {Davidson/1996}.

Of special relevance to us is Gelfand's Theorem \ref {Theorem 2.1.10/Murphy/1990} which asserts that there is an
equivalence between the category of locally compact Hausdorff (\subjex {LCH}{LCH topological space} for short)
topological spaces, with proper\fn {A map between topological spaces is said to be \subjex {proper}{proper map} if the
inverse image of every compact set is compact.}  continuous maps, on the one hand, and the category of abelian
C*-algebras, with \subjex {\nonDegHomo }{non-degenerate *-homomorphism}\fn {A *-homomorphism $\varphi $ from a
C*-algebra $A$ to another C*-algebra $B$, or perhaps even into the multiplier algebra $\Mult (B)$, is said to be
\"{{\nonDegHomo }} if $B = \clspan {\varphi (A)B}$ (brackets denoting closed linear span).  By taking adjoints, this is
the same as saying that $B = \clspan {B\varphi (A)}$.  } *-homomorphisms, on the other hand.  This equivalence is
implemented by the contravariant functor $$ X \rightsquigarrow C_0(X), $$ where $C_0(X)$ refers to the C*-algebra formed
by all continuous complex valued functions $f$ defined on $X$, vanishing\fn {We say that a map f \subjex {vanishes at
$\infty $}{vanish at $\infty $}, if for every real number $\varepsilon >0,$ the set $\{x\in X:|f(x)|\geq \varepsilon \}$
is compact.}  at $\infty $.

If $X$ and $Y$ are LCH spaces then any proper continuous map $h:X\to Y$ induces a {\nonDegHomo } *-homomorphism $$ \phi
_h: f\in C_0(Y)\mapsto f\circ h\in C_0(X), $$ and conversely, any {\nonDegHomo } *-homomorphism from $C_0(Y)$ to
$C_0(X)$ is induced, as above, by a unique proper continuous map from $X$ to $Y$.  Moreover, $\phi _h$ is an isomorphism
if and only if $h$ is a homeomorphism.

Ideals (always assumed to be norm-closed and two-sided) in C*-algebras are automatically self-adjoint.  Every ideal in a
C*-algebra is both {\nonDegAlg } and idempotent, and hence the conclusion of \cite {Condforextraassoc} holds for them.

If $X$ is a LCH space then there is a one-to-one correspondence between open subsets of $X$ and ideals of $C_0(X)$ given
as follows: to an open set $U\subseteq X$ we attach the ideal given by $$ C_0(U):= \{f\in C_0(X): f = 0 \hbox { \ on \ }
X\setminus U\}.  $$

The reader should be aware that any open set $U\subseteq X$ may also be seen as a LCH space in itself, so this notation
has a potential risk of confusion since, besides the above meaning of $C_0(U),$ one could also think of $C_0(U)$ as the
set of all continuous complex valued functions \"{defined on} $U,$ and vanishing at $\infty $.

However the two meanings of $C_0(U)$ give rise to naturally isomorphic C*-algebras, the isomorphism taking any function
defined on $U$ to its extension to the whole of $X,$ declared zero outside of $U$.  The very slight distinction between
the two interpretations of this notation will fortunately not cause us any problems.

Recall from \cite {DefineParSym} that $\IX $ denotes the inverse semigroup formed by all partial symmetries of a set
$X.$

\definition Given a C*-algebra $A,$ we will say that a partial symmetry $\phi \in {\cal I}(A)$ is a \subj {partial
automorphism} of $A,$ if the domain and range of $\phi $ are closed two-sided ideals of $A,$ and $\phi $ is a
*-isomorphism from its domain to its range.  We will denote by $\IA $ the collection of all partial automorphisms of
$A$.  It is evident that $\IA $ is an inverse sub-semigroup of ${\cal I}(A)$.

Given any partial homeomorphism of a LCH space $X,$ say $h:U\to V$, where $U$ and $V$ are open subsets of $X$, the map
$$ \phi _h: C_0(V)\to C_0(U) $$ is a *-isomorphism between ideals of $C_0(X),$ and hence may be seen as a partial
automorphism of $C_0(X)$.  It follows from what was said above that the correspondence $$ h\in \IXt \ \mapsto \ \phi
_{h\inv }\in {\rm pAut}\big (C_0(X)\big ) \equationmark ComaprePaSpaceAlg $$ is a semigroup isomorphism.

\definition \label DefCstarPa A \subj {C*-algebraic partial action} of the group $G$ on the C*-algebra $A$ is a partial
action $\Th = \big (\{\D g\}_{g\in G}, \{\th g\}_{g\in G}\big )$ on the underlying set $A,$ such that each $\D g$ is a
\"{closed} two-sided ideal of $A,$ and each $\th g$ is a *-isomorphism from $\Di g$ to $\D g$.  By a \subj {C*-algebraic
partial dynamical system} we shall mean a partial dynamical system $$ \big (A,\ G,\ \{\D g\}_{g\in G},\ \{\th g\}_{g\in
G}\big ) $$ where $A$ is a C*-algebra and $\big (\{\D g\}_{g\in G}, \{\th g\}_{g\in G}\big )$ is a C*-algebraic partial
action of $G$ on $A$.

When it is understood that we are working in the category of C*-algebras and there is no chance for confusion we will
drop the adjective \"{C*-algebraic} and simply say \"{partial action} or \"{partial dynamical system}.

As an immediate consequence of \cite {PAasPR} we have:

\state Proposition \label PAasPRCstar Let $G$ be a group, $A$ be a C*-algebra, and $$ \Th :G\to \IA $$ be a map.  Then
$\Th $ is a C*-algebraic partial action of $G$ on $A$ if and only if conditions (i--iv) of \cite {PAasPR} are fulfilled.

Putting together \cite {PAasPRTop}, \cite {PAasPRCstar} and \cite {ComaprePaSpaceAlg}, one concludes:

\state Corollary \label TopoCstarActions If $G$ is a group and $X$ is a LCH space, then \cite {ComaprePaSpaceAlg}
induces a natural equivalence between topological partial actions of\/ $G$ on $X$ and C*-algebraic partial actions of
$G$ on $C_0(X)$.

\fix We now fix, for the time being, a C*-algebraic partial action $$ \Th = \big (\{\D g\}_{g\in G},\ \{\th g\}_{g\in
G}\big ) $$ of a group $G$ on a C*-algebra $A$.

Since $\Th $ is in particular a *-algebraic partial action, we may apply the construction of the crossed product
described in \cite {DefCrossProd} to $\Th $.  However the resulting algebra, which we will temporarily denote by $$
A\art \Th G, $$ will most certainly not be a C*-algebra, so we will modify the construction a bit in order to stay in
the category of C*-algebras.  Meanwhile we observe that $A\art \Th G$ is an associative algebra by \cite
{TheoremOnAssoc}, as well as a *-algebra by \cite {StarInCrossProd}.

\definition A \subj {C*-seminorm} on a complex *-algebra $B$ is a seminorm $p:B \to \R _+$, such that, for all $a,b\in
B,$ one has that \izitem \zitem $p(ab) \leq p(a)p(b)$, \zitem $p(a^*) = p(a)$, \zitem $p(a^*a) = p(a)^2$.

\bigskip

If $B$ is a C*-algebra and $p$ is a C*-seminorm on $B,$ it is well known that $$ p(b) \leq \Vert b\Vert , \equationmark
CStarSeminormDominated $$ for all $b\in B$.

\bigskip

\state Proposition \label CPAdmissible Let $p$ be a C*-seminorm on $A\art \Th G$. Then, for every $ a = \soma {g\in
G}a_g\delta _g $ in $ A\art \Th G, $ one has that $$ p(a) \leq \soma {g\in G}\Vert a_g\Vert .  $$

\Proof Notice that $A\delta _1$ is isomorphic to the C*-algebra $A,$ so by what was said above we have that $ p(a\delta
_1)\leq \Vert a\Vert , $ for all $a\in A$.  We then have that $$ p(a_g\delta _g)^2 = p\big ( (a_g\delta _g) (a_g\delta
_g)^*\big ) \={LotsAFormulas} p(a_ga_g^*\delta _1) \leq \Vert a_ga_g^*\Vert = \Vert a_g\Vert ^2, $$ so the statement
follows from the triangle inequality. \endProof

\bigskip

Let us therefore define a seminorm on $ A\art \Th G, $ by $$ \maxnorm a = {\rm sup} \big \{p(a) : p \hbox { is a
C*-seminorm on } A\art \Th G\big \}.  \equationmark DefineMaxNorm $$ By \cite {CPAdmissible} we see that $\maxnorm a$ is
always finite and it is not hard to see that $ \maxnorm \ponto $ is a C*-seminorm on $A\art \Th G$ (we will later prove
that it is in fact a norm).

\definition \label DefineCStarCP The \subj {C*-algebraic crossed product} of a C*-algebra $A$ by a group $G$ under a
C*-algebraic partial action $\Th = \big (\{\D g\}_{g\in G}, \{\th g\}_{g\in G}\big )$ is the C*-algebra $A\rt \Th G$
obtained by completing $A\art \Th G$ relative to the seminorm $\maxnorm \ponto $ defined above.

The process of completing a semi-normed space involves first modding out the subspace formed by all vectors of zero
length.  However, as already mentioned, $\maxnorm \ponto $ will be shown to be a norm on $A\art \Th G$, so the modding
out part will be seen to be unnecessary.

\definition \label RedefineDelta From now on, for any $a\in \D g$, we will let $a\delta ^{alg}_g$ denote the element of
$A\art \Th G$ we have so far been denoting by $a\delta _g$, while we will reserve the notation $a\delta _g$ for the
canonical image of $a\delta ^{alg}_g$ in $A\rt \Th G$.

As we will be mostly working with the C*-algebraic crossed product, the notation $a\delta ^{alg}_g$ will only rarely be
used in the sequel.

\definition \label DefineStdInclusion We will denote by $$ \iota : A \to A\rt \Th G, $$ the mapping defined by $\iota
(a) = a\delta _1$, for every $a\in A$.

As already mentioned, we will later prove that the natural map from $A\art \Th G$ to $A\rt \Th G$ is injective and
consequently $\iota $ will be seen to be injective as well.

\medskip The following is a useful device in producing *-homomorphisms defined on crossed product algebras:

\state Proposition \label UniversalityOfEnvelop Let $B$ be a C*-algebra and let $$ \varphi _0: A\art \Th G \to B $$ be a
*-homomor\-phism.  Then there exists a unique *-homomorphism $\varphi $ from $A\rt \Th G$ to $B$, such that the diagram
\vskip -0.2cm \null \hfill \beginpicture \setcoordinatesystem units <0.0030truecm, -0.0030truecm> point at 0 0 \put
{$A\art \Th G$} at 000 000 \put {$B$} at 1000 000 \arrow <0.11cm> [0.3,1.2] from 400 0 to 800 000 \put {$\varphi _0$} at
600 -100 \put {$A\rt \Th G$} at 000 700 \arrow <0.11cm> [0.3,1.2] from 0 200 to 0 500 \arrow <0.11cm> [0.3,1.2] from 350
550 to 800 200 \put {$\varphi $} at 650 480 \endpicture \hfill \null \vskip 0.3cm \noindent commutes, where the vertical
arrow is the canonical mapping arising from the completion process.

\Proof It is enough to notice that $p(x) := \Vert \varphi _0(x)\Vert $ defines a C*-seminorm on $A\art \Th G$, which is
therefore bounded by $\maxnorm \ponto $.  Thus $\varphi _0$ is continuous for the latter, and hence extends to the
completion.  \endProof

\nrem As already mentioned, partial actions on C*-algebras were introduced in \ref {Exel/1994a} for the case of the
group of integers, and in \ref {McClanahan/1995} for general groups.  Although not covered by this book, the notion of
\"{continuous} partial actions of \"{topological groups} on C*-algebras, \"{twisted} by a \"{cocycle} or otherwise, has
also been considered \ref {Exel/1997a}.

\chapter Partial isometries

\def \TT {\eubox {T}} \def \SS {\eubox {S}} \def \isol #1{\, #1\, } \def \ini #1{\isol {#1^*#1}} \def \fin #1{\isol
{#1#1^*}} \def \sam {self-adjoint multiplicative}

When working in the category of C*-algebras, we will will often consider *-partial representations $$ u: G \to A
\equationmark SomePartRep $$ of a given group $G$ in a given C*-algebra $A$, according to Definition \cite {DefPR}.  In
other words, the definition of *-partial representations given in \cite {DefPR} needs no further adaptation to the
C*-algebraic case, considering that C*-algebras are special cases of *-algebras.  Incidentally, many *-partial
representations studied in this book will take place in $\Lin (H)$, the C*-algebra of all bounded linear operators on a
Hilbert space $H$.

Given a *-partial representation $\Pr $, as in \cite {SomePartRep}, observe that, by \cite {StarPrepCond} and \cite
{PrgIsPiso}, one has $$ \pr g\pr g^*\pr g = \pr g \for g\in G.  \equationmark PrepFormedByPisos $$

Elements satisfying this equation are crucial for the present work, so we shall dedicate this entire chapter to their
study.  We begin by giving them a well deserved name:

\definition Let $A$ be a *-algebra.  \izitem \zitem An element $s$ in $A$ is said to be a \subj {partial isometry}, if
$ss^*s=s$.  \zitem An element $p$ in $A$ is said to be a \subj {projection}, if $p = p^* = p^2$.

A useful characterization of partial isometries in terms of projections is as follows:

\state Proposition \label CharacPisoViaProj Let $A$ be a C*-algebra and let $s\in A$.  Then the following are
equivalent: \izitem \zitem $s$ is a partial isometry, \zitem $s^*s$ is a projection, \zitem $ss^*$ is a projection.

\Proof The implications (i)$\Rightarrow $(ii \& iii) are evident.  On the other hand, we have $$ (ss^*s-s)^*(ss^*s-s) =
s^*ss^*ss^*s-s^*ss^*s -s^*ss^*s+s^*s \$= (s^*s)^3-2(s^*s)^2 + s^*s.  $$ Thus, if $s^*s$ is a projection, the above
vanishes and hence $$ \Vert ss^*s-s\Vert ^2 = \Vert (ss^*s-s)^*(ss^*s-s)\Vert = 0, $$ showing that $ss^*s=s$, and hence
that $s$ is a partial isometry.  This proves that (ii)$\Rightarrow $(i).  The proof of (iii)$\Rightarrow $(i) is
obtained by replacing $s$ by $s^*$.  \endProof

Given a projection $p$, we have by \cite {DefineCstarAlg/iii} that $$ \Vert p\Vert ^2 = \Vert p^*p\Vert = \Vert p^2\Vert
= \Vert p\Vert , $$ so, unless $p=0$, we have that $\Vert p\Vert =1$.  If $s$ is a partial isometry, then $s^*s$ is a
projection by \cite {CharacPisoViaProj/ii}, so unless $s=0$, we have $$ \Vert s\Vert ^2 = \Vert s^*s\Vert =1.  $$ This
shows that nonzero projections, as well as nonzero partial isometries have norm exactly 1.

\definition Let $A$ be a C*-algebra and let $s$ be a partial isometry in $A$.  Then the projections $s^*s$ and $ss^*$
are called the \subjex {initial}{initial projection} and \subjex {final}{final projection} projections of $s$,
respectively.

Given a Hilbert space $H$, it is well known that a bounded linear operator $s$ in $\Lin (H)$ is a partial isometry if
and only if $s$ is isometric when restricted to the orthogonal complement of its kernel, a space which is known as the
\subj {initial space} of $s$.  On the other hand, the range of $s$ is known as its \subj {final space}.  It is easy to
see that the range of the initial projection of $s$ coincides with its initial space, and similarly for the final
projection and the final space.

If we denote the initial space of $s$ by $H_0$ and its final space by $H_1$, then the effect of applying $s$ to a vector
$\xi \in H$ consists in projecting $\xi $ orthogonally onto $H_0$, followed by the application of an isometric linear
transformation from $H_0$ onto $H_1$ (namely the restriction of $s$ to $H_0$).

The adjoint $s^*$ of a partial isometry is easily seen to be a partial isometry, while the roles of the initial and
final spaces of $s$ are interchanged with those of $s^*$.

We thus see that, given a *-partial representation $\Pr $ of a group $G$ in a C*-algebra, each $\pr g$ is a partial
isometry and the element $$ \e g = \pr g\pri g = \pr g \pr g^*, $$ which has already played an important role, is
nothing but the final projection of $\pr g$.  On the other hand, the initial projection of $\pr g$ is clearly $\ei g$.

\def \Gen {\Lambda }

Let us now discuss some simple facts about partial isometries and projections in C*-algebras.

\state Lemma \label Idempot Let $A$ be a C*-algebra and let $p\in A$ be such that $p^2=p$, and $\Vert p\Vert \leq
1$. Then $p=p^*$.

\Proof By \ref {Theorem 1.7.3/Arveson/1981} we may assume that $A$ is a closed *-subalgebra of operators on some Hilbert
space $H$.  In that case, notice that for every $\xi \in p(H)^\perp $ and every $\lambda \in \R $, we have $$ |1+\lambda
|\ \Vert p(\xi )\Vert = \Vert p(\xi ) + \lambda p(\xi )\Vert = \Vert p\big (\xi +\lambda p(\xi )\big )\Vert \leq \Vert
\xi +\lambda p(\xi )\Vert , $$ so, by Pythagoras Theorem, $$ (1+\lambda )^2 \Vert p(\xi )\Vert ^2 \leq \Vert \xi \Vert
^2 +\lambda ^2 \Vert p(\xi )\Vert ^2, $$ which is easily seen to imply that $$ (1+2\lambda ) \Vert p(\xi )\Vert ^2 \leq
\Vert \xi \Vert ^2, $$ and since $\lambda $ is arbitrary, we must have that $p(\xi )=0$.  This says that $p$ vanishes on
$p(H)^\perp $ and, since $p$ is the identity on $p(H)$, it must coincide with the orthogonal projection onto $p(H)$.
Hence $p=p^*$.  \endProof

\state Lemma \label ComutProj Let $p$ and $q$ be projections in a C*-algebra $A$.  Then $pq$ is idempotent if and only
if $p$ and $q$ commute.

\Proof If $pq$ is idempotent, then, since $\|pq\|\leq 1$, we have, by \cite {Idempot}, that $pq = (pq)^* = qp$.  The
converse is trivial.  \endProof

The product of two partial isometries is not always a partial isometry.  However we have:

\state Proposition \label ProdPIso Let $s$ and $t$ be partial isometries in a C*-algebra $A$.  Then $st$ is a partial
isometry if and only if $s^*s$ and $tt^*$ commute.

\Proof By definition $st$ is a partial isometry if and only if $$ st(st)^*st = st \iff stt^*s^*st = st \$\iff
s^*stt^*s^*stt^* = s^*stt^* \iff (s^*stt^*)^2 = s^*stt^*, $$ which, by \cite {ComutProj}, is equivalent to the
commutativity of $s^*s$ and $tt^*$.  \endProof

If we are given a set $S$ of partial isometries in a C*-algebra $A$, we may always consider the multiplicative semigroup
$\langle S\rangle $ generated by $S$.  However, as seen above, unless there is enough commutativity among range and
source projections, it is likely that $\langle S\rangle $ will include elements which are not partial isometries.

\definition A set $S$ of partial isometries in a C*-algebra $A$ is said to be \subjex {tame}{tame set of partial
isometries} if the multiplicative sub-semigroup of $A$ generated by $S \cup S^*$, henceforth denoted by $\langle S \cup
S^*\rangle $, consists exclusively of partial isometries.

Given a tame set $S$, we then have that $\langle S \cup S^*\rangle $ is a multiplicative semigroup formed by partial
isometries, which is easily seen to be a self-adjoint set.

\state Proposition \label SemigroupIsInverse Let $\SS $ be a {\sam } sub-semigroup of a C*-algebra consisting of partial
isometries.  Then $\SS $ is an inverse semigroup.

\Proof Given $s$ in $\SS $, we must prove that any element $t$ in $\SS $ such that $sts=s$, and $tst=t$, necessarily
coincides $s^*$.

Observe that both $ts$ and $st$ are idempotent elements with norm no bigger than 1.  So $ts$ and $st$ are self-adjoint
by \cite {Idempot}.  Therefore $$ ts = (ts)^* = (tss^*s)^* = (s^*s)^* (ts)^* = s^*sts = s^*s, $$ and $$ st = (st)^* =
(ss^*st)^* = (st)^* (ss^*)^* = stss^* = ss^*.  $$ Hence $$ t = tst = tss^*st = s^*ss^*ss^* = s^*.  \omitDoubleDollar
\endProof

As a consequence of the above result and our discussion just before it, we have:

\state Corollary \label TameIsISG For any tame set $S$ of partial isometries in a C*-algebra, one has that $\langle
S\cup S^*\rangle $ is an inverse semigroup.

Recall from \cite {MainTheEgCommute} that the projections $\e g$ arising from a partial representation always commute.
This is the basis for the close relationship between tame sets of partial isometries and *-partial representations, as
we shall now see.

\state Proposition \label PrepIsTame Let $\Pr $ be a *-partial representation of a group $G$ in a unital C*-algebra $A$.
Then the range of $\Pr $ is a tame set of partial isometries.

\Proof Let $s$ be an element in the multiplicative semigroup generated by the range of $\Pr $ (which is a self-adjoint
set).  Thus $s=\pr {g_1}\cdots \pr {g_n}$, for suitable elements $g_1,\ldots ,g_n$ of $G$.  By induction on $n$ we then
have $$ \def \nlin {\hfill \stake {14pt}\cr } \def \uline #1{\ \underline {\stake {5pt}#1}} \matrix { ss^*s & = & \pr
{g_1}\cdots \pr {g_{n-1}} \e {g_n} \pr {g_{n-1}}^*\cdots \pr {g_1}^* \pr {g_1}\cdots \pr {g_n} \={PrepComutRel} \nlin &
= & \e {g_1\cdots g_n} \uline {\pr {g_1}\cdots \pr {g_{n-1}}} \uline {\pr {g_{n-1}}^*\cdots \pr {g_1}^*} \uline {\pr
{g_1}\cdots \pr {g_{n-1}}} \ \pr {g_n} = \nlin & = & \e {g_1\cdots g_n} \uline {\pr {g_1}\cdots \pr {g_{n-1}}} \pr {g_n}
= \nlin & = & \pr {g_1}\cdots \pr {g_{n-1}} \e {g_n} \pr {g_n} = \nlin & = & \pr {g_1}\cdots \pr {g_{n-1}} \pr {g_n} =
\nlin & = & s.\hfill } $$ So $s$ is a partial isometry.  \endProof

Roughly speaking, the following result includes a converse of the previous one.

\state Proposition \label PisoSet Let $S = \{s_\lambda \}_{\lambda \in \Gen }$ be a family of partial isometries in a
unital C*-algebra $A$, and denote by $\F $ the free group on the index set $\Gen $.  Then the following are equivalent:
\izitem \zitem There exists a semi-saturated partial representation $\Pr $ of\/ $\F $ in $A$, such that $\pr \lambda
=s_\lambda $, for every $\lambda \in \Gen $.  \zitem There exists a partial representation $\Pr $ of\/ $\F $ in $A$,
such that $\pr \lambda =s_\lambda $, for every $\lambda \in \Gen $.  \zitem $S$ is tame.

\Proof (i) $\Rightarrow $ (ii): Obvious.

\medskip \noindent (ii) $\Rightarrow $ (iii): Follows from \cite {PrepIsTame}.

\medskip \noindent (iii) $\Rightarrow $ (i): For all $\lambda \in \Gen $, we define $\pr {\lambda }=s_\lambda $, and
$\pr {\lambda \inv }=s_\lambda ^*$.  If $g= g_1\cdots g_n$, with $g_i\in \Gen \cup \Gen \inv $, is in reduced form, we
put $$ \pr {g}=\pr {g_1}\cdots \pr {g_n}.  $$

This defines a map $\Pr : \F \to A$, which we claim is a semi-saturated partial representation.  Adopting the usual
convention that the reduced form of the unit group element of $\F $ is the empty string, and also that a product
involving zero factors equals one, we see that $\Pr $ satisfies \cite {PrOnOne}.  The easy verification of \cite
{StarPrepCond} is left to the reader.

Before concluding the verification of the remaining axioms in \cite {DefPR}, we observe that the condition for
semi-saturatedness, namely that $ \pr {g}\pr {h} = \pr {gh}, $ whenever $g$ and $h$ satisfy $|gh| = |g| + |h|$, is
evidently satisfied simply because, in this case, the reduced form $gh$ is the juxtaposition of the corresponding
reduced forms of $g$ and $h$.

In order to prove \cite {ObserveRight}, namely $$ \pr {g} \pr {h} \pri {h} = \pr {gh} \pri {h} \for g,h \in \F ,
\equationmark TheAxiomToBeProved $$ we use induction on $|g| + |h|$.  If either $|g|$ or $|h|$ is zero, there is nothing
to prove.  So, assuming that $|g|,|h|\geq 1$, we may write \def \tilg {a} \def \tilh {b} $$ g = \tilg \lambda \and h =
\mu \tilh , $$ where $\tilg , \tilh \in \F $, $\lambda ,\mu \in \Gen \cup \Gen \inv $ and, moreover, $|g|=|\tilg |+1$
and $|h|=|\tilh |+1$.

In case $\lambda \inv \neq \mu $, we have $|gh|=|g| + |h|$, so $\pr {gh} = \pr {g} \pr {h}$, as seen above.  If, on the
other hand, $\lambda \inv = \mu $, we have $$ \pr {g} \pr {h} \pri {h} = \pr {\tilg \lambda } \pr {\lambda \inv \tilh }
\pr {\tilh \inv \lambda } = \pr {\tilg } \pr {\lambda } \pri {\lambda } \pr {\tilh } \pri {\tilh } \pr {\lambda } =
\cdots $$ By (iii) and the induction hypothesis, we conclude that the above equals $$ \cdots = \pr {\tilg } \pr {\tilh }
\pri {\tilh } \pr {\lambda } \pri {\lambda } \pr {\lambda } = \pr {\tilg \tilh } \pri {\tilh } \pr {\lambda } \$= \pr
{gh} \pr {\tilh \inv \lambda } = \pr {gh} \pri {h}, $$ proving \cite {TheAxiomToBeProved}.  As already remarked, \cite
{StarPrepCond} and \cite {ObserveRight} together imply \cite {ObserveLeft}, so $\Pr $ is indeed a semi-saturated partial
representation and the proof is concluded.  \endProof

The study of partial isometries in a C*-algebra resembles the theory of inverse semigroups in the sense that every
\"{tame} set of partial isometries is contained in an inverse semigroup by \cite {TameIsISG}.  Thus, as long as we are
focusing on partial isometries lying in a single tame set, we may apply many of the tools of the theory of inverse
semigroups.

This should be compared to quantum theory in the sense that, when we are working with a set of \"{commuting}
self-adjoint operators, we are allowed to apply results from function theory, since our set of operators generates a
commutative C*-algebra which, by Gelfand's Theorem, is necessarily of the form $C_0(X)$, for some locally compact
Hausdorff space $X$.  On the down side, if our operators do not commute, function theory becomes unavailable and we must
face true quantum phenomena.

In what follows we will develop some elementary results about partial isometries in a C*-algebra which are not always
explicitly required to lie in a tame set.  Since wild (as opposed to tame) sets of partial isometries are very hard to
handle, our guiding principle will be to stay as close as possible to the theory of inverse semigroups.

We begin with a result supporting a subsequent definition of an order relation among partial isometries.  This should be
compared with \cite {DefinePOinISG}.

\state Proposition \label OrderPisos Let $s$ and $t$ be partial isometries in a C*-algebra $A$.  Then the following are
equivalent \izitem \zitem $ts^*s = s$, \zitem $ss^*t = s$.  \medskip \noindent In this case $s^*s\leq t^*t$, and
$ss^*\leq tt^*$.

\Proof Given that (i) holds, we have $$ ts^*s = s = ss^*s = ts^*s(ts^*s)^*ts^*s = ts^*st^*ts^*s.  $$ Multiplying on the
left by $t^*$ gives $$ t^*ts^*s = t^*ts^*st^*ts^*s, $$ so, if we let $p= t^*t$, and $q=s^*s$, we see that $pq = pqpq$,
so $pq$ is idempotent.  It then follows from \cite {ComutProj} that $p$ and $q$ commute.  Consequently $$ t^*s =
t^*ts^*s = pq, $$ so $t^*s$ is self-adjoint, and then $$ t^*s = (t^*s)^* = s^*t.  \eqno
{(\diamond )} $$

From (i) we also deduce that $$ st^* = ts^*st^*, $$ so $st^*$ is self-adjoint as well, and hence $$ st^* = (st^*)^* =
ts^*.   \eqno {(\star )} $$

Focusing on (ii) we then have $$ ss^*t \={osXt} st^*s \={ostX} ts^*s \explica {=}{i} s.  $$

The converse is verified by applying the part that we have already proved for $s':=s^*$ and $t':=t^*$.

With respect to the final sentence in the statement, we have $$ s^*s \explica =i (ts^*s)^*ts^*s = s^*st^*ts^*s =
s^*st^*t, $$ where the last equality is a consequence of the commutativity of $p$ and $q$.  So $s^*s\leq t^*t$, and one
similarly proves that $ss^*\leq tt^*$.  \endProof

The above result should be compared with \cite {DefinePOinISG}, hence motivating the following:

\definition \label DefineOrderPisos Given two partial isometries $s$ and $t$ in a C*-algebra, we will say that $s$ is
\"{dominated} by $t$, or that $s\preceq t$, if the equivalent conditions of \cite {OrderPisos} are satisfied.

It is elementary to verify that ``$\preceq $'' is a reflexive and antisymmetric relation.  We also have:

\state Proposition The order ``$\preceq $'' defined above is transitive.

\Proof Suppose that $r$, $s$ and $t$ are partial isometries with $r\preceq s\preceq t$.  Then $r = sr^*r$ and $s =
ts^*s$.  Therefore $$ tr^*r = t(sr^*r)^*sr^*r = tr^*rs^*sr^*r \={OrderPisos} ts^*sr^*r = sr^*r = r.  $$ So $r\preceq t$.
\endProof

The following provides a useful alternative characterization of the order relation ``$\preceq $''.

\state Proposition \label AltCharOrderPisos Given two partial isometries $s$ and $t$ in a C*-algebra, the following are
equivalent: \izitem \zitem $s\preceq t$, \zitem $ts^* = ss^*$, \zitem $s^*t = s^*s$.

\Proof (i)$\Rightarrow $(ii): $$ ts^* = ts^*ss^*= ss^*.  $$

\medskip \noindent (ii)$\Rightarrow $(i): $$ ts^*s = ss^*s = s.  $$ The proof that (i)$\Leftrightarrow $(iii) follows
along similar lines.  \endProof

Let us now prove invariance of ``$\preceq $'' under multiplication:

\state Lemma \label PrecInvarMult Let $s_1$, $s_2$, $t_1$ and $t_2$ be partial isometries in a C*-algebra, such that
$s_1\preceq s_2$, and $t_1\preceq t_2$.  If $s_i^*s_i$ commutes with $t_it_i^*$, for all $i=1,2$, then $s_1t_1\preceq
s_2t_2$.

\Proof Recall that each $s_it_i$ is a partial isometry by \cite {ProdPIso}.  We have $$ s_2t_2 (s_1t_1)^* = s_2t_2
t_1^*s_1^* \={AltCharOrderPisos/ii} s_2t_1 t_1^*s_1^*s_1s_1^* \$= s_2s_1^*s_1t_1 t_1^*s_1^* = s_1t_1 t_1^*s_1^* =
s_1t_1(s_1t_1)^*.  $$ This verifies \cite {AltCharOrderPisos/ii}, so $s_1t_1\preceq s_2t_2$, as desired.  \endProof

If a net $\{s_i\}_i$ of partial isometries on a Hilbert space strongly converges to a partial isometry $s$, then the
final projections of the $s_i$ might not converge to the final projection of $s$.  An example of this is obtained by
taking $s_n = (s^*)^n$, where $s$ is the unilateral shift on $\ell ^2({\bf N})$.  In this case the $s_n$ converge
strongly to zero, but the final projections of the $s_n$ all coincide with the identity operator.  On the bright side we
have:

\state Lemma \label IncreasingPisos Let $\{s_i\}_{i\in I}$ be an \underbar {increasing} net of partial isometries on a
Hilbert space $H$.  Then \izitem \zitem $\{s_i\}_{i\in I}$ strongly converges to a partial isometry $s$, \zitem
$\{s_i^*\}_{i\in I}$ strongly converges to $s^*$, \zitem $\{s_i^*s_i\}_{i\in I}$ strongly converges to $s^*s$, \zitem
$\{s_is_i^*\}_{i\in I}$ strongly converges to $ss^*$.

\Proof We will first show that $$ \exists \, \lim _is_i(\xi ), \equationmark MyLimClaim $$ for all $\xi $ in $H$.

By the last sentence of \cite {OrderPisos} we have that the corresponding initial projections, say $$ e_i = s_i^*s_i, $$
form an increasing net in the usual order of projections.  Letting $H_i$ be the range of $e_i$, also known as the
initial space of $s_i$, we then see that the $H_i$ form an increasing family of subspaces of $H,$ so $$ H_0:= \med
\bigcup _{i\in I} H_i $$ is a linear subspace of $H$.  We will next prove \cite {MyLimClaim} for every $\xi $ in $H_0$.
Given such $\xi $, let $i$ be such that $\xi \in H_i$.  Consequently, for all $j\geq i$, we have $$ s_j(\xi ) =
s_je_i(\xi ) = s_js_i^*s_i(\xi ) = s_i(\xi ), $$ so we see that the net $\{s_j(\xi )\}_j$ is eventually constant, hence
convergent.  Observing that our net is uniformly bounded, we then have that \cite {MyLimClaim} also holds for all $\xi $
in the closure of $H_0$.  On the other hand, if $\xi \in H_0^\perp $, then $e_i(\xi )=0$, for all $i$, hence $$ s_i(\xi
) = s_ie_i(\xi ) = 0, $$ so \cite {MyLimClaim} is verified for $\xi $ in $H_0^\perp $ as well, hence also for all $\xi $
in $H$. So we may define $$ s(\xi ) = \lim _is_i(\xi ) \for \xi \in H, $$ and it is easy to see that $s$ is isometric on
$\overline H_0$, while $s$ vanishes on $H_0^\perp $, so $s$ is a partial isometry. This proves (i).

The initial space of $s$ is easily seen to coincide with $\overline H_0$, whence $s^*s$ is the orthogonal projection
onto $\overline H_0$.  On the other hand, it is clear that the orthogonal projection onto $\overline H_0$ is the strong
limit of the $e_i$, so $$ s^*s = \lim _ie_i = \lim _is_i^*s_i, $$ hence (iii) follows.

The order ``$\preceq $'' being evidently invariant under conjugation, we have that $\{s_i^*\}_{i\in I}$ is an increasing
net of partial isometries, hence by the above reasoning it strongly converges to some partial isometry $t$, and moreover
the corresponding net of initial projections $\{s_is_i^*\}_{i\in I}$ (sic) converges to $t^*t$.

Unfortunately the operation of conjugation is not strongly continuous, but it is well known to be weakly continuous,
hence $s_i^*\convrg i s^*$ weakly.  Since the weak limit is unique, we deduce that $t=s^*$, from where (ii) and (iv)
follow.  \endProof

In the following we present another important relation involving partial isometries.

\definition \label DefCompat Let $A$ be a C*-algebra.  Given two partial isometries $s$ and $t$ in $A$, we will say that
$s$ and $t$ are \subjex {compatible}{compatible partial isometries}, if $$ st^*t=ts^*s \and tt^*s = ss^*t.  $$ If $S$ is
a subset of $A$ consisting of partial isometries, we will say that $S$ is a \"{compatible set} when the elements in $S$
are pairwise compatible.

We will soon give a geometric interpretation of this concept, but let us first prove a useful result.

\state Proposition \label CondForCompat Given partial isometries $s$ and $t$ in a C*-algebra $A$, the following are
equivalent: \iaitem \aitem $s$ and $t$ are compatible, \aitem $st^*$ and $s^*t$ are positive elements of $A$.  \medskip
\noindent In this case one also has that \izitem \zitem the final projections $ss^*$ and $tt^*$ commute, \zitem the
initial projections $s^*s$ and $t^*t$ commute, \zitem $st^* = ts^* = ss^* tt^*$, \zitem $s^*t = t^*s = s^*s t^*t$.
\zitem $ts^*s = tt^*s$, and consequently all of the four terms involved in the definition of compatibility coincide.

\Proof Assuming that $s$ and $t$ are compatible, we have $$ st^* = st^*tt^* = ts^*st^* \geq 0, $$ and $$ s^*t = s^*ss^*t
= s^*tt^*s \geq 0, $$ proving (b).  Conversely, assuming (b), observe that, since positive elements are necessarily
self-adjoint, we have $$ st^* = (st^*)^* = ts^*,  \eqno {(\star )} $$ and $$ s^*t = (s^*t)^*
= t^*s,  \eqno {(\diamond )} $$ thus verifying the first identities in (iii) and (iv).
Therefore $$ tt^*ss^* \={sXt} ts^*ts^* \={stX} st^*st^* \={sXt} ss^*tt^*, $$ and $$ t^*ts^*s \={stX} t^*st^*s \={sXt}
s^*ts^*t \={stX} s^*st^*t, $$ proving (i) and (ii).

Repeating an earlier calculation we have $$ (st^*)^2 = st^*st^* \={sXt} ss^*tt^*.  \eqno {(\dagger )} $$ This implies
that $(st^*)^2$ is a projection, whose spectrum is therefore contained in $\{0,1\}$.  By the Spectral Mapping Theorem we
have $$ \big (\sigma (st^*)\big )^2 = \sigma \big ((st^*)^2\big )\subseteq \{0,1\}, $$ so the spectrum of $st^*$ is
contained in $\{-1,0,1\}$, but since $st^*$ is assumed to be positive, its spectrum must in fact be a subset of
$\{0,1\}$.  Consequently $st^*$ is a projection, and then the last identity in (iii) follows from $(\dagger )$.  On the
other hand, $$ (s^*t)^2 = s^*ts^*t \={stX} s^*st^*t, $$ so $(s^*t)^2$ is also a projection and the same reasoning
adopted above leads to the proof of the last identity in (iv).

We may now prove the first of the two conditions in \cite {DefCompat}, namely $$ ts^*s \={stX} st^*s \={sXt} ss^*t =
ss^*tt^*t \={sXt} st^*st^*t = st^*t.  $$

The second condition in \cite {DefCompat} is proved in a similar way, so we deduce that $s$ and $t$ are compatible, as
desired.  Finally, with respect to (v), we have $$ ts^*s \explica ={iii} ss^*tt^*s \explica ={i} tt^*ss^*s = tt^*s.
\omitDoubleDollar \endProof

\nostate \label GeoInterpr As promised, let us give a geometric interpretation for the notion of compatibility of
partial isometries.  For this let $S$ and $T$ be compatible partially isometric linear operators on a Hilbert space $H$.
By \cite {CondForCompat} we have that the initial projections of $S$ and $T$ commute, so we may decompose $H$ as an
orthogonal direct sum $$ H = K\oplus H_S\oplus H_T\oplus L, $$ such that $K\oplus H_S$ is the initial space of $S$, and
$K\oplus H_T$ is the initial space of $T$.  $K$ is therefore the intersection of the initial spaces of $S$ and $T$, and
the orthogonal projection onto $K$ is thus the product of the initial projections of $S$ and $T$, namely $S^*S T^*T$.
If $k$ is a vector in $K$ we then have that $$ S(k) = ST^*T(k) \={DefCompat} TS^*S(k) = T(k), $$ so $S=T$ on $K$.  We
may then define an operator $S\vee T$ on $H$ by $$ \def \quad {\kern 4pt}\pilar {20pt} \stake {16pt} \matrix { (S\vee
T)(k,x_S,x_T,l) & = & S(k) + S(x_S) + T(x_T)\hfill \cr \pilar {18pt} & = & T(k) + S(x_S) + T(x_T),} $$ for all $k$ in
$K$, $x_S$ in $H_S$, $x_T$ in $H_T$, and $l$ in $L$, and it may be proved that $S\vee T$ is a partial isometry which
coincides with $S$ on the initial space of $S$, and with $T$ on the initial space of $T$.

In our next result we will generalize this idea for partial isometries in any C*-algebra.

\state Proposition \label JoinPisos Let $A$ be a C*-algebra and let $s$ and $t$ be partial isometries in $A$.  Then the
following are equivalent \iaitem \aitem $s$ and $t$ are compatible \aitem $s^*s$ and $t^*t$ commute, $ss^*$ and $tt^*$
commute, and there exists a partial isometry dominating both $s$ and $t$.  \medskip \noindent In this case, defining
$$\def \quad {\kern 4pt}\pilar {20pt} \stake {16pt} \matrix {u &=& s+t-st^*t \hfill \cr \pilar {12pt} &=& s+t-ts^*s,\cr
} $$ we have that $u$ is a partial isometry such that: \izitem \zitem $s\preceq u$, and $t\preceq u$, \zitem if $v$ is a
partial isometry such that $s\preceq v$, and $t\preceq v$, then $u\preceq v$, \zitem the initial projection of $u$
coincides with\fn {If $p$ and $q$ are commuting projections in an algebra $A$, one denotes by $p\vee q = p+q-pq$.  It is
well known that $p\vee q$ is again a projection, which is the least upper bound of $p$ and $q$ among the projections in
$A$.}  $\ s^*s \vee t^*t$, \zitem the final projection of $u$ coincides with \ $ss^*\!\vee tt^*$.

\Proof Supposing that $s$ and $t$ are compatible, we compute $$ u^*u = (s^*+t^*-t^*ts^*)(s+t-st^*t) \$=
s^*s+s^*t-s^*st^*t+ t^*s+t^*t-t^*st^*t- t^*ts^*s-t^*ts^*t+t^*ts^*st^*t \$= s^*s+s^*t-s^*st^*t+ t^*s+t^*t-s^*tt^*t
\phantom {\ -t^*ts^*s} -t^*tt^*s \phantom {\ \ +t^*ts^*st^*t}\$= s^*s\phantom {\ +s^*t}-s^*st^*t\phantom {\ \ +
t^*s}+t^*t\phantom {\ -s^*tt^*t} \phantom {\ -t^*ts^*s} \phantom {\ -t^*tt^*s} \phantom {\ \ +t^*ts^*st^*t}\$= s^*s\vee
t^*t, \equationmark SourceU $$ where we have used \cite {CondForCompat} to conclude that $s^*s$ and $t^*t$ commute.
This proves that $u^*u$ is a projection, and by \cite {CharacPisoViaProj} we deduce that $u$ is a partial isometry.
This also proves (iii), and the proof of (iv) is done along similar lines.  In order to prove (i) we compute $$ us^*s =
(s+t-ts^*s)s^*s = ss^*s+ts^*s-ts^*s = s, $$ while $$ ut^*t = (s+t-st^*t)t^*t = st^*t + tt^*t -st^*t =t, $$ so
$s,t\preceq u$.  Notice that this also proves that (a)$\Rightarrow $(b).

Next suppose that $v$ is a partial isometry dominating $s$ and $t$.  Then $$ vu^*u \={SourceU} v(s^*s\vee t^*t) = vs^*s
+ vt^*t - vs^*st^*t = s + t - st^*t = u, $$ so $u\preceq v$, taking care of (ii).

In order to prove that (b)$\Rightarrow $(a), assume that the commutativity conditions in (b) hold, and let $v$ be a
partial isometry dominating both $s$ and $t$. Then $$ st^*t = us^*st^*t = ut^*t s^*s = ts^*s, $$ while $$ tt^*s =
tt^*ss^*v = ss^*tt^*v = ss^*t.  \omitDoubleDollar \endProof

The result above shows that, when two partial isometries are compatible, their least upper bound exists. This motivates
and justifies the introduction of the following notation:

\definition Given a set $S$ of partial isometries in a C*-algebra $A$, suppose that there exists a partial isometry $u$
that dominates every element of $S$, and such that $u\preceq v$, for every other partial isometry $v$ dominating all
elements of $S$.  Observing that such a $u$ is necessarily unique, we denote it by $ \mx S.  $ If $S$ is a two-element
set, say $S=\{s,t\}$, and $\mx S$ exists, we also denote $\mx S$ by $s\vee t$.

Of course this is nothing but the usual notion of least upper bounds, meaningful in any ordered set.  We have spelled it
out just for emphasis.

Given two compatible partial isometries $s$ and $t$ in a C*-algebra $A$, observe that \cite {JoinPisos} implies that
$s\vee t$ exists, and moreover $$\def \quad {\kern 4pt}\pilar {20pt} \stake {16pt} \matrix {s\vee t &=& s+t-st^*t \hfill
\cr \pilar {12pt} &=& s+t-ts^*s.\cr } $$

Projections being special cases of partial isometries, the above notions may also be applied to the former.  Given two
projections $p$ and $q$ in a C*-algebra, notice that $p\leq q$ if and only if $p\preceq q$.  Also, $p$ and $q$ are
compatible if and only if they commute.  In this case the two meanings of the expression $p\vee q$ so far defined are
easily seen to coincide.

\medskip Here is a sort of a distributivity property mixing compatibility of partial isometries and the notion of least
upper bounds just mentioned:

\state Proposition \label KindOfTrans Let $r$, $s$ and $t$ be compatible partial isometries in a C*-algebra.  Then $r$
is compatible with $s\vee t$.

\Proof We have $$ r(s\vee t)^*(s\vee t) \={JoinPisos/iii} r(s^*s\vee t^*t) = rs^*s + rt^*t - rs^*st^*t \$= sr^*r + tr^*r
- sr^*rt^*t = (s + t - st^*t)r^*r = (s\vee t)r^*r.  $$ On the other hand $$ (s\vee t)(s\vee t)^*r \={JoinPisos/iv}
(ss^*\vee tt^*)r = ss^*r + tt^*r - ss^*tt^*r \$= rr^*s + rr^*t - ss^*rr^*t = rr^*(s + t - ss^*t) = rr^*(s + t - st^*s)
\$= rr^*(s + t - ts^*s) = rr^*(s\vee t).  \omitDoubleDollar \endProof

\state Proposition \label ExistMaxCompat Let $S$ be a compatible set of partial isometries in a C*-algebra $A$.  In
addition we assume that either $A$ is a von Neumann algebra, or that $S$ is finite.  Then \izitem \zitem $\mx S$ exists,
\zitem the initial \resp {final} projection of\/ $\mx S$ coincides with the least upper bound of the initial \resp
{final} projections of the members of $S$, \zitem every partial isometry in $A$ which is compatible with the members of
$S$, is also compatible with $\mx S$.

\Proof \def \uf {_{\scriptscriptstyle F}} Let us first deal with the case of a finite $S$, and so we assume that
$S=\{s_1,s_2,\ldots ,s_n\}$.  Our proof will be by induction on $n$.

In case $n=0$, then $S$ is the empty set and $\mx S=0$.  If $n=1$, then $\mx S=s_1$, and in both cases (ii--iii) hold
trivially.

In case $n=2$, both the existence of $\mx S$ and (ii) follow from $\cite {JoinPisos}$, while (iii) follows from \cite
{KindOfTrans}.

Assuming now that $n>2$, let $$ S'=\{s_1,s_2,\ldots ,s_{n-1}\}.  $$

By the induction hypothesis we have that $\mx S'$ exists, and it is compatible with every partial isometry $r$, which in
turn is compatible with the members of $S'$.  This evidently includes $s_n$.  We then apply the already verified case
$n=2$ to the set $$ S'' = \{\mx S', s_n\}.  $$ It is then elementary to prove that $$ \mx S'' = (\mx S')\vee s_n = \mx
S.  $$ By induction we have that $$ (\mx S')^* (\mx S') = \mx \{s^*s: s\in S'\},$$ and $$(\mx S') (\mx S')^* = \mx
\{ss^*: s\in S'\}, $$ so (ii) follows from \cite {JoinPisos/iii\&iv}

If the partial isometry $r$ is compatible with all of the members of $S$, then it is also compatible with $\mx S'$, by
induction, and obviously also with $s_n$.  Hence $r$ is compatible with $\mx S$ by \cite {KindOfTrans}.

We next tackle the case that $S$ is infinite, assuming that $A$ is a von Neumann algebra of operators on a Hilbert space
$H$.  For each finite set $F\subseteq S$, let us denote by $$ s\uf =\mx F.  $$ We may then view $\{s\uf \}\uf $ as an
increasing net, indexed by the set of all finite subsets of $S$, ordered by inclusion. By \cite {IncreasingPisos} the
strong limit of this net exists, and we shall denote it by $u$.

Given $s\in S$, and a finite set $F\subseteq S$, with $s\in F$, we have that $s\preceq s\uf $, whence, for all $\xi $ in
$H$, we have $$ s(\xi ) = s\uf s^*s(\xi ).  $$ Taking the limit as $F\to \infty $, we deduce that $ s(\xi ) = us^*s(\xi
), $ so $s\preceq u$, hence $u$ is an upper bound for $S$.

In order to prove that $u$ is in fact the least upper bound, let $v$ be an upper bound for $S$.  In particular $v$ is
also an upper bound for all finite $F\subseteq S$, whence $s\uf \preceq v$, which is to say that for all $\xi $ in $H$
one has $$ s\uf (\xi ) = v s\uf ^*s\uf (\xi ).  $$ Taking the limit as $F\to \infty $, we have $$ u(\xi ) = \lim \uf
s\uf (\xi ) = \lim \uf vs\uf ^*s_f(\xi ) \={IncreasingPisos/iii} vu^*u(\xi ), $$ so $u\preceq v$.

Point (ii) follows easily from the first part of the proof and \cite {IncreasingPisos/iii\&iv}.  With respect to (iii),
if the partial isometry $r$ is compatible with all of the members of $S$, then $r$ is compatible with $s\uf $, for all
finite $F\subseteq S$, by the first part of the proof.  Thus $$ s\uf r^*r=rs\uf ^*s\uf \and rr^*s\uf = s\uf s\uf ^*r.
$$ Again taking the limit as $F\to \infty $, we deduce from \cite {IncreasingPisos} that $$ ur^*r=ru^*u \and rr^*u =
uu^*r, $$ so $u$ is compatible with $r$.  \endProof

While the ideas of the above proof are still hot, let us observe that we have also proven:

\state Proposition \label MaxIsLimit Let $S$ be a compatible set of partial isometries in a von Neumann algebra $A$.
Then $\mx S$ is the strong limit of the net $\{\mx F\}_F$, where $F$ ranges in the directed set consisting of all finite
subsets of $S$.

\state Proposition \label PrdMx Let $A$ be a C*-algebra \resp {von-Neumann algebra} and let $S$ and $T$ be compatible
sets of partial isometries in $A$ such that $s^*s$ commutes with $tt^*$ for every $s$ in $S$, and every $t$ in $T$.  In
the C*-case we moreover assume that $S$ and $T$ are finite.  Then $$ ST := \{st: s\in S, \ t\in T\} $$ is a compatible
set of partial isometries and $$ \mx (ST) = (\mx S)(\mx T).  $$

\Proof That $ST$ consists of partial isometries is a consequence of \cite {ProdPIso}.

We will next prove that $ST$ is a compatible set. For this choose any pair of elements in $ST$, say $s_1t_1$ and
$s_2t_2$, where $s_1,s_2\in S$, and $t_1,t_2\in T$.  Notice that the set $$ \{s_1^*s_1,\ s_2^*s_2,\ t_1t_1^*,\
t_2t_2^*\} $$ is commutative because the $s_i^*s_i$ commute with the $t_jt_j^*$ by hypothesis, while the $t_jt_j^*$
commute amongst themselves by \cite {CondForCompat/i}, and the $s_i^*s_i$ commute with each other by \cite
{CondForCompat/ii}.  We then have $$ s_1t_1(s_2t_2)^* = s_1t_1t_2^*s_2^* \={CondForCompat/iii} s_1\fin {t_1}\fin
{t_2}\ini {s_2}s_2^* \$= s_1\ini {s_2}\fin {t_1}\fin {t_2}s_2^* = s_2\ini {s_1}\fin {t_1}\fin {t_2}s_2^* \geq 0.  $$ We
also have $$ (s_1t_1)^*s_2t_2 = t_1^* s_1^* s_2 t_2 \={CondForCompat/iv} t_1^*\fin {t_1} \ini {s_1}\ini {s_2} t_2 \$=
t_1^*\ini {s_1}\ini {s_2} \fin {t_1} t_2 = t_1^*\ini {s_1}\ini {s_2} \fin {t_2} t_1 \geq 0.  $$ By \cite {CondForCompat}
we then have that $s_1t_1$ is compatible with $s_2t_2$, proving that $ST$ is indeed a compatible set.  Employing \cite
{ExistMaxCompat/ii} we have $$ \ini {(\mx S)} = \mx \{\ini s: s\in S\} \and \fin {(\mx T)} = \mx \{\fin t: t\in T\}, $$
and hence the initial projection of $\mx S$ commutes with the final projection of $\mx T$, so we deduce from \cite
{ProdPIso} that $(\mx S)(\mx T)$ is a partial isometry.

We will next prove that $$ \mx (ST) = (\mx S)(\mx T), \equationmark MaxProdFinite $$ under the assumption that $S$ and
$T$ are finite sets.  We begin by analyzing several possibilities for the number of elements in $S$ and $T$.  For
example, when $S=\{s_1,s_2\}$, and $T=\{t\}$, we have $$ \mx (ST) = \mx \{s_1t,s_2t\} = (s_1t)\vee (s_2t) = s_1t + s_2t
- s_1t(s_2t)^*s_2t \$= s_1t + s_2t - s_1tt^*s_2^*s_2t = s_1t + s_2t - s_1s_2^*s_2t \$= (s_1 + s_2 - s_1s_2^*s_2)t = (s_1
\vee s_2)t = (\mx S)(\mx T).  $$

By induction it then easily follow follows that \cite {MaxProdFinite} holds when $S$ has an arbitrary finite number of
elements and $T$ a singleton.  Suppose now that $S=\{s\}$, and $T=\{t_1,t_2\}$.  Then $$ \mx (ST) = (st_1)\vee (st_2)
\={CondForCompat/v} st_1 + st_2 - st_1(st_1)^*st_2 \$= st_1 + st_2 - st_1t_1^*s^*st_2 = st_1 + st_2 - st_1t_1^*t_2 \$=
s(t_1 + t_2 - t_1t_1^*t_2) = s(t_1 \vee t_2) = (\mx S)(\mx T).  $$ Again by induction it follows that \cite
{MaxProdFinite} holds for $S$ a singleton and $T$ finite.

We will now prove \cite {MaxProdFinite} in the general finite case using induction on the number of elements of $T$.
Writing $ T=\{t_1,t_2,\ldots ,t_n\}, $ consider the set $T' =\{t_2,\ldots ,t_n\}$, and observe that $$ (\mx S) (\mx T) =
(\mx S) \big (t_1 \vee (\mx T')\big ) = (\mx S)t_1 \ \vee \ (\mx S)(\mx T') \$= (\mx St_1) \ \vee \ (\mx ST') = \mx
(St_1\cup ST') = \mx ST.  $$ Finally, let us tackle the general infinite case, assuming that $A$ is a von-Neumann
algebra.  For this recall from \cite {MaxIsLimit} that $$ \mx S = \lim _{F\to \infty } \mx F \and \mx T = \lim _{G\to
\infty } \mx G, $$ where $F$ and $G$ range in the directed sets consisting of all finite subsets of $S$ and $T$,
respectively.  Since multiplication is doubly continuous on bounded sets for the the strong operator topology, we have
$$ (\mx S) (\mx T) = (\lim _{F\to \infty } \mx F)(\lim _{G\to \infty } \mx G) \$= \lim _{F,G\to \infty } (\mx F)(\mx G)
= \lim _{F,G\to \infty } (\mx FG).  \equationmark msSmxT $$

Again by \cite {MaxIsLimit} we have that $\mx (ST)$ is the strong limit of the net $$ \{\mx H\}_{\buildrel H\subseteq ST
\over {\scriptscriptstyle H \hbox { \sixrm finite}}}, \equationmark NetInProd $$ \medskip \noindent and we observe that
the collection formed by all $H=FG$, where $F$ and $G$ are finite subsets of $S$ and $T$, respectively, is co-final in
the collection of finite subsets of $ST$.  This gives rise to a subnet of \cite {NetInProd}, which therefore also
converges to $\mx (ST)$.  In other words, the last limit in \cite {msSmxT} coincides with $\mx (ST)$, thus proving \cite
{MaxProdFinite} in the general case.  \endProof

\definition Let $S$ be a set of partial isometries in a C*-algebra $A$.  \izitem \zitem We will say that $S$ is \subj
{finitely-$\vee $-closed} if, whenever $T$ is a \"{finite} compatible subset of $S$, one has that $\mx T$ lies in $S$.
\zitem In case $A$ is a von-Neumann algebra, we will say that $S$ is \subjex {$\vee $-closed}{closed under $\vee $, or
$\vee $-closed} if the condition above holds for every compatible subset $T\subseteq S$, regardless of whether it is
finite of infinite.

The set of all partial isometries in a C*-algebra $A$ is evidently finitely-$\vee $-closed.  Moreover, the intersection
of an arbitrary family of finitely-$\vee $-closed sets is again finitely-$\vee $-closed.  Thus we may speak of the \subj
{finite-$\vee $-closure} of a given set $S$ of partial isometries in $A$, namely the intersection of all finitely-$\vee
$-closed sets of partial isometries containing $S$.

All that was said in the above paragraph clearly remains true if we remove all references to finiteness, as long as we
take $A$ to be a von-Neumann algebra.  In particular we may also speak of the \subjex {$\vee $-closure}{closure with
respect to $\vee $} of a set of partial isometries in a von-Neumann algebra.

\state Lemma \label ClosureISG Let $A$ be a C*-algebra \resp {von-Neumann algebra} and let $S$ be a {\sam }
sub-semigroup of $A$ consisting of partial isometries.  Then the finite-$\vee $-closure \resp {$\vee $-closure} of $S$
is also a {\sam } sub-semigroup of $A$, hence an inverse semigroup by \cite {SemigroupIsInverse}.

\Proof Letting $$ \SS = \big \{\mx T: T \hbox { is a (finite) compatible subset of } S\big \}, $$ we will first prove
that $\SS $ is (finitely-)$\vee $-closed.  For this let us be given a (finite) compatible subset of $\SS $, say $$ \TT =
\{\mx T_j: j\in J\}, $$ where each $T_j$ is a (finite) compatible subset of $S$, and $J$ is a (finite) set of indices.
We then claim that $$ T:= \cup _jT_j $$ is a compatible subset of $S$.  To see this, we must prove that any two elements
$t$ and $s$ in $T$ are compatible, and we will do this by verifying \cite {JoinPisos/b}.  Initially notice that, thanks
to $S$ being a {\sam } semigroup, both $s^*t$ and $st^*$ are partial isometries.  So we deduce from \cite {ProdPIso}
that $ss^*$ commutes with $tt^*$, and that $s^*s$ commutes with $t^*t$.

On the other hand, since $\TT $ is a compatible set, we may use \cite {ExistMaxCompat} to obtain a partial isometry $u$
in $A$ with $$ \mx T_j \preceq u \for j\in J.  $$ It is then evident that $u$ dominates every element of every $T_j$,
and consequently $u$ dominates every element of $T$.  In particular, $u$ dominates $s$ and $t$, thus verifying all of
the conditions in \cite {JoinPisos/b}, and hence we conclude that $s$ and $t$ are compatible.  This proves that $T$ is a
compatible set.

It is now easy to see that $$ \mx \, \TT = \medvee _{j\in J} \big (\mx T_j\big ) = \mx T \in \SS .  $$ thus showing that
$\SS $ is (finitely-)$\vee $-closed\fn {It is worth observing that the assumption that $S$ is a {\sam } semigroup, hence
tame, was used above in an essential way in order to prove that $\hbox {\eulersmall S}$ is $\vee $-closed.  In order to
get to the $\vee $-closure of an arbitrary (wild) set of partial isometries it might be necessary to iterate the above
construction more than once.}.

In order to prove that $\SS $ is closed under multiplication, we pick two elements in $\SS $, say $\mx T$ and $\mx U$,
where $T$ and $U$ are (finite) compatible subset of $S$.  Recalling that $S$ is assumed to be closed under
multiplication, we have that $TU\subseteq S$, so clearly $TU$ consists of partial isometries, which is equivalent to
saying that the initial projections of the members of $T$ commute with the final projections of the members of $U$.
From \cite {PrdMx} we then conclude that $$ (\mx T)(\mx U) = \mx (TU) \in \SS , $$ proving that $\SS $ is closed under
multiplication.

In order to prove that $\SS $ is also closed under adjoints, pick a generic element $\mx T\in \SS $, where $T$ is a
(finite) compatible subset of $S$.  Then it is easy to see that $T^*$ is also a compatible set, and $$ (\mx T)^* = \mx
(T^*) \in \SS .  \omitDoubleDollar \endProof

An easy consequence of the above result is in order:

\state Corollary \label ClosureTame The finite-$\vee $-closure \resp {$\vee $-closure} of any tame set of partial
isometries in a C*-algebra \resp {von-Neumann algebra} is tame.

\Proof Given a tame set $T$ of partial isometries in a C*-algebra \resp {von-Neumann algebra} $A$, let $S$ be the
multiplicative sub-semigroup of $A$ generated by $T\cup T^*$.  Then $S$ satisfies the hypothesis of \cite {ClosureISG},
hence the finite-$\vee $-closure \resp {$\vee $-closure} of $S$, which we denote by $\SS $, is a {\sam } semigroup of
partial isometries, therefore necessarily tame.

The finite-$\vee $-closure \resp {$\vee $-closure} of $T$ is then clearly a subset of $\SS $, hence it is also tame.
\endProof

\nrem Partial isometries have been a topic of interest for a long time.  Among other things it appears in the polar
decomposition of bounded operators.  A thorough study of single partial isometries is to be found in \ref
{HalmosMcLaughlin/1963}.

Proposition \cite {ProdPIso} has been reproved in the literature many times.  The first reference for this result we
know of is \ref {Erdelyi/1968}.

The notion of a tame set of partial isometries and its relationship to partial representations, as discussed in
Proposition \cite {PisoSet}, first appeared in \ref {Exel/2000}.  The order relation among partial isometries defined in
\cite {DefineOrderPisos} has been considered in \ref {HalmosMcLaughlin/1963} in the case of operators on Hilbert's
space.  As mentioned above, it is inspired by the usual order relation on inverse semigroups.

\syschapter {Covariant representations of C*-algebraic dynamical systems}{Covariant representations of C*-systems}

In \cite {CovarTopRepInAlg/iii} we have seen that a covariant representation of an algebraic partial dynamical system
gives rise to a representation of the crossed product.  In the paragraph following \cite {CovarTopRepInAlg} we have also
noticed that a converse of this result is easily obtained in the case of unital ideals.  In the present chapter we will
prove a similar converse for C*-algebras, despite the fact that ideals are not always unital.

\fix Let us therefore fix a C*-algebraic partial dynamical system $$ \big (A,\ G,\ \{\D g\}_{g\in G},\ \{\th g\}_{g\in
G}\big ), $$ for the duration of this chapter.

The appropriate definition of covariant representations of a C*-algebraic partial dynamical system is identical to
Definition \cite {DefineCovarRep}, with the understanding that the target algebra $B$ mentioned there will always be a
C*-algebra.  In respect to this, it is noteworthy that every *-homomorphism between C*-algebras is automatically
continuous.

Quite often the target algebra for our covariant representations will be taken to be the C*-algebra $\Lin (H)$
consisting of all bounded linear operator on a Hilbert space $H$.

\state Proposition \label RepFromCovarInCStar Given a covariant representation $(\pi ,\Pr )$ of a C*-alge\-braic partial
dynamical system $ \big (A, G, \{\D g\}_{g\in G}, \{\th g\}_{g\in G}\big ) $ in a unital C*-algebra $B$, there exists a
unique *-homomorphism $$ \pi \times \Pr : A\rt \Th G\to B, $$ such that $ (\pi \times \Pr ) (a\delta _g) = \pi (a)\pr g,
$ for all $g$ in $G$, and all $a$ in $\D g$.

\Proof Let us use $\pi \times ^{alg}\Pr $ to denote the *-homomorphism from $A\art \Th G$ to $B$, provided by \cite
{FromCovarToRep}.  The conclusion then follows immediately by applying \cite {UniversalityOfEnvelop} to $\pi \times
^{alg}\Pr $.  \endProof

The following is the main result of this chapter:

\state Theorem \label AllRepsComeFromCovar Given a C*-algebraic partial action $$ \Th = \big (\{\D g\}_{g\in G},\ \{\th
g\}_{g\in G}\big ) $$ of a group $G$ on a C*-algebra $A$, and a {\nonDegRep }\fn {\def \Linsmall {\hbox {\rsfootnote
L}\,}\def \Ksmall {\hbox {\rsfootnote K}\,}A *-representation $\rho :B\to \Linsmall (H)$, of a C*-algebra $B$ on a
Hilbert space $H$ is said to be \subjex {\nonDegRep }{non-degenerate representation} if $\clspan {\rho (B)H}=H$ (closed
linear span).  The reader should however be warned that the notion of non-degeneracy for *-representations, as defined
here, is different from the notion of non-degeneracy for *-homomorphisms used earlier.  For example, denoting by
$\Ksmall (H)$ the algebra of all compact operators on a Hilbert space $H$, notice that the inclusion of $\Ksmall (H)$
into $\Linsmall (H)$ is {\nonDegRep } as a *-representation, but it is not {\nonDegHomo } as a *-homomorphism.
Fortunately this somewhat imprecise terminology, which incidentally is used throughout the modern literature, will not
cause any confusion.}  *-representation $$ \rho : A\rt \Th G \to \Lin (H), $$ where $H$ is a Hilbert space, there exists
a unique covariant representation $(\pi ,\Pr )$ of\/ $\Th $ in $\Lin (H)$, such that: \izitem \zitem $\pi $ is a
{\nonDegRep } representation of $A$, \zitem $ \pr g \pri g $ is the orthogonal projection onto $\clspan {\pi (\D g)H}$
(brackets meaning closed linear span), \zitem $\rho = \pi \times \Pr $.

\Proof Denote by $\pi $ the representation of $A$ on $H$ given by $$ \pi (a) = \rho (a\delta _1) \for a\in A.  $$ To see
that $\pi $ is {\nonDegRep }, observe that, since $\rho $ is {\nonDegRep } and since $A\art \Th G$ is dense in $A\rt \Th
G$, the restriction of $\rho $ to the former is {\nonDegRep }.  Thus, the set $$ X := \big \{\rho (a\delta _g)\eta :
g\in G,\ a\in \D g,\ \eta \in H\big \} $$ spans a dense subset of $H$.  Given any element $\xi $ in $X$, say $\xi =\rho
(a\delta _g)\eta $, as above, use the Cohen-Hewitt Theorem \ref {32.22/HewittRoss/1970} to write $a=bc$, with $b,c\in \D
g$.  Then $$ \xi =\rho (a\delta _g)\eta \={LotsAFormulas} \rho (b\delta _1 c \delta _g)\eta = \rho (b\delta _1)\rho (c
\delta _g)\eta = \pi (b)\rho (c \delta _g)\eta \in \pi (A)H.  $$ This shows that $X\subseteq \pi (A)H$, so $\pi (A)H$
also spans a dense subset of $H$, meaning that $\pi $ is {\nonDegRep }, as desired.  For each $g$ in $G$, let $$ H_g =
\clspan {\pi (\D g)H}.  $$ We should remark that, again by the Cohen-Hewit Theorem, every element $\xi $ in $H_g$ may be
written as $$ \xi =\pi (a)\eta , $$ for some $a$ in $\D g$, and $\eta $ in $H$, which means that $$ H_g = \pi (\D g)H,
\equationmark CHforHg $$ even without taking closed linear span.

Letting $\{v_i\}_{i\in I}$ be an approximate identity\fn {An \subj {approximate identity} for a C*-algebra $B$ is a net
$\{v_i\}_{i\in I}\subseteq B$ of positive elements, with $\Vert v_i\Vert \leq 1$, such that $b = \lim _{i\to \infty }
bv_i = \lim _{i\to \infty } v_ib$, for every $b$ in $B$.  Every C*-algebra (and hence also every ideal in a C*-algebra)
is known to admit an approximate identity.}  for $\D g$, we claim that the orthogonal projection onto $H_g$, which we
henceforth denote by $\e g$, is given by $$ \e g(\xi ) = \lim _{i\to \infty } \pi (v_i)\xi \for \xi \in H.
\equationmark ProjFromAproxUnit $$ In fact, given $\xi $ in $H$, write $\xi =\xi _1+\xi _2$, with $\xi _1$ in $ H_g$,
and $\xi _2$ in $H_g^\perp $.  For every $a$ in $\D g$, and every $\eta $ in $H$, observe that $$ \langle \pi (a)\xi
_2,\eta \rangle = \langle \xi _2,\underbrace {\pi (a^*)\eta }_{\in H_g}\rangle = 0, $$ from where we see that $\pi
(a)\xi _2=0$.  In particular $\pi (v_i)\xi _2=0$, for all $i\in I$.  On the other hand, let us use \cite {CHforHg} in
order to write $\xi _1=\pi (a)\eta $, where $a\in \D g$, and $\eta \in H$.  Then $$ \lim _{i\to \infty } \pi (v_i)\xi =
\lim _{i\to \infty } \pi (v_i)\big (\pi (a)\eta + \xi _2\big ) = \lim _{i\to \infty } \pi (v_ia)\eta \$= \pi \big (\lim
_{i\to \infty } v_ia\big )\eta = \pi (a)\eta = \xi _1 = \e g(\xi ), $$ proving \cite {ProjFromAproxUnit}.

It is clear that each $H_g$ is invariant under $\pi $, so one may easily prove that $\e g$ commutes with $\pi (a)$, for
every $a$ in $A$.  We further claim that, given another group element $h\in G$, one has $$ \e g \e h = \e h\e g.
\equationmark TheEgCommute $$ Indeed, given $\xi $ in $H$, observe that $$ \e h\big (\e g(\xi )\big ) = \e h\big (\lim
_{i\to \infty } \pi (v_i)\xi \big ) = \lim _{i\to \infty } \e h\big (\pi (v_i)\xi \big ) \$= \lim _{i\to \infty } \pi
(v_i)\big (\e h(\xi )\big ) = \e g\big (\e h(\xi )\big ).  $$

Another fact we will need later is $$ H_g\cap H_h = \clspan {\pi (\D g\cap \D h)H} \for g,h\in G.  \equationmark
IntersectionsOfH $$ In order to prove it, we note that the inclusion ``$\supseteq $'' is evident, while the reverse
inclusion may be proved as follows: pick any $\xi $ in $H_g\cap H_h$, and emphasizing the role of $h$, write $\xi =\pi
(a)\eta $, with $a$ in $\D h$, and $\eta $ in $H$, by \cite {CHforHg}.  Once more employing our approximate identity
$\{v_i\}_{i\in I}$ for $\D g$, we then have $$ \xi = \e g(\xi ) = \e g\big (\pi (a)\eta \big ) = \lim _{i\to \infty }
\pi (v_i)\pi (a)\eta \$= \lim _{i\to \infty } \pi (v_ia)\eta \in \clspan {\pi (\D g\cap \D h)H}, $$ proving that $$
H_g\cap H_h \ \subseteq \ \clspan {\pi (\D g\cap \D h)H}.  $$ As already observed, Cohen-Hewit implies that $\pi (\D
g\cap \D h)H$ (without taking closed linear span) is a closed linear subspace of $H$, so \cite {IntersectionsOfH} is
verified.

In order to kick-start the construction of the partial representation $\Pr $ referred to in the statement, let us now
construct, for any given $g$ in $G$, an isometric linear operator $ v_g:H_{g\inv } \to H, $ satisfying $$ v_g\big (\pi
(a)\xi \big ) = \rho \big (\th g(a)\delta _g\big )\xi \for a\in \Di g \for \xi \in H.  $$ With this goal in mind we
claim that, for all $a_1,\ldots ,a_n\in \Di g$, and all $\xi _1,\ldots ,\xi _n\in H$, one has that $$ \Big \Vert \usoma
{i=1}n \pi (a_i)\xi _i\Big \Vert = \Big \Vert \usoma {i=1}n \rho \big (\th g(a_i)\delta _g\big )\xi _i\Big \Vert .
\equationmark EqualityOfNorms $$ Indeed, starting from the right-hand-side above, we have $$ \Big \Vert \usoma {i=1}n
\rho \big (\th g(a_i)\delta _g\big )\xi _i \Big \Vert ^2 = \Big \langle \usoma {i,j=1}n \rho \Big (\big (\th
g(a_j)\delta _g\big )^* \big (\th g(a_i)\delta _g\big )\Big ) \xi _i,\xi _j \Big \rangle \={LotsAFormulas} $$$$= \Big
\langle \usoma {i,j=1}n \pi (a_j^* a_i) \xi _i,\xi _j \Big \rangle = \Big \langle \usoma {i=1}n \pi (a_i) \xi _i, \usoma
{j=1}n \pi (a_j) \xi _j \Big \rangle = \Big \Vert \usoma {i=1}n \pi (a_i) \xi _i\Big \Vert ^2, $$ proving \cite
{EqualityOfNorms}.  This said, the correspondence $$ \usoma {i=1}n \pi (a_i)\xi _i \mapsto \usoma {i=1}n \rho \big (\th
g(a_i)\delta _g\big )\xi _i, $$ may now be shown to provide a well defined isometric linear map on the linear span of
$\pi (\Di g)H$, which we have seen coincides with $H_{g\inv }$ by \cite {CHforHg}, providing an isometric linear
operator $v_g$ on $H_{g\inv }$ satisfying the desired properties.

Next we show that the range of $v_g$ is precisely $H_g$.  For this, given $a\in \Di g$, write $a=bc$, with $b,c\in \Di
g$, by the Cohen-Hewitt Theorem.  So, for all $\xi \in H$, we have that $$ v_g\big (\pi (a)\xi \big ) = \rho \big (\th
g(bc)\delta _g\big )\xi = \rho \big (\th g(b)\delta _1\big ) \ \rho \big (\th g(c)\delta _g\big )\ \xi \in \pi \big (\th
g(b)\big ) H \subseteq H_g, $$ proving that $v_g(H_{g\inv })\subseteq H_g$.  On the other hand, given any $a$ in $\D g$
and $\xi $ in $H$, use Cohen-Hewit again to write $\thi g(a) = bc$, with $b,c\in \Di g$.  Then $$ a\delta _1= \th
g(bc)\delta _1 \={LotsAFormulas} \big (\th g(b)\delta _g\big )\ (c\delta _{g\inv }), $$ so $$ \pi (a)\xi = \rho (a\delta
_1)\xi = \rho \big (\th g(b)\delta _g\big ) \underbrace {\rho (c\delta _{g\inv })\xi }_{\eta } = v_g\big (\pi (b)\eta
\big ) \in v_g(H_{g\inv }), $$ where $\eta $ is as indicated.  It follows that $H_g \subseteq v_g(H_{g\inv })$, thus
showing that $v_g$ is indeed an isometric operator from $H_{g\inv }$ onto $H_g$.

Taking one step closer to our desired partial representation, for each $g$ in $G$, we let $$ \pr g: H \to H $$ be the
linear operator defined on the whole of $H$ by extending $v_g$ to be zero on the orthogonal complement of $H_{g\inv }$.
It is then clear that $\pr g$ is a partial isometry in $\Lin (H)$, with initial space $H_{g\inv }$ and final space
$H_g$.  Consequently the orthogonal projection onto $H_g$ satisfies $$ \e g = \pr g \pr g^*.  \equationmark
EgIsProjection $$

If $a$ is in $\Di g$ we have by definition that $ \pr g \pi (a)\xi = \rho \big (\th g(a)\delta _g\big )\xi , $ for all
$\xi $ in $H$, which means that $$ \pr g \pi (a) = \rho \big (\th g(a)\delta _g\big ) \for a\in \Di g.  \equationmark
UgPiaSemXi $$

Having shown that $\pi $ is {\nonDegRep }, it is evident that $H_1=H$, and it is easy to see that $\pr 1$ is the
identity operator on $H$, as required by \cite {PrOnOne}.

We next claim that $$ \pr g \big (\pri g (\xi )\big )= \xi \for g\in G\for \xi \in H_g.  \equationmark UgStarIsInverse
$$ Assuming, as we may, that $\xi =\pi (a)\eta $, with $a\in \D g$, and $\eta \in H$, write $\thi g(a)=bc$, with $b,c\in
\Di g$, by Cohen-Hewit, and observe that $$ \pri g (\xi ) = \pri g \big (\pi (a)\eta \big ) = \rho \big (\thi g(a)\delta
_{g\inv }\big )\eta = \rho (b\delta _1\ c\delta _{g\inv })\eta \$= \pi (b) \underbrace {\rho (c\delta _{g\inv })\eta
}_{\zeta } = \pi (b) \zeta .  $$ Therefore the left-hand-side of \cite {UgStarIsInverse} equals $$ \pr g \big (\pri g
(\xi )\big )= \pr g \big (\pi (b) \zeta \big ) = \rho \big (\th g(b)\delta _g\big )\zeta = \rho \big (\th g(b)\delta
_g\big )\ \rho (c\delta _{g\inv })\eta \$= \rho \big (\th g(bc)\delta _1\big )\eta = \pi (a)\eta = \xi , $$ proving
\cite {UgStarIsInverse}, from where we deduce that $$ \pri g = \pr g^* \for g\in G, \equationmark RespectsAdjoint $$
hence verifying \cite {StarPrepCond} and, in view of \cite {EgIsProjection}, this also proves condition (ii) in the
statement.

We next focus on proving \cite {ObserveRight}, namely $$ \pr {gh}\pri h = \pr g \pr h \pri h, $$ for all $g,h\in G$.  As
a first step we claim that $$ \pr {gh}\pri h (\eta ) = \pr g (\eta ) \for \eta \in H_ h\cap H_{g\inv }.  \equationmark
SomeClaimizinha $$ In order to do this we use \cite {IntersectionsOfH} to write $ \eta =\pi (a)\xi , $ with $a$ in $\D
h\cap \Di g$.  Observing that $$ \thi h(a) \in \Di h\cap \D {h\inv g\inv }, $$ we may write $\th {h\inv }(a) = bc$, with
$b$ and $c$ in $\Di h\cap \D {h\inv g\inv }$.  Therefore $$ \pri h(\eta ) = \pri h\big (\pi (a)\xi \big ) = \rho \big
(\th {h\inv }(a)\delta _{h\inv }\big )\xi = \rho (bc\delta _{h\inv })\xi \$= \rho (b\delta _1)\,\rho (c\delta _{h\inv
})\xi = \pi (b)\,\rho (c\delta _{h\inv })\xi .  $$ So the left-hand-side of \cite {SomeClaimizinha} equals $$ \pr
{gh}\big (\pri h (\eta )\big ) = \pr {gh}\big (\pi (b)\,\rho (c\delta _{h\inv })\xi \big ) = \rho \big (\th
{gh}(b)\delta _{gh}\big ) \,\rho (c\delta _{h\inv })\xi \$= \rho \big (\th {gh}(bc)\delta _g\big )\xi = \rho \big (\th
g(a)\delta _g\big )\xi = \pr g \big (\pi (a)\xi \big ) = \pr g \big (\eta \big ), $$ taking care of \cite
{SomeClaimizinha}.

Observing that the right-hand-side of \cite {ObserveRight} equals $\pr g \e h$, our next claim is that $$ \pr g \e h =
\pr {gh}\pri h \ei g \for g,h\in G.  \equationmark SomeOtherClaimizinha $$ In order to prove this we have $$ \pr g \e h
= \pr g \ei g\e h \={SomeClaimizinha} \pr {gh}\pri h \ei g\e h \={TheEgCommute} $$$$ = \pr {gh}\pri h \e h\ei g = \pr
{gh}\pri h \ei g, $$ proving \cite {SomeOtherClaimizinha}.  Taking adjoints in \cite {SomeOtherClaimizinha}, using \cite
{RespectsAdjoint}, and changing variables appropriately leads to $$ \e h\pr g = \e g \pr h \pr {h\inv g}.  \equationmark
SomeHelp $$

Focusing now on the left-hand-side of \cite {ObserveRight}, observe that $$ \pr {gh} \pri h = \pr {gh} \ei {(gh)} \pri h
\={SomeHelp} \pr {gh} \ei h \pri {(gh)} \pr g.  $$ Since $\pr g = \pr g \ei g$, the right-hand-side above is unaffected
by right multiplication by $\ei g$, and so is the left-hand-side.  Thus $$ \pr {gh} \pri h = \pr {gh} \pri h \ei g
\={SomeOtherClaimizinha} \pr g \e h = \pr g \pr h \pri h, $$ proving \cite {ObserveRight}.  It has already been observed
that \cite {ObserveLeft} is implied by \cite {ObserveRight} and \cite {StarPrepCond}, so the verification that $\Pr $ is
a partial representation is complete.  We now verify that $(\pi ,\Pr )$ is a covariant representation, which is to say
that $$ \pr g \pi (a) \pri g = \pi \big (\th g(a)\big ), \equationmark TheCovariance $$ for any given $g$ in $G$, and
$a$ in $\Di g$.  To prove it, let $\xi \in H$, and let us first suppose that $\xi $ is in $H_g^\perp $.  Then, evidently
$ \pri g(\xi ) =0, $ so the operator on the left-hand-side of \cite {TheCovariance} vanishes on $\xi $.  We claim that
the same is true with respect to the operator on the right-hand-side.  Indeed, given any $\eta $ in $H$, we have that $$
\big \langle \pi \big (\th g(a)\big )\xi ,\eta \big \rangle = \big \langle \xi ,\underbrace {\pi \big (\th g(a^*)\big
)\eta }_{\in H_g}\big \rangle = 0.  $$ Since $\eta $ is arbitrary, we have that $\pi \big (\th g(a)\big )\xi =0$, as
claimed.

Suppose now that $\xi $ is in $H_g$.  Then we may write $\xi =\pi (b)\eta $, for $b\in \D g$, and $\eta \in H$, by \cite
{CHforHg}, and $$ \pr g \pi (a) \pri g(\xi ) = \pr g \pi (a) \pri g\big (\pi (b)\eta \big ) \={UgPiaSemXi} $$$$ = \rho
\big (\th g(a)\delta _g\big ) \rho \big (\thi g(b)\delta _{g\inv }\big )\eta = \rho \Big (\th g(a)\delta _g \ \thi
g(b)\delta _{g\inv }\Big )\eta \={LotsAFormulas} $$$$ = \rho \big (\th g(a)b\delta _1\big )\eta = \pi \big (\th g(a)\big
) \pi (b)\eta = \pi \big (\th g(a)\big ) \xi .  $$ This concludes the proof of \cite {TheCovariance}, and hence that
$(\pi ,\Pr )$ is a covariant representation as needed.  In order to show point (iii) in the statement, it clearly
suffices to prove that $$ (\pi \times u)(a\delta _g) = \rho (a\delta _g), \equationmark RecoverRho $$ for every $g$ in
$G$, and $a\in \D g$, which we will now do.  Writing $a=bc$, with $b$ and $c$ in $\D g$, we have $$ \rho (a\delta _g)
\={LotsAFormulas} \rho \big (b\delta _g\, \thi g(c)\delta _1\big ) = \rho (b\delta _g)\pi \big (\thi g(c)\big )
\={TheCovariance} $$$$ = \rho (b\delta _g)\pri g \pi \big (c)\pr g \={UgPiaSemXi} \rho (b\delta _g)\rho \big (\thi
g(c)\delta _{g\inv }\big )\pr g \={LotsAFormulas} $$$$ = \rho (bc\delta _1)\pr g = \pi (a) \pr g = (\pi \times \Pr
)(a\delta _g), $$ so \cite {RecoverRho} is verified.

We are therefore left with the only remaining task of proving uniqueness of the $(\pi ,\Pr )$.  With this goal in mind,
assume that $(\pi ',\Pr ')$ is another covariant representation of our system in $\Lin (H)$, such that $\pi '\times \Pr
' = \rho $.  For every $a$ in $A$ we then have that $\pi '(a) = \rho (a\delta _1) = \pi (a)$, so $\pi $ and $\pi '$ must
coincide.

Given $g$ in $ G$, and $\xi $ in $H_{g\inv }$, write $\xi =\pi (a)\eta $, with $a$ in $\Di g$, and $\eta $ in $H$, by
\cite {CHforHg}.  Then $$ \Pr '_g(\xi ) = \Pr '_g\pi (a)\eta = \Pr '_g\pi '(a)\eta = \big (\pi '(a^*)\Pr '_{g\inv }\big
)^*\eta = \big ((\pi '\times \Pr ')(a^*\delta _{g\inv })\big )^*\eta \$= \big ((\pi \times \Pr )(a^*\delta _{g\inv
})\big )^*\eta = \pr g\pi (a)\eta = \pr g(\xi ), $$ which says that $\Pr '_g$ coincides with $\pr g$ on $H_{g\inv }$.
By (ii) we have that both $\Pr '_g$ and $\pr g$ vanish on $H_{g\inv }^\perp $, so $\Pr '_g$ coincides with $\pr g$ on
the whole of $H$.  This proves that $\Pr '=\Pr $, and hence completes the proof of the uniqueness part.  \endProof

\nrem Theorem \cite {AllRepsComeFromCovar} is due do McClanahan \ref {Proposition 2.8/McClanahan/1995}.  The above proof
is inspired by \ref {Theorem 1.3/ExelLacaQuigg/1997}\fn {Please note that \ref {ExelLacaQuigg/1997} is the preprint
version of \ref {ExelLacaQuigg/2002}.}, except that we have avoided the use of approximate identities whenever possible,
basing the arguments on the Cohen-Hewit Theorem instead.  Even though the proof presented here turned out to be a bit
long, we have chosen it because we believe a possible generalization of this result to algebraic partial actions will
most likely use arguments based on idempotency rather than convergence of limits.

\chapter Partial representations subject to relations

\label PrepCondChap \def \G {{\cal G}} \def \ur {^\Rel } \def \BZero {A} \def \Cstarpar {C^*\underpar (G)} \def \segt
#1#2{\overline {#1#2}} \def \conv {\hbox {\sixrm conv}}

In the present chapter we will describe one of the most efficient ways to show a given C*-algebra to have a partial
crossed product description.

In order to motivate our method, suppose that $B$ is a unital C*-algebra defined in terms of generators and relations.
Suppose also that the given relations imply that the generators are partial isometries or even, as it often happens,
explicitly state this fact.  Let us moreover suppose that these partial isometries form a tame set so, by \cite
{PisoSet}, we may find a *-partial representation $\Pr $ of a group $G$ (quite often a free group) whose range also
generates $B$.  Under all of these favorable circumstances $B$ may then be given another presentation in which the set
of generators turns out to be the range of a *-partial representation, and we may then consider the partial dynamical
system described in \cite {SystemFromPrep}.

Having the same goal as Theorem \cite {GradedAlgsAsCrossProd}, although with different hypothesis, and in a different
category, the main result of this chapter is intended to specify conditions under which the homomorphism given in \cite
{MapToBzero} is in fact an isomorphism.  Another major objective of this chapter is to give a very concrete picture of
the spectrum of the commutative algebra involved in this system.

Given a group $G$, we will be dealing with universal\fn {See \ref {Blackadar/1985} for a definition of universal
C*-algebras given by generators and relations.}  unital C*-algebras on the set of generators $$ \G = \{\pr g: g\in G\},
$$ subject to a given set of relations.  In all instances below this set of relations will split as $$ \Rel '\cup \Rel ,
$$ where $\Rel '$ consists precisely of relations \cite {DefPR/i--iv}, which is to say that the correspondence $g\mapsto
\pr g$ is a *-partial representation of $G$ in $B$.  With respect to the remaining set of relations, namely $\Rel $, we
will always make the assumption that it consists of algebraic relations involving only the $\e g$ (defined to be $\e g =
\pr g\pri g$).  By this we mean relations of the form $$ p(\e {g_1},\e {g_2}, \ldots ,\e {g_n}) =0, \equationmark
PolyRel $$ where $p$ is a complex polynomial in $n$ variables, and the $g_i$ are in $G$.  Recall that the $\e g$ commute
by \cite {MainTheEgCommute}, so it is OK to apply a polynomial in $n$ commuting variables to them.

\definition \label DefineCstarParRel Given a set $\Rel $ of relations of the form \cite {PolyRel}, we will denote by
$\CstarparRel G\Rel $ the universal unital C*-algebra generated by the set $\G = \{\pr g: g\in G\}$, subject to the set
of relations $\Rel '\cup \Rel $, where $\Rel '$ consists of the relations \cite {DefPR/i--iv}.

Notice that, since the set $\Rel '$ is always involved by default, we have decided not to emphasize it in the notation
introduced above.  The subscript ``$\underpar $'' should be enough to remind us that the axioms for *-partial
representations, namely $\Rel '$, is also being taken into account.

\definition Let $v$ be a *-partial representation of a group $G$ in a unital C*-algebra $B$, and let $\Rel $ be a set of
relations of the form \cite {PolyRel}.  We will say that $v$ \"{satisfies} $\Rel $ if, for every relation ``$ p(\e
{g_1},\e {g_2}, \ldots ,\e {g_n}) = 0 $'' in $\Rel $, one has that $$ p(v_{g_1}v_{g_1\inv },v_{g_2}v_{g_2\inv }, \ldots
,v_{g_n}v_{g_n\inv }) =0.  $$

The universal property of $\CstarparRel G\Rel $ may then be expressed as follows.

\state Proposition \label UniversalForPreps For every *-partial representation $v$ of $G$ in a unital C*-algebra $B$,
satisfying $\Rel $, there exists a unique *-homomorphism $$ \varphi :\CstarparRel G\Rel \to B, $$ such that $\varphi
(\pr g) = v_g$, for every $g\in G$.

\Proof The proof follows immediately from the universality of $\CstarparRel G\Rel $, once we realize that, besides
satisfying $\Rel $, the $v_g$ also satisfy $\Rel '$, as a consequence of it being a *-partial representation.  \endProof

As an example of relations of the form \cite {PolyRel} one is allowed to express that a certain $\pr g$ is an isometry,
since this reads as ``$ (\pr g)^* \pr g =1 $'', and may therefore be expressed as ``$ \ei g - 1 = 0.  $'' As another
example, one may express that two partial isometries $\pr g$ and $\pr h$ have orthogonal ranges by writing down the
relation ``$ \e g\e h =0 $''.  On the other hand we are ruling out relations such as ``$\pr g+\pr h = \pr k$'', since
these are not of the above form.

We should however remark that certain relations which do not immediately fit under \cite {PolyRel} may sometimes be
given an equivalent formulation within that framework, an example of which we would now like to present.

\state Proposition \label CondForSSat Let $\Pr $ be a *-partial representation of the group $G$ in a unital C*-algebra.
Given two elements $g$ and $h$ in $G$, the following are equivalent \izitem \zitem $\pr {gh} = \pr g \pr h$, \zitem $\e
{gh} = \e {gh}\e g$, \medskip \noindent where $\e g = \pr g\pri g$, as usual.

\Proof By the C*-identity \cite {DefineCstarAlg/iii}, we have $$ \Vert \pr {gh} - \pr g\pr h \Vert ^2 = \Vert (\pr {gh}
- \pr g\pr h)(\pr {gh} - \pr g\pr h)^*\Vert \$= \Vert \pr {gh} \pr {h\inv g\inv } - \pr {gh} \pri h \pri g - \pr g\pr h
\pr {h\inv g\inv } + \pr g\pr h \pri h \pri g \Vert \={ObserveLeft} $$$$ = \Vert \e {gh} - \pr {gh} \pri h \pri g\Vert .
$$ In addition, we have $$ \pr {gh}\pr {h\inv }\pr {g\inv } = \pr {gh}\pr {h\inv }\pr {g\inv } \pr g \pr {g\inv } = \pr
{gh}\pr {h\inv g\inv }\pr g \pr {g\inv } = \e {gh}\e g, $$ which plugged above gives $$ \Vert \pr {gh} - \pr g\pr h
\Vert ^2 = \Vert \e {gh} - \e {gh}\e g\Vert .  $$ from where the statement follows immediately.  \endProof

As an important application of this idea, one may phrase the fact that a *-partial representation is semi-saturated (see
Definition \cite {DefineSemiSatPRept}), relative to a given length function $\ell $ on $G$, by requiring it to satisfy
the set of relations $$ \e {gh} - \e {gh}\e g = 0, $$ for all $g$ and $h$ in $G$ such that $\ell (gh)=\ell (g)+\ell
(h)$.

\medskip One of the main goals in this chapter is to prove that $\CstarparRel G\Rel $ is isomorphic to a partial crossed
product of the form $$ C(\OR )\rt {\Th _\Rel } G.  $$

The partial dynamical system involved in this result will be presented as a restriction of the partial Bernoulli action,
introduced in \cite {PartialBernoulli}, to a closed invariant subset (see \cite {RestrToInvar}).  The first step towards
this goal will therefore be to describe the subset $\OR \subseteq \OuG $ mentioned above.

Let us consider, for each $g$ in $G$, the mapping $$ \ee g: \OG \to \{0,1\}, \equationmark DefineEEg $$ defined by $$
\ee g(\omega ) = \omega _g = \bool {g\in \omega } \for \omega \in \OG , $$ (see \cite {VectorsAndSets}).  Seeing $\OG $
as a product space, the $\ee g$ are precisely the standard projections, hence continuous functions.

Consider a relation $$ p(\e {g_1},\e {g_2}, \ldots ,\e {g_n}) =0, $$ of the form \cite {PolyRel}.  If we replace each
$\e g$ by $\ee g$ in the left-hand-side above we get $$ p(\ee {g_1},\ee {g_2}, \ldots ,\ee {g_n}), $$ which may be
interpreted as a complex valued function on $\OG $, namely $$ \omega \in \OG \ \mapsto \ p\big (\ee {g_1}(\omega ), \ee
{g_2}(\omega ),\ldots , \ee {g_n}(\omega )\big ) \in \C .  \equationmark DefineFforRel $$

\definition \label DefineOmegaR Given a set $\Rel $ of relations of the form \cite {PolyRel}, we will let $\FRel $ be
the set of functions on $\OG $ of the form \cite {DefineFforRel}, obtained via the above substitution procedure from
each relation in $\Rel $.  The \subjex {spectrum}{spectrum of a set of relations} of $\Rel $, denoted by $\OR $, is then
defined to be the subset of $\OuG $ defined by $$ \OR = \big \{ \omega \in \OuG : f(g\inv \omega ) = 0,\ \forall f\in
\FRel ,\ \forall g\in \omega \big \}.  $$

\state Proposition \label ORIsCompactInvar For any set $\Rel $ of relations of the form \cite {PolyRel}, one has that
$\OR $ is a compact subset of\/ $\OuG $, which is moreover invariant under the partial Bernoulli action $\beta $ defined
in \cite {PartialBernoulli}.

\Proof In order to prove that $\OR $ is closed in $\OuG $, let $\{\omega _i\}_{i\in I}$ be a net in $\OR $ converging to
some $\omega $ in $\OuG $.  Given $g\in \omega $, observe that $\ee g(\omega )=1$, so $\omega $ lies in the open set
$\ee g\inv (\{1\})$.  We may therefore assume, without loss of generality, that every $\omega _i$ lie in $\ee g\inv
(\{1\})$.  This is to say that $g\in \omega _i$, for every $i\in I$, and since $\omega _i$ is in $\OR $, we deduce that
$$ f(g\inv \omega _i) =0 \for i\in I \for f\in \FRel .  $$

Observing that the correspondence $$ \nu \in \OG \mapsto f(g\inv \nu )\in \C $$ is a continuous mapping, we conclude
that $$ f(g\inv \omega ) = \lim _{i\to \infty } f(g\inv \omega _i) = 0, $$ proving that $\omega \in \OR $, and hence
that $\OR $ is closed.  Since $\OuG $ is compact, then so is $\OR $.

Recall that the partial Bernoulli action $\beta = \big (\{\D g\}_{g\in G}, \{\beta _g\}_{g\in G}\big )$ is given by $$
\D g = \{\omega \in \OuG \ : \ g\in \omega \} \and \beta _g(\omega ) = g\omega \for \omega \in \Di g.  $$

In order to prove invariance, we must show that $\beta _g(\OR \cap \Di g)\subseteq \OR $, so pick $\omega $ in $\OR \cap
\Di g$, and let $f\in \FRel $, and $h\in \beta _g(\omega )$ be given.  Then $ g\inv h\in \omega , $ and since $\omega
\in \OR $, we have that $$ 0 = f\big ((g\inv h)\inv \omega \big ) = f(h\inv g\omega ) = f\big (h\inv \beta _g(\omega
)\big ), $$ proving that $\beta _g(\omega )$ is in $\OR $.  \endProof

\fix From now on we will fix an arbitrary set $\Rel $ of relations of the form \cite {PolyRel}.  We may then use \cite
{RestrToInvar} to restrict the partial Bernoulli action to $\OR $, obtaining a partial action $$ \Th _\Rel = \big (\{\D
g\ur \}_{g\in G}, \{\th g\ur \}_{g\in G}\big ), \equationmark RelPAction $$ where $$ \D g\ur = \OR \cap \D g = \{\omega
\in \OR : g\in \omega \}, \equationmark DescriptBernouliRelDoms $$ and $$ \th g\ur (\omega ) = g\omega \for \omega \in
\D {g\inv }.  $$

It should be noticed that, since each $\D g$ is open in $\OuG $, one has that $\D g\ur $ is open in $\OR $.  This, plus
the obvious fact that each $\th g\ur $ is continuous, means that $\Th _\Rel $ is a topological partial action of $G$ on
$\OR $.

Observe that $\D g\ur $ is also closed in $\OR $ by \cite {DescriptBernouliRelDoms}, so \cite {ClopenGlobalizes} applies
and hence we see that $\Th _\Rel $ admits a Hausdorff globalization.  In fact we may describe the globalization of $\Th
_\Rel $ without resorting to \cite {ClopenGlobalizes} as follows: define $$ X_\Rel = \big \{ \omega \in \OG : f(g\inv
\omega ) = 0,\ \forall f\in \FRel ,\ \forall g\in \omega \big \}, $$ that is, $X_\Rel $ is defined as in \cite
{DefineOmegaR}, except that we have replaced ``$\omega \in \OuG $'' by ``$\omega \in \OG $''.  In other words, $$ \OR =
X_\Rel \cap \OuG .  $$

By staring at the above definition of $X_\Rel $, it is clear that $X_\Rel $ is invariant under the global Bernoulli
action $\eta $ described in \cite {BernoulliAction}.

Notice that $\omega =\varnothing $, namely the empty set, is a member of $X_\Rel $ since the membership condition above
is vacuously verified for $\varnothing $.  If $\omega $ is any element of $X_\Rel $ other than $\varnothing $, choose
$g$ in $\omega $, and notice that $1 \in g\inv \omega = \eta _{g\inv }(\omega )$, so $$ \eta _{g\inv }(\omega )\in
X_\Rel \cap \OuG = \OR , $$ whence $$ \omega = \eta _g\big (\eta _{g\inv }(\omega )\big ) \in \eta _g(\OR ).  $$

With this it is easy to see that the orbit of $\OR $ under $\eta $ is given by $$ \bigcup _{g\in G} \eta _g(\OR ) =
X_R\setmenos \{\varnothing \}.  $$

This proves the following:

\state Proposition The restriction of the global Bernoulli action $\eta $ to the invariant subset $X_R\setmenos
\{\varnothing \}$ coincides with the globalization of the partial action $\Th _\Rel $ introduced in \cite {RelPAction}.

Let us now turn our topological partial action into a C*-algebraic one.  For this we invoke \cite {TopoCstarActions} to
obtain a partial action of $G$ on $C(\OR )$, which, by abuse of language, we will also denote by $\Th _\Rel $.

As seen in our discussion after the definition of the partial Bernoulli action \cite {PartialBernoulli}, $\D g$ is a
compact space for every $g$ in $G$, and hence so is $\D g\ur $.  The corresponding ideal of $C(\OR )$, namely $$ C_0(\D
g\ur ) = C(\D g\ur ), $$ is therefore a unital ideal.  In particular, we are entitled to use \cite {PRinCrossProd} in
order to get a *-partial representation $$ v: G \to C(\OR )\rt {\Th _\Rel } G, \equationmark ThePartialRepOnCrosPod $$
defined by $$ v_g = 1_g\delta _g \for g\in G, $$ where $1_g$ denotes the unit\fn {We could also use the heavier notation
$1_g^\Rel $ for the unit of $C(\D g\ur )$, but this will be unnecessary.}  of $C(\D g\ur )$.  Viewed within $C(\OR )$,
observe that $1_g$ is the characteristic function of $\D g\ur $.  By the description of $\D g\ur $ given in \cite
{DescriptBernouliRelDoms}, we may then write $$ 1_g(\omega ) = \bool {g\in \omega } \for \omega \in \OR , $$ which is to
say that $$ 1_g = \ee g|_{\OR }.  \equationmark OnegIsEg $$

\state Lemma \label CanonicRepSatisfies The *-partial representation $v$ defined in \cite {ThePartialRepOnCrosPod}
satisfies $\Rel $.

\Proof Given a relation ``$ p(\e {g_1},\e {g_2}, \ldots ,\e {g_n}) = 0 $'' in $\Rel $, we must therefore prove that $$
p(v_{g_1}v_{g_1\inv },v_{g_2}v_{g_2\inv }, \ldots ,v_{g_n}v_{g_n\inv }) =0.  \eqno (\dagger ) $$ For any $g$ in $G$,
notice that $$ v_gv_{g\inv } = (1_g\delta _g) (1_{g\inv }\delta _{g\inv }) \={LotsAFormulas} 1_g\delta _1, $$ so, under
the usual identification of $C(\OR )$ as a subalgebra of $C(\OR )\rt {\Th _\Rel } G$ provided by \cite
{EmbeddIntoCrossProd}, we may write $$ v_gv_{g\inv } = 1_g = \ee g|_{\OR }.  $$ Therefore we see that the left-hand-side
of $(\dagger )$ equals the restriction of $$ p(\ee {g_1}, \ee {g_2}, \ldots ,\ee {g_n}), $$ to $\OR $.  This is
precisely one of the functions in $\FRel $ which, by the very definition of $\OR $ (with $g=1$) vanishes identically on
$\OR $.  This verifies $(\dagger )$ and hence completes the proof.  \endProof

With this we arrive at the main result of this chapter:

\state Theorem \label CstarRelAsCP Let $G$ be a group and let $\Rel $ be a set of relations of the form \cite {PolyRel}.
Then there exists a *-isomorphism $$ \varphi :\CstarparRel G\Rel \to C(\OR )\rt {\Th _\Rel } G, $$ such that $$ \varphi
(\pr g) = 1_g\delta _g \for g\in G, $$ where $1_g$ denotes the characteristic function of $\D g\ur $.

\Proof By \cite {CanonicRepSatisfies} and \cite {UniversalForPreps} there exists a *-homomorphism $\varphi $ satisfying
the conditions in the statement, and all we must do is prove $\varphi $ to be an isomorphism.

We begin by building a covariant representation of $\Th _\Rel $ in $\CstarparRel G\Rel $.  For this, let $\BZero $ be
the closed *-subalgebra of $\CstarparRel G\Rel $ generated by the set $\{\e g: g\in G\}$.  Noticing that $A$ is
commutative by \cite {TheEgCommute}, we denote by $\hat \BZero $ its spectrum.

Consider the mapping $$ h: \gamma \in \hat \BZero \mapsto \big (\gamma (\e g)\big )_{g\in G} \in \OG .  $$ Since each
$\e g$ is idempotent, we have that $\gamma (\e g)\in \{0,1\}$, so the above is well defined and the reader will not have
any difficulty in proving that $h$ is continuous.  We next claim that the range of $h$ is contained in $\OR $.  Picking
$\gamma \in \hat \BZero $, we must then prove that $ \omega :=h(\gamma )\in \OR .  $ Evidently $$ \bool {1\in \omega } =
\omega _1 =\gamma (\e 1) = \gamma (1)=1, $$ so $\omega \in \OuG $.  Let $f$ be a function on $\OG $ of the form \cite
{DefineFforRel}, associated to a relation in $\Rel $.  Picking any $g\in \omega $, we must therefore prove that $$
f(g\inv \omega ) = 0.  $$ With $p,g_1,\ldots ,g_n$, as in \cite {DefineFforRel}, observe that $$ f(g\inv \omega ) =
p\big (\ee {g_1}(g\inv \omega ),\ldots , \ee {g_n}(g\inv \omega )\big ).  \equationmark FginvOmega $$ In order to get a
better grasp on the meaning of the above, notice that for any $i=1,\ldots ,n$, one has $$ \ee {g_i}(g\inv \omega ) =
\bool {g_i\in g\inv \omega } = \bool {gg_i\in \omega } = \omega _{gg_i} = \gamma (\e {gg_i}) = \cdots $$ Since $g\in
\omega $, we see that $\omega _g=1$, which is to say that $\gamma (\e g)=1$, so the above equals $$ \cdots = \gamma (\e
g) \gamma (\e {gg_i}) = \gamma (\pr g \pri g\e {gg_i} ) \={PrepComutRel} \gamma (\pr g \e {g_i} \pri g\big ) = \gamma
_g(\e {g_i}\big ), $$ where $\gamma _g: \BZero \to {\bf C}$ is defined by $$ \gamma _g(b) =\gamma (\pr g b \pri g\big )
\for b\in \BZero .  $$

It is easy to see that $\gamma _g$ is a homomorphism, which therefore commutes with polynomials.  So by \cite
{FginvOmega} we have $$ f(g\inv \omega ) = p\big (\gamma _g(\e {g_1}), \ldots , \gamma _g(\e {g_n})\big ) = \gamma
_g\big (p\big (\e {g_1}, \ldots , \e {g_n})\big ) = 0.  $$ This proves that $h$ indeed maps $\hat \BZero $ into $\OR $,
so by dualization we obtain a unital *-homomorphism $$ \hat h : f\in C(\OR ) \mapsto f\circ h\in C(\hat \BZero ).  $$

Observe that $\hat h(1_g)$ coincides with the Gelfand transform $\hat \e g$ of $\e g$, because, for every $\gamma \in
\hat \BZero $, one has $$ \hat h(1_ g)(\gamma ) = 1_g\big (h(\gamma )\big ) \={OnegIsEg} \ee g\big (h(\gamma )\big ) =
\gamma (\e g) = \hat \e g(\gamma ).  $$

Identifying $C(\hat \BZero )$ with $\BZero $ via the Gelfand transform, we then have that $\hat h$ is a *-homomor\-phism
$$ \hat h : C(\OR ) \to \BZero , $$ such that $ \hat h(1_ g) = \e g, $ for all $g$ in $G$.  We will now prove that the
pair $(\hat h, \Pr )$ is a covariant representation of $\Th _\Rel $ in $\CstarparRel G\Rel $, which is to say that $$
\pr g \hat h(f) \pri g = \hat h\big (\th g\ur (f)\big ), \equationmark ThisIsCovarHere $$ for all $f\in C(\Di g\ur )$,
and all $g\in G$.

Using the \pilar {12pt} Stone-Weiestrass Theorem one may show that the functions of the form $f = 1_ {g\inv }1_ h$, with
$h$ in $G$, generate $C(\Di g\ur )$, as a C*-algebra, so it is enough to prove \cite {ThisIsCovarHere} only for such
functions.  In this case the right-hand-side of \cite {ThisIsCovarHere} becomes $$ \hat h\big (\th g\ur (f)\big ) = \hat
h\big (\th g\ur (1_ {g\inv }1_ h)\big ) \={IdealUnits} \hat h(1_g 1_{gh}) = \e g \e {gh}.  $$ The left-hand-side of
\cite {ThisIsCovarHere}, on the other hand, is given by $$ \pr g \hat h(f) \pri g = \pr g \hat h(1_ {g\inv }1_ h) \pri g
= \pr g \e {g\inv }\e h \pri g \$= \pr g \e h \pri g \={PrepComutRel} \pr g \pri g \e {gh} = \e g \e {gh}, $$ so \cite
{ThisIsCovarHere} is seen to hold, as claimed, whence $(\hat h,\Pr )$ is indeed a covariant representation.  From \cite
{RepFromCovarInCStar} we then obtain a *-homomorphism $$ \hat h\times \Pr : C(\OR )\rt ? G \to \CstarparRel G\Rel , $$
satisfying $$ (\hat h\times \Pr )(1_g\delta _g) = \hat h(1_g)\pr g = \e g \pr g = \pr g \pri g \pr g = \pr g, $$ for all
$g$ in $G$.  Observing that $\varphi $ sends $\pr g$ to $1_g\delta _g$, and that $\CstarparRel G\Rel $ is generated by
the $\pr g$, one easily proves that $(\hat h\times \Pr )\circ \varphi $ is the identity on $\CstarparRel G\Rel $.

To prove that the composition $\varphi \circ (\hat h\times \Pr )$ also coincides with the identity it suffices to prove
that $C(\OR )\rt ?  G$ is generated by the $1_g\delta _g$.

We have already seen that each $\D g\ur $ is generated by the elements $1_g1_h$, for $h$ ranging in $G$, so $ C(\OR )\rt
? G $ is generated by the elements of the form $$ 1_h1_g\delta _g = (1_h\delta _e)(1_g\delta _g) = (1_h\delta
_h)(1_{h\inv }\delta _{h\inv })(1_g\delta _g).  $$ Thus the crossed product is indeed generated by the $1_g\delta _g$,
as claimed, so $\varphi \circ (\hat h\times \Pr )$ is the identity map, whence $\varphi $ and $(\hat h\times \Pr )$ are
each other's inverse, proving that $\varphi $ is an isomorphism and hence concluding the proof.  \endProof

In this and in later chapters we will present several applications of \cite {CstarRelAsCP}, describing certain classes
of C*-algebras as partial crossed products.  In the first such application we will prove a version of \cite
{ParGpAlgAsCP} to the context of C*-algebras, so let us begin by adapting the definition of $\Kpar $ given in \cite
{DefineKPar} to our context.

\definition \label DefineCStarPar Given a group $G$, the \subj {partial group C*-algebra} of $G$, denoted $\Cstarpar $,
is defined to be the algebra $\CstarparRel G\Rel $, where $\Rel $ is the empty set of relations.

Thus, by \cite {UniversalForPreps} we see that $\Cstarpar $ is the universal unital C*-algebra for *-partial
representations of $G$ (without any further requirements), a fact that may be considered a version of \cite {UnivPRep}
to C*-algebras.

Notice that if $\Rel $ is the empty set of relations, then $\OR =\OuG $ so, as an immediate consequence of \cite
{CstarRelAsCP} we have:

\state Corollary \label CstarParAsCP For every group\/ $G$, one has that $\Cstarpar $ is *-isomorphic to the crossed
product of $ C(\OuG ) $ by $G$, relative to the partial Bernoulli action $\beta $ defined in \cite {PartialBernoulli},
under an isomorphism which sends each canonical generating partial isometry $\pr g$ in $\Cstarpar $ to $1_g\delta _g$,
where $1_g$ denotes the characteristic function of $\D g$.

It is not hard to see that $C(\OuG )$ is the universal unital C*-algebra generated by a set $ {\cal E} = \{\e g: g\in
G\}, $ subject to the relations stating that the $\e g$ are commuting self-adjoint idempotents, and that $\e 1=1$.  In
other words, $C(\OuG )$ is the analogue of the algebra $\Apar $ introduced in \cite {DefineApar} within the category of
C*-algebras.  Therefore \cite {CstarParAsCP} may be seen as the analogue of \cite {ParGpAlgAsCP} in the present
category.

\bigskip

As our second example, let us consider the universal C*-algebra for semi-saturated partial representations.  For this
let us suppose that $G$ is a group equipped with a length function $\ell $.  One may then define a \subj {pseudo-metric}
on $G$ by setting $$ d(g,h) = \ell (g\inv h) \for g,h\in G.  $$ By pseudo-metric we simply mean that $d$ is a
nonnegative real valued function satisfying the triangle inequality: $$ d(g,h) \leq d(g,k) + d(k,h) \for g,h,k\in G.  $$

Although we will not need it here, one may clearly relate certain properties of $\ell $ with the remaining axioms for
metric spaces.

In Euclidean space we know that the triangle inequality becomes an equality if and only if the three points involved are
suitably placed in a straight line.  This motivates the following:

\definition \label DefineConvex Let $G$ be a group equipped with a length function $\ell $.  \iaitem \aitem Given $g$
and $h$ in $G$, the \subjex {segment}{segment in a group} joining $g$ and $h$ is the subset of $G$ given by $$ \segt gh
= \big \{x\in G: d(g,h) = d(g,x) + d(x,h)\big \}.  $$ \aitem A subset $\omega \subseteq G$ is said to be \subjex
{convex}{convex subset of a group} if $\segt gh \subseteq \omega $, for all $g$ and $h$ in $\omega $.

The above notion of convexity and the notion of semi-saturated partial representations, as rephrased by \cite
{CondForSSat}, may be bridged as follows:

\state Proposition \label SemiSatVsConvex Let $G$ be a group and let $\SSat $ be the set consisting of one relation of
the form $$ \e {gh} - \e {gh}\e g=0, $$ for each pair of elements $g$ and $h$ in $G$ satisfying $\ell (gh)=\ell (g)+\ell
(h)$.  Then $$ \ORSat = \big \{\omega \in \OuG : \omega \hbox { is convex}\big \}.  $$

\Proof Notice that the corresponding set ${\cal F}_{\SSat }$ of functions on $\OG $ is formed by the functions $$
f(\omega ) = \bool {gh\in \omega } - \bool {gh\in \omega }\bool {g\in \omega } , $$ for all $g,h\in G$, such that $\ell
(gh)=\ell (g)+\ell (h)$.  For each such function, each $\omega \in \OuG $, and each $k\in \omega $, notice that $$
f(k\inv \omega ) = \bool {gh\in k\inv \omega } - \bool {gh\in k\inv \omega }\bool {g\in k\inv \omega } \$= \bool {kgh\in
\omega } - \bool {kgh\in \omega }\bool {kg\in \omega }, $$ Letting $\Phi $ and $\Psi $ denote the formulas ``$kgh\in
\omega $'' and ``$kg\in \omega $'', respectively, the above may be shortened to $$ f(k\inv \omega ) = \bool {\Phi } -
\bool {\Phi }\bool {\Psi }.  $$

We next claim that the Boolean value of the formula ``$\Phi \Rightarrow \Psi $'' is given by $1- f(k\inv \omega )$.  In
fact, using well known properties of Boolean operators we have $$ \def \cra {\cr \pilar {12pt}} \eqalign { \bool {\Phi
\Rightarrow \Psi } & = \bool {\neg \Phi \vee \Psi } \cra & = \bool {\neg \Phi } +\bool {\Psi } - \bool {\neg \Phi }
\bool {\Psi } \cra & = 1-\bool {\Phi } +\bool {\Psi } - (1-\bool {\Phi })\bool {\Psi } \cra & = 1-\bool {\Phi } +\bool
{\Psi } - \bool {\Psi } +\bool {\Phi }\bool {\Psi } \cra & = 1-\bool {\Phi } +\bool {\Phi }\bool {\Psi }\cra & = 1-
f(k\inv \omega ).  } $$

Therefore, to say that $\omega $ is in $\ORSat $ is the same as saying that for all $g$ and $h$ in $G$ with $\ell
(gh)=\ell (g)+\ell (h)$, one has that the Boolean value of ``$\Phi \Rightarrow \Psi $'' is equal to 1, obviously meaning
that $\Phi $ implies $\Psi $, so $$ (\forall k\in \omega )\quad kgh\in \omega \Rightarrow kg\in \omega .  $$

One may easily show that the most general situation in which an element $x$ lies in a segment $\segt kz$ is when $x=kg$,
and $z=kgh$, where $g$ and $h$ are elements in $G$ such that $\ell (gh) =\ell (g) + \ell (h)$.  Thus the above condition
for $\omega $ to lie in $\ORSat $ is precisely saying that $\omega $ is convex.  \endProof

The following is a direct consequence of \cite {CstarRelAsCP}:

\state Corollary Let $G$ be a group equipped with a length function $\ell $.  Then the set $$ \Omega _{\conv } = \big
\{\omega \in \OuG : \omega \hbox { is convex}\big \} $$ is a closed subspace of\/ $\OuG $, invariant under the partial
Bernoulli action.  Moreover, denoting by $\Th $ the corresponding restricted partial action, one has that $ C(\Omega
_{\conv })\Rt {\Th } G $ is the universal C*-algebra for semi-saturated partial representations in the following sense:
\izitem \zitem The map $$ g\in G \mapsto 1_g\delta _g\in C(\Omega _{\conv })\rt {\Th } G, $$ where $1_g$ is the
characteristic function of the set $ \{\omega \in \Omega _{\conv }: g\in \omega \}, $ is a semi-sa\-tu\-rated *-partial
representation.  \zitem Given any semi-saturated *-partial representation of\/ $G$ in a unital C*-algebra $B$, there
exists a unique *-homomorphism $$ \varphi : C(\Omega _{\conv })\rt {\Th } G \to B, $$ such that $\varphi (1_g\delta _g)
= \pr g$.

\nrem Most results of this chapter are taken from \ref {ExelLacaQuigg/2002}.

\syschapter {Hilbert modules and \RME lence}{Hilbert modules and MR~equivalence}

 \font \bigmath = cmr12 \def \bigstar {{\raise -3pt \hbox {\kern -3pt \bigmath *}}}

One of the main tools to study partial crossed products in the realm of C*-algebras is the theory of Hilbert modules.
In this short chapter we will therefore outline some of the main results from that theory which we will need in the
sequel.  The reader is referred to \ref {JensenThomsen/1991}, \ref {BrownMigoShen/1994} and \ref {Lance/1995} for
careful treatments of this important subject.

Hilbert modules are also crucial in defining the concept of \RME lence between C*-algebras and between C*-algebraic
dynamical systems, which we will also attempt to briefly discuss, while referring the reader to \ref {Rieffel/1974} for
a more extensive treatment and examples.

\definition \label DefineHilbMod Let $A$ be a C*-algebra.  By a \subj {right pre-Hilbert $A$-module} we mean a complex
vector space $M $, equipped with the structure of a right $A$-module as well as an \subj {$A$-valued inner-product} $$
\langle \ponto ,\ponto \rangle : M \times M \to A, $$ satisfying \izitem \zitem $\langle \xi ,\lambda \eta +\eta
'\rangle = \lambda \langle \xi ,\eta \rangle + \langle \xi ,\eta '\rangle $, \zitem $\langle \xi ,\xi \rangle \geq 0$,
\zitem $\langle \xi ,\xi \rangle =0 \Rightarrow \xi =0$, \zitem $\langle \xi ,\eta a\rangle =\langle \xi ,\eta \rangle
a$, \zitem $\langle \xi ,\eta \rangle =\langle \eta ,\xi \rangle ^*$, \medskip \noindent for every $\xi ,\eta ,\eta '\in
M $, $\lambda \in {\bf C}$, and $a\in A$.

Notice that it follows from \cite {DefineHilbMod/i} and \cite {DefineHilbMod/v} that $\langle \ponto ,\ponto \rangle $
is conjugate-linear in the first variable.  Also, from \cite {DefineHilbMod/iv} and \cite {DefineHilbMod/v}, one has
that $$ \langle a\xi ,\eta \rangle =a^*\langle \xi ,\eta \rangle \for \xi ,\eta ,\in M \for a\in A.  $$

It is a well known fact \ref {Proposition 2.3/Paschke/1973} that the expression $$ \Vert \xi \Vert _2=\Vert \langle \xi
,\xi \rangle \Vert ^{1/2} \for \xi \in M , $$ defines a norm on $M $, and we then say that $M $ is a \subj {right
Hilbert $A$-module} if $M $ is complete relative to this norm.  When $M $ is not complete, the Banach space completion
of $M $ may be shown to carry the structure of a right Hilbert $A$-module extending the one of $M$.

Contrary to the Hilbert space case, a bounded $A$-linear operator $T$ between two Hilbert modules $M$ and $N$, does not
necessarily admit an adjoint.  By an adjoint of $T$ we mean an operator $T^*:N\to M$ satisfying the familiar property $$
\langle T(\xi ),\eta \rangle = \langle \xi ,T^*(\eta )\rangle \for \xi \in M\for \eta \in N .  $$

For a counter-example, let $A$ be a C*-algebra and let us view $A$ as a Hilbert module over itself with inner-product
defined by $$ \langle a,b\rangle = a^*b \for a,b\in A.  $$

Every ideal $J\ideal A$ is a sub-Hilbert module of $A$ and the inclusion $\iota :J\to A$ is an isometric $A$-linear map,
hence bounded.  However the adjoint of $\iota $ may fail to exist under certain circumstances.  Suppose, for instance,
that $A$ is unital and $J$ is not.  If $\iota ^*$ exists then for every $x$ in $J$ we have $$ x^* = x^*1 = \langle \iota
(x),1\rangle = \langle x,\iota ^*(1)\rangle = x^*\iota ^*(1).  $$ This implies that $\iota ^*(1)$ is a unit for $J$,
contradicting our assumptions.  Therefore $\iota ^*$ does not exist.

When studying Hilbert module one therefore usually restrict attention to the \subjex {adjointable operators}
{adjointable operator}, meaning the operators which happen to have an adjoint.

The set of all adjointable operators on a Hilbert module $M$ is denoted $\AdjOp (M)$.  One may prove without much
difficulty that $\AdjOp (M)$ is a C*-algebra with respect to the composition of operators, the adjoint defined above,
and the operator norm.

Given a C*-algebra $A$, one may likewise define the concept of \subj {left pre-Hilbert $A$-module}, the only differences
relative to Definition \cite {DefineHilbMod} being that $M$ is now assumed to be a left $A$-module and axioms \cite
{DefineHilbMod/i\&iv} become \Item {(i')} $\langle \lambda \xi +\xi ',\eta \rangle = \lambda \langle \xi ,\eta \rangle +
\langle \xi ',\eta \rangle $, \Item {(iv')} $\langle a\xi ,\eta \rangle =a\langle \xi ,\eta \rangle $, \medskip
\noindent for every $\xi ,\xi ',\eta $ in $M $, $\lambda $ in ${\bf C}$, and $a$ in $ A$.  The notion of \subj {left
Hilbert $A$-module} is defined in the obvious way.

\state Lemma \label ApproxIdentityHMod Let $A$ be a C*-algebra and let $M$ be a right \resp {left} Hilbert $A$-module.
If $\{v_i\}_i$ is an approximate identity for $A$, then for every $\xi $ in $M$, one has that $\xi = \lim _i \xi v_i$
\resp {$\xi = \lim _i v_i\xi $}.

\Proof We have $$ \Vert \xi -\xi v_i\Vert ^2 = \Vert \langle \xi -\xi v_i,\xi -\xi v_i\rangle \Vert \$= \Vert \langle
\xi ,\xi \rangle -\langle \xi ,\xi \rangle v_i -v_i^*\langle \xi ,\xi \rangle + v_i^*\langle \xi ,\xi \rangle v_i\Vert
\$\leq \Vert \langle \xi ,\xi \rangle -\langle \xi ,\xi \rangle v_i\Vert + \Vert v_i\Vert \Vert \langle \xi ,\xi \rangle
v_i - \langle \xi ,\xi \rangle \Vert \convrg i 0.  $$ A similar argument proves the left-handed version.  \endProof

As a consequence of the previous result, we see that $M=\clspan {MA}$ (closed linear span), but we can in fact get a
slightly more precise result:

\state Lemma \label Roots Let $A$ be a C*-algebra and let $M$ be a right \resp {left} Hilbert $A$-module. Then, for
every $\xi $ in $M$ one has that $$ \xi = \lim _{n\to \infty }\xi \langle \xi ,\xi \rangle ^{1/n} \quad (\hbox {resp.\ }
\xi = \lim _{n\to \infty }\langle \xi ,\xi \rangle ^{1/n}\xi ).  $$ Consequently $M = \Clspan {M\langle M,M\rangle }$
\resp {$M = \Clspan {\langle M,M\rangle M}$}.

\Proof We have $$ \big \langle \xi - \xi \langle \xi ,\xi \rangle ^{1/n}, \xi - \xi \langle \xi ,\xi \rangle ^{1/n}\big
\rangle \$= \langle \xi ,\xi \rangle - \langle \xi ,\xi \rangle \langle \xi ,\xi \rangle ^{1/n} - \langle \xi ,\xi
\rangle ^{1/n}\langle \xi ,\xi \rangle + \langle \xi ,\xi \rangle ^{1/n}\langle \xi ,\xi \rangle \langle \xi ,\xi
\rangle ^{1/n} \$= \langle \xi ,\xi \rangle - 2\langle \xi ,\xi \rangle ^{1+1/n} + \langle \xi ,\xi \rangle ^{1+2/n} \
\>{n\to \infty } 0.  $$ Therefore $\Vert \xi - \xi \langle \xi ,\xi \rangle ^{1/n}\Vert \>{n\to \infty } 0$, concluding
the proof in the right module case, while a similar argument proves the left-handed version.  \endProof

It is interesting to notice that every right Hilbert $A$-module $M$ is automatically a left Hilbert module over $\AdjOp
(M)$ as follows: given $T$ in $\AdjOp (M)$ and $\xi $ in $M$, we evidently define the multiplication of $T$ by $\xi $ by
$$ T\ponto \xi := T(\xi ), $$ thus providing the left module structure.  As for the $\AdjOp (M)$-valued inner-product,
given $\xi $ and $\eta $ in $M$, consider the mapping $$ \Omega _{\xi ,\eta }: \zeta \in M \mapsto \xi \langle \eta
,\zeta \rangle \in M.  $$

It is easy to see that $\Omega _{\xi ,\eta }$ is an adjointable operator with $\Omega _{\xi ,\eta }^*=\Omega _{\xi ,\eta
}$.  We thus define the $\AdjOp (M)$-valued inner-product on $M$ by $$ \ip {\AdjOp (M)} \xi \eta = \Omega _{\xi ,\eta }
\for \xi ,\eta \in M.  $$

Observe that the two inner-products defined on $M$ satisfy the compatibility condition $$ \ip {\AdjOp (M)}\xi \eta \zeta
= \xi \ip B \eta \zeta \for \xi ,\eta ,\zeta \in M, $$ which motivates our next concept.

\definition \label DefineHBimod Given C*-algebras $A$ and $B$, by a \subj {Hilbert $A$-$B$-bi\-module} we mean a left
Hilbert $A$-module $M$, which is also equipped with the structure of a right Hilbert $B$-module such that, denoting by
$\ip A\ponto \ponto $ and $\ip B\ponto \ponto $ the $A$-valued and $B$-valued inner-products, respectively, one has that
\izitem \zitem $(a\xi )b = a(\xi b)$.  \zitem $\ip A \xi \eta \zeta = \xi \ip B \eta \zeta $, \medskip \noindent for all
$a\in A$, $b\in B$, and $\xi ,\eta ,\zeta \in M$.

It is possible to show \ref {Corollary 1.11/BrownMigoShen/1994} that the norms originating from the two inner-products
on a Hilbert bimodule agree, meaning that $$ \Vert \ip A \xi \xi \Vert = \Vert \ip B \xi \xi \Vert \for \xi \in M.  $$

\definition Let $M$ be a Hilbert $A$-$B$-bimodule.  \iaitem \aitem We say that $M$ is \subjex {right \resp {left}
full}{full Hilbert module}, if the linear span of the range of $\ip A \ponto \ponto $ \resp {$\ip B\ponto \ponto $} is
dense in $A$ \resp {$B$}.  \aitem If $M$ is both left and right full, we say that $M$ is an \subj {imprimitivity
bimodule}.  \aitem If there exists an imprimitivity $A$-$B$-bimodule, we say that $A$ and $B$ are \subjex {\RME
lent}{\RME lent C*-algebras}.

A rather common situation in which a \RME lence takes place is as follows: suppose that $B$ is a C*-algebra and $A$ is a
closed *-sub\-algebra of $B$.

Recall that $A$ is said to be a \subj {hereditary subalgebra} if $$ ABA\subseteq A.  $$ This is equivalent to saying
that any element $b$ in $B$, such that $0\leq b\leq a$, for some $a$ in $A$, necessarily satisfies $b\in A$.

On the other hand, $A$ is said to be a \subj {full subalgebra} if the ideal generated by $A$ in $B$ coincides with $B$,
namely $$ \clspan {BAB} = B.  $$

If $A$ is a hereditary subalgebra of $B$, the right ideal of $B$ generated by $A$, namely $M=\clspan {AB}$, becomes a
Hilbert $A$-$B$-bimodule with inner-products defined by $$ \ip A xy = xy^* \and \ip B xy = x^*y \for x,y\in M.  $$

Since $A\subseteq M$, it is evident that $M$ is left-full.  If we moreover suppose that $A$ is a full subalgebra in the
above sense, then $M$ is also right-full, hence an imprimitivity bimodule, whose existence tells us that $A$ and $B$ are
\RME lent.  Thus:

\state Proposition \label FullHeredMorita If $A$ is a full, hereditary, closed *-subalgebra of a C*-algebra $B$, then
$A$ and $B$ are \RME lent.

\RME lent C*-algebras share numerous interesting properties relating to representation theory and $K$-theory \ref
{Rieffel/1974}, \ref {BrownGreenRieffel/1977}.  One of the most striking results in this field states that if $A$ and
$B$ are separable\fn {In fact this result is known to hold under less stringent conditions, namely if $A$ and $B$
possess strictly positive elements.  However most of the applications we have in mind are to the separable case, so we
will not bother to deal with non-separable algebras here, unless the separability condition is irrelevant.}
C*-algebras, then they are \RME lent if and only if they are \subj {stably isomorphic}, meaning that $$ \Kp \otimes
A\simeq \Kp \otimes B, $$ where $\Kp $ is the algebra of compact operators on a separable infinite dimensional Hilbert
space (see \ref {BrownGreenRieffel/1977} for more details).

The concept of \RME lence has an important counterpart for C*-algebraic partial dynamical systems:

\definition \label RmorPDsys Let $G$ be a group and suppose that for each $k=1,2$ we are given a C*-algebraic partial
dynamical system $$ \Th ^k = \big (A^k,\ G,\ \{A_g^k\}_{g\in G},\ \{\th g^k\}_{g\in G}\big ).  $$ We will say that $\Th
^1$ and $\Th ^2$ are \subjex {\RME lent}{\RME lent partial actions} if there exists a Hilbert $A^1$-$A^2$-bimodule $M$,
and a (set-theoretical) partial action $$ \gamma = \big (\{M_g\}_{g\in G},\ \{\gamma _g\}_{g\in G}\big ) $$ of $G$ on
$M$, such that for all $k=1,2$, and all $g$ in $G$, one has that \izitem \zitem $M_g$ is a norm closed,
sub-$A^1$-$A^2$-bimodule of $M$, \zitem $A_g^k = \clspan {\ip {A^k}{M_g}{M_g}}$ (closed linear span), \zitemBismark
FullIdealsInMor \zitem $\gamma _g$ is a complex linear map, \zitem $\bip {A^k}{\gamma _g(\xi )}{\gamma _g(\eta )} = \th
g^k\big (\ip {A^k}\xi \eta \big )$, for all $\xi ,\eta \in M_{g\inv }$. \zitemBismark TringCovar \medskip \noindent In
this case we say that $$ \gamma = \big (M,\ G,\ \{M_g\}_{g\in G},\ \{\gamma _g\}_{g\in G}\big ) $$ is an \subj
{imprimitivity system} for $\Th ^1$ and $\Th ^2$.

Speaking of \cite {TringCovar}, observe that $\ip {A^k}\xi \eta $ lies in the domain of $\th g^k$ by \cite
{FullIdealsInMor}.

Using \cite {TringCovar}, and the fact that C*-algebra automorphisms are necessarily isometric, one sees that the
$\gamma _g$ must be isometric as well.

Notice that each $M_g$ may be seen as a full Hilbert $A_g^1$-$A_g^2$-bimodule, so necessarily $A_g^1$ and $A_g^2$ are
\RME lent.  This applies in particular to $g=1$, so $A^1$ and $A^2$ must be \RME lent as well.

\state Proposition \label NotAdjoinGama If $\Th ^1$ and $\Th ^2$ are \RME lent partial dynamical systems, and $\gamma $
is an imprimitivity system for $\Th ^1$ and $\Th ^2$, as in \cite {RmorPDsys}, then, given $\xi $ in $M_{g\inv }$, one
has that $$ \gamma _g(a\xi ) = \th g^1(a)\gamma _g(\xi ) \for a\in \Di g^1, $$ and $$ \gamma _g(\xi b) = \gamma _g(\xi
)\, \th g^2(b) \for b\in \Di g^2.  $$

\Proof Focusing on the first assertion, let $\eta \in M_{g\inv }$.  Then by \cite {TringCovar} we have $$ \ip
{A^1}{\gamma _g(a\xi )}{\gamma _g(\eta )} = \th g^1\big (\ip {A^1}{a\xi }\eta \big ) = \th g^1(a)\th g^1\big (\ip
{A^1}\xi \eta \big ) \$= \th g^1(a) \ip {A^1}{\gamma _g(\xi )}{\gamma _g(\eta )} = \bip {A^1}{\th g^1(a) \gamma _g(\xi
)}{\gamma _g(\eta )}.  $$

Since $\eta $ is arbitrary, and since $\gamma _g(\eta )$ can take on any value in $M_g$, we conclude that $$ \gamma
_g(a\xi )=\th g^1(a) \gamma _g(\xi ).  $$ The second assertion is proved similarly.  \endProof

In view of the above result, it is not reasonable to expect the $\gamma _g$ to be bi-module maps (hence not adjointable
either) relative to any of the available Hilbert module structures.  However, if we see each $M_g$ as a \"{ternary
C*-ring} \ref {Zettl/1983} under the ternary operation $$ \{\xi ,\eta ,\zeta \}_g:= \xi \ip {A^2} \eta \zeta , $$ it
easily follows from \cite {NotAdjoinGama} that $\gamma _g$ is an isomorphism of ternary C*-rings.

\bigskip One of our long term goals, unfortunately not to be achieved too soon, will be proving that \RME lent partial
actions lead to \RME lent crossed products.  For this it is important to introduce the concept of linking algebra.

\medskip \fix So, suppose for the time being that $A$ and $B$ are C*-algebras, and that $M$ is a Hilbert
$A$-$B$-bimodule.

We first define the \subj {adjoint Hilbert bimodule} as follows.  Let $M^*$ be any set admitting a bijective function $$
\xi \in M\mapsto \xi ^*\in M^*.  $$ We then define the structures of vector space, left $B$-module, and right $A$-module
on $M^*$ by $$ \xi ^*+\lambda \eta ^* = (\xi +\bar \lambda \eta )^*, $$$$ b\xi ^* = (\xi b^*)^*, $$$$ \xi ^*a = (a^*\xi
)^*, $$ for all $\xi $ and $\eta $ in $M$, $\lambda $ in ${\bf C}$, $a$ in $A$ and $b$ in $B$.  We further define $A$-
and $B$-valued inner-products on $M^*$ by $$ \ip A{\xi ^*}{\eta ^*} = \ip A{\xi }{\eta } \and \ip B{\xi ^*}{\eta ^*} =
\ip B{\xi }{\eta }, $$ for all $\xi $ and $\eta $ in $M$.

The reader might be annoyed by the fact that, unlike module structures, the inner-products defined on $M^*$ have not
changed relative to the original ones.  However this turns out to be the only sensible choice, and one will quickly be
convinced of this by checking the appropriate axioms, which the reader is urged to do.  Once this is done $M^*$ becomes
a Hilbert $B$-$A$-bimodule.

We then write $$ L = \pmatrix {A & M \cr M^* & B}, $$ simply meaning the cartesian product $A \times M \times M^* \times
B$, denoted in a slightly unusual way.

We make $L$ a complex-vector space in the obvious way, and a *-algebra by introducing multiplication and adjoint
operations as follows: $$ \pmatrix {a_1 & \xi _1 \cr \eta ^*_1 & b_1} \pmatrix {a_2 & \xi _2 \cr \eta ^*_2 & b_2} =
\pmatrix {a_1 a_2 + \ip A{\xi _1}{\eta _2} & a_1\xi _2 + \xi _1b_2 \cr \eta ^*_1a_2 + b_1\eta _2^*& \ip B {\eta _1}{\xi
_2} + b_1b_2}, $$ and $$ \pmatrix {a & \xi \cr \eta ^* & b}^\bigstar = \pmatrix {a^* & \eta \cr \xi ^* & b^*}.  $$

In order to give $L$ a norm, we define a representation $\pi _B$ of $L$ as adjointable operators on the right Hilbert
$B$-module $M\oplus B$ (where $B$ is seen as a right Hilbert $B$-module in the obvious way) by $$ \pi _B\pmatrix {a_1 &
\xi _1 \cr \eta ^*_1 & b_1} \pmatrix {\xi _2 \cr b_2} = \pmatrix {a_1\xi _2 + \xi _1b_2 \cr \ip B {\eta _1}{\xi _2} +
b_1b_2}, $$ and another representation $\pi _A$ on the right Hilbert $A$-module $A\oplus M^*$, by $$ \pi _A \pmatrix
{a_1 & \xi _1 \cr \eta ^*_1 & b_1} \pmatrix {a_2 \cr \eta ^*_2 } = \pmatrix {a_1 a_2 + \ip A{\xi _1}{\eta _2} \cr \eta
^*_1a_2 + b_1\eta _2^*}.  $$

We then define a norm on $L$ by $$ \Vert c\Vert = \max \big \{\Vert \pi _A(c)\Vert ,\Vert \pi _B(c)\Vert \big \} \for
c\in L, $$ and one may then prove that $L$ becomes a C*-algebra with the above structure \ref {Proposition
2.3/BrownMigoShen/1994}.

\definition Given C*-algebras $A$ and $B$, and a Hilbert $A$-$B$-bimo\-dule $M$, the \subj {linking algebra} of $M$ is
the C*-algebra $L$ described above.

The linking algebra relative to the imprimitivity bimodule implementing a \RME lence between partial actions also
carries a partial action, as we will now show:

\state Proposition \label LinkPaction Let $G$ be a group and suppose we are given C*-alge\-braic partial dynamical
systems $$ \alpha = \big (A,G,\{A_g\}_{g\in G},\{\alpha _g\}_{g\in G}\big ) \and \beta = \big (B,G,\{B_g\}_{g\in
G},\{\beta _g\}_{g\in G}\big ).  $$ Suppose, moreover, that $\alpha $ and $\beta $ are \RME lent partial actions, and
that $$ \gamma = \big (M,\ G,\ \{M_g\}_{g\in G},\ \{\gamma _g\}_{g\in G}\big ) $$ is an imprimitivity system for $\alpha
$ and $\beta $.  Letting $L$ be the linking algebra of $M$, for each $g$ in $G$, we have that \izitem \zitem the subset
$L_g$ of $L$ defined by $$ L_g = \pmatrix {A_g & M_g \cr M^*_g & B_g} $$ is a closed two-sided ideal, \zitem the mapping
$$ \lambda _g: \pmatrix {a & \xi \cr \eta ^* & b} \in L_{g\inv }\ \longmapsto \ \pmatrix {\alpha _g(a) & \gamma _g(\xi )
\cr \gamma _g(\eta )^* & \beta _g(b)} \in L_g.  $$ is a *-isomorphism, \zitem the pair $$ \lambda = \big (\{L_g\}_{g\in
G},\ \{\lambda _g\}_{g\in G}\big ), $$ is a C*-algebraic partial action of $G$ on $L$.

\Proof We first claim that if $\xi \in M$, and $\eta \in M_g$, then $$ \ip A\xi \eta \in A_g \and \ip B\xi \eta \in B_g.
\equationmark FunnyIdealPropOne $$ In order to see this, let $\{v_i\}_i$ be an approximate identity for $B_g$.  Viewing
$M_g$ as a right Hilbert $B_g$-module, we get from \cite {ApproxIdentityHMod} that $\eta = \lim _i \eta v_i$, so $$ \ip
B\xi \eta = \ip B\xi {\lim _i \eta v_i}= \lim _i \ip B\xi { \eta }v_i \in B_g.  $$ Using a similar reasoning one shows
that $\ip A\xi \eta \in A_g$.

Of course, if it is $\xi $, rather than $\eta $, which belongs to $M_g$, then \cite {FunnyIdealPropOne} still hold (by
taking adjoints).

We will need another fact of a similar nature, that is, if $a$ is in $A_g$, $b$ is in $B_g$, and $\xi $ is in $M$, then
$$ a\xi \in M_g \and \xi b\in M_g.  \equationmark FunnyIdealPropTwo $$ To prove it, notice that by \cite
{FullIdealsInMor}, we may assume that $b = \ip B\eta \zeta $, for some $\eta $ and $\zeta $ in $M_g$.  In this case we
have $$ \xi b = \xi \ip B\eta \zeta \={DefineHBimod/ii} \ip A\xi \eta \zeta \in M_g.  $$ The proof that $a\xi \in M_g$
is similar.

We therefore have that, in any one of the expressions: $$ ab,\ a\xi ,\ \xi b,\ \ip A\xi \eta \ \hbox {and}\ \ip B\xi
\eta , $$ where $a$ is in $A$, $b$ is in $B$, and $\xi $ and $\eta $ are in $M$, if one of the two terms involved
belongs to a set named with a subscript ``$g$'', then so does the whole expression.

Keeping this principle in mind, and staring at the definition of the multiplication operation in $L$ for a while, one
then sees that $L_g$ is in fact a two-sided ideal of $L$.

In \ref {Proposition 2.3/BrownMigoShen/1994} it is proven that the norm topology of $L$ coincides with the product
topology, when $L$ is viewed as $A \times M \times M^* \times B$.  Thus, since $A_g$ is closed in $A$, $B_g$ is closed
in $B$ and $M_g$ is closed in $M$, we see that $L_g$ is closed in $L$.

The proof that each $\gamma _g$ is a *-isomorphism is an easy consequence of \cite {TringCovar} and \cite
{NotAdjoinGama}, and the corresponding property for $\alpha _g$ and $\beta _g$.

Finally, the last point is easily verified by observing that $\lambda $ is the direct sum of four partial actions.
\endProof

This result will later be used to prove that \RME lent partial actions lead to \RME lent crossed products.  The strategy
for doing so will be to form the crossed product of $L$ by $G$, and then make sense of the expression $$ L\rt \lambda G
= \pmatrix {A\rt \alpha G & M\rt \gamma G \cr (M\rt \gamma G)^* & B\rt \beta G}, $$ so $M\rt \gamma G$ will be seen to
be an imprimitivity bimodule implementing a \RME lence between $A\rt \alpha G$ and $B\rt \beta G$.  The biggest hurdle
we will face, and the reason we need to further develop our theory, is showing $A\rt \alpha G$ and $B\rt \beta G$ to be
isomorphic to subalgebras of $L\rt \lambda G$.  Even though we may easily find *-homomorphisms $$ A\rt \alpha G \to L\rt
\lambda G \and B\rt \beta G \to L\rt \lambda G, $$ proving these to be injective requires some extra work.  As soon as
we have the appropriate tools, we will return to this point and we will prove the following:

\state Theorem \label FutureRME If $$ \alpha = \big (A,G,\{A_g\}_{g\in G},\{\alpha _g\}_{g\in G}\big ) \and \beta = \big
(B,G,\{B_g\}_{g\in G},\{\beta _g\}_{g\in G}\big ) $$ are \RME lent C*-algebraic partial dynamical systems, then $A\rt
\alpha G$ and $B\rt \beta G$ are \RME lent C*-algebras.

\nrem Inspired by Kaplansky's C*-modules, Hilbert modules were introduced by Paschke in \ref {Paschke/1973}.  The notion
of \RME lence, sometimes also referred to as \subj {strong Morita equivalence}, was introduced by Rieffel in \ref
{Rieffel/1974}, adapting the purely algebraic concept of Morita equivalence \ref {Morita/1958} to operator algebras.

\RME lence of actions of groups on C*-algebras were considered independently in \ref {CurtoMurphyWilliams/1984} and \ref
{Combes/1984} and, for the case of partial actions, in \ref {Abadie/1999} and \ref {Abadie/2003}, where a version of
Theorem \cite {FutureRME} for reduced crossed products appeared.  When we are ready to prove \cite {FutureRME}, we will
take care of both the reduced and full versions.

\partpage {II}{FELL BUNDLES}

\chapter Fell bundles

\def \Bu #1{B^\Pr _{#1}}

If one wishes to understand the structure of a given algebra, a traditional method is to try to decompose it in the
direct sum of ideals.  However, when working with simple algebras, such a decomposition is evidently unavailable.

An alternative method is to try to find a grading (see \cite {AlgGradedALg}) of our algebra, in which case we will have
decomposed it in its \"{homogeneous spaces}.  These are of course subspaces of the given algebra but they can also be
seen as separate entities.  In other words, starting from a graded algebra we may see the collection of its homogeneous
spaces as the parts we are left with after \"{disassembling} the algebra along its grading.  The study of the separate
pieces might then hopefully shed some light on the structure of our algebra.

\def \.{.\kern 2pt} From now on we will employ the concept of Fell bundles, introduced by J\.M\.G\.Fell under the name
of \"{C*-algebraic bundles}, in order to deal with disassembled C*-algebras\fn {It is important to stress that, over
topological groups, Fell bundles carry in its topology a lot more information than we will be using here.  However,
since all of our groups are discrete, this extra topological data will not be relevant for us here.}.  However it is
crucial to note that, in the category of C*-algebras, the concept of a grading, to be defined shortly, only requires the
direct sum of the homogeneous spaces to be \"{dense}, so the process of passing from a graded C*-algebra to its
corresponding Fell bundle, that is, the collection formed by its homogeneous spaces, involves a significant loss of
information: there is no straightforward process to reassemble a graded C*-algebra from its parts!  In other words,
there are examples of non-isomorphic graded C*-algebras whose associated Fell bundles are indistinguishable.

While this can be considered a weakness of the method, it often helps to organize one's tasks in two broad groups.  In
order to understand the structure of a graded C*-algebra one should therefore attempt to separately understand:
\smallskip \itemitem {(1)} the structure of its associated Fell bundle, and \smallskip \itemitem {(2)} the way in which
the various parts are pieced together.

\smallskip \noindent Breaking up assignments along these lines frequently makes a lot of sense because the two tasks
often require very different tool sets.

In this and the forthcoming chapters we plan to provide both sets of tools for dealing with graded C*-algebras: (1)
understanding the structure of Fell bundles will be done by proving that, under suitable hypothesis, every Fell bundle
arises from a partial dynamical system, while (2) understanding how to piece together the homogeneous spaces to
reassemble the algebra will be done via the theory of cross-sectional algebras and amenability.

\definition \label DefineFellBundle A \subj {Fell bundle} (also known as a \subj {C*-algebraic bundle}) over a group $G$
is a collection $$ \Bun = \{B_g\}_{g\in G} $$ of Banach spaces, each of which is called a \subj {fiber}.  In addition,
the \subj {total space} of $\Bun $, namely the disjoint union of all the $B_g$'s, which we also denote by $\Bun $, by
abuse of language, is equipped with a multiplication operation and an involution $$ \cdot :\Bun \times \Bun \to \Bun ,
\qquad \qquad * : \Bun \to \Bun , $$ satisfying the following properties for all $g$ and $h$ in $G$, and all $b$ and $c$
in $\Bun $: \iaitem \aitem $B_gB_h\subseteq B_{gh}$, \aitem multiplication is bi-linear from $B_g\times B_h$ to
$B_{gh}$, \aitem multiplication on $\Bun \,$ is associative, \aitem $\Vert bc\Vert \leq \Vert b\Vert \Vert c\Vert $,
\aitem $(B_g)^* \subseteq B_{g\inv }$, \aitem involution is conjugate-linear from $B_g$ to $B_{g\inv }$, \aitem $(bc)^*
= c^* b^*$, \aitem $b^{**}=b$, \aitem $\Vert b^*\Vert = \Vert b\Vert $, \aitem $\Vert b^*b\Vert = \Vert b\Vert ^2$,
\aitem $b^*b \geq 0$ in $B_1$. \aitemBismark PositivityInFell

\medskip Axioms (a--d) above define what is known as a \subj {Banach algebraic bundle}.  Adding (e--i) gives the
definition of a \subj {Banach *-algebraic bundle}.

Observe that axioms (a--j) imply that $B_1$ is a C*-algebra with the restricted operations.  We will often refer to
$B_1$ as the \subj {unit fiber algebra}.

With respect to axiom \cite {PositivityInFellLoc}, notice that the reference to positivity there is to be taken with
respect to the standard order relation in $B_1$, seen as a C*-algebra.  Should one prefer to avoid any reference to this
order relation, an alternative formulation of \cite {PositivityInFellLoc} is to require that for each $b$ in $\Bun $,
there exists some $a$ in $B_1$, such that $b^*b=a^*a$.

\medskip As already hinted upon, a concept which is closely related to Fell bundles is the notion of graded C*-algebras,
which is not equivalent, and hence should be distinguished from its purely algebraic counterpart \cite {AlgGradedALg}.

\definition \label DefineGradedCstar Let $B$ be a C*-algebra and $G$ be a group.  We say that a linearly independent
collection $ \{B_g\}_{g\in G} $ of closed subspaces of $B$ is a \"{C*-grading} for $B$, if $\bigoplus _{g\in G} B_g$ is
dense in $B$, and for every $g$ and $h$ in $G$, one has that \zitemno =0 \zitem $B_g B_h \subseteq B_{gh}$, \zitem
$B_g^* \subseteq B_{g\inv }$.  \medskip \noindent In this case we say that $B$ is a \subj {$G$-graded C*-algebra} and
each $B_g$ is called a \subjex {grading subspace}{grading subspace in a C*-algebra}.

The reason why the above is not a special case of \cite {AlgGradedALg} is that here the direct sum of the $B_g$'s is
only required to be dense in $B$.

\nostate \label Disassemble If one is given a $G$-graded C*-algebra $B$, then the collection of all grading subspaces
will clearly form a Fell bundle with the norm, the multiplication operation and the adjoint operation borrowed from $B$.

We will see that, conversely, every Fell bundle may be obtained from a $G$-graded C*-algebra $B$, as above, although $B$
is not uniquely determined since there might be many ways to complete the direct sum of the $B_g$'s.

One of the reasons why one is interested in studying Fell bundles rather than graded C*-algebras is to avoid getting
distracted by the problems related to the different completions mentioned above.

\nostate \label GroupBundle An important example of Fell bundles is given by the \subj {group bundle} $$ \Bun = {\bf
C}\times G, $$ where $G$ is any group.  Each fiber $B_g$ is defined to be ${\bf C}\times \{g\}$, with usual linear and
norm structure, and with operations $$ (\lambda ,g) (\mu ,h) = (\lambda \mu ,gh) \and (\lambda ,g)^* = (\bar \lambda
,g\inv ), $$ for all $\lambda ,\mu \in {\bf C}$, and $g,h\in G$.

\nostate One of the main reasons we are interested in Fell bundles is because they may be built from C*-algebraic
partial dynamical systems.  In order to describe this construction, let us fix a C*-algebraic partial action $$ \Th =
\big (\{\D g\}_{g\in G},\ \{\th g\}_{g\in G}\big ) $$ of a group $G$ on a C*-algebra $A$.  We will begin the
construction of our Fell bundle by defining its total space to be $$ \Bun = \{(b,g) \in A\times G: b\in \D g\}.  $$

Inspired by the notation introduced in \cite {Polynomial}, we will write $b\delta _g$ to refer to $(b,g)$, whenever
$b\in \D g$.  We may then identify the fibers of our bundle as $$ B_g = \{b\delta _g: b\in \D g\}.  $$

The linear structure and the norm on each $B_g$ is borrowed from $\D g$, while the multiplication operation is defined
exactly as in \cite {ProdInParCrossProd}, namely $$ (a \delta _g) (b \delta _h) = \th g \big (\thi g(a)b\big ) \delta
_{gh} \for a\in \D g \for b\in \D h.  $$ We finally borrow the definition of the involution from \cite
{StarInCrossProd}, namely $$ (a\delta _g)^* = \thi g(a^*)\delta _{g\inv } \for g\in G \for a\in \D g.  $$

\state Proposition \label BundleFromPA Given a C*-algebraic partial action $$ \Th = \big (\{\D g\}_{g\in G},\ \{\th
g\}_{g\in G}\big ) $$ of a group $G$ on a C*-algebra $A$, one has that $\Bun $, with the above operations, is a Fell
bundle over $G$, henceforth called the \subj {semi-direct product bundle} relative to $\Th $.

\Proof The task of checking the axioms in \cite {DefineFellBundle} is mostly routine.  We restrict ourselves to a few
comments.  With respect to associativity, we observe that it follows by the same argument used in \cite
{TheoremOnAssoc}.  Axiom (d), (i) and (j), namely the axioms referring to the norm structure, all follow easily from the
fact that partial automorphisms of C*-algebras are isometric.  Axioms (g) and (h) may be proved following the ideas used
in the proof of \cite {StarInCrossProd}, while (k) easily follows from \cite {LotsLefInnProd} and the fact that partial
automorphisms preserve positivity.  \endProof

\nostate \label BunFromPrep Another important class of examples of Fell bundles is given in terms of partial
representations.  In order to describe these, let $$ \Pr : G \to A $$ be a *-partial representation of a given group $G$
in a unital C*-algebra $A$.  For each $g$ in $G$, consider the closed linear subspace $\Bu g\subseteq A$, spanned by the
elements of the form $$ \pr {h_1}\pr {h_2}\cdots \pr {h_n}, $$ where $n\geq 1$ is any integer, the $h_i$ are in $G$, and
$$ h_1\cdots h_n = g.  $$

It is elementary to check that $$ (\Bu g)^*=\Bu {g\inv } \and \Bu g\Bu h\subseteq \Bu {gh}, $$ for all $g$ and $h$ in
$G$, so the collection $$ \Bun ^\Pr = \{\Bu g\}_{g\in G} $$ is seen to be a Fell bundle with the operations borrowed
from $A$.

\state Proposition \label BundleFromPR Let $\Pr $ be a *-partial representation of the group $G$ in the unital
C*-algebra $A$.  For each $g$ in $G$, let $ \e g = \pr g \pri g, $ as usual.  Then, regarding the Fell bundle $\Bun ^\Pr
$ defined above, we have: \izitem \zitem the unit fiber algebra $\Bu 1$ coincides with the C*-algebra generated by the
set $\{\e g:g\in G\}$, \zitem $\Bu 1$ is abelian, \zitem for each $g$ in $G$, one has that $\pr g \Bu 1 = \Bu g = \Bu
1\pr g$.

\Proof Given $h_1,\ldots ,h_n$ in $G$, we claim that $$ \pr {h_1}\pr {h_2}\cdots \pr {h_n} = \e {p_1}\e {p_2}\cdots \e
{p_n}\pr {p_n}, \equationmark ProdPrepProdIdemp $$ where $$ p_k = h_1\cdots h_k \for 1\leq k\leq n.  $$

If $n=1$, we have $$ \pr {h_1}= \pr {h_1}\pr {h_1\inv }\pr {h_1} = \e {h_1}\pr {h_1} = \e {p_1}\pr {p_1}, $$ proving the
claim in this case.  Assuming that $n\geq 2$, we have $$ \pr {h_1}\pr {h_2} \pr {h_3}\cdots \pr {h_n} = \pr {h_1} \pr
{h_1\inv } \pr {h_1} \pr {h_2} \pr {h_3}\cdots \pr {h_n} \={ObserveLeft} $$$$ = \pr {h_1} \pr {h_1\inv } \pr {h_1h_2}
\pr {h_3}\cdots \pr {h_n} = \e {h_1} \pr {h_1h_2} \pr {h_3}\cdots \pr {h_n}, $$ and the claim follows by induction.

If we moreover assume that $h_1\cdots h_n = 1$, so that $\pr {h_1}\pr {h_2}\cdots \pr {h_n}$ is an arbitrary generator
of $\Bu 1$, we deduce from the claim that $$ \pr {h_1}\pr {h_2}\cdots \pr {h_n} = \e {h_1}\e {h_1h_2}\cdots \e
{h_1h_2\cdots h_n}, $$ proving (i), and then (ii) follows from \cite {MainTheEgCommute}.

Let us next prove (iii).  Observing that $\pr g$ lies in $\Bu g$, it is immediate that $\pr g \Bu 1 \subseteq \Bu g$.
On the other hand, notice that if $h_1\cdots h_n = g$, and $b$ is defined by $ b=\pr {h_1}\pr {h_2}\cdots \pr {h_n}, $
then $$ b\ei g = b\pri g \pr g \={ProdPrepProdIdemp} \e {p_1}\e {p_2}\cdots \e {p_n}\pr {p_n} \pri g \pr g \={PrgIsPiso}
$$$$ = \e {p_1}\e {p_2}\cdots \e {p_n}\pr {p_n} = b.  $$

Since $\Bu g$ is the closed linear span of the set of elements $b$, as above, we then see that $$ b\ei g = b \for b\in
\Bu g.  $$ For each such $b$, we then have that $$ b = b\ei g = b \pri g \pr g \in \Bu 1 \pr g, $$ because $b \pri g$ is
clearly in $\Bu 1$.  This proves that $\Bu g = \Bu 1 \pr g$, and the remaining statement in (iii) follows by taking
adjoints.  \endProof

It is interesting to notice that the Fell bundle associated to a given *-partial representation, as discussed above, may
also be described as a semi-direct product bundle.  In fact, one may adapt \cite {SystemFromPrep} to the context of
C*-algebras, obtaining a partial action of $G$ on $\Bu 1$, whose associated semi-direct product bundle is isomorphic to
$\Bun ^\Pr $.  We leave the verification of these statements as an easy exercise.

\medskip \fix From now on we fix an arbitrary Fell bundle $\Bun = \{B_g\}_{g\in G}$.

\state Lemma \label ApproximateIdentity If $\{v_i\}_i$ is an approximate identity for $B_1$, then, for every $b$ in any
$B_g$, one has that $$ b = \lim _i bv_i = \lim _i v_i b.  $$

\Proof Noticing that $B_g$ is a Hilbert $B_1$-$B_1$-bimodule in a natural way, the result follows from \cite
{ApproxIdentityHMod}.  \endProof

Using square brackets to denote closed linear span, notice that for any two fibers $B_g$ and $B_h$ of $\Bun $, we have
that $\clspan {B_gB_h}$ is a closed linear subspace of $B_{gh}$.  Questions related to whether or not $\clspan {B_gB_h}$
coincide with $B_{gh}$ will have a special significance for us so let us introduce two concepts which are based on this.
Recall from \cite {DefineSemiSatPAct} and \cite {DefineSemiSatPRept} that the term \"{semi-saturated} has already been
defined both in the context of partial actions and of partial representation.  We will now extend it to the context of
Fell bundles.

\definition \label DefineSemiSatFellBun Let $\Bun = \{B_g\}_{g\in G}$ be a Fell bundle.  \iaitem \aitem We say that
$\Bun $ is \subjex {saturated}{saturated Fell bundle} if $\clspan {B_gB_h} = B_{gh}$, for every $g$ and $h$ in $G$.
\aitem If $G$ is moreover equipped with a length function $\ell $, we say that $\Bun $ is \subjex
{semi-saturated}{semi-saturated Fell bundle} (with respect to the given length function $\ell $) if, for all $g$ and $h$
in $G$ satisfying $\ell (gh)=\ell (g)+\ell (h)$, one has that $\clspan {B_gB_h} = B_{gh}$.

As already seen, the unit fiber $B_1$ of a Fell bundle is always a C*-algebra, and it is easy to see that for each $g$
in $G$, one has that $\clspan {B_gB_{g\inv }}$ is an ideal of $B_1$.

\state Lemma \label ApproximateIdentityTwo Let $\Bun = \{B_g\}_{g\in G}$ be a Fell bundle, and let $g\in G$.  Given an
approximate identity $\{u_i\}_i$ for $\clspan {B_gB_{g\inv }}$, and another approximate identity $\{v_i\}_i$ for
$\clspan {B_{g\inv }B_g}$, one has that $$ b = \lim _i u_ib = \lim _i bv_i \for b\in B_g.  $$

\Proof This is an immediate consequence of \cite {ApproxIdentityHMod}, once we notice that $B_g$ is a Hilbert $\clspan
{B_gB_{g\inv }}$-\null $\clspan {B_{g\inv }B_g}$-\null bimodule.  \endProof

\state Lemma \label FellIsTRO For every $g$ in $G$, one has that $\clspan {B_gB_{g\inv }B_g} = B_g$.

\Proof Given $b$ in $B_g$, choose an approximate identity $\{u_i\}_i$ for $\clspan {B_gB_{g\inv }}$.  Then $$ b = \lim
_i u_ib \in \clspan {B_gB_{g\inv }B_g}, $$ proving that $B_g \subseteq \clspan {B_gB_{g\inv }B_g}$.  The reverse
inclusion is obvious.  \endProof

We may use this to give a simpler characterization of saturated Fell bundles.

\state Proposition A necessary and sufficient condition for a Fell bundle $\Bun = \{B_g\}_{g\in G}$ to be saturated is
that $$ \clspan {B_gB_{g\inv }} = B_1 \for g\in G.  $$

\Proof That the condition is necessary is obvious.  Conversely, notice that for all $g$ and $h$ in $G$, we have $$
B_{gh} \={FellIsTRO} \clspan {B_{gh} B_{(gh)\inv } B_{gh}} \subseteq \clspan {B_1B_{gh}} = \clspan {B_gB_{g\inv }B_{gh}}
\$\subseteq \clspan {B_g B_h}\subseteq B_{gh}.  $$ Therefore equality holds throughout, proving that $\Bun $ is
saturated.  \endProof

We may easily characterize saturatedness and semi-saturatedness for semi-direct product bundles:

\state Proposition Let $ \Th = \big (\{\D g\}_{g\in G},\ \{\th g\}_{g\in G}\big ) $ be a C*-algebraic partial action of
the group $G$ on the C*-algebra $A$, and let $\Bun $ be the associated semi-direct product bundle.  \izitem \zitem Then
$\Bun $ is saturated if and only if $\Th $ is a global action.  \zitem If $G$ is moreover equipped with a length
function, then $\Bun $ is a semi-saturated Fell bundle (according to \cite {DefineSemiSatFellBun/b}) if and only if $\Th
$ is a semi-saturated partial action (according to \cite {DefineSemiSatPAct}).

\Proof Given $g$ in $G$, we have $$ \clspan {B_gB_{g\inv }} = \clspan {(\D g\delta _g)(\D {g\inv }\delta _{g\inv })}
\={LotsAFormulas} $$$$= \clspan {\D g\th g(\D {g\inv })\delta _1} = \D g\delta _1.  \equationmark BgBgDg $$ So, to say
that $\clspan {B_gB_{g\inv }} = B_1$ is equivalent to saying that $D_g=A$.  This proves (i).

Focusing on (ii), pick $g$ and $h$ in $G$, with $\ell (gh) = \ell (g) + \ell (h)$.  Reasoning as in \cite
{DomainCompos}, it is easy to see that the range of $\th g\circ \th h$ is given by $\D g \cap \D {gh}$.  Since $\th
g\circ \th h$ is a restriction of the bijective mapping $\th {gh}$, it is clear that these two maps coincide if and only
if they have the same range, that is, $$ \th g\circ \th h =\th {gh} \iff \D g \cap \D {gh} = \D {gh} \iff \D {gh}
\subseteq \D g.  $$ Assuming that $\Bun $ is semi-saturated, we have that $B_{gh}=\clspan {B_g B_h}$, so $$ \D
{gh}\delta _1 \={BgBgDg} \clspan {B_{gh}B_{h\inv g\inv }} = \clspan {B_g B_h B_{h\inv g\inv }} \subseteq \clspan {B_g
B_{g\inv }} \={BgBgDg} \D g\delta _1, $$ so $\D {gh} \subseteq \D g$, whence $\th g\circ \th h =\th {gh}$, as seen
above, proving $\Th $ to be semi-saturated.

Conversely, assuming that $\Th $ is semi-saturated, we have $$ B_{gh} \={FellIsTRO} \clspan { B_{gh} B_{h\inv g\inv }
B_{gh}} \={BgBgDg} \Clspan { (\D {gh} \delta _1) B_{gh}} \$\subseteq \Clspan { (\D g \delta _1) B_{gh}} \={BgBgDg}
\clspan { B_g B_{g\inv } B_{gh}} \subseteq \clspan { B_g B_h}, $$ which says that $\Bun $ is semi-saturated.  \endProof

Our next result relates saturatedness for partial representations and for Fell bundles.

\state Proposition Given a *-partial representation $\Pr $ of a group $G$ in a nonzero unital C*-algebra $A$, consider
its associated Fell bundle $\Bun ^\Pr $, as described in \cite {BunFromPrep}.  Then: \izitem \zitem $\Bun ^\Pr $ is a
saturated Fell bundle if and only if $\Pr $ is a unitary group representation.  In this case $\Bun ^\Pr $ is isomorphic
to the group bundle ${\bf C}\times G$.  \zitem $\Bun ^\Pr $ is a semi-saturated Fell bundle if and only if $\Pr $ is a
semi-saturated partial representation.

\Proof Assuming that $\Bun ^\Pr $ is saturated, let $g$ be in $G$.  Then $$ \Bu 1 = \clspan {\Bu g \Bu {g\inv }}
\={BundleFromPR/iii} \clspan {\Bu 1 \pr g \pri g \Bu 1} = \clspan {\e g \Bu 1}.  $$

This is to say that the ideal generated by the idempotent element $\e g$ coincides with the unital algebra $\Bu 1$, but
this may only happen if $\e g=1$.  Since $g$ is arbitrary, this easily implies that $\Pr $ is a unitary group
representation.

Conversely, supposing that $\Pr $ is a unitary group representation, notice that whenever $ h_1\cdots h_n = g, $ one has
that $ \pr {h_1}\pr {h_2}\cdots \pr {h_n} = \pr g, $ which is to say that $\Bu g = {\bf C}\pr g$.  It is then easy to
see that $\Bun ^\Pr $ is isomorphic to the group bundle ${\bf C}\times G$, whence a saturated Fell bundle.

With respect to (ii), pick $g$ and $h$ in $G$, with $\ell (gh) = \ell (g) + \ell (h)$.  Then, assuming $\Pr $ to be
semi-saturated, we have that $\pr g \pr h =\pr {gh}$, so $$ \clspan {\Bu g \Bu h} = \clspan {\Bu 1\pr g \pr h \Bu 1} =
\clspan {\Bu 1\pr {gh} \Bu 1} = \clspan {\Bu 1\Bu {gh}} \={ApproximateIdentity} \Bu {gh}, $$ so $\Bun ^\Pr $ is seen to
be semi-saturated.

In order to prove the converse, given any $g$ in $G$, we claim that $$ \e g b = b \for g\in \Bu g.  $$ To see this,
notice that, if $b\in \Bu g$, we may write $b = \pr g a$, for some $a$ in $\Bu 1$, by \cite {BundleFromPR/iii}, so $$ \e
g b = \e g \pr g a = \pr g \pri g \pr g a = \pr g a = b, $$ proving the claim.  If $\Bun ^\Pr $ is semi-saturated, and
still under the assumption that $\ell (gh) = \ell (g) + \ell (h)$, we then immediately see that $$ \e g b = b \for g\in
\Bu {gh}, $$ because $\Bu {gh} = \clspan {\Bu g\Bu h}$.  In particular, since $\pr {gh}$ is in $\Bu {gh}$, we deduce
that $$ \e g \e {gh} = \e g \pr {gh} \pri {gh} = \pr {gh} \pri {gh} = \e {gh}.  $$

Using \cite {CondForSSat} we then conclude that $\pr {gh} = \pr g \pr h$, thus verifying that $\Pr $ is semi-saturated.
\endProof

\medskip \fix Let us now return to the above situation in which $\Bun = \{B_g\}_{g\in G}$ denotes an arbitrary, fixed
Fell bundle.

\medskip Every element $c$ of $B_1$ defines \"{left and right multiplication operators} $$ L_c: b\in B_g \mapsto cb \in
B_g \and R_c: b\in B_g \mapsto bc \in B_g.  $$ We will next show that the above may be extended to multipliers of $B_1$.

\state Proposition \label MultOnBun Given any $m$ in the multiplier algebra $\Mult (B_1)$, there is a unique pair of
bounded linear maps $$ L_m,R_m: B_g \to B_g, $$ such that, for all $a$ in $B_1$, and all $b$ in $B_g$, one has that
\izitem \zitem $L_m(ba) = L_m(b)a$, \zitem $R_m(ab) = aR_m(b)$, \zitem $a L_m(b) = (am)b$, \zitem $R_m(b)a = b(ma)$.
\medskip \noindent Moreover, if $\{u_i\}_i$ is any approximate identity for $B_1$, then $$ L_m(b) = \lim _i (u_im)b \and
R_m(b) = \lim _i b(mu_i) \for b\in B_g.  $$

\Proof Let $\{v_i\}_i$ be an approximate identity for $B_1$.  Addressing uniqueness, notice that if $b$ is in $B_g$, we
have that $$ L_m(b) \={ApproximateIdentity} \lim _i v_iL_m(b) \explica ={iii} \lim _i (v_im)b.  $$ One similarly proves
that $R_m$ is unique.  As for existence, we claim that the limits $$ \lim _i (v_im)b \and \lim _i b(mv_i) $$ exist.
Checking the Cauchy condition relative to the first limit we have, for all $i$ and $j$, that $$ \Vert (v_im)b -
(v_jm)b\Vert ^2 = \Vert \big ((v_im)b - (v_jm)b\big )\big (b^*(m^*v_i) - b^*(m^*v_j)\big )\Vert \$= \Vert
(v_im)bb^*(m^*v_i) - (v_im)bb^*(m^*v_j) - (v_jm)bb^*(m^*v_i) + (v_jm)bb^*(m^*v_j)\Vert \$\leq \Vert v_im(bb^*)m^*v_i -
v_im(bb^*)m^*v_j\Vert + \Vert v_jm(bb^*)m^*v_i - v_jm(bb^*)m^*v_j\Vert \$\leq \Vert v_i\Vert \Vert m(bb^*)m^*v_i -
m(bb^*)m^*v_j\Vert \${\ +\ } \Vert v_j\Vert \Vert m(bb^*)m^*v_i - m(bb^*)m^*v_j\Vert \convrg {i,j} 0.  $$ This proves
the existence of the first limit, while the second limit is shown to exist by a similar argument.  We may then define
$L_m$ and $R_m$ as in the last sentence of the statement, and the proofs of (i--iv) are now mostly routine, so we
restrict ourselves to proving (iii), leaving the proofs of the other points to the reader.  Given $a$ in $B_1$, and $b$
in $B_g$, the left-hand-side of (iii) equals $$ a L_m(b) = a\lim _i (v_im)b = \lim _i (av_im)b = (am)b.
\omitDoubleDollar \endProof

The following result gives a concrete way to describe the \"{fiber-wise multipliers} given by \cite {MultOnBun} in the
case of semi-direct product bundles:

\state Proposition \label MultOnSDPBun Let $\Bun $ be the semi-direct product bundle associated to a given C*-algebraic
partial dynamical system $$ \Th = \big (A,\ G,\ \{\D g\}_{g\in G},\ \{\th g\}_{g\in G}\big ).  $$ Then for every
multiplier $m\in \Mult (A)$, and every $a$ in any $\D g$, one has that $$ L_m(a\delta _g) = (ma)\delta _g \and
R_m(a\delta _g) = \th g\big (\thi g(a)m\big )\delta _g.  $$

\Proof Picking an approximate identity $\{v_i\}_i$ for $A$, we have that $\{v_i\delta _1\}_i$ is an approximate identity
for $B_1=A\delta _1$, so for any $a$ in $\D g$, we have that $$ L_m(a\delta _g) = \lim _i (v_im\delta _1)a\delta _g =
\lim _i (v_ima)\delta _g = (ma)\delta _g, $$ proving the first identity.  As for the second one, we have $$ R_m(a\delta
_g) = \lim _i (a\delta _g)(mv_i\delta _1) = \lim _i \th g\big (\thi g(a)mv_i\big )\delta _g \$= \th g\big (\thi
g(a)m\big )\delta _g.  \omitDoubleDollar \endProof

The reader is urged to compare the above formula for $R_m(a\delta _g)$ with the multiplication $(a\delta _g)(m\delta
_1)$, should $m$ be an element of $A$.

\medskip We will now present the first process of assembling a C*-algebra from a given Fell bundle.  This will produce
the biggest possible outcome.

\medskip \fix Recall that $\Bun = \{B_g\}_{g\in G}$ denotes a Fell bundle which has been fixed near the beginning of
this chapter.

\definition By a \subjex {section}{section of a Fell bundle} of $\Bun $ we shall mean any function $y $ from $G$ to the
total space of $\Bun $, such that $y_g \in B_g$, for every $g$ in $G$.  We will moreover denote by $\CCB $ the
collection of all finitely supported sections.  Given two sections $y $ and $z $ in $\CCB $, we define their \subj
{convolution product} by $$ (y \star z )_g = \soma {h\in G} y_hz_{h\inv g} \for g\in G.  $$ We also define an \"{adjoint
operation} by $$ y ^*_g = (y_{g\inv })^* \for g\in G.  $$

\medskip With respect to the convolution product above, notice that for every $h$ in $G$ we have $$ y_h z_{h\inv g} \in
B_hB_{h\inv g}\subseteq B_g, $$ so all of the summands lie in the same vector space, hence the sum is well defined as an
element of $B_g$.

We leave it for the reader to prove the following easy fact:

\state Proposition With the above operations one has that $\CCB $ is an associative *-algebra.

As an example, notice that for the semi-direct product bundle $\Bun $ described in \cite {BundleFromPA}, we have that
$\CCB $ is isomorphic to $A\art \Th G$.

\definition \label DefineRepBundle A *-representation of a Fell bundle $\Bun = \{B_g\}_{g\in G}$ in a *-algebra $C$ is a
collection $\pi =\{\pi _g\}_{g\in G}$ of linear maps $$ \pi _g: B_g \to C, $$ such that \izitem \zitem $\pi _g(b)\pi
_h(c) = \pi _{gh}(bc)$, \zitem $\pi _g(b)^* = \pi _{g\inv }(b^*)$, \medskip \noindent for all $g,h\in G$, and all $b\in
B_g$, and $c\in B_h$.

Given a Fell bundle $\Bun $, for each $g$ in $G$ we will denote by $j_ g$ the natural inclusion of $B_g$ in $\CCB $,
namely the map $$ j_ g:B_g\to \CCB $$ defined by $$ j_g(b)|_h = \left \{\matrix { b, & \hbox { if } h=g, \cr \pilar
{12pt} 0, & \hbox { otherwise.}  }\right .  \equationmark DefineJg $$

\state Proposition \label FormulasForJ The collection of maps $j = \{j_ g\}_{g\in G}$ is a *-representa\-tion of $\Bun $
in $\CCB $.

\Proof Left to the reader. \endProof

We would now like to construct a C*-algebra from $\CCB $ in the same way we built the C*-algebraic partial crossed
product from its algebraic counterpart.  For this we need to obtain a bound for C*-seminorms on $\CCB $, much like we
did in \cite {CPAdmissible}.

\state Proposition \label CCBAdmissible Let $p$ be a C*-seminorm on $\CCB $. Then, for every $y $ in $\CCB $, one has
that $$ p(y ) \leq \soma {g\in G}\Vert y_g \Vert .  $$

\Proof Using \cite {FormulasForJ} we have for every $b\in B_g$, that $$ p\big (j_g(b)\big )^2 = p\big (j_g(b)^*
j_g(b)\big ) = p\big (j_{g\inv }(b^*) j_g(b)\big )= p\big (j_1(b^*b)\big ).  $$

Observe that the composition $p\circ j_ 1$ is clearly a C*-seminorm on $B_1$.  So, as already mentioned in \cite
{CStarSeminormDominated}, we have $$ p\big (j_1(b^*b)\big )\leq \Vert b^*b\Vert = \Vert b\Vert ^2.  $$ The conclusion is
thus that $$ p\big (j_g(b)\big ) \leq \Vert b\Vert .  $$ For a general $y $ in $\CCB $, we may write $ y = \soma {g\in
G}j_g(y_g), $ so $$ p(y ) \leq \soma {g\in G}p\big (j_g(y_g)\big ) \leq \soma {g\in G}\Vert y_g \Vert .  $$ This
concludes the proof.  \endProof

\bigskip Given $y $ in $\CCB $, define $$ \maxnorm y = \mathop {\rm sup}\limits _p\ p(y ), \equationmark
DefMaxNormInFell $$ where $p$ ranges in the collection of all C*-seminorms on $\CCB $.  By the above result we see that
$\maxnorm \ponto $ is finite for every $y$, and hence it is a well defined C*-seminorm on $\CCB $.

\definition \label DefineFullCrossSect The \subj {cross sectional C*-algebra} of $\Bun $, denoted $\CB $, is the
C*-algebra obtained by taking the quotient of $\CCB $ by the ideal consisting of the elements of norm zero and then
completing it.

We will soon see that $\maxnorm \ponto $ is in fact a norm, and hence that $\CCB $ may be embedded as a dense subalgebra
of $\CB $.  Meanwhile we will call the canonical mapping arising from the completion process by $$ \kappa : \CCB \to \CB
.  $$

Since $\kappa $ is a *-homomorphism, and since $j$ is a *-representation, the following is evident:

\state Proposition \label DefineHatJota For each $g$ in $G$, let us denote by $\hj _g$ the composition $$ B_g \> {j_g}
\CCB \> \kappa \CB .  $$ Then $\hj = \{\hj _ g\}_{g\in G}$ is a *-representation of $\Bun $ in $\CB $, henceforth called
the \subjex {universal representation}{universal representation of a Fell bundle}.

We have already seen that $B_1$ is always a C*-algebra and it is easy to see that $\hj _1$ is necessarily a
*-homomorphism.  In fact we have:

\state Proposition \label BOneEssential For every Fell bundle $\Bun $, the map $$ \hj _1:B_1\to \CB $$ is a {\nonDegHomo
} *-homomorphism.

\Proof Follows immediately from \cite {ApproximateIdentity}.  \endProof

The result stated below is simply the recognition that the C*-algebraic partial crossed product defined in \cite
{DefineCStarCP} is a special case of the cross sectional C*-algebra.

\state Proposition \label CPIsCrossSectAlg Let $$ \Th = \big (\{\D g\}_{g\in G},\ \{\th g\}_{g\in G}\big ) $$ be a
C*-algebraic partial action of a group $G$ on a C*-algebra $A$, and let $\Bun $ be the corresponding semi-direct product
bundle.  Then $A\rt \Th G$ is naturally isomorphic to $\CB $.

In the case of the group bundle ${\bf C}\times G$, the cross sectional C*-algebra is denoted $C^*(G)$ and is called the
\subj {group C*-algebra} of $G$.

\medskip

The cross sectional C*-algebra of a Fell bundle possesses the following universal property with respect to bundle
representations:

\state Proposition \label FromRepBunToRelAlg Let $\pi =\{\pi _g\}_{g\in G}$ be a *-representation of the Fell bundle
$\Bun $ in a C*-algebra $C$.  Then there exists a unique *-homomorphism $$ \varphi :\CB \to C, $$ such that $$ \varphi
\big (\hj _ g(b)\big ) = \pi _g(b) \for g\in G\for b\in B_g.  $$ We will say that $\varphi $ is the \subj {integrated
form} of $\pi $.

\Proof Given $\pi $, define $ \varphi _0: \CCB \to C, $ by $$ \varphi _0(y) = \soma {g\in G} \pi _g(y_g ) \for y\in \CCB
.  $$ It is a routine matter to prove that $\varphi _0$ is a *-homomorphism, whence the formula $$ p(y ) = \Vert \varphi
_0(y )\Vert \for y \in \CCB , $$ defines a C*-seminorm on $\CCB $, which is therefore dominated by $\maxnorm \ponto $.
This says that $\varphi _0$ is continuous and so it extends to the completion, providing the required map.  \endProof

\nrem Fell bundles were introduced in the late 60's by J.~M.~G. Fell in \ref {FellDoran/1969}, under the name
\"{C*-algebraic bundles}, setting the stage for a far reaching generalization of Harmonic Analysis.  See \ref
{FellDoran/1988} for a very extensive treatment of this and other important concepts.

Fell bundles over topological groups are defined in a different fashion since the topology of the group must also be
taken into account.  However, since we have chosen to work only with discrete group, we do not need to worry about any
extra topological conditions.

\chapter Reduced cross-sectional algebras

The construction of the cross-sectional C*-algebra for a Fell bundle given in \cite {DefineFullCrossSect} shares with
\cite {DefineCStarCP}, as well as many other universal constructions, a norm defined via a supremum (see \cite
{DefMaxNormInFell} and \cite {DefineMaxNorm}).  While we have so far verified that these maximal norms are finite, we
have not worried about them being zero!  In other words, we are facing the potentially tragic situation in which the
cross-sectional algebra of nontrivial Fell bundles may turn out to be the zero algebra!  By \cite {CPIsCrossSectAlg}
this would also affect C*-algebraic partial crossed products.

In order to rule out this undesirable situation, given a Fell bundle $\Bun $, we must provide nontrivial C*-seminorms on
$\CCB $, and one of the best ways to do so is through representation theory.  That is, given any representation $\rho $
of $\CCB $ on a Hilbert space, we may define a seminorm by $p(y) = \Vert \rho (y)\Vert $, which will be nontrivial as
long as $\rho $ is nontrivial.

However, instead of representing $\CCB $ on a Hilbert space, our representation will be on a Hilbert module over the
unit fiber algebra $B_1$, which we now set out to construct.  As a byproduct we will describe our second process of
assembling a C*-algebra from a given Fell bundle.

\medskip \fix From now on we fix an arbitrary Fell bundle $\Bun = \{B_g\}_{g\in G}$.

\medskip By \cite {FormulasForJ} we have that $j_ 1$ is a *-homomorphism through which we may view $B_1$ as a subalgebra
of $\CCB $.  This makes $\CCB $ a right $B_1$-module in a standard way and we will now introduced a $B_1$-valued
inner-product on $\CCB $ as follows $$ \langle y, z\rangle = \soma {g\in G} (y_g)^*z_g \for y ,z \in \CCB .  $$

The easy verification that this is indeed an inner-product is left to the reader.  Once this is done we have that $\CCB
$ is a right pre-Hilbert $B_1$-module.

\definition \label DefineEllTwo We shall denote by $\ell ^2(\Bun )$ the right Hilbert $B_1$-module obtained by
completing $\CCB $ under the norm $\Vert \ponto \Vert _2$ arising from the inner-product defined above.

Seeing $\CCB $ as a dense subspace of $\ell ^2(\Bun )$, we may view $j_g$ as a map from $B_g$ to $\ell ^2(\Bun )$, and
it is interesting to remark that this is an isometry, since, for every $b$ in any $B_g$, one has $$ \Vert j_g(b)\Vert
_2^2 = \Vert \langle j_g(b),j_g(b)\rangle \Vert = \Vert b^*b\Vert = \Vert b\Vert ^2.  $$

Our next step will be the construction of a representation of $\Bun $ within $\AdjOp \big (\ell ^2(\Bun )\big )$.  The
following technical fact will be useful in that direction:

\state Lemma \label OrderNorm Given $g,h\in G$, $b\in B_g$, and $c\in B_h$, one has that $$ c^*b^*bc\leq \Vert b\Vert
^2c^*c.  $$

\Proof It is well known that a positive element $h$ in a C*-algebra $A$ satisfies $h\leq \Vert h\Vert $, in the sense
that $v^*(\Vert h\Vert -h)v\geq 0$, for every $v$ in $A$.  Here, we either temporarily work in the unitization of $A$,
or we simply write $v^*(\Vert h\Vert -h)v$ to actually mean $\Vert h\Vert v^*v - v^*hv$.

Picking an approximate identity $\{v_i\}_i$ for $B_1$, we then have for all $i$ that, $$ v_i^*(\Vert b\Vert ^2 -
b^*b)v_i\geq 0.  $$ So, there exists some $a_i$ in $B_1$ such that $$ v_i^*(\Vert b\Vert ^2 - b^*b)v_i = a_i^*a_i.  $$

Applying \cite {PositivityInFell} for $a_ic\in B_h$, and using \cite {ApproximateIdentity}, we get $$ 0\leq (a_ic)^*a_ic
= c^*a_i^*a_ic = c^*v_i^*(\Vert b\Vert ^2 - b^*b)v_i c \ \convrg {i} \ c^*(\Vert b\Vert ^2 - b^*b)c \$= \Vert b\Vert
^2c^*c - c^*b^*bc, $$ concluding the proof.  \endProof

\state Proposition \label PresentLambdaG Given $b$ in any $B_g$, the operator $$ \lambda _g(b) : \ y \ \in \ \CCB \
\longmapsto \ j_g(b)\star y \ \in \ \CCB $$ is continuous relative to the Hilbert module norm on $\CCB $, and hence
extends to a bounded operator on $\ell ^2(\Bun )$, still denoted by $\lambda _g(b)$ by abuse of language, such that
$\Vert \lambda _g(b)\Vert = \Vert b\Vert $, and which moreover satisfies $$ \lambda _g(b)\big (j_h(c)\big ) = j_{gh}(bc)
\for h\in G \for c\in B_h.  $$

\Proof Notice that for every $y$ in $\CCB $ and every $h$ in $G$, the formula for the convolution product gives $$
\lambda _g(b)y |_h = b y_{g\inv h}.  $$ Therefore, if $c$ is in $B_h$, then for every $k$ in $G$, we have that $$
\lambda _g(b)\big (j_h(c)\big )|_k = b j_h(c)_{g\inv k} = \delta _{h,g\inv k}bc = \delta _{gh,k}bc = j_{gh}(bc)_k, $$
proving the last assertion in the statement.  Addressing the boundedness of $\lambda _g(b)$, given $y$ in $\CCB $,
notice that $$ \langle \lambda _g(b)y ,\lambda _g(b)y \rangle = \soma {h\in G} (y_{g\inv h})^* b^* b y_{g\inv h} \$=
\soma {h\in G} (y_h)^* b^* b y_h \_\leq {OrderNorm} \Vert b\Vert ^2\soma {h\in G} (y_h)^* y_h = \Vert b\Vert ^2\langle y
,y \rangle .  $$ This implies that $\Vert \lambda _g(b)y \Vert _2\leq \Vert b\Vert \Vert y \Vert _2$, from where we see
that $\lambda _g(b)$ is bounded and $\Vert \lambda _g(b)\Vert \leq \Vert b\Vert $.

In order to prove the reverse inequality let $\xi $ be the element of $\ell ^2(\Bun )$ obtaining by mapping $b^*$ into
$\ell ^2(\Bun )$ via $j_{g\inv }$.  Recalling that $j_{g\inv }$ is an isometry, we have that $\Vert \xi \Vert _2 = \Vert
b\Vert $.  We then have $$ \lambda _g(b)\xi = \lambda _g(b)j_{g\inv }(b^*) = j_g(b)\star j_{g\inv }(b^*) = j_1(bb^*), $$
so $$ \Vert b\Vert ^2 = \Vert b^*b\Vert = \Vert j_1(bb^*)\Vert _2 = \Vert \lambda _g(b)\xi \Vert _2 \leq \Vert \lambda
_g(b)\Vert \Vert \xi \Vert _2 = \Vert \lambda _g(b)\Vert \Vert b\Vert .  $$ This proves that $\Vert \lambda _g(b)\Vert
\geq \Vert b\Vert $, completing the proof.  \endProof

We may now finally introduce the first nontrivial representation of our bundle.

\state Proposition \label DefineRegularRep The collection of maps $$ \lambda =\{\lambda _g\}_{g\in G} $$ introduced
above is a representation of $\Bun $ in $\AdjOp \big (\ell ^2(\Bun )\big )$, henceforth called the \subjex {regular
representation}{regular representation of a Fell bundle} of $\Bun $.

\Proof For $b$ in any $B_g$, and $y ,z $ in $\CCB $, we have $$ \langle \lambda _g(b)y ,z \rangle = \soma {h\in G}
(y_{g\inv h})^*b^*z_h = \soma {h\in G} (y_h)^*b^*z_{gh} = \langle y ,\lambda _{g\inv }(b^*)z \rangle .  $$ Since $\CCB $
is dense in $\ell ^2(\Bun )$, we conclude from the above that $$ \langle \lambda _g(b)\xi ,\eta \rangle = \langle \xi
,\lambda _{g\inv }(b^*)\eta \rangle \for \xi ,\eta \in \ell ^2(\Bun ), $$ and hence that $\lambda _g(b)$ is an
adjointable operator, with $\lambda _g(b)^*=\lambda _{g\inv }(b^*)$.  This also proves \cite {DefineRepBundle/ii}.
Given $b_g$ in $B_g$, $b_h$ in $B_h$, and $y $ in $\CCB $, we have $$ \lambda _g(b_g)\big (\lambda _h(b_h)y \big )=
j_g(b_g)\star j_h(b_h)\star y \={FormulasForJ} $$$$= j_{gh}(b_gb_h)\star y = \lambda _{gh}(b_gb_h)y.  $$ Again because
$\CCB $ is dense in $\ell ^2(\Bun )$, we get $\lambda _g(b_g)\lambda _h(b_h) = \lambda _{gh}(b_gb_h)$, proving \cite
{DefineRepBundle/i} and concluding the proof.  \endProof

By \cite {FromRepBunToRelAlg}, the regular representation of $\Bun $ gives rise to a *-homomor\-phism $$ \Lambda
:C^*(\Bun ) \to \AdjOp \big (\ell ^2(\Bun )\big ), \equationmark RegRepCross $$ satisfying $\Lambda \circ \hj _g =
\lambda _g$, for any $g$ in $G$.

\definition \label DefineRedCrossSect The *-homomorphism $\Lambda $ defined above will be called the \subj {regular
representation of $C^*(\Bun )$}, and its range, which we will denote by $\CrB $, will be called the \subj {reduced cross
sectional C*-algebra} of $\Bun $.

\medskip We thus see that $\CrB $ is isomorphic to the quotient of $C^*(\Bun )$ by the null space of $\Lambda $.

\beginpicture \setcoordinatesystem units <0.0020truecm, 0.0020truecm> point at 0 0 \setplotarea x from -3600 to 4000, y
from -000 to 750 \put {$B_g \>{j_g} \CCB \> \kappa C^*(\Bun ) \trepa \Lambda \twoheadrightarrow \CrB .$} at 000 000
\setquadratic \plot -1650 180 -125 700 1300 150 / \arrow <0.11cm> [0.3,1.2] from 1290 158 to 1300 150 \plot -1550 150
-775 400 000 150 / \arrow <0.11cm> [0.3,1.2] from -10 156 to 000 150 \put {$\hj _g$} at -300 500 \put {$\lambda _g$} at
300 850 \advseqnumbering \put {(\current ) \deflabel {RegRepProp}{\current }} at 3000 300 \endpicture

\bigskip

In the special case of the group bundle ${\bf C}\times G$, the reduced cross sectional C*-algebra is denoted $C^*\lred
(G)$, and it is called the \subj {reduced group C*-algebra} of $G$.

\medskip

The consequences of the existence of the regular representation are numerous.  We begin with a technical remark which
should be seen as a generalization of Fourier coefficients to Fell bundles.

\state Lemma \label Fourier For each $g$ in $G$, there exists a \subjex {contractive}{contractive linear map}\fn {A
linear map $T$ is said to be \"{contractive} if $\Vert T(x)\Vert \leq \Vert x\Vert $, for all $x$ in its domain.} linear
mapping $$ E_g:\CrB \to B_g, $$ such that, for every $h$ in $G$, and every $b$ in $B_h$, one has that $$ E_g\big
(\lambda _h(b)\big ) = \left \{\matrix { b, & \hbox {if } g=h, \cr 0, & \hbox {if } g\neq h.  } \right .  $$ For any
given $z$ in $\CrB $ we will refer to $E_g(z)$ as the $g^{th}$ \subj {Fourier coefficient} of $z$.

\Proof The map $$ P_g: y\in \CCB \to y_g\in B_g, $$ is easily seen to be continuous for the Hilbert module norm, and
hence it extends to a continuous mapping, also denoted $P_g$, from $\ell ^2(\Bun )$ to $B_g$.

Letting $\{v_i\}_i$ be an approximate identity for $B_1$, we then claim that the limit $$ \lim _iP_g \Big (z \big (
j_1(v_i)\big )\Big ) \equationmark LimitForEg $$ exists for every $ z $ in $\CrB $.  In order to prove our claim, let us
first check it for $ z =\lambda _h(b)$, where $h\in G$, and $b\in B_h$.  In this case, for every $i$, we have $$ P_g
\Big (z \big ( j_1(v_i)\big )\Big ) = P_g \Big (\lambda _h(b) \big (j_1(v_i)\big )\Big ) = P_g \big (j_h(b)\star
j_1(v_i)\big ) \$= P_g \big (j_h(bv_i)\big ) = \delta _{g,h}bv_i, $$ where $\delta _{g,h}$ is the Kronecker symbol.
Thus the limit in \cite {LimitForEg} exists and coincides with $\delta _{g,h}b$, by \cite {ApproximateIdentity}.  It
follows that the limit also exists for every $ z $ in $$ \soma {h\in G} \lambda _h(B_h) = \Lambda \big (\kappa (\CCB
)\big ), $$ which is clearly dense in $\CrB $.  Since the operators involved in \cite {LimitForEg} are uniformly
bounded, we conclude that the above limit indeed exists for all $ z $ in $\CrB $.  We may thus define $$ E_g( z ) = \lim
_i P_g \Big (z \big ( j_1(v_i)\big )\Big ) \for z \in \CrB , $$ observing that the calculation above gives $ E_g\big
(\lambda _h(b)) = \delta _{g,h} b, $ for all $b$ in $B_h$.

Finally, notice that since $\Vert P_g\Vert =1$, and $\Vert v_i\Vert \leq 1$, then $$ \Vert E_g( z )\Vert = \lim _i \Big
\Vert P_g \Big (z \big ( j_1(v_i)\big )\Big ) \Big \Vert \leq \Vert z \Vert .  \omitDoubleDollar \endProof

The more technical aspects having been taken care of, we may now reap the benefits:

\state Proposition \label FinallyInjective Let $\Bun $ be a Fell bundle.  Then: \izitem \zitem the map $\kappa :\CCB \to
C^*(\Bun )$ is injective, \zitem the map $\Lambda \circ \kappa :\CCB \to \CrB $ is injective, \zitem there exists a
C*-seminorm on $\CCB $ that is a norm, \zitem for every $g$ in $G$, the map $\hj _g:B_g\to C^*(\Bun )$ is isometric,
\zitem for every $g$ in $G$, the map $\lambda _g:B_g\to \CrB $ is isometric, \zitem $C^*(\Bun )$ is a graded C*-algebra,
with grading subspaces $\hj _g(B_g)$, \zitem $\CrB $ is a graded C*-algebra, with grading subspaces $\lambda _g(B_g)$.

\Proof (ii) Given $y\in \CCB $ such that $\Lambda \big (\kappa (y)\big )=0$, write $y = \soma {h\in G}j_h(y_h)$.  Then,
for every $g$ in $G$, we have $$ 0 = E_g\big (\Lambda (\kappa (y))\big ) = \soma {h\in G}E_g\Big (\Lambda \big (\kappa
(j_h(y_h))\big )\Big ) \={RegRepProp} \soma {h\in G}E_g\big (\lambda _h(y_h)\big ) \={Fourier} y_g, $$ hence $y=0$.

\medskip \noindent (i) Follows immediately from (ii).

\medskip \noindent (iii) In view of (i), it is enough to take $p(y) := \maxnorm y = \Vert \kappa (y)\Vert $.

\medskip \noindent (iv) Given $b$ in any $B_g$, on the one hand we have that $$ \Vert \hj _g(b)\Vert = \Vert \kappa \big
(j_g(b)\big )\Vert = \maxnorm {j_g(b)} \_\leq {CCBAdmissible}\Vert b\Vert , $$ and on the other $$ \Vert \hj _g(b)\Vert
\geq \Vert \Lambda \big (\hj _g(b)\big )\Vert \={RegRepProp} \Vert \lambda _g(b)\Vert \={PresentLambdaG} \Vert b\Vert .
$$

\medskip \noindent (v) This is just \cite {PresentLambdaG}.  It is included here for comparison purposes only.

\medskip \noindent (vii) By (v) we have that each $\lambda _g(B_g)$ is a closed subspace of $\CrB $.  To show that they
are independent, suppose that $$ \soma {g\in G} z _g =0, $$ with $ z _g\in \lambda _g(B_g)$, the sum having only
finitely many nonzero terms.  For each $g$, write $ z _g=\lambda _g(y_g)$, with $y_g\in B_g$.  We may then view $y$ as
an element of $\CCB $ and we have $$ \Lambda \big (\kappa (y)\big ) = \soma {g\in G}\Lambda \big (\kappa \big
(j_g(y_g)\big )\big ) = \soma {g\in G}\lambda _g(y_g) = \soma {g\in G} z _g =0, $$ whence $y=0$ by (ii), and
consequently $ z _g = \lambda _g(y_g)=0$, for all $g$ in $G$.  This shows that the $\lambda _g(B_g)$ form an independent
collection of subspaces.  It is then easy to see that $$ \bigoplus _{g\in G} \lambda _g(B_g) = \Lambda \big (\kappa
(\CCB )\big ), \equationmark CCDenseInCred $$ which is clearly dense in $\CrB $.  In order to check the remaining
conditions \cite {DefineGradedCstar/i--ii}, observe that $\lambda $ is a representation by \cite {DefineRegularRep}, so
for all $g,h\in G$, we have that $$ \lambda _g(B_g)\lambda _h(B_g) = \lambda _{gh}(B_{gh}), $$ and $$ \lambda _g(B_g)^*
= \lambda _{g\inv }(B_g^*) = \lambda _{g\inv }(B_{g\inv }).  $$

Point (vi) may be proved in a similar fashion and so it is left for the reader.  \endProof

Recall from \cite {CPIsCrossSectAlg} that, given a C*-algebraic partial action of a group $G$ on a C*-algebra $A$, the
crossed product $A\rt \Th G$ is isomorphic to the cross sectional C*-algebra of the semi-direct product bundle.
Inspired by this we may give the following:

\definition \label DefineRedCP The \subj {reduced crossed product} of a C*-algebra $A$ by a partial action $\Th $ of a
group $G$, denoted $A\redrt \Th G$, is the reduced cross sectional algebra of the corresponding semi-direct product
bundle.

Applying \cite {FinallyInjective} to semi-direct product bundles we obtain:

\state Corollary \label PartCPIsGraded Let $ \Th $ be a C*-algebraic partial action of a group $G$ on a C*-algebra $A$.
Then: \izitem \zitem the map $\kappa :A\art \Th G\to A\rt \Th G$ is injective, \zitem the map $\Lambda \circ \kappa
:A\art \Th G\to A\redrt \Th G$ is injective, \zitem there exists a C*-seminorm on $A\art \Th G$ that is a norm, \zitem
for every $g$ in $G$, and all $a\in \D g$, one has that $\Vert a\delta _g\Vert = \Vert a\Vert $, \zitem for every $g$ in
$G$, and all $a\in \D g$, one has that $\Vert \Lambda (a\delta _g)\Vert = \Vert a\Vert $, \zitem $A\rt \Th G$ is a
graded C*-algebra, with grading subspaces $\D g\delta _g$, \zitem $A\redrt \Th G$ is a graded C*-algebra, with grading
subspaces $\Lambda (\D g\delta _g)$.

We have already briefly indicated that the $E_g$ of \cite {Fourier} are to be thought of as Fourier coefficients.  In
fact, when the bundle in question is the group bundle over ${\bf Z}$, namely ${\bf C}\times {\bf Z}$, then the reduced
cross sectional algebra is isomorphic to the algebra of all continuous functions on the unit circle.  Moreover, for any
such function $f$, one has that $E_n(f)$ coincides with the Fourier coefficient $\hat f(n)$.

We will now present a few results about Fell bundles which are motivated by Fourier analysis and which emphasize further
similarities between these theories.  Our first such result should be compared with the well known fact that the matrix
of a multiplication operator on $L^2(S^1)$, relative to the standard basis, is a Laurent matrix, that is, a doubly
infinite matrix with constant diagonals.

\state Proposition \label MatrixCoef Given $ z \in \CrB $, one has that $$ \big \langle j_g(b), z j_h(c)\big \rangle =
b^*E_{gh\inv }( z )c, $$ for all $g,h\in G$, $b\in B_g$, and $c\in B_h$.

\Proof It is clearly enough to prove the result for $ z =\soma {k\in G} \lambda _k(y_k)$, with each $y_k\in B_k$, the
sum having only finitely many nonzero terms.  In this case we have $$ \big \langle j_g(b), z j_h(c)\big \rangle = \Big
\langle j_g(b),\soma {k\in G} \lambda _k(y_k)j_h(c)\Big \rangle \$= \Big \langle j_g(b),\soma {k\in G} j_k(y_k)\star
j_h(c)\Big \rangle = \Big \langle j_g(b),\soma {k\in G} j_{kh}(y_kc)\Big \rangle \$= b^*y_{gh\inv }c \={Fourier}
b^*E_{gh\inv }( z )c.  \omitDoubleDollar \endProof

Another well known result in Classical Harmonic Analysis states that, if all of the Fourier coefficients of a continuous
function $f$ vanish, then $f=0$.  Our next result should be considered as a generalization of this fact.

\state Proposition \label EOneIsFaith Given $ z \in \CrB $, the following are equivalent: \izitem \zitem $E_1(z^*z)=0$,
\zitem $E_g(z)=0$, for every $g$ in $G$, \zitem $\,z=0$.

\Proof (i) $\Rightarrow $ (iii).  Given $b$ in any $B_g$, we have $$ \langle z j_g(b), z j_g(b)\rangle = \langle j_g(b),
z ^* z j_g(b)\rangle \={MatrixCoef} b^*E_1( z ^* z )b =0, $$ whence $ z j_g(b) =0$.  Since the $j_g(b)$ span a dense
subspace of $\ell ^2(\Bun )$, we have $z =0$.

\medskip \noindent (ii) $\Rightarrow $ (iii).  Given $b$ in any $B_g$, and $c$ in any $B_h$, we have $$ \big \langle
j_g(b), z j_h(c)\big \rangle = b^*E_{gh\inv }( z )c = 0.  $$ Since the $j_g(b_g)$ span a dense subspace of $\ell ^2(\Bun
)$, we see that $ z j_h(c)=0$.  As above, this gives $z =0$.

\medskip \noindent

\medskip \noindent (i) $\Leftarrow $ (iii) $\Rightarrow $ (ii). Obvious.  \endProof

\medskip A very important fact in Classical Harmonic Analysis is Parseval's identity, which asserts that, for a
continuous function $f$ on the unit circle ${\bf T}$ (in fact also for more general functions), one has that $$ \int
_{\bf T}|f(z)|^2\,dz = \soma {n\in {\bf Z}}|\hat f(n)|^2.  $$

We will now present a generalization of this to Fell bundles by seeing the right-hand side above in terms of our version
of Fourier coefficients, while the left-hand-side is interpreted based on the fact that, in the classical case, the
integral of a function coincides with its zeroth Fourier coefficient.

We begin with a preparatory result, the second part of which is a generalization of yet another important fact in
Harmonic Analysis.

\state Lemma \label SomaEmLDois \iaitem \aitem Let $y\in \CCB $, and let $ z = \soma {g\in G}\lambda _g(y_g).  $ Then $$
\soma {g\in G}E_g(z)^*E_g(z) = E_1(z^*z).  $$ \aitem (Bessel's inequality) For every $z\in \CrB $, and for every finite
subset $K\subseteq G$, one has that $$ \soma {g\in K}E_g(z)^*E_g(z) \leq E_1(z^*z).  $$

\Proof Initially observe that if $z$ is as in (a), then by \cite {Fourier} one has that $E_g(z) = y_g$, so the sum in
(a) is actually a finite sum.  We then have $$ E_1(z^*z) = E_1\Big (\soma {g,h\in G}\lambda _g(y_g)^*\lambda _h(y_h)\Big
) \={DefineRepBundle} E_1\Big (\soma {g,h\in G}\lambda _{g\inv h}(y_g^*y_h)\Big ) \={Fourier} $$$$ = \soma {g\in
G}y_g^*y_g = \soma {g\in G}E_g(z)^*E_g(z), $$ proving (a).  In order to prove (b), assume first that $z$ is as in (a).
Then it is clear that for every finite set $K\subseteq G$, one has that $$ E_1(z^*z) \trepa {(\rm a)} = \soma {g\in
G}E_g(z)^*E_g(z) \geq \soma {g\in K}E_g(z)^*E_g(z).  $$ Since the set of elements $z$ considered above is dense in $\CrB
$, and since both the left- and the right-hand expressions above represent continuous functions of the variable $z$, the
result follows by taking limits.  \endProof

We now wish to show that the conclusion in \cite {SomaEmLDois/a} in fact holds for every $z$ in $\CrB $.  Evidently the
sum might now have infinitely many nonzero terms, so we need to discuss the kind of convergence to be expected.  We will
see that unconditional, rather than absolute convergence is the appropriate notion to be used here.

Recall that a series $\soma {i\in I}x_i$ in a Banach space is said to be \subj {unconditionally convergent} with sum $s$
if, for every $\varepsilon >0$, there exists a finite set $F_0\subseteq I$, such that for every finite set $F\subseteq
I$, with $F_0\subseteq F$, one has that $$ \normsum {\kern 2pt s- \soma {i\in F}x_i}< \varepsilon .  $$

This is equivalent to saying that $s$ is the limit of the net $\big \{\soma {i\in F}x_i\big \}_F$, indexed by the
directed set formed by all finite subsets $F\subseteq I$, ordered by inclusion.

\state Proposition \label ParsevalId (Parseval's identity) For every $z$ in $\CrB $, one has that $$ \soma {g\in
G}E_g(z)^*E_g(z) = E_1(z^*z), $$ where the series converges unconditionally.

\Proof Given a finite subset $K\subseteq G$, write $K = \{g_1,g_2,\ldots ,g_n\}$, and let us consider the mapping $$
E_K: z\in \CrB \mapsto \pmatrix { E_{g_1}(z) & 0 & \ldots & 0 \cr \pilar {12pt} E_{g_2}(z) & 0 & \ldots & 0 \cr \vdots &
\vdots & \ddots & \vdots \cr E_{g_n}(z) & 0 & \ldots & 0 } \in M_n\big (\CrB \big ), $$ where, as usual, we are
identifying each $B_g$ as a subspace of $\CrB $ via the regular representation.  It is evident that $E_K$ is a linear
mapping, and we claim that it is contractive.  In fact, given $z$ in $\CrB $, we have that $$ \Vert E_K(z)\Vert ^2 =
\Vert E_K(z)^*E_K(z)\Vert = \normsum {\soma {g\in K}E_g(z)^*E_g(z)} \_\leq {SomaEmLDois/b} $$$$ \leq \Vert
E_1(z^*z)\Vert \_\leq {Fourier} \Vert z^*z\Vert = \Vert z\Vert ^2, $$ proving our claim.  Given $z,z_0\in \CrB $, we
then have that $$ \Vert E_K(z)^*E_K(z) - E_K(z_0)^*E_K(z_0) \Vert \leq \Vert z\Vert \Vert z-z_0\Vert + \Vert z-z_0\Vert
\Vert z_0\Vert , $$ from where one sees that the collection of (non-linear) functions $$ \Big \{E_K(\ponto )^*E_K(\ponto
) : K\subseteq G,\ |K|<\infty \Big \}, $$ is equicontinuous.  Observe that each function in the above set takes values
in the $(1,1)$ upper left corner of our matrix algebra.  Moreover, for every $z$ in $\CrB $, we have that $$ \big
(E_K(z)^*E_K(z)\big )_{1,1} = \soma {g\in K}E_g(z)^*E_g(z).  $$

By \cite {SomaEmLDois/a} we see that, as $K{\uparrow } G$, the expression above converges to $E_1(z^*z)$ on the dense
set formed by the finite sums, as in \cite {SomaEmLDois/a}.  Being equicontinuous, it therefore converges to
$E_1(z^*z)$, for every $z$ in $\CrB $.  \endProof

\nrem The construction of the reduced cross sectional algebra of a Fell bundle was introduced in \ref {Exel/1997b} in
analogy with the reduced crossed products defined for C*-dynamical systems, and McClanahan's reduced partial crossed
products \ref {Section 3/McClanahan/1995}.  Our approach involves representations on Hilbert modules, rather than
Hilbert spaces, allowing us to focus on a single representation at the expense of a bit more abstraction.

\chapter Fell's absorption principle

The classical version of Fell's absorption principle\fn {Apparently Fell has not attempted to generalize this for Fell
bundles.}  \ref {Theorem 2.5.5/BrownOzawa/2008} states that the tensor product of any unitary group representation by
the left-regular representation is weakly-contained in the left-regular representation.  The main goal of this chapter
is to prove a version of this powerful principle for Fell bundles.

In order to set up the context let us first introduce our notation for the left regular representation of the group $G$:
we will denote by $\{e_g\}_{g\in G}$ the canonical basis of $\ell ^2(G)$, and by $\lreg $ the regular representation of
$G$ on $\ell ^2(G)$, so that $$ \lreg _g(e_h) = e_{gh} \for g,h\in G.  $$

Our slightly unusual choice of notation ``$\lreg $'', as opposed to ``$\lambda $'', is due to a potential conflict
between this and our previous notation for the regular representation of a Fell bundle introduced in \cite
{DefineRegularRep}.

Given a Fell bundle $$ \Bun = \{B_g\}_{g\in G} $$ and a representation $\pi =\{\pi _g\}_{g\in G}$ of $\Bun $ on a
Hilbert space $H$, let us consider another representation of $\Bun $, this time on $H\otimes \ell ^2(G)$, by putting $$
(\pi _g\otimes \lreg )(b) = \pi _g(b)\otimes \lreg _g \for b\in B_g.  $$ It is an easy exercise to check that $$ \pi
\otimes \lreg := \{\pi _g\otimes \lreg _g \}_{g\in G} \equationmark TensorProdRep $$ is indeed a *-representation of
$\Bun $.  The integrated form of $\pi \otimes \lreg $, according to \cite {FromRepBunToRelAlg}, is then a
*-representation $$ \varphi : C^*(\Bun ) \to \Lin \big (H\otimes \ell ^2(G)\big ), \equationmark IntegrTensorProdRep $$
satisfying $$ \varphi \big (\hj _g(b)\big ) = \pi _g(b)\otimes \lreg _g \for g\in G \for b\in B_g.  \equationmark
RepFromPiTensorV $$

\state Proposition \label FellAbsorption (Fell's absorption principle for Fell bundles).  Let $\pi $ be a representation
of the Fell bundle $\Bun $ on a Hilbert space $H$.  Also let $\varphi $ be the integrated form of the representation
$\pi \otimes \lreg $ described above.  Then $\varphi $ vanishes on the kernel of the regular representation $\Lambda $,
and hence factors through $\CrB $, providing a representation $\psi $ of the latter,

\null \hfill \beginpicture \setcoordinatesystem units <0.0020truecm, -0.0020truecm> point at 0 0 \put {$C^*(\Bun )$} at
-900 000 \put {$\Lin \big (H\otimes \ell ^2(G)\big )$} at 1200 000 \put {$\CrB $} at 000 900 \arrow <0.11cm> [0.3,1.2]
from -380 000 to 300 000 \put {$\varphi $} at 000 -200 \arrow <0.11cm> [0.3,1.2] from -680 250 to -230 650 \put
{$\Lambda $} at -630 550 \arrow <0.11cm> [0.3,1.2] from 230 650 to 680 250 \put {$\psi $} at 630 550 \endpicture \hfill
\null

\bigskip \noindent such that $\psi \circ \Lambda = \varphi $.  In addition, if $\pi _1$ is faithful, then $\psi $ is
also faithful.

\Proof Consider the linear mapping $F$ from $\Lin \big (H\otimes \ell ^2(G)\big )$ into itself sending each bounded
operator to its diagonal relative to the decomposition $$ H\otimes \ell ^2(G) = \med {\bigoplus }_{g\in G} H\otimes e_g.
$$ To be precise, $$ F(T) = \soma {g\in G} P_gTP_g \for T\in \Lin \big (H\otimes \ell ^2(G)\big ), $$ where each $P_g$
is the orthogonal projection onto $H\otimes e_g$, and the sum is interpreted in the strong operator topology.  It is a
well known fact that $F$ is a \subj {conditional expectation}\null \fn {A \"{conditional expectation} is a linear
mapping $F$ from a C*-algebra $B$ onto a closed *-subalgebra $A$ which is positive, idempotent, contractive and an
$A$-bimodule map.  We moreover say that $F$ is \"{faithful} if $F(b^*b)=0\Rightarrow b=0$.}  onto the algebra of
diagonal operators, which is moreover faithful.  We next claim that $$ F\big (\varphi (y)\big ) = \pi _1\Big (E_1\big
(\Lambda (y)\big )\Big )\otimes 1 \for y\in C^*(\Bun ).  \equationmark CompareTwoCondexp $$ Since both sides represent
continuous linear maps on $C^*(\Bun )$, it is enough to check it for $ y = \hj _g(b), $ where $b\in B_g$.  In this case
we have $$ F\big (\varphi (y)\big ) = F\big (\varphi \big (\hj _g(b)\big )\big ) \={RepFromPiTensorV} F(\pi _g(b)\otimes
\lreg _g\big ) = \delta _{g,1}\pi _1(b)\otimes 1.  $$ In the last step above we have used the fact that $$ \big (\pi
_g(b)\otimes \lreg _g\big ) (H\otimes e_h)\subseteq H\otimes e_{gh} \for h\in G, $$ so that $\pi _g(b)\otimes \lreg _g$
has zero diagonal entries when $g\neq 1$, while it is a diagonal operator when $g=1$.  Focusing now on the
right-hand-side of \cite {CompareTwoCondexp}, we have $$ \pi _1\Big (E_1\big (\Lambda (y)\big )\Big ) = \pi _1\Big
(E_1\big (\Lambda (\hj _g(b))\big )\Big ) \={RegRepProp} \pi _1\Big (E_1\big (\lambda _g(b)\big )\Big ) \={Fourier}
\delta _{g,1}\pi _1(b), $$ hence proving \cite {CompareTwoCondexp}.

Given $y$ in $C^*(\Bun )$ such that $\Lambda (y)=0$, we then have that $ \Lambda (y^*y)=0, $ so \cite
{CompareTwoCondexp} gives $$ 0 = F\big (\varphi (y^*y)\big ) = F\big (\varphi (y)^*\varphi (y)\big ), $$ which implies
that $\varphi (y)=0$, since $F$ is faithful.  This proves that $\varphi $ vanishes on the kernel of $\Lambda $.

Assuming that $\pi _1$ is faithful, we will prove that the kernel of $\Lambda $ in fact coincides with the kernel of
$\varphi $, from where the last sentence of the statement will follow.  We have already seen above that $\Ker (\Lambda
)\subseteq \Ker (\varphi )$, so we need only prove the reverse inclusion.  Assuming that $\varphi (y)=0$, we deduce from
\cite {CompareTwoCondexp} that $ \pi _1\Big (E_1\big (\Lambda (y^*y)\big )\Big )=0, $ and since $\pi _1$ is faithful,
also that $$ 0 = E_1\big (\Lambda (y^*y)\big ) = E_1\big (\Lambda (y)^*\Lambda (y)\big ).  $$ From \cite {EOneIsFaith}
we then get $\Lambda (y)=0$, so $\Ker (\varphi )\subseteq \Ker (\Lambda )$.  \endProof

Recall from \ref {Section 3.3/BrownOzawa/2008} that if $A$ and $B$ are C*-algebras, then the algebraic tensor product
$A\odot B$ may have many different C*-norms, the smaller of which, usually denoted by $\Vert \ponto \Vert \lmin $, is
called the \"{minimal} or \"{spatial tensor norm}, while the biggest one, usually denoted by $\Vert \ponto \Vert \lmax
$, is called the \"{maximal} norm.

Accordingly, we denote by $A\tmin B$ and $A\tmax B$, the completions of $A\odot B$ under $\Vert \ponto \Vert \lmin $ and
$\Vert \ponto \Vert \lmax $, respectively.  These are often referred to as the \"{minimal} and \"{maximal} tensor
products of $A$ and $B$.

If $A$ and $B$ are faithfully represented on Hilbert spaces $H$ and $K$, respectively, then the natural representation
of $A\odot B$ on $H\otimes K$ is isometric for the minimal norm, and this is precisely the reason why this norm is also
called the spatial norm.  See \ref {BrownOzawa/2008} for a comprehensive study of tensor products of C*-algebras.

In our next result we use the minimal tensor product to give a slightly different, but sometimes more useful way to
state Fell's absorption principle:

\state Corollary \label FellAbsorptionThree Let $\Bun = \{B_g\}_{g\in G}$ be a Fell bundle and let $\pi = \{\pi
_g\}_{g\in G}$ be a *-representation of $\Bun $ in a C*-algebra $A$.  Then there is a *-homomor\-phism $$ \psi : \CrB
\to A\tmin C^*\lred (G), $$ such that $$ \psi \big (\lambda _g(b)\big ) = \pi _g(b)\otimes \lreg _g \for g\in G \for
b\in B_g.  $$ Moreover, if $\pi _1$ is faithful, then so is $\psi $.

\Proof Supposing, without loss of generality, that $A$ is an algebra of operators on a Hilbert space $H$, we may view
$\pi $ as a representation of $\Bun $ on $H$.  Letting $\varphi $ and $\psi $ be as in \cite {FellAbsorption}, we have
for all $b$ in any $B_g$, that $$ \psi \big (\lambda _g(b)\big ) = \psi \big (\Lambda (\hj _g(b))\big ) = \varphi \big
(\hj _g(b)\big ) = \pi _g(b)\otimes \lreg _g, $$ as required.  Observing that $\Cr G$ is precisely the subalgebra of
operators on $\ell ^2(G)$ generated by the range of the left-regular representation $\lreg $, we see from the above
computation that $\psi \big (\lambda _g(b)\big )$ lies in the spatial tensor product of $A$ by $\Cr G$, which, by \ref
{Theorem 3.3.11/BrownOzawa/2008}, is isomorphic to $A\tmin C^*\lred (G)$.  We may then view $\psi $ as a map from $\CrB
$ to $A\tmin C^*\lred (G)$, as needed.  Finally, assuming that $\pi _1$ is faithful, \cite {FellAbsorption} implies that
$\psi $ is also faithful.  \endProof

Fell's absorption principle may be used to produce a variety of maps simultaneously involving the full and the reduced
cross-sectional C*-algebras, such as the following:

\state Proposition \label FullMiniAbsorpt Given a Fell bundle $\Bun = \{B_g\}_{g\in G}$, there exists an injective
*-homomorphism $$ \tau :\CrB \to \CB \tmin C^*\lred (G), $$ such that $$ \tau \big (\lambda _g(b)\big ) = \hj
_g(b)\otimes \lreg _g \for g\in G \for b\in B_g.  $$

\Proof Apply \cite {FellAbsorptionThree} to $\pi =\hj $.  \endProof

A similar result is as follows:

\state Proposition \label MiniAbsorp Given a Fell bundle $\Bun = \{B_g\}_{g\in G}$, there exists an injective
*-homomorphism $$ \sigma : \CrB \to \CrB \tmin \Cr G, $$ such that $$ \sigma \big (\lambda _g(b)\big ) = \lambda
_g(b)\otimes \lreg _g \for g\in G \for b\in B_g.  $$

\Proof Apply \cite {FellAbsorptionThree} to $\pi =\lambda $.  \endProof

A very interesting question arises when one replaces the minimal tensor product by the maximal one in \cite
{MiniAbsorp}.  To be precise, given a Fell bundle $\Bun = \{B_g\}_{g\in G}$, for each $g$ in $G$, consider the mapping
$$ \mu _g: B_g \to \CrB \tm \Cr G, $$ given by $$ \mu _g(b) = \lambda _g(b)\otimes \lreg _g \for b\in B_g.  $$

It is easy to see that $\mu = \{\mu _g\}_{g\in G}$ is a *-representation of $\Bun $, the integrated form of which is
then a *-homomorphism $$ \weird : \CB \to \CrB \tm \Cr G, \equationmark WeirdMap $$ satisfying $$ \weird (\hj _g\big
(b)\big ) = \lambda _g(b)\otimes \lreg _g \for b\in B_g.  $$

Observe that the right-hand side of the above expression is identical to the expression defining $\sigma $ in \cite
{MiniAbsorp}, although there we used the minimal tensor product, while here we are using the maximal one.  In any case,
since the above formula for $\weird $ involves not one, but two regular representations, one is left wandering whether
or not it factors through the reduced cross-sectional algebra.  The answer is a bit surprising: we will prove below that
$\weird $ is injective!

\state Proposition \label CrazyFact Given any Fell bundle $\Bun = \{B_g\}_{g\in G}$, there exists an injective
*-homomorphism $$ \weird : \CB \to \CrB \tmax \Cr G, $$ such that $$ \weird (\hj _g\big (b)\big ) = \lambda _g(b)\otimes
\lreg _g \for b\in B_g.  $$

\Proof The proof will of course consist in verifying that the map $\weird $ introduced in \cite {WeirdMap} is injective.
Choose a faithful representation $\rho $ of $C^*(\Bun )$ on a Hilbert space $H$, and apply \cite {FellAbsorptionThree}
to $\rho \circ \hj $, obtaining a faithful representation $$ \psi :\CrB \to \Lin (H) \tmin \Cr G \subseteq \Lin \big
(H\otimes \ell ^2(G)\big ), $$ such that $$ \psi (\lambda _g(b)) = \rho (\hj _g(b))\otimes \lreg _g \for b\in B_g.
\equationmark CharPsio $$

\def \rreg {r^{\scriptscriptstyle G}}

Let $\rreg $ be the \"{right regular representation} of $G$ on $\ell ^2(G)$, given on the standard orthonormal basis
$\{e_g\}_{g\in G}$ of $\ell ^2(G)$ by $$ \rreg _g(e_h) = e_{hg\inv } \for g,h\in G.  $$ We claim that there exist a
*-representation $$ \tilde \rreg : \Cr G \to \Lin \big (H\*\ell ^2(G)\big ), $$ such that $$ \tilde \rreg (\lreg _g) =
1\*\rreg _g \for g\in G.  $$ To see this, consider the unitary operator $U$ on $H\*\ell ^2(G)$ given by $$ U(\xi \otimes
e_g) = \xi \otimes e_{g\inv } \for \xi \in H \for g\in G.  $$ Then, for every $g,h\in G$, and $\xi \in H$, we have $$
U^*(1\otimes \lreg _g)U(\xi \otimes e_h) = U^*(id\otimes \lreg _g)(\xi \otimes e_{h\inv }) = U^*(\xi \otimes e_{gh\inv
}) \$= \xi \otimes e_{hg\inv } = (1\otimes \rreg _g)(\xi \otimes e_h), $$ so $ U^*(1\otimes \lreg _g)U= 1\otimes \rreg
_g, $ for all $g$ in $G$, and the desired representation may then be defined by $$ \tilde \rreg (x) = U^*(1\otimes x)U
\for x\in \Cr G.  $$

It is easy to see that the range of $\psi $ commutes with the range of $\tilde \rreg $, so by the universal property of
the maximal norm \ref {Theorem 3.3.7/BrownOzawa/2008}, we get a representation $$ \psi \times \tilde \rreg : \CrB \tm
\Cr G \to \Lin \big (H\*\ell ^2(G)\big ), $$ such that $$ (\psi \times \tilde \rreg ) (x\otimes y) = \psi (x)\tilde
\rreg (y) \for x\in \CrB \for y\in \Cr G.  $$

Consider the map $\varphi $ defined to be the composition $$ \varphi : C^*(\Bun ) \> \weird \CrB \tm \Cr G \>{\psi
\times \tilde \rreg } \Lin \big (H\*\ell ^2(G)\big ), $$ and observe that, for $b$ in any $B_g$, we have that $$ \varphi
\big (\hj _g(b)\big )= (\psi \times \tilde \rreg )\Big (\weird \big (\hj _g(b)\big )\Big ) = (\psi \times \tilde \rreg
)\big ( \lambda _g(b)\otimes \lreg _g \big ) \={CharPsio} $$$$ = \Big (\rho \big (\hj _g(b)\big )\*\lreg _g\Big )
(1\*\rreg _g) = \rho \big (\hj _g(b)\big )\*\lreg _g\rreg _g.  $$ Notice that the above operator, when applied to a
vector of the form $\xi \*e_1$, with $\xi \in H$, produces $$ \varphi \big (\hj _g(b)\big ) \calcat {\xi \*e_1} = \rho
\big (\hj _g(b)\big )\xi \* e_1.  $$ Since the $\hj _g(b)$ span a dense subset of $\CB $, we conclude that $$ \varphi
\big (y\big ) \calcat {\xi \*e_1} = \rho \big (y\big )\xi \* e_1 \for y\in C^*(\Bun ) \for \xi \in H.  $$

Therefore, assuming that $\weird (y)=0$, for some $y$ in $C^*(\Bun ) $, we also have $\varphi (y)=0$, whence $\rho
(y)\xi =0$, for all $\xi $ in $H$, so $\rho (y)=0$.  Since $\rho $ was supposed to be injective on $C^*(\Bun )$, we
deduce that $y=0$.  This proves that $\weird $ is injective.  \endProof

\nrem Proposition \cite {CrazyFact} first appeared in \ref {Theorem 6.2/AraExelKatsura/2011}.  It is based on an idea
verbally suggested to me by Eberhard Kirchberg at the CRM in Barcelona, in 2011.

\chapter Graded C*-algebras

We have seen in \cite {FinallyInjective/vi\&vii} that a Fell bundle $\Bun $ gives rise to two graded C*-algebras, namely
$$ \CB \and \CrB .  $$ It is not hard to see that the Fell bundles obtained by disassembling these two algebras, as in
\cite {Disassemble}, are both isomorphic to the original Fell bundle $\Bun $.  Since there are situations in which $\CB
$ and $\CrB $ are not isomorphic\fn {An example is the group bundle over a non-amenable group.}, we see that the
correspondence between Fell bundles and graded C*-algebras is not a perfect one.  The goal of this chapter is thus to
study this rather delicate relationship.

Given a graded C*-algebra $B$ with associated Fell bundle $\Bun = \{B_g\}_{g\in G}$, by \cite {FromRepBunToRelAlg} we
see that there is a unique C*-algebra epimorphism $$ \varphi : C^*(\Bun ) \to B, $$ which is the identity on each $B_g$
(identified both with a subspace of $C^*(\Bun )$ and of $B$ in the natural way).  This says that $C^*(\Bun )$ is, in a
sense, the biggest C*-algebra whose associated Fell bundle is $\Bun $.  Our next result will show that the reduced cross
sectional algebra is on the other extreme of the range.  It is also a very economical way to show a C*-algebra to be
graded.

\state Theorem \label Economical Let $B$ be a C*-algebra and assume that for each $g$ in a group $G$, there is
associated a closed linear subspace $B_g \subseteq B$ such that, for all $g$ and $h$ in $G$, one has \zitemno =0 \zitem
$B_g B_h \subseteq B_{gh}$, \zitem $B_g^* = B_{g\inv }$, \zitem $\soma {g\in G} B_g$ is dense in $B$.  \medskip
\noindent Assume, in addition, that there is a bounded linear map $$ F : B \to B_1, $$ such that $F$ is the identity map
on $B_1$, and that $F$ vanishes on each $B_g$, for $g\neq 1$.  Then \iaitem \aitem The $B_g$ are independent and hence
$\{B_g\}_{g\in G}$ is a grading for $B$.  \aitem $F$ is a conditional expectation onto $B_1$.  \aitem If $\Bun $ denotes
the associated Fell bundle, then there exists a C*-algebra epimorphism $$ \psi : B \to \CrB , $$ such that $\psi (b) =
\lambda _g(b)$, for each $g$ in $G$, and each $b$ in $B_g$.

\Proof If $x=\soma {g\in G} b_g$ is a finite sum\fn {By \"{finite sum} we actually mean that the set of indices $g$, for
which $b_g\neq 0$, is finite.}  with $b_g \in B_g$, then $$ x^*x = \soma {g,h\in G } b_g^* b_h = \soma {k\in G } \Big (
\soma {g\in G} b_g^* b_{gk} \Big ).  $$ Observing that $b_g^* b_{gk}$ lies in $B_k$, we see that $$ F(x^*x) = \soma
{g\in G} b_g^* b_g.  $$ Therefore, if $x=0$, then each $b_g = 0$, which shows the independence of the $B_g$'s.  This
also shows that $F$ is positive.

Given $a$ in $B_1$, it is easy to see that $$ F(ax) = aF(x) \and F(xa) = F(x)a \for x\in B, $$ by first checking on
finite sums, as above.  So, apart from contractivity, (b) is proven.

Define a right pre-Hilbert $B_1$-module structure on $B$ via the $B_1$-valued inner product $$ \langle x,y\rangle =
F(x^*y) \for x,y\in B.  $$ For $b,x \in B$, using the positivity of $F$, we have that $$ \langle bx,bx\rangle =
F(x^*b^*bx) \leq \|b\|^2 F(x^*x) = \|b\|^2 \langle x,x\rangle .  $$ So, the left multiplication operators $$ L_b : x \in
B \mapsto bx \in B $$ are bounded and hence extend to the completion $X$ of $B$ (after moding out vectors of norm zero).
It is then easy to show that $$ L : b \in B \mapsto L_b \in \AdjOp (X) $$ is a C*-algebra homomorphism.  Let $$ x=\soma
{g\in G} b_g \and y=\soma {g\in G } c_g $$ be finite sums with $b_g,c_g\in B_g$, and regard both $x$ and $y$ as elements
of $X$.  We then have $$ \langle x,y\rangle = F\Big (\soma {g,h\in G } b_g^* c_h\Big ) = \soma {g\in G} b_g^* c_g = \Big
\langle \soma {g\in G} j_g(b_g), \soma {g\in G} j_g(c_g) \Big \rangle , $$ where the last inner product is that of $\ell
^2(\Bun )$.  This is the key ingredient in showing that the formula $$ U\Big (\soma {g\in G} b_g\Big ) = \soma {g\in G}
j_g(b_g) $$ may be used to define an isometry of Hilbert $B_1$-modules $$ U : X \to \ell ^2(\Bun ).  $$ Here it is
important to remark that the continuity of $F$ ensures that the set of finite sums $\soma {g\in G} b_g$ is not only
dense in $B$ but also in $X$.  For $b$ in $B_g$ and $c$ in $B_h$ we have $$ U L_{b}(c) = U(b c) = j_{gh}(b c) = j_g(b)
\star j_h(c) = \lambda _g(b) j_h(c) = \lambda _g(b) U (c).  $$ Since the finite sums $\soma {h\in G } c_h$ are dense in
$X$, as observed above, we conclude that $$ U L_{b} U^* = \lambda _g(b).  $$ This implies that, for all $b$ in any
$B_g$, the operator $U L_b U^*$ belongs to $\CrB $, and hence the same holds for an arbitrary $b$ in $B$.  This defines
a map $$ \psi : b \in B \mapsto U L_b U^* \in \CrB , $$ which satisfies the requirements in (c).

The only remaining task, the proof of the contractivity of $F$, now follows easily because $F = E_1\circ \psi $, where
$E_1$ is given by \cite {Fourier}.  \endProof

The map $\psi $, above, should be thought of as a generalized regular representation of $B$.

From now on we will be mostly interested in graded algebras possessing a conditional expectation, so we make the
following:

\definition A C*-grading $\{B_g\}_{g\in G}$ on the C*-algebra $B$ is said to be a \subj {topological grading} if there
exists a (necessarily unique) conditional expectation of $B$ onto $B_1$, vanishing on all $B_g$, for $g\neq 1$, as in
\cite {Economical}.

\state Proposition \label CBCrBTopGraded If $\Bun $ is a Fell bundle, then both $\CB $ and $\CrB $ are topologically
graded C*-algebras.

\Proof With respect to $\CrB $, notice that the mapping $E_1:\CrB \to B_1$ given by \cite {Fourier} is a conditional
expectation by \cite {Economical}, clearly satisfying all of the required conditions.

The statement for $\CB $ follows immediately by taking the conditional expectation $$ \CB \trepa \Lambda \longrightarrow
\CrB \>{E_1} B_1.  \omitDoubleDollar \endProof

Not all graded C*-algebras are topologically graded, as the following example shows: suppose that $X$ is a closed subset
of the unit circle $\bf T$, and let $z$ be the standard generator of the algebra $C(X)$, namely the inclusion function
$X\hookrightarrow {\bf C}$.  Letting $$ B_n ={\bf C}z^n \for n\in {\bf Z}, $$ it is easy to see that $B_nB_m\subseteq
B_{n+m}$, and $B_n^*=B_{-n}$.  In order to decide whether or not the $B_n$ are linearly independent subspaces of $C(X)$,
suppose that a finite sum $ \soma {n\in {\bf Z}}b_n = 0, $ where each $b_n\in B_n$. Write $b_n=a_nz^n$, for some scalar
$a_n$, so that $$ \soma {n\in {\bf Z}}b_n = \soma {n\in {\bf Z}}a_nz^n=: p, $$ which is to say that $p$ is a
trigonometric polynomial.  The condition that we need, namely $$ p=0 \quad \Rightarrow \quad (\forall n\in {\bf Z}, \
a_n=0) $$ does not always hold, as $X$ could be contained in the set of roots of $p$, but if we assume that $X$ is
infinite, the above implication is clearly true.  In this case the $B_n$ are linearly independent and so $\{B_n\}_{n\in
{\bf Z}}$ is a grading for $C(X)$.

We claim that when $X$ is a proper infinite subset of $\bf T$, we do not have a topological grading.  To see this,
assume the contrary, which is to say that there exists a \"{continuous} conditional expectation $ F: C(X) \to {\bf C}, $
such that $$ F\Big (\soma {n\in {\bf Z}}a_nz^n\Big ) = a_0, $$ for every trigonometric polynomial $\soma {n\in {\bf
Z}}a_nz^n$.  Equivalently $$ F\big (p\restr X\big ) = {1\over 2\pi }\int _0^{2\pi }\kern -4pt p(e^{it})\,dt,
\equationmark IntegralFormula $$ for every trigonometric polynomial $p$, but since these are dense in $C({\bf T})$, the
above integral formula holds for every $p$ in $C({\bf T})$.

Since we are assuming that $X$ is not the whole of $\bf T$, we may use Tietze's extension Theorem to obtain a nonzero,
positive $p\in C({\bf T})$, vanishing on $X$.  Plugging $p$ into \cite {IntegralFormula} will then bring about a
contradiction, thus proving that our grading does not admit a conditional expectation.

\bigskip

Recalling our discussion about $C^*(\Bun )$ being the biggest graded algebra for a given Fell bundle, we now see that
$\CrB $ is the smallest such, at least among topologically graded algebras.

Summarizing we obtain the following important consequence of \cite {Economical/c} and the remark preceding the statement
of \cite {Economical}:

\state Theorem \label AllGradedInBetween Let $B$ be a topologically graded C*-algebra with grading $\Bun = \{B_g\}_{g\in
G}$.  Then there exists a commutative diagram of surjective *-homomorphism: \vskip -10pt \null \hfill \beginpicture
\setcoordinatesystem units <0.0020truecm, -0.0020truecm> point at 0 0 \put {$C^*(\Bun )$} at -900 000 \put {$\CrB $} at
900 000 \put {$B$} at 000 800 \arrow <0.11cm> [0.3,1.2] from -380 000 to 300 000 \put {$\Lambda $} at 000 -200 \arrow
<0.11cm> [0.3,1.2] from -680 250 to -230 650 \put {$\varphi $} at -630 500 \arrow <0.11cm> [0.3,1.2] from 230 650 to 680
250 \put {$\psi $} at 630 500 \endpicture \hfill \null

\bigskip Another important consequence of \cite {Economical} is the existence of \"{Fourier coefficients} in
topologically graded algebras as follows:

\state Corollary \label Ft Let $B$ be a topologically graded C*-algebra with grading $\{B_g\}_{g\in G}$.  Then, for
every $g$ in $G$, there exists a contractive linear map $$ F_g : B \to B_g, $$ such that, for all finite sums $x=\soma
{g\in G} b_g$, with $b_g \in B_g$, one has $F_g(x) = b_g$.  Moreover, given $h\in G$, and $b\in B_h$, one has that $$
F_g(bx) = bF_{h\inv g}(x) \and F_g(xb) = F_{gh\inv }(x)b, $$ for any $x$ in $B$.

\Proof For the existence it is enough to define $F_g = E_g\circ \psi $, where $E_g$ is as in \cite {Fourier}, and $\psi
$ is given by \cite {Economical/c}.  In order to prove the last two identities, one easily checks them for finite sums
of the form $x=\soma {k\in G} b_k$, with $b_k\in B_k$, and the result then follows by density of the set formed by such
elements.  \endProof

A measure of how much bigger a graded C*-algebra is, relative to the corresponding reduced cross sectional algebra, is
evidently the null space of the regular representation.  The next result gives a characterization of this null space.

\state Proposition \label KerLRR Given a topologically graded C*-algebra $B$ with conditional expectation $F$, let $\psi
$ be its regular representation, as in \cite {Economical/c}. Then $$ \hbox {Ker}(\psi ) = \{b\in B : F(b^* b) =0 \}.  $$

\Proof Observing that $F = E_1\circ \psi $, we have $$ F(b^* b) = E(\psi (b)^* \psi (b)), $$ from where we see that
$F(b^* b) = 0$, if and only if $\psi (b) = 0$, by \cite {EOneIsFaith}.  \endProof

This can be employed to give a useful characterization of $\CrB $ among graded algebras:

\state Proposition \label Charact Suppose we are given a topologically graded C*-alge\-bra $B$ with grading $\Bun =
\{B_g\}_{g\in G}$, and conditional expectation $F$.  If $F$ is \subjex {faithful}{faithful conditional expectation}, in
the sense that $$ F(b^*b)=0 \ \Rightarrow \ b=0 \for b\in B, $$ then $B$ is canonically isomorphic to $\CrB $.

\Proof The regular representation $\psi $ of \cite {Economical/c} will be an isomorphism by \cite {KerLRR}.  \endProof

\nrem Theorems \cite {Economical} and \cite {AllGradedInBetween} were first proved in \ref {Exel/1997b}.

\chapter Amenability for Fell bundles

Theorem \cite {AllGradedInBetween} tells us that the topologically graded C*-algebras whose associated Fell bundles
coincide with a given Fell bundle $\Bun $ are to be found among the quotients of $\CB $ by ideals contained in the
kernel of the regular representation $\Lambda $.

It is therefore crucial to understand the kernel of $\Lambda $ and, in particular, to determine conditions under which
$\Lambda $ is injective.  In the case of the group bundle over $G$, it is well known that the injectivity of $\Lambda $
is equivalent to the amenability of $G$.

In this chapter we will extend to Fell bundles some of the techniques pertaining the the theory of amenable groups,
including the important notion of amenability of group actions introduced by Anantharaman-Delaroche \ref
{AnantharamanDelaroche/1987}.  We begin with some terminology.

\state {Definition} \label DefineAmenaBun A Fell bundle $\Bun $ is said to be \subjex {amenable} {amenable Fell bundle}
if the regular representation $$ \Lambda : \CB \arw \CrB $$ is injective.  \stateend

Amenability, when applied to the context of groups, may be characterized by a large number of equivalent conditions \ref
{Greenleaf/1969}, such as the equality of the full and reduced group C*-algebras \ref {Theorem 7.3.9/Pedersen/1979}.
Many of these conditions may be rephrased for Fell bundles although they are not always known to be equivalent.  One
therefore has many possibly inequivalent choices when defining the term \"{amenable} in the context of Fell bundles.
Definition \cite {DefineAmenaBun} is just one among many alternatives.

An immediate consequence of \cite {AllGradedInBetween} is:

\stateRef {Proposition}{Democracy} Let $\Bun $ be an amenable Fell bundle.  Then all topologically graded C*-algebras
whose associated Fell bundles coincide with $\Bun $ are isomorphic to each other.  \stateend

The next technical result is intended to pave the way for the introduction of sufficient conditions for the amenability
of a given Fell bundle.  The idea is to produce \"{wrong way maps} which will later be used to prove the injectivity of
$\Lambda $.

\state Lemma \label WrongWay Let $\Bun = \{B_g\}_{g\in G}$ be a Fell bundle.  Given a finitely supported function $ a :
G \to B_1, $ the formula $$ V(z) = \soma {g,h\in \G } \hj _g\big (a(gh)^* E_g(z) a(h)\big ) \for z\in \CrB ,
\equationmark FormulaDefiningVInStatement $$ gives a well defined completely positive linear map from $\CrB $ to $\CB $,
such that $$ \| V \| \leq \normsum {\soma {g\in \G } a(g)^* a(g)}.  $$ Identifying each $B_g$ as a subspace of $\CrB $
or $\CB $, as appropriate, one moreover has $$ V (b) = \soma {h\in \G } a(gh)^* b a(h) \for b\in B_g.  $$

\Proof Let $$ \rho :\CB \to \Lin (H) $$ be any faithful *-representation of $\CB $ on a Hilbert space $H$.  For each $g$
in $G$, let $\pi _g = \rho \circ \hj _g$, so that $ \pi =\{\pi _g\}_{g\in G} $ is a *-representation of $\Bun $ on $H$.
Consider the representation $\pi \otimes \lreg $ of $\Bun $ on $H\otimes \ell ^2(G)$ described in \cite {TensorProdRep},
and let $$ \varphi : C^*(\Bun ) \to \Lin \big (H\otimes \ell ^2(G)\big ) $$ be the integrated form of $\pi \otimes \lreg
$.  By \cite {FellAbsorption}, we have that $\varphi $ vanishes on $\Ker (\Lambda )$, so it factors through $\CrB $,
giving a *-homomorphism $\psi $ such that the diagram \vskip -10pt \beginpicture \setcoordinatesystem units
<0.0020truecm, -0.0020truecm> point at 0 0 \setplotarea x from -4000 to 4000, y from -500 to 2900 \put {$C^*(\Bun )$} at
-900 000 \put {$\Lin \big (H\otimes \ell ^2(G)\big )$} at 1150 000 \put {$\CrB $} at 000 900 \arrow <0.11cm> [0.3,1.2]
from -380 000 to 300 000 \put {$\varphi $} at 000 -200 \arrow <0.11cm> [0.3,1.2] from -680 250 to -230 650 \put
{$\Lambda $} at -630 500 \arrow <0.11cm> [0.3,1.2] from 230 650 to 680 250 \put {$\psi $} at 630 540 \endpicture

\bigskip \noindent commutes.  We next consider the operator \def \A {A} $$ \A : H \arw H \* \ell ^2(\G ) $$ given by $$
\A (\xi ) = \soma {g\in \G } \alpha (g) \xi \* e_g \for \xi \in H, $$ where $\alpha (g) = \pi _1\big (a(g)\big )$.  In
order to compute the norm of $\A $, let $\xi \in H$.  Then $$ \Vert \A (\xi )\Vert ^2 = \soma {g\in \G } \| \alpha (g)
\xi \|^2 = \soma {g\in \G } \langle \alpha (g)^* \alpha (g) \xi , \xi \rangle = \Big \langle \soma {g\in \G } \alpha
(g)^* \alpha (g) \xi , \xi \Big \rangle \$\leq \normsum {\soma {g\in \G } \alpha (g)^* \alpha (g)} \| \xi \|^2, $$ from
where we deduce that $$ \| \A \| \leq \normsum {\soma {g\in \G } \alpha (g)^* \alpha (g)} ^{1/2} \leq \normsum {\soma
{g\in \G } a(g)^* a(g)} ^{1/2}.  \equationmark NormOfA $$

Considering the completely positive map $$ W : T \in \Lin \big (H\* \ell ^2(\G )\big ) \longmapsto \A ^* T \A \in \Lin
(H), $$ observe that for $b$ in any $B_g$, and every $\xi $ in $H$, one has $$ W\big (\pi _g(b) \* \lreg _g\big ) \xi =
\A ^* \big (\pi _g(b) \* \lreg _g\big ) \Big ( \soma {h\in \G } \alpha (h) \xi \* e_h \Big ) \$= \A ^* \Big ( \soma
{h\in \G } \pi _g(b) \alpha (h) \xi \* e_{gh} \Big ) = \soma {h\in \G } \alpha (gh)^* \pi _g(b) \alpha (h) \xi , $$
where we have used the fact that $ \A ^*(\eta \* e_k) = \alpha (k)^* \eta , $ which the reader will have no difficulty
in proving.  Summarizing, we have that $$ W\big (\pi _g(b) \* \lreg _g\big ) = \soma {h\in \G } \alpha (gh)^* \pi _g(b)
\alpha (h).  $$

If we then define $ V $ to be the composition $$ V\ :\ \CrB \> \psi \Lin \big (H\otimes \ell ^2(G)\big ) \> W \Lin (H),
$$ we have for $b$ in any $B_g$ that $$ V\big (\lambda _g(b)\big ) = W\big (\psi (\lambda _g(b))\big ) = W\big (\varphi
(\hj _g(b))\big ) \={FromRepBunToRelAlg} W\big (\pi _g(b) \* \lreg _g\big ) \$= \soma {h\in \G } \alpha (gh)^* \pi _g(b)
\alpha (h) = \rho \Big (\soma {h\in \G } \hj _g\big (a(gh)^* b a(h)\big )\Big ).  $$

Recalling that $\CrB $ is generated by the $\lambda _g(B_g)$, we see that the range of $V$ is contained in the range of
$\rho $, which is isomorphic to $\CB $, since we have chosen $\rho $ to be faithful.  Identifying $\CB $ with its image
under $\rho $, we may then view $V$ as a map $$ V : \CrB \arw \CB , $$ such that $$ V\big (\lambda _g(b)\big ) = \soma
{h\in \G } \hj _g\big (a(gh)^* b a(h)\big ) \for b\in B_g.  \equationmark MainIdForV $$

Observe that, with the appropriate identifications $$ B_g \simeq \lambda _g(B_g) \subseteq \CrB \and B_g \simeq \hj
_g(B_g) \subseteq \CB , $$ we have therefore proven the last identity in the statement.  For a finite sum $$ z =\soma
{g\in G} \lambda _g(b_g), $$ with each $b_g$ in $B_g$, we then have by \cite {MainIdForV} that $$ V(z) = \soma {g\in G}
V\big (\lambda _g(b_g)\big ) = \soma {g,h\in \G } \hj _g\big (a(gh)^* b_g a(h)\big ) \$= \soma {g,h\in \G } \hj _g\big
(a(gh)^* E_g(z) a(h)\big ), $$ proving \cite {FormulaDefiningVInStatement} for all $z$ of the above form.  Observe that
if $(g,h)$ is a pair of group elements for which the corresponding term in the last sum above is nonzero, then both $h$
and $gh$ must lie in the finite support of $a$, which we henceforth denoted by $A$.  For each such pair we then have $$
g = (gh)h\inv \in AA\inv , $$ so $$ (g,h) \in AA\inv \times A.  $$ Consequently the sum in \cite
{FormulaDefiningVInStatement} is finite, hence representing a continuous function of $z$.  Since the set of $z$'s
considered above is dense, \cite {FormulaDefiningVInStatement} is proven.

In order to estimate the norm of $V$ we have $$ \| V \| \leq \| W \| \leq \Vert A\Vert ^2 \_\leq {NormOfA} \normsum
{\soma {g\in \G } a(g)^* a(g)}, $$ as desired.  \endProof

This puts us in position to describe an important concept.

\definition \label AP Let $\Bun = \{B_g\}_{g\in G}$ be a Fell bundle.  By a \subj {\aproxNet } for $\Bun $ we mean a net
$\{ a_i \}_{i\in I}$ of finitely supported functions $$ a_i: \G \arw B_1, $$ such that \izitem \zitem $ \ds \sup _{i\in
I} \normsum {\soma {g\in \G } a_i(g)^* a_i(g)} < \infty , $ \zitem $ \ds \lim _{i \rightarrow \infty } \soma {h\in \G }
a_i(gh)^* b a_i(h)=b, $ \ for all $b$ in each $B_g$.  \medskip \noindent If a {\aproxNet } exists, we will say that
$\Bun $ has the \subj {approximation property}.  \stateend

Sometimes one might have difficulty in checking \cite {AP/ii} for every single $b$ in $B_g$, so the following might come
in handy:

\state Proposition \label LazyApproxProp Let $\Bun = \{B_g\}_{g\in G}$ be a Fell bundle and let $\{ a_i \}_{i\in I}$ be
a {\aproxNet } for $\Bun $, except that \cite {AP/ii} is only known to hold for $b$ in a total\fn {A \"{total set} in a
normed space is a set spanning a dense subspace.}  subset of each $B_g$.  Then \cite {AP/ii} does hold for all $b$ in
$B_g$.  Consequently $\{ a_i \}_{i\in I}$ is a {\aproxNet } and $\Bun $ has the approximation property.

\Proof For each $i$ in $I$, and each $g$ in $G$, consider the map $$ T_i \ : \ b\in B_g \ \mapsto \ \soma {h\in \G }
a_i(gh)^* b a_i(h) \ \in \ B_g.  $$ By hypothesis we have that the $T_i$ converge to the identity map of $B_g$,
pointwise on a dense set.  To get the conclusion it is then enough to establish that the $T_i$ are uniformly bounded.

Letting $V_i$ be the map provided by \cite {WrongWay} in terms of $a_i$, notice that $T$ may be seen as a restriction of
$V_i$, once $B_g$ is identified with $\lambda _g(B_g)\subseteq \CrB $, and with $\hj _g(B_g)\subseteq \CB $.  Since
these are isometric copies of $B_g$ by \cite {FinallyInjective/iv--v}, we deduce that $\Vert T_i\Vert \leq \Vert
V_i\Vert $, from where the desired uniform boundedness follows, hence concluding the proof.  \endProof

The relationship between the approximation property and the amenability of Fell bundles is the main result of this
chapter:

\stateRef {Theorem}{Main} If a Fell bundle $\Bun $ has the approximation property, then it is amenable.  \stateend

\Proof Given a {\aproxNet } $\{a_i\}_{i\in I}$ as in \cite {AP}, let us consider, for each $i$ in $I$, the map $$ V _i :
\CrB \arw \CB $$ provided by \cite {WrongWay}, relative to $a_i$.  Define a map $\Phi _i$ from $\CB $ to itself to be
the composition $\Phi _i = V _i \circ \Lambda $, and observe that, by hypothesis we have $$ \lim _{i \rightarrow \infty
} \Phi _i(b) = b, $$ for every $b$ in each $B_g$.  Because the $b$'s span a dense subspace of $\CB $, and because the
$\Phi _i$'s are uniformly bounded, we conclude that $\lim _{i \rightarrow \infty } \Phi _i(y) = y$, for every $y$ in
$\CB $.

Now, if $y$ is in the kernel of the regular representation, that is, if $y$ is in $\CB $ and $\Lambda (y)= 0$, then $$ y
= \lim _{i \rightarrow \infty } \Phi _i(y) = \lim _{i \rightarrow \infty } V _i(\Lambda (y)) = 0, $$ which proves that
$\Lambda $ is injective.  \endProof

A somewhat trivial example of this situation is that of bundles over amenable groups.

\state Theorem \label BundleAmenaGroup Let $\Bun $ be a Fell bundle over the amenable group $\G $. Then $\Bun $
satisfies the approximation property and hence is amenable.

\Proof According to \ref {7.3.8/Pedersen/1979}, there exists a net $\{f_i\}_{i\in I}$ in the unit sphere of $\ell ^2(G)$
such that $$ \lim _{i \rightarrow \infty } \soma {h\in \G } \overline { f_i(gh) } f_i(h) = 1 \for g\in \G .  $$ By
truncating the $f_i$ to larger and larger subsets of $G$ we may assume them to be finitely supported.  Letting
$\{u_j\}_{j\in J}$ be an approximate identity for $B_1$, one checks that the doubly indexed net $$ \{a_{i,j}
\}_{(i,j)\in I\times J} $$ defined by $a_{i,j}(g) = f_i(g) u_j$ is a {\aproxNet } for $\Bun $.  \endProof

There is another proof that $\Bun $ is amenable when $G$ is amenable based on Fell's absorption principle, as follows:
if $G$ is amenable then the trivial representation is weakly contained in the left-regular representation, which means
that there is a *-homomorphism $$ \varepsilon :\Cr G \to {\bf C} $$ such that $\varepsilon (\lreg _ g) = 1$, for all $g$
in $G$.  Taking into account that $C^*\lred (G)$ is nuclear, and using the map $\tau $ given by \cite {FullMiniAbsorpt},
one has that the map $\varphi $ defined by the composition $$ \varphi \ : \ \CrB \> \tau \CB \tmin C^*\lred (G) \simeq
\CB \tmax C^*\lred (G) \> {id\otimes \varepsilon } \CB $$ satisfies $$ \varphi \big (\lambda _g(b)\big ) = \hj _g(b)
\for g\in G \for b\in B_g, $$ whence $\varphi \circ \Lambda $ is the identity of $\CB $, and hence $\Lambda $ is
injective.

\medskip One of the major questions we are usually confronted with in Classical Harmonic Analysis is whether or not a
function on the unit circle may be reconstructed from its Fourier coefficients.  The same question may be stated for any
topologically graded algebra:

\state Question \label DoesFourierConverge \rm Let $B$ be a topologically graded C*-algebra and let $F_g$ be the
\"{Fourier coefficient operators} given by \cite {Ft}.  Given an element $x$ in $B$, can we reconstruct $x$ if we are
given the value of $F_g(x)$ for every $g$ in $G$?

There are two possible interpretations of this question, depending on what one means by \"{to reconstruct}.  On the one
hand, \cite {EOneIsFaith} implies that if $x$ and $y$ are two elements of $\CrB $ having the same Fourier coefficients,
then $x=y$.  This means that there is at most one $x$ in $\CrB $ with a given \"{Fourier transform}, but it gives us no
algorithm to actually produce such an $x$.  Of course, if $x$ is a finite sum as in \cite {Ft}, then we see that $$ x =
\soma {g\in \G } F_g(x), $$ the sum having only finitely many nonzero terms.  Any attempt at generalizing this beyond
the easy case of finite sums will necessarily require a careful analysis of the convergence of the above series.  The
fact that, over the group bundle ${\bf C}\times {\bf Z}$, this reduces to the usual question of convergence of Fourier
series, is a stark warning that question \cite {DoesFourierConverge} does not have a straightforward answer.

In the classical case of functions on the unit circle, the lack of convergence of Fourier series is circumvented by
Cesaro sums, an analogue of which we shall now discuss.

\stateRef {Definition}{SumProc} Let $B$ be a topologically graded C*-algebra with grading $\{B_g\}_{g\in G}$, and let
$F_g$ be the Fourier coefficient operators given by \cite {Ft}.  A completely positive bounded linear map $S :B \arw B$,
is said to be a \subj {summation process} if, \zitemno = 0 \zitem $S\circ F_g=0$, for all but finitely many $g\in G$,
\zitem $S(x) = \soma {g\in \G } S(F_g(x))$, for all $x$ in $B$.  \stateend

For example, choosing a finite subset $K\subseteq G$, one obtains a summation process by defining $$ S(x) = \soma {g\in
K} F_g(x) \for x\in B.  $$

The following result answers question \cite {DoesFourierConverge} in the affirmative for Fell bundles possessing the
approximation property:

\state Proposition \label Cesaro Suppose that we are given a topologically graded C*-algebra $B$ whose grading $\Bun =
\{B_g\}_{g\in G}$ admits a {\aproxNet } $\{ a_i \}_{i\in I}$ (and hence satisfies the approximation property).  Then the
maps $S_i:B\to B$, given by $$ S_i(x) = \soma {g,h\in \G } a_i(gh)^* F_g(x) a_i(h) \for x\in B, $$ form a bounded net
$\{S_i\}_{i\in I}$ of summation processes converging pointwise to the identity map of $B$.

\Proof Since $\Bun $ is amenable we may invoke \cite {Democracy} to assume without loss of generality that $B=\CrB $.
By \cite {WrongWay} we have that the $S_i$ indeed form a well defined uniformly bounded net of completely positive
linear maps.  Moreover, as in the proof of \cite {Main}, one checks that the $S_i$ converge pointwise to the identity
map.

We will now prove that the $S_i$ are summation processes.  With respect to \cite {SumProc/i}, suppose that $S_i\circ
F_g\neq 0$.  Then there is some $b\in B_g$, such that $$ 0 \neq S_i(b) = \soma {h\in \G } a_i(gh)^* b a_i(h).  $$ In
particular there must be at least one $h$ for which both $a_i(gh)$ and $a_i(h)$ are nonzero.  If $A$ is the support of
$a_i$, then this implies that $gh,h\in A$, whence $$ g = (gh)h\inv \in AA\inv .  $$ Thus $S_i\circ F_g=0$, for all $g$
not belonging to the finite set $AA\inv $, proving \cite {SumProc/i}.

If $x$ is a finite sum $ x=\soma {g\in \G } b_g, $ with each $b_g$ in $B_g$, then clearly $$ x=\soma {g\in \G } F_g(x).
$$ This implies \cite {SumProc/ii} for such an $x$ and hence also for all $x$ in $B$, by continuity of both sides of
\cite {SumProc/ii}, now that we know that the sum involved is finite.  \endProof

Given a C*-algebraic partial dynamical system $$ \Th = \big (A,\ G,\ \{\D g\}_{g\in G},\ \{\th g\}_{g\in G}\big ) $$ we
may consider the associated semi-direct product bundle $\Bun $, and we could try to determine conditions on $\Th $ which
would say that $\Bun $ satisfies the approximation property.  Unraveling the definition of the approximation property
for the special case of semi-direct product bundles, the condition is the existence of a net $\{ a_i \}_{i\in I}$ of
finitely supported functions $$ a_i: \G \arw A, $$ which is bounded in the sense that $$ \ds \sup _{i\in I} \normsum
{\soma {g\in \G } a_i(g)^* a_i(g)} < \infty , $$ and such that $$ \lim _{i \rightarrow \infty } \soma {h\in \G }
a_i(gh)^* \th g\big (c a_i(h)\big ) = \th g(c) \for c\in \Di g.  $$ One arrives at this by simply plugging $b = \th
g(c)\delta _g$ in \cite {AP/ii}, where $c$ is in $\Di g$.

\bigskip In \ref {Zimmer/1978} Zimmer introduced a notion of amenability of global group actions on measure spaces (see
also \ref {Zimmer/1977}), which was later adapted by Anan\-thara\-man-Delaroche \ref {AnantharamanDelaroche/1987} to
actions on C*-algebras.  A slightly different version of this concept, taken from \ref {Definition
4.3.1/BrownOzawa/2008}, is as follows:

\definition \label ClaireAmena Let $\beta $ be a global action of a discrete group $G$ on a unital C*-algebra $A$.  Then
$\beta $ is said to be an \subj {amenable action} if there exists a sequence $\{T_i\}_{i\in {\bf N}}$ of finitely
supported functions $$ T_i:G\to {\cal Z}(A), $$ where ${\cal Z}(A)$ stands for the center of $A$, such that \izitem
\zitem $T_i(g)\geq 0$, for all $i\in {\bf N}$, and all $g\in G$, \zitem $\soma {g\in G} T_i(g)^2 = 1$, \zitem $\lim
_{i\to \infty }\Vert T_i-T^g_i\Vert _2 =0$, for all $g$ in $G$, where $$ T^g_i(h)=\beta _g\big (T_i(g\inv h)\big ) \for
h\in G, $$ and $\Vert T\Vert _2$ is defined for all finitely supported functions $T:G\to A$ by $$ \Vert T\Vert _2 =
\normsum {\soma {g\in G} T(g)^*T(g)}^{^{1\over 2}}.  $$

We would now like to relate the above notion with the approximation property for Fell bundles.

\state Proposition \label AmenaActApprox Let $\beta $ be a global amenable action of a group $G$ on a unital C*-algebra
$A$.  Then the corresponding semi-direct product bundle has the approximation property.

\Proof Let $\Bun = \{B_g\}_{g\in G}$ be the semi-direct product bundle for $\beta $, so that each $B_g = A\delta _g$.
As usual we will identify $B_1$ with $A$ in the obvious way.

Choosing $\{T_i\}_{i\in {\bf N}}$ as in \cite {ClaireAmena}, let us consider each $T_i$ as a $B_1$-valued function, and
let us prove that $\{T_i\}_{i\in {\bf N}}$ satisfies the conditions of \cite {AP}.

Observing that \cite {AP/i} follows immediately from \cite {ClaireAmena/ii}, we concentrate our efforts on proving \cite
{AP/ii}.  For this, choose any element $b$ in any $B_g$, so that necessarily $b = c\delta _g$, for some $c$ in $A$.
Then $$ \soma {h\in \G } T_i(gh)^* b T_i(h) = \soma {h\in \G } T_i(gh) c\delta _g T_i(h) = \soma {h\in \G } T_i(gh)
c\beta _g\big (T_i(h)\big ) \delta _g \$= \soma {k\in \G } T_i(k) \beta _g\big (T_i(g\inv k)\big ) = \soma {k\in \G }
T_i(k) T^g_i(k) \ b \convrg i b, $$ where the last step follows from \ref {Lemma 4.3.2/BrownOzawa/2008}.  This verifies
\cite {AP/ii}, so the proof is concluded.  \endProof

\medskip We would now like to point out a couple of results from \ref {Exel/2000}, whose proofs are perhaps a bit too
long to be included here, so we shall restrict ourselves to giving their statements.  The following is taken from \ref
{Theorems 4.1 \& 6.3/Exel/2000}:

\state Theorem \label ApproxBundle \def \Gen {\hbox {\rs G}}\null Let $\Gen $ be any set and let $\F $ be the free group
on $\Gen $, equipped with the usual length function.  \izitem \zitem Given a semi-saturated *-partial representation
$\Pr $ of\/ $\F $ in a C*-algebra $A$, which is also \subjex {orthogonal}{orthogonal partial representation}, in the
sense that $$ \pr g^*\pr h=0, $$ for any two distinct $g,h\in \Gen $, then the associated Fell bundle $\Bun ^\Pr $,
introduced in \cite {BunFromPrep}, satisfies the approximation property and hence is amenable.  \zitem Given a
semi-saturated Fell bundle $\Bun = \{B_g\}_{g\in G}$ over $\F $, which is also \subjex {orthogonal}{orthogonal Fell
bundle}, in the sense that $$ B_g^*B_h=\{0\}, $$ for any two distinct $g,h\in \Gen $, then $\Bun $ is amenable.

\nrem The definition of the approximation property given in \cite {AP} is inspired by the equivalent conditions of \ref
{Th\'eor\`eme 3.3/AnantharamanDelaroche/1987}.  The question of whether the converse of Theorem \cite {Main} holds is a
most delicate one.  In the special case when $\Bun $ is the semi-direct product bundle relative to a global action of an
exact group on an abelian C*-algebra, a converse of Theorem \cite {Main} has been proven in \ref {Matsumura/2012}.
Another open question of interest is whether \cite {AmenaActApprox} admits a converse.

\chapter Functoriality for Fell bundles

\def \quo {\Bun /\Jbun }

\def \smallp {{\scriptscriptstyle P}} \def \lpb {\ell ^2_P(\Bun )} \def \lp {\lambda ^\smallp } \def \Lp {\Lambda
^{\scriptscriptstyle \!P}} \def \jp {j^\smallp } \def \V {V} \def \funcfull {} \def \funcred #1{#1\lred }

So far we have been studying individual Fell bundles, but we may also see them as objects of a category.  To do so we
will introduce a notion of \"{morphism} between Fell bundles as well as the notion of a \"{\subFellBun dle}.

One of the most disturbing questions we will face is whether or not the cross sectional algebra of a \subFellBun dle is
a subalgebra of the cross sectional algebra of the ambient Fell bundle.  Strange as it might seem, the natural map
between these algebras is not always injective!  Attempting to dodge this anomaly we will see that under special
hypothesis (heredity or the existence of a conditional expectation) sub-bundles lead to bona fide subalgebras.

\definition \label DefineMorphBun Given Fell bundles $\Abun = \{A_g\}_{g\in G}$ and $\Bun = \{B_g\}_{g\in G}$, a \subjex
{morphism}{morphism between Fell bundles} from $\Abun $ to $\Bun $ is a collection $ \varphi = \{\varphi _g\}_{g\in G} $
of linear maps $$ \varphi _g: A_g \to B_g, $$ such that \izitem \zitem $\varphi _g(a)\varphi _h(b) = \varphi _{gh}(ab)$,
and \zitem $\varphi _g(a)^* = \varphi _{g\inv }(a^*)$, \medskip \noindent for all $g$ and $h$ in $G$, and for all $a$ in
$A_g$ and $b$ in $A_h$.

Observe that if $\varphi $ is as above then $\varphi _1$ is evidently a *-homomorphism from $A_1$ to $B_1$, hence
$\varphi _1$ is necessarily continuous.  Also, given $a$ in any $A_g$, we have that $$ \Vert \varphi _g(a)\Vert ^2 =
\Vert \varphi _g(a)^*\varphi _g(a)\Vert = \Vert \varphi _{g\inv }(a^*)\varphi _g(a)\Vert \$= \Vert \varphi _1(a^*a)\Vert
\leq \Vert a^*a\Vert = \Vert a\Vert ^2, $$ so indeed all of the $\varphi _g$ are necessarily continuous.

If every $\varphi _g$ is bijective, it is easy to see that $$ \varphi \inv := \{\varphi _g\inv \}_{g\in G} $$ is also a
morphism.  In this case we say that $\varphi $ is an \"{isomorphism}.  If an isomorphism exists between the Fell bundles
$\Abun $ and $\Bun $, we say that $\Abun $ and $\Bun $ are \subj {isomorphic Fell bundles}.

Morphisms between Fell bundles lead to *-homomorphisms between the corresponding cross sectional C*-algebras as we will
now show.

\state Proposition \label FunctorialityFull Let $\varphi = \{\varphi _g\}_{g\in G}$ be a morphism from the Fell bundle
$\Abun = \{A_g\}_{g\in G}$ to the Fell bundle $\Bun = \{B_g\}_{g\in G}$.  Then there exists a unique *-homomorphism $$
\funcfull \varphi : C^*(\Abun ) \to \CB , $$ (denoted by $\funcfull \varphi $ by abuse of language) such that $$
\funcfull \varphi \big (\hj ^A_g(a)\big ) = \hj ^B_g\big (\varphi _g(a)\big ), $$ for all $a$ in any $A_g$, where we
denote by $\hj ^A$ and $\hj ^B$ the universal representations of $\Abun $ and $\Bun $, respectively.

\Proof For each $g$ in $G$, consider the composition $$ \pi _g: A_g \> {\varphi _g} B_g \> {\hj ^B_g} \CB .  $$ It is
then routine to check that $\pi = \{\pi _g\}_{g\in G}$ is a *-representation of $\Abun $ in $\CB $, so the result
follows from \cite {FromRepBunToRelAlg}.  \endProof

As above, we will now show that morphisms between Fell bundles lead to *-homomorphisms between the corresponding reduced
cross sectional C*-algebras:

\state Proposition \label FunctorialityRed Let $\varphi = \{\varphi _g\}_{g\in G}$ be a morphism from the Fell bundle
$\Abun = \{A_g\}_{g\in G}$ to the Fell bundle $\Bun = \{B_g\}_{g\in G}$.  Then there exists a unique *-homomorphism $$
\funcred \varphi : \Cr \Abun \to \CrB , $$ such that $$ \funcred \varphi \big (\lambda ^A_g(a)\big ) = \lambda ^B_g\big
(\varphi _g(a)\big ), $$ for all $a$ in any $A_g$, where we denote by $\lambda ^A$ and $\lambda ^B$ the regular
representations of $\Abun $ and $\Bun $, respectively.  Moreover, if $\varphi _1$ is injective, then so is $\funcred
\varphi $.

\Proof \def \aa {'}\def \bb {} Consider the diagram

\null \hfill \beginpicture \setcoordinatesystem units <0.0020truecm, -0.0025truecm> point at 1000 1000 \put {$C^*(\Abun
)$} at -1000 0000 \arrow <0.11cm> [0.3,1.2] from -1000 250 to -1000 600 \put {$\Lambda \aa $} at -1250 400 \put {$\Cr
{\Abun }$} at -1000 800 \arrow <0.11cm> [0.3,1.2] from -1000 1050 to -1000 1400 \put {$E\aa _1$} at -1250 1200 \put
{$A_1$} at -1000 1650

\put {$C^*(\Bun )$} at 1000 0000 \arrow <0.11cm> [0.3,1.2] from 1000 250 to 1000 600 \put {$\Lambda \bb $} at 1250 400
\put {$\Cr {\Bun }$} at 1000 800 \arrow <0.11cm> [0.3,1.2] from 1000 1050 to 1000 1400 \put {$E\bb _1$} at 1250 1200
\put {$B_1$} at 1000 1650 \arrow <0.11cm> [0.3,1.2] from -550 000 to 550 000 \put {$\funcfull \varphi $} at 0 -180
\arrow <0.11cm> [0.3,1.2] from -700 1650 to 700 1650 \put {$\varphi _1$} at 000 1500 \arrow <0.11cm> [0.3,1.2] from -600
180 to 500 700 \put {$\psi $} at -100 600 \endpicture \hfill \null

\bigskip \noindent where $\funcfull \varphi $ is given by \cite {FunctorialityFull}, $\Lambda \aa $ and $\Lambda \bb $
are the regular representations of $\Abun $ and $\Bun $, respectively, and $E\aa _1$ and $E\bb _1$ are the conditional
expectations given by \cite {Fourier} for $\Abun $ and $\Bun $, respectively.  By checking on $C_c(\Abun )$, it is then
easy to see that the outside of the diagram is commutative.

Defining $\psi =\Lambda \bb \circ \funcfull \varphi $, we then claim that the null space of $\psi $ contains the null
space of $\Lambda \aa $.  To see this, notice that for any $y\in C^*(\Abun )$, we have that $$ \Lambda \aa (a) = 0 \
\_\Longleftrightarrow {EOneIsFaith} \ E_1\aa \big (\Lambda \aa (a^*a)\big ) = 0 \ \Longrightarrow \ \varphi _1\Big
(E_1\aa \big (\Lambda \aa (a^*a)\big )\Big ) = 0 \ \Longleftrightarrow $$$$ \Longleftrightarrow \ E_1\bb \big (\psi
(a^*a)\big ) = 0 \ \_\Longleftrightarrow {EOneIsFaith} \ \psi (a) = 0.  $$

This shows our claim, so we see that $\psi $ factors through $\Cr {\Abun }$, producing the map $\funcred \varphi $
mentioned in the statement and one now easily checks that it satisfies all of the required conditions.

In case $\varphi _1$ is injective, the only non-reversible arrow ``$\Longrightarrow $'' above may be replaced by a
reversible one, from where it follows that the null space of $\psi $ in fact coincides with the null space of $\Lambda
\aa $, which is to say that the factored map, namely $\funcred \varphi $, is injective.  \endProof

\medskip Morphisms between Fell bundles have properties somewhat similar to morphisms between C*-algebras, as we shall
now discuss.

\state Proposition \label MorpLikeCStar Let $\varphi = \{\varphi _g\}_{g\in G}$ be a morphism from the Fell bundle
$\Abun = \{A_g\}_{g\in G}$ to the Fell bundle $\Bun = \{B_g\}_{g\in G}$.  \iaitem \aitem If $\varphi _1$ is injective,
then all of the $\varphi _g$ are isometric.  \aitem If the range of each $\varphi _g$ is dense in $B_g$, then all of the
$\varphi _g$ are in fact surjective.

\Proof Since $A_1$ and $B_1$ are C*-algebras, under the hypothesis in (a), we have that $\varphi _1$ is in fact
isometric.  Given $a$ in any $A_g$, we then have that $$ \Vert \varphi _g(a)\Vert ^2 = \Vert \varphi _g(a)^*\varphi
_g(a)\Vert = \Vert \varphi _1(a^*a)\Vert = \Vert a^*a\Vert = \Vert a\Vert ^2.  $$ This proves that $\varphi _g$ is
isometric.

Supposing now that all of the $\varphi _g$ have dense range, it is easy to see that the *-homomorphism $ \funcred
\varphi : \Cr \Abun \to \CrB , $ given by \cite {FunctorialityRed}, has dense range.  Since $\funcred \varphi $ is a
*-homomorphism between C*-algebras, we deduce that $\funcred \varphi $ is surjective.  Thus, given $b$ in any $B_g$,
there exists some $y$ in $\Cr \Abun $ such that $$ \funcred \varphi (y) = \lambda ^B_g(b).  $$ Observe that $\funcred
\varphi $ satisfies $$ \varphi _g\big (E_g^A(x)\big ) = E_g^B\big ( \funcred \varphi (x)\big ) \for g\in G\for x\in \Cr
\Abun , $$ where $E_g^A$ and $E_g^B$ are given by \cite {Fourier}.  This may be easily checked by first assuming that
$x=\lambda _h(a)$, for $a$ in any $A_h$.  Thus $$ \varphi _g\big (E_g^A(y)\big ) = E_g^B\big ( \funcred \varphi (y)\big
) = E_g^B\big (\lambda ^B_g(b)\big ) = b, $$ proving that $\varphi _g$ is surjective.  \endProof

\definition \label DefineSubBun Let $\Bun = \{B_g\}_{g\in G}$ be a Fell bundle.  By a \subj {\subFellBun dle} of $\Bun $
we mean a collection $\Abun = \{A_g\}_{g\in G}$ of closed subspaces $A_g\subseteq B_g$, such that \izitem \zitem
$A_gA_h\subseteq A_{gh}$, and \zitem $A_g^*\subseteq A_{g\inv }$, \medskip \noindent for all $g$ and $h$ in $G$.

It is evident that a \subFellBun dle is itself a Fell bundle with the restricted operations, and moreover the inclusion
from $\Abun $ into $\Bun $ is a morphism.  Thus, given a \subFellBun dle $\Abun $ of a Fell bundle $\Bun $, from \cite
{FunctorialityFull} and \cite {FunctorialityRed} we get *-homomorphisms $$ C^*(\Abun ) \> \iota \CB \and \Cr \Abun \>
{\funcred \iota } \CrB , \equationmark SubBunMaps $$ the second of which is injective by the last sentence of \cite
{FunctorialityRed}.

This of course raises the question as to whether $\iota $ is also injective but, surprisingly, this is not always the
case.  We would now like to discuss a rather anomalous phenomenon leading up to the failure of injectivity for $\iota $.

Recall that a C*-algebra $C$ is said to be \subjex {exact}{exact C*-algebra} if, whenever $$ 0 \to J \to A \to B \to 0
$$ is an exact sequence of C*-algebras, the corresponding sequence $$ 0 \to J\tmin C \to A\tmin C \to B\tmin C \to 0 $$
is also exact.  Likewise, a group $G$ is said to be \subjex {exact}{exact group} if $\Cr G$ is an exact C*-algebra.  See
\ref {Chapters 2 and 5/BrownOzawa/2008} for more details.

All amenable groups are exact, but there are many non-amenable exact groups, such as the free group on two or more
generators \ref {Proposition 5.1.8/BrownOzawa/2008}.  The first examples of non-exact groups were found by Gromov in
\ref {Gromov/2003}, shattering a previously held belief that all groups were exact.

\state Proposition \label GoodSubBunOnlyAmenable Let $G$ be an exact group.  Then the following are equivalent: \izitem
\zitem for every Fell bundle $\Bun $ over $G$, and every \subFellBun dle $\Abun $ of $\Bun $, the map $\iota $ of \cite
{SubBunMaps} is injective, \zitem $G$ is amenable.

\Proof By \ref {Theorem 5.1.7/BrownOzawa/2008} there exists a compact space $X$ and a global amenable action $\beta $ of
$G$ on $X$.  Letting $\Bun $ be the semi-direct product bundle for the corresponding action of $G$ on $C(X)$, we have
that $\Bun $ satisfies the approximation property by \cite {AmenaActApprox}, and hence $\CB $ is naturally isomorphic to
$\CrB $ via the regular representation $\Lambda $, by \cite {Main}.

Since $\beta $ is a global action, and $C(X)$ is a unital algebra, $1\delta _g$ is an element of $B_g$, for each $g$ in
$G$, and we may then consider the \subFellBun dle $\Abun = \{A_g\}_{g\in G}$ of $\Bun $ given by $$ A_g = {\bf C}\delta
_g \for g\in G.  $$ By (i) we then have that the map $$ \iota : C^*(\Abun ) \to \CB $$ given by \cite {SubBunMaps} is
injective.  Therefore the composition $$ C^*(\Abun ) \>\iota \CB \> \Lambda \CrB $$ is also injective.  Letting $\tau
=E_1\circ \Lambda \circ \iota $, where $E_1$ is given by \cite {Fourier}, it is easy to see that $\tau $ takes values in
the subspace ${\bf C}\delta _1\subseteq B_1$, so we may view $\tau $ as a complex linear functional on $C^*(\Abun )$,
which is faithful in the sense that $$ \tau (a^*a)=0 \imply a=0 \for a\in C^*(\Abun ), $$ because $E_1$ is a faithful
conditional expectation by \cite {EOneIsFaith}.

Recall that $C^*(\Abun )$ is a topologically graded algebra by \cite {CBCrBTopGraded}, and it is easy to see that $\tau
$ is its canonical conditional expectation.  Therefore it follows from \cite {KerLRR} that the regular representation $$
\Lambda ^\Abun :C^*(\Abun )\to \Cr \Abun $$ is injective.  Since the full \resp {reduced} cross sectional C*-algebra of
$\Abun $ is nothing but the full \resp {reduced} group C*-algebra of $G$, we conclude from \ref {Theorem
2.6.8/BrownOzawa/2008} that $G$ is amenable.

The converse implication is easily proven by observing that for both $\Bun $ and $\Abun $, the reduced and full cross
sectional algebras coincide by \cite {BundleAmenaGroup}, so the result follows from the last sentence of \cite
{FunctorialityRed}.  \endProof

\nostate \label InjectiveNot Given any non-amenable exact group $G$, such as the free group on two generators, the
result above implies that there is a Fell bundle $\Bun $ over $G$, and a \subFellBun dle $\Abun $ of $\Bun $, for which
the canonical mapping $\iota $ of \cite {SubBunMaps} is non-injective.

\medskip The reader might have noticed that the last sentence of \cite {FunctorialityRed}, relating to the injectivity
of $\funcred \varphi $, has no counterpart in \cite {FunctorialityFull}.  In fact, as the above analysis shows, it is
impossible to add to \cite {FunctorialityFull} a sentence similar to the last one of \cite {FunctorialityRed}.

\medskip From our perspective, one of the most important examples of \subFellBun dles is as follows:

\state Proposition \label SubFellFromInvar Let $ \beta = \big (B, G, \{B_g\}_{g\in G}, \{\beta _g\}_{g\in G}\big ) $ be
a C*-algebraic partial dynamical system, and let $A$ be a $\beta $-invariant, closed *-sub\-algebra of $B$.  Also let $
\alpha = \big (\{A_g\}_{g\in G},\{\alpha _g\}_{g\in G}\big ) $ be the restriction of $\beta $ to $A$, as defined in
\cite {RestrToInvar}.  Letting $\Abun $ and $\Bun $ be the semi-direct product bundles relative to $\alpha $ and $\beta
$, respectively, one has that $\Abun $ is a \subFellBun dle of $\Bun $ (see below for the case when $A$ is not
invariant).

\Proof Recall that $A_g$ is defined by $$ A_g = A\cap B_g, $$ so evidently $A_g\delta _g$ is a closed subspace of
$B_g\delta _g$.  We leave it for the reader to verify that the multiplication and adjoint operations of $\Abun $ are
precisely the restrictions of the corresponding operations on $\Bun $, thus proving that $\Abun $ is indeed a
\subFellBun dle of $\Bun $, as desired.  \endProof

\definition \label HeredIdealBun Let $\Bun = \{B_g\}_{g\in G}$ be a Fell bundle and let us be given a \subFellBun dle
$\Abun = \{A_g\}_{g\in G}$ of $\Bun $.  \iaitem \aitem We say that $\Abun $ is a \subj {hereditary \subFellBun dle} of
$\Bun $ if $$ A_gB_hA_k\subseteq A_{ghk} \for g,h,k\in G.  $$ \aitem We say that $\Abun $ is an \subjex {ideal}{ideal of
a Fell bundle} of $\Bun $ if $$ A_gB_h\subseteq A_{gh} \and B_gA_h\subseteq A_{gh} \for g,h\in G.  $$

It is evident that every ideal is a hereditary \subFellBun dle.

\medskip A very important source of examples of hereditary \subFellBun dles is given by the process of restriction of
global actions to (not necessarily invariant) ideals, as we will now show:

\state Proposition \label SubFellFromIdeal Let $\beta $ be a global C*-algebraic action of a group $G$ on a C*-algebra
$B$.  Given a closed two-sided ideal $A$ of $B$, let $$ \alpha = \big (\{\D g\}_{g\in G},\{\alpha _g\}_{g\in G}\big ) $$
be the restriction of $\beta $ to $A$, as defined in \cite {DefineRestriction}.  Denoting by $\Abun $ and $\Bun $ the
semi-direct product bundles relative to $\alpha $ and $\beta $, respectively, one has that $\Abun $ is a hereditary\fn
{One should not be misled into thinking that $\Abun $ is an ideal in $\Bun $, despite the fact that $A$ is an ideal in
$B$.}  \subFellBun dle of $\Bun $.

\Proof Recall that $\D g$ is defined by $$ \D g = A\cap \beta _g(A), $$ so $\D g$ is a closed subspace of $B$, and we
may then see $\D g\delta _g$ is a closed subspace of $B\delta _g$.  We leave it for the reader to verify that the
multiplication and adjoint operations of $\Abun $ are precisely the restrictions of the corresponding operations on
$\Bun $, so that $\Abun $ is indeed a \subFellBun dle of $\Bun $, as desired.

In order to prove that $\Abun $ is hereditary, suppose that we are given $g$, $h$ and $k$ in $G$, $a$ in $\D g$, $b$ in
$B$, and $c$ in $\D k$.  Then $$ (a\delta _g)(b\delta _h)(c\delta _k) = a\beta _g\big (b)\beta _{gh}(c)\delta _{ghk}.
\equationmark VerifHered $$

Observing that each $\D g$ is an ideal in $B$, notice that the coefficient of $\delta _{ghk}$ above satisfies $$ a\beta
_g\big (b)\beta _{gh}(c) \in \D g \cap \beta _{gh}(\D k) = A \cap \beta _g(A) \cap \beta _{gh}\big (A \cap \beta
_k(A)\big ) \$= A \cap \beta _g(A) \cap \beta _{gh}(A) \cap \beta _{ghk}(A)\subseteq A \cap \beta _{ghk}(A) = \D {ghk},
$$ so the element mentioned in \cite {VerifHered} in fact lies in $D_{ghk}\delta _{ghk}$, proving $\Abun $ to be
hereditary.  \endProof

\state Proposition \label HeredIdeal Let $\Bun $ be a Fell bundle and let $\Abun $ be a hereditary \sFb dle \resp
{ideal} of $\Bun $.  Then $\Cr \Abun $ is a hereditary subalgebra \resp {ideal} of $\CrB $.

\Proof We first observe that $\Cr \Abun $ may indeed be seen as a subalgebra of $\CrB $ by \cite {FunctorialityRed}.  It
is then easy to see that $$ C_c(\Abun ) C_c(\Bun ) C_c(\Abun ) \subseteq C_c(\Abun ), $$ in the hereditary case, and $$
C_c(\Abun ) C_c(\Bun ) \subseteq C_c(\Abun ) \and C_c(\Bun ) C_c(\Abun ) \subseteq C_c(\Abun ), $$ in the ideal case,
from where the result follows.  \endProof

As seen in \cite {InjectiveNot}, the canonical map $\iota $ of \cite {SubBunMaps} is not always injective. However under
more restrictive conditions such pathologies disappear:

\state Theorem \label FullInjHered Let $\Bun $ be a Fell bundle and let $\Abun $ be a hereditary \subFellBun dle.  Then
the canonical map $$ \iota : C^*(\Abun ) \to \CB $$ of \cite {SubBunMaps} is injective, so $C^*(\Abun )$ is naturally
isomorphic to the range on $\iota $, which is a hereditary *-subalgebra of $\CB $.

\Proof Consider the diagram

\null \hfill \beginpicture \setcoordinatesystem units <0.0020truecm, -0.0025truecm> point at 1000 1000 \put {$C^*(\Abun
)$} at -1000 0000 \arrow <0.11cm> [0.3,1.2] from -550 000 to 100 000 \put {$\weird \,'$} at -200 -180 \put {$\Cr \Abun
\tmax \Cr G$} at 1230 0000 \put {$C^*(\Bun )$} at -1000 1000 \arrow <0.11cm> [0.3,1.2] from -550 1000 to 100 1000 \put
{$\weird $} at -200 820 \put {$\CrB \tmax \Cr G$} at 1250 1000 \arrow <0.11cm> [0.3,1.2] from -1000 240 to -1000 760
\put {$\iota $} at -1200 450 \arrow <0.11cm> [0.3,1.2] from 1100 240 to 1100 760 \put {$\funcred \iota \otimes id$} at
1550 450 \endpicture \hfill \null

\bigskip \noindent where $\weird \,'$ and $\weird $ are given by \cite {CrazyFact} for the respective bundles.  By
checking on the dense set $C_c(\Abun )\subseteq C^*(\Abun )$, it is easy to see that this is a commutative diagram.

Employing \ref {Proposition 3.6.4/BrownOzawa/2008} we conclude that $\funcred \iota \otimes id$ is injective because
$\funcred \iota $ is injective by the last sentence of \cite {FunctorialityRed}, and $\Cr \Abun $ is a hereditary
subalgebra of $\CrB $, by \cite {HeredIdeal}.  Since $\weird \,'$ and $\weird $ are also injective, we conclude that
$\iota $ is injective.  That $C^*(\Abun )$ is a hereditary subalgebra may be proved easily by the same argument used in
the proof of \cite {HeredIdeal}.  \endProof

Given an ideal $\Jbun = \{J_g\}_{g\in G}$ of a Fell bundle $\Bun = \{B_g\}_{g\in G}$, consider, for each $g$ in $G$, the
quotient space $B_g/J_g$.  It is clear that the operations of $\Bun $ drop to the quotient, giving multiplication
operations $$ \cdot : {B_g \over J_g} \times {B_h\over J_h}\ \to \ {B_{gh} \over J_{gh}} $$ and involutions $$ * : {B_g
\over J_g}\ \to \ {B_{g\inv }\over J_{g\inv }}.  $$

One then easily verifies that the collection $$ \quo := \{B_g/J_g\}_{g\in G} $$ is a Fell bundle with the above
operations.

\definition Given an ideal $\Jbun $ of a Fell bundle $\Bun $, the Fell bundle $\quo $ constructed above will be called
the quotient of $\Bun $ by $\Jbun $.

Let us now compare the cross sectional C*-algebra of the quotient Fell bundle with the quotient of the corresponding
cross sectional C*-algebras:

\state Proposition \label FullExactSeq If $\Jbun $ is an ideal in the Fell bundle $\Bun $, then $C^*(\Jbun )$ is an
ideal in $\CB $.  Moreover the quotient of $\CB $ by $C^*(\Jbun )$ is isomorphic to $C^*(\quo )$, thus giving an exact
sequence of C*-algebras $$ 0 \longrightarrow C^*(\Jbun ) \> \iota C^*(\Bun ) \> q C^*(\quo ) \longrightarrow 0.  $$

\Proof Since ideals are necessarily hereditary \sFb dles, \cite {FullInjHered} applies and hence $C^*(\Jbun )$ is
naturally isomorphic to a subalgebra of $\CB $, which the reader may easily prove to be an ideal using the same argument
employed in the proof of \cite {HeredIdeal}.

For each $g$ in $G$, let us denote by $$ q_g:B_g\to B_g/J_g $$ the quotient map.  It is then obvious that
$q=\{q_g\}_{g\in G}$ is a morphism from $\Bun $ to $\quo $, which therefore induces via \cite {FunctorialityFull} a
clearly surjective *-homo\-morphism, $$ q:\CB \to C^*(\quo ), $$ still denoted by $q$ by abuse of language, such that $$
q(\hj _g(b)) = \hj _g(b+J_g) \for g\in G\for b\in B_g, $$ where we rely on the context to determine which universal
representation $\hj $ is meant in each case.  Notice that $q$ vanishes on $C^*(\Jbun )$, so $$ C^*(\Jbun )\subseteq \Ker
(q).  $$ The proof will then be complete once we show that these sets are in fact equal.  For each $g$ in $G$, define $$
\pi ^0_g:b\in B_g \mapsto \hj _g(b) + C^*(\Jbun ) \in \CB /C^*(\Jbun ).  $$ Since $\pi ^0_g$ clearly vanishes on $J_g$,
it factors through $B_g/J_g$, giving a linear mapping $$ \pi _g:B_g/J_g \to \CB /C^*(\Jbun ).  $$ One may now check that
$\{\pi _g\}_{g\in G}$ is a representation of $\quo $, the integrated form of which, according to \cite
{FromRepBunToRelAlg}, is a *-homomorphism $$ \varphi :C^*(\quo ) \to \CB /C^*(\Jbun ), $$ such that $$ \varphi \big (\hj
_g(b+J_g)\big ) = \hj _g(b) + C^*(\Jbun ) \for g\in G \for b\in B_g, $$ where again the context should be enough to
determine the appropriate versions of $\hj $.  The composition $$ C^*(\Bun ) \> {q} C^*(\quo ) \> \varphi \CB /C^*(\Jbun
) $$ thus sends each $\hj _g(b)$ to $\hj _g(b) + C^*(\Jbun )$, which means that it is precisely the quotient map modulo
$C^*(\Jbun )$.  This immediately implies that $$ \Ker (q)\subseteq C^*(\Jbun ), $$ thus concluding the proof.  \endProof

\state Remark \label LeComLe \rm It is interesting to observe that \"{reduced} cross sectional algebras behave well with
respect to \"{\subFellBun dles} by \cite {SubBunMaps}, while \"{full} cross sectional algebras behave well with respect
to \"{exact sequences} by \cite {FullExactSeq}.

Interchanging the roles of \"{full} and \"{reduced}, we have already obtained a partial result regarding \"{\subFellBun
dles} and \"{full} cross sectional algebras in \cite {FullInjHered}, so we will now discuss the behavior of \"{exact
sequences} under \"{reduced} cross sectional algebras.  We begin with a technical result:

\state Lemma \label LemaDoPsi Given a Fell bundle $\Bun = \{B_g\}_{g\in G}$, there exists a bounded completely positive
linear map $$ \Psi : C^*(\Bun )\tmin \Cr G \to \CrB , $$ such that, for every $g$ and $h$ in $G$, and every $b$ in $B_g$
one has that $$ \Psi \big (\hj _g(b)\*\lreg _h\big ) = \left \{ \matrix { \lambda _g(b), & \hbox { if $g=h$,} \cr \pilar
{14pt} 0, & \hbox { if $g\neq h$.}}\right .  $$ Consequently $\Psi \circ \tau $ coincides with the identity of $\CrB $,
where $\tau $ is given by \cite {FullMiniAbsorpt}.

\Proof Let $\pi $ be any *-representation of $\Bun $ on a Hilbert space $H$ such that $\pi _1$ is faithful.  This can
easily be obtained by composing a faithful representation of $\CB $ with the universal representation $\hj $, for
example.  Also let $$ \varphi : \CB \to \Lin \big (H\otimes \ell ^2(G)\big ) $$ be built from $\pi $, as in \cite
{IntegrTensorProdRep}, so that $$ \varphi \big (\hj _g(b)\big ) = \pi _g(b)\otimes \lreg _g \for g\in G \for b\in B_g.
$$ Recalling that $\Cr G$ is the closed *-subalgebra of $\Lin \big (\ell ^2(G)\big )$ generated by the range of the
left-regular representation of $G$, consider the representation $$ \varphi \otimes id: \CB \tmin \Cr G \to \Lin \big
(H\otimes \ell ^2(G)\otimes \ell ^2(G)\big ), $$ characterized by the fact that $$ (\varphi \otimes id)(b\otimes \lreg
_h) = \pi _g(b)\otimes \lreg _g\otimes \lreg _h, $$ for any $g$ and $h$ in $G$, and any $b$ in $B_g$.  Let $K$ be the
closed subspace of $H\otimes \ell ^2(G)\otimes \ell ^2(G)$, given by $$ K = \segura {\bigoplus }_{g\in G} H\*e_g\*e_g,
$$ and let $P$ be the orthogonal projection onto $K$.  Also consider the completely positive map $$ \Psi : C^*(\Bun
)\tmin \Cr G \to \Lin (K) $$ given by $\Psi (x) = P\Big ( (\varphi \otimes id) (x) \Big )P$.

There is an obvious isometric isomorphism between $K$ and $H\*\ell _2(G)$ under which a vector of the form $\xi
\*e_g\*e_g$ is mapped to $\xi \*e_g$.  If we identify $\Lin (K)$ with $\Lin \big (H\*\ell _2(G)\big )$ under this map,
one sees that, $$ \Psi (\hj _g(b)\*\lreg _h) = \left \{ \matrix { \pi _g(b)\*\lreg _g, & \hbox { if $g=h$,} \cr \pilar
{14pt} 0, & \hbox { if $g\neq h$,}}\right .  $$ for every $g$ and $h$ in $G$, and every $b$ in $B_g$.

The range of $\Psi $ is therefore the closed linear span of the set of all $\pi _g(b)\*\lreg _g$, for $b$ in $B_g$.
This is also the range of $\varphi $, which in turn is the same as the range of the map $\psi $ of \cite
{FellAbsorption}.

Since $\psi $ is a faithful representation by the last sentence of \cite {FellAbsorption}, given that $\pi _1$ is
faithful by construction, we may then view $\Psi $ as taking values in $\CrB $, thus providing the desired map.
\endProof

Let us now discuss an important consequence of the exactness of the base group to Fell bundles:

\state Theorem \label ReduExactSeq Let $G$ be an exact group, and let $\Bun $ be a Fell bundle over $G$.  If $\Jbun $ is
an ideal in $\Bun $, then the quotient of $\CrB $ by $\Cr \Jbun $ is isomorphic to $\Cr \quo $, thus yielding an exact
sequence of C*-algebras $$ 0 \longrightarrow \Cr \Jbun \> {\iota \lred } \Cr \Bun \> {q\lred } \Cr \quo \longrightarrow
0.  $$

\Proof Recall from \cite {HeredIdeal} that $\iota \lred $ is a natural isomorphism from $\Cr \Jbun $ onto an ideal in
$\Cr \Bun $.  On the other hand, the map $q\lred $ referred to above is the one given by \cite {FunctorialityRed} in
terms of the quotient morphism $q =\{q_g\}_{g\in G}$.

As in \cite {FullExactSeq}, it is immediate to check that $q\lred $ is surjective and that the range of $\iota \lred $
is contained in the kernel of $q\lred $.  Therefore the only point requiring our attention is the proof that the range
of $\iota \lred $ contains the kernel of $q\lred $.  So, let us pick $z$ in $\CrB $, such that $$ q\lred (z) =0.  $$
Temporarily turning to full cross sectional algebras, recall from \cite {FullExactSeq} that $$ 0 \longrightarrow
C^*(\Jbun ) \> \iota C^*(\Bun ) \> q C^*(\quo ) \longrightarrow 0 $$ is an exact sequence.  Since $\Cr G$ is an exact
C*-algebra by hypothesis, the middle row of the diagram below is exact at $C^*(\Bun )\tmin \Cr G$.

$$ \def \flxc {\big \uparrow } \def \tmin {\!\mathop {\otimes }\limits \lmin \!}  \def \hold {\vrule height 16pt depth
8pt width 0pt} \matrix { \Cr \Jbun & \>{\ds \funcred \iota } & \CrB &\cr \Psi _1 \flxc \qquad &&\Psi _2 \flxc \qquad
\hold \cr C^*(\Jbun )\tmin \Cr G & \>{\ds \iota \*1} & C^*(\Bun )\tmin \Cr G & \>{\ds q\*1}& C^*(\quo )\tmin \Cr G\cr
&&\tau _2 \flxc \qquad &&\tau _3 \flxc \quad \hold \cr && \CrB & \>{\ds q\lred } & \Cr {\quo } } $$

\medskip \noindent In this diagram we have also marked the maps $\Psi _1$ and $\Psi _2$ given by \cite {LemaDoPsi} for
$\Jbun $ and $\Bun $, respectively, as well as the maps $\tau _2$ and $\tau _3$ given by \cite {FullMiniAbsorpt} for
$\Bun $ and $\quo $.  The usual method of checking on the appropriate dense subalgebras easily shows that the diagram is
commutative.

Since $q\lred (z)=0$, we have that $(q\*1)\big (\tau _2(z)\big )=0$, so by exactness there exists $x$ in $C^*(\Jbun
)\tmin \Cr G$ such that $$ (\iota \otimes 1)(x) = \tau _2(z).  $$ Therefore $$ \iota \lred \big (\psi _1(x)\big ) = \psi
_2\big ((\iota \otimes 1)(x)\big ) = \psi _2\big (\tau _2(z)\big ) \={LemaDoPsi} z, $$ proving $z$ to be in the range of
$\iota \lred $, as desired.  \endProof

The above result should be seen from the point of view of \cite {LeComLe}.  That is, even though \"{reduced} cross
sectional algebras and \"{exact sequences} are not the best friends, under the assumption that the group is exact, we
get a satisfactory result.

We will now go back to discussing another slightly disgruntled relationship, namely that of \"{full} cross sectional
algebras and \"{\subFellBun dles}, a relationship that has already appeared in \cite {FullInjHered}, where heredity
proved to be the crucial hypothesis.  Instead of heredity we will now work under the existence of conditional
expectations, as defined below.

\definition \label DefineConExBun Let $\Bun = \{B_g\}_{g\in G}$ be a Fell bundle and $\Abun = \{A_g\}_{g\in G}$ be a
\subFellBun dle.  A \subjex {conditional expectation}{conditional expectation in Fell bundles} from $\Bun $ to $\Abun $
is a collection of maps $$ P = \{P_g\}_{g\in G}, $$ where each $$ P_g: B_g \rightarrow B_g $$ is a bounded, idempotent
linear mapping, with range equal to $A_g$, such that $P_1$ is a conditional expectation from $B_1$ to $A_1$ and, for
every $g$ and $h$ in $G$, every $b$ in $B_g$, and every $c$ in $B_h$, one has that \izitem \zitem $P_g(b)^*=P_{g\inv
}(b^*)$, \zitem $P_{gh}(bc)=bP_h(c)$, provided $b\in A_g$, \zitem $P_{gh}(bc)=P_g(b)c$, provided $c\in A_h$.

Observe that \cite {DefineConExBun/iii} easily follows from \cite {DefineConExBun/i--ii} by taking adjoints.

\medskip \fix From now on, let us assume that we are given a Fell bundle $\Bun = \{B_g\}_{g\in G}$, and a \subFellBun
dle $\Abun = \{A_g\}_{g\in G}$, admitting a conditional expectation $$ P = \{P_g\}_{g\in G}.  $$

It is our next goal to show that there exists a conditional expectation from $\CrB $ to $\Cr \Abun $ extending each
$P_g$.  This result will later be used in conjunction with \cite {CrazyFact} to give us the desired embedability result
for full crossed sectional algebras.

The method we will adopt uses a version of the Jones-Watatani basic construction \ref {Watatani/1990}.  As the first
step we will construct a right Hilbert $A_1$-module from $\CCB $, vaguely resembling the $B_1$-module $\ell ^2(\Bun )$
defined in \cite {DefineEllTwo}.  Given $y$ and $z$ in $\CCB $, define $$ \ip Pyz = \soma {g\in G} P_1\big
((y_g)^*z_g\big ).  $$

We have already seen that $\CCB $ has the structure of a right $B_1$-module via the standard inclusion $j_1:B_1\to \CCB
$.  Restricting this module structure to $A_1$, we may view $\CCB $ as a right $A_1$-module.  It is then easy to prove
that $\ip P\ponto \ponto $ is an $A_1$-valued pre-inner product, the positivity following from \cite {PositivityInFell}
and the assumption that $P_1$ is positive.

\definition We shall denote by $\lpb $ the right Hilbert $A_1$-module obtained by completing $\CCB $ under the semi-norm
$\Vert \ponto \Vert _{2,\smallp }$ arising from the inner-product defined above (after modding out the subspace of
vectors of length zero).  For each $b$ in any $B_g$, we will denote by $\jp _g(b)$ the canonical image of $j_g(b)$ in
$\lpb $.

The following result is a version of \cite {PresentLambdaG} to the present situation:

\state Proposition \label PresentLambdaGP Given $b$ in any $B_g$, the operator $$ \lambda _g(b) : \ y \ \in \ \CCB \
\longmapsto \ j_g(b)\star y \ \in \ \CCB $$ is bounded relative to $\Vert \ponto \Vert _{2,\smallp }$ and hence extends
to a bounded operator on $\lpb $, which we will denote by $\lp _g(b)$, such that $\Vert \lp _g(b)\Vert \leq \Vert b\Vert
$, and which moreover satisfies $$ \lp _g(b)\big (\jp _h(c)\big ) = \jp _{gh}(bc) \for h\in G \for c\in B_h.  $$

\Proof The last assertion follows from the corresponding identity proved in \cite {PresentLambdaG}.  Addressing the
boundedness of $\lambda _g(b)$, given $y$ in $\CCB $, notice that $$ \bip P {\lambda _g(b)y}{\lambda _g(b)y} = \soma
{h\in G} P_1\big ((y_{g\inv h})^* b^* b y_{g\inv h}\big ) \$= \soma {h\in G} P_1\big ((y_h)^* b^* b y_h\big ) \_\leq
{OrderNorm} \Vert b\Vert ^2\soma {h\in G} P_1\big ((y_h)^* y_h\big ) = \Vert b\Vert ^2\ip Py y, $$ from where the result
follows.  \endProof

Much like we did in \cite {DefineRegularRep}, one may now prove that the collection of maps $\lp = \{\lp _g\}_{g\in G}$
is a representation of $\Bun $ as adjointable operators on $\lpb $, and hence gives rise to an integrated form, which we
denote by $$ \Lp : \CB \to \AdjOp \big (\lpb \big ), $$ and which in turn is characterized by the fact that $$ \Lp \big
(\hj _g(b)\big ) = \lp _g(b) \for g\in G \for b\in B_g.  \equationmark CharacLp $$

\state Proposition \label LpVanishes The map $\Lp $ defined above vanishes on the kernel of the regular representation
$\Lambda $ of $\Bun $.

\Proof Mimicking \cite {MatrixCoef} one may prove that $$ \bip P {\jp _g(b)} {\Lp (z)\jp _h(c)} = P_1\Big (b^*E_{gh\inv
}\big ( \Lambda (z) \big )c\Big ), $$ for all $z\in \CB $, $g,h\in G$, $b\in B_g$, and $c\in B_h$.  Thus, if $z$ is in
the kernel of $\Lambda $, we have that the left-hand-side above vanishes identically, and this easily implies that $\Lp
(z)=0$.  \endProof

As a consequence, we see that $\Lp $ factors through $\CrB $, producing a representation $$ \Lp \lred : \CrB \to \AdjOp
\big (\lpb \big ), $$ such that $$ \Lp \lred \big (\lambda _g(b)\big ) = \lp _g(b) \for g\in G \for b\in B_g.
\equationmark CharacLpRed $$

In order to proceed, we need to prove the following key inequality:

\state Lemma \label DesigCondEx For $b$ in any $B_g$, one has that $$ P_g(b)^*P_g(b)\leq P_1(b^*b).  $$

\Proof Observing that $b-P_g(b)$ belongs to $B_g$, and hence by \cite {PositivityInFell} one has $$ \big (b-P_g(b)\big
)^*\big (b-P_g(b)\big )\geq 0 $$ in $B_1$, we deduce from the positivity of $P_1$ that $$ 0 \leq P_1\Big (\big
(b-P_g(b)\big )^*\big (b-P_g(b)\big )\Big ) \$= P_1\Big (b^*b-b^*P_g(b) - P_g(b)^*b + P_g(b)^*P_g(b)\Big ) \$=
P_1(b^*b)-P_{g\inv }(b^*)P_g(b) - P_g(b)^*P_g(b) + P_g(b)^*P_g(b) \$= P_1(b^*b)- P_g(b)^*P_g(b), $$ from where the
conclusion follows.  \endProof

Observe that $C_c(\Abun )$ may be seen as a right $A_1$-sub-module of $\CCB $ in a natural way.  Moreover, the inclusion
of the former into the latter is clearly isometric for the usual 2-norm on $C_c(\Abun ) \subseteq \ell ^2(\Abun )$, and
the norm $\Vert \ponto \Vert _{2,\smallp }$ on $\CCB $.  Consequently this map extends to an isometric
right-$A_1$-linear map $$ \V : \ell ^2(\Abun )\to \lpb .  $$

Adopting the policy of using single quotes for denoting the relevant maps for $\Abun $, let us denote by $$ j'_g: A_g
\to C_c(\Abun ) $$ the maps given by \cite {DefineJg}.  We may then characterize $\V $ by the fact that $$ \V \big
(j'_g(a)\big ) = \jp _g(a) \for a\in A_g.  \equationmark CharacV $$

Since bounded linear maps on Hilbert modules are not necessarily adjointable, we need a bit of work to provide an
adjoint for $\V $.

\state Proposition \label AdjInclusion The mapping $$ p: \CCB \to C_c(\Abun ) \subseteq \ell ^2(\Abun ), $$ defined by
$$ p(y)_g = P_g(y_g) \for y\in \CCB \for g\in G, $$ is bounded with respect to $\Vert \ponto \Vert _{2,\smallp }$, and
hence extends to a bounded operator from $\lpb $ to $\ell ^2(\Abun )$, still denoted by $p$, by abuse of language.
Moreover $p$ is the adjoint of the map $\V $ defined above.

\Proof Given $y$ in $\CCB $, one has $$ \bip {}{p(y)}{p(y)} = \soma {g\in G} P_g(y_g)^*P_g(y_g) \_\leq {DesigCondEx}
\soma {g\in G} P_1\big ((y_g)^*y_g\big ) = \ip P yy.  $$ This proves the boundedness of $p$.  In order to prove the last
assertion in the statement, let $y\in C_c(\Abun )$ and $z\in \CCB $.  Then $$ \ip {} y{p(z)} = \soma {h\in G}
y_g^*P_g(z_g) = \soma {h\in G} P_1(y_g^*z_g) = \ip P {\V (y)}z.  \omitDoubleDollar \endProof

Since $p = \V ^*$, we will from now on dispense with the notation ``$p$'', using ``$\V ^*$'' instead.

As already mentioned, $\V $ is an isometry, so $\V ^*\V $ is the identity operator on $\ell ^2(\Abun )$, while $\V \V
^*$ is a projection in $\AdjOp (\lpb )$, whose range is clearly the canonical copy of $\ell ^2(\Abun )$ within $\lpb $.

The following result describes what happens when we compress elements of the range of $\Lp $ down to $\ell ^2(\Abun )$
via $\V $.  In its statement we will continue with the trend of using single quotes when denoting the relevant maps for
$\Abun $.

\state Proposition \label HappensCompress For $b$ in any $B_g$, one has that $$ \V ^* \lp _g(b)\V = \lambda '_g\big
(P_g(b)\big ).  $$

\Proof Given $h$ in $G$, and $a$ in $A_h$ we have that $$ \Big ( \V ^* \lp _g(b)\V \Big )\modulogrande \,_{j'_h(a)}
\={CharacV} \V ^*\lp _g(b) \big (\jp _h(a)\big ) \={PresentLambdaGP} \V ^*\jp _{gh}(ba) \={AdjInclusion} $$$$ =
j'_{gh}\big (P_{gh}(ba)\big ) = j'_{gh}\big (P_g(b)a\big ) \={PresentLambdaG} \lambda '_g\big (P_g(b)\big )j'_h(a).  $$
Since the elements $j'_h(a)$ considered above span a dense subspace of $\ell ^2(\Abun )$, the proof is concluded.
\endProof

With this we may now prove our next main result:

\state Theorem \label CondExpOnRedAlg Given a Fell bundle $\Bun = \{B_g\}_{g\in G}$, and a \subFellBun dle $\Abun =
\{A_g\}_{g\in G}$ admitting a conditional expectation $ P = \{P_g\}_{g\in G}, $ there exists a conditional expectation
$$ E:\CrB \to \Cr \Abun , $$ such that $$ E\big (\lambda _g(b)\big )=\lambda '_g\big (P_g(b)\big ) \for g\in G \for b\in
B_g, $$ where $\lambda '$ and $\lambda $ denote the regular representations of $\Abun $ and $\Bun $, respectively.

\Proof Define $$ E: \CrB \to \AdjOp \big (\ell ^2(\Abun )\big ), $$ by $$ E(z) = \V ^*\Lp \lred (z)\V \for z\in \CrB .
$$

In case $z=\lambda _g(b)$, for some $b$ in $B_g$, notice that $$ E\big (\lambda _g(b)\big ) = \V ^*\Lp \lred \big
(\lambda _g(b)\big )\V \={CharacLpRed} \V ^*\lp _g(b)\V \={HappensCompress} \lambda '_g\big (P_g(b)\big ), $$ proving
the last assertion in the statement.  Since the $\lambda _g(b)$ span $\CrB $, this also proves that the range of $E$ is
contained in $\Cr \Abun $, so we may view $E$ as a $\Cr \Abun $-valued map, as required.

Conditional expectations are meant to be maps from an algebra to a subalgebra.  So it is important for us to view $\Cr
\Abun $ as a subalgebra of $\CrB $, and we of course do it through \cite {SubBunMaps}.  Technically this means that we
identify $\lambda '_g(a)$ and $\lambda _g(a)$, for $a$ in any $A_g$.  In this case, since $P_g(a)=a$, we deduce from the
equation just proved that $E\big (\lambda _g(a)\big )=\lambda _g(a)$, so $E$ is the identity on $\Cr \Abun $, and the
reader may now easily prove that $E$ is in fact a conditional expectation.  \endProof

We may now return to the discussion of \"{full} cross sectional algebras versus \"{\subFellBun dles}:

\state Theorem \label FullInjCondExp Let $\Bun $ be a Fell bundle and let $\Abun $ be a \subFellBun dle.  If $\Abun $
admits a conditional expectation, then the canonical map $$ \iota : C^*(\Abun ) \to \CB $$ of \cite {SubBunMaps} is
injective.

\Proof The proof is essentially identical to the proof of \cite {FullInjHered}, with the exception that, in order to
conclude that $\funcred \iota \otimes id$ is injective, instead of invoking \ref {Proposition 3.6.4/BrownOzawa/2008} we
rely on \ref {Proposition 3.6.6/BrownOzawa/2008}, observing that condition \ref {3.6.6.1/BrownOzawa/2008} is fulfilled
by the conditional expectation $E$ provided by \cite {CondExpOnRedAlg}.  \endProof

Let us now study amenability of \subFellBun dles.

\state Proposition Let $\Bun $ be an amenable Fell bundle and let $\Abun $ be a \subFellBun dle of $\Bun $.  If $\Abun $
is either hereditary or admits a conditional expectation, then $\Abun $ is also amenable.

\Proof Consider the diagram $$ \matrix { C^*(\Abun ) & \>{\ \iota \ } & \CB \cr \Lambda ' \big \downarrow && \big
\downarrow \Lambda \pilar {15pt} \stake {8pt} \cr \Cr {\Abun } & \>{\iota \lred } & \CrB \cr } $$ where $\Lambda '$ and
$\Lambda $ are the regular representations of $\Abun $ and $\Bun $, respectively, and $\iota $ and $\iota \lred $ are as
in \cite {SubBunMaps}.  It is elementary to check that this is a commutative diagram.

Employing either \cite {FullInjHered} or \cite {FullInjCondExp}, as appropriate, we have that $\iota $ is injective.
Since $\Bun $ is amenable, $\Lambda $ is injective as well, so we conclude that $\Lambda '$ is injective, whence $\Abun
$ is amenable.  \endProof

Using the method employed in the proof of \cite {GoodSubBunOnlyAmenable} it is possible to produce an amenable Fell
bundle containing a non-amenable \subFellBun dle.  Therefore the hypothesis that $\Abun $ is either hereditary or admits
a conditional expectation is crucial for the validity of the result above.

\medskip Let us now study the approximation property for \subFellBun dles.

\state Proposition \label SubBunHasAP Let $\Bun = \{B_g\}_{g\in G}$ be a Fell bundle satisfying the approximation
property, and let $\Abun = \{A_g\}_{g\in G}$ be a \subFellBun dle of $\Bun $.  If $A_1$ is a hereditary subalgebra of
$B_1$, then $\Abun $ also satisfies the approximation property.

\Proof Let $\{ a_i \}_{i\in I}$ be a {\aproxNet } for $\Bun $, as defined in \cite {AP}, and let $\{v_j\}_{j\in J}$ be
an approximate identity for $A_1$.  For each $(i,j)$ in $I\times J$, consider the function $$ c_{i,j}: g\in G \mapsto
v_ja_i(g)v_j \in A_1.  $$ Notice that $v_ja_i(g)v_j$ indeed lies in $A_1$ because $A_1B_1A_1\subseteq A_1$, given that
$A_1$ is supposed to be hereditary.

Considering $I\times J$ as a directed set with coordinate-wise order, we claim that $ \{c_{i,j}\}_{(i,j)\in I\times J} $
is a {\aproxNet } for $\Abun $.

Since each $a_i$ is a finitely supported function on $G$, it is clear that so are the $c_{i,j}$.  In order to prove
condition \cite {AP/i}, suppose we are given $(i,j)\in I\times J$.  Then, taking into account that $v_j^*v_j \leq 1$,
one has that $$ \normsum {\soma {g\in \G } c_{i,j}(g)^* c_{i,j}(g)} = \normsum {\soma {g\in \G } v_j^*a_i(g)^*v_j^*v_j
a_i(g) v_j} \$\leq \normsum {\soma {g\in \G } v_j^*a_i(g)^* a_i(g) v_j} \leq \normsum {\soma {g\in \G } a_i(g)^*
a_i(g)}, $$ so \cite {AP/i} follows from the corresponding property of the {\aproxNet } $\{ a_i \}_{i\in I}$.  In order
to prove \cite {AP/ii}, given $b$ in any $A_g$, we must prove that $$ \soma {h\in \G } c_{i,j}(gh)^* b c_{i,j}(h)
\convrg {(i,j)} b.  \equationmark TargetAmenaSub $$

For each $i$ in $I$, let us consider the linear operator $W_i$ defined on $A_g$ by $$ W_i(b) = \soma {h\in \G }
a_i(gh)^* b a_i(h) \for b\in A_g.  $$ Since this may be seen as the restriction of the map $V$ of \cite {WrongWay} to
$A_g$, we have that $$ \Vert W_i\Vert \leq \normsum {\soma {h\in \G } a_i(h)^*a_i(h)} \leq M, $$ where $M$ may be chosen
independently of $i$.  Focusing on the left-hand-side of \cite {TargetAmenaSub}, we have $$ \soma {h\in \G }
c_{i,j}(gh)^* b c_{i,j}(h) = \soma {h\in \G } v_j^*a_i(gh)^* v_j^*b v_ja_i(h)v_j = v_j^*W_i(v_j^*b v_j)v_j, $$ so \cite
{TargetAmenaSub} translates into $$ v_j^*W_i(v_j^*b v_j)v_j\convrg {(i,j)} b.  \equationmark TargetAmenaSubTwo $$

Observing that the terms of an approximate identity have norm at most 1, we have that $$ \Vert b - v_j^*W_i(v_j^*b
v_j)v_j\Vert \$\leq \Vert b - v_j^*b v_j\Vert + \Vert v_j^*b v_j - v_j^*W_i(b)v_j\Vert + \Vert v_j^*W_i(b)v_j -
v_j^*W_i(v_j^*b v_j)v_j \Vert \$\leq \Vert b - v_j^*b v_j\Vert + \Vert b - W_i(b)\Vert + \Vert W_i\Vert \Vert b-v_j^*b
v_j \Vert .  $$

By \cite {ApproximateIdentity} we have that $b v_j \convrg j b$, so \cite {TargetAmenaSubTwo} follows.  This verifies
that $\{c_{i,j}\}_{(i,j)\in I\times J}$ is indeed a {\aproxNet } for $\Abun $, as desired.  \endProof

Observe that the last hypothesis of the above result, namely that $A_1$ is a hereditary subalgebra of $B_1$, is
evidently satisfied if $\Abun $ is a hereditary \sFb dle of $\Bun $.  It might also be worth pointing out that the fact
that $A_1$ is a hereditary subalgebra of $B_1$ does not imply that $\Abun $ is a hereditary \subFellBun dle of $\Bun $.
A counter example may be easily obtained by choosing any Fell bundle $\Bun $ and selecting $A_1=B_1$, while $A_g=\{0\}$,
for all $g\neq 1$.  In this case there is no reason for $A_1B_gA_1$ to be contained in $A_g$.

After spending a little effort, we have not been able to determine if the approximation property passes to \sFb dles
under the existence of a conditional expectation, but we believe there is a good chance this is true.

\medskip With respect to amenability of quotient Fell bundles we have the following:

\state Proposition \label QuoAmena Suppose we are given Fell bundles $\Bun = \{B_g\}_{g\in G}$ and $\Abun =
\{A_g\}_{g\in G}$, as well as a morphism $ \varphi = \{\varphi _g\}_{g\in G} $ from $\Bun $ to $\Abun $, such that
$\varphi _g$ is onto $A_g$ for every $g$.  \izitem \zitem If $\Bun $ is amenable and $G$ is an exact group, then $\Abun
$ is also amenable.  \zitem If $\Bun $ satisfies the approximation property, then so does $\Abun $.

\Proof Given a {\aproxNet } $\{a_i\}_{i\in I}$ for $\Bun $, it is elementary to check that $\{\varphi _1(a_i)\}_{i\in
I}$ is a {\aproxNet } for $\Abun $.  Thus (ii) follows.

With respect to (i), for each $g$ in $G$, let $J_g$ be the kernel of $\varphi _g$.  It is then clear that $\Jbun =
\{J_g\}_{g\in G}$ is an ideal of $\Bun $, in the sense of \cite {HeredIdealBun/b}, and one has that $\Bun /\Jbun $ is
isomorphic to $\Abun $.

Consider the diagram $$ \def \flxc {\quad \big \downarrow } \matrix { 0 & \longrightarrow & C^*(\Jbun ) &
\longrightarrow & C^*(\Bun ) & \longrightarrow & C^*(\Abun ) & \longrightarrow & 0 \cr && \flxc \Lambda _{\!\Jbun } &&
\flxc \Lambda _\Bun && \flxc \Lambda _\Abun \pilar {15pt} \stake {9pt}\cr 0 & \longrightarrow & \Cr {\Jbun } &
\longrightarrow & \Cr {\Bun } & \longrightarrow & \Cr {\Abun } & \longrightarrow & 0 \cr } $$ whose rows are exact by
\cite {FullExactSeq} and \cite {ReduExactSeq}, and the vertical arrows indicate the various regular representations.  It
is immediate to check that this is commutative.  Moreover, $\Lambda _\Bun $ is one-to-one by hypothesis, and $\Lambda
_{\!\Jbun }$ is onto by definition of reduced cross sectional algebras, so one may prove that $\Lambda _\Abun $ is
injective by a standard diagram chase, similar to one used to prove the Five Lemma.  \endProof

\nrem Theorem \cite {ReduExactSeq} was proved in \ref {Exel/2002}.  It would be interesting to decide if \cite
{QuoAmena/i} holds without the hypothesis that $G$ is exact.

\chapter Functoriality for partial actions

\label FunctPASec

As already mentioned, some of the most important examples of Fell bundles arise from partial dynamical systems.  In this
chapter we will therefore present some consequences of the results proved in the above chapter to semi-direct product
bundles.  We will also study the relationship between the semi-direct product bundle for a partial action $\beta $ and
the corresponding bundle for a restriction of $\beta $.

\smallskip \fix So, let us now fix two C*-algebraic partial dynamical systems $$ \alpha = \big (A,G,\{A_g\}_{g\in
G},\{\alpha _g\}_{g\in G}\big ) $$ and $$ \beta = \big (B,G,\{B_g\}_{g\in G},\{\beta _g\}_{g\in G}\big ).  $$

We will denote the semi-direct product bundles for $\alpha $ and $\beta $ by $$ \Abun = \{A_g\delta _g\}_{g\in G} \and
\Bun = \{B_g\delta _g\}_{g\in G}, $$ respectively.

Recall from \cite {DefineEquivar} that a *-homomorphism $\varphi :A\to B$ is said to be $G$-equivariant provided
$\varphi (A_g) \subseteq B_g$, and $\varphi \big (\alpha _g(a)\big ) = \beta _g\big (\varphi (a)\big )$, for $a$ in any
$A_ g$.

\state Proposition Given a $G$-equivariant *-homomorphism $\varphi :A\to B$, for each $g$ in $G$ consider the mapping $$
\varphi _g:A_g\delta _g \to B_g\delta _g, $$ given by $$ \varphi _g(a\delta _g) = \varphi (a)\delta _g \for a\in A_g.
$$ Then $\{\varphi _g\}_{g\in G}$ is a morphism from $\Abun $ to $\Bun $.

\Proof Left to the reader.  \endProof

As an immediate consequence of \cite {FunctorialityFull} and \cite {FunctorialityRed}, we have:

\state Proposition \label FunctorPact Given a $G$-equivariant *-homomorphism $\varphi :A\to B$, there are
*-homomorphisms $$ \Phi : A\rt \alpha G \to B\rt \beta G \and \Phi \lred : A\redrt \alpha G \to B\redrt \beta G $$
sending $a\delta _g$ to $\varphi (a)\delta _g$, for $a$ in any $A_g$, with the appropriate interpretation of $\delta _g$
in each case.

In case $\varphi $ is injective one may see $\Abun $ as a \subFellBun dle of $\Bun $.  Two of the most important
examples of this are given by \cite {SubFellFromInvar} and \cite {SubFellFromIdeal}.

\state Proposition \label RestrGivesSubBun Let $ \beta = \big (B, G, \{B_g\}_{g\in G}, \{\beta _g\}_{g\in G}\big ) $ be
a C*-algebraic partial dynamical system, and let $A$ be a closed *-sub\-algebra of $B$. Suppose that either: \izitem
\zitem $A$ is invariant under $\beta $, or \zitem $\beta $ is a global action and $A$ is an ideal of $B$.  \medskip
\noindent In either case, let $ \alpha = \big (\{A_g\}_{g\in G},\{\alpha _g\}_{g\in G}\big ) $ be the restriction of
$\beta $ to $A$ (defined in \cite {RestrToInvar} in the first case, or in \cite {DefineRestriction} in the second case).
Then we have a natural inclusion $$ A\redrt \alpha G \subseteq B\redrt \beta G.  $$ In situation (ii) one moreover has
that $A\redrt \alpha G$ is a hereditary subalgebra of $B\redrt \beta G$.

\Proof We have that $\Abun $ is a \subFellBun dle of $\Bun $ by \cite {SubFellFromInvar} or \cite {SubFellFromIdeal}.
That $A\redrt \alpha G$ is naturally a subalgebra of $B\redrt \beta G$ then follows from the last assertion in \cite
{FunctorialityRed}.  Under hypothesis (ii), we have by \cite {SubFellFromIdeal} that $\Abun $ is a hereditary
\subFellBun dle of $\Bun $ so the last sentence in the statement follows from \cite {HeredIdeal}.  \endProof

For the case of full crossed products we have the following:

\state Proposition \label RestrGivesFullCP Let $\beta $ be a global action of a group $G$ on a C*-algebra $B$, and let
$A$ be an ideal of $B$.  Considering the partial action $\alpha $ of $G$ on $A$ given by restriction, as defined in
\cite {DefineRestriction}, the canonical mapping $$ \iota :A\rt \alpha G \to B\rt \beta G $$ is injective. Moreover the
range of $\iota $ is a hereditary *-subalgebra of $B\rt \beta G$.

\Proof We have seen in \cite {SubFellFromIdeal} that the semi-direct product bundle for $\alpha $ is hereditary in the
one for $\beta $.  So the conclusion follows immediately from \cite {FullInjHered}.  \endProof

In order to take advantage of \cite {FullInjCondExp}, namely the injectivity of the canonical mapping under the
existence of conditional expectations, we must understand the relationship between partial actions and conditional
expectations.  This is the goal of our next result.

\state Proposition Suppose we are given a C*-algebraic partial dynamical system $$ \beta = \big (B,G,\{B_g\}_{g\in G},
\{\beta _g\}_{g\in G}\big ), $$ and let $A$ be a closed *-subalgebra of $B$.  Suppose moreover that there is a
$G$-equivariant conditional expectation $F$ from $B$ onto $A$.  Then \izitem \zitem $A$ is invariant under $\beta $,
\zitem the semi-direct product bundle for the action $\alpha $ obtained by restricting $\beta $ to $A$, henceforth
denoted by $\Abun $, seen as a \subFellBun dle of the semi-direct product bundle for $\beta $, henceforth denoted by
$\Bun $, admits a conditional expectation $P = \{P_g\}_{g\in G}$, where $$ P_g(b\delta _g) = F(b)\delta _g \for b\in
B_g.  $$ \zitem the canonical mapping $$ \iota :A\rt \alpha G \to B\rt \beta G $$ is injective.

\Proof Given $g$ in $G$, we have that $$ \beta _g(A\cap B_{g\inv }) = \beta _g\big (F(A\cap B_{g\inv })\big )
\={DefineEquivar/ii} F\big (\beta _g(A\cap B_{g\inv })\big ) \subseteq A, $$ proving that $A$ is indeed invariant under
$\beta $.

Notice that if $b$ is in any $B_g$, then $F(b)\in B_g\cap A=:A_g$, by invariance, so indeed $F(b)\delta _g$ lies in
$A_g\delta _g$, as needed.  To prove \cite {DefineConExBun/i}, we compute $$ P_g(b\delta _g)^* = \big (F(b)\delta _g\big
)^* = \beta _{g\inv }\big (F(b^*)\big )\delta _{g\inv } \$= F\big (\beta _{g\inv }(b^*)\big )\delta _{g\inv } = P_{g\inv
}\big (\beta _{g\inv }(b^*)\delta _{g\inv }\big ) = P_{g\inv }\big ((b\delta _g)^*\big ).  $$ If we now take $c$ in some
$A_h$, we have $$ P_{gh}\big ((b\delta _g)(c\delta _h)\big ) = P_{gh}\Big (\beta _g\big (\beta _{g\inv }(b)c\big )\delta
_{gh}\Big ) = F\Big (\beta _g\big (\beta _{g\inv }(b)c\big )\Big )\delta _{gh} \$= \beta _g\Big (F\big (\beta _{g\inv
}(b)c\big )\Big )\delta _{gh} = \beta _g\Big (F\big (\beta _{g\inv }(b)\big )c\Big )\delta _{gh} = \beta _g\Big (\beta
_{g\inv }\big (F(b)\big )c\Big )\delta _{gh} \$= \big (F(b)\delta _g\big )\big (c\delta _h\big ) = P_g\big (b\delta
_g\big )\big (c\delta _h\big ), $$ proving \cite {DefineConExBun/iii} and then \cite {DefineConExBun/ii} follows by
taking adjoints.  The fact that $P_g(B_g\delta _g)= A\delta _g$, or equivalently that $F(B_g)=A_g$, may be easily proved
based on the $G$-equivariance of $F$.  Point (iii) now follows from \cite {FullInjCondExp}.  \endProof

Let us now study the behavior of crossed products relative to ideals and quotients.

\medskip \fix We therefore fix, for the time being, a C*-algebraic partial action $$ \alpha = \big (\{A_g\}_{g\in G},
\{\alpha _g\}_{g\in G}\big ) $$ of a group $G$ on a C*-algebra $A$, and an $\alpha $-invariant closed two-sided ideal
$J\ideal A$.  The restriction of $\alpha $ to $J$, defined according to \cite {RestrToInvar}, will henceforth be denoted
by $$ \tau = \big (\{J_g\}_{g\in G}, \{\tau _g\}_{g\in G}\big ).  $$

Denoting by $B$ the quotient C*-algebra, we obtain the exact sequence $$ 0 \to J \> i A \>q B \to 0, \equationmark
SomeExactSequence $$ where $\iota $ denoted the inclusion, and $q$ the quotient mapping.  Letting $$ B_g = q(A_g) \for
g\in G, $$ it is then easy to see that $B_g$ is a closed two-sided ideal in $B$ and that $B_g$ is *-isomorphic to
$A_g/(A_g\cap J)$.  Moreover, since each $\alpha _g$ maps $A_{g\inv }$ to $A_g$ in such a way that $$ \alpha _g(A_{g\inv
}\cap J) \subseteq A_g\cap J $$ by invariance of $J$, we see that $\alpha _g$ drops to the quotient providing a
*-homomorphism $$ \beta _g: B_{g\inv } \to B_g, $$ such that $$ \beta _g\big (q(a)\big ) = q\big (\alpha _g(a)\big )
\for g\in G \for a\in A_{g\inv }.  $$

\state Proposition \label QuotientPAct One has that $$ \beta := \big (\{B_g\}_{g\in G}, \{\beta _g\}_{g\in G}\big ) $$
is a partial action of $G$ on $B$.

\Proof Since \cite {PAIdentity} is trivially true, it suffices to verify \cite {SecondPADef/i--ii}.  Given $g,h\in G$,
we claim that $$ B_g\cap B_h = q(A_g\cap A_h).  \equationmark QuotientIntersec $$ By the Cohen-Hewitt Theorem every
ideal in a C*-algebra is idempotent.  In particular, given any two ideals $I$ and $J$, one has $$ I\cap J = (I\cap
J)(I\cap J) \subseteq IJ.  $$ Since $IJ$ is obviously contained in $I\cap J$, we deduce that $ I\cap J = IJ.  $ Thus,
for every $g,h\in G$, we have $$ B_g\cap B_h = B_gB_h = q(A_g)q(A_h) = q(A_gA_h) = q(A_g\cap A_h), $$ proving \cite
{QuotientIntersec}.  In order to prove \cite {PAIntersecContains}, we then compute $$ \beta _g\big (B_{g\inv } \cap
B_h\big ) = \beta _g\big (q(A_{g\inv } \cap A_h)\big ) \$= q\big (\alpha _g(A_{g\inv }\cap A_h)\big ) \subseteq q\big
(A_{gh}) = B_{gh}.  $$

Regarding \cite {PACompos}, pick $x$ in $B_{h\inv }\cap B_{(gh)\inv }$, and write $x = q(a)$, with $a$ in $A_{h\inv
}\cap A_{(gh)\inv }$, by \cite {QuotientIntersec}.  Then $$ \beta _g\big (\beta _h(x)\big ) = \beta _g\big (\beta
_h(q(a))\big ) = \beta _g\big (q(\alpha _h(a))\big ) = q\big (\alpha _g(\alpha _h(a))\big ) \$= q\big (\alpha
_{gh}(a)\big ) = \beta _{gh}\big (q(a)\big ) = \beta _{gh}(x).  $$ This completes the proof.  \endProof

Regarding the corresponding semi-direct product bundles we have:

\state Proposition \label QuoBundleandPA Let $ \alpha = \big (\{A_g\}_{g\in G}, \{\alpha _g\}_{g\in G}\big ) $ be a
C*-algebraic partial action of a group $G$ on a C*-algebra $A$, and let $J\ideal A$ be an $\alpha $-invariant closed
two-sided ideal.  Also let $\Jbun $, $\Abun $ and $\Bun $ be the semi-direct product bundles relative to the partial
actions $\tau $, $\alpha $ and $\beta $, above.  Then $\Jbun $ is an ideal in $\Abun $, and $\Bun $ is naturally
isomorphic to the quotient Fell bundle $\Abun /\Jbun $.

\Proof Left for the reader.  \endProof

We thus have the following consequences of our study of ideals in Fell bundles in the previous chapter:

\state Theorem \label ExactSeqCross Let $$ 0 \to J \> i A \> q B \to 0, $$ be an exact sequence of C*-algebras and let $
\alpha = \big (\{A_g\}_{g\in G}, \{\alpha _g\}_{g\in G}\big ) $ be a partial action of $G$ on $A$, relative to which $J$
is invariant.  Then the corresponding sequence $$ 0 \to J\rt \tau G \to A\rt \alpha G \to B\rt \beta G \to 0 $$ is also
exact.  Moreover, if \izitem \zitem $G$ is an exact group, or \zitem the semi-direct product bundle associated to
$\alpha $ satisfies the approximation property, \medskip \noindent then the sequence $$ 0 \to J\redrt \tau G \to A\redrt
\alpha G \to B\redrt \beta G \to 0 $$ is exact as well.

\Proof The first assertion is a direct application of \cite {FullExactSeq}.

Under hypothesis (i), the second assertion follows immediately from \cite {ReduExactSeq}, so we need only prove the
second assertion under hypothesis (ii).

Let $\Jbun $, $\Abun $ and $\Bun $ be the semi-direct product bundles relative to the partial actions $\tau $, $\alpha $
and $\beta $, as in \cite {QuoBundleandPA}.

By hypothesis we have that $\Abun $ satisfies the approximation property and hence so does $\Jbun $, by \cite
{SubBunHasAP}, as well as $\Bun $, by \cite {QuoAmena/ii}.  Therefore all three Fell bundles in sight are amenable by
\cite {Main}, so the above sequence of reduced crossed products coincides with the corresponding one for full crossed
products whose exactness has already been verified.  \endProof

This result has a useful application to the study of ideals in the crossed product generated by subsets of the
coefficient algebra.  We state it only for full crossed products since we only envisage applications of it in this case.

\state Proposition \label QuotientOfCP Let $$ \alpha = \big (A,G,\{A_g\}_{g\in G},\{\alpha _g\}_{g\in G}\big ) $$ be a
C*-algebraic partial dynamical system.  Given any subset $W \subseteq A$, \izitem \zitem let $K$ be the ideal of $A\rt
\alpha G$ generated by $\iota (W )$, where $\iota $ is the map defined in \cite {DefineStdInclusion}, and \zitem let $J$
be the smallest\fn {It is easy to see that an arbitrary intersection of invariant ideals is again invariant, so the
smallest invariant ideal always exist.}  $\alpha $-invariant ideal of $A$ containing $W $.  \medskip \noindent Then $K$
coincides with $J\rt \tau G $ (or rather, its canonical image within $A\rt \alpha G$).  In addition, there exists a
*-isomorphism $$ \varphi : {A\rt \alpha G \over K} \to \Big ({A \over J}\Big )\rt \beta G, $$ such that $$ \varphi
(a\delta _g + K) = (a+J)\delta _g \for g\in G \for a\in A_g.  $$

\Proof We first claim that $\iota \inv (K)$ is an $\alpha $-invariant ideal.  While it is evident that $\iota \inv (K)$
is an ideal, we still need to prove it to be $\alpha $-invariant.  Given $a\in \iota \inv (K)\cap A_{g\inv }$, we must
check that $\alpha _g(a)\in \iota \inv (K)$, which is to say that $\iota \big (\alpha _g(a)\big ) \in K$.  Using
Cohen-Hewitt, write $a = bc$, where both $b$ and $c$ lie in $\iota \inv (K)\cap A_{g\inv }$.  Then $$ \iota \big (\alpha
_g(a)\big ) = \alpha _g(a) \delta _1 = \alpha _g(bc) \delta _1 \={LotsAFormulas} \big (\alpha _g(b)\delta _g\big ) (c
\delta _1) = \big (\alpha _g(b)\delta _g\big ) \iota (c) \in K.  $$

So $\iota \inv (K)$ is indeed an invariant ideal, which evidently contains $W$.  Since $J$ is the smallest among such
ideals, we have that $J\subseteq \iota \inv (K)$, which is to say that $\iota (J)\subseteq K$, whence $$ \iota (W
)\subseteq \iota (J)\subseteq K.  \equationmark SetInJInK $$

This implies that $K$ coincides with the ideal of $A\rt \alpha G$ generated by $\iota (J)$.  On the other hand notice
that, by \cite {BOneEssential}, the ideal generated by $\iota (J)$ is $J\rt \tau G$, thus proving that $K= J\rt \tau G$.
The second and last assertion in the statement is now an immediate consequence of \cite {ExactSeqCross}.  \endProof

\chapter Ideals in graded algebras

\label IdealsChapter

Let $B$ be a graded C*-algebra with grading $\{B_g\}_{g\in G}$.  If $J$ is an ideal (always assumed to be closed and
two-sided) in $B$, there might be no relationship between $J$ and the grading of $B$.  It is even possible that $J\cap
B_g$ is trivial for every $g$ in $G$.

For example, let $\Bun $ be the group bundle over ${\bf Z}$, so that $\CB $ is isomorphic to $C({\bf T})$, where ${\bf
T}$ denotes the unit circle.  Fixing $z_0\in {\bf T}$, the ideal $$ J = \{f\in C({\bf T}): f(z_0)=0\} $$ has trivial
intersection with every homogeneous subspace $ B_n = {\bf C}z^n.  $ This is because a nonzero element in $B_n$ is
invertible, and hence cannot belong to any proper ideal.

\medskip The purpose of this chapter is thus to study the relationship between ideals in graded algebras and the grading
itself.  For this we shall temporarily

\fix fix a graded C*-algebra $B$, with grading $\{B_g\}_{g\in G}$, and a closed two-sided ideal $J\ideal B$.

\state Proposition The closed two-sided ideal of $B$ generated by $J\cap B_1$ coincides with the closure of \ $
\bigoplus _{g\in G} J\cap B_g.  $

\Proof Given $g$ and $h$ in $G$, notice that $$ (J\cap B_g)B_h \subseteq (JB_h)\cap (B_gB_h) \subseteq J\cap B_{gh}.  $$
Therefore we see that $$ K:= \overline {\pilar {9pt}\medoplus _{g\in G} J\cap B_g} $$ is invariant under right
multiplication by elements of $B_h$, and a similar reasoning shows invariance under left multiplication as well.  Since
$K$ is closed by definition, we then deduce that $K$ is a two-sided ideal.  Consequently, noticing that $J\cap B_1$ is
contained in $K$, we have that $$ \langle J \cap B_1\rangle \subseteq K, $$ where the angle brackets above indicate the
closed two-sided ideal generated by $J \cap B_1$.

In order to prove the reverse inclusion, given $g$ in $G$, and $x\in J\cap B_g$, notice that $x^*x\in J\cap B_1$.  Using
\cite {Roots} (which is stated for Hilbert modules, and hence also holds for C*-algebras) we then have $$ x = \lim
x(x^*x)^{1/n} \in \langle J \cap B_1\rangle , $$ We therefore conclude that $$ J\cap B_g\subseteq \langle J \cap
B_1\rangle , $$ from where it follows that $ K \subseteq \langle J \cap B_1\rangle , $ as desired.  \endProof

This justifies the introduction of the following concept:

\definition We shall say that $J$ is an \subj {induced ideal} (sometimes also called a \subj {graded ideal}) provided
any one of the following equivalent conditions hold: \iaitem \aitem $J$ coincides with the ideal generated by $J \cap
B_1$, \aitem $\bigoplus _{g\in G} J\cap B_g$ is dense in $J$.

For topologically graded algebras there is a lot more to be said, so we shall assume from now on that $B$ is
topologically graded.  Recall from \cite {Ft} that in this case $B$ admits \"{Fourier coefficient operators} $$ F_g: B
\to B_g \for g\in G, $$ such that $$ F_g(b) = \delta _{g,h}b \for g,h\in G \for b\in B_h.  $$ Given an arbitrary ideal
$J\ideal B$, let us consider the following subsets of $B$: $$ \matrix { J' & = & \langle J \cap B_1\rangle , \hfill \cr
\vrule height 12pt depth 8pt width0pt J'' & = & \left \{ x\in B: F_g(x) \in J, \ \forall g\in \G \right \},\cr J''' & =
& \left \{ x\in B: F_1(x^* x) \in J \right \},\hfill } \equationmark ThreeIdeals $$ where the angle brackets in the
definition of $J'$ are again supposed to mean the closed two-sided ideal generated.

\state Proposition \label InducedI Given any ideal $J$ in a topologically graded C*-algebra $B$, one has that the sets
$J'$, $J''$ and $J'''$ defined above are closed two-sided ideals in $B$, and moreover $$ J'\subseteq J''=J'''.  $$

\Proof It is evident that these are closed subspaces of $B$, and moreover that $J'$ is an ideal.

In order to prove that $J''$ is an ideal, let $b\in B$, and $x\in J''$, and let us prove that $bx\in J''$.  Since the
$B_g$ span a dense subspace of $B$, we may assume that $b\in B_h$, for some $h\in G$.  Then $$ F_g(bx) \={Ft} bF_{h\inv
g}(x) \in J, $$ proving that $bx\in J''$, and hence that $J''$ is a left ideal.  One similarly proves that $J''$ is a
right ideal.  Since we will eventually prove that $J'''=J'',$ we skip the proof that $J'''$ is an ideal for now.

Observing that $J\cap B_1$ is contained in $J''$, and since we now know that $J''$ is an ideal, the ideal generated by
$J \cap B_1$, namely $J'$, is also contained in $J''$.

Given any $x\in B$, notice that by Parseval's identity \cite {ParsevalId}, we have that $$ \soma {g\in G}F_g(z)^*F_g(z)
= F_1(z^*z).  $$ In fact \cite {ParsevalId} refers to the $E_g$, but since $F_g=E_g\circ \psi $ (see the proof of \cite
{Ft}), our identity follows easily from \cite {ParsevalId}.  Since ideals are hereditary, we then have that $$
F_1(z^*z)\in J \iff F_g(z)^*F_g(z)\in J,\ \forall g\in G, $$ and we notice that the condition in the right-hand side
above is also equivalent to $F_g(z)\in J$.  This proves that $J''=J'''$.  \endProof

Having seen how the ideals defined in \cite {ThreeIdeals} relate to each other, let us also discuss how do they relate
to $J$, itself. It is elementary to see that $J$ always contains $J'$, but the relationship between $J$ and $J''$ is not
straightforward.  We will see below that $J'=J''$ under certain conditions, in which case it will follow that
$J''\subseteq J$.  However there are examples in which $J''$ is not a subset of $J$, and in fact it may occur that, on
the contrary, $J$ is a proper subset of $J''$.

This is the case, for example, if $B$ is the full group C*-algebra of a non-amenable group $G$, and $J=\{0\}$.  One may
then prove that $J''$ is the kernel of the regular representation, hence $J''$ is strictly larger than $J$.

\def \tensorG #1{#1\tmin \Cr G} \def \tmap #1{#1\otimes 1}

One might suspect that the culprit for this anomaly is the failure of faithfulness of the standard conditional
expectation on $C^*(G)$, but examples may also be found in topologically graded C*-algebras with faithful conditional
expectations.  Take, for example, a group $G$ and a C*-algebra $B$.  One may then prove that $\tensorG B$ is isomorphic
to the reduced crossed product of $B$ by $G$ under the trivial action, so the former is a topologically $G$-graded
C*-algebra with faithful conditional expectation by \cite {EOneIsFaith}.

Assuming that $G$ is a non-exact group, one may find a short exact sequence of C*-algebras $$ 0 \to J \>i B \>\pi A \to
0, $$ for which $$ 0 \to \tensorG J \>{\tmap i} \tensorG B \>{\tmap \pi } \tensorG A \to 0, $$ is not exact.  It is well
known that the only place where exactness may fail is at the mid point of this sequence, meaning that the range of
$\tmap i$ is properly contained in the kernel of $\tmap \pi $.  Letting $J$ be the range of $\tmap i$, one may prove
that $J''$ is the kernel of $\tmap \pi $, whence $J$ is properly contained in $J''$, as claimed.  Consequently we have
that $$ J' \subseteq J \subsetneq J'', \equationmark FourierNorIndu $$ so this also produces an example in which $J''$
is strictly larger than $J'$.

\state Proposition \label InducedII Let $B$ be a topologically graded C*-algebra with grading $\{B_g\}_{g\in G}$, and
assume that the associated Fell bundle has the approximation property.  Then, for every ideal $J\ideal B$, one has that
the ideals $J'$ and $J''$ defined in \cite {ThreeIdeals} are equal.

\Proof It clearly suffices to prove that $J''\subseteq J'$.  Given $x\in J''$, we have by definition that each
$F_g(x)\in J$, and we claim that $F_g(x)\in J'$.  In order to see this, notice that $$ F_g(x)^*F_g(x)\in J\cap B_1
\subseteq J', $$ so we have that $ F_g(x)^*F_g(x) \equiv 0 \hbox { (mod $J'$)}.  $ Since $B/J'$ is a C*-algebra, we have
that $ F_g(x) \equiv 0 \hbox { (mod $J'$)}, $ as well, meaning that $F_g(x)\in J'$, thus proving our claim.

Let $\{a_i\}_i$ be a {\aproxNet } for $\Bun $, and let $\{S_i\}_i$ be the net of summation processes provided by \cite
{Cesaro}.  A glimpse at the formula defining $S_i$ is enough to convince ourselves that $S_i(x)$ is also in $J'$, hence
also $$ x= \lim _iS_i(x) \in J'.  \omitDoubleDollar \endProof

There is another situation in which we may guarantee the coincidence of the ideals $J'$ and $J''$.

\state Theorem \label InducedIII Let $G$ be a discrete group and let $B$ be a topologically $G$-graded C*-algebra.
Suppose that $G$ is exact and that the standard conditional expectation $F:B\to B_1$ is faithful.  Then for every ideal
$J$ of $B$ one has that the ideals $J'$ and $J''$ defined in \cite {ThreeIdeals} coincide.

\Proof Denote by $\Bun = \{B_g\}_{g\in G}$ the underlying Fell bundle and note that $B$ is isomorphic to $\CrB $ by
\cite {Charact}.  For each $g$ in $G$, let $J_g=J\cap B_g$ so that $\Jbun :=\{J_g\}_{g\in G}$ is an ideal in $\Bun $.
Employing \cite {ReduExactSeq} we have that the sequence $$ 0 \to \Cr {\Jbun } \> {\iota \lred } B \> {q \lred } \Cr
{\Bun /\Jbun } \to 0 \equationmark SeExataComB $$ is exact.  We next claim that $$ J' = \iota \lred \big (\Cr {\Jbun
}\big ) \and J'' = \Ker (q \lred ).  \equationmark TwoIdealIdentities $$ Given $g$ in $G$, notice that $$ J_g^*J_g
\subseteq B_1\cap J \subseteq J', $$ so by \cite {FellIsTRO} one has that $$ J_g = \clspan {J_g J_g^* J_g} \subseteq
\clspan {J_gJ'} \subseteq J', $$ so $\iota \lred \big (\Cr {\Jbun }\big ) \subseteq J'$.  Since the reverse inclusion is
evident, we have proven the first identity in \cite {TwoIdealIdentities}.

On the other hand, denoting by $E$ the faithful standard conditional expectation of $\Cr {\Bun /\Jbun }$, it is easy to
see that $ E\circ q \lred = q \lred \circ F, $ so for any $b$ in $B$ we have that $$ q \lred (b) = 0 \iff E\big (q \lred
(b^*b)\big ) = 0 \iff q \lred \big (F(b^*b)\big ) = 0 \iff F(b^*b)\in J_1, $$ where the last step is a consequence of
the fact that $F(b^*b)$ is in $B_1$, and that the behavior of $q$ on $B_1$ is that given by \cite {FunctorialityRed}.
This shows that $\Ker (q \lred ) = J''$, concluding the verification of \cite {TwoIdealIdentities}.

Since the sequence \cite {SeExataComB} is exact, the proof follows.  \endProof

\definition \label InduFourier Let $B$ be a topologically graded C*-algebra with grading $\{B_g\}_{g\in G}$ and Fourier
coefficient operators $F_g$.  We will say that a closed, two-sided ideal $J\ideal B$ is a \subj {Fourier ideal}, if
$F_g(J)\subseteq J$, for every $g$ in $G$.

Thus $J$ is a Fourier ideal if and only if $J\subseteq J''$, while $J$ is induced if and only if $J=J'$.  We may thus
reinterpret \cite {InducedI}, \cite {InducedII} and \cite {InducedIII} as follows.

\state Proposition \label InduIsFourier Let $B$ be a topologically graded C*-algebra.  Then every induced ideal of $B$
is a Fourier ideal.  Moreover, the converse holds if either \izitem \zitem the associated Fell bundle has the
approximation property, or \zitem $G$ is exact and the standard conditional expectation of $B$ onto $B_1$ is faithful.

For an example of a Fourier ideal which is not induced, see \cite {FourierNorIndu}.

\medskip The next result is stated for Fourier ideals but, because of the reasoning above, it also holds for induced
ones.

\state Proposition \label Ideals Let $B$ be a topologically graded C*-algebra with grading $\{B_g\}_{g\in G}$, and let
$J$ be a Fourier ideal of $B$. Then $B/J$ is topologically graded by the spaces $q(B_g)$, where $q$ is the quotient map.

\Proof Since $J$ is invariant under each Fourier coefficient operator $F_g$, we have that $F_g$ passes to the quotient
giving a well defined bounded map on $B/J$, namely $$ \tilde F_g (x+J) = F_g(x) + J \for x\in B.  $$

Notice that $$ q(B_g)= \tilde F_g(B/J) = \Ker (id - \tilde F_g), $$ the last step holding thanks to the fact that
$\tilde F_g$ is idempotent.  As a consequence we deduce that $q(B_g)$ is a closed subspace of $B/J$.

It is now immediate to verify that the collection $\{ q(B_g) \}_{g\in \G }$ satisfies \cite {Economical/i--iii}, and
that $\tilde F_1$ fills in the rest of the hypothesis there to allow us to conclude that this is in fact a topological
grading for $B/J$.  \endProof

All of the above results have their versions in the setting of partial crossed product algebras, since these are graded
algebras.  The next simple result, which we will use later, has no counterpart for graded algebras since its conclusion
explicitly mentions the partial action.

\state Proposition \label IntersInvar Let $ \Th = \big (A,G,\{\D g\}_{g\in G},\{\th g\}_{g\in G}\big ) $ be a
C*-algebraic partial dynamical system.  Given an ideal $J$ of either $A\rt \Th G$ or $A\redrt \Th G$, let $K=J\cap A$.
Then $K$ is a $\Th $-invariant ideal of $A$.

\Proof Of course we are identifying $A$ with its copy $A\delta _1$ in the crossed product algebra.  Regardless of
whether we are working with the full or reduced crossed product, the proof is the same: given $a$ in $K\cap \Di g$,
choose an approximate identity $\{v_i\}_i$ for $\D g$.  Then $$ J \ni (v_i\delta _g)(a\delta _1)(v_i\delta _g)^*
\={LotsAFormulas} v_i\th g(a)v_i^*\delta _1 \convrg {i}\th g(a)\delta _1, $$ so $\th g(a)$ is in $K$, proving the
statement.  \endProof

\nrem The motivation for this chapter comes from Nica's work on induced ideals of algebras of Wiener-Hopf operators \ref
{Section 6/Nica/1992}, which in turn is inspired by Str\u atil\u a and Voiculescu's work on AF-algebras \ref
{StratilaVoiculescu/1975}.  Propositions \cite {InducedI} and \cite {InducedII} have been proven in \ref {Exel/1997b},
while Theorem \cite {InducedIII} is from \ref {Theorem 5.1/Exel/2002p}\fn {Please note that \ref {Exel/2002p} is the
preprint version of \ref {Exel/2002}.}.

\chapter \PreFellBun dles

\def \half {^{1\over 2}}

In this short chapter we will develop some algebraic aspects of the theory of Fell bundles.  We will begin by
introducing the notion of an \"{algebraic} Fell bundle and the main question we shall analyze is whether or not they
admit a norm with which one may obtain a classical Fell bundle.

The motivation for this is the study of tensor products of Fell bundles by C*-algebras which we will do in the next
chapter.

\definition An \subj {algebraic Fell bundle} over a group $G$ is a collection $$ \Cun = \{C_g\}_{g\in G} $$ of complex
vector spaces, such that the disjoint union of all the $C_g$'s, which we also denote by $\Cun $, by abuse of language,
is moreover equipped with a multiplication operation and an involution $$ \cdot :\Cun \times \Cun \to \Cun , \qquad
\qquad * : \Cun \to \Cun , $$ satisfying the following properties for all $g$ and $h$ in $G$, and all $b$ and $c$ in
$\Cun $: \iaitem \aitem $C_gC_h\subseteq C_{gh}$, \aitem multiplication is bi-linear from $C_g\times C_h$ to $C_{gh}$,
\aitem multiplication on $\Cun \,$ is associative, \advance \aitemno by 1 \medskip \aitem $(C_g)^* \subseteq C_{g\inv
}$, \aitem involution is conjugate-linear from $C_g$ to $C_{g\inv }$, \aitem $(bc)^* = c^* b^*$, \aitem $b^{**}=b$.

Compared to \cite {DefineFellBundle} observe that axioms (d) and (i--k) are missing since these refer to norms which are
absent in the present case.  Our goal here is to furnish a norm on each $C_g$ with respect to which the completions form
a Fell bundle.

When $G=\{1\}$, an algebraic Fell bundle consists of a single *-algebra and the reader is probably aware that providing
a C*-norm on a *-algebra is not a problem with a straightforward solution.  In the case of a general group $G$, the
above axioms imply that $C_1$ is a *-algebra so we should therefore expect our construction of norms on the $C_g$ to
require at least an initial choice of a C*-norm on $C_1$.

If we are indeed given a C*-norm on $C_1$, we may consider the completion $$ B_1 := \overline {C}_1 $$ with respect to
this norm, which will therefore be a C*-algebra.  Given $c$ in any $C_g$, we may view $c^*c$ as an element of $B_1$,
however without further hypothesis there is no reason why $c^*c$ is positive in $B_1$.  Without this positivity
condition there is clearly no hope of turning our algebraic Fell bundle into a Fell bundle.

Still assuming that $c\in C_g$, consider the mapping $$ a\in C_1 \mapsto c^*ac \in C_1.  $$ Again there seems to be no
reason why this is continuous with respect to the given norm on $C_1$ and again this continuity is a necessary condition
for us to proceed.

\definition \label DefinePreFell A \subj {\preFellBun dle} over a group $G$ is an algebraic Fell bundle $$ \Cun =
\{C_g\}_{g\in G} $$ equipped with a C*-norm $\Vert \ponto \Vert $ on $C_1$ such that, for every $c$ in any $C_g$, one
has that \izitem \zitem $c^*c$ is a positive element in the C*-algebra $B_1$ obtained by completing $C_1$ relative to
the norm given above, \zitem the mapping $$ Ad_c:a\in C_1 \mapsto c^*ac\in C_1, $$ is continuous with respect to $\Vert
\ponto \Vert $.

\medskip Should one prefer to avoid any reference to the completed algebra $B_1$ in \cite {DefinePreFell/i}, above, one
could instead require that for every $c\in C_g$, there exists a sequence $\{a_n\}_n\subseteq C_1$, such that $\Vert
a_n^*a_n - c^*c\Vert \convrg i 0$.

\medskip \fix From now on we will fix a \preFellBun dle $\Cun = \{C_g\}_{g\in G}$, and we will look for suitable norms
$\Vert \ponto \Vert _g$ on the other $C_g$'s with the intent of obtaining a Fell bundle.  In view of \cite
{DefineFellBundle/j}, we have no choice but to define $$ \Vert c\Vert _g = \Vert c^*c\Vert \half \for c\in C_g, $$ and
our task is then to check the remaining axioms.  We begin with some technical results.

\state Proposition \label LemmaStarInvNorm For every $c$ in any $C_g$, one has that $\Vert cc^*\Vert = \Vert c^*c\Vert
$.

\Proof Recall that a self-adjoint element $x$ in a C*-algebra satisfies $\Vert x^n\Vert = \Vert x\Vert ^n$, for every
$n\in {\bf N}$.  This identity also holds for our C*-norm on $C_1$, as it follows from the corresponding fact applied to
the completed algebra $B_1$.  So $$ \Vert cc^*\Vert ^{n+1} = \Vert (cc^*)^{n+1}\Vert = \Vert c(c^*c)^nc^*\Vert = \Vert
Ad_c\big ((c^*c)^n\big )\Vert \$\leq \Vert Ad_c\Vert \Vert (c^*c)^n\Vert \leq \Vert Ad_c\Vert \Vert c^*c\Vert ^n.  $$
Taking $n^{th}$ root leads to $$ \Vert cc^*\Vert ^{n+1\over n} \leq \Vert Ad_c\Vert ^{1\over n}\Vert c^*c\Vert , $$ and
when $n\to \infty $, we get $$ \Vert cc^*\Vert \leq \Vert c^*c\Vert .  $$ The reverse inequality follows by replacing
$c$ with $c^*$.  \endProof

Let us now prove a technical result, reminiscent of \cite {ApproximateIdentity}.

\state Lemma \label CVCapprox Let $\{v_i\}_i$ be an approximate identity for $C_1$.  Then, for every $c$ in any $C_g$,
one has that $$ \lim _ic^*v_ic = c^*c.  $$

\Proof Using that $\Vert \ponto \Vert $ is a C*-norm, we have $$ \Vert c^*v_ic - c^*c\Vert ^2 = \Vert (c^*v_ic -
c^*c)^*(c^*v_ic - c^*c)\Vert \$= \Vert c^*v_i^*cc^*v_ic - c^*v_i^*cc^*c - c^*cc^*v_ic + c^*cc^*c\Vert \$= \Vert
Ad_c(v_i^*cc^*v_i - v_i^*cc^* - cc^*v_i + cc^*)\Vert \convrg i 0.  \omitDoubleDollar \endProof

With this we may prove a version of \cite {OrderNorm}.

\state Lemma \label OrderNormPre \def \pos {B_1\null _+}\null Given $c$ in any $C_g$, one has that $$ c^*xc\leq \Vert
x\Vert c^*c \for x\in C_1\cap \pos , $$ where $\pos $ denotes the set of positive elements of the C*-algebra $B_1$
obtained by completing $C_1$, and the order relation ``$\leq $'' is that of $B_1$.

\Proof \def \pos {B_1\null _+}\null Given $c$ in $C_g$, we claim that $$ x\in C_1\cap \pos \imply c^*xc\in \pos .
\equationmark ConjugaPosit $$ In fact, for each $x$ in $C_1\cap \pos $, write $x=y^*y$, with $y\in B_1$, and choose a
sequence $\{z_n\}_n\subseteq C_1$, converging to $y$.  Then $$ c^*xc = Ad_c\big (\lim _{n\to \infty } z_n^*z_n\big )
\={DefinePreFell/ii} \lim _{n\to \infty } Ad_c(z_n^*z_n) = \lim _{n\to \infty } c^*z_n^*z_nc \_\in {DefinePreFell/i}\pos
, $$ proving \cite {ConjugaPosit}.  Next choose an approximate identity $\{v_i\}_i$ for $B_1$, contained in $C_1$. Then,
for $x\in C_1\cap \pos $, and for any $i$ we have that $$ v_i^*(\Vert x\Vert - x)v_i \in C_1\cap \pos , $$ whence by
\cite {ConjugaPosit} we have $$ \pos \ni c^*v_i^*(\Vert x\Vert - x)v_ic = \Vert x\Vert c^*v_i^*v_ic - c^*v_i^*xv_ic.
\equationmark SomePositiveSeq $$ Observing that $\{v_i^*v_i\}_i$ is also an approximate identity for $B_1$, we deduce
from \cite {CVCapprox} that $$ c^*v_i^*v_ic \convrg i c^*c, $$ and from \cite {DefinePreFell/ii} we have $$
c^*v_i^*xv_ic\convrg i c^*xc.  $$

Taking the limit as $i\to \infty $ in \cite {SomePositiveSeq} we finally obtain $$ \Vert x\Vert c^*c - c^*xc\in \pos .
\omitDoubleDollar \endProof

We now have all of the necessary tools in order to prove the main result of this chapter:

\state Theorem \label NormsOnTPBundle Given a \preFellBun dle $\Cun = \{C_g\}_{g\in G}$, there is a unique family of
seminorms $\Vert \ponto \Vert _g$ on the $C_g$'s, such that $\Vert \ponto \Vert _1$ is the given norm on $C_1$, and for
$g,h\in G$, $b\in C_g$, and $c\in C_h$, one has that \smallskip \aitemno = 2 \aitem $\Vert bc\Vert _{gh} \leq \Vert
b\Vert _g \Vert c\Vert _h$, \aitemno = 7 \aitem $\Vert b^*\Vert _{g\inv } = \Vert b\Vert _g$, \aitem $\Vert b^*b\Vert _1
= \Vert b\Vert _g^2$.

\Proof Define $\Vert \ponto \Vert _1$ to coincide with the given norm on $C_1$, and for all $g\neq 1$, and all $b\in
C_g$, set $$ \Vert b\Vert _g = \Vert b^*b\Vert \half .  $$ If $b\in C_g$ and $c\in C_h$, we have $$ (bc)^*bc = c^*b^*bc
\_\leq {OrderNormPre} \Vert b^*b\Vert c^*c, $$ whence $$ \Vert bc\Vert _{gh}^2 = \Vert (bc)^*bc\Vert \leq \Vert
b^*b\Vert \Vert c^*c\Vert = \Vert b\Vert ^2_g\Vert c\Vert ^2_h, $$ proving (d). Notice that (i) follows from \cite
{LemmaStarInvNorm}, while (j) follows by definition when $g\neq 1$, and otherwise from the fact that the norm on $C_1$
is assumed to be a C*-norm.

To conclude we prove the triangle inequality: given $b,c\in C_g$, we have $$ \Vert b+c\Vert ^2_g = \Vert
(b+c)^*(b+c)\Vert = \Vert b^*b+b^*c+c^*b+c^*c\Vert \explica \leq {d} $$ $$ \leq \Vert b^*\Vert _{g\inv }\Vert b\Vert
_g+\Vert b^*\Vert _{g\inv }\Vert c\Vert _g+\Vert c^*\Vert _{g\inv }\Vert b\Vert _g+\Vert c^*\Vert _{g\inv }\Vert c\Vert
_g \explica ={i} $$$$ = \pilar {15pt} \Vert b\Vert _g\Vert b\Vert _g+\Vert b\Vert _g\Vert c\Vert _g+\Vert c\Vert _g\Vert
b\Vert _g+\Vert c\Vert _g\Vert c\Vert _g = (\Vert b\Vert _g+\Vert c\Vert _g)^2.  $$ This concludes the proof of the
existence part, while uniqueness follows easily from (j).  \endProof

Given a \preFellBun dle $\Cun = \{C_g\}_{g\in G}$, one may then consider the completion of each $C_g$ under $\Vert
\ponto \Vert _g$, say $$ B_g = \overline {C}_g, $$ and, after extending the operations by continuity, the verification
of axioms \cite {DefineFellBundle/i--k} becomes routine, especially since they are already known to hold on dense sets.
We thus obtain a Fell bundle $\Bun = \{B_g\}_{g\in G}$.

\definition \label DefinFellCompletion The Fell bundle $\Bun $, obtained as above from a given \preFellBun dle $\Cun $,
will be called the \subjex {completion}{completion of a \preFellBun dle} of $\Cun $.

The question of extending representations from a \preFellBun dle to its completion is an important one which we will
discuss next.

\state Proposition \label ExtendToFellComplet Let $\Cun = \{C_g\}_{g\in G}$ be a \preFellBun dle and let $A$ be a
C*-algebra.  Suppose we are given a collection of linear maps $\pi =\{\pi _g\}_{g\in G}$, where $$ \pi _g: C_g \to A, $$
such that \izitem \zitem $\pi _g(b)\pi _h(c) = \pi _{gh}(bc)$, \zitem $\pi _g(b)^* = \pi _{g\inv }(b^*)$, \medskip
\noindent for all $g,h\in G$, and all $b\in C_g$, and $c\in C_h$.  Suppose, moreover, that $\pi _1$ is continuous
relative to the norm on $C_1$.  Then each $\pi _g$ is continuous relative to the norm $\Vert \ponto \Vert _g$ on $C_g$,
and hence it extends to a bounded linear map $$ \tilde \pi _g : B_g \to A, $$ where $\Bun = \{B_g\}_{g\in G}$ is the
Fell bundle completion of\/ $\Cun $.  In addition the collection of maps $\tilde \pi =\{\tilde \pi _g\}_{g\in G}$ is a
representation of $\Bun $ in $A$.

\Proof For $c$ in any $C_g$ we have $$ \Vert \pi _g(c)\Vert ^2 = \Vert \pi _g(c)^*\pi _g(c)\Vert \explica ={i\&ii} \Vert
\pi _1(c^*c)\Vert \leq \Vert c^*c\Vert _1 = \Vert c\Vert ^2_g, $$ proving $\pi _g$ to be continuous.  The remaining
statements are verified in a routine way.  \endProof

\chapter Tensor products of Fell bundles

\def \tautp {\otimes _\tau }

Tensor products form a very important part of the theory of C*-algebras and, given the very close relationship between
C*-algebras and Fell bundles, no treatment of Fell bundles is complete without a careful study of their tensor products.

The most general form of this theory would start by considering two Fell bundles $\Bun $ and $\Bun '$, over two groups
$G$ and $G'$, respectively, and one would then attempt to construct a Fell bundle $\Bun \otimes \Bun '$ over the group
$G\times G'$.  However, given the applications we have in mind, we will restrict ourselves to the special case in which
$G'$ is a trivial group.  In other words, we will restrict our study to tensor products of Fell bundles by single
C*-algebras.

Among the main aspects of tensor products we plan to analyze is the relationship between the corresponding versions of
spatial and maximal norms.

\medskip

\fix Let us fix a C*-algebra $A$ and a Fell bundle $ \Bun = \{B_g\}_{g\in G}.  $ For each $g$ in $G$, let us consider
the vector space tensor product $$ C_g := A\odot B_g.  $$

Given $g$ and $h$ in $G$, it is easy to see that there exists a bi-linear operation $C_g \times C_h \to C_{gh}$, such
that $$ (a\otimes b)(a'\otimes b') = (aa')\otimes (bb') \for a,a'\in A \for b\in B_g \for b'\in B_h.  $$ Likewise one
may show the existence of a conjugate-linear map $^*$ from $C_g$ to $C_{g\inv }$, such that $$ (a\otimes b)^* =
a^*\otimes b^* \for a\in A \for b\in B_g.  $$

One then easily checks that the collection $$ A \odot \Bun = \{C_g\}_{g\in G} \equationmark DefineTPAlgFb $$ is an
algebraic Fell bundle with the above operations.

The next result is intended to aid the verification of \cite {DefinePreFell/i}, once we have a C*-norm on $C_1$.

\state Proposition Given $c$ in any $C_g$, there exist $x_1,\ldots ,x_n\in C_1$, such that $$ c^*c=\somain x_i^*x_i.  $$

\Proof Writing $c=\somain a_i\otimes b_i$, with $a_i\in A$, and $b_i\in B_g$, we claim that the $n\times n$ matrix $$ m
= \{b_i^*b_j\}_{i,j} \in M_n(B_1), $$ is positive.  To prove it notice that, viewed as a matrix over the cross sectional
C*-algebra of $\Bun $, we have that $m=y^*y$, where $y$ is the row matrix $$ y = [b_1\ b_2\ \ldots \ b_n].  $$

Thus $m$ is positive as a matrix over $C^*(\Bun )$, but since $B_1$ is a subalgebra of $C^*(\Bun )$ by \cite
{FinallyInjective/iv}, and since the coefficients of $m$ lie in $B_1$, we have that $m$ is positive as a matrix over
$B_1$.  So there exists $h$ in $M_n(B_1)$, such that $m = h^*h$, which translates into $$ b_i^*b_j = \usoma {k=1}n
(h^*)_{i,k}h_{k,j} = \usoma {k=1}n (h_{k,i})^*h_{k,j} \for i,j=1,\ldots ,n.  $$ We then have $$ c^*c = \Big (\usoma
{i=1}na_i\otimes b_i\Big )^* \Big (\usoma {j=1}na_j\otimes b_j\Big ) = \usoma {i,j=1}n a_i^*a_j\otimes b_i^*b_j \$=
\usoma {i,j,k=1}n a_i^*a_j\otimes (h_{ki})^*h_{kj} = \usoma {k=1}n \ \usoma {i,j=1}n (a_i\otimes h_{ki})^*(a_j\otimes
h_{kj}) \$= \usoma {k=1}n \Big (\usoma {i=1}n a_i\otimes h_{ki}\Big )^* \Big (\usoma {j=1}n a_j\otimes h_{kj}\Big ) =
\usoma {k=1}n x_k^*x_k, $$ where $x_k = \usoma {i=1}n a_i\otimes h_{ki}$.  \endProof

We then see that axiom \cite {DefinePreFell/i} will hold for any choice of C*-norm on $C_1$.

\state Proposition \label MaxMinPreFell Let \ $\Vert \ponto \Vert \lmax $ \ and \ $\Vert \ponto \Vert \lmin $ \ be the
maximal and minimal C*-norms on the algebraic tensor product $A\odot B_1$, respectively.  Then $A\odot \Bun $ is a
\preFellBun dle with either one of these norms.

\Proof After the remark in the paragraph just before the statement, it is now enough to prove that, for every $a$ and
$a'$ in $A$, and every $b$ and $b'$ in any $B_g$, the mapping $$ \rho \ : \ t\in C_1 \mapsto (a\otimes b)^*t(a'\otimes
b')\in C_1, $$ is continuous with respect to both the maximal and minimal C*-norms.  If $x\in A$, and $y\in B_1$, notice
that $$ \rho (x\otimes y) = a^*xa' \otimes b^*yb', $$ thus $\rho $ is seen to be the tensor product of the maps $\varphi
$ and $\psi $ given by $$ \varphi :x\in A \mapsto a^*xa'\in A \and \psi :y\in B_1 \mapsto b^*yb'\in B_1.  $$

By the polarization identity we may write $\varphi $ as a linear combination of four maps of the form $$ \varphi _i:x\in
A \mapsto a_i^*xa_i\in A, $$ with $a_i\in A$, and likewise $\psi $ may be written as a linear combination of four maps
of the form $$ \psi _j:y\in B_1 \mapsto b_j^*yb_j\in B_1, $$ with $b_j\in B_g$.  Consequently $\rho $ may be written as
a linear combination of the sixteen maps $\varphi _i\odot \psi _j$, and it is therefore enough to show that these are
continuous.  Since both the $\varphi _i$ and the $\psi _j$ are completely positive maps, the result follows from \ref
{Theorem 3.5.3/BrownOzawa/2008}.  \endProof

\definition \label DefineFellTP Given a C*-algebra $A$ and a Fell bundle $\Bun = \{B_g\}_{g\in G}$, let $\tau $ be any
C*-norm on $A\odot B_1$, with respect to which $A\odot \Bun $ is a \preFellBun dle.  Then its completion, according to
\cite {DefinFellCompletion}, will be denoted by $A\tautp \Bun $.  For each $g$ in $G$, we will denote the corresponding
fiber of $A\tautp \Bun $ by $A\tautp B_g$.  In the case of the maximal and minimal C*-norms on $A\odot B_1$, we will
respectively denote the corresponding completions by $$ A\tmax \Bun \and A\tmin \Bun , $$ with fibers $$ A\tmax B_g \and
A\tmin B_g, $$ for each $g$ in $G$.

Employing similar methods one could also study the tensor product of a Fell bundle $ \Bun = \{B_g\}_{g\in G} $ over a
group $G$, by a Fell bundle $ \Cun = \{C_h\}_{h\in H} $ over a group $H$, obtaining a Fell bundle $$ \Bun \otimes \Cun =
\{\overline {\pilar {9pt}B_g\otimes C_h}\}_{(g,h)\in G\times H} $$ over the group $G\times H$.  This will of course
require an appropriate choice of norm on $B_1\odot C_1$.  We will however not pursue these ideas here.

\state Proposition \label TensorWitMap Given a Fell bundle $ \Bun = \{B_g\}_{g\in G}, $ let $C$ and $D$ be C*-algebras
and suppose that $\tau $ and $\nu $ are C*-norms on $C\odot B_1$ and $D\odot B_1$, respectively, making $C\odot \Bun $
and $D\odot \Bun $ into \preFellBun dles.  Suppose moreover that $\varphi :C\to D$ is a *-homomorphism such that $$
\varphi \otimes id: C\odot B_1 \to D\odot B_1 $$ is continuous relative to $\tau $ and $\nu $.  Then there exists a
morphism $\tilde \varphi = \{\varphi _g\}_{g\in G}$ from $C\otimes _\tau \Bun $ to $D\otimes _\nu \Bun $, such that $$
\varphi _g(c\otimes b) = \varphi (c)\otimes b, $$ for every $b$ in any $B_g$.

\Proof Given $g$ in $G$, define $$ \varphi ^0_g : C\odot B_g \to D\odot B_g $$ by $\varphi ^0_g = \varphi \otimes id$.
It is then easy to see that the collection of maps $$ \widetilde {\varphi ^0} = \{\varphi ^0_g\}_{g\in G} $$ satisfy
\cite {DefineMorphBun/i\&ii}, relative to the algebraic Fell bundle structures of $C\odot B_g$ and $D\odot B_g$ given in
\cite {DefineTPAlgFb}.  We moreover claim that each $\varphi ^0_g$ is continuous relative to the norms $\Vert \ponto
\Vert ^\tau _g$ and $\Vert \ponto \Vert ^\nu _g$, given by \cite {NormsOnTPBundle} on $C\odot B_g$ and $D\odot B_g$,
respectively.

In order to verify the claim, notice that the case $g=1$ is granted by hypothesis.  For an arbitrary $g$ in $G$, pick
any $x$ in $C\odot B_g$, and observe that $$ \big (\Vert \varphi ^0_g(x)\Vert _g^\nu \big )^2 = \Vert \varphi
^0_g(x)^*\varphi ^0_g(x)\Vert _1^\nu = \Vert \varphi ^0_1(x^*x)\Vert _1^\nu \leq \Vert x^*x\Vert _1^\tau = \big (\Vert
x\Vert _g^\tau \big )^2, $$ where the inequality above is a consequence of the already verified case $g=1$.  We may
therefore extend each $\varphi ^0_g$ to a continuous mapping $$ \varphi _g : C\otimes _\tau B_g \to D\otimes _\nu B_g,
$$ and the reader may then easily verify that the resulting collection of maps $\tilde \varphi = \{\varphi _g\}_{g\in
G}$ satisfies the required conditions.  \endProof

Recall that any Fell bundle admits a \"{universal} representation in its full cross sectional C*-algebra by \cite
{DefineHatJota}, as well a \"{regular} representation in its reduced cross sectional C*-algebra by \cite
{DefineRegularRep}.  In the case of $A\tmax \Bun $ and $A\tmin \Bun $, we believe it is best not to introduce any
special notation for these representations, instead relying on \cite {FinallyInjective/iv--v}, which allows us to
identify the fibers of our bundles as subspaces of both the full and reduced cross sectional C*-algebras.

In principle this has a high risk of confusion, since when $a$ is in $A$ and $b$ is in some $B_g$, the expression
$a\otimes b$ may be interpreted as elements of eleven different spaces, namely $A\odot B_g$, $A\tmax B_g$, $A\tmin B_g$,
$C^*(A\tmax \Bun )$, $C^*(A\tmin \Bun )$, $C^*\lred (A\tmax \Bun )$, $C\lred ^*(A\tmin \Bun )$, besides the maximal and
minimal tensor products of $A$ by either $C^*(\Bun )$ or $\CrB $.  Fortunately, as we will see, the context will always
suffice to determine the correct interpretation of $a\otimes b$.

\state Proposition \label CrossNorms For every $a$ in $A$, and $b$ in any $B_g$, one has that $$ \Vert a\otimes b\Vert
=\Vert a\Vert \Vert b\Vert , $$ for any one of the various interpretations of $a\otimes b$, except, of course, for
$A\odot B_g$ which is a space devoid of a norm.

\Proof We begin by treating the interpretation of $a\otimes b$ in $A\tmax B_g$.  By \cite {NormsOnTPBundle/j} we have $$
\Vert a\otimes b\Vert ^2_g = \Vert (a\otimes b)^*(a\otimes b)\Vert _1 = \Vert a^*a\otimes b^*b\Vert _1 = \Vert
a^*a\otimes b^*b\Vert \lmax \$= \Vert a^*a\Vert \Vert b^*b\Vert = \Vert a\Vert ^2\Vert b\Vert ^2, $$ where, in the
penultimate step above, we have used that the norm on $A\tmax B_1$ is a cross-norm (i.e.,~satisfies the identity in the
statement).

A very similar argument proves the result in the case of $A\tmin B_g$, and then \cite {FinallyInjective/iv} takes care
of $C^*(A\tmax \Bun )$ and $C^*(A\tmin \Bun )$, while \cite {FinallyInjective/v} does it for $C^*\lred (A\tmax \Bun )$
and $C\lred ^*(A\tmin \Bun )$.

Again by \cite {FinallyInjective/iv--v}, we know that the norm of $b$, as interpreted within $C^*(\Bun )$ or $\CrB $,
coincide with its norm as an element of $B_g$.  Thus, the result relative to the maximal or minimal tensor products of
$A$ by either $C^*(\Bun )$ or $\CrB $ follow, since both the maximal or minimal norms are cross-norms.  \endProof

\state Theorem \label TensorMaxFull Given a C*-algebra $A$ and a Fell bundle $\Bun $, one has that $$ C^*(A\tmax \Bun )
\simeq A\tmax C^*(\Bun ), $$ via an isomorphism which sends $a\otimes b \to a\otimes \hj _g(b)$, for $b$ in any $B_g$,
and all $a$ in $A$.

\Proof For each $g$ in $G$, define $$ \pi _g: A\odot B_g \to A\tmax C^*(\Bun ) $$ by $$ \pi _g(a\otimes b) = a\otimes
\hj _g(b) \for a\in A \for b\in B_g, $$ where $\hj $ is the universal representation of $\Bun $ in $C^*(\Bun )$.

It is clear that the $\pi _g$ satisfy \cite {ExtendToFellComplet/i--ii}.  In addition, by the universal property of the
maximal norm \ref {Theorem 3.3.7/BrownOzawa/2008}, we have that $\pi _1$ is continuous for $\Vert \ponto \Vert \lmax $.
It then follows from \cite {ExtendToFellComplet} that the $\pi _g$ extend to a representation $\tilde \pi =\{\tilde \pi
_g\}_{g\in G}$ of $A\tmax \Bun $ \/ in $A\tmax C^*(\Bun )$.  We then conclude from \cite {FromRepBunToRelAlg} that there
exists a unique *-homomorphism $$ \varphi : C^*(A\tmax \Bun ) \to A\tmax C^*(\Bun ) \equationmark FromCrossToTensor $$
such that $\varphi $ composed with the universal representation of $A\tmax \Bun $ in $C^*(A\tmax \Bun )$ yields $\tilde
\pi $.  Employing the identifications discussed in the paragraph before \cite {CrossNorms}, this means that $$ \varphi
(a\otimes b) = a\otimes \hj _g(b) \for a\in A \for g\in G \for b\in B_g.  $$

In order to prove that $\varphi $ is an isomorphism we will construct an inverse for it.  This will involve
understanding a representation of $\Bun $ in the multiplier algebra of $C^*(A\tmax \Bun )$, which we will now describe.
Given $b$ in any $B_g$, we claim that there are bounded linear operators $$ L_b, R_b: C^*(A\tmax \Bun ) \to C^*(A\tmax
\Bun ), $$ such that $$ L_b(a\otimes c) = a\otimes bc \and R_b(a\otimes c) = a\otimes cb, $$ for any $a$ in $A$, and any
$c$ in the total space of $\Bun $.

To see this let $\{v_i\}_i$ be an approximate identity for $A$, and consider the operators $L_{v_i\otimes b}$ on
$C^*(A\tmax \Bun )$ given by left-multiplication by $v_i\otimes b$.  We then have that $$ L_{v_i\otimes b}(a\otimes c) =
v_ia\otimes bc \convrg i a\otimes bc.  $$ Note that this convergence is guaranteed by \cite {CrossNorms}.

Therefore the net $\{L_{v_i\otimes b}\}_i$ converges pointwise on $C_c(A\tmax \Bun )$ and, being a uniformly bounded
net, it actually converges pointwise everywhere to a bounded operator which we denote by $L_b$, and which clearly
satisfies the required conditions.  The existence of $R_b$ is proved similarly, and it is then easy to see that the pair
$(L_b,R_b)$ is a multiplier of $C^*(A\tmax \Bun )$.  Moreover the mappings $$ \mu _g : \ b\ \in \ B_g\ \mapsto \
(L_b,R_b) \ \in \ \Mult \big (C^*(A\tmax \Bun )\big ), $$ once put together, form a representation $$ \mu =\{\mu
_g\}_{g\in G} $$ of $\Bun $ in $\Mult \big (C^*(A\tmax \Bun )\big )$.  Feeding this representation into \cite
{FromRepBunToRelAlg} provides a *-homomorphism $$ \psi :C^*(\Bun ) \to \Mult \big (C^*(A\tmax \Bun )\big ), $$ such that
$\psi \circ \hj _g=\mu _g$, for all $g$ in $G$.  In particular $$ \psi \big (\hj _g(b)\big )(a\otimes c) = a\otimes bc
\for b\in B_g \for c\in \Bun \for a\in A , $$ where $\Bun $ is being used here to denote total space.

We now prove the existence of *-homomorphism $$ \chi :A \to \Mult \big (C^*(A\tmax \Bun )\big ), $$ such that $$ \chi
(a)(a'\otimes c) = aa'\otimes c \for a,a'\in A \for c\in \Bun .  $$ Since the reasoning is similar to the above, we
limit ourselves to sketching it.  Fixing $a$ in $A$ and an approximate identity $\{v_i\}_i$ for $B_1$, one has that the
left-multiplication operators $L_{a\otimes v_i}$ satisfy $$ L_{a\otimes v_i}(a'\otimes c) = aa'\otimes v_ic \convrg i
aa'\otimes c \for a'\in A \for c\in \Bun , $$ by \cite {ApproximateIdentity}.  The strong limit of the $L_{a\otimes
v_i}$ therefore exists and gives the first coordinate of a multiplier pair $(L_a,R_a)$ which we set to be $\chi (a)$.

Once the existence of $\chi $ is established, observing that the ranges of $\psi $ and $\chi $ commute, we may employ
once again the universal property of the maximal norm \ref {Theorem 3.3.7/BrownOzawa/2008}, obtaining a *-homomorphism $
\chi \times \psi $ from $A\tmax C^*(\Bun )$ to the above multiplier algebra, such that $$ (\chi \times \psi )(a\otimes
y) = \chi (a)\psi (y) \for a\in A \for y\in C^*(\Bun ).  $$ In particular, for $a,a'\in A$, $b\in B_g$, and $c\in \Bun $
(total space), we have $$ \Big ((\chi \times \psi )\big (a\otimes \hj _g(b)\big )\Big )(a'\otimes c) = aa'\otimes bc.
$$ This implies that $$ (\chi \times \psi )\big (a\otimes \hj _g(b)\big ) = a\otimes b, $$ so the range of $\chi \times
\psi $ is actually contained in $C^*(A\tmax \Bun )$ (or rather, in its canonical copy inside the multiplier algebra).
It is now easy to see that $\chi \times \psi $ is the inverse of the map $\varphi $ in \cite {FromCrossToTensor}, so the
proof is concluded.  \endProof

The following is the reduced/minimal version of the result above.

\state Theorem \label TensorMinRed Given a C*-algebra $A$ and a Fell bundle $\Bun $, one has that $$ C^*\lred (A\tmin
\Bun ) \simeq A\tmin \CrB , $$ via an isomorphism sending $a\otimes b \to a\otimes \lambda _g(b)$, for $b$ in any $B_g$,
and all $a$ in $A$.

\Proof The first few steps of this proof are very similar to the proof of \cite {TensorMaxFull}, but significant
differences will appear along the way.  For each $g$ in $G$, define $$ \pi _g: A\odot B_g \to A\tmin \CrB $$ by $$ \pi
_g(a\otimes b) = a\otimes \lambda _g(b) \for a\in A \for b\in B_g, $$ where $\lambda $ is the regular representation of
$\Bun $ in $\CrB $ given in \cite {DefineRegularRep}.

It is clear that the $\pi _g$ satisfy \cite {ExtendToFellComplet/i--ii}.  In addition, by \ref {Proposition
3.6.1/BrownOzawa/2008}, we have that $\pi _1$ is isometric for the minimal norm on $A\odot B_1$, and hence also
continuous.  It then follows from \cite {ExtendToFellComplet} that the $\pi _g$ extend to a representation $\tilde \pi
=\{\tilde \pi _g\}_{g\in G}$ of $A\tmin \Bun $ in $A\tmin \CrB $.  Note that $\tilde \pi _1$ is then isometric on
$A\tmin B_1$.

Let us denote the integrated form of $\tilde \pi $ by $$ \rho : C^*(A\tmin \Bun ) \to A\tmin \CrB $$ so that by \cite
{FromRepBunToRelAlg} we have $$ \rho (a\otimes b) = a\otimes \lambda _g(b) \for a\in A \for b\in B_g.  $$

We next apply \cite {FellAbsorptionThree} for $\tilde \pi $, obtaining a *-homomorphism $$ \psi : \BigCr {A\tmin \Bun }
\to A\tmin \CrB \tmin \Cr G, $$ such that for all $a\otimes b$ in any $A\tmin B_g$, one has $$ \psi (a\otimes b) =
a\otimes \lambda _g(b)\otimes \lreg _g.  $$

Since $\tilde \pi _1$ is faithful, we have by \cite {FellAbsorptionThree} that $\psi $ is faithful.  Using the above
maps, in addition to the map $\sigma $ provided by \cite {MiniAbsorp}, we build the diagram $$ \matrix { C^*(A\tmin \Bun
) & \>{\ds \rho } & A\tmin \CrB \hfill \cr \cr \Lambda \ \Big \downarrow && \qquad \Big \downarrow \ id\otimes \sigma
\hfill \cr \cr \Cr {A\tmin \Bun } & \> {\ds \psi } & A\tmin \CrB \tmin \Cr G } $$ which the reader may easily prove to
be commutative.

Since $\sigma $ is faithful, we have by \ref {Proposition 3.6.1/BrownOzawa/2008}, that $id\otimes \sigma $ is faithful.
We have also seen above that $\psi $ is faithful, so one deduces that $\rho $ and $\Lambda $ must have the same null
spaces.  Consequently $\rho $ factors through the kernel of $\Lambda $, producing the required isomorphism.  \endProof

Let us now study the approximation property for tensor products of C*-algebras by Fell bundles.

\state Proposition \label TensorAproxProp Let $A$ be a C*-algebra and $\Bun = \{B_g\}_{g\in G}$ be a Fell bundle.  Also
let $\tau $ be any C*-norm on $A\odot B_1$, with respect to which $A\odot \Bun $ is a \preFellBun dle.  If $\Bun $
satisfies the approximation property described in \cite {AP}, then so does $A\tautp \Bun $.

\Proof Let $\{v_i\}_{i\in I}$ be an approximate identity for $A$, and choose a {\aproxNet } $\{ a_j \}_{j\in J}$ for
$\Bun $.  Considering $I\times J$ as an ordered set with coordinate-wise order relation, it is clear that $I\times J$ is
a directed set.  For each $(i,j)\in I\times J$, let $$ \alpha _{i,j} : G \to A\tautp B_1, $$ be defined by $$ \alpha
_{i,j}(g) = v_i\otimes a_j(g) \for g\in G.  $$

We then claim that $\alpha $ is a {\aproxNet } for $A\tautp \Bun $.  Indeed, given $(i,j)\in I\times J$, we have $$
\normsum {\soma {g\in \G } \alpha _i(g)^* \alpha _i(g)} = \normsum {\soma {g\in \G } v_i^*v_i\otimes a_i(g)^* a_i(g)}
\$= \Vert v_i^*v_i\Vert \,\normsum {\soma {g\in \G } a_i(g)^* a_i(g)}, $$ where we have denoted $\tau $ by the usual
norm symbol, noting that, as well as any C*-norm on $A\odot B_1$, $\tau $ is a cross norm \ref {Lemma
3.4.10/BrownOzawa/2008}.  This proves that the $\alpha _i$ satisfy \cite {AP/i}.

Given $a$ in $A$, and $b$ in any $B_g$, we have for all $(i,j)\in I\times J$, that $$ \soma {h\in \G } \alpha _i(gh)^*
(a\otimes b) \alpha _i(h) = \soma {h\in \G } \big (v_i\otimes a_j(gh)\big )^* (a\otimes b) \big (v_i\otimes a_j(h)\big )
\$= \soma {h\in \G } v_i^*av_i\otimes a_j(gh)^* b a_j(h) = v_i^*av_i\otimes \soma {h\in \G } a_j(gh)^* b a_j(h) \convrg
{i,j} a\otimes b.  $$ The conclusion now follows from \cite {LazyApproxProp}.  \endProof

Recall from \ref {Theorem 3.8.7/BrownOzawa/2008} that a C*-algebra $A$ is \subjex {nuclear}{nuclear C*-algebra} if, for
every C*-algebra $B$, there is a unique C*-norm on $A\odot B$.  This is clearly the same as saying that the maximal and
minimal tensor norms on $A\odot B$ coincide.

The next result gives sufficient conditions for the nuclearity of cross sectional C*-algebras.

\state Proposition \label AproxImplyNuc Suppose that the Fell bundle $\Bun $ satisfies the approximation property and
that $B_1$ is nuclear.  Then $C^*(\Bun )$, which is necessarily isomorphic to $\CrB $, is also nuclear.

\Proof We must check that, for any C*-algebra $A$, there is only one C*-norm on $A\odot C^*(\Bun )$.  This is equivalent
to showing that the canonical map from the maximal to the minimal tensor product, shown in the first row of the diagram
below, is an isomorphism.

\beginpicture \setcoordinatesystem units <0.0020truecm, -0.0020truecm> point at 0 0 \setplotarea x from -3300 to 4000, y
from -400 to 1000 \put {$A \tmax C^*(\Bun )$} at -2000 000 \put {$A \tmin C^*(\Bun )$} at 2000 000 \put {$C^*\big (A
\tmax \Bun \big )$} at -2000 1000 \put {$C^*(A \tmin \Bun )$} at 000 1000 \put {$C^*\lred (A \tmin \Bun )$} at 2100 1000
\arrow <0.11cm> [0.3,1.2] from -1100 000 to 1100 000 \arrow <0.11cm> [0.3,1.2] from -2000 750 to -2000 250 \arrow
<0.11cm> [0.3,1.2] from 2000 750 to 2000 250 \put {$=$} at -1000 1000 \arrow <0.11cm> [0.3,1.2] from 770 1000 to 1200
1000 \put {$\Lambda $} at 950 800 \endpicture

\bigskip Consider the isomorphisms of \cite {TensorMaxFull} and \cite {TensorMinRed} in place of the left and right-hand
vertical arrows above, respectively.  Observe, in addition, that since $B_1$ is nuclear by hypothesis, the maximal and
minimal norms on $A\odot B_1$ coincide, so that $A \tmax \Bun $ and $A \tmin \Bun $ are in fact equal.  Moreover $A
\tmin \Bun $ satisfies the approximation property by \cite {TensorAproxProp}, and hence is amenable by \cite {Main}, so
its regular representation, marked as $\Lambda $ in the above diagram, is an isomorphism.  Finally, since all maps above
are essentially the identity on the various dense copies of $A\odot \CCB $, we deduce that the diagram commutes, which
implies that the arrow in the first row of the diagram is an isomorphism, as desired.  \endProof

The following is a a partial converse of the above result:

\state Theorem \label NuclearAmenable If the reduced cross-sectional C*-algebra of a Fell bundle $\Bun $ is nuclear,
then $\Bun $ is amenable.

\Proof Consider the diagram $$ \matrix { C^*(\Bun ) &\>{\weird } & \CrB \tm \Cr G \cr \pilar {20pt} \Lambda \ \Big
\downarrow && \quad \Big \downarrow \ q\cr \pilar {20pt} \CrB &\> {\ds \sigma }& \CrB \*\lmin \Cr G \cr } $$ where
$\sigma $ is provided by \cite {MiniAbsorp}, $\weird $ \ by \cite {CrazyFact}, and $q$ is the natural map from the
maximal to the minimal tensor product.  By checking on elements of the form $\hj _g(b)$, it is easy to see that the
diagram commutes.  Assuming that $\CrB $ is nuclear, we have that $q$ is injective by \ref {3.6.12/BrownOzawa/2008}, and
hence $\Lambda $ is injective as well.  \endProof

Our next result gives sufficient conditions for the reduced cross sectional C*-algebra of a Fell bundle to be exact.

\state Proposition \label ExactReducedCrossSec Let $\Bun = \{B_g\}_{g\in G}$ be a Fell bundle over an exact group $G$,
such that $B_1$ is an exact C*-algebra.  Then $\CrB $ is exact.

\Proof Let $$ 0 \to J \> \iota A \> \pi Q \to 0 $$ be an exact sequence of C*-algebras.  We need to prove that $$ 0 \to
J\tmin \CrB \> {\iota \*id} A\tmin \CrB \> {\pi \*id} Q\tmin \CrB \to 0 \equationmark MustProveExact $$ is also exact.

Considering the Fell bundles $J\tmin \Bun $, $A\tmin \Bun $, and $Q\tmin \Bun $, introduced in \cite {DefineFellTP}, we
will next show that $J\tmin \Bun $ is naturally isomorphic to an ideal in $A\tmin \Bun $, and the corresponding quotient
is isomorphic to $Q\tmin \Bun $.  By \ref {Theorem 3.5.3/BrownOzawa/2008} we have that both $$ \iota \otimes id : J\odot
B_1 \to A\odot B_1 \and \pi \otimes id : A\odot B_1 \to Q\odot B_1 $$ are continuous for the minimal tensor product
norms, so we may employ \cite {TensorWitMap} to conclude that there are morphisms $$ \tilde \iota = \{\iota _g\}_{g\in
G} \and \tilde \pi = \{\pi _g\}_{g\in G} $$ from $J\tmin \CrB $ to $A\tmin \CrB $, and from there to $Q\tmin \CrB $,
respectively.

By \ref {Theorem 3.6.1/BrownOzawa/2008} we have that $\iota \otimes id$, also known as $\iota _1$, is injective on
$J\tmin B_1$, so all the $\iota _g$ are isometric by \cite {MorpLikeCStar/a}.  We may therefore identify $J\tmin \CrB $
with a \subFellBun dle of $A\tmin \CrB $.

By first checking the conditions in \cite {HeredIdealBun/b} for the corresponding dense algebraic tensor products it is
easy to see that $J\tmin \CrB $ is in fact an ideal in $A\tmin \CrB $.

It is evident that each $\pi _g$ vanishes on $J\tmin B_g$, and we claim that the null space of $\pi _g$ is precisely
$J\tmin B_g$.  For the special case $g=1$ this follows from the exactness of the sequence $$ 0 \to J\tmin B_1 \> {\iota
\*id} A\tmin B_1 \> {\pi \*id} Q\tmin B_1 \to 0, $$ since $B_1$ is assumed to be an exact C*-algebra.  Given an
arbitrary $g$ in $G$, and given any $x$ in $A\tmin B_g$, with $\pi _g(x)=0$, we have that $$ 0 = \pi _g(x)^*\pi _g(x) =
\pi _1(x^*x), $$ from where we conclude that $x^*x$ lies in $J\tmin B_1$.  Then $$ x \={Roots} \lim _{n\to \infty
}x(x^*x)^{1/n} \in \Clspan {(A\tmin B_g)(J\tmin B_1)} \subseteq J\tmin B_g.  $$ This shows that indeed $J\tmin B_g =
\Ker (\pi _g)$.

Observe that the range of each $\pi _g$ contains the dense subspace $Q\odot B_g$, so $\pi _g$ is surjective by \cite
{MorpLikeCStar/b}.  Together with the above conclusion about the kernel of $\pi _g$ we then conclude that $Q\tmin B_g$
is isomorphic to the quotient of $A\tmin B_g$ by $J\tmin B_g$.  In other words, $Q\tmin \Bun $ is isomorphic to the
quotient Fell bundle $$ \Big (A\tmin \Bun \Big ) \Big / \Big (J\tmin \Bun \Big ) .  $$

We may then invoke \cite {ReduExactSeq} to obtain the exact sequence of C*-algebras $$ 0 \to \CstarRed \big (J\tmin \Bun
\big ) \to \CstarRed \big (A\tmin \Bun \big ) \to \CstarRed \big (Q\tmin \Bun \big ) \to 0.  $$ The isomorphisms
provided by \cite {TensorMinRed} may now be used to transform this sequence into \cite {MustProveExact}, which is
therefore also an exact sequence.  This concludes the proof.  \endProof

\nrem This chapter is based on \ref {Sections 5 and 6/AraExelKatsura/2011}.  In particular, Theorem \cite
{NuclearAmenable} first appeared in \ref {Theorem 6.4/AraExelKatsura/2011} and Proposition \cite {ExactReducedCrossSec}
is from \ref {Proposition 5.2/AraExelKatsura/2011}.

\chapter Smash product

\def \ShortBarSum #1{\kern 3pt \overline {\kern -5pt \pilar {9pt} #1 \kern -2pt}\kern 2pt}

The title of this chapter, as well as the terminology for the main concept to be introduced here, is taken from the
theory of Hopf algebras, where one defines the \"{smash product} of an algebra by the \"{co-action} of a Hopf algebra.

In the special case of the Hopf C*-algebra $C^*(G)$, where $G$ is a discrete group, Quigg has shown that co-actions
correspond to Fell bundles \ref {Quigg/1996}.  Incidentally, our study of Fell bundles may therefore be seen as a
special case of the theory of co-actions of Hopf C*-algebras.  Having decided to concentrate on Fell bundles, rather
than more general co-actions, we will likewise introduce the notion of smash product, below, in an ad-hoc way, avoiding
the whole apparatus of Hopf algebras.  Our limited view of smash products will nevertheless have many important
applications in the sequel.

Given a group $G$, we will denote the algebra of all compact operators on $\ell ^2(G)$ by $\Kp \big (\ell ^2(G)\big )$.
Given $g$ in $G$, we will let $e_g$ be the canonical basis vector of $\ell ^2(G)$, and for each $g,h\in G$, we will
denote by $\E gh$ the rank-one operator on $\ell ^2(G)$ defined by $$ \E gh(\xi ) = \langle \xi ,e_h\rangle e_g \for \xi
\in \ell ^2(G), $$ so that $$ \E gh(e_k) = \delta _{h,k}e_g \for k\in G.  $$ It is well known that $\Kp \big (\ell
^2(G)\big )$ is the closed linear span of the set formed by all the $\E gh$.

Since $\Kp \big (\ell ^2(G)\big )$ is a nuclear algebra, the minimal tensor product of a C*-algebra $A$ by $\Kp \big
(\ell ^2(G)\big )$ is isomorphic to its maximal version, so we will use the symbol ``$\otimes $'' in $$ A\otimes \Kp
\big (\ell ^2(G)\big ), $$ indistinctly meaning ``$\tmin $'' or ``$\tmax $''.

\state Proposition \label IntroSmashZero If $\Bun = \{B_g\}_{g\in G}$ is a Fell bundle, then the subset of\/ $\CB
\otimes \Kp \big (\ell ^2(G)\big )$ given by $$ \smshz = \soma {g,h\in G}B_{g\inv h}\otimes \E gh $$ is a *-subalgebra.

\Proof It is enough to notice that, for every $g,h,k,l\in G$, one has that $$ (B_{g\inv h}\otimes \E gh) (B_{k\inv
l}\otimes \E kl ) \subseteq \delta _{h,k} (B_{g\inv h}B_{h\inv l}\otimes \E gl ) \subseteq B_{g\inv l}\otimes \E gl , $$
and also that $$ (B_{g\inv h}\otimes \E gh)^* = B_{h\inv g}\otimes \E hg .  \omitDoubleDollar \endProof

\definition \label IntroSmash The \subj {smash product}\fn {One may prove that $C^*(\hbox {\rsfootnote B})$ admits a
co-action of the Hopf algebra $C^*(G)$, and a construction based on the usual Hopf algebra concept of smash product
leads to the notion presently being introduced.}  of the Fell bundle $\Bun $ \ by $G$, denoted $\smsh $, is the closed
*-subalgebra of $\CB \otimes \Kp \big (\ell ^2(G)\big )$ given by the closure of $\smshz $.

\state Remark \rm If $B$ is any graded C*-algebra whose grading coincides with our Fell bundle $\Bun $, we could define
a variant of $\smshz $ by looking at $$ \soma {g,h\in G}B_{g\inv h}\otimes \E gh $$ as a subalgebra of $B\otimes \Kp
\big (\ell ^2(G)\big )$, as opposed to $\CB \otimes \Kp \big (\ell ^2(G)\big )$.  Its closure within $B\otimes \Kp \big
(\ell ^2(G)\big )$ could then be taken as an alternative definition of $\smsh $.  However it is not difficult to see
that this alternative smash product is isomorphic to the one defined in \cite {IntroSmash}.  The reason is that for all
finite subsets $F\subseteq G$, the set $$ S_F:= \soma {g,h\in F}B_{g\inv h}\otimes \E gh, $$ (no closure) seen within
$B\otimes \Kp \big (\ell ^2(G)\big )$, is a closed *-subalgebra whose isomorphism class clearly does not depend on the
way $\Bun $ is represented in $B$.  Since the smash product (any variant of it) is the inductive limit of the $S_F$, as
$F$ ranges over the finite subsets of $G$, we then see that the smash product itself does not depend on the graded
algebra $B$.  Our choice of $\CB $ in the definition of $\smsh $, above, is therefore arbitrary.

Faithfully representing $\CB $ on a Hilbert space $H$, we have that $$ \CB \otimes \Kp \big (\ell ^2(G)\big ) \subseteq
\Lin \big (H\otimes \ell ^2(G)\big ), \equationmark RepresOfSmash $$ (where we think of the tensor product sign in the
left-hand-side above as the spatial tensor product).  Hence we may also view $\smsh $ as an algebra of operators on
$H\otimes \ell ^2(G)$.

Observe that we have for all $g,h,k\in G$, that $$ (\uE gg) (B_{h\inv k}\otimes \E hk) = \delta _{gh}(B_{h\inv k}\otimes
\E hk) \subseteq \smsh , \equationmark EggIsMult $$ from where we see that $(\uE gg) (\smsh ) \subseteq \smsh $, and
similarly one has that $(\smsh )(\uE gg) \subseteq \smsh $.  In other words, $\uE gg$ is a multiplier of $\smsh $.

It should be noticed, however, that the same is not true for $\uE gh$, when $g\neq h$.

\state Proposition \label IntroWgh Given any $w$ in $\smsh $, and given $g$ and $h$ in $G$, there is a unique
$w_{g,h}\in B_{g\inv h}$ such that $$ (\uE gg) w (\uE hh) = w_{g,h}\otimes \E gh.  $$

\Proof This is obvious for $w$ in $\smshz $.  So the result follows from the density of the latter in $\smsh $.
\endProof

We will soon be dealing with numerous ideals in $\smsh $.  In preparation for this we present the following simple
criterion for deciding when an element of $\smsh $ belongs to a given ideal.

\state Proposition \label IdealMembership Let $J$ be a closed subspace of $\smsh $.  Given $w$ in $\smsh $, consider the
statements: \izitem \zitem $w_{g,h}\otimes \E gh\in J$, for all $g$ and $h$ in $G$, \zitem $w\in J$.  \medskip \noindent
Then (i) implies (ii).  In case $J$ satisfies $$ (1\otimes \Kp )J (1\otimes \Kp )\subseteq J, $$ in particular if $J$ is
a two-sided ideal, the converse also holds.

\Proof For every finite subset $F\subseteq G$, let $$ P_F = \soma {g\in F}\uE gg.  $$ We then claim that, $$ w = \lim
_{F\uparrow G} P_FwP_F \for w\in \smsh , $$ where we think of the collection of all finite subsets $F\subseteq G$ as a
directed set relative to set inclusion.  The reason for the claim is that it clearly holds for $w\in \smshz $, which is
dense in $\smsh $, while the $P_F$ are uniformly bounded.

Observe that for $F$ as above, we have $$ P_FwP_F = \soma {g,h\in F}(\uE gg)w(\uE hh) \={IntroWgh} \soma {g,h\in
F}w_{g,h}\otimes \E gh.  $$

Thus, assuming that $w$ satisfies (i), we deduce that $P_FwP_F$ is in $J$, and hence also that $w$ is in $J$, because
$J$ is closed under taking limits.

The last assertion in the statement is obvious, observing that closed two-sided ideals are invariant under left and
right multiplication by multipliers of the ambient algebra.  \endProof

\definition The \subj {restricted smash product} of the Fell bundle $\Bun $ by $G$, denoted $\smshid $, is defined to be
$$ \smshid = \ShortBarSum {\soma {g,h\in G}\clspan {B_{g\inv } B_h}\otimes \E gh}, $$ where the closure is taken within
$\CB \otimes \Kp \big (\ell ^2(G)\big )$.

Observe that $\smshid \subseteq \smsh $, because $\clspan {B_{g\inv } B_h} \subseteq B_{g\inv h}$.  Moreover, unless
$\Bun $ is saturated, this is a proper inclusion.

\state Proposition One has that $\smshid $ is a closed two-sided ideal of $\smsh $.

\Proof Given $g,h,k,l\in G$, notice that $$ (B_{g\inv h}\otimes \E gh) (B_{k\inv } B_l\otimes \E kl ) \subseteq \delta
_{h,k} (B_{g\inv h}B_{h\inv } B_l\otimes \E gl ) \$\subseteq B_{g\inv } B_l\otimes \E gl \subseteq \smshid , $$ from
where we deduce that $\smshid $ is a left ideal of $\smsh $, and a similar reasoning proves it to be a right ideal as
well.  \endProof

We may now use our membership criterion \cite {IdealMembership} to characterize elements of $\smshid $ as follows:

\state Proposition \label MemberSmIdeal Let $w$ be in $\smsh $. Then the following are equivalent: \izitem \zitem
$w_{g,h}\in [B_{g\inv } B_h]$ (closed linear span), for every $g$ and $h$ in $G$, \zitem $w\in \smshid $.

\Proof (i) $\Rightarrow $ (ii)\enspace Follows from \cite {IdealMembership}.

\medskip \noindent (ii) $\Rightarrow $ (i)\enspace Evident.  \endProof

Recalling that $\lreg $ denotes the left-regular representation of $G$ on $\ell ^2(G)$, notice that $$ \lreg _g\E hk =
\E {gh}k \and \E hk\lreg _g = \E h{g\inv k} \for g,k,h\in G.  $$ Consequently, $$ \lreg _g\E hk \lreg _{g\inv } = \E
{gh}{gk}, $$ from where we see that $$ (1\otimes \lreg _g) (B_{h\inv k}\otimes \E hk) (1\otimes \lreg _{g\inv }) =
B_{h\inv k}\otimes \E {gh}{gk} \$= B_{(gh)\inv (gk)}\otimes \E {gh}{gk} \ \subseteq \ \smsh .  $$ This implies that
$\smsh $ is invariant under conjugation by $1\otimes \lreg _g$, so we may define an action $\rEta $ of $G$ on $\smsh $
by $$ \rEta _g:w\in \smsh \mapsto (1\otimes \lreg _g)w(1\otimes \lreg _{g\inv }) \in \smsh .  \equationmark
DefineEtaAction $$

Notice that the ideal $\smshid $ is not necessarily invariant under $\rEta $, since $$ (1\otimes \lreg _g) (B_{h\inv
}B_k\otimes \E hk) (1\otimes \lreg _{g\inv }) = B_{h\inv }B_k\otimes \E {gh}{gk}, $$ and there is no reason why $$
B_{h\inv }B_k \trepa ?\subseteq \clspan {B_{h\inv g\inv }B_{gk}}, $$ even though both of these are subset of $B_{h\inv
k}$.  In case $\Bun $ is saturated then $\clspan {B_{h\inv }B_k}$ and $\clspan {B_{h\inv g\inv }B_{gk}}$ both coincide
with $B_{h\inv k}$, and hence the above inclusion would hold, but in general it does not.

Regardless of this lack of invariance, we may still restrict $\rEta $ to a \"{partial action} $\rTh $ of $G$ on $\smshid
$, according to \cite {DefineRestriction}.

\definition \label RegPartActn Let $\Bun $ be a Fell bundle.  \iaitem \aitem The global action $\rEta $ of $G$ on $\smsh
$, defined above, will be called the \subj {dual global action} for $\Bun $.  \aitem The partial action $\rTh $ of $G$,
obtained by restricting $\rEta $ to $\smshid $, will be called the \subj {dual partial action} for $\Bun $.

\medskip The justification for this terminology is as follows: if $\Bun $ is the Fell bundle formed by the spectral
subspaces for an action $\Th $ of a compact abelian group $K$ on a C*-algebra $B$, one may show that the crossed product
$B\Rt \Th K$ is isomorphic to $\smsh $, and that dual action of $\widehat K$ on $B\Rt \Th K$ is equivalent to the
\"{dual global action} defined above.

\medskip The dual partial action will play a central role from now on, so it pays to describe the ideals involved:

\state Proposition \label DomainStdAct For each $g$ in $G$, let $ E_g = \rEta _g(\smshid )\cap (\smshid ).  $ Then $$
E_g = \ShortBarSum {\soma {h,k\in G} \clspan {B_{h\inv } D_gB_k}\otimes \E hk}, $$ where $D_g = \clspan {B_gB_{g\inv
}}$.

\Proof Given $w$ in $E_g$, in particular $w$ is in $\smshid $, so we have by \cite {MemberSmIdeal} that $w_{h,k}\in
\clspan {B_{h\inv }B_k}$, for all $h$ and $k$ in $G$.  Since $w$ is also in $\rEta _g(\smshid )$, then $ y:=\rEta
_{g\inv }(w) \in \smshid , $ so $$ \clspan {B_{h\inv g}B_{g\inv k}} \ni y_{g\inv h,g\inv k} = w_{h,k}.  $$

Choosing approximate identities $\{u_i\}_i$ and $\{v_j\}_j$, for the ideals $\clspan {B_{h\inv }B_h}$ and $\clspan
{B_{k\inv }B_k}$, respectively, we have by \cite {ApproximateIdentityTwo} that $$ w_{g,h} = \lim _{i,j} u_iw_{g,h}v_j
\in \clspan {B_{h\inv }B_h B_{h\inv g}B_{g\inv k} B_{k\inv }B_k} \$\subseteq \clspan {B_{h\inv }B_g B_{g\inv } B_k} =
\clspan {B_{h\inv }D_g B_k}.  $$ That $w$ is in the set on the right-hand side in the statement then follows from \cite
{IdealMembership}, therefore proving ``$\subseteq $''.

Conversely, given $g,h,k\in G$, notice that $D_g \subseteq B_1$, so $D_gB_k\subseteq B_k$, whence $$ B_{h\inv }
D_gB_k\otimes \E hk \subseteq B_{h\inv } B_k\otimes \E hk \subseteq \smshid .  $$ On the other hand we also have that $$
\clspan {B_{h\inv } D_g B_k}= \clspan {B_{h\inv } B_gB_{g\inv } B_k} \subseteq \clspan {B_{h\inv g} B_{g\inv k}}, $$ so
$ B_{h\inv } D_g B_k\otimes \E {g\inv h}{g\inv k} $ is also contained in $\smshid $.  We thus have that $$ B_{h\inv }
D_g B_k\otimes \E hk \$= (1\otimes \lreg _g)(B_{h\inv } D_g B_k\otimes \E {g\inv h}{g\inv k})(1\otimes \lreg _{g\inv })
\subseteq \rEta _g(\smshid ), $$ proving that each $B_{h\inv } D_gB_k\otimes \E hk$ is contained in $E_g$.  This shows
the remaining inclusion ``$\supseteq $'', and hence the proof is concluded.  \endProof

Since $\rTh $ is a restriction of $\rEta $, it is interesting to ask what exactly is the globalization of $\rTh $.  The
answer could not be other than the corresponding dual global action!

\state Proposition \label GlobaGloba The dual global action for a Fell bundle is the globalization of the corresponding
dual partial action.

\Proof Calling our bundle $\Bun = \{B_g\}_{g\in G}$, all we must do is prove that $$ \soma {g\in G}\rEta _g(\smshid ) $$
is dense in $\smsh $. For each $g$ and $h$ in $G$, observe that $B_{g\inv h}\otimes \E 1{g\inv h}$ is contained in
$\smshid $.  Moreover $$ \rEta _g(B_{g\inv h}\otimes \E 1{g\inv h}) = B_{g\inv h}\otimes \E gh, $$ so $B_{g\inv
h}\otimes \E gh$ is contained in the orbit of $\smshid $ under $\rEta $, from where the proof follows.  \endProof

\nrem Algebras resembling the restricted smash product in the context of partial actions first appeared in \ref
{Exel/1994a}.  In the case of partial actions of continuous groups, the smash product was used by Abadie in \ref
{Abadie/1999} and \ref {Abadie/2003}.

\syschapter {Stable Fell bundles as partial crossed products}{Stable Fell bundles as crossed products}

As we have seen in \cite {PartCPIsGraded/vi\&vii}, a partial crossed product is always a graded C*-algebra.  In this
chapter we will present one of the most important results of the theory of partial actions, proving a converse of the
above statement under quite broad hypotheses.

The most general form of such a converse is unfortunately not true, meaning that not all graded C*-algebras are partial
crossed products.  This may be seen from the following example: consider the ${\bf Z}$-grading of $M_3({\bf C})$ given
by $$ B_0 = \left [\matrix { \bigstar & 0 & 0 \cr 0 & \bigstar & \bigstar \cr 0 & \bigstar & \bigstar } \right ], \quad
B_1 = \left [\matrix { 0 & 0 & 0 \cr \bigstar & 0 & 0 \cr \bigstar & 0 & 0 } \right ], \quad B_{-1} = \left [\matrix { 0
& \bigstar & \bigstar \cr 0 & 0 & 0 \cr 0 & 0 & 0 } \right ], $$ while $B_n=\{0\}$, for all $n\in {\bf Z}\setminus
\{0,\pm 1\}$.  In order to see that this grading does not arise from a partial action, let us argue by contradiction and
suppose that there is a partial action $$ \Th = \big (\{\D n\}_{n\in {\bf Z}}, \{\th n\}_{n\in {\bf Z}}\big ) $$ of
${\bf Z}$ on $B_0$ whose associated semi-direct product bundle coincides with the Fell bundle given by the above grading
of $M_3({\bf C})$.  In this case, one would have $$ \clspan {B_{-1}B_1} \={LotsRigInnProd} D_{-1}\ \simeq \ D_{1}\ =\
\clspan {B_1B_{-1}}, \equationmark FalsoIso $$ (brackets meaning closed linear span) but notice that $$ \clspan
{B_{-1}B_{1}} = \left [\matrix { \bigstar & 0 & 0 \cr 0 & 0 & 0 \cr 0 & 0 & 0 } \right ] \and \clspan {B_{1}B_{-1}} =
\left [\matrix { 0 & 0 & 0 \cr 0 & \bigstar & \bigstar \cr 0 & \bigstar & \bigstar } \right ], $$ which are not
isomorphic algebras, hence bringing about a contradiction.

Further analyzing this contradiction, notice that even though $\clspan {B_{-1}B_{1}}$ and $\clspan {B_{1}B_{-1}}$ fail
to be isomorphic, they are \RME lent, given that $B_1$ is an imprimitivity bimodule.  More generally, if $B$ is any
graded C*-algebra with grading $\{B_g\}_{g\in G}$, then for every $g$ in $G$, one has that $$ D_g:= \clspan {B_gB_{g\inv
}} $$ is an ideal in $B_1$.

Assuming that the grading arises from a partial action, as in \cite {PartCPIsGraded}, then $D_{g\inv }$ is isomorphic to
$D_g$ by the reasoning used in \cite {FalsoIso}.  But, regardless of this assumption, $D_{g\inv }$ is always \RME lent
to $D_g$, with $B_g$ playing the role of the imprimitivity bimodule.  Thus, considering \RME lence as a weak form of
isomorphism, we see that the rudiments of a partial action are present in any graded C*-algebra.  Furthermore, in case
$B_1$ is a separable stable C*-algebra, then by \ref {BrownGreenRieffel/1977} the above \RME lence actually implies that
$D_{g\inv }$ is isomorphic to $D_g$, so we may choose, for each $g$, an isomorphism $$ \th g:D_{g\inv }\to D_g, $$
getting us even closer to obtaining a partial action.

In this chapter we will carefully explore these ideas in order to prove that every separable Fell bundle with stable
unit fiber algebra is isomorphic to the semi-direct product bundle for a suitable partial action of the base group on
$B_1$.

\bigskip

Recall that a C*-algebra $A$ is said to be \subjex {stable}{stable C*-algebra} if $A$ is isomorphic to the tensor
product of some other C*-algebra $B$ by the algebra $\Kp $ of all compact operators on a separable, infinite dimensional
Hilbert space.  In symbols $$ A \simeq B\otimes \Kp .  $$ We again refrain from specifying either ``$\tmin $'' or
``$\tmax $'', since $\Kp $ is a nuclear C*-algebra, whence the minimal and maximal norms are equal.

Since $\Kp \otimes \Kp \simeq \Kp $, in case $A$ is stable we have $$ A \simeq B\otimes \Kp \simeq B\otimes \Kp \otimes
\Kp \simeq A\otimes \Kp , $$ which is to say that the algebra $B$, referred to above, might as well be taken to be $A$
itself.

We will now present a useful criterion for the stability of C*-algebras.

\state Lemma \label CharStable A C*-algebra $A$ is stable if and only if there exists a {\nonDegHomo }
*-homo\-mor\-phism from $\Kp $ to the multiplier algebra $ \Mult (A).  $

\Proof Assuming that $A$ is stable, write $A=B\otimes \Kp $, where $B$ is a C*-algebra.  Supposing without loss of
generality that $B$ is faithfully represented on a Hilbert space $H$, we may view $B\otimes \Kp $ as an algebra of
operators on $H\otimes \ell ^2$.  Defining $$ \gamma : k\in \Kp \mapsto 1\otimes k\in \Lin (H\otimes \ell ^2), $$ it is
easy to see that the range of $\gamma $ is contained in the algebra of multipliers of $B\otimes \Kp $ and that, seen as
a map from $\Kp $ to $\Mult (B\otimes \Kp )=\Mult (A)$, one has that $\gamma $ is {\nonDegHomo }.

Conversely, given a {\nonDegHomo } *-homo\-mor\-phism from $\Kp $ to $\Mult (A)$, let $\{\E ij\}_{i,j\in {\bf N}}$ be
the standard matrix units of $\Kp =\Kp (\ell ^2)$, and let $$ B = \gamma (\E 11)A\gamma (\E 11).  $$ For each $n\in {\bf
N}$ one may easily prove that the map $$ \varphi _n\ :\ \{b_{i,j}\}_{i,j=1}^n \ \in \ M_n(B) \ \longmapsto \usoma
{i,j=1}n \gamma (\E i1) b_{i,j} \gamma (\E 1j) \ \in \ A, $$ is an injective *-homomorphism.  After also checking that
these maps are compatible with the usual inductive limit structure of $$ B\otimes \Kp \simeq \lim _{\longrightarrow }
M_n(B), $$ we obtain a *-homomorphism $$ \varphi : B\otimes \Kp \to A, $$ which is easily seen to be injective.  In
order to prove that $\varphi $ is also surjective, first notice that $$ a = \lim _{n\to \infty }\ \usoma {i,j=1}n \gamma
(\E ii) a \gamma (\E jj) \for a\in A, \equationmark AproxNDeg $$ since $\gamma $ is assumed to be {\nonDegHomo }.  Given
$a\in A$, and setting $$ b_{i,j} = \gamma (\E 1i) a \gamma (\E j1) \for i,j\leq n, $$ one has that $b:=
\{b_{i,j}\}_{i,j=1}^n$ is an element of $M_n(B)$, and that $$ \varphi _n(b) = \usoma {i,j=1}n\gamma (\E i1) \gamma (\E
1i) a \gamma (\E j1) \gamma (\E 1j) = \usoma {i,j=1}n\gamma (\E ii) a \gamma (\E jj) \>{(\AproxNDeg )} a, $$ proving
that $\varphi $ is surjective.  Therefore $A\simeq B\otimes \Kp $, and hence $A$ is stable.  \endProof

Our first use of this tool will be in proving the stability of $\smshid $ under appropriate hypotheses.

\state Lemma \label StableFell Let $\Bun = \{B_g\}_{g\in G}$ be a Fell bundle over a countable group $G$.  If $B_1$ is
stable, then so is $\smshid $.

\Proof Recall that if $\CB $ is faithfully represented on a Hilbert space $H$, then $$ \smshid \ \subseteq \ \CB \otimes
\Kp \ \subseteq \ \Lin \big (H\otimes \ell ^2(G)\big ), $$ where $\Kp = \Kp \big (\ell ^2(G)\big )$.  Considering the
representation of $B_1$ on $H\otimes \ell ^2(G)$ given by $$ \pi (b) = b\otimes 1 \for b\in B_1, $$ it is easy to see
that the range of $\pi $ is contained in the multiplier algebra of $\smshid $, so we may view $\pi $ as a *-homomorphism
$$ \pi :B_1\to \Mult (\smshid ), $$ (which should not be confused with another rather canonical mapping from $B_1$ to
$\smshid $, namely $b\to b\otimes \E 11$).

We claim that $\pi $ is {\nonDegHomo }.  In order to see this, let $\{v_i\}_{i\in I}$ be an approximate identity for
$B_1$.  Then, for every $g,h\in G$, and every $b\in B_{g\inv }B_h$, one has that $$ \pi (v_i)(b\otimes \E gh) =
(v_i\otimes 1)(b\otimes \E gh) = v_ib\otimes \E gh \convrg i b\otimes \E gh, $$ by \cite {ApproximateIdentity}.  This
shows that $\pi $ is {\nonDegHomo }, as desired.

It is well known that if $C$ and $D$ are C*-algebras, and $\varphi :C\to \Mult (D)$ is a given {\nonDegHomo }
*-homomorphism, then $\varphi $ admits a unique extension to a unital *-homomorphism $\tilde \varphi :\Mult (C)\to \Mult
(D)$.  Applying this to the present situation, let $$ \tilde \pi :\Mult (B_1)\to \Mult (\smshid ) $$ be the extension of
$\pi $ thus obtained.  Since $B_1$ is stable, we may use \cite {CharStable} to get a {\nonDegHomo } *-homomorphism
$\gamma :\Kp \to \Mult (B_1)$, and we may then consider the composition $$ \Kp \> \gamma \Mult (B_1) \>{\tilde \pi }
\Mult (\smshid ), $$ which we denote by $\sigma $.  Notice that $$ \clspan {\sigma (\Kp ) (\smshid )} = \clspan {\sigma
(\Kp ) \pi (B_1)(\smshid )} = \clspan {\tilde \pi \big (\gamma (\Kp )\big )\pi (B_1)(\smshid )} \$= \clspan {\tilde \pi
\big (\gamma (\Kp )B_1\big )(\smshid )} = \clspan {\pi (B_1)(\smshid )} = \smshid , $$ so $\sigma $ is {\nonDegHomo },
and hence $\smshid $ is stable.  \endProof

The reader could use the above method to show that, if $B_1$ is stable, then $\CB $ and $\CrB $ are also stable.

The following is a slight improvement on \ref {Lemma 2.5/Brown/1977}, using ideas from the proof of \ref {Theorem
3.4/BrownGreenRieffel/1977}:

\state Lemma \label SoupedUpBGR Let A be a separable stable C*-algebra. If $p$ is a \subj {full projection}\fn {A
projection $p$ in $\Mult (A)$ is said to be \"{full}, when $A=\clspan {ApA}$ (closed linear span).} in $\Mult (A)$ such
that $pAp$ is also stable, then there is $v\in \Mult (A)$ such that $v^*v=1$, and $vv^* = p$.

\Proof \def \va {v_{_A}}\null \def \vb {v_{_B}}\null Let $e$ be any minimal projection in $\Kp $, for example $e=\E 11$.
We first claim that there exists $\va \in \Mult (A\otimes \Kp )$ such that $$ \va ^*\va =1\otimes 1\and \va \va
^*=1\otimes e.  \equationmark ConditionsForV $$

To see this, recall that any two separable infinite dimensional Hilbert spaces are isometrically isomorphic to each
other, so there is $u$ in $\Mult (\Kp \otimes \Kp )$ such that $u^*u=1\otimes 1$, and $uu^* = 1\otimes e$.  Letting
$\varphi :A\otimes \Kp \to A$ be any *-isomorphism, notice that $$ \varphi \otimes id:A\otimes \Kp \otimes \Kp \to
A\otimes \Kp $$ is also a *-isomorphism, which therefore extends to the respective multiplier algebras.  Defining $\va =
(\varphi \otimes id)(1\otimes u)$, observe that $$ \va ^*\va = (\varphi \otimes id)\big ((1\otimes u)^*(1\otimes u)\big
) = (\varphi \otimes id)(1\otimes u^*u) \$= (\varphi \otimes id)(1\otimes 1\otimes 1) = 1\otimes 1, $$ while $$ \va \va
^* = (\varphi \otimes id)\big ((1\otimes u)(1\otimes u)^*\big ) = (\varphi \otimes id)(1\otimes uu^*) \$= (\varphi
\otimes id)(1\otimes 1\otimes e) = 1\otimes e, $$ so the claim is proved.

Since $B:= pAp$ is also supposed to be stable, the same argument applies to produce $\vb \in \Mult (B\otimes \Kp )$,
such that $$ \vb ^*\vb =p\otimes 1 \and \vb \vb ^*=p\otimes e, $$ where the slight difference between the above
equations and \cite {ConditionsForV} is due to the fact that the unit of $\Mult (B)$ is called $p$.

We now observe that separable C*-algebras possess strictly positive elements, so we may apply \ref {Lemma
2.5/Brown/1977} to show the existence of an element $w\in \Mult (A\otimes \Kp )$ such that $w^*w=1\otimes 1$, and $ww^*
= p\otimes 1$.  We then define $u= \vb w\va ^*$, $$ \matrix { 1\otimes 1 & \> {\ds \va } & 1\otimes e \cr w\big
\downarrow \vrule height 18pt depth 8pt width 0pt \cr p\otimes 1 & \> {\ds \vb } & p\otimes e. } $$

Strictly speaking the product $\vb w$ makes no sense, since $\vb $ belongs to $\Mult (B\otimes \Kp )$, while $w$ is in
$\Mult (A\otimes \Kp )$.  However, since $$ B\otimes \Kp = (p\otimes 1)(A\otimes \Kp )(p\otimes 1), $$ we may naturally
embed $\Mult (B\otimes \Kp )$ in $\Mult (A\otimes \Kp )$.  A simple calculation then shows that $$ u^*u = 1\otimes e
\and uu^* = p\otimes e.  $$ Therefore $u$ is a partial isometry and we claim that $(1\otimes e)u(1\otimes e) = u$.  In
fact, we have $$ (1\otimes e)u(1\otimes e) = (1\otimes e)uu^*u = (1\otimes e)(p\otimes e)u \$= (p\otimes e)u = uu^*u =
u.  $$

The claim proved, and since $e$ is a minimal projection, we see that $u$ must be of the form $$ u = v\otimes e, $$ where
$v\in \Mult (A)$ satisfies $v^*v=1$, and $vv^*=p$.  \endProof

From now on our study of Fell bundles will rely on \cite {SoupedUpBGR}, so we will have to restrict our attention to
Fell bundles satisfying suitable separability conditions.

\definition A Fell bundle $\Bun = \{B_g\}_{g\in G}$ is said to be \subjex {separable}{separable Fell bundle} if \izitem
\zitem $G$ is a countable group, \zitem every $B_g$ is a separable Banach space.

\medskip \fix From now on we will fix a separable Fell bundle $\Bun = \{B_g\}_{g\in G}$.  It is then easy to see that
most constructions originating from $\Bun $, such as $\CB $, $\CrB $, $\smsh $ and $\smshid $, lead to separable
C*-algebras.

\medskip Recall from \cite {EggIsMult} that, for every $g$ in $G$, one has that $\uE gg$ is a multiplier of $\smsh $,
and hence also of $\smshid $.  By \cite {MemberSmIdeal} the corresponding \"{corner} of $\smshid $ is then $$ (\uE
gg)(\smshid )(\uE gg) = \clspan {B_{g\inv }B_g}\otimes \E gg, $$ and in particular $$ (\uE 11)(\smshid )(\uE 11) =
B_1\otimes \E 11.  \equationmark BOneIsCorner $$

\state Proposition $\uE 11$ is a a full projection in $\Mult (\smshid )$.

\Proof Given $g$ in $G$, notice that $B_g\otimes \E 1g$ is contained in $\smshid $.  The result then follows immediately
from the identity $$ (B_g\otimes \E 1g)^*(\uE 11)(B_h\otimes \E 1h) = B_{g\inv }B_h\otimes \E gh.  \omitDoubleDollar
\endProof

\fix From now on we will assume, in addition, that $B_1$ is stable.  We then have from \cite {StableFell} that $\smshid
$ is a separable stable C*-algebra, and in its multiplier algebra one finds the full projection $\uE 11$, whose
associated corner is the stable algebra $B_1\otimes \E 11$.  We are then precisely in the situation of \cite
{SoupedUpBGR}, so we conclude that there is a $v\in \Mult (\smshid )$, such that $$ v^*v=1\otimes 1\and vv^*=\uE 11.
\equationmark VisBorn $$

In particular the mapping $$ \Phi _1: b\in B_1 \mapsto v^*(b\otimes \E 11)v\in \smshid \equationmark ConjugaPorV $$ is
an isomorphism.  It is our goal to prove that $\Phi _1$ is but one ingredient of an isomorphism $$ \Phi = \{\Phi
_g\}_{g\in G} $$ between $\Bun $ and the semi-direct product bundle arising from the dual partial action of $G$ on
$\smshid $ introduced in \cite {RegPartActn}.  We will eventually define each $\Phi _g$ by the formula $$ \def \phig
{v^*(b\otimes \E 1g)\Gamma _g(v)} \Phi _g(b) = \phig \delta _g \for b\in B_g, $$ and in preparation for this we first
prove the following technical facts:

\state Lemma \label TechForBunIso For every $g$ in $G$, one has that \izitem \zitem $B_g\otimes \E 1g\subseteq E_g$,
\zitem if $y\in E_g$, then $\rthi g(vy)v^* \in B_g\otimes \E {g\iinv }1$.

\Proof By \cite {DomainStdAct} we have that $$ B_{h\inv } D_gB_k\otimes \E hk \subseteq E_g, $$ for every $h$ and $k$ in
$G$.  So, plugging in $h=1$, and $k=g$, and using \cite {FellIsTRO}, we see that (i) follows.

In order to prove (ii), recall that $\rTh $ is the restriction of $\Gamma $ to $\smshid $, and that, in turn, $\Gamma $
is the adjoint action relative to the unitary representation $1\otimes \lreg $, as defined in \cite {DefineEtaAction}.
Therefore $\Gamma $ may be seen as an action of $G$ on the whole algebra of bounded operators on $H\otimes \ell ^2(G)$
(see \cite {RepresOfSmash}).

Given $y$ in $E_g$, let $x= \rthi g(vy)v^*$, so $$ x = \rth {g\inv }(vv^*vy)v^*vv^* \={VisBorn} \Gamma _{g\inv }(\uE
11)\rth {g\inv }(vy)v^*(\uE 11) \$= (\uE {g\inv }{g\inv })x(\uE 11) \={IntroWgh} x_{g\inv ,1}\otimes \E {g\inv }1.  $$
Using \cite {MemberSmIdeal}, we have that $x_{g\inv ,1}\in B_g$, proving (ii).  \endProof

We may now prove the main result of this chapter:

\state Theorem \label ClassifyBundles Let $\Bun $ be a separable Fell bundle whose unit fiber algebra $B_1$ is stable.
Then there exists a C*-algebraic partial action of $G$ of $B_1$ whose associated semi-direct product bundle is
isomorphic to $\Bun $.

\Proof \def \phig {v^*(b\otimes \E 1g)\Gamma _g(v)}\null As already mentioned we will prove that $\Bun $ is isomorphic
to the semi-direct product bundle for the dual partial action of $G$ on $\smshid $.  One may then transfer this action
over to $B_1$ via the isomorphism of \cite {ConjugaPorV}, arriving at the conclusion in the precise form stated above.

Working with the multiplier $v$ of $\smshid $, as in \cite {VisBorn}, recall that $v$ also acts as a multiplier of any
ideal of $\smshid $, such as $E_g$ and $E_{g\inv }$.  Consequently $\Gamma _g(v)$ may be seen as a multiplier of $\Gamma
_g(E_{g\inv })= E_g$, so we see that\fn {Observe that if $u=(L_u,R_u)$ and $v=(L_v,R_v)$ are multipliers of $E_g$, then
for every $a$ in $E_g$, the expression $uav$ stands for either $L_u(R_v(a))$ or $R_v(L_u(a))$, which coincide with each
other thanks to \cite {Condforextraassoc}, and to the fact that $E_g$, being a C*-algebra, is both {\nonDegAlg } and
idempotent.}  $$ v^*E_g\Gamma _g(v)\subseteq E_g.  $$

Using \cite {TechForBunIso/i}, we then have that $$ \varphi _g: b\in B_g \mapsto \phig \in E_g, $$ is a well defined map
for each $g$ in $G$.  On the other hand, by \cite {TechForBunIso/ii} we see that there is a map $$ \psi _g:E_g\to B_g,
$$ such that $$ \psi _g(y)\otimes \E {g\iinv }1 = \rthi g(vy)v^* \for y\in E_g.  \equationmark CharacPsiG $$

We will now prove that $\varphi _g$ and $\psi _g$ are each other's inverse.  For this, notice that if $y\in E_g$, then
$$ \varphi _g\big (\psi _g(y)\big ) = v^*(\psi _g(y)\otimes \E 1g)\Gamma _g(v) \$= v^*\Gamma _g\big ((\psi _g(y)\otimes
\E {g\iinv }1)v\big ) \={CharacPsiG} v^*\Gamma _g\big ( \rthi g(vy)v^* v\big ) \={VisBorn} y.  $$ On the other hand,
given $b$ in $B_g$, we have $$ \psi _g\big (\varphi _g(b)\big ) \otimes \E {g\iinv }1 = \rthi g(v\varphi _g(b))v^* =
\rthi g\big (v \phig \big )v^* \={VisBorn} $$$$= \rthi g\big ((\uE 11)(b\otimes \E 1g)\Gamma _g(v) \big )v^* = \rthi
g\big ((b\otimes \E 1g)\Gamma _g(v) \big )v^* \$= (b\otimes \E {g\iinv }1)v v^* = (b\otimes \E {g\iinv }1)(\uE 11) =
b\otimes \E {g\iinv }1, $$ from where we see that $\psi _g\big (\varphi _g(b)\big ) =b$, whence $\psi _g$ is indeed the
inverse of $\varphi _g$.  Since both of these maps are contractive, then both are in fact isometries.  We then define $$
\Phi _g: B_g \to E_g\delta _g, $$ by $$ \Phi _g(b) = \varphi _g(b)\delta _g \for b\in B_g, $$ and we will prove that the
collection of maps $\Phi = \{\Phi _g\}_{g\in G}$ gives an isomorphism from $\Bun $ to the semi-direct product bundle
relative to $\rTh $.

Since we already know that the $\Phi _g$ are bijective, it now suffices to prove that $\Phi $ is a morphism of Fell
bundles.  We begin by proving \cite {DefineMorphBun/i}.  Given $b\in B_g$, and $c\in B_h$, let $ x = \varphi _g(b), $
and $ y = \varphi _h(c).  $ We then have $$ \Phi _g(b) \Phi _h(c) = (x\delta _g)(y\delta _h) = \rth g\big (\rthi
g(x)y\big )\delta _{gh} = x\Gamma _g(y)\delta _{gh} \$= v^*(b\otimes \E 1g)\Gamma _g(v) \Gamma _g\big (v^*(c\otimes \E
1h)\Gamma _h(v)\big )\delta _{gh} \$= v^*(b\otimes \E 1g)\Gamma _g(\uE 11)(c\otimes \E g{gh})\Gamma _{gh}(v)\delta _{gh}
\$= v^*(b\otimes \E 1g)(1\otimes \E gg)(c\otimes \E g{gh})\Gamma _{gh}(v)\delta _{gh} \$= v^*(bc\otimes \E 1{gh})\Gamma
_{gh}(v)\delta _{gh} = \Phi _{gh}(bc).  $$ Referring to \cite {DefineMorphBun/ii}, pick $b$ in any $B_g$, and notice
that $$ \Phi _g(b)^* = \big (\varphi _g(b)\delta _g\big )^* = \rthi g\big (\varphi _g(b)^*\big )\delta _{g\inv } \$=
\Gamma _{g\inv }\big (\Gamma _g(v^*)(b^*\otimes \E g1)v\big )\delta _{g\inv } = v^*(b^*\otimes \E 1{g\inv })\Gamma
_{g\inv }(v)\delta _{g\inv } \$= \Phi _{g\inv }(b^*).  $$ This concludes the proof.  \endProof

As an application of this to graded algebras we present the following:

\state Corollary Suppose we are given a countable group $G$, and a separable, topologically $G$-graded C*-alge\-bra $B$.
Suppose moreover that $B_1$ is stable and the canonical conditional expectation onto $B_1$ is faithful.  Then there
exists a partial action of\/ $G$ on $B_1$ such that $$ B\simeq B_1\redrt \null G.  $$

\Proof All of the conditions of \cite {ClassifyBundles} are clearly fulfilled for the associated Fell bundle $\Bun $, so
we may assume that $\Bun $ is the semi-direct product bundle for a partial action of $G$ on $B_1$.  By \cite {Charact}
we then have that $B$ is isomorphic to the reduced cross-sectional C*-algebra of $\Bun $, also known as $B_1\redrt \null
G$.  \endProof

\nrem The fact that every separable, stable Fell bundle arises from a partial crossed product is known for almost 20
years \ref {Theorem 7.3/Exel/1997a}, except that the partial action might be twisted by a cocycle.  On the other hand,
the so called Packer-Raeburn trick \ref {Theorem 3.4/PackerRaeburn/1989} asserts that, up to stabilization, every
twisted global action is \"{exterior equivalent} to a genuine (untwisted) action.  Therefore it has been widely
suspected that the cocycle in \ref {Theorem 7.3/Exel/1997a} could be eliminated.  Theorem \cite {ClassifyBundles} does
precisely that.  It has been proven by Sehnem in her Masters Thesis \ref {Sehnem/2014}.  A purely algebraic version of
\ref {Theorem 7.3/Exel/1997a} may be found in \ref {Theorem 8.5/DokuchaevExelSimon/2008}.

\chapter Globalization in the C*-context

\def \mat #1,#2; #3,#4;{\pmatrix {#1 & #2 \cr #3 & #4}} \def \Prim {\hbox {Prim}} \def \up #1{^{\scriptscriptstyle #1}}

The question of globalization is one of the richest parts of the theory of partial actions.  Recall that partial actions
on sets always have a unique globalization \cite {GlobaSets}, and so do partial actions on topological spaces \cite
{TopoGlob}, although not always on a Hausdorff space \cite {TopoGlobIffClosed}.  In the case of algebraic partial
actions, existence \cite {NoGlobalInAlg} and uniqueness \cite {NotUniqueGloba} may fail, except when the corresponding
ideals are unital \cite {GlobaUnital}.

With such a track record, when it comes to C*-algebras we clearly shouldn't expect a smooth ride.  In order to be able
to properly discuss globalization for C*-algebraic partial actions we will need a large part of the material developed
so far and it is for this reason that the present chapter has been postponed until now.

The concept of restriction, as defined in \cite {DefineRestriction}, is perfectly suitable for the category of
C*-algebras, requiring no further adaptation.  In other words, if $\eta $ is a global C*-algebraic action of a group $G$
on an algebra $B$, and if we are given a closed two-sided ideal $A\ideal B$, the restriction of $\eta $ to $A$ gives a
bona fide C*-algebraic partial action of $G$ on $A$.  However the concept of globalization given in \cite
{DefineAlgGloba} requires some fine tuning if it is to be of any use when working with C*-algebras.  The difference
between \cite {DefineAlgGloba} and the following definition is essentially the occurrence of the word \"{closure} below.

\definition \label DefineCstarGloba Let $\eta $ be a C*-algebraic global action of a group $G$ on an algebra $B$, and
let $A$ be a closed two-sided ideal of $B$.  Also let $\Th $ be the partial action obtained by restricting $\eta $ to
$A$.  If $B$ is the \"{closure} of $$ \soma {g\in G}\eta _g(A), $$ we will say that $\eta $ is a \subj {C*-algebraic
globalization} of $\Th $.

It is not difficult to adapt the example given in \cite {NoGlobalInAlg} to produce a C*-algebraic partial action
admitting no globalization.  However, example \cite {NotUniqueGloba}, exploiting a somewhat grotesque algebraic
structure (identically zero product) to show the lack of uniqueness for globalizations of algebraic partial actions, has
no counterpart for C*-algebras: our next result shows that if a C*-algebraic partial action admits a globalization, then
that globalization is necessarily unique.

\state Proposition \label UniqueGlobaCstar Let $ \Th = (\{\D g\}_{g\in G},\ \{\th g\}_{g\in G}) $ be a C*-algebraic
partial action of the group $G$ on the algebra $A$, and suppose that for each $k=1,2$, we are given a globalization
$\eta \up k$ of\/ $\Th $, acting on a C*-algebra $B\up k$.  Then there exists an equivariant *-isomorphism $$ \varphi
:B\up 1 \to B\up 2 $$ which is the identity on the respective copies of $A$ within $B\up 1$ and $B\up 2$.

\Proof As a first step we claim that, if $a$ and $b$ are in $A$, then $$ \eta \up 1_g(a)b=\eta \up 2_g(a)b \for g\in G.
\equationmark ProdInASame $$

Choosing an approximate identity $\{v_i\}_i$ for $\Di g$, notice that $\{\th g(v_i)\}_i$ is an approximate identity for
$\D g$.  Also, since $$ \eta \up k_g(a)b \in \eta \up k_g(A)\cap A = \D g \for k=1,2, $$ we have $$ \eta \up k_g(a)b =
\lim _{i\to \infty } \th g(v_i)\eta \up k_g(a)b = \lim _{i\to \infty } \eta \up k_g(v_ia)b = \lim _{i\to \infty } \th
g(v_ia)b, $$ where the last step is justified by the fact that $v_ia$ is in $\Di g$.  Since the right-hand-side above
does not depend on $k$, the claim is proved.

Recall that each $B\up k$ is the closure of the sum of the ideals $\eta \up k_g(A)$, for $g$ in $G$.  This suggests the
following slightly more general situation: suppose we are given a C*-algebra $B$, which is the closure of the sum of a
family $\{J_i\}_{i\in I}$ of closed two-sided ideals.  Then every $b$ in $B$ acts as a multiplier of each $J_i$ by left
and right multiplication, thus providing a canonical *-homomorphism $$ \mu _i:B\to \Mult (J_i).  $$

Setting $$ \textstyle \mu = \prod _{i\in I}\mu _i : B\to \prod _{i\in I}\Mult (J_i), $$ (here the product on the
right-hand-side is defined to be the C*-algebra formed by all \"{bounded} families $x = (x_i)_{i\in I}$, with each $x_i$
in $\Mult (J_i)$, equipped with coordinate-wise operations and the supremum norm), we claim that $\mu $ is injective.
In fact, if $b\in B$ is such that $\mu _i(b)=0$, for all $i\in I$, then $$ bx=0 \for x \in \textstyle \bigcup _{i\in I}
J_i.  $$ Since the set of all such $x$'s span a dense subspace of $B$, we conclude that $b=0$.  Being injective, $\mu $
is also necessarily isometric, so for each $b$ in $B$, one has that $$ \Vert b\Vert = \Vert \mu (b)\Vert = \sup _{i\in
I} \Vert \mu _i(b)\Vert = \sup _{i\in I}\sup _{\scriptstyle x\in J_i \atop \scriptstyle \Vert x\Vert \leq 1}\Vert
bx\Vert .  $$

Returning to the above setting, given $a_1,\ldots ,a_n\in A$, and $g_1,\ldots ,g_n\in G$, we may then compute the norm
of the element $$ b = \usoma {i=1}n\eta \up 1_{g_i}(a_i)\in B\up 1, $$ as follows: $$ \Vert b\Vert = \sup _{h\in G}\sup
_{\scriptstyle a\in A \atop \scriptstyle \Vert a\Vert \leq 1}\Vert b\eta \up 1_h(a)\Vert = \sup _{h,a}\normsum {\usoma
{i=1}n\eta \up 1_{g_i}(a_i) \eta \up 1_h(a)} \$= \sup _{h,a}\normsum {\eta \up 1_h \Big (\usoma {i=1}n\eta \up 1_{h\inv
g_i}(a_i) a\Big )} = \sup _{h,a}\normsum {\usoma {i=1}n\eta \up 1_{h\inv g_i}(a_i) a} \={ProdInASame} $$$$ = \sup
_{h,a}\normsum {\usoma {i=1}n\eta \up 2_{h\inv g_i}(a_i) a} = \ \cdots \ = \normsum {\usoma {i=1}n\eta \up
2_{g_i}(a_i)}, $$ where the ellipsis indicates the reversal of our computations with the superscript ``1'' replaced by
``2''.  This implies that the correspondence $$ \usoma {i=1}n\eta \up 1_{g_i}(a_i) \mapsto \usoma {i=1}n\eta \up
2_{g_i}(a_i), \equationmark FromOneToOtherGlob $$ is well defined and extends to give an isometric linear mapping
$\varphi :B\up 1 \to B\up 2$, which we claim to satisfy all of the required conditions.

We begin with the verification that $\varphi $ is multiplicative, for which it clearly suffices to check that $$ \varphi
\big (\eta \up 1_g(a) \eta \up 1_h(b)\big ) = \eta \up 2_g(a) \eta \up 2_h(b) \for g,h\in G \for a,b\in A.  $$ We have
$$ \varphi \big (\eta \up 1_g(a) \eta \up 1_h(b)\big ) = \varphi \Big (\eta \up 1_h\big (\eta \up 1_{h\inv g}(a)b\big
)\Big ) \={ProdInASame} $$$$ = \varphi \Big (\eta \up 1_h\big (\eta \up 2_{h\inv g}(a)b\big )\Big )
\={FromOneToOtherGlob} \eta \up 2_h\big (\eta \up 2_{h\inv g}(a)b\big ) = \eta \up 2_g(a) \eta \up 2_h(b).  $$

The easy verification that $\varphi $ is equivariant and preserves the star operation is left to the reader.  \endProof

Another important aspect of globalization is that it does not mix the commutative and the non-commutative worlds:

\state Proposition \label GlobaComata Let $\eta $ be a globalization of a C*-algebraic partial action $\Th $.  If $\Th $
acts on a commutative algebra, then so does $\eta $.

\Proof Let $A$ and $B$ be the algebras where $\Th $ and $\eta $ act, respectively, so that $A$ is a commutative ideal in
$B$.  We first claim that $A$ is contained in the center of $B$.

To see this pick $a$ in $A$, and $b$ in $B$.  Using Cohen-Hewitt, we may write $a=a_1a_2$, with $a_1$ and $a_2$ in $A$.
We then have $$ ab = (a_1a_2)b = a_1(a_2b) = (a_2b)a_1 = a_2(ba_1) = (ba_1)a_2 = ba, $$ proving that $a$ is in the
center of $B$.  This argument in fact shows that any commutative idempotent ideal must be central.

Since $B$ is generated by the translates of $A$, it is clearly enough to prove that $$ \eta _g(a)\eta _h(b)= \eta _h(b)
\eta _g(a), $$ for any $a$ and $b$ in $A$, and for any $g$ and $h$ in $G$, but this is easily proven with the following
computation: $$ \eta _g(a)\eta _h(b)= \eta _g\big (a\eta _{g\inv h}(b)\big ) = \eta _g\big (\eta _{g\inv h}(b)a\big ) =
\eta _{h}(b)\eta _g(a\big ).  \omitDoubleDollar \endProof

Since partial actions on commutative C*-algebras correspond bijectively to partial actions on LCH (locally compact
Hausdorff) spaces by \cite {TopoCstarActions}, we may correlate the globalization questions for commutative C*-algebras
on the one hand, and for topological spaces, on the other:

\state Proposition \label ComutGloba Let $\Th $ be a partial action of a group $G$ on a LCH space $X$, and denote by
$\Th '$ the partial action of\/ $G$ on $C_0(X)$ corresponding to $\Th $ via \cite {TopoCstarActions}.  Then a necessary
and sufficient condition for $\Th '$ to admit a globalization is that the globalization of $\Th $ provided by \cite
{TopoGlob} take place on a Hausdorff space.

\Proof If $(\eta ,Y)$ is a globalization of $\Th $ and $Y$ is Hausdorff, it is easy to see that the corresponding action
$\eta '$ of $G$ on $C_0(Y)$ is a globalization for $\Th '$.

Conversely, if we are given a globalization $(\eta ',B)$ for $\Th '$, we have by \cite {GlobaComata} that $B$ is
commutative, hence $B\simeq C_0(Y)$, where $Y$ is the spectrum of $B$.  Denoting by $\eta $ the global action of $G$ on
$Y$ corresponding to $\eta '$ under \cite {TopoCstarActions}, it is easy to see that $\eta $ is the globalization for
$\Th $.  Being the spectrum of a commutative C*-algebra, $Y$ is Hausdorff.  \endProof

This result may easily be used to produce examples of C*-algebraic partial actions not admitting a globalization: just
take a partial action on a LCH space $X$ whose graph is not closed, so that its globalization will be non-Hausdorff by
\cite {TopoGlobIffClosed}.  The corresponding partial action on $C_0(X)$ will therefore admit no globalization by \cite
{ComutGloba}.

Having seen that the existence question for globalization of C*-algebraic partial actions often has a negative answer,
one might try to relax the question itself, and look for \RME lent actions admitting a globalization.  The following
main result shows that this can always be attained.

\state Theorem \label ExistMorGlob Every C*-algebraic partial action is \RME lent to one admitting a globalization.
More precisely, every C*-algebraic partial action is \RME lent to the dual partial action $\rTh $ on the restricted
smash product for the corresponding semi-direct product bundle (which admits a globalization by \cite {GlobaGloba}).

\Proof Given a C*-algebraic partial action $$ \Th = (\{\D g\}_{g\in G},\ \{\th g\}_{g\in G}) $$ of the group $G$ on the
C*-algebra $A$, denote by $\Bun $ its semi-direct product bundle.  Our task is therefore to construct a Hilbert
$A$-$\smshid $-bimodule $M$, and a partial action $\gamma $ of $G$ on $M$, satisfying the conditions required by \cite
{RmorPDsys}.  We begin by letting $M$ be the subspace of $\smshid $ given by $$ M = \barsum {\soma {h\in G}B_h\otimes \E
1h}.  $$

If the reader is used to thinking of elements of $\CB \otimes \Kp $ as infinite matrices with entries in $\CB $, then
$M$ should be thought of as the space of all row matrices in $\smshid $.  Observe that $M$ is invariant under
left-multiplication by $B_1\otimes \E 11$, so, upon identifying $A$ with $B_1\otimes \E 11$ via $$ a\in A \mapsto
a\delta _1\otimes \E 11\in B_1\otimes \E 11, \equationmark IdentifyABOne $$ we may use the multiplicative structure of
$\smshid $ to define the left $A$-module structure.  Given $h$, $k$ and $l$ in $G$, notice that $$ (B_h\otimes \E 1h)
(B_{k\inv } B_l\otimes \E kl)= \delta _{h,k}(B_hB_{h\inv } B_l\otimes \E 1l) \subseteq B_l\otimes \E 1l \subseteq M, $$
from where we see that $M$ is a right ideal in $\smshid $, hence also a right $\smshid $-module.

Given $\xi $ and $\eta $ in $M$, it is easy to see that $\xi \eta ^*$ (operated as elements of the C*-algebra $\smshid
$) lies in $B_1\otimes \E 11$.  So, under the identification \cite {IdentifyABOne}, we may view $\xi \eta ^*$ as an
element of $A$.  In other words, we define the $A$-valued inner-product of $\xi $ and $\eta $ to be the unique element
$\ip A\xi \eta \in A$, such that $$ \xi \eta ^* = \ip A\xi \eta \delta _1\otimes \E 11.  \equationmark FirstInProd $$

For the $\smshid $-valued inner-product we simply set $$ \ip {\smshid }\xi \eta = \xi ^*\eta \for \xi ,\eta \in M.  $$
Proving $M$ to be a Hilbert $A$-$\smshid $-bimodule is now entirely routine.  In order to construct the partial action,
we let for each $g$ in $G$, $$ M_g = \barsum {\soma {h\in G}\clspan {B_gB_{g\inv }B_h}\otimes \E 1h}.  $$

It is easy to see that $(B_1\otimes \E 11)M_g\subseteq M_g$, so that $M_g$ is a left $A$-sub-module.  In order to prove
it to be a right $\smshid $-sub-module, notice that for each $h$ and $k$ in $G$, we have $$ (B_gB_{g\inv }B_h\otimes \E
1h) (B_{h\inv }B_k\otimes \E hk) \subseteq B_gB_{g\inv }B_k\otimes \E 1k \subseteq M_g, $$ thus concluding the
verification of \cite {RmorPDsys/i}.

Recall that in \cite {DomainStdAct} we used the notation $D_g$ to mean $\clspan {B_gB_{g\inv }}$, while here $D_g$ is
supposed to mean the range of $\th g$.  Insisting on the present use of $D_g$, observe that $B_g=\D g\delta _g$, hence
$$ \clspan {B_gB_{g\inv }} = \clspan {(\D g\delta _g)(\Di g\delta _{g\inv })} \={LotsPAreceGloba} \clspan {\D g\th g(\Di
g)\delta _1} = \D g\delta _1, $$ so the current meaning of $D_g$ is compatible with the one used in \cite
{DomainStdAct}, up to the usual identification of $\D g$ with $\D g\delta _1$.

Given $g$ and $h$ in $G$ observe that $$ \clspan {B_gB_{g\inv }B_h} = \clspan {\D g\D h\delta _h} = (\D g\cap \D
h)\delta _h, $$ so $M_g$ may be alternatively described as $$ M_g = \barsum {\soma {h\in G} (\D g\cap \D h)\delta
_h\otimes \E 1h}.  \equationmark DefineMg $$

Speaking of \cite {FullIdealsInMor}, notice that for each $h$ in $G$, we have that $$ \Big ((\D g\cap \D h)\delta _h\Big
)\Big ((\D g\cap \D h)\delta _h\Big )^* \={LotsRigInnProd} (\D g\cap \D h)\delta _1, $$ from where we see that $\clspan
{\ip A{M_g}{M_g}}=\D g$.  On the other hand, given $h$ and $k$ in $G$, we have $$ \clspan {\big (B_gB_{g\inv }B_h\otimes
\E 1h\big )^* \big (B_gB_{g\inv }B_k\otimes \E 1k\big )} \$= \clspan {B_{h\inv }B_gB_{g\inv }B_gB_{g\inv }B_k\otimes \E
hk} = \clspan {B_{h\inv }D_gB_k\otimes \E hk}, $$ so we deduce from \cite {DomainStdAct} that $ \clspan {\ip \smshid
{M_g}{M_g}} = E_g, $ as desired.

We next claim that for each $g$ in $G$, there exists a bounded linear mapping $$ \gamma _g:M_{g\inv }\to M_g, $$ such
that $$ \gamma _g(a\delta _h\otimes \E 1h) = \th g(a)\delta _{gh}\otimes \E 1{gh}, \equationmark DefineGammaG $$ for all
$g$ and $h$ in $G$, and for all $a$ in $\D {g\inv }\cap \D h$.  In order to see this, first notice that if $a$ is as
above, then $$ \th g(a)\in \th g(\D {g\inv }\cap \D h) \_\subseteq {PAIntersecContains} \D g\cap \D {gh}, $$ from where
we see that the right-hand-side of \cite {DefineGammaG} indeed lies in $M_g$, as needed.  Denoting by $M_g^0$ the dense
subspace of $M_g$ defined as in \cite {DefineMg}, but without taking closures, it is therefore immediate to prove the
existence of a map $\gamma _g^0:M_{g\inv }^0\to M_g^0$, satisfying \cite {DefineGammaG}.  Given $\xi $ and $\eta $ in
$M_{g\inv }^0$, write them as finite sums $$ \xi = \soma {h\in G}a_h\delta _h\otimes \E 1h \and \eta = \soma {h\in
G}b_h\delta _h\otimes \E 1h, $$ with each $a_h$ and $b_h$ in $\D {g\inv }\cap \D h$, and notice that $$ \gamma _g^0(\xi
) \gamma _g^0(\eta )^* = \Big (\soma {h\in G}\th g(a_h)\delta _{gh}\otimes \E 1{gh}\Big )\Big (\soma {k\in G}\th
g(b_k)\delta _{gk}\otimes \E 1{gk}\Big )^* \={LotsRigInnProd} $$$$ = \soma {h\in G}\th g(a_hb_h^*)\delta _1\otimes \E
11.  $$ In view of \cite {FirstInProd}, this shows that $$ \ip A{\gamma _g^0(\xi )}{\gamma _g^0(\eta )} = \soma {h\in
G}\th g(a_hb_h^*) = \th g\big (\ip A\xi \eta \big ), $$ from where we deduce two important facts: firstly, that $\gamma
_g^0$ is an isometry, so it extends to a bounded linear map $\gamma _g$ from $M_{g\inv }$ to $M_g$, proving our claim
and, secondly, that $\gamma _g$ satisfies \cite {TringCovar} relative to the $A$-valued inner-product.  Let us now also
prove \cite {TringCovar} for the $\smshid $-valued inner-product, so we begin by taking $\xi $ and $\eta $ in $M_{g\inv
}$ of the form $$ \xi = a\delta _h\otimes \E 1h \and \eta = b\delta _k\otimes \E 1k, $$ with $a$ in $\D {g\inv }\cap \D
h$, and $b$ in $\D {g\inv }\cap \D k$.  So $$ \ip \smshid {\gamma _g(\xi )}{\gamma _g(\eta )} = \big (\th g(a)\delta
_{gh}\otimes \E 1{gh}\big )^*\big (\th g(b)\delta _{gk}\otimes \E 1{gk}\big ) \={LotsLefInnProd} $$$$ = \th {h\inv g\inv
}\big (\th g(a^*b)\big )\delta _{h\inv k}\otimes \E {gh}{gk} \={PACompos} \thi h(a^*b)\delta _{h\inv k}\otimes \E
{gh}{gk} \$= \rth g\Big (\thi h(a^*b)\delta _{h\inv k}\otimes \E hk\Big ) = \rth g\Big (\ip \smshid \xi \eta \Big ).  $$
Since the set of elements $\xi $ and $\eta $ considered above span a dense subspace of $M_{g\inv }$, we have proven
\cite {TringCovar}.

To conclude we must now prove that $\gamma $ is indeed a partial action on $M$.  Since \cite {PAIdentity} is elementary,
we prove only \cite {PAExtend}.  In order to do this, we first observe that, as a consequence of \cite
{IdealMembership}, given $w$ in $\smshid $, a necessary and sufficient condition for $w$ to be in $M_g$ is that $$
w_{1,k}\in (\D g\cap \D k)\delta _k \and w_{h,k}=0, \equationmark CondForWinMg $$ for all $h$ and $k$ in $G$, with
$h\neq 1$.  This said, we claim that, for all $g$ and $h$ in $G$, one has that $$ \gamma _h\inv (M_h\cap M_{g\inv })
\subseteq M_{(gh)\inv }.  \equationmark ClaimDomain $$ Notice that the set in the left-hand-side above is precisely the
domain of $\gamma _{g}\circ \gamma _{h}$ by \cite {DomainCompos}.  We thus pick an element $\xi $ in this set, so that
$$ \xi \in M_{h\inv } \and \gamma _{h}(\xi )\in M_{g\inv }.  $$ Letting $a_k\in \Di {h}\cap \D k$ be defined by $\xi
_{1,k}=a_k\delta _k$, observe that $$ \gamma _{h}(\xi )_{1,hk} = \th {h}(a_k)\delta _{hk}, $$ as a quick look at \cite
{DefineGammaG} will reveal.  Since $\gamma _{h}(\xi )\in M_{g\inv }$, we deduce from \cite {CondForWinMg} that $$ \gamma
_{h}(\xi )_{1,hk} \in (\Di {g}\cap \D {hk}) \delta _{hk}, $$ which implies that $$ \th {h}(a_k)\in \Di {g}\cap \D {hk}
\quad \cap \D {h}, $$ so $$ a_k \in \thi {h}(\Di {g}\cap \D {hk}\cap \D {h}) \_\subseteq {PAIntersecContains} \D {h\inv
g\inv }\cap \D k\cap \D {h\inv }.  $$ Therefore $$ \xi _{1,k}=a_k\delta _k\in (\D {h\inv g\inv }\cap \D k)\delta _k, $$
whence $\xi $ is in $M_{(gh)\inv }$, by \cite {CondForWinMg}, proving claim \cite {ClaimDomain}.  In particular, this
shows that the domain of $\gamma _{g}\circ \gamma _{h}$ is contained in the domain of $\gamma _{gh}$, and it is now easy
to see that $\gamma _{gh}(\xi )=\gamma _{g}\big (\gamma _{h}(\xi )\big )$, for all $\xi $ in the domain of $\gamma
_{g}\circ \gamma _{h}$, thus verifying \cite {PAExtend}, and hence finishing the proof.  \endProof

Starting with an arbitrary C*-algebraic partial dynamical system $$ \Th = (A,G,\{\D g\}_{g\in G},\ \{\th g\}_{g\in G}),
$$ we are then led to considering two other actions: the dual partial action $\rTh $ \cite {RegPartActn} for the
semi-direct product bundle, and its globalization $\rEta $ \cite {GlobaGloba}.

Each of these actions comes with its full and reduced crossed products:

\advseqnumbering \label Tabela \bigskip \begingroup \def \tabrule {\noalign {\hrule }} \def \poste {\pilar {15pt}\stake
{8pt}} \offinterlineskip \centerline {\vbox {\halign { \strut \vrule \hfil \quad #\quad \hfil &\vrule \hfil \quad #\quad
\hfil &\vrule \hfil \quad #\quad \hfil \vrule \cr \tabrule $\Th $ & $\rTh $ & $\rEta $ \pilar {14pt}\cr \eightrm Given
action & \eightrm Dual action & \eightrm Globalization \pilar {0pt}\stake {8pt}\cr \tabrule $A\rt \Th G$ & $(\smshid
)\rt \rTh G$ & $(\smsh )\rt \rEta G$\poste \cr \tabrule $\quad A\redrt \Th G\quad $ & $(\smshid )\redrt \rTh G$ &
$(\smsh )\redrt \rEta G$ \poste \cr \tabrule }}} \endgroup

\centerline {\eightrm \current . Table. Derived actions and crossed products.}

\bigskip We already know that $\Th $ and $\rTh $ are \RME lent as partial actions.  Even though $\rTh $ and $\rEta $ are
not necessarily \RME lent\fn {An example of a partial action not \RME lent to its globalization is given right after
\cite {GradedAlgsAsCrossProd}.  It is easy to see that the globalization of this partial action is the action of ${\bf
Z}$ on itself by translation.  These actions are not \RME lent because \RME lence among commutative algebras is
tantamount to isomorphism.}, we will now show that in each row of the above table all of the three algebras are \RME
lent to each other.  Our next two results are designed to justify this claim.

\state Theorem\fn {The statement in part (ii) of this result has been anticipated in \cite {FutureRME}.}  \label
EqPAEquivCP Let $G$ be a group and suppose we are given \RME lent C*-algebraic partial dynamical systems $$ \alpha =
\big (A,G,\{A_g\}_{g\in G},\{\alpha _g\}_{g\in G}\big ) \and \beta = \big (B,G,\{B_g\}_{g\in G},\{\beta _g\}_{g\in
G}\big ).  $$ Then \izitem \zitem $A\redrt \Th G$ and $B\redrt \eta G$ are \RME lent, and \zitem $A\rt \Th G$ and $B\rt
\eta G$ are \RME lent.

\Proof \def \L {{\rsbox {L}\,}} Let $$ \gamma = \big (M,G,\{M_g\}_{g\in G},\{\gamma _g\}_{g\in G}\big ) $$ be an
imprimitivity system for $\alpha $ and $\beta $, and let $L$ be the linking algebra of $M$.  Consider the partial action
$$ \lambda = \big (\{L_g\}_{g\in G},\ \{\lambda _g\}_{g\in G}\big ), $$ of $G$ on $L$ given by \cite {LinkPaction}.
Identified with the upper left-hand corner of $L$, it is clear that $A$ is a $\lambda $-invariant subalgebra, so we may
use case (i) of \cite {RestrGivesSubBun} to view $A\redrt \alpha G$ as a subalgebra of $L\redrt \lambda G$.

Unless $A$ is a unital algebra, there is no canonical way to view the $2\times 2$ matrix $$ \E 11 = \mat 1, 0; 0, 0; $$
as an element of $L$, but, even when $A$ is non-unital, the formal left- or right-multiplication of $\E 11$ by elements
of $L$ is easily seen to define a multiplier of $L$.

By \cite {BOneEssential}, which is stated for the full cross-sectional algebra, but which evidently also holds for the
reduced one, we have that the inclusion of $L$ in $L\redrt \lambda G$ is a {\nonDegHomo } *-homomorphism.  We may then
extend it to a *-homomorphism $$ \Mult (L) \to \Mult \big (L\redrt \lambda G\big ), $$ and we will denote the image of
$\E 11$ under this map by $\E 11\delta _1$.  Given $g$ in $G$, and $$ x = \mat a, m; n^*, b; \in L_g, $$ notice that $$
(\E 11\delta _1)(x\delta _g)(\E 11\delta _1) \={MultOnSDPBun} \lambda _ g\big (\lambda _{g\inv }(\E 11x)\E 11\big
)\delta _g \$= \lambda _ g\left ( \mat \alpha _{g\inv }(a), \gamma _{g\inv }(m); 0, 0; \E 11\right )\delta _g = \mat a,
0; 0, 0; \delta _g.  \equationmark CondexpLinkAlg $$ With this it is easy to see that $$ (\E 11\delta _1)(L\redrt
\lambda G)(\E 11\delta _1) = A\redrt \alpha G, \equationmark AHeredInL $$ from where it follows that $A\redrt \alpha G$
is a hereditary subalgebra of $L\redrt \lambda G$.

We next claim that $A$ is a full subalgebra.  To see this, we first check that the closed two-sided ideal generated by
$A\redrt \alpha G$ coincides with $$ \Clspan {(L\redrt \lambda G)(\E 11\delta _1) (L\redrt \lambda G)}.  $$ Temporarily
denoting \def \upperrt #1{\buildrel \rtimes \over #1} \def \upperrt #1{#1^{\!\scriptscriptstyle \rtimes }}\null $$
\upperrt L:= L\redrt \lambda G \and \upperrt A:=A\redrt \alpha G, $$ and denoting $\E 11\delta _1$ by $p$, notice that
$\clspan {\upperrt L p\upperrt L}$ is an ideal in $\upperrt L$, hence idempotent, so $$ \clspan {\upperrt L p\upperrt L}
= \clspan {\upperrt L p\upperrt L \upperrt Lp\upperrt L} = \clspan {\upperrt L p\upperrt L p\upperrt L} \={AHeredInL}
\clspan {\upperrt L\upperrt A \upperrt L}, $$ proving our claim.  Given $g\in G$, $$ x = \mat a, m; n^*, b; \in L_1 = L
\and x'= \mat a', m'; n'^*, b'; \in L_g, $$ we have $$ (x\delta _1)(\E 11\delta _1)(x'\delta _g) \={MultOnSDPBun} x\E
11x'\delta _g = \mat aa', am'; n^*a',\ip Bn{m'};\delta _g.  $$ This implies that the ideal generated by $A\redrt \alpha
G$ contains $$ \mat \clspan {AA_g}, \clspan {AM_g}; \Clspan {(A_gM)^*}, \pilar {12pt}\Clspan {\ip B M{M_g}};\delta _g.
$$ Observing that $$ \clspan {AA_g} = A_g,\qquad \clspan {AM_g} = M_g, $$ $$ M_g \={ApproxIdentityHMod} \clspan {\ip
A{M_g}{M_g}M_g} \subseteq \clspan {A_gM}, $$ and $$ B_g \={FullIdealsInMor} \Clspan {\ip B {M_g}{M_g}} \subseteq \Clspan
{\ip B M{M_g}}, $$ we therefore conclude that the ideal generated by $A\redrt \alpha G$ contains $L_g\delta _g$, for
every $g$, and consequently also $L\redrt \lambda G$.

The grand conclusion is that $A\redrt \alpha G$ is a full hereditary subalgebra of $L\redrt \lambda G$, so these are
\RME lent C*-algebras by \cite {FullHeredMorita}.

We may now rerun the whole argument above in order to prove $B\redrt \beta G$ to be \RME lent to $L\redrt \lambda G$.
Since \RME lence is well known to be an equivalence relation, the proof of (i) is concluded.

Focusing now on (ii), let us consider the semi-direct product bundles $\Abun $ and $\L $ associated to the partial
actions $\alpha $ and $\lambda $, respectively.  We claim that $\Abun $ and $\L $ satisfy the hypothesis of \cite
{FullInjCondExp}, namely that there exists a conditional expectation from $\L $ onto $\Abun $.

Since each fiber of $\L $ is faithfully represented in $\Cr \L $ by \cite {FinallyInjective/ii}, we will produce the
conditional expectation needed by working within $\Cr \L $, or equivalently, within $L\redrt \lambda G$.  We then claim
that the mapping $$ E: x\in L\redrt \lambda G \mapsto (\E 11\delta _1) x(\E 11\delta _1)\in A\redrt \alpha G, $$ sends
each fiber $L_g\delta _g$ to the corresponding $A_g\delta _g$, and that the restrictions $P_g$ thus defined form a
conditional expectation $P = \{P_g\}_{g\in G}$ from $\L $ to $\Abun $.

An expression for $P_g$ may be obtained from our previous calculation \cite {CondexpLinkAlg}, namely $$ P_g\left (\mat
a, m; n^*, b; \delta _g\right ) = \mat a, 0; 0, 0; \delta _g.  $$ The verification that the $P_g$ satisfy \cite
{DefineConExBun/i--iii} is now straightforward and is left as an exercise.

As claimed, the hypothesis of \cite {FullInjCondExp} is now verified, so we conclude that the canonical embedding of
full cross-sectional algebras is a monomorphism, which we interpret as $$ A\rt \alpha G \subseteq L\rt \lambda G.  $$
The proof of (i) now generalizes ipsis literis, after replacing reduced crossed products by their full versions.
\endProof

Returning to the context of table \cite {Tabela}, recall that $\Th $ is \RME lent to $\rTh $ by \cite {ExistMorGlob}.
Thus, by the above result, in each one of the last two rows of table \cite {Tabela}, the algebras in the first and
second columns are \RME lent.

In order to include the algebras in the third column, we present the following:

\state Theorem \label GlobaMorita Let $$ \Th = \big (A,\ G,\ \{\D g\}_{g\in G},\ \{\th g\}_{g\in G}\big ) $$ be a
C*-algebraic partial dynamical system admitting a globalization $\eta $, acting on a C*-algebra $B$.  Then: \izitem
\zitem $A\redrt \Th G$ is isomorphic to a full hereditary subalgebra of $B\redrt \eta G$, in a natural way, and hence
$A\redrt \Th G$ and $B\redrt \eta G$ are \RME lent.  \zitem $A\rt \Th G$ is isomorphic to a full hereditary subalgebra
of $B\rt \eta G$, in a natural way, and hence $A\rt \Th G$ and $B\rt \eta G$ are \RME lent.

\Proof By \cite {RestrGivesSubBun} we have that $A\redrt \Th G$ is a hereditary subalgebra of $B\redrt \eta G$, so let
us now prove that $A\redrt \Th G$ is a full subalgebra.  This is to say that the closed two-sided ideal generated by
$A\redrt \Th G$, say $$ J:=\Clspan {(B\redrt \eta G)(A\redrt \Th G)(B\redrt \eta G)}, $$ coincides with $B\redrt \eta
G$.  For this pick $g$ and $h$ in $G$, and notice that $$ J\supseteq \Clspan {(B\delta _h)(A\delta _1)(B\delta _{h\inv
g})} = \Clspan {B\eta _h\big (AB)\delta _g} = \eta _h\big (A)\delta _g, $$ so $J$ also contains $$ \Big (\,\barsum
{\soma {h\in G}\eta _h\big (A)}\,\Big )\delta _g = B\delta _g, $$ for every $g$ in $G$, from where it follows that
$J=B\redrt \eta G$, proving that $A\redrt \Th G$ is indeed full.

As for the second point, we have by \cite {RestrGivesFullCP} that $A\rt \Th G$ is naturally isomorphic to a hereditary
subalgebra of $B\rt \eta G$, and the above proof that $A\redrt \Th G$ is a full subalgebra of $B\redrt \eta G$
generalizes ipsis literis to show that $A\rt \Th G$ is a full subalgebra of $B\rt \eta G$.

The statements about \RME lence now follow immediately from \cite {FullHeredMorita}.  \endProof

Given a partial dynamical system $$ \Th = (A,G,\{\D g\}_{g\in G},\ \{\th g\}_{g\in G}), $$ we may summarize the main
results obtained in this chapter as follows:

\begingroup \parindent 14pt \parskip 4pt \def \bit {\item {$\bullet $}}

\bit When a globalization for $\Th $ exists, it is unique. See \cite {UniqueGlobaCstar}.

\bit Regardless of whether or not $\Th $ admit a globalization, there always exists a partial dynamical system $\rTh $,
\RME lent to $\Th $, which admits a globalization $\rEta $. See \cite {ExistMorGlob}.

\bit The reduced crossed products relative to the partial dynamical systems $\Th $, $\rTh $, and $\rEta $ are all \RME
lent to each other. See \cite {EqPAEquivCP/i} and \cite {GlobaMorita/i}.

\bit The full crossed products relative to the partial dynamical systems $\Th $, $\rTh $, and $\rEta $ are all \RME lent
to each other.  See \cite {EqPAEquivCP/ii} and \cite {GlobaMorita/ii}.

\endgroup

\bigskip There are a few other important results related to the globalization of C*-algebraic partial actions which we
would now like to mention without proofs.

When a partial action $\Th $ is \RME lent to another partial action $\rTh $, which, in turn, admits a globalization
$\rEta $, one says that $\rEta $ is a \subj {\RiefMor globalization} of $\Th $.

\medskip The following is \ref {Proposition 6.3/Abadie/2003}:

\state Theorem \label UniqueMoritaGloba Any two \RiefMor globalizations of a given C*-al\-gebraic partial action are
\RME lent to each other.

Another very interesting result proved by Abadie regards partial actions on abelian C*-algebras.  In order to describe
it, let $\Th $ be a partial action of a group $G$ on a LCH space $X$ and denote by $\Th '$ the partial action of $G$ on
$C_0(X)$ corresponding to $\Th $ via \cite {TopoCstarActions}.

By \cite {TopoGlob}, we have that $\Th $ always admits a globalization, say $\eta $, acting on a space $Y$.  However $Y$
may be non-Hausdorff, and in this case $\Th '$ admits no globalization by \cite {ComutGloba}.  So let $\beta $ be a
\RiefMor globalization of $\Th '$, acting on a C*-algebra $B$, which always exists by \cite {ExistMorGlob}.

The action $\beta $ evidently induces an action $\hat \beta $ on the primitive ideal space $\Prim (B)$, and because any
two \RiefMor globalizations are equivalent by \cite {UniqueMoritaGloba}, $\hat \beta $ does not depend on the specific
choice of $\beta $, being thus intrinsically associated to $\Th $.

The next result is essentially \ref {Proposition 7.4/Abadie/2003}:

\state Theorem Let $\Th $ be a topological partial action of a group $G$ on a LCH space $X$, and let $\beta $, acting on
$B$, be a \RiefMor globalization of the corresponding partial action on $C_0(X)$.  Then the action $\hat \beta $ of $G$
on $\Prim (B)$ is the (possibly non-Hausdorff) globalization of $\Th $.

\nrem Most of the results in this chapter are extracted from Abadie's PhD Thesis \ref {Abadie/1999}, \ref {Abadie/2003},
where in fact the more general case of partial actions of locally compact groups is considered.  See also \ref
{AbadieMarti/2009} for a study of the relationship between the amenability of a partial action and that of its Morita
enveloping action.

\chapter Topologically free partial actions

\def \Cz #1{C_0(\D {#1})}

In chapter \cite {IdealsChapter} we have studied induced and Fourier ideals in topologically graded algebras.  These may
be considered as the ideals which respect the given grading.  In the present chapter we will study conditions which
imply that all ideals somehow take the grading into account, such as having a nonzero intersection with the unit fiber
algebra or, in the best-case scenario, are induced ideals.

The graded algebra at the center of the stage here will be the crossed product of an abelian C*-algebra $A$ by a given
group $G$, and the main conditions we will impose relate to the corresponding partial action of $G$ on the spectrum of
$A$, namely that there are not too many fixed points.  As an application we will give conditions for a partial crossed
product to be simple.

\definition Let $\Th = \big (\{\D g\}_{g\in G}, \{\th g\}_{g\in G}\big )$ be a topological partial action of a group $G$
on a locally compact Hausdorff space $X$, and let $g\in G$.  \iaitem \aitem A \subj {fixed point} for $g$ is any element
$x$ in $\Di g$, such that $\th g(x)=x$.  \aitem The set of all fixed points for $g$ will be denoted by $F_g$.  \aitem We
say that $\Th $ is \subjex {free}{free partial action} if $F_g$ is empty for every $g\neq 1$.  \aitem We say that $\Th $
is \subjex {topologically free}{topologically free partial action} if the \"{interior} of $F_g$ is empty for every
$g\neq 1$.

Observe that $F_g =F_{g\inv }$, for all $g$ in $G$.  Also notice that $F_g$ is a closed subset of $\Di g$ (and hence
also of $\D g$ by the previous remark) in the relative topology, but it is not necessarily closed in $X$.  Incidentally,
here is a useful result of general topology:

\state Lemma \label CloInOp Let $X$ be a topological space, and suppose we are given subsets $F\subseteq D\subseteq X$,
such that $F$ is closed relative to $D$, and $D$ is open in $X$. If the interior of $F$ is empty, then the interior of
$\bar F$ is also empty, where $\bar F$ denotes the closure of $F$ relative to $X$.

\Proof We should observe that, since $D$ is open, the interior of $F$ relative to $D$ is the same as the interior of $F$
relative to $X$.

It suffices to prove that, if $V$ is an open subset of $X$, with $V\subseteq \bar F$, then $V=\emptyset $.  Notice that
for each such $V$, we have $$ V\cap D \subseteq \bar F \cap D = F, $$ the last equality being the expression that $F$ is
closed in $D$.  Since $F$ has no interior, we deduce that $V\cap D = \emptyset $.

It is well known that, when two open subsets are disjoint, each one is necessarily disjoint from the other's closure,
whence $$ \emptyset = V \cap \bar D \supseteq V \cap \bar F = V.  \omitDoubleDollar \endProof

Another topological fact we will need is the following baby version of Baire's category Theorem, which incidentally
holds in any topological space:

\state Lemma \label BabyBaire Let $X$ be a topological space and let $F_1, F_2, \ldots , F_n$ be nowhere dense\fn {A
subset of a topological space is said to be \"{nowhere dense} when its closure has empty interior.}  subsets of $X$.
Then $F_1 \cup F_2 \cup \ldots \cup F_n$ is nowhere dense.

\Proof Replacing each $F_i$ by its closure, we may clearly assume them to be closed, and hence so is their union.

In order to prove the statement we must show that every open set $$ V\subseteq F_1 \cup F_2 \cup \ldots \cup F_n $$ is
necessarily empty.  Given such a $V$, notice that, $$ W:= V \setminus (F_1 \cup F_2 \cup \ldots \cup F_{n-1}) $$ is an
open set contained in $F_n$, whence $W$ is empty by hypothesis.  Thus $ V \subseteq F_1 \cup F_2 \cup \ldots \cup
F_{n-1}, $ and the conclusion follows by induction.  \endProof

\medskip \fix From now on we will fix a group $G$, a locally compact Hausdorff space $X$ and a topological partial
action $$ \Th = \big (\{\D g\}_{g\in G}, \{\th g\}_{g\in G}\big ) $$ of $G$ on $X$.  Our attention will be focused on
the partial action of $G$ on $C_0(X)$ induced by $\Th $ via \cite {TopoCstarActions}, which we will henceforth denote by
$\alpha $.  More precisely $$ \alpha = \Big (\big \{\Cz g\big \}_{g\in G}, \{\alpha _ g\}_{g\in G}\Big ), $$ where each
$\alpha _g$ is given by $$ \alpha _g : f\in C_0(\D {g\inv }) \mapsto f\circ \thi g \in \Cz g.  $$

Above all, we are interested in the reduced crossed product $$ C_0(X)\redrt \alpha G, $$ so, given $f$ in $\Cz g$, we
will always interpret the expression ``$f\delta _g$'' as an element of the above crossed product algebra.

The following technical Lemma is the key tool in the proof of our main result below.  It is intended to \"{shut out}
certain elements in the grading subspaces $ \Cz g\delta _g, $ with $g\neq 1$, by \"{compressing} them away with positive
elements lying in the unit fiber algebra.

\state Lemma \label hlemma Given \iaitem \aitem $g\in G$, with $g\neq 1$, \aitem $f \in \Cz g$, \aitem $x_0\in
X\setminus F_g$, and \aitem $\varepsilon >0$, \medskip \noindent there exists $h$ in $C_0(X)$ such that \izitem \zitem
$h(x_0) = 1$, \zitem $ 0 \leq h \leq 1$, and \zitem $\Vert (h \delta _1) ( f \delta _g) (h\delta _1) \Vert \leq
\varepsilon $.

\Proof We separate the proof in two cases, according to whether or not $x_0$ lies in $\D g$.  If $x_0$ is not in $\D g$,
let $$ K= \{x\in \D g: |f(x)| \geq \varepsilon \}.  $$ Then $K$ is a compact subset of $\D g$ and, since $x_0$ is not in
$K$, we may use Urysohn's Lemma obtaining a continuous, real valued function on $X$ such that $$ 0 \leq h \leq 1,\quad
h(x_0) = 1 \and h|_K = 0.  $$

By temporarily working in the one-point compactification of $X$, and replacing $K$ by $K\cup \{\infty \}$, we may also
assume that $h$ vanishes at infinity, meaning that $h$ is in $C_0(X)$.

Since $f$ is bounded by $\varepsilon $ outside of $K$, it follows that $\Vert h f \Vert \leq \varepsilon $, from where
(iii) easily follows, hence concluding the proof in the present case.

If $x_0$ is in $\D g$, then $\theta _{g\inv }(x_0)$ is defined and not equal to $x_0$.  We may then take disjoint open
sets $V_1$ and $V_2$ such that $$ x_0 \in V_1 \and \theta _{g\inv }(x_0) \in V_2, $$ and we may clearly assume that $V_1
\subseteq \D g$, and $V_2 \subseteq \D {g\inv }$.  We may further shrink $V_1$ by replacing it with $$ V_1' = V_1 \cap
\theta _g(V_2), $$ (notice that $x_0$ is still in $V'_1$), so that $\theta _{g\inv } (V_1') \subseteq V_2$, and then $$
\theta _{g\inv } (V_1') \cap V_1' = \emptyset .  \equationmark VoneVtwoDisj $$ Using Urysohn's Lemma again, pick $h$ in
$C_0(X)$ such that $$ 0 \leq h \leq 1, \ h(x_0) = 1 \and h|_{X\setminus V_1'} = 0.  $$ We then have $$ (h \delta _1) ( f
\delta _g) (h\delta _1) = \alpha _g\big ( \alpha _{g\inv } (h f) h \big ) \delta _g =0, $$ because $\alpha _{g\inv } (h
f)$ is supported $\theta _{g\inv }(V_1')$, while $h$ is supported in $V_1'$, and these are disjoint sets by \cite
{VoneVtwoDisj}.  This verifies (iii) in the final case, and the proof is thus concluded.  \endProof

The following result, together with its consequences to be presented below, is gist of this chapter:

\state Theorem \label topfree Let $\Th $ be a topologically free partial action of a group $G$ on a locally compact
Hausdorff space $X$.  Then any nonzero closed two-sided ideal $$ J \ideal C_0(X) \redrt \alpha G $$ has a nonzero
intersection with $C_0(X)$.

\Proof We should notice that the last occurrence of $C_0(X)$, above, stands for its copy $C_0(X)\delta _1$ within
$C_0(X) \redrt \alpha G.$

Denoting the conditional expectation provided by \cite {Fourier} by $$ E_1: C_0(X) \redrt \alpha G \to C_0(X), $$ we
claim that, for every $z$ in $C_0(X) \redrt \alpha G$, and for every $\varepsilon >0$, there exists $h \in C_0(X)$, such
that \izitem \zitem $ 0 \leq h \leq 1$, \zitem $\Vert h E_1(z) h \Vert \geq \Vert E_1(z) \Vert - \varepsilon $, and
\zitem $\Vert (h\delta _1) z (h\delta _1) - h E_1(z) h\delta _1 \Vert \leq \varepsilon $.

\medskip Assume first that $z$ is a linear combination of the form $$ z=z_1\delta _1 + \soma {g\in T} z_g \delta _g,
\equationmark SpecialZ $$ where $T$ is a finite subset of $G$, with $1\not \in T$, in which case $E_1(z) = z_1$.  Let $$
V = \{x \in X: |z_1(x)| > \Vert z_1\Vert - \varepsilon \}, $$ which is clearly open and nonempty.  By the topological
freeness hypothesis and by \cite {CloInOp}, each $F_g$ is nowhere dense in $X$.  Furthermore, by \cite {BabyBaire} one
has that $\bigcup _{g\in T} F_g$ is likewise nowhere dense, hence there exists some $$ x_0 \in V \setminus \Big (\segura
{\bigcup }_{g\in T}F_g\Big ).  $$

For each $g$ in $T$ we may then apply \cite {hlemma}, obtaining an $h_g$ in $C_0(X)$ satisfying $$ h_g(x_0) = 1, \quad 0
\leq h_g \leq 1 \and \Vert (h_g\delta _1) (z_g \delta _g) (h_g\delta _1) \Vert \leq {\varepsilon \over |T|}.  $$

Here we are tacitly assuming that $|T|>0$, and we leave it for the reader to treat the trivial case in which $|T|=0$.

We will now show that $ h := \segura {\prod }_{g\in T} h_g, $ satisfies conditions (i--iii), above.  Noticing that (i)
is immediate, we prove (ii) by observing that $x_0$ is in $V$, so $$ \Vert h E_1(z) h \Vert = \Vert h z_1 h \Vert \geq
|z_1(x_0)| > \Vert z_1 \Vert - \varepsilon = \Vert E_1(z) \Vert - \varepsilon .  $$ As for (iii), we have $$ \Vert
(h\delta _1) z (h\delta _1) - h E_1(z) h\delta _1 \Vert = \Vert (h\delta _1) \big (z-E_1(z)\delta _1\big ) (h\delta
_1)\Vert \$= \normsum {\soma {g \in T} (h\delta _1) (z_g \delta _g) (h\delta _1) } \leq \soma {g \in T} \Vert (h\delta
_1) (z_g \delta _g) (h\delta _1) \Vert \$\leq \soma {g \in T} \Vert (h_g\delta _1) (z_g \delta _g) (h_g\delta _1) \Vert
\leq \varepsilon .  $$ This proves (i--iii) under special case \cite {SpecialZ}, but since the elements of that form are
dense in $C_0(X)\redrt \alpha G$, a standard approximation argument gives the general case.

The claim verified, let us now address the statement.  Arguing by contradiction we suppose that $J$ is a nonzero ideal
such that $$ J \cap C_0(X) =\{0\}.  $$

Pick a nonzero element $y$ in $J$, and let $ z = y^*y.  $ Using the claim, for each positive $\varepsilon $ we choose
$h$ in $C_0(X)$, satisfying (i--iii) above.  Let $$ q \ : \ C_0(X) \redrt \alpha G \ \to \ {C_0(X) \redrt \alpha G \over
J} $$ be the quotient map.  Since $z$ lies in $J$, we have that $$ q\big ((h\delta _1) z (h\delta _1)\big )=0, $$ whence
$$ \big \Vert q \big ( h E_1(z) h\delta _1\big ) \big \Vert = \big \Vert q \big ( (h\delta _1) z (h\delta _1) - h E_1(z)
h\delta _1\big ) \big \Vert \$\leq \big \Vert (h\delta _1) z (h\delta _1) - h E_1(z) h\delta _1 \big \Vert \explica \leq
{iii} \varepsilon .  $$

Since $J\cap C_0(X)=\{0\}$, we deduce that $q$ is injective, hence isometric, on $C_0(X)$.  So $$ \varepsilon \geq \Vert
h E_1(z) h \Vert \explica \geq {ii} \Vert E_1(z) \Vert - \varepsilon , $$ from where we see that $\Vert E_1(z) \Vert
\leq 2 \varepsilon $, and since $\varepsilon $ is arbitrary, we have $$ 0 = E_1(z) = E_1(y^*y).  $$ It then results from
\cite {EOneIsFaith} that $y=0$, a contradiction.  This concludes the proof.  \endProof

A useful consequence to the representation theory of reduced crossed products is as follows:

\state Corollary \label CoroTopFree Let $\Th $ be a topologically free partial action of a group $G$ on a locally
compact Hausdorff space $X$.  Then a *-representa\-tion of $C_0(X) \redrt \alpha G$ is faithful if and only if it is
faithful on $C_0(X)$.

\Proof Apply \cite {topfree} to the kernel of the given representation.  \endProof

\definition We will say that $\theta $ is a \subj {minimal partial action} if there are no $\theta $-invariant closed
subsets of $X$, other than $X$ and the empty set.

The complement of an invariant set is invariant too, so minimality is equivalent to the absence of nontrivial \"{open}
invariant subsets.

\state Corollary \label TofFreeMinimal If, in addition to the conditions of \cite {topfree}, $\Th $ is minimal, then
$C_0(X) \redrt \alpha G$ is a simple\fn {A C*-algebra is said to be simple if there are no nontrivial closed two-sided
ideals.} C*-algebra.

\Proof Let $J$ be a nonzero, closed two-sided ideal of $C_0(X) \redrt \alpha G$.  Employing \cite {topfree}, we have
that $$ K := J\cap C_0(X) \neq \{0\}.  \equationmark KNonZero $$ Since $K$ is an ideal in $C_0(X)$, there is an open
subset $U\subseteq X$, such that $K = C_0(U)$.  Using \cite {IntersInvar} we have that $K$ is $\alpha $-invariant, and
it is easy to see that this implies $U$ to be $\Th $-invariant.  By minimality of $\Th $, one has that either
$U=\emptyset $, or $U=X$, but under the first case we would have $K=\{0\}$, contradicting \cite {KNonZero}.  So $U=X$,
and hence $K=C_0(X)$, meaning that $C_0(X)\subseteq J$.

By \cite {BOneEssential}, which is stated for the full cross-sectional algebra, but which evidently also holds for the
reduced one, we have that $C_0(X)$ (the unit fiber algebra in the semi-direct product bundle) generates $C_0(X) \redrt
\alpha G$, as an ideal, whence $J=C_0(X) \redrt \alpha G$, and the proof is concluded.  \endProof

The upshot of this result is that when only two $\Th $-invariant open subsets exist in $X$, namely $\emptyset $ and the
whole space, then only two closed two-sided ideals exist in $C_0(X) \redrt \alpha G$, namely $\{0\}$ and the whole
algebra.

This raises the question as to whether a correspondence may still be found between invariant open subsets and closed
two-sided ideals, in case these exist in greater numbers.  In order to find such a correspondence we need a bit more
than topological freeness.

\state Theorem \label ClassIdeal Let $\Th $ be a topological partial action of a group $G$ on a locally compact
Hausdorff space $X$, such that either: \izitem \zitem $G$ is an exact group, or \zitem the semi-direct product bundle
associated to $\Th $ satisfies the approximation property.  \medskip \noindent In addition we suppose that $\Th $, as
well as the restriction of\/ $\Th $ to any closed $\Th $-invariant subset of $X$, is topologically free.  Then there is
a one-to-one correspondence between $\Th $-invariant open subsets $U\subseteq X$, and closed two-sided ideals in $C_0(X)
\redrt \alpha G$, given by $$ U\mapsto C_0(U)\redrt \alpha G.  $$

\Proof Given a $\Th $-invariant open subset $U\subseteq X$, it is easy to see that $C_0(U)$ is an $\alpha $-invariant
ideal in $C_0(X)$.  By \cite {ExactSeqCross}, one then has that $C_0(U)\redrt \alpha G$ is an ideal in $C_0(X)\redrt
\alpha G$.

By first checking on elements of the algebraic crossed product, it is easy to see that the conditional expectation $E_1$
of \cite {Fourier} satisfies $$ E_1\big (C_0(U)\redrt \alpha G\big ) = C_0(U)\delta _1.  $$ From this it follows that
our correspondence is injective, and we are then left with the task of proving that any ideal $$ J\ideal C_0(X)\redrt
\alpha G $$ is of the above form.  Given $J$, let $U\subseteq X$ be the unique open subset of $X$ such that $$ J\cap
C_0(X) = C_0(U).  $$

By \cite {IntersInvar} we have that $C_0(U)$ is $\alpha $-invariant, whence $U$ is $\Th $-invariant and hence so is $$ F
:= X\setminus U.  $$

It is a well known fact that the map $$ q: f \in C_0(X) \mapsto f|_F\in C_0(F) $$ passes to the quotient modulo
$C_0(U)$, leading up to a *-isomorphism $$ C_0(X)/C_0(U) \simeq C_0(F).  $$

Observe that $C_0(X)/C_0(U)$ carries a quotient partial action of $G$ as proved in \cite {QuotientPAct}, while the
restriction of $\Th $ to $F$ induces a partial action of $G$ on $C_0(F)$, via \cite {TopoCstarActions}.  The isomorphism
mentioned above may then be easily shown to be $G$ equivariant.  So, employing \cite {ExactSeqCross} we get the
following exact sequence of C*-algebras and *-homomorphisms $$ 0\ \to \ C_0(U)\redrt \null G \ \> {\iota \lred } \
C_0(X)\redrt \alpha G \ \>{q\lred } \ C_0(F)\redrt \beta G \ \to \ 0.  $$

By \cite {BOneEssential} notice that $C_0(U)\redrt \null G$ coincides with the ideal generated by $C_0(U)$ within
$C_0(X)\redrt \alpha G$.  Since $C_0(U)\subseteq J$, we then deduce that $$ C_0(U)\redrt \null G\subseteq J.
\equationmark ContJota $$ We next claim that $$ q\lred (J) \cap C_0(F) = \{0\}.  \equationmark ClaimqJCzF $$ In order to
see this, let $z$ be in $q\lred (J) \cap C_0(F)$, and choose $y$ in $J$ such that $q\lred (y)=z$.  Since $C_0(F)=q\lred
\big (C_0(X)\big )$, we may also choose $f$ in $C_0(X)$ such that $q\lred (f) = z$.  Therefore $q\lred (f-y)=0$, hence
we see that $$ f-y \in \Ker (q\lred ) =C_0(U)\redrt \null G\subseteq J, $$ from where it follows that $$ f \in J\cap
C_0(X) = C_0(U), $$ so $$ 0 = q\lred (f) = z.  $$

This proves \cite {ClaimqJCzF}, and since the restriction of $\Th $ to $F$ is topologically free by hypothesis, we may
use \cite {topfree} to deduce that $q\lred (J) = \{0\}$, which is to say that $$ J \subseteq \Ker (q\lred )
=C_0(U)\redrt \null G $$ so the inclusion in \cite {ContJota} is in fact an equality of sets.  This concludes the proof.
\endProof

As a consequence we have:

\state Corollary Under the hypotheses of Theorem \cite {ClassIdeal}, every ideal of \/ $C_0(X) \redrt \alpha G$ is
induced, and hence also a Fourier ideal.

\Proof In view of \cite {ClassIdeal} it is enough to prove that the ideal $$ J =C_0(U) \redrt \alpha G $$ is induced for
every $\Th $-invariant open subset $U\subseteq X$.  Given any such $U$, we need to prove that $J$ is generated by $J\cap
C_0(X)\delta _1$.  However, since it is clear that $$ C_0(U)\delta _1\subseteq J\cap C_0(X)\delta _1, $$ we deduce that
$$ J = C_0(U) \redrt \alpha G \={BOneEssential} \langle C_0(U)\delta _1\rangle \subseteq \langle J\cap C_0(X)\delta
_1\rangle \subseteq J, $$ so $J$ is indeed an induced ideal.  By \cite {InduIsFourier} we then also have that $J$ is a
Fourier ideal.  \endProof

\nrem Theorem \cite {topfree} first appeared in \ref {Theorem 2.6/ExelLacaQuigg/1997}.  It is a direct generalization of
\ref {Theorem 4.1/KawamuraTomiyama/1990} for partial actions.  See also \ref {Theorem 4.4/Exel/2011}.

A groupoid version of Theorem \cite {ClassIdeal} is stated in \ref {Proposition 4.5/Renault/1980}, although there is a
missing hypothesis in it, without which one might have the bad behavior discussed in \ref {Remark 4.10/Renault/1991}.
The correct statement applies for amenable groupoids and is to be found in \ref {Corollary 4.9/Renault/1991}.

Since a notion of \"{exact groupoid} does not seem to exist, Theorem \cite {ClassIdeal} under hypothesis (i) does not
appear to have a groupoid counterpart.

Many authors have studied results similar to the ones we proved above for global actions of discrete groups on
non-abelian C*-algebras.  See, for example, \ref {Elliott/1980}, \ref {KawamuraTomiyama/1990} and \ref
{ArchboldSpielberg/1993}.

\partpage {III}{APPLICATIONS}

\chapter Dilating partial representations

\def \Cstarpar {C^*\underpar (G)} \def \t {\tilde }

With this chapter we start a series of applications of the theory developed so far.

Recall from \cite {CompressPrep} that if $v$ is a unitary group representation of $G$, and if $p$ is a self-adjoint
idempotent such that $v_gpv_{g\inv }$ commutes with $p$, for every $g$ in $G$, then $$ \pr g := pv_gp \for g\in G $$
defines a partial representation of $G$.

One may think of $\Pr $ as a restriction of $v$, a process akin to restricting a partial action to a not necessarily
invariant subset.  In this chapter we will study a situation in which one may revert this process, proving that partial
representations of a group on a Hilbert space may always be obtained from unitary representations via the above process.

\definition In the context of \cite {CompressPrep}, we will say that $\Pr $ is the \subjex {restriction}{restriction of
a unitary representation} of $v$ and, conversely, that $v$ is a \subjex {dilation}{dilation of a partial representation}
of \Pr .

In order to prepare for the proof of our main dilation result, we will now prove a related dilation Lemma, applicable to
*-representations of C*-algebras.

\state Lemma \label ExtendRep Let $B$ be a C*-algebra and $A$ be a closed *-subalgebra of $B$.  Given a *-representation
$\pi $ of $A$ on a Hilbert space $H$, there exists a representation $\t \pi $ of $B$ on another Hilbert space $\t H$,
and a (not necessarily surjective) isometric linear operator $j :H\to \t H$, such that $$ j \pi (a) = \t \pi (a)j \for
a\in A.  $$

\Proof By splitting $\pi $ in the direct sum of cyclic representations (perhaps including an identically zero
sub-representation), we may assume, without loss of generality, that $\pi $ is cyclic.  So, let $\xi $ be a cyclic
vector for $\pi $.  Defining $$ \varphi (a) = \langle \pi (a)\xi ,\xi \rangle \for a\in A, $$ we get a positive linear
functional on $A$.  Let $\t \varphi $ be a positive linear functional on $B$ extending $\varphi $ \ref {Proposition
3.1.6/Pedersen/1979}, and let $\t \pi $ be the GNS representation of $B$ associated to $\t \varphi $, acting on a
Hilbert space $\t H$, with cyclic vector $\t \xi $.

We then claim that there exists an isometric linear operator $j:H\to \t H$, such that $$ j \big (\pi (a)\xi \big ) = \t
\pi (a)\t \xi \for a\in A.  $$ To see this, first define $j $ on the dense subspace $\pi (A)\xi \subseteq H$ by the
above formula.  This is well defined and isometric because for every $a$ in $A$, we have $$ \Vert \t \pi (a)\t \xi \Vert
^2 = \langle \t \pi (a)\t \xi ,\t \pi (a)\t \xi \rangle = \langle \t \pi (a^*a)\t \xi ,\t \xi \rangle = \t \varphi
(a^*a) \$= \varphi (a^*a) = \langle \pi (a^*a)\xi ,\xi \rangle = \Vert \pi (a)\xi \Vert ^2.  $$

Extending $j $ by continuity to the whole of $H$, we clearly obtain the claimed operator.  So, in order to conclude the
proof, it is enough to verify the identity in the statement.

It $\eta $ is a vector in $H$ of the form $\eta =\pi (c)\xi $, with $c$ in $A$, then, for every $a$ in $A$ we have $$
j\big (\pi (a)\eta \big ) = j\big (\pi (a)\pi (c)\xi \big ) = j \big (\pi (ac)\xi \big ) = \t \pi (ac)\t \xi \$= \t \pi
(a)\t \pi (c)\t \xi = \t \pi (a)j\big (\pi (c)\xi \big )= \t \pi (a)j(\eta ).  $$

Since the set of $\eta $'s considered above is dense in $H$, the proof is concluded.  \endProof

We should notice that, in the context of the above result, the range of $j$ is easily seen to be invariant under $\t \pi
(A)$.  So the projection onto the range of $j$, namely $jj^*$, commutes with $\t \pi (A)$.  Consequently, if we let
$\rho $ be the restriction of $\t \pi $ to $A$, then the sub-representation of $\rho $ determined by the range of $j$ is
clearly equivalent to $\pi $, the equivalence being implemented by $j$.  Another useful remark is that, since $j$ is an
isometry, we have that $j^*j$ coincides with the identity on $H$, so $$ j^*\t \pi (a)j = j^*j\pi (a) = \pi (a), $$ for
every $a$ in $A$.

We may now prove our main result on dilation of partial representations.

\state Theorem \label AbadieOnDilation Let $\Pr $ be a \"{partial} representation of a group $G$ on a Hilbert space $H$.
Then there exists a \"{unitary} group representation $\t \Pr $ of $G$ on a Hilbert space $\t H$, and an isometric linear
operator $j : H\to \t H$, such that $$ \pr g = j^*\t \pr gj \for g\in G.  $$ In addition, is $p$ denotes the projection
of $\t H$ onto the range of $j$, namely $p=jj^*$, then $p$ commutes with $\t \pri g p \t \pr g$ for every $g$ in $G$.

\Proof \def \can #1{w_{#1}} We should first observe that, if we identify $H$ with its image via $j$, hence thinking of
$j$ as the inclusion map, then, in line with \cite {CompressPrep}, the formula displayed in the statement becomes $$ \pr
g = p\t \pr g p.  $$ However, in this proof we will not enforce the identification of $H$ as a subspace of $\t H$.

Taking $\Rel $ to be the empty set of relations in \cite {UniversalForPreps} (see also \cite {DefineCStarPar} and \cite
{UnivPRep}), we conclude that there exists a *-representation $\pi $ of $\Cstarpar $ on $H$ such that $$ \pi (\can g) =
\pr g \for g\in G, $$ where the $\can g$ denote the canonical generators of $\Cstarpar $ (please notice that these are
denoted by $\pr g$ in \cite {UniversalForPreps}).  In particular $$ \pi (1) = \pi (\can 1) = \pr 1 = 1, $$ so $\pi $ is
a {\nonDegRep } representation.

Recall from \cite {CstarParAsCP} that $\Cstarpar $ is isomorphic to the crossed product algebra relative to the partial
Bernoulli action $\beta $ of $G$ on $\OuG $.  By definition, this is the restriction of the global Bernoulli action
$\eta $, introduced in \cite {DefineGlobalBernouli}, to $\OuG $.  So we conclude from \cite {RestrGivesFullCP} that the
canonical mapping $$ \iota : \Cstarpar = C(\OuG )\rt \beta G \ \longrightarrow \ C\big (\OG \big )\rt \beta G $$ is
injective, hence we may view $\Cstarpar $ as a subalgebra of $C\big (\OG \big )\rt \beta G$.

Notice that $\OuG $ is a clopen subset of $\OG $, whence the characteristic function of the former, which we denote by
$1_{\OuG }$, lies in $C\big (\OG \big )$.

Given that $\can g$ identifies with $1_g\delta _g$ by \cite {CstarParAsCP}, the following identity, to be proved in the
sequel, will be useful: $$ \iota (1_g\delta _g) = 1_{\OuG }\t \delta _g 1_{\OuG }, \equationmark AbadieRelation $$ where
$\t \delta _g$ denoted the canonical unitary element of $C\big (\OG \big )\rt \beta G$ implementing the global Bernoulli
action.  To see this we compute $$ 1_{\OuG }\t \delta _g 1_{\OuG } = 1_{\OuG } 1_{\eta _g(\OuG )}\t \delta _g = 1_{\OuG
\cap \eta _g(\OuG )}\t \delta _g \={FirstDescriptBernouliDoms} 1_g\t \delta _g = \iota (1_g\delta _g), $$ proving \cite
{AbadieRelation}.

We next choose a representation $\t \pi $ of $C\big (\OG \big )\rt \beta G$ on a Hilbert \break space $\t H$, and an
isometric linear operator $j : H\to \t H$, satisfying the conditions in \cite {ExtendRep} relative to $\pi $.  In
particular, the formula $$ \t \pr g = \t \pi (\t \delta _g) \for g\in G, $$ defines a unitary representation $\t \Pr $
of $G$ on $\t H$, which we will prove to satisfy all of the conditions in the statement.

As a first step, insisting on the fact that $\can g$ identifies with $1_g\delta _g$ by \cite {CstarParAsCP}, we have $$
\pr g = \pi (\can g) = \pi (1_g\delta _g) = j^* \t \pi \big (\iota (1_g\t \delta _g)\big ) j \={AbadieRelation} $$$$ =
j^* \t \pi (1_{\OuG }\t \delta _g 1_{\OuG }) j = j^* \t \pi (1_{\OuG })\t \pr g \t \pi (1_{\OuG }) j = \ldots
\equationmark ComputingPrg $$ We would now like to get rid of the two occurrences of term $\t \pi (1_{\OuG })$ above,
and to justify this we observe that $$ \t \pi (1_{\OuG }) j \={ExtendRep} j \pi (1_{\OuG }) = j, $$ since $\pi $ is a
{\nonDegRep } representation and $1_{\OuG }$ represents the unit of $C\big (\OuG \big )\rt \beta G$.  Using this, and
the identity obtained by taking stars on both sides, we deduce from \cite {ComputingPrg} that $$ \pr g = j^* \t \pr g j.
$$

In order to prove the last sentence in the statement, we compute $$ \t \pri g p \t \pr g p = \t \pri g jj^* \t \pr g
jj^*= \t \pri g j \pr g j^*= \t \pri g j \pi (1_g\delta _g) j^* \={ExtendRep} $$$$ = \t \pri g \t \pi \big (\iota
(1_g\delta _g)\big ) j j^*= \t \pi \big (\t \delta _{g\inv } 1_g\t \delta _g) j j^*= \t \pi \big (1_{g\inv }) j j^*= j
\pi \big (1_{g\inv }) j^*.  $$ This shows that the element with which we started the above computation is self-adjoint,
so $$ \t \pri g p \t \pr g p = (\t \pri g p \t \pr g p)^* = p \t \pri g p \t \pr g, $$ proving the desired
commutativity.  \endProof

\nrem We believe that Lemma \cite {ExtendRep} is well known among specialists, but we have not been able to find a
reference for it in the literature.  Theorem \cite {AbadieOnDilation} was first proved in \ref {Proposition
3.3/Abadie/2003} for the case of amenable groups.

\chapter Semigroups of isometries

\def \D {{\bf D}} \def \T {{\bf T}} \def \t {\tilde }

The Toeplitz algebra, namely the C*-algebra generated by a non-unitary isometry \ref {Coburn/1967}, is among the most
historically significant examples of C*-algebras.  After Coburn's pioneering example, many authors have studied
C*-algebras generated by sets of isometries, sometimes also including partial isometries.  In many such examples, the
generating isometries are parametrized by a semigroup.  The present chapter is therefore dedicated to studying
semigroups of isometries and their relationship to partial group representations.

Recall that the term \"{semigroup} in Mathematics always refers to a set $P$, equipped with an associative operation
(often denoted multiplicatively).

\definition \label SGIso Let $A$ be a unital C*-algebra.  A \subj {semigroup of isometries} in $A$, based on a semigroup
$P$, is a map $v:P\to A$, such that \izitem \zitem $v_p^*v_p=1$, \zitem $v_pv_q=v_{pq}$, \medskip \noindent for every
$p,q\in P$.

Notice that if $P$ has a unit, say 1, then $$ v_1 = v_1^*v_1v_1 = v_1^*v_1=1, \equationmark SGUnital $$ so $v $ must
necessarily send 1 to 1.

\nostate \label SGIPowerIsometry A simple example of a semigroup of isometries is obtained by taking any isometric
linear operator $S$ on a Hilbert space $H$, and defining $$ v :n\in {\bf N} \mapsto S^n\in \Lin (H).  $$

\nostate \label NCrosN To describe an example based on the semigroup ${\bf N}\times {\bf N}$, let us fix two isometric
linear operator $S$ and $T$ on the same Hilbert space $H$, and let us assume that $S$ and $T$ commute.  Then, defining
$$ v : (n,m)\in {\bf N}\times {\bf N} \mapsto S^n T^m \in \Lin (H), $$ we again get a semigroup of isometries on $H$.

\definition A semigroup $P$ is said to be left-cancellative provided $$ pm=pn \imply m=n, $$ for every $m$, $n$ and $p$
in $P$.

Given a left-cancellative semigroup $P$, there is a somewhat canonical example of semigroup of isometries which we would
now like to present.  Consider the Hilbert space $\ell ^2(P)$, with its usual orthonormal basis $\{e_n\}_{n\in P}$.

For each $p$ in $P$, consider the bounded linear operator $\lambda _p$ on $\ell ^2(P)$, specified by $$ \lambda _p(e_n)
= e_{pn} \for n\in P.  $$

Notice that $\lambda _p$ is an isometric operator\fn {If $P$ did not possess this property, in particular if infinitely
many $n$'s were mapped to the same element of $P$ under left multiplication by $p$, there would be no \"{bounded} linear
operator sending the $e_n$ to $e_{pn}$.}  because the map $n\mapsto pn$ is injective as a consequence of $P$ being
left-cancellative.  It is now easy to check that $\lambda $ is a semigroup of isometries in $\Lin \big (\ell ^2(P)\big
)$.

\definition \label RegularSGP We shall refer to the above $\lambda $ as the \subj {regular semigroup of isometries} of
$P$.

If $P$ is left-cancellative, observe that $\lambda _p\lambda _p^*$ is the orthogonal projection onto $$ \lambda _p\big
(\ell ^2(P)\big ) = \barspan \{e_{pn}: n\in P\}.  \equationmark RangeLambdaP $$ In particular it is easy to see that
$\lambda _p\lambda _p^*$ commutes with $\lambda _q\lambda _q^*$, for every $p$ and $q$ in $P$.

The commutativity of range projections is not always granted, even for the case of the commutative semigroup ${\bf
N}\times {\bf N}$, as the following example shows.  Let $\varphi $ be a continuous, complex valued function on the
closed unit disk $$ \D := \{z\in {\bf C}: |z|\leq 1\}, $$ which is holomorphic on the interior of $\D $.  Suppose also
that $\varphi $ is \"{unimodular}, in the sense that $|\varphi (z)|=1$, for every $$ z\in \T := \{z\in {\bf C}: |z|=1\}.
$$ It is then easy to see that the operator $$ U_\varphi :\xi \in L^2(\T ) \mapsto \varphi \xi \in L^2(\T ) $$ ($\varphi
\xi $ referring to pointwise product) is unitary.

Denoting by $H^2$ the \subj {Hardy space}\fn {The Hardy space is the closed subspace of $L^2(\T )$ spanned by the set of
all functions of the form $z\mapsto z^n$, with $n\geq 0$.}  of Classical Harmonic Analysis, one may prove that
$U_\varphi (H^2)\subseteq H^2$ (this is where we need $\varphi $ to be holomorphic on the interior of the disk), so the
restriction of $U_\varphi $ to $H^2$, henceforth denoted $$ T_\varphi : H^2\to H^2, $$ is an isometry.  This is often
referred to in the literature as the \subj {Toeplitz operator} with symbol $\varphi $.

Denoting by $P$ the \subj {Hardy projection}, namely the orthogonal projection from $L^2(\T )$ to $H^2$, the invariance
of $H^2$ under $U_\varphi $ may be expressed by the formula $$ PU_\varphi P = U_\varphi P, $$ and the Toeplitz operator
may be alternatively defined by $$ T_\varphi = PU_\varphi P.  $$

A useful consequence of this description is the following formula for the adjoint of $T_\varphi $: $$ T_\varphi ^* =
PU_\varphi ^*P = PU_{\bar \varphi }P.  \equationmark AdjToep $$

Given two unimodular holomorphic functions $\varphi $ and $\psi $, as above, it is easy to see that $T_\varphi $ and
$T_\psi $ commute.  So, in the context of \cite {NCrosN}, we get the semigroup of isometries $$ v : (n,m)\in {\bf
N}\times {\bf N} \mapsto T_\varphi ^n T_\psi ^m \in \Lin (H^2).  \equationmark SGIFromToep $$

In order to present our intended example illustrating that the range projections need not commute, let us fix the
following two unimodular holomorphic functions on the unit disk: $$ \varphi (z) = z \and \psi (z) = {z-a \over 1-\bar a
z} \for z\in \T , $$ where $a$ is a fixed complex number with $|a|<1$.

Incidentally, $\psi $ is usually called a \subj {Blaschke factor}, and it is an easy exercise to prove that it is indeed
a unimodular holomorphic function on the unit disk (the denominator vanishes only for $z=1/\bar a$, which lies outside
the disk).

Our next goal is to analyze the commutativity of the range projections $$ Q':= v_{(1,0)} v_{(1,0)}^* \and Q'':=
v_{(0,1)} v_{(0,1)}^*.  $$ Upon identifying the Hardy space with $\ell ^2({\bf N})$, according to the usual orthonormal
basis $\{e_n\}_{n\in {\bf N}}$ of $H^2$ given by $$ e_n(z)=z^n \for z\in {\bf N}, $$ we have that $v_{(1,0)}$, also
known as $T_\varphi $, becomes the \subj {unilateral shift}, in the sense that $$ T_\varphi (e_n) = e_{n+1} \for n\in
{\bf N}.  $$

The range projection $$ Q' = T_\varphi T_\varphi ^*, $$ is then the orthogonal projection onto the subspace of $H^2$
spanned by the $e_n$, with $n\in {\bf N}\setmenos \{0\}$ and, in particular, we have that $$ Q'(e_0)=0.  \equationmark
QprimeVanishEzero $$

Assuming that $Q'$ and $Q''$ commute, we have for all $j>0$ that $$ \langle T_\psi T_\psi ^*(e_0),e_j\rangle = \langle
Q''(e_0),Q'(e_j)\rangle = \langle Q'Q''(e_0),e_j\rangle \$= \langle Q''Q'(e_0),e_j\rangle \={QprimeVanishEzero} 0.  $$
It follows that all but the zeroth coefficient of $T_\psi T_\psi ^*(e_0)$ relative to the canonical basis vanish, so $$
T_\psi T_\psi ^*(e_0)=ce_0, $$ for some constant $c\in {\bf C}$.  On the other hand $$ T_\psi T_\psi ^*(e_0) \={AdjToep}
PU_\psi PPU_{\bar \psi }P(e_0) = \psi P\bar \psi .  $$

Notice that, since $\psi $ is holomorphic, its Taylor series around zero coincides with its Fourier expansion on the
unit circle: $$ \psi (z) \sim \somalim {n=0}^\infty {\psi ^{(n)}(0)\over n!} z^n = \somalim {n=0}^\infty \hat \psi (n)
z^n.  $$ In particular all of the Fourier coefficients $\widehat {(\bar \psi \,)}(n)$ vanish, for $n>0$, which implies
that $P\bar \psi = \bar \psi (0)e_0$.  Comparing the two expressions for $T_\psi T_\psi ^*(e_0)$ obtained above we
deduce that $$ ce_0 = T_\psi T_\psi ^*(e_0)= \psi P\bar \psi = \psi \bar \psi (0)e_0.  $$

This equation entails two possibilities: in case $c\neq 0$, one has that $\psi $ is equal to the constant $c\bar \psi
(0)\inv $, almost everywhere on $\T $, which is clearly not the case.  The second possibility, namely that $c=0$, in
turn implies that $\psi \bar \psi (0)=0$, whence $\psi (0)=0$.  This proves the following:

\state Proposition \label RangeToeplitzCommute Let $\varphi $ and $\psi $ be given by $$ \varphi (z) = z \and \psi (z) =
{z-a \over 1-\bar a z} \for z\in \T , $$ where $a$ is a fixed complex number with $|a|<1$.  Then the corresponding
Toeplitz operators $T_\varphi $ and $T_\psi $ are commuting isometric linear operators on Hardy's space.  Moreover their
range projections $T_\varphi T_\varphi ^*$ and $T_\psi T_\psi ^*$ commute if and only if $a=0$.

\nostate \label BadSGI The above result allows for the construction of an example of the situation mentioned above:
choosing $\psi $ to be any Blaschke factor with $a\neq 0$, the semigroup of isometries presented in \cite {SGIFromToep}
will be such that the final projections of $v_{(1,0)}$ and $v_{(0,1)}$ do not commute.

This is in stark contrast with \cite {MainTheEgCommute}, according to which the final projections of the partial
isometries involved in a partial group representation always commute!  We must therefore recognize that the theory of
partial representations of groups is unsuitable for dealing with such badly behaved examples.

We will therefore restrict our attention to special cases of semigroup of isometries which can be effectively studied
via the theory of partial representations.

If we are to find a partial representation of a group somewhere behind a semigroup of isometries, we'd better make sure
there is a group around.  So we will assume from now on that the semigroup $P$ is a sub-semigroup of a group $G$.  If we
are moreover given a semigroup of isometries $$ v :P\to A, $$ where $A$ is a unital C*-algebra, we will now discuss
conditions under which one may extend $v $ to a partial representation of $G$.

Assume for a while that the problem has been solved, i.e., that a partial representation $\Pr $ of $G$ in $A$ has been
found such that $\pr n = v_n$, for all $n$ in $P$.  Since each $v_n$ is an isometry, it is clearly left-invertible.
Therefore, from \cite {RegardDisregard} it follows that, for all $m,n\in P$, $$ \pr {m\inv n} = \pr {m\inv } \pr n = \pr
m^* \pr n = v_m^* v_n.  \equationmark PrepDetermined $$ In other words, if $g\in P\inv P$, then $\pr g$ may be recovered
from $v $.

\definition A sub-semigroup $P$ of a group $G$ is called an \subj {Ore sub-semigroup}\fn {An abstract semigroup $P$ may
be embedded in a group $G$, such that $G=P\inv P$, if and only if $P$ is cancellative and for every $m$ and $n$ in $P$,
there exist $x$ and $y$ in $P$, such that $xa=yb$.  This is the content of Ore's Theorem. In view of this, some authors
prefer to say that $P$ is an \"{Ore semigroup}, rather than an \"{Ore sub-semigroup of $G$}.}, provided it satisfies
$G=P\inv P$.

\state Lemma \label PrWellDefined Let $G$ be a group and let $P$ be an Ore sub-semigroup of $G$.  Also let $A$ be a
unital C*-algebra and $ v :P\to A $ be a semigroup of isometries.  Then, for every $m, n, p, q \in P$, one has that $$
m\inv n = p\inv q \imply v_m^*v_n = v_p^*v_q.  $$

\Proof Given that $m\inv n = p\inv q$, we have that $pm\inv = qn\inv $.  By hypothesis we may choose $r$ and $s$ in $G$,
such that $$ pm\inv = qn\inv = r\inv s.  $$ Therefore $ rp = sm$ and $rq = sn, $ so $$ v_m = v_s^*v_sv_m = v_s^*v_{sm} =
v_s^*v_{rp} = v_s^*v_rv_p, $$ whence $$ v_m^*v_n = v_p^* v_r^* v_sv_n = v_p^* v_r^* v_{sn} = v_p^* v_r^* v_{rq} = v_p^*
v_r^* v_rv_q = v_p^* v_q.  \omitDoubleDollar \endProof

We thus see that the condition that $G=P\inv P$ gives a canonical way to extend a semigroup of isometries $v $ on $P$ to
a map defined on $G$.  Simply take any $g$ in $G$, decompose it as $g = m\inv n$, and define $\pr g = v_m^*v_n$.  By the
above result, $\pr g$ does not depend on the chosen decomposition of $g$.

As already mentioned, we are interested in extending $v $ to a partial representation of $G$, but it is by no means
clear that the map $\Pr $ thus described will be a partial representation.  In particular, if the range projections
$v_nv_n^*$ do not commute with each other, then \cite {MainTheEgCommute} implies that $\Pr $ cannot be a partial
representation.

\state Theorem \label OreExtension Let $G$ be a group and let $P$ be an Ore sub-semigroup of $G$.  Also let $A$ be a
unital C*-algebra and $ v :P\to A $ be a semigroup of isometries.  Then the following are equivalent: \izitem \zitem
$v_mv_m^*$ commutes with $v_nv_n^*$, for all $n,m$ in $P$, \zitem there exists a *-partial representation $\Pr $ of $G$
in $A$ extending $v $.  \medskip \noindent In this case, the partial representation $\Pr $ referred to in (ii) is
unique.

\Proof The implication (ii)$\Rightarrow $(i) is an easy consequence of \cite {MainTheEgCommute}, so let us focus on the
converse.  Given $g$ in $G$, write $g = m\inv n$, with $m,n\in P$, and set $\pr g = v_m^*v_n$.  By \cite {PrWellDefined}
we see that $\Pr $ is a well defined function on $G$, and we will next prove it to be a *-partial representation.

Upon replacing $P$ with $P\cup \{1\}$, we may assume that $1\in P$, and as seen in \cite {SGUnital}, we have that
$v_1=1$, so \cite {PrOnOne} is verified.  The proof of \cite {StarPrepCond} is also elementary, so let us prove \cite
{ObserveRight}.  Given $g,h\in G$, write $$ g = m\inv n \and h = p\inv q, $$ with $m, n, p, q \in P$.  Using the
hypothesis, pick $r$ and $s$ in $P$ such that $$ np\inv = r\inv s, $$ so $rn = sp$.  Letting $$ m' = rm,\quad
n'=rn,\quad p'=sp \and q'=sq, $$ notice that $$ g = {m'}\inv n' \and h = {p'}\inv q', $$ but now $n'=p'$.  This says
that, in our initial choice of $m$, $n$, $p$ and $q$, we could have taken $n=p$.  With this extra assumption (and
removing all diacritics) we then have that $gh = m\inv q$, and $$ \pr g \pr h \pri h = v_m^*\ v_p v_p^*\ v_qv_q^*\ v_p =
v_m^*\ v_qv_q^*\ v_p v_p^*\ v_p \$= v_m^*v_qv_q^*v_p = \pr {gh}\pri h.  $$ This proves \cite {ObserveRight} and, as
already noticed in \cite {LeaveOutOneAxiom}, axiom \cite {ObserveLeft} must also hold.

The uniqueness of $\Pr $ now follows immediately from \cite {RegardDisregard}, as observed in \cite {PrepDetermined}.
\endProof

Motivated by the result above, we will now study semigroups of isometries which are not necessarily based on Ore
sub-semigroups, but which may be extended to a partial representation.  So let us introduce the following terminology.

\definition \label DefineExtendable Let $G$ be a group and let $P$ be a sub-semigroup of $G$.  Also let $A$ be a unital
C*-algebra and $ v :P\to A $ be a semigroup of isometries.  We will say that $v $ is \subjex {extendable}{extendable
semigroup}, if there exists a *-partial representation $\Pr $ of $G$, such that $\pr n = v_n$, for every $n\in P$.

In the case of Ore semigroups, as we have seen in \cite {OreExtension}, a necessary and sufficient condition for
extendability is the commutativity of range projections.  Also, in the general case it easy to see that this condition
is still necessary but we have unfortunately not been able to determine an elegant set of sufficient conditions for
extendability.

In fact the concept of extendability elicits many questions that we will not even attempt to answer.  For example, it is
not clear whether the extended partial representation is unique, even under the assumption that $G$ is generated by $P$.

An example of a non-extendable semigroup of isometries is clearly presented by \cite {SGIFromToep}, provided $\varphi $
and $\psi $ are as in \cite {RangeToeplitzCommute}, with $a\neq 0$.

\medskip The following is a dilation result for extendable semigroups of isometries:

\state Theorem \label DilationExtendabSG Let $G$ be a group and let $P$ be a sub-semigroup of $G$.  Also let $$ v :P\to
\Lin (H) $$ be an extendable semigroup of isometries, where $H$ is a Hilbert space.  Then there exists a Hilbert space
$\tilde H$ containing $H$, and a unitary group representation $\tilde u$ of $G$ on $\tilde H$ such that, for every $n$
in $P$, one has that: \izitem \zitem $\tilde u_n(H)\subseteq H$, \zitem $v_n$ coincides with the restriction of $\tilde
u_n$ to $H$, \zitem denoting by $p$ the orthogonal projection from $\tilde H$ to $H$, one has that $p$ commutes with
$\tilde u_{g\inv }p\tilde u_g$, for every $g$ in $G$.

\Proof Since $v $ is extendable, we may pick a *-partial representation $\Pr $ of $G$ on $H$ extending $v $.  Let
$\tilde \Pr $, $\tilde H$ and $j$ be obtained by applying \cite {AbadieOnDilation} to $\Pr $.  Viewing $H$ as a subspace
of $\tilde H$ via $j$, we then have that (iii) is verified.  Given $n$ in $P$, we have by \cite {AbadieOnDilation} that
$$ j^*\tilde u_nj = \pr n = v_n, $$ and, taking into account that $v_n$ is an isometry, we conclude that $$ 1 = v_n^*v_n
= (j^*\tilde u_nj)^*j^*\tilde u_nj = j^*\tilde u_{n\inv }j j^*\tilde u_nj.  $$ Multiplying this on the left by $j$ and
on the right by $j^*$, we get $$ p = jj^* = j1j^* = jj^*\tilde u_{n\inv }j j^*\tilde u_njj^* = p\tilde u_{n\inv }p\tilde
u_np = \ldots $$ By the already proved point (iii), the above equals $$ \ldots = \tilde u_{n\inv }p\tilde u_npp = \tilde
u_{n\inv }p\tilde u_np.  $$ This proves that $p=\tilde u_{n\inv }p\tilde u_np$, and if we now multiply this on the left
by $\tilde u_n$ we arrive at $$ \tilde u_np=p\tilde u_np, $$ which says that the range of $p$, namely $H$, is invariant
under $\tilde u_n$, hence proving (i).  As for (ii), we have for all $n$ in $P$ that $$ v_n = \pr n = j^*\tilde u_nj =
\tilde u_n|_H, $$ where the last equality is justified by our identification of $H$ as a subspace of $\tilde H$, and the
invariance of the former under $\tilde u_n$.  \endProof

\definition Let $P$ be a semigroup (not necessarily contained in a group) and let $$ v :P\to \Lin (H) $$ be a semigroup
of isometries, where $H$ is a Hilbert space.  By a \subjex {unitary dilation}{unitary semigroup dilation} of $v $ we
mean a semigroup of isometries $$ w:P\to \Lin (\tilde H) $$ on a Hilbert space $\tilde H$ containing $H$, such that, for
every $n$ in $P$, $w_n$ is a \"{unitary operator}, leaving $H$ invariant, and whose restriction to $H$ coincides with
$v_n$.

As a consequence of \cite {DilationExtendabSG}, we see that every extendable semigroup of isometries admits a dilation.
In particular, every semigroup of isometries based on an Ore semigroup whose final projections commute, is extendable by
\cite {OreExtension}, and hence admits a dilation.  This should be compared with \ref {Theorem 1.4/Laca/2000}, where a
similar result has been proven without the hypothesis of the commutativity of final projections, but also without
conclusion \cite {DilationExtendabSG/iii}.

\nrem The study of isometries and semigroups thereof dates back at least to the 1960's and have involved numerous
authors, such as Arveson \ref {Arveson/1991}, Burdak \ref {Burdak/2004}, Coburn \ref {Coburn/1967}, Douglas \ref
{Douglas/1969}, \ref {Douglas/1972}, Hor\'ak and M\umlaut uller \ref {HorakMuller/1989}, Laca \ref {Laca/2000}, Nica
\ref {Nica/1992}, Phillips and Raeburn \ref {PhillipsRaeburn/1993}, Popovici \ref {Popovici/2004} and S\l oci\'nski \ref
{Slocinski/1980} among others.

\chapter Quasi-lattice ordered groups

\def \t {\tilde }

Besides the case of Ore semigroups, we will now present another important example of extendable semigroup of isometries.

Let us begin by noticing that if $P$ is a sub-semigroup of a group $G$, such that $P\cap P\inv = \{1\}$, then one may
define a left-invariant partial order relation on $G$ by $$ g\leq h \iff g\inv h\in P, $$ for all $g$ and $h$ in $G$.
Conversely, if ``$\leq $'' is a left-invariant partial order relation on $G$, then $$ P:= \{g\in G: g\geq 1\} $$ is a
sub-semigroup with $P\cap P\inv = \{1\}$.  So there is a one-to-one correspondence between such sub-semigroups and
left-invariant partial order relations on $G$.

\definition Let $G$ be a group with a distinguished sub-semigroup $P$, such that $P\cap P\inv = \{1\}$, and denote by
``$\leq $'' the corresponding order relation on $G$.  Given a nonempty subset $A\subseteq G$, \izitem \zitem we say that
an element $k\in G$ is an \subj {upper bound} for $A$ and, if $k\geq g$, for all $g$ in $A$, \zitem we say that an
element $m\in G$ is a \subj {least upper bound} for $A$, if $m$ is an upper bound for $A$, and $m\leq k$, whenever $k$
is an upper bound for $A$, \zitem we say that an element \underbar {$m\in P$} is a \subjex {least upper bound \underbar
{in $P$}}{least upper bound in $P$} for $A$, if $m$ is an upper bound for $A$, and $m\leq k$, whenever $k$ is an upper
bound for $A$, with \underbar {$k\in P$}, \zitem \zitemBismark QLO we say that $(G,P)$ is a \subj {quasi-lattice}, or
that $G$ is \"{quasi-lattice ordered}, if every nonempty finite subset $A\subseteq G$, admitting an upper bound in $P$,
necessarily admits a least upper bound in $P$, \goodbreak \zitem \zitemBismark WQLO we say that $(G,P)$ is a \subj {weak
quasi-lattice}, or that $G$ is \"{weakly quasi-lattice ordered}, if every nonempty finite subset $A\subseteq P$,
admitting an upper bound\fn {Notice that, if $A\subseteq P$, then any upper bound of $A$ is necessarily in $P$.},
necessarily admits a least upper bound.

\bigskip The first two definitions above are, of course, part of the elementary theory of ordered sets.  They are listed
just to provide a backdrop for the others.  Moreover, given the subtleties involved in the last two definitions, it is
perhaps worth repeating them in more technical terms.  Given a nonempty subset $A\subseteq G$, let us denote by $$
A^\uparrow = \{k\in G: k \hbox { is an upper bound for } A\} \equationmark DefUpBdSet $$ Then \iaitem \aitem $(G,P)$ is
a quasi-lattice if and only if, for every nonempty finite subset $A\subseteq G$, such that $A^\uparrow \cap P\neq
\emptyset $, there exists a smallest element in $A^\uparrow \cap P$.  \aitem $(G,P)$ is a weak quasi-lattice if and only
if, for every nonempty finite subset $A\subseteq P$, such that $A^\uparrow \neq \emptyset $, there exists a smallest
element in $A^\uparrow $.

\bigskip Notice that the notion of weak quasi-lattice may be defined using only the algebraic structure of $P$, while
checking that a pair $(G,P)$ is a quasi-lattice requires a knowledge of the whole group $G$.  In other words, the notion
of weak quasi-lattice easily generalizes for semigroups which are not necessarily contained in a group, while the notion
of quasi-lattice must necessarily refer to a pair $(G,P)$, as opposed to a single semigroup $P$.

It is clear that every subset $A\subseteq P$ which admits a least upper bound, also admits a least upper bound in
$P$. Therefore every quasi-lattice is necessarily a weak quasi-lattice.

\medskip If a set $A$ admits a least upper bound $m$, then $m$ is clearly unique and we denote it by $$ m = \mx A.  $$
In addition, if $A$ is a two element set, say $A=\{g,h\}$, we also denote by $$ g\vee h = \mx \{g,h\}.  $$

The notion of quasi-lattice ordered group was introduced by Nica in \ref {Nica/1992}, where it was noticed that every
weak quasi-lattice is a quasi-lattice, provided condition \cite {QLO} holds for sets $A$ with a single element.  Nica
denoted this special case of \cite {QLO} as (QL1), which he phrased as follows:

\bigskip \itemitem {(QL1)} Any $g$ in $PP\inv $ (these are precisely the elements having an upper bound in $P$) has a
least upper bound in $P$.

\bigskip Ten years after the publication of Nica's paper, Crisp and Laca \ref {Lemma 7/CrispLaca/2002} found what might
well be called ``Columbus' egg" of quasi-lattices, realizing that in fact (QL1) alone is a sufficient condition for
$(G,P)$ to be a quasi-lattice (see below).

\medskip In hindsight, the following simple characterizations should hopefully clean up our act.  We suggest that the
reader take this as alternative definitions for the concepts involved, should the above perhaps too extensive discussion
have blurred the big picture:

\state Proposition Let $G$ be a group and let $P\subseteq G$ be a sub-semigroup such that $P\cap P\inv = \{1\}$.  Then
\izitem \zitem $(G,P)$ is a weak quasi-lattice, if and only if $m\vee n$ exists for every $m$ and $n$ in $P$ admitting a
common upper bound, \zitem $(G,P)$ is a quasi-lattice, if and only if any $g$ in $PP\inv $ has a least upper bound in
$P$.

\Proof Point (i) is easily proven by induction.

As for (ii), assume that $(G,P)$ is a quasi-lattice and let $g\in PP\inv $.  Writing $g=mn\inv $, with $n$ and $m$ in
$P$, it is clear that $g\leq m$, so $\{g\}$ admits an upper bound in $P$, hence also a least upper bound in $P$.

Conversely, suppose that any $g$ in $PP\inv $ has a least upper bound in $P$.  We will first prove that $(G,P)$ is a
weak quasi-lattice by verifying condition (i).  So suppose that $n$ and $m$ are in $P$, and that $p$ is an upper bound
of $n$ and $m$.  Then $$\def \quad {\kern 4pt} \matrix { 1 &=& n\inv n &\leq & n\inv p & \and \cr \pilar {16pt} && n\inv
m &\leq & n\inv p,} $$ so $n\inv m$ admits an upper bound in $P$, namely $n\inv p$.  By hypothesis there exists a least
upper bound in $P$ for $n\inv m$, say $q$.  Then $q$ is a least upper bound for the set $\{1,n\inv m\}$, and by
left-invariance of the order relation, $nq$ is a least upper bound for $\{n,m\}$, as desired.

We will now prove that $(G,P)$ is a quasi-lattice.  So let $A$ be a nonempty finite subset of $G$ admitting an upper
bound in $P$, say $p$.  By hypothesis we know that, for each $g$ in $A$, there exists the least upper bound of $g$ in
$P$, which we denote by $n_g$.  Then it is clear that $n_g\leq p$, for all $g$ in $A$, so $p$ is also an upper bound for
$A' = \{n_g:g\in A\}$.  Since $A'\subseteq P$, the weak quasi-lattice property of $(G,P)$ may be used to ensure that the
least upper bound in $P$ of $A'$ exists, but this is clearly also the least upper bound in $P$ of $A$.  \endProof

One of the main examples of quasi-lattice ordered groups is obtained by taking $\F _n$ to be the free group on an
arbitrary (finite or infinite) number $n$ of generators, and taking ${\bf P}_n$ to be the smallest sub-semigroup of $\F
_n$ containing 1 and its free generators.

An example of a weak quasi-lattice which is not a quasi-lattice was recently obtained by Scarparo \ref {Scarparo/2014},
who also generalized some important results on quasi-lattice ordered groups to the context of weak quasi-lattices.
Scarparo's example is as follows: take $\F _2$ to be the free group on a set of two generators, say $\{a,b\}$, and let
$$ {\bf P}':= b{\bf P}_2\cup \{1\}, $$ where ${\bf P}_2$ was defined above. Then ${\bf P}'$ is clearly a sub-semigroup
of $\F _2$ and the pair $(\F _2,{\bf P}')$ is not a quasi-lattice although it is a weak quasi-lattice.

The reason why it is not a quasi-lattice is that the singleton $\{ba\inv b\inv \}$ is bounded in ${\bf P}'$ by $b$ and
$ba$, but it does not admit a least upper bound in ${\bf P}'$, as the reader may easily verify.

On the other hand $(\F _2,{\bf P}')$ is a weak quasi-lattice because ${\bf P}'$ is isomorphic to ${\bf P}_\infty $, as a
semigroup, via an isomorphism defined using the set $$ \{ba^n: n\in \N \}, $$ as a set of free generators for ${\bf
P}'$.  Since $(\F _\infty ,{\bf P}_\infty )$ is a quasi-lattice, as seen above, it is necessarily also a weak
quasi-lattice.

Another interesting consequence of Scarparo's example is that a semigroup $P$ may be seen as a sub-semigroup of two
non-isomorphic groups $G_1$ and $G_2$, in such a way that both $(G_1,P)$ and $(G_2,P)$ are weak quasi-lattices.  Namely,
take $$ (G_1,P) = (\F _2,{\bf P}') \and (G_2,P) = (\F _\infty ,{\bf P}_\infty ), $$ were we are identifying ${\bf P}'$
and ${\bf P}_\infty $, as already mentioned.

\medskip Whenever possible we will strive to work with weak quasi-lattices, not only because it provides for stronger
results, but also due to the fact that weak quasi-lattices may be defined intrinsically, independently of the group
containing them.

\nostate \label NicasIdea Nica's theory of semigroups of isometries based on quasi-lattices is rooted on the following
simple idea.  Let $(G,P)$ be a weak quasi-lattice and consider the regular semigroup of isometries $$ \lambda : P \to
\Lin (\ell ^2(P)) $$ defined in \cite {RegularSGP}.  Reinterpreting \cite {RangeLambdaP}, we may say that $\lambda
_p\lambda _p^*$ is the orthogonal projection onto $$ \barspan \{e_n: n\geq p\}.  \equationmark RangeLambdaPReinterp $$
Thus, given $p$ and $q$ in $P$, the product of the commuting projections $\lambda _p\lambda _p^*$ and $\lambda _q\lambda
_q^*$ coincides with the orthogonal projection onto the closed linear span of the set $$ \{e_n: n\geq p\} \cap \{e_n:
n\geq q\} = \{e_n: n\ \geq \ p\vee q\}, $$ assuming of course that $p\vee q$ exists.  When $p\vee q$ does not exist,
then by hypothesis $p$ and $q$ have no upper bound whatsoever, so the above intersection of sets is empty.  This
motivates the following:

\definition \label DefNCC Let $(G,P)$ be a weak quasi-lattice and $A$ be a C*-algebra.  We shall say that a given
semigroup of isometries $v :P\to A$ satisfies \subj {Nica's covariance condition}, \subj {NCC} for short, if for every
$m$ and $n$ in $P$, one has that $$ v_mv_m^* \ v_nv_n^* = \left \{\matrix { v_{m\vee n}v_{m\vee n}^*\hfill ,& \hbox {if
$m\vee n$ exists}, \cr \pilar {12pt} \hfill 0\hfill ,& \hbox {otherwise}.\hfill } \right .  $$

As briefly indicated above, in the case of a weak quasi-lattice, the regular semigroup of isometries satisfies NCC.

\state Lemma \label TechNCCTame Let $(G,P)$ be a weak quasi-lattice, $A$ be a C*-algebra, and $v :P\to A$ be a semigroup
of isometries satisfying NCC.  Given $m,n\in P$, one has that: \izitem \zitem if $m$ and $n$ have no common upper bound,
then $v_m^*v_n=0$, \zitem if $m$ and $n$ have a common upper bound, then there are $x$ and $y$ in $P$ such that $mx =
ny$, and $v_m^*v_n=v_xv_y^*$.

\Proof Under condition (i) we have $$ v_m^*v_n = v_m^*\underbrace {v_mv_m^*v_nv_n^*}_{(\DefNCC )}v_n = 0.  $$

On the other hand, if $m$ and $n$ have a common upper bound, then the least upper bound exists, and there are $x$ and
$y$ in $P$ such that $$ m\vee n = mx = ny.  $$ So $$ v_m^*v_n = v_m^*\underbrace {v_mv_m^*v_nv_n^*}_{(\DefNCC )}v_n =
v_m^*v_{m\vee n}v_{m\vee n}^*v_n = v_m^*v_{mx}v_{ny}^*v_n = v_xv_y^*.  \omitDoubleDollar \endProof

The main goal of this chapter is to prove extendability of a given semigroups of isometries $v$ satisfying NCC.  After
this is achieved, by \cite {PrepIsTame} we will have that the range of $v$ is a tame set.  However the strategy we will
adopt will require knowing in advance that the range of $v$ is tame, so we need to prove this as an intermediate step.

\state Proposition \label RangeNCCTame Let $(G,P)$ be a weak quasi-lattice, $A$ be a C*-algebra, and $v :P\to A$ be a
semigroup of isometries satisfying NCC.  Then, \Zitem the set $\{v_mv_n^*: m,n\in P \}\cup \{0\}$ is a multiplicative
semigroup, \zitem the range of $v$ is a tame set of partial isometries.

\Proof Given $m$, $n$, $p$ and $q$ in $P$, let us suppose that $n$ and $p$ have no common upper bound.  Then $$
v_mv_n^*v_pv_q^* \={TechNCCTame/i} 0.  $$ On the other hand, if $n$ and $p$ have a common upper bound, then by \cite
{TechNCCTame/ii} we may write $ v_n^*v_p = v_xv_y^*, $ for suitable $x$ and $y$, so $$ v_mv_n^*v_pv_q^* = v_m v_xv_y^*
v_q^* = v_{mx}v_{qy}^*.  $$ This proves (i).

In order to prove (ii), let $S$ be the semigroup referred to in (i).  Since $P$ contains 1 by hypothesis, and since
$v_1=1$ by \cite {SGUnital}, we see that both the range of $v$ and its adjoint are contained in $S$.  Therefore we see
that $S$ coincides with the multiplicative semigroup generated by the range of $v$ and its adjoint.  To prove that the
range of $v$ is tame, it therefore suffices to verify that every element in $S$ is a partial isometry.  Given $m$ and
$n$ in $P$, we have that $$ (v_mv_n^*) (v_mv_n^*)^* = v_mv_n^*v_nv_m^* = v_mv_m^*, $$ which is clearly a projection, and
hence $v_mv_n^*$ is a partial isometry by \cite {CharacPisoViaProj}.  This concludes the proof of (ii).  \endProof

We will now present a crucial technical result designed to help proving the extendability of semigroups of isometries
satisfying NCC.  It should be compared to \cite {PrWellDefined}.

\state Lemma \label CompaIndepRep Let $(G,P)$ be a weak quasi-lattice, $A$ be a unital C*-algebra and let $$ v :P\to A,
$$ be a semigroup of isometries satisfying NCC.  Given $m$, $n$, $p$ and $q$ in $P$ such that $mn\inv = pq\inv $, one
has that $v_mv_n^*$ and $v_pv_q^*$ are compatible partial isometries, as defined in \cite {DefCompat}.

\Proof We will soon see that the question of upper bounds for the set $\{n,q\}$ is a crucial one, so let us begin by
analyzing it.  Suppose first that $n$ and $q$ admit an upper bound $k$.  So there are $x$ and $y$ in $P$ such that $$ k=
nx = qy.  $$ Then $$ (mx) (nx)\inv = mn\inv = pq\inv = (py) (qy)\inv , $$ so $mx = py$.  Defining $$ h= mx = py, $$ it
is then clear that $h$ is an upper bound for $\{m,p\}$.  Denoting by $ g= mn\inv = pq\inv , $ notice that $$ gk = mn\inv
nx = mx = h, $$ so this provides for a well defined function $$ k\in \{n,q\}^\uparrow \ \mapsto \ gk\in \{m,p\}^\uparrow
, $$ where we are using the notation introduced in \cite {DefUpBdSet} for the set of upper bounds.  Given the symmetric
roles played by $m$, $n$, $p$ and $q$, it is clear that $$ h\in \{m,p\}^\uparrow \ \mapsto \ g\inv h\in \{n,q\}^\uparrow
$$ is also well defined.  Therefore $\{n,q\}$ admits an upper bound if and only if $\{m,p\}$ does, and in this case the
left-invariance of our order relation implies that $$ g(n\vee q) = m\vee p.  \equationmark CorrespondMax $$

Returning to our semigroup of isometries, and letting $$ s=v_mv_n^*, $$ we have $$ ss^* = v_mv_n^*v_nv_m^* = v_mv_m^*,
$$ so we see that $s$ is indeed a partial isometry, sharing final projection with $v_m$.  It is also easy to see that
the initial projection of $s$ coincides with the final projection of $v_n$.  Likewise, if we let $$ t=v_pv_q^*, $$ then
the initial and final projections of $t$ are $v_qv_q^*$ and $v_pv_p^*$, respectively.

Recalling that our task is to prove that $$ st^*t=ts^*s \and tt^*s = ss^*t, \equationmark CompatGoal $$ let us assume
first that $\{n,q\}$ admits no upper bound.  In this case, NCC implies that $v_nv_n^*$ is orthogonal to $v_qv_q^*$,
whence $v_n^*v_q=0$.  So, $$ st^*t = v_mv_n^* v_qv_q^* =0.  $$ On the other hand, $$ ts^*s = v_pv_q^* v_nv_n^* =
v_p(v_n^* v_q)^*v_n^* = 0, $$ proving the first equation in \cite {CompatGoal}.  In order to prove the second one,
observe that, by the reasoning at the beginning of this proof, $\{m,p\}$ admits no upper bound either, so $v_p^*v_m=0$,
whence $$ tt^*s = v_pv_p^*v_mv_n^* = 0 = v_mv_m^*v_pv_q^* = ss^*t, $$ completing the proof of \cite {CompatGoal} under
the assumption that $\{n,q\}$ admits no upper bound.  Let us therefore assume that $\{n,q\}$ admits an upper bound, and
hence by hypothesis also a least upper bound.  We may then find $x$ and $y$ in $P$ such that $$ n\vee q=nx=qy.  $$
Therefore $$ st^*t = v_mv_n^* v_qv_q^* = v_mv_n^* v_n v_n^* v_qv_q^* = v_mv_n^* v_{n\vee q} v_{n\vee q}^* = v_mv_n^*
v_{nx} v_{nx}^* \$= v_mv_n^* v_nv_x v_{nx}^* = v_m v_x v_{nx}^* = v_{mx} v_{nx}^*.  $$ By \cite {CorrespondMax} we have
that $m\vee p$ exists and $$ \def \quad {\kern 4pt}\pilar {16pt} \stake {14pt} \matrix { m\vee p = g(n\vee q) &=& mn\inv
nx &=& mx & =\cr \pilar {12pt} &=& pq\inv qy &=& py.}  $$ Therefore $$ ts^*s = v_pv_q^*v_nv_n^* =
v_pv_q^*v_qv_q^*v_nv_n^* = v_pv_q^*v_{q\vee n}v_{q\vee n}^* = v_pv_q^*v_{qy}v_{nx}^* \$= v_pv_q^*v_qv_yv_{nx}^* =
v_{py}v_{nx}^* = v_{mx}v_{nx}^*, $$ thus proving that $ st^*t = ts^*s, $ which is the first equation in \cite
{CompatGoal}.  The second one may now be verified similarly, or by applying what has already been proved to $$ s':= s^*
=v_nv_m^* \and t':= t^* =v_qv_p^*.  $$ This concludes the proof.  \endProof

The following is the main result of this chapter.  It was first proved in \ref {QuiggRaeburn/1997} for the case of
quasi-lattices.  Our plan to prove it for weak quasi-lattices will require that the range of our semigroup of isometries
lies in a von Neumann algebra.

\state Theorem \label QLExtendable Let $(G,P)$ be a weak quasi-lattice, and let $A$ be a unital C*-algebra.  Suppose
moreover that either \izitem \zitem $A$ is a von Neumann algebra, or \zitem $(G,P)$ is a quasi-lattice.  \medskip
\noindent Then every semigroup of isometries $ v : P\to A $ which satisfies NCC is extendable.

\Proof Given $g$ in $G\setmenos PP\inv $, we define $$ \pr g = 0.  $$ On the other hand, given $g$ in $PP\inv $,
consider the collection of partial isometries $$ T_g = \{v_mv_n^*: m,n\in P,\ mn\inv = g\}.  $$

By \cite {CompaIndepRep} we have that the elements of $T_g$ are mutually compatible and we would now like to argue that
$\mx T_g$ exists.  Under the hypothesis that $A$ is a von Neumann algebra we may simply use \cite {ExistMaxCompat}, so
let us prove the existence of $\mx T_g$ under (ii).  In this case we will actually prove that $T_g$ contains a maximum
element.

Observing that $g$ may be written as $g=pq\inv $, for some $p,q\in P$, notice that we then have that $g\preceq p$, so we
see that $g$ admits an upper bound in $P$.  Using the quasi-lattice property, the least upper bound in $P$ of $\{g\}$,
here denoted by $\sigma (g)$, exists.  In particular $ g\preceq \sigma (g), $ so $$ \tau (g):= g\inv \sigma (g)\in P.
$$ We may then write $$ g = \sigma (g)\tau (g)\inv , \equationmark MostEfficientDec $$ which might be thought of as the
\"{most efficient way} of writing $g$.  Evidently $ v_{\sigma (g)} v_{\tau (g)}^* $ lies in $T_g$, and we claim that
this element dominates every other element of $T_g$.  In fact, given any $m,n\in P$ such that $g=mn\inv $, we have that
$g\preceq m$.  Since $\sigma (g)$ is the least upper bound in $P$ of $g$, it follows that $\sigma (g)\preceq m$ as well,
so there exists $x$ in $P$ such that $m =\sigma (g)x$.  Consequently $$ n = g\inv m = \tau (g)\sigma (g)\inv \sigma
(g)x= \tau (g)x.  $$

To prove our claim that $v_{\sigma (g)} v_{\tau (g)}^*$ is the biggest element in $T_g$, we must prove that
$v_mv_n^*\preceq v_{\sigma (g)} v_{\tau (g)}^*$, so we compute $$ v_{\sigma (g)} v_{\tau (g)}^*(v_mv_n^*)^*v_mv_n^* =
v_{\sigma (g)} v_{\tau (g)}^*v_nv_m^*v_mv_n^* = v_{\sigma (g)} v_{\tau (g)}^*v_nv_n^* \$= v_{\sigma (g)} v_{\tau
(g)}^*v_{\tau (g)}v_xv_n^* = v_{\sigma (g)x}v_n^* = v_mv_n^*.  $$ This proves that $v_{\sigma (g)} v_{\tau (g)}^*$ is
the biggest element in $T_g$, so it is obvious that $\mx T_g$ exists.

Having proven the existence of $\mx T_g$ under either (i) or (ii), we may define $$ \pr g = \mx T_g, $$ and we will now
set out to prove that $\Pr $ is a partial representation, the strategy being to apply \cite {characPRViaIneq}.  We first
notice that, thanks to \cite {RangeNCCTame}, the range of $v$ is a tame set, so $$ S := \big \langle v (P)\cup v
(P)^*\big \rangle , $$ meaning the multiplicative sub-semigroup of $A$ generated by $v (P)\cup v (P)^*$, is an inverse
semigroup by \cite {TameIsISG}.

Notice that by construction, under hypothesis (ii), the values of $\Pr $ lie in $S$ itself, while under (i), the values
of $\Pr $ lie in the $\vee $-closure of $S$, which is an inverse semigroup by \cite {ClosureISG}.

In either case we may view $\Pr $ as a map from $G$ to some inverse semigroup, which is a necessary condition for the
application of \cite {characPRViaIneq}.  In fact we need $\Pr $ to take values in a \"{unital} inverse semigroup, so we
actually have to throw in the unit of $A$.

We will now check conditions \cite {characPRViaIneq/i--iii}, but since the verification of \cite {characPRViaIneq/i\&ii}
is trivial, our task is simply to show that $$ \pr g \pr h \preceq \pr {gh}, \equationmark PrgPrhLessPrgh $$ for all $g$
and $h$ in $G$.  If either $g$ or $h$ lie outside $PP\inv $, the left-hand-side of \cite {PrgPrhLessPrgh} vanishes, so
there is nothing to do.  We therefore assume that both $g$ and $h$ lie in $PP\inv $, and in this case we claim that $\pr
{gh}$ dominates $T_gT_h$.  To see this, choose elements in $T_g$ and $T_h$, respectively of the form $v_mv_n^*$ and
$v_pv_q^*$, where $g=mn\inv $ and $h=pq\inv $.  We must then prove that $$ v_mv_n^* v_pv_q^* \preceq \pr {gh}.
\equationmark VmnpqLessPrgh $$

If $\{n,p\}$ admits no upper bound then \cite {TechNCCTame/i} implies that $v_n^*v_q=0$, so the left-hand-side of \cite
{VmnpqLessPrgh} vanishes and, again, there is nothing to do.  On the other hand, if $\{n,p\}$ admits an upper bound, we
have by \cite {TechNCCTame/ii} that there are $x$ and $y$ in $P$ such that $nx=py$, and $ v_n^*v_p = v_{x}v_{y}^*, $
whence $$ v_mv_n^* v_pv_q^* = v_m v_{x}v_{y}^*v_q^* = v_{mx}v_{qy}^*.  \equationmark VmVnVmxVqy $$ Observe that $$
mx(qy)\inv = mxy\inv q\inv = m n\inv p q\inv = gh, $$ so $ v_{mx}v_{qy}^*\in T_{gh}, $ and we then deduce from \cite
{VmVnVmxVqy} that $$ v_mv_n^* v_pv_q^* \preceq \mx T_{gh} = \pr {gh}.  $$

This proves \cite {VmnpqLessPrgh}, and hence that $\pr {gh}$ indeed dominates $T_gT_h$.  Consequently $$ \pr g\pr h =
(\mx T_g) (\mx T_h) \={PrdMx} \mx (T_gT_h) \preceq \pr {gh}, $$ taking care of \cite {PrgPrhLessPrgh}, whence allowing
us to apply \cite {characPRViaIneq}, proving that $\Pr $ is a partial representation, as desired.

Let us now prove that $\Pr $ extends $v$.  For this observe that if $p$ is in $P$, then $$ \pr p = \mx T_p = \mx
\{v_mv_n^*: m,n\in P,\ mn\inv =p \}.  $$ Given $m$ and $n$ in $P$ such that $mn\inv =p$, we have that $m = pn$, and we
claim that $$ v_mv_n^*\preceq v_p.  $$ In fact, $$ v_p(v_mv_n^*)^*(v_mv_n^*) = v_pv_nv_m^*v_mv_n^* = v_{pn}v_n^* =
v_mv_n^*.  $$

This proves the claim, so we see that $v_p$ indeed dominates $T_p$.  Since $v_p$ moreover belongs to $T_p$, we deduce
that $$ v_p = \mx T_p = \pr p, $$ proving that $v$ indeed extends $\Pr $, which is to say that $v$ is extendable.  This
concludes the proof.  \endProof

In the above proof, we did not worry about the uniqueness of the extended partial representation $\Pr $.  However this
will later become relevant, so let us discuss it now.  In doing so we will be content with the quasi-lattice case,
leaving the study of uniqueness in the case of weak quasi-lattices as an outstanding problem.

Again referring to the strategy adopted in the above proof, notice that our first step was to set $\pr g=0$, whenever
$g$ is not in $PP\inv $.  Evidently this is equivalent to saying that $$ \e g = \pr g\pri g =0.  $$ On the other hand,
in the quasi-lattice case, we saw that every $g$ in $PP\inv $ has a \"{most efficient} decomposition, namely \cite
{MostEfficientDec}, in which case recall that we have defined $$ \pr g = \mx T_g = v_{\sigma (g)} v_{\tau (g)}^*, $$
whence $$ \e g = \pr g\pri g = v_{\sigma (g)} v_{\tau (g)}^* v_{\tau (g)} v_{\sigma (g)}^* = v_{\sigma (g)} v_{\sigma
(g)}^* \$= \pr {\sigma (g)} \pri {\sigma (g)} = \e {\sigma (g)} = \e {g\vee 1}.  $$

\state Proposition \label UniqueExtOfSGIso Let $(G,P)$ be a quasi-lattice, and let $A$ be a unital C*-algebra.  Then
every semigroup of isometries $ v : P\to A $ satisfying NCC admits a \underbar {unique} extension to a partial
representation $ \Pr $ of $G$ in $A$, satisfying $$ \e g = \left \{\matrix { \e {g\vee 1} ,& \hbox {if $g\in PP\inv $},
\cr \pilar {12pt} \hfill 0\hfill ,& \hbox {otherwise},\hfill } \right .  \equationmark ConditionExtMCCIso $$ for all $g$
in $G$, where $\e g = \pr g\pri g$, as usual.

\Proof The existence follows from \cite {QLExtendable} and our discussion above.  Regarding uniqueness, let $\Pr $ be an
extension of $v$ satisfying the above properties.  Then, for all $g$ in $G\setmenos PP\inv $, one has that $$ \pr g =
\pr g \pri g \pr g = \e g \pr g = 0 \pr g = 0.  $$

On the other hand, if $g$ is in $PP\inv $, then $g$ admits an upper bound in $P$, so $$ m:= g\vee 1 $$ exists by
hypothesis, so we may write $g = mn\inv $, where $n = g\inv m\in P$.  By hypothesis we then have that $ \e g = \e m, $
so, $$ \pr g = \pr g \pri g \pr g = \e g \pr g = \e m \pr {mn\inv } = \pr m \pri m \pr {mn\inv } = \pr m \pri n = v_m
v_n^*.  $$ This proves that $\Pr $ is unique.  \endProof

Although apparently innocuous, condition \cite {ConditionExtMCCIso} has very interesting consequences for our purposes:

\state Proposition \label RestrictToIsomSg Let $(G,P)$ be a quasi-lattice, and let $A$ be a unital C*-algebra.  Given
any partial representation $\Pr : G\to A$ satisfying \cite {ConditionExtMCCIso}, the restriction of $\Pr $ to $P$ is a
semigroup of isometries satisfying NCC.

\Proof Given any $n$ in $P$, notice that $n\inv \leq 1$, so $n\inv \vee 1=1$.  Consequently $$ 1 = \e 1 = \e {n\inv \vee
1} = \e {n\inv } = \pr n^* \pr n, $$ from where we conclude that $\pr n$ is an isometry.  In particular each $\pr n$ is
left-invertible, so by \cite {RegardDisregard/iii} we have that $$ \pr m\pr n = \pr {mn} \for m,n\in P.  $$ This shows
that $\Pr |_P$ is a semigroup of isometries.

To conclude, let us now prove NCC.  For this, suppose first that $m$ and $n$ are given elements of $P$ not possessing a
common upper bound.  Then we claim that $m\inv n\not \in PP\inv $.  To see this, assume by contradiction that $$ m\inv n
= pq\inv , $$ where $p,q\in P$.  Then $$ m,n\leq nq = mp, $$ contradicting our assumption.  By hypothesis we then have
$\e {m\inv n}=0$, so also $\pr {m\inv n}=0$.  Focusing on \cite {DefNCC} we then compute $$ \pr m \pr m^* \ \pr n \pr
n^* = \pr m \pri m \pr n \pri n \={ObserveRight} \pr m \pr {m\inv n} \pri n = 0.  $$

Suppose now that $m$ and $n$ are elements of $P$ possessing a common upper bound.  Then $m\vee n$ exists and there are
$x$ and $y$ in $P$ such that $$ m\vee n = mx = ny.  $$ Moreover, given the left invariance of the order relation on $G$,
we have $$ 1\vee (m\inv n) = m\inv (m\vee n) = x, $$ whence, by hypothesis we have $ \e {m\inv n} = \e x, $ and then $$
\pr m^* \ \pr n \={RegardDisregard/iii} \pr {m\inv n} = \e {m\inv n}\pr {m\inv n} = \e x\pr {m\inv n} \$= \pr x \pri x
\pr {m\inv n} = \pr x \pr {x\inv m\inv n} = \pr x \pri y.  $$ Therefore, $$ \pr m \pr m^* \ \pr n \pr n^* = \pr m \pr x
\pri y \pri n = \pr {mx} \pri {(ny)} = \pr {m\vee n} \pr {m\vee n}^*, $$ as desired.  \endProof

\syschapter {C*-algebras generated by semi\-groups of isometries}{C*-algebras generated by semigroups}

\def \cstarab {C^*_{ab}} \def \Isom {{\rsbox P}\,} \def \OIsom {\Omega _\Isom }

Given a semigroup $P$, we may consider the universal C*-algebra $C^*(P)$ generated by a set $$ \{w_n:n\in P\}, $$
subject to the requirement that the correspondence $$ n\mapsto w_n $$ be a semigroup of isometries.  The universal
property of $C^*(P)$ may then be phrased as follows: for every unital C*-algebra $A$, and every semigroup of isometries
$v:P\to A$, there exists a unique *-homomorphism $\varphi :C^*(P)\to A$, such that the diagram

\null \hfill \beginpicture \setcoordinatesystem units <0.0020truecm, -0.0025truecm> point at 0 0 \put {$P$} at 000 000
\put {$A$} at 1000 000 \arrow <0.11cm> [0.3,1.2] from 200 0 to 800 000 \put {$v$} at 500 -100 \put {$C^*(P)$} at 000 700
\arrow <0.11cm> [0.3,1.2] from 0 200 to 0 500 \put {$w$} at -150 350 \arrow <0.11cm> [0.3,1.2] from 300 550 to 800 200
\put {$\varphi $} at 650 480 \endpicture \hfill \null

\vskip 0.5cm \noindent commutes.  So, to study the representation theory of $C^*(P)$ is equivalent to studying all
semigroups of isometries based on $P$.

In the case of ${\bf N}{\times }{\bf N}$, observe that given any semigroup of isometries $v:{\bf N}{\times }{\bf N} \to
A$, one has that $v_{(1,0)}$ and $v_{(0,1)}$ are commuting isometries.  Conversely, given any pair of commuting
isometries, we obtain by \cite {NCrosN} a semigroup of isometries based on ${\bf N}{\times }{\bf N}$.  Thus one may
argue that, to fully understand $C^*({\bf N}{\times }{\bf N})$ is to understand all commuting pairs of isometries, and
in fact many authors have dedicated a significant amount of energy in this endeavor.

In order to give the reader a feeling for the difficulty of the problem at hand, we point out that many results in the
literature on this subject have the following format: given any pair of commuting isometries $S$ and $T$ on a Hilbert
space $H$, there is an orthogonal decomposition $$ H=H_1\oplus \cdots \oplus H_n $$ in invariant subspaces, in such a
way that the restrictions of $S$ and $T$ to the $H_i$ may be described in more or less concrete terms.  For example,
\izitem \zitem $S$ and $T$ are unitary on $H_1$, \zitem $S|_{H_2}$ is unitary and $T|_{H_2}$ is a multiple of the
unilateral shift, \zitem $S|_{H_3}$ is a multiple of the unilateral shift and $T|_{H_3}$ is unitary, etc.

\medskip However, as one approaches the last space $H_n$, the available descriptions often become more and more
technical and sometimes include an \"{evanescent} space, about which little can be said.  We refer the reader to \ref
{Slocinski/1980} and \ref {HorakMuller/1989} for some of the earliest results, and to \ref {Popovici/2004} and \ref
{Burdak/2004} for what we believe represents the state of the art in this hard subject.

The difficulties presented by the theory of semigroups of isometries based on the otherwise nicely behaved semigroup
${\bf N}{\times }{\bf N}$ is perhaps a sign that insurmountable obstacles lie ahead of the corresponding theory for more
general semigroups.  We will therefore restrict our study to a nicer class of semigroups of isometries, namely the
extendable ones.

Given a group $G$ and a sub-semigroup $P\subseteq G$, recall from \cite {DefineExtendable} that a semigroup of
isometries $$ v: P\to A $$ is extendable if there exists a partial representation $\Pr $ of $G$ in $A$ such that $\pr n
= v_n$, for every $n$ in $P$.  Thus, to study extendable semigroups of isometries based on $P$ is the same as studying
partial representations of $G$ which ascribe isometries to the elements of $P$.

\definition \label DefineCstarGP Let $G$ be a group and let $P\subseteq G$ be a sub-semigroup.  We will denote\fn {We
should warn the reader that our choice of notation here is not standard.  In particular it conflicts with the notation
adopted in \ref {Nica/1992} for a C*-algebra which we will later study under the notation $\NicaAlg $.}  by $C^*(G,P)$
the universal unital C*-algebra generated by a set $ \{\pr g: g\in G\} $ subject to relations \cite {DefPR/i--iv}, in
addition to $$ \pr n^*\pr n = 1\for n\in P.  $$

Observe that this is a special case of \cite {DefineCstarParRel}, in the sense that $C^*(G,P)$ coincides with
$\CstarparRel G\Isom $, where $\Isom $ consists of the above set of relations stating that $\pr n$ is an isometry for
every $n$ in $P$.

As a consequence of \cite {CstarRelAsCP}, we therefore have that $C^*(G,P)$ may be described as the partial crossed
product algebra $C(\OIsom )\rt {\Th \ur } G$, and we will now give a precise description of the space $\OIsom $.

Our relations, if written in the form \cite {PolyRel}, become $$ \ei n -1 =0 \for n\in P, $$ so the set $\FRel $
mentioned in \cite {DefineOmegaR} is formed by the functions $$ \omega \in \OG \ \mapsto \ \bool {n\inv \in \omega } -1
\in \C .  $$ Still according to \cite {DefineOmegaR}, we have that $\OIsom $ consists of all $\omega \in \OuG $ such
that $$ \bool {n\inv \in g\inv \omega } -1 =0 \for n\in P \for g\in \omega .  $$

The above condition for an element $\omega \in \OuG $ to be in $\OIsom $ may be interpreted as $$ g\in \omega \imply
gn\inv \in \omega \for n\in P.  \equationmark SpecIsom $$

Since $h\leq g$ if and only if $h = gn\inv $, for some $n$ in $P$, we see that \cite {SpecIsom} precisely expresses that
$\omega $ is hereditary\fn {Recall that a subset $S$ of an ordered set $X$ is said to be hereditary if, whenever $x$ and
$y$ are elements of $X$, such that $x\leq y\in S$, one has that $x\in S$.}.  The above analysis and \cite {CstarRelAsCP}
therefore give us the following:

\state Theorem \label CstarGPasCP Let $G$ be a group and let $P\subseteq G$ be a sub-semigroup.  Consider the closed
subspace $\OIsom $ of\/ $\OuG $ formed by the hereditary elements, equipped with the partial action of $G$ obtained by
restricting the partial Bernoulli action.  Then there is a natural *-isomorphism between $C^*(G,P)$ and the partial
crossed product $C(\OIsom )\rt {\Th \ur } G$.  Under this isomorphism, each $\pr g$ corresponds to $1_g\delta _g$, where
$1_g$ denotes the characteristic function of $$ \D g^\Isom = \{\omega \in \OIsom : g\in \omega \}.  $$

One may argue that the above result is not fully satisfactory since the definition of $C^*(G,P)$ puts too much emphasis
on the group $G$, while one might actually only be interested in the semigroup $P$.  Further research is therefore
needed in order to answer a few outstanding questions such as:

\state Question \rm \label QuestionOne Can extendability of a semigroup of isometries $v$ be characterized via algebraic
relations involving only the isometries $v_n$, for $n$ in $P$?

Recall that \cite {OreExtension} provides an answer to this question in the case of Ore semigroups.

We have already mentioned that the commutativity of range projections is a necessary condition for extendability of a
semigroup of isometries, so any answer to \cite {QuestionOne} is likely to include the commutativity of range
projections.  This also motivates the following:

\definition \label DefineCstarAbP Given a semigroup $P$, we denote by $\cstarab (P)$ the universal unital C*-algebra
generated by a set $\{w_n:n\in P\}$, subject to the relations \izitem \zitem $w_n^*w_n = 1$, \zitem $w_mw_n = w_{mn}$,
\zitem $w_mw_m^*$ commutes with $w_nw_n^*$, \medskip \noindent for all $n,m$ in $P$.

Another relevant open question is as follows:

\state Question \label QuestionTwo \rm Given a sub-semigroup $P$ of a group $G$, can the subalgebra of $C^*(G,P)$
generated by $\{\pr p:p\in P\}$ be concretely described?

In the above context, observe that, given a group element $g\in P\inv P$, and writing $g=m\inv n$, with $m,n\in P$, one
has by \cite {RegardDisregard} that $$ \pr g = \pr {m\inv n} = \pri m \pr n, $$ since $\pr n$ is an isometry, and hence
left-invertible.  The image of $P\inv P$ under $\Pr $ is therefore contained in the subalgebra referred to in \cite
{QuestionTwo}.  If $G$ coincides with $P\inv P$, that is, if $P$ is an Ore sub-semigroup of $G$, the answer to question
\cite {QuestionTwo} is therefore that the subalgebra mentioned there coincides with the whole of $C^*(G,P)$.

Being able to satisfactorily answer the above questions for Ore sub-semigroups, we may give a nice description of
$\cstarab (P)$, as follows:

\state Theorem Let $P$ be an Ore sub-semigroup of a group $G$.  Then $\cstarab (P)$ is naturally isomorphic to
$C^*(G,P)$ and consequently there is a *-isomorphism $$ \varphi : \cstarab (P) \to C(\OIsom )\rt {\Th \ur } G, $$ such
that $$ \varphi (w_n) = 1_n\delta _n \for n\in P, $$ where $\OIsom $ and $1_n$ are as in \cite {CstarGPasCP}.

\Proof We will prove that there is an isomorphism $$ \psi : \cstarab (P) \to C^*(G,P), \equationmark WhatPhiIs $$ such
that $$ \psi (w_n) = \pr n \for p\in P.  \equationmark WhatPhiDoes $$ The conclusion will then follow immediately from
\cite {CstarGPasCP}.  As a first step, observe that the correspondence $$ n\in P \mapsto \pr n\in C^*(G,P) $$ is a
semigroup of isometries, therefore evidently satisfying \cite {DefineCstarAbP/i\&ii}, and which also satisfies \cite
{DefineCstarAbP/iii} thanks to \cite {MainTheEgCommute}.  The universal property of $\cstarab (P)$ may therefore be
invoked to prove the existence of a *-homomor\-phism $\psi $, as in \cite {WhatPhiIs}, satisfying \cite {WhatPhiDoes},
and we need only prove that $\psi $ is an isomorphism.

Noticing that the correspondence $$ n\in P\mapsto w_n\in \cstarab (P) $$ is a semigroup of isometries satisfying \cite
{OreExtension/i}, we have by \cite {OreExtension} that there exists a partial representation $$ \tilde w: G\to \cstarab
(P), $$ such that $\tilde w_n = w_n$, for every $n$ in $P$.  Since the $\tilde w_n$ are isometries, the universal
property of $C^*(G,P)$ implies that there exists a *-homomorphism $$ \gamma :C^*(G,P)\to \cstarab (P), $$ such that $
\gamma (\pr n) = \tilde w_n, $ for all $n$ in $P$.  We then conclude that, for all $n\in P$, one has that $$ \gamma \big
(\psi (w_n)\big ) \={WhatPhiDoes} \gamma (\pr n) = \tilde w_n =w_n, $$ from where we see that $\gamma \circ \psi $ is
the identity on $\cstarab (P)$.  In particular this shows that $\psi $ is one-to-one.  In order to prove that $\psi $ is
onto $C^*(G,P)$, it is enough to show that the standard generating set $\{\pr g: g\in G\}$ of $C^*(G,P)$ lies in the
range of $\psi $.  For this, given $g$ in $G$, write $g=m\inv n$, with $m$ and $n$ in $P$.  Then $$ \pr g = \pr {m\inv
n} \={RegardDisregard} \pri m \pr n = (\pr m)^* \pr n = \psi (w_ m)^* \psi (w_n) \$= \psi (w_ m^*w_n) \in \psi \big
(\cstarab (P)\big ).  $$ This proves that $\psi $ is an isomorphism, as desired.  \endProof

Since the range of a partial representation is always a tame set by \cite {PrepIsTame}, the above result implies that
the isometries canonically generating $\cstarab (P)$ form a tame set in case $P$ is an Ore sub-semigroup.  However the
following seems to be open:

\state Question \rm Given a semigroup $P$, is the subset $\{w_n:n\in P\}$ of $\cstarab (P)$ tame?  If not, what is the
smallest set of relations one can add to the definition of $\cstarab (P)$ to make the answer affirmative?

\chapter {Wiener-Hopf C*-algebras}

\def \CParNica {\CstarparRel G{\NicaRel }}

In this chapter we continue the study initiated above of C*-algebras associated to semigroups of isometries.  One of the
most important among these is the Wiener-Hopf C*-algebra associated to a quasi-lattice ordered group $(G,P)$, so we
dedicate the present chapter, in its entirety, to the study of this example.

We will initially concentrate on the study of a C*-algebra introduced by Nica in \ref {Nica/1992}, which we shall denote
by $\NicaAlg $, and which should be seen as the \"{full} version of the actual Wiener-Hopf algebra to be defined later.

\fix We will now fix, for the entire duration of this chapter, a group $G$, and a sub-semigroup $P\subseteq G$, such
that $(G,P)$ is a quasi-lattice.

\definition \label DefineNCCAlg We will denote by $\NicaAlg $ the universal unital C*-algebra for semigroups of
isometries based on $P$ satisfying NCC. \ More precisely, $\NicaAlg $ is the universal unital C*-algebra generated by a
set $\{v_n: n\in P\}$, subject to the relations below for all $n$ and $m$ in $P$: \izitem \zitem $v_n^*v_n = 1$, \zitem
$v_mv_n = v_{mn}$, \zitem $ v_mv_m^* \ v_nv_n^* = \left \{\matrix { v_{m\vee n}v_{m\vee n}^*\hfill ,& \hbox { if $m\vee
n$ exists}, \cr \pilar {12pt} \hfill 0\hfill ,& \hbox { otherwise}.\hfill } \right .  $

Some authors denote the algebra introduced above by $C^*(G,P)$, but we have chosen $\NicaAlg $ in order to avoid
conflict with the notation introduced in \cite {DefineCstarGP}.  Our notation is also intended to emphasize that
$\NicaAlg $ depends only on the algebraic structure of $P$, rather than on $G$.

As we have seen in \cite {UniqueExtOfSGIso} and \cite {RestrictToIsomSg}, there exists a one-to-one correspondence
between semigroup of isometries satisfying NCC and partial representations of $G$ satisfying \cite {ConditionExtMCCIso}.

\state Proposition \label NicaIsPrepCond Let $\NicaRel $ be the set of relations \cite {ConditionExtMCCIso}, for all $g$
in $G$, and consider the universal C*-algebra $\CParNica $ for partial representations of $G$ satisfying $\NicaRel $, as
defined in \cite {DefineCstarParRel}.  Then there is a *-isomorphism $$ \varphi : \NicaAlg \to \CParNica , $$ such that
$\varphi (v_n) = \pr n$, for all $n$ in $P$, where we denote the canonical partial representation of $G$ in $\CParNica $
by $\Pr $.

\Proof Since $\Pr $ satisfies \cite {ConditionExtMCCIso} by construction, \cite {RestrictToIsomSg} implies that the
restriction of $\Pr $ to $P$ is a semigroup of isometries satisfying NCC.  The universal property of $\NicaAlg $ then
provides for a unital *-homomorphism $$ \varphi : \NicaAlg \to \CParNica , $$ such that $$ \varphi (v_n) = \pr n \for
n\in P.  \equationmark PhiVnPrn $$

On the other hand, by \cite {UniqueExtOfSGIso}, the universal semigroup of isometries $v:P\to \NicaAlg $ extends to a
partial representation $\tilde v:G\to \NicaAlg $ satisfying \cite {ConditionExtMCCIso}. So, again by universality, there
is a unital *-homomorphism $$ \psi : \CParNica \to \NicaAlg , $$ such that $$ \psi (\pr g) = \tilde v_g \for g\in G.  $$

Given any $n\in P$, we have that $$ \psi \big (\varphi (v_n)\big ) = \psi \big (\pr n\big ) = \tilde v_n = v_n, $$ from
where we see that $\psi \circ \varphi $ is the identity mapping on $\NicaAlg $.

Notice that the correspondence $$ n\in P \mapsto \varphi (v_n) \in \CParNica $$ is a semigroup of isometries satisfying
NCC, since this is obtained by composing $v$ with a *-homomorphism.  We could now apply \cite {UniqueExtOfSGIso} to
prove that it extends to a partial representation of $G$, except that we already know two extensions of it, the first
one being $$ g\in G \mapsto \varphi (\tilde v_g)\in \CParNica , $$ and the second one being simply $\Pr $, as one may
easily deduce from \cite {PhiVnPrn}.  Clearly both of these extensions satisfy \cite {ConditionExtMCCIso}, so the
uniqueness part of \cite {UniqueExtOfSGIso} implies that $$ \varphi (\tilde v_g) = \pr g \for g\in G.  $$ Consequently
we have $$ \varphi \big (\psi (\pr g)\big ) = \varphi (\tilde v_g) = \pr g, $$ so $\varphi \circ \psi $ coincides with
the identity mapping of $\CParNica $ on a generating set, whence $\varphi \circ \psi $ is the identity.  This proves
that $\varphi $ is indeed an isomorphism, as desired.  \endProof

Using \cite {CstarRelAsCP} we may then describe $\CParNica $, and hence also $\NicaAlg $, as a partial crossed product.
But, in order to state a more concrete result, let us first give a description of the spectrum of $\NicaRel $ which is
more geometrical than the one given in \cite {DefineOmegaR}.

\state Proposition \label SpecNicaHeredDir The spectrum of the set of relations $\NicaRel $ is given by $$ \global
\advance \fnctr by 1 \ON = \big \{\omega \in \OuG : \omega \hbox { is hereditary and directed\kern 1pt}^{\number \fnctr
}\big \}.  \footnote {}{\kern -9.5pt $^{\number \fnctr }$\kern 5pt \eightpoint Recall that a subset $S$ of an ordered
set $X$ is said to be directed if, for all $x,y\in S$, there exists $z$ in $S$, with $z\geq x,y$.\par \vskip -10pt} $$

\Proof We will first prove ``$\subseteq $'', so we pick any $\omega \in \ON $.  For every $n$ in $P$, we have that
$n\inv \vee 1=1$, so $\NicaRel $ includes the relation ``$\ei n = \e 1$''.  Consequently, for every $g\in \omega $, we
have by the definition of $\ON $, that $$ \bool {n\inv \in g\inv \omega } = \bool {1\in g\inv \omega } = \bool {g\in
\omega } = 1, $$ which is to say that $ gn\inv \in \omega , $ thus proving that $\omega $ is hereditary.

We next claim that $$ \omega \subseteq PP\inv .  \equationmark OmegaInPPinv $$ To see this let $g\in G\setmenos PP\inv
$.  Then $\NicaRel $ includes the relation ``$\e g = 0$'' and, again by the definition of $\ON $, we have that $[g\in
\omega ] = 0$, so $g\notin \omega $, proving the claim.

We will next prove that $\omega $ is directed.  For this let $g,h\in \omega $, so by \cite {OmegaInPPinv} we may write
$g=mn\inv $ and $h = pq\inv $.  Assuming, without loss of generality, that $m=g\vee 1$ and $p = h\vee 1$, we have that
$$ 1 = \bool {g\in \omega } = \bool {g\vee 1\in \omega } = \bool {m\in \omega }, $$ so $m\in \omega $, and similarly
$p\in \omega $.  Observing that $$ g=mn\inv \leq m \and h = pq\inv \leq p, $$ we see that it is enough to find a common
upper bound for $m$ and $p$ in $\omega $, since this will automatically be a common upper bound for $g$ and $h$.

Recalling that $\ON $ is invariant under the partial Bernoulli action, and since $m\in \omega $, we have that $m\inv
\omega \in \ON $.  So the conclusion reached in \cite {OmegaInPPinv} also applies to $m\inv \omega $, whence $$ m\inv p
\in m\inv \omega \subseteq PP\inv , $$ and therefore $(m\inv p)\vee 1$ exists.  The relation ``$\e {m\inv p} = \e
{(m\inv p)\vee 1}$'' is then among the relations in $\NicaRel $, so for all $k$ in $\omega $ one has $$ \bool {m\inv p
\in k\inv \omega } = \bool {(m\inv p){\vee }1 \in k\inv \omega }.  $$ Plugging in $k=m$, left-hand-side above evaluates
to 1, so the right-hand-side does too, meaning that $$ (m\inv p)\vee 1 \in m\inv \omega , $$ whence $$ \omega \ni m\big
((m\inv p)\vee 1\big ) = p\vee m, $$ where the last step above, including the existence of $p\vee m$, is a consequence
of the left-invariance of the order relation on $G$.  This proves that $m$ and $p$ have a common upper bound in $\omega
$, thus proving that $\omega $ is directed.

Let us now prove the inclusion ``$\supseteq $'' between the sets mentioned in the statement, so we pick an hereditary
and directed $\omega \in \OuG $.  We will first prove that $\omega $ is contained in $PP\inv $.  For this pick $g\in
\omega $ and observe that, since 1 is also in $\omega $, there exists an upper bound $m$ for $\{g,1\}$ in $\omega $.
Thus $m$ is in $P$ and, since $g\leq m$, we deduce that $n:= g\inv m$ also lies in $P$.  So $$ g = mn\inv \in PP\inv .
$$ This verifies our claim that $\omega \subseteq PP\inv $.

In order to prove that $\omega \in \ON $, we must show that $f(h\inv \omega )=0$, for all $h$ in $\omega $, and for all
$f$ in $\FNRel $ (see \cite {DefineOmegaR} for the definition of $\FNRel $).

According to \cite {ConditionExtMCCIso} we need to consider two cases.  Firstly, when $g$ is not in $PP\inv $, the
corresponding relation ``$\e g=0$'' in $\NicaRel $ leads to the function $$ f(x) = \bool {g\in x} \for x\in \OuG , $$ so
we must prove that $$ g\not \in h\inv \omega \for h\in \omega .  \equationmark SoWeMustProveThat $$

By the left-invariance of the order relation on $G$, for every $h$ in $\omega $, we have that $h\inv \omega $ is also
hereditary and directed.  Moreover, if $h\in \omega $, then $h\inv \omega $ also lies in $\OuG $, so \cite
{OmegaInPPinv} holds for $h\inv \omega $, just like it does for $\omega $.  So $$ h\inv \omega \subseteq PP\inv , $$
from where \cite {SoWeMustProveThat} follows immediately.

The second case to be considered is when $g\in PP\inv $, in which case $g\vee 1$ exists and we are led to the function
$f$ in $\FNRel $ given by $$ f(\xi ) = \bool {g\in \xi } - \bool {g{\vee }1\in \xi } \for \xi \in \OuG .  $$ Given $h$
in $\omega $ we must then show that $f(h\inv \omega )=0$, which translates into $$ \bool {g\in h\inv \omega } = \bool
{g{\vee }1\in h\inv \omega }, $$ or, equivalently $$ \def \quad {\ } \matrix { hg\in \omega \iff & h(g \vee 1) & \in
\omega , \cr &\pilar {10pt} \vertequal \cr & (hg)\vee h\stake {10pt}} $$ where the vertical equal sign is a consequence
of left-invariance.  Using that $\omega $ is hereditary and directed, and that $h$ is in $\omega $, this may now be
easily verified.  \endProof

Interpreting \cite {CstarRelAsCP} in the present case, and using \cite {NicaIsPrepCond}, we have the following concrete
description of $\NicaAlg $ as a partial crossed product.

\state Theorem \label NicaAsCP Given a quasi-lattice $(G,P)$, consider the subset $\ON $ of $2^G$ formed by the
hereditary directed subsets of $G$ containing 1.  Then $\ON $ is compact and invariant under the partial Bernoulli
action.  In addition there is a *-isomorphism $$ \psi :\NicaAlg \to C(\ON )\rt \null G, $$ such that $\psi (v_n) =
1_n\delta _n$, where $1_n$ denotes the characteristic function of the set $\{\omega \in \ON : n\in \omega \}$.

Observe that if $\omega $ is a hereditary and directed subset of $G$, then $\omega \cap P$ is a hereditary and directed
subset of $P$.  We will next explore this correspondence in order to obtain another description of $\ON $.

\state Proposition \label HerDirInP The map $$ \omega \in \ON \mapsto \omega \cap P \in 2^P $$ is a homeomorphism from
$\ON $ onto the subset of\/ $2^P$ formed by all nonempty, hereditary, directed subsets of $P$.

\Proof We leave it for the reader to verify the continuity of our map.

Given any $\omega $ in $\ON $, let $\xi =\omega \cap P$.  Using that $\omega $ is hereditary and directed by \cite
{SpecNicaHeredDir}, it is easy to see that $$ \omega = \xi P\inv = \big \{g\in G: g\leq n, \hbox { for some } n\in \xi
\big \}, $$ from where it follows that the map in the statement is injective.

It is obvious that $\omega \cap P$ is a hereditary directed subset of $P$ for every $\omega $ in $\ON $, and conversely,
for every nonempty, hereditary, directed subset $\xi \subseteq P$, one has that $$ \omega _\xi := \xi P\inv $$ lies in
$\ON $, and $\xi = \omega _\xi \cap P$. This proves that the range of our map is as described in the statement.  Since
$\ON $ is compact and $2^P$ is Hausdorff, it follows that our map is a homeomorphism onto its range.  \endProof

The above result can easily be used to transfer the partial Bernoulli action on $\ON $ to a partial action of $G$ on the
set of all nonempty, hereditary, directed subsets of $P$, thus providing yet another model for the partial dynamical
system leading up to $\NicaAlg $ in \cite {NicaAsCP}.

\medskip Given any $m$ in $P$, notice that the set $$ mP\inv = \{mn\inv : n\in P\} = \{g\in G: g\leq m\} \equationmark
DefmPInv $$ is a hereditary and directed subset of $G$ containing 1, hence $mP\inv $ is an element of $\ON $.

Given any $\omega \in \ON $, we may then consider the net $$ \{mP\inv \}_{m\in \omega \cap P}, \equationmark SelfNet $$
where $\omega \cap P$, carrying the order induced from $G$, is a directed set by \cite {HerDirInP}, hence suitable for
playing the role of the index set for a net.

\state Proposition \label ConvergenceInON For every $\omega $ in $\ON $, the net described in \cite {SelfNet} converges
to $\omega $.

\Proof Given $\omega $ in $\ON $ and a neighborhood $V$ of $\omega $, by definition of the product topology there are $$
g_1,g_2,\ldots ,g_r;\ h_1,h_2,\ldots ,h_s\in G, $$ such that the \"{basic neighborhood} $$ U := \{\eta \in \ON : g_i\in
\eta ,\ h_j\not \in \eta ,\ \forall i\leq r,\ \forall j\leq s\} $$ satisfies $$ \omega \in U\subseteq V.  $$

In particular the $g_i$ lie in $\omega $, so we may use the fact that $\omega $ is directed to find some $m_0$ in
$\omega $ such that $m_0\geq g_i$, for all $i\leq r$.  Upon assuming without loss of generality that 1 is among the
$g_i$, we have that $m_0\in \omega \cap P$.

We then claim that, for all $m$ in $\omega \cap P$ with $m\geq m_0$, one has that $mP\inv \in U$, which is to say that
\iaitem \aitem $g_i\in mP\inv $, for all $i\leq r$, and \aitem $h_j\not \in mP\inv $, for all $j\leq s$.

\medskip Point (a) is an obvious consequence of the fact that each $$ g_i\leq m_0\leq m.  $$ With respect to (b), assume
by way of contradiction, that some $h_j\in mP\inv $.  Then $h_j\leq m$, whence $h_j\in \omega $, because $m\in \omega $,
and $\omega $ is hereditary.  This is a contradiction with the fact that $\omega $ belongs to $U$, hence verifying (b).
This concludes the proof.  \endProof

Recall from \cite {NicasIdea} that the regular representation of $P$ on $\ell ^2(P)$ satisfies NCC.  Thus, by
universality, there is a *-representation of $\NicaAlg $ on $\ell ^2(P)$ consistent with the regular representation.
Our next goal will be to study the range of such a representation and, in particular, to show that it is isomorphic to
the \"{reduced} crossed product relative to the same partial dynamical system whose associated \"{full} crossed product
was shown to coincide with $\NicaAlg $ in \cite {NicaAsCP}.

\definition \label DefineWH The closed *-algebra of bounded operators on $\ell ^2(P)$ generated by the range of $\lambda
$ is called the \subjex {Wiener-Hopf}{Wiener-Hopf algebra} algebra of $P$, and it is denoted by $\WH $.

In \ref {Nica/1992}, Nica denotes the C*-algebra defined above by ${\cal W}(G,P)$, but since the above definition is
given only in terms of the semigroup $P$, we prefer not to explicitly mention the ambient group $G$.

\medskip We will next describe $\WH $ as a reduced partial crossed product.

\state Theorem \label WienerReduced There is a *-isomorphism $$ \chi :\WH \to C(\ON )\redrt \null G, $$ such that $\chi
(\lambda _n) = 1_n\delta _n$, where $1_n$ is as in \cite {NicaAsCP}.

\Proof The universal property of $\NicaAlg $ implies that there exists a *-homo\-mor\-phism $$ \sigma :\NicaAlg \to \WH
, $$ such that $\sigma (v_n) = \lambda _n$, for all $n$ in $P$, which is clearly onto.

Considering the *-isomorphism $\psi $ given in \cite {NicaAsCP} we may define a surjective *-homomorphism $ \rho =
\sigma \circ \psi \inv , $ as indicated in the diagram: $$ \matrix { C(\ON )\rt \null G & \trepa {\textstyle \rho
}\longrightarrow & \WH \cr \pilar {16pt} \stake {8pt} \hfill _{\textstyle \psi } \nwarrow \kern -5pt&& \kern -10pt
\nearrow _{_{\textstyle \sigma }} \hfill \cr &\NicaAlg } $$

We will now prove that the kernel of $\rho $ coincides with the kernel of the regular representation $$ \Lambda :C(\ON
)\rt \null G \to C(\ON )\redrt \null G, $$ from where the statement will follow.

Consider the conditional expectation $F$ from $\Lin \big (\ell ^2(P)\big )$ onto the subalgebra formed by the diagonal
operators relative to the standard orthonormal basis of $\ell ^2(P)$.  To be precise, $$ F(T) = \soma {n\in P}
e_{n,n}Te_{n,n} \for T\in \Lin \big (\ell ^2(P)\big ), $$ where the sum is interpreted in the strong operator topology,
and each $e_{n,n}$ is the orthogonal projection onto the one-dimensional subspace generated by the corresponding basis
vector.  It is a well known fact that $F$ is a faithful conditional expectation onto the algebra of diagonal operators.

We will also consider the conditional expectation $E$, given by the composition

\vskip -10pt \null \hfill \beginpicture \setcoordinatesystem units <0.0020truecm, 0.0025truecm> point at 0 0
\setplotarea x from 0000 to 1000, y from -350 to 450 \put {$ C(\ON )\rt \null G \ \trepa {\textstyle \Lambda
}\longrightarrow \ C(\ON )\redrt \null G \ \trepa {\textstyle E_1}\longrightarrow \ C(\ON ), $} at 000 000 \setquadratic
\plot -1150 150 75 400 1300 150 / \arrow <0.11cm> [0.3,1.2] from 1290 154 to 1300 150 \put {$E$} at 75 550 \endpicture
\hfill \null

\noindent where $E_1$ is the faithful conditional expectation given by \cite {Fourier}.  We then claim that the diagram
$$ \matrix { C(\ON )\rt \null G & \trepa {\textstyle \rho }\longrightarrow & \WH \cr \pilar {16pt} \stake {8pt} E\big
\downarrow && \big \downarrow F\cr C(\ON ) &\trepa {\textstyle \rho }\longrightarrow & \Lin \big (\ell ^2(P)\big ) } $$
commutes, where the bottom arrow is to be interpreted as the restriction of $\rho $ to the canonical copy of $C(\ON )$
within the crossed product.

By \cite {RangeNCCTame/i} we have that $\NicaAlg $ is the closed linear span of the set of elements of the form
$v_mv_n^*$, as $m$ and $n$ range in $P$.  Using \cite {NicaAsCP} we therefore have that $C(\ON )\rt \null G$ is spanned
by the elements $y$ of the form $$ y = \psi (v_mv_n^*) = (1_m\delta _m) (1_n\delta _n)^* = \beta _m(1_{m\inv }1_{n\inv
})\delta _{mn\inv } \$= 1_m1_{mn\inv }\delta _{mn\inv }, \equationmark EltYofForm $$ and the commutativity of the above
diagram will follow once we check that $F\rho (y) = \rho E(y)$, for every $y$ of the above form.

Assuming that $m\neq n$, one has that $mn\inv \neq 1$, and $E(y)=0$.  On the other hand, $$ \rho (y) = \sigma \circ \psi
\inv \!\!\circ \psi (v_mv_n^*) = \sigma (v_mv_n^*) = \lambda _m\lambda _n^*.  $$

Because $m\neq n$, one has that the matrix of $\lambda _m\lambda _n^*$ has no nonzero diagonal entry, which is to say
that $$ F\big (\rho (y)\big )= 0 = \rho \big (E(y)\big ), $$ as desired.

Assuming now that $m=n$, we have that $y$ lies in $C(\ON )$ (which we identify with $C(\ON )\delta _1$ as usual), so
$E(y)=y$.  On the other hand $$ \rho (y) = \lambda _m\lambda _m^*, $$ which is a diagonal operator by \cite
{RangeLambdaP}, whence $$ F\big (\rho (y)\big ) = \rho (y) = \rho \big (E(y)\big ), $$ proving our diagram to be
commutative, as claimed.

We will eventually like to show that $\rho $ is one-to-one on $C(\ON )$.  With this goal in mind we next claim that, for
every $k$ in $P$, and every $f$ in $C(\ON )$, one has that the $k^{th}$ diagonal entry of $\rho (f)$, here denoted by
$\rho (f)_{k,k}$, is given by $$ \rho (f)_{k,k} = f(kP\inv ), \equationmark DenseCharsDiag $$ where we are seeing
$kP\inv $ as an element of $\ON $, as we already did in \cite {SelfNet}.

Let us first assume that $f=1_m$, where $m$ is in $P$.  In this case, identifying $1_m$ with $1_m\delta _1$, as usual,
we have that $$ 1_m\delta _1 \={EltYofForm} \psi (v_mv_m^*), $$ so $\rho (1_m) = \lambda _m\lambda _m^*$, whence for
every $k$ in $P$, the corresponding diagonal entry of $\rho (1_m)$ is given by $$ \rho (1_m)_{k,k} = (\lambda _m\lambda
_m^*)_{k,k} \={RangeLambdaPReinterp} \bool {k\geq m} = \bool {m\in kP\inv } = 1_m(kP\inv ), $$ proving \cite
{DenseCharsDiag} for the particular case of $f=1_m$.  To prove the general case of this identity, it is enough to verify
that the $1_m$ generate $C(\ON )$ as a C*-algebra, which we do as follows: by Stone-Weierstrass we have that $$ \{1_g:
g\in G\} \equationmark GenerateWeiers $$ generates $C(\ON )$ as a C*-algebra.  On the other hand, since each $\omega $
in $\ON $ is hereditary and directed, and since $1\in \omega $, we have that $$ g\in \omega \iff g\vee 1\in \omega \for
g\in G.  $$ Consequently $1_g=1_{g\vee 1}$, when $g\vee 1$ exists, and $1_g=0$, otherwise.  Except for the identically
zero function we then have that \cite {GenerateWeiers} coincides with $$ \{1_n: n\in P\}, $$ which is therefore also a
generating set for $C(\ON )$, concluding the proof of \cite {DenseCharsDiag}.

Still aiming at the proof that the kernel of $\rho $ coincides with the kernel of the regular representation of $C(\ON
)\rt \null G$, we will next prove that $\rho $ is injective on $C(\ON )$.

For this, suppose that $\rho (f)=0$, for some $f$ in $C(\ON )$.  From \cite {DenseCharsDiag} we conclude that $f$
vanishes on every $kP\inv $ in $\ON $.  However, the set formed by these is dense in $\ON $ by \cite {ConvergenceInON},
so $f=0$.

With this we may now describe the null space of $\rho $, as follows: given any $y$ in $C(\ON )\rt \null G$, and using
that $F$ and $E_1$ are faithful and $\rho $ is injective on $C(\ON )$, we have that $$ \rho (y) = 0 \iff F\big (\rho
(y)^*\rho (y)\big ) = 0 \iff \rho \big (E(y^*y)\big ) = 0 \$\iff E(y^*y) = 0 \iff E_1\big (\Lambda (y^*y)\big ) = 0 \iff
\Lambda (y) =0.  $$

This shows that $\rho $ and $\Lambda $ share kernels, as claimed, so $\rho $ factors through $\Lambda $, producing a
*-isomorphism $$ C(\ON )\redrt \null G \ \trepa {\textstyle \tilde \rho }\longrightarrow \ \WH , $$ whose inverse
satisfies all of the required conditions.  \endProof

We may now use the results on amenability of Fell bundles to obtain conditions under which $\NicaAlg $ coincides with
the Wiener-Hopf algebra:

\state Theorem Let $(G,P)$ be a quasi-lattice, where $G$ is an amenable group.  Then $\NicaAlg $ is naturally isomorphic
to $\WH $.

\Proof Follows immediately from \cite {NicaAsCP}, \cite {WienerReduced}, and \cite {BundleAmenaGroup}.  \endProof

We dedicate the remainder of this chapter to prove a useful characterization of faithful representations of the
Wiener-Hopf algebra.

\state Proposition \label GPTopFree Let $(G,P)$ be a quasi-lattice.  Then the partial Ber\-noulli action restricted to
$\ON $ is topologically free.

\Proof Let $g\in G\setmenos \{1\}$ and assume by contradiction that there is an open subset $V$ of the domain of $\beta
_g$, where $\beta $ refers to the partial Bernoulli action, formed by fixed points for $\beta _g$.

By \cite {ConvergenceInON} the collection of all $mP\inv $, as $m$ range in $P$, forms a dense subset of $\ON $, so
there is some $m$ in $P$ such that $mP\inv \in V$.  Therefore $$ mP\inv = \beta _g(mP\inv ) = gmP\inv .  $$

Since $m$ is the maximum element in $mP\inv $, and $gm$ is the maximum element in $gmP\inv $, we deduce that $gm=m$, and
hence that $g=1$, a contradiction.  This concludes the proof.  \endProof

In view of \cite {WienerReduced} and the result above, Corollary \cite {CoroTopFree} applies to give a characterization
of faithful representations of the Wiener-Hopf algebra as those whose restriction to $C(\ON )$ are faithful.  However, a
more careful analysis will lead us to an even more precise result.  But before that we need the following auxiliary
result.

\state Lemma \label InvarContainsThis Let $(G,P)$ be a quasi-lattice and let $W$ be a subset of $\ON $ which is
nonempty, open and invariant.  Then there are $p_1,p_2,\ldots ,p_k\in P\setmenos \{1\}$, such that $W$ contains the
subset $$ \{\omega \in \ON : p_i\notin \omega , \hbox { for all } i=1,\ldots ,k\}.  $$

\Proof As before, we denote the partial Bernoulli action of $G$ on $\ON $ by $\beta $, and for each $g$ in $G$, we
denote the range of $\beta _g$ by $\D g^\NicaRel $.

Regarding the elements $mP\inv $ of \ $\ON $ described in \cite {DefmPInv}, we claim that $1P\inv $ lies in $W$.  To see
this observe that, for every $g$ in $G$, one has $$ 1P\inv \in \D g^\NicaRel \ \_\Leftrightarrow {DescriptBernouliDoms}\
g\in 1P\inv \iff g\inv \in P.  $$ The orbit of $1P\inv $ under $\beta $ is therefore given by $$ \{\beta _n(1P\inv ) :
n\in P\} = \{nP\inv : n\in P\}, $$ which is a dense subset of $\ON $ by \cite {ConvergenceInON}.  Since $W$ is nonempty,
it follows that $W$ contains $nP\inv $ for some $n\in P$. So, using the invariance of $W$, we deduce that $$ 1P\inv =
\beta _{n\inv }(nP\inv ) \in \beta _{n\inv }\big (W\cap \D n^\NicaRel \big )\subseteq W, $$ proving the claim.  We may
therefore pick $$ g_1,g_2,\ldots ,g_r;\ h_1,h_2,\ldots ,h_s\in G, $$ such that the basic neighborhood $$ U := \{\omega
\in \ON : g_i\in \omega ,\ h_j\not \in \omega ,\ \forall i\leq r,\ \forall j\leq s\} $$ satisfies $$ 1P\inv \in
U\subseteq W.  $$

Consequently, for each $i\leq r$, one has that $g_i\in 1P\inv $, so $g_i$ has the form $g_i=n_i\inv $, for some $n_i$ in
$P$.  However notice that $ n_i\inv \in \omega , $ for every single $\omega $ in $\ON $, because $$ n_i\inv \leq 1\_\in
{DescriptBernouliDoms} \omega , $$ and $\omega $ is hereditary by \cite {SpecNicaHeredDir}.  So the condition ``$g_i\in
\omega $'', appearing in the definition of $U$ above, is innocuous and hence may be omitted, meaning that $$ U =
\{\omega \in \ON : h_j\not \in \omega ,\ \forall j\leq s\}.  $$

Since every $\omega $ in $\ON $ is directed by \cite {SpecNicaHeredDir}, one may easily see that $\omega \subseteq
PP\inv $ (this was in fact explicitly proved in \cite {OmegaInPPinv}).  So, if any given $h_j$ is not in $PP\inv $, the
condition ``$h_j\not \in \omega $'' is true for any $\omega $ in $\ON $, and hence may also be eliminated from the
definition of $U$.  We may therefore assume, without loss of generality, that the $h_j$ all lie in $PP\inv $.

Recall that such elements are precisely the ones which admit an upper bound in $P$, and hence $h_j\vee 1$ exist, for
every $j$.  Again because every $\omega $ in $\ON $ is hereditary and directed, it is easy to see that $$ h_j\in \omega
\iff h_j\vee 1\in \omega \for \omega \in \ON .  $$ We then conclude that the condition ``$h_j\not \in \omega $'' in the
definition of $U$ may be replaced by ``$h_j\vee 1\notin \omega $'', allowing us to assume without loss of generality
that $h_j\in P$.  Moreover we must have $h_j\neq 1$, for every $j$, since otherwise $U=\emptyset $.  This concludes the
proof.  \endProof

The following is the characterization of faithful representations of the Wiener-Hopf algebra announced earlier.

\state Theorem \label FaithfulWiener Let $(G,P)$ be a quasi-lattice and let $\pi $ be a representation of\/ $\WH $ on a
Hilbert space.  Then $\pi $ is faithful if and only if, given any $p_1,p_2,\ldots ,p_k\in P\setmenos \{1\}$, one has
that $$ (1-V_1V_1^*) (1-V_2V_2^*) \cdots (1-V_kV_k^*) \neq 0, $$ where each $V_i = \pi (\lambda _{p_i})$.

\Proof Recall from \cite {DefineWH} that $\WH $ is the closed *-algebra of operators on $\ell ^2(P)$ generated by the
range of the regular semigroup of isometries of $P$, namely $$ \{\lambda _p:p\in P\}.  $$ It is also interesting to
point out that, apart from $\lambda _1$, which is the identity operator, all of the other $\lambda _p$'s are \"{proper}
isometries, in the sense that $$ 1-\lambda _p\lambda _p^*\neq 0 \for p\in P\setmenos \{1\}.  $$

Denoting by $e_1$ the first element of the canonical basis of $\ell ^2(P)$, observe that, for each $p$ in $P\setmenos
\{1\}$, one has that $e_1$ is orthogonal to the range of $\lambda _p\lambda _p^*$ by \cite {RangeLambdaPReinterp}, so
that in fact we have $$ (1-\lambda _p\lambda _p^*)(e_1)=e_1.  \equationmark EOneNotKilled $$

Therefore, given $p_1,p_2,\ldots ,p_k\in P\setmenos \{1\}$, as above, we have that $$ (1-\lambda _{p_1}\lambda _{p_1}^*)
(1-\lambda _{p_2}\lambda _{p_2}^*) \cdots (1-\lambda _{p_k}\lambda _{p_k}^*) \neq 0, $$ because this operator sends
$e_1$ to itself, as one may easily verify by successive applications of \cite {EOneNotKilled}.

If $\pi $ is a faithful representation of $\WH $, it therefore does not vanish on the above operator, hence proving the
``only if'' part of the statement.

In order to prove the converse, let us assume by contradiction that $\pi $ satisfies the condition in the statement and
yet it is not faithful.  Using the isomorphism $$ \chi : \WH \to C(\ON )\redrt \null G $$ provided by \cite
{WienerReduced}, we have that $$ \rho := \pi \circ \chi \inv $$ is a representation of the above reduced crossed
product, which is likewise non-faithful.

Denoting by $J$ the null space of $\rho $, we then have that $J$ is a nontrivial ideal of $C(\ON )\redrt \null G$.  So,
putting together \cite {GPTopFree} and \cite {topfree}, we deduce that $$ K:= J\cap C(\ON ) \neq \{0\}, $$ and moreover
$K$ is invariant by \cite {IntersInvar}.  Writing $K=C_0(W)$, where $W$ is an open subset of $\ON $, we then have that
$W$ is nonempty and invariant.

We may then invoke \cite {InvarContainsThis} to produce $p_1,p_2,\ldots ,p_k$ in $P\setmenos \{1\}$ such that $$ U:=
\{\omega \in \ON : p_i\notin \omega , \hbox { for all } i=1,\ldots ,k\} \subseteq W.  $$ It follows that any function in
$C(\ON )$, whose support is contained in $U$, necessarily lies in the null space of $\rho $.

Denoting, as usual, the characteristic function of the set $$ \{\omega \in \ON : p\in \omega \} $$ by $1_p$, observe
that the characteristic function of $U$, here denoted by $1_U$, is given by $$ 1_U= (1-1_{p_1}) (1-1_{p_2}) \ldots
(1-1_{p_k}).  $$ Therefore $$ 0 = \rho (1_U) = \medprod _{i=1}^k\rho (1-1_{p_i}) \={LotsRigInnProd} \medprod
_{i=1}^k\rho \big (1-(1_{p_i}\delta _{p_i})(1_{p_i}\delta _{p_i})^*\big ) \={WienerReduced} $$$$ = \medprod _{i=1}^k\rho
\big (1-\chi (\lambda _{p_i})\chi (\lambda _{p_i})^*\big ) = \medprod _{i=1}^k\big (1-\pi (\lambda _{p_i})\pi (\lambda
_{p_i})^*\big ) = \medprod _{i=1}^k\big (1-V_iV_i^*\big ), $$ where the $V_i$ are as in the statement.  This is a
contradiction, and hence the proof is concluded.  \endProof

\nrem Wiener-Hopf operators were first introduced by Norbert Wiener and Eberhard Hopf in \ref {WienerHopf/1931}, who
studied them from the point of view of integral equations.  In the context of quasi-lattice ordered groups, the
Wiener-Hopf C*-algebra $\WH $ and the algebra $\NicaAlg $ (with a different notation) was first studied by Nica in \ref
{Nica/1992}, where he observed that both $\WH $ and $\NicaAlg $ have a ``crossed product type structure'', emphasizing
the role of the abelian subalgebra ${\cal D}=C(\ON )$ \ref {Section 6.2/Nica/1992}.  We believe the above description of
these algebras as reduced and full partial crossed products vindicates Nica's suspicion in a telling way.  Theorem \cite
{FaithfulWiener} is due to Laca and Raeburn \ref {Proposition 2.3(3)/LacaRaeburn/1996}.

\chapter The Toeplitz C*-algebra of a graph

As we have already mentioned, C*-algebras generated by partial isometries have played a very important role in the
theory.  Besides the C*-algebras associated to semigroups of isometries treated in the previous chapter, some of the
most prevalent examples are the Cuntz-Krieger \ref {CuntzKrieger/1980}, Exel-Laca \ref {ExelLaca/1999} and the graph
C*-algebras \ref {Raeburn/2005}.

The class of Exel-Laca algebras includes the Cuntz-Krieger algebras, but the relationship between the former and graph
algebras is not as straightforward.  When a certain matrix which parametrizes Exel-Laca algebras have the property that
two distinct rows are either equal or orthogonal (a well known property of the incidence matrix of a graph), Exel-Laca
algebras are easily seen to produce all graph algebras, except for graphs containing vertices which are not the range of
any edge (known as sources).  However, up to \RME lence, the class of Exel-Laca algebras coincides with the class of
graph algebras, as recently proved by Katsura, Muhly, Sims and Tomforde \ref {KatsuraMuhlySimsTomforde/2010}.

Both Exel-Laca and graph C*-algebras may be described as partial crossed products, and in fact these algebras have been
among the motivating examples for the development of the theory of partial actions.  Since the partial crossed product
description of graph C*-algebras is a bit easier than that of Exel-Laca algebras, we will dedicate the remaining
chapters of this book to studying graph C*-algebras from the point of view of partial actions.

\medskip \fix From now on we will fix a \subjex {graph}{graph as a triple} $$ E=(E^0,E^1,r,d), $$ where $E^0$ is the set
of \subjex {vertices}{vertices in a graph}, $E^1$ is the set of \subjex {edges}{edges in a graph}, and $$ r,d :E^1 \to
E^0 $$ are the \subjex {range}{range map in a graph} and \subjex {domain}{domain, or source map in a graph} (or
\"{source}) maps, respectively.

\definition \label DefineGraphAlg The \subjex {Toeplitz C*-algebra}{Toeplitz C*-algebra of a graph} of $E$, denoted $\TE
$, is the universal C*-algebra generated by a set of mutually orthogonal projections $$ \Gen _0 = \{\p v : v\in E^0\},
$$ and a set of partial isometries $$ \Gen _1 = \{\s a : a \in E^1\}, $$ subject to the relations: \izitem \zitem $\s a
^* \s b = \bool {a =b }\ \p {d(a )}$, where the brackets correspond to Boolean value, \zitemBismark CKOrthogonal \zitem
$\s a \s a ^*\leq \p {r(a )}$, \zitemBismark CKLessThan \medskip \noindent for all $a $ and $b $ in $E^1$.  The \subjex
{graph C*-algebra of $E$}{graph C*-algebra}, denoted $C^*(E)$, is likewise defined, where in addition to the above
relations, one has \zitem $\p v = \soma {r(a )=v}\s a \s a ^*$, \zitemBismark CKSum \smallskip \noindent for all $v\in
E^0$, such that $r\inv (v)$ is finite and nonempty.

The first requirement for a C*-algebra generated by partial isometries to be tractable with our methods is that the
partial isometries involved form a tame set, so our first medium term goal will be to prove this fact.  Since $C^*(E)$
is derived from the Toeplitz algebra, we will initially concentrate our attention on $\TE $.

As usual, a (finite) \subjex {path}{path in a graph} in $E$ is defined to be a sequence $ \alpha = \alpha _1 \ldots
\alpha _n, $ where $n\geq 1$, and the $\alpha _i$ are edges in $E$, such that $d(\alpha _i) = r(\alpha _{i+1})$, for all
$i = 1,\ldots ,n-1$.  This is the usual convention when treating graphs from a categorical point of view, in which
functions compose from right to left.

The length of $\alpha $, denoted $|\alpha |$, is the number $n$ of edges in it.  As a special case we will also consider
each vertex of $E$ as a \subj {path of length zero}.

The \subjex {source}{source of a path} of a path $\alpha $ of length $n>0$ is defined by $ d (\alpha ) = d(\alpha _n), $
and its \subjex {range}{range of a path} is defined by $ r(\alpha ) = r (\alpha _1).  $ In the special case that $\alpha
$ is a path of length zero, hence consisting of a single vertex, the range and source of $\alpha $ are both defined to
be that vertex.

The set of all paths of length $n$ will be denoted by $E^n$ (this being consistent with the notations for $E^0$ and
$E^1$ already in use), and the set of all (finite) paths will be denoted by $\FinPath $.

By an \subj {infinite path} in $E$ we shall mean an infinite sequence $$ \alpha = \alpha _1\alpha _2\ldots $$ where each
$\alpha _i$ is an edge in $E$, and $d(\alpha _i) = r(\alpha _{i+1})$, for all $i\geq 1$.  The set of all infinite paths
will be denoted by $\InfPath $.  We define the range of an infinite path $\alpha $ to be $r(\alpha )=r(\alpha _1)$, and
the length of such a path to be $|\alpha |=\infty $.  There is no sensible notion of domain for an infinite path.

The set formed by all finite or infinite paths in $E$ will be denoted by $\Path $, namely $$ \Path = \FinPath \cup
\InfPath .  $$

\definition Given paths $\alpha \in \FinPath $ and $\beta \in \Path $, with $d(\alpha ) = r(\beta )$, we will denote by
$\alpha \beta $ the path defined as follows: \izitem \zitem if $|\alpha |>0$, and $|\beta |>0$, then $\alpha \beta $ is
simply the concatenation of $\alpha $ and $\beta $, \zitem if $|\alpha |=0$, and $|\beta |>0$, then $\alpha \beta =\beta
$, \zitem if $|\alpha |>0$, and $|\beta |=0$, then $\alpha \beta =\alpha $, \zitem if $|\alpha |=|\beta |=0$, in which
case $\alpha $ necessarily coincides with $\beta $ (because $\alpha = d(\alpha ) = r(\beta ) = \beta $), then $\alpha
\beta $ is defined to be either $\alpha $ or $\beta $.

Thus we see that paths of length zero get absorbed at either end of path multiplication except, of course, when both
path being multiplied have length zero, when only one of them disappears. One should note, however, that even though
paths of length zero do not show up in the resulting product, they play an important role in determining whether or not
the multiplication is defined.  In particular, if $\alpha $ is a path of length zero consisting of a single vertex $v$,
then $\alpha \beta $ is not defined, hence forbidden, unless $r(\beta )=v$.

This should be compared with the operation of composition of morphisms in a category, where vertices play the role of
identities.

\definition Given two paths $\alpha \in \FinPath $, and $\beta \in \Path $, we will say that $\alpha $ is a \subjex
{prefix}{prefix of a path} of $\beta $ if there exists a path $\gamma \in \Path $, such that $d(\alpha ) = r(\gamma )$,
and $\beta =\alpha \gamma $.

As an example, notice that for every path $\alpha $, one has that $r(\alpha )\alpha =\alpha $, so $r(\alpha )$ is a
prefix of $\alpha $.  The following is an alternate way to describe this concept:

\state Proposition \label DefinePrefix Given two paths $\alpha \in \FinPath $, and $\beta \in \Path $, one has that
$\alpha $ is a prefix of $\beta $ if and only if \izitem \zitem $r(\alpha ) = r(\beta )$, \zitem $|\alpha | \leq |\beta
|$, \zitem $\alpha _i=\beta _i$, for all $i=1,\ldots ,|\alpha |$.

Notice that when $\alpha $ and $\beta $ are paths of nonzero length, in order to check that $\alpha $ is a prefix of
$\beta $ it is enough to verify \cite {DefinePrefix/ii--iii}, since \cite {DefinePrefix/i}, follows from \cite
{DefinePrefix/iii}, with $i=1$.

On the other hand, if $|\alpha |=0$, then $\alpha $ is a prefix of $\beta $ if and only if \cite {DefinePrefix/i} holds.

Given a finite path $\alpha = \alpha _1\ldots \alpha _n$ of length $n>0$, we will denote by $\s \alpha $ the element of
$\TE $ given by $$ \s \alpha = \s {\alpha _1}\ldots \s {\alpha _n}.  $$ In the special case that $v$ is a vertex in
$E^0$, and $\alpha $ is the path of length zero given by $\alpha =v$, we will let $$ \s \alpha = \p v.  $$

\state Lemma For every $\alpha \in \FinPath $, one has \izitem \zitem $\s \alpha = \p {r(\alpha )}\s \alpha $,
\zitemBismark ProjTimesS \zitem $\s \alpha = \s \alpha \p {d(\alpha )}$, \zitemBismark STimesProj \zitem $\s \alpha \s
\alpha ^*\leq \p {r(\alpha )}$.  \zitemBismark RanLessProj

\Proof Observing that the case $|\alpha |=0$ is trivial, we assume that $|\alpha |\geq 1$.  Notice that, for every edge
$a \in E^1$, one has $$ \p {r(a )}\s a = \p {r(a )}\s a \s a ^*\s a \={CKLessThan} \s a \s a ^*\s a = \s a .  $$
Applying this to the leading edge of $\alpha $, we deduce \cite {ProjTimesSLoc}.  On the other hand, $$ \s a = \s a \s a
^*\s a \={CKOrthogonal} \s a \p {d(a )}, $$ which, applied to the trailing edge of $\alpha $, provides \cite
{STimesProjLoc}.

Regarding \cite {RanLessProjLoc}, we have $$ \s \alpha \s \alpha ^* \={ProjTimesSLoc} \p {r(\alpha )}\s \alpha \, \s
\alpha ^* \p {r(\alpha )} \leq \|\s \alpha \|^2 \p {r(\alpha )} \leq \p {r(\alpha )}.  \eqno {(\dagger )} $$ In the
first inequality above we have used the identity $$ abb^*a^*\leq \|b\|^2 aa^*, $$ known to hold in any C*-algebra, while
the second inequality in $(\dagger )$ holds because $\s \alpha $ is a product of partial isometries, each of which has
norm no bigger than 1, hence $\|\s \alpha \|\leq 1$.  \endProof

\state Remark \label NotPath \rm If $\alpha = \alpha _1\ldots \alpha _n$ is a random sequence of edges, not necessarily
forming a path, we may still define $\s \alpha $, as above.  However, unless $\alpha $ is a path, we will have $\s
\alpha =0$.  The reason is as follows: if there exists $i$ with $d(\alpha _i) \neq r(\alpha _{i+1})$, then $$ \s {\alpha
_i} \s {\alpha _{i+1}} = \s {\alpha _i} \p {d(\alpha _i)} \p {r(\alpha _{i+1})} \s {\alpha _{i+1}} = 0, $$ because the
vertex projections $\p v$ are pairwise orthogonal by construction.

\state Lemma \label DomainForPath If $\alpha $ is a finite path in $E$, then $\s \alpha ^*\s \alpha =\p {d(\alpha )}$.

\Proof This is evident if $|\alpha |=0$, while \cite {DefineGraphAlg/i} gives the case $|\alpha |=1$.  If $|\alpha |\geq
2$, then $$ \s {\alpha _1}^*\s {\alpha _1}\s {\alpha _2} = \p {d(\alpha _1)}\s {\alpha _2} = \p {r(\alpha _2)}\s {\alpha
_2} \={ProjTimesS} \s {\alpha _2}, \eqno {(\star )} $$ so $$ \s \alpha ^*\s \alpha = \s {\alpha _n}^*\ldots \s {\alpha
_2}^*\underbrace {\s {\alpha _1}^*\s {\alpha _1}\s {\alpha _2}}_{(\star )}\ldots \s {\alpha _n} = \s {\alpha _n}^*\ldots
\s {\alpha _2}^*\s {\alpha _2}\ldots \s {\alpha _n}, $$ and the result follows by induction.  \endProof

\state Lemma \label MustPrefix Given $\alpha $ and $\beta $ in $\FinPath $ such that $\s \alpha ^*\s \beta \neq 0$, then
either $\alpha $ is a prefix of $\beta $, or vice versa.

\Proof Since $ \s \beta ^*\s \alpha = (\s \alpha ^*\s \beta )^* \neq 0, $ we see that $\alpha $ and $\beta $ play
symmetric roles, so we may assume, without loss of generality, that $|\alpha |\leq |\beta |$.

Observe that $$ 0\neq \s \alpha ^*\s \beta \={ProjTimesS} \s \alpha ^*\p {r(\alpha )}\p {r(\beta )}\s \beta , $$ which
yields $\p {r(\alpha )}\p {r(\beta )}\neq 0$, and consequently $r(\alpha )=r(\beta )$.  If $|\alpha |=0$, the proof is
thus concluded, so we suppose that $|\alpha |>0$, writing $\alpha =\alpha _1\alpha _2\ldots \alpha _n$, and $\beta
=\beta _1\beta _2\ldots \beta _m$.

Notice that $\alpha _1=\beta _1$, since otherwise we would have $\s {\alpha _1}^*\s {\beta _1} = 0$ by \cite
{DefineGraphAlg/i} , which would imply $\s {\alpha }^*\s {\beta } = 0$.  Writing $v=d(\beta _1)=r(\beta _2)$, we then
have $$ 0\neq \s \alpha ^*\s \beta = \s {\alpha _n}^*\ldots \s {\alpha _2}^*\s {\alpha _1}^*\s {\beta _1}\s {\beta
_2}\ldots \s {\beta _m} \={CKOrthogonal} $$$$ = \s {\alpha _n}^*\ldots \s {\alpha _2}^*\p v\s {\beta _2}\ldots \s {\beta
_m} \={ProjTimesS} \s {\alpha _n}^*\ldots \s {\alpha _2}^*\s {\beta _2}\ldots \s {\beta _m}, $$ and the proof again
follows by induction.  \endProof

Our goal of proving the generating set of $\TE $ to be a tame set of partial isometries depends on an understanding of
the semigroup they generate.  With the above preparations we are now able to describe this semigroup.

\state Proposition \label DescriptSemiGroup Let $\Gen $ be the subset of\/ $\TE $ given by $$ \Gen = \Gen _0 \cup \Gen
_1 = \{\p v: v\in E^0\} \cup \{\s a : a \in E^1\}.  $$ Then the multiplicative sub-semigroup of\/ $\TE $ generated by
$\Gen \cup \Gen \,^* \cup \{0\}$ coincides with $ \{\s \alpha \s \beta ^*: \alpha ,\beta \in \FinPath \} \cup \{0\}.  $

\Proof We will initially prove that the set $M$ defined by $$ M = \{\s \alpha \s \beta ^*: \alpha ,\beta \in \FinPath \}
\cup \{0\}, $$ is closed under multiplication.  So pick $\alpha $, $\beta $, $\mu $ and $\nu $ in $\FinPath $, and let
us show that $$ \s \alpha \s \beta ^*\s \mu \s \nu ^*\in M.  \equationmark ClosedMulti $$

If $\s \beta ^*\s \mu =0$, then our claim follows trivially.  Otherwise $\s \beta ^*\s \mu \neq 0$, and \cite
{MustPrefix} implies that $\beta $ is a prefix of $\mu $, or vice versa.  Assuming without loss of generality that
$\beta $ is a prefix of $\mu $, we have that $\s \mu =\s \beta \s \xi $, for some path $\xi $, whence $$ \s \alpha \s
\beta ^*\s \mu \s \nu ^* = \s \alpha \s \beta ^*\s \beta \s \xi \s \nu ^* \={DomainForPath} \s \alpha \p {d(\beta )}\s
\xi \s \nu ^* \={ProjTimesS} \s \alpha \p {d(\beta )}\p {r(\xi )}\s \xi \s \nu ^* = \cdots $$ If $d(\beta )\neq r(\xi
)$, the above vanishes, in which case \cite {ClosedMulti} is proved, or else $d(\beta )=r(\xi )$, whence $\p {d(\beta
)}\p {r(\xi )} = \p {r(\xi )}$, and the above equals $$ \cdots = \s \alpha \p {r(\xi )}\s \xi \s \nu ^* = \s \alpha \s
\xi \s \nu ^* = \s {\alpha \xi }\s \nu ^* \in M.  $$

Strictly speaking, $\alpha \xi $ might not be a path, but in this case $\s {\alpha \xi }=0$, by \cite {NotPath}, so $\s
{\alpha \xi }\s \nu ^*$ lies in $M$, as well.

This concludes the proof that $M$ is a semigroup, as claimed.  Using angle brackets ``$\langle \kern 4pt \rangle $'' to
denote generated semigroup, we clearly have $$ \Gen \, \cup \Gen \, ^* \cup \{0\} \subseteq M \subseteq \big \langle
\Gen \, \cup \Gen \, ^* \cup \{0\} \big \rangle , $$ from where the proof follows easily.  \endProof

Speaking of the elements $\s \alpha \s \beta ^*$ forming the above semigroup $M$, notice that when $d(\alpha )\neq
d(\beta )$, one has $$ \s \alpha \s \beta ^* \={STimesProj} \s \alpha \p {d(\alpha )}\p {d(\beta )}\s \beta ^* = 0.  $$
So, we may alternatively describe $M$ as $$ M = \{\s \alpha \s \beta ^*: \alpha ,\beta \in \FinPath ,\ d(\alpha )\neq
d(\beta )\} \cup \{0\}.  $$

Having a concrete description of our semigroup, we may now face our first main objective:

\state Theorem For every graph $E$, the set of standard generators of $\TE $, namely $ \Gen = \Gen _0 \cup \Gen _1, $ is
a tame set of partial isometries.

\Proof By \cite {DescriptSemiGroup}, it suffices to prove that $\s \alpha \s \beta ^*$ is a partial isometry for every
$\alpha ,\beta \in \FinPath $.  We have $$ (\s \alpha \s \beta ^*) {(\s \alpha \s \beta ^*)\pilar {9pt}}^* (\s \alpha \s
\beta ^*) = \s \alpha \s \beta ^* \s \beta \s \alpha ^* \s \alpha \s \beta ^* \={DomainForPath} \s \alpha \p {d(\beta
)}\p {d(\alpha )}\s \beta ^* \$= \s \alpha \p {d(\alpha )}\p {d(\beta )}\s \beta ^* \={STimesProj} \s \alpha \s \beta
^*.  $$ This shows that $\s \alpha \s \beta ^*$ is a partial isometry, and hence that $\Gen $ is tame.  \endProof

It follows from the above that the set $\{\s a : a \in E^1\}$ is also tame.  So, letting $\F $ be the free group on the
set $E^1$, we may invoke \cite {PisoSet} to conclude that there exists a unique semi-saturated *-partial representation
$$ \Pr : \F \to \TTE , \equationmark PRepInTE $$ such that $ \pr a = \s a , $ for all $a \in E^1$, where $\TTE $ is the
algebra obtained from $\TE $ by adding a unit to it, even if it already has one\fn {It is easy to see that $\TE $ is
unital if and only if $E^0$ is finite, in which case the unit of $\TE $ is given by the sum of all the $\p v$, with $v$
ranging in $E^0$.}.

As usual let us denote by $\e g =\pr g \pri g$.

\state Proposition \label RelationsForTE The partial representation described in \cite {PRepInTE} satisfies, for every
$a $ and $b $ in $E^1$, \izitem \zitem $d(a )=d(b ) \imply \ei a = \ei b $, \zitem $d(a )\neq d(b ) \imply \ei a \ei b =
0$, \zitem $r(a )=d(b ) \imply \e a \leq \ei b $, \zitem $a \neq b \phantom {d()d()} \imply \e a \e b = 0$.

\Proof Noticing that $$ \ei a = \pri a \pr a = \s a ^*\s a = \p {d(a )}, $$ points (i) and (ii) follow at once.  As for
(iii) we have $$ \e a = \pr a \pri a = \s a \s a ^* \_\leq {CKLessThan} \p {r(a )} = \p {d(b )} = \ei b .  $$ Finally,
if $a \neq b $, then $$ \e a \e b = \pr a \pri a \pr b \pri b = \s a \s a ^* \s b \s b ^* \={CKOrthogonal} 0.
\omitDoubleDollar \endProof

In order to avoid technical details, we will now restrict ourselves to studying graphs for which $\TE $ is generated by
the $\s a $ alone, a property that follows from the absence of sinks, as defined below:

\definition \label DefineSinkSourceRegular A vertex $v\in E^0$ is said to be a: \izitem \zitem \subjex {sink}{sink
vertex}, if $d\inv (v)$ is the empty set, \zitem \subjex {source}{source vertex}, if $r\inv (v)$ is the empty set,
\zitem \subj {regular vertex}, if $r\inv (v)$ is finite and nonempty.

If the vertex $v$ is not a sink, then there exists some edge $a $ in $E^1$ such that $d(a )=v$.  Therefore \cite
{CKOrthogonal} implies that $$ \s a ^* \s a = \p {d(a )} = \p v, $$ whence $\p v$ lies in the algebra generated by the
$\s a $.

\medskip \fix From now on we will concentrate on graphs without sinks, so we suppose throughout that $E$ has no sinks.

In this case, as seen above, all of the $\p v$ may be expressed in terms of the $\s a $, which is to say that $$ \Gen _1
= \{\s a : a \in E^1\} $$ generates $\TE $, as a C*-algebra.

\state Theorem \label UniversalPropTTE Let $E$ be a graph without sinks, and let $\F $ be the free group on the set
$E^1$, equipped with the usual word-length function.  Also let \izitem \zitem $\Rel '$ be the set consisting of
relations \cite {DefPR/i--iv}, \zitem $\SSat $ be the set of relations described in \cite {SemiSatVsConvex}, and \zitem
$\Rel _E$ be the set consisting of relations \cite {RelationsForTE/i--iv}, for each $a $ and $b $ in $E^1$.  (By this we
mean that we take the relation on the right-hand-side of ``$\,\Rightarrow $'', whenever the condition in the
left-hand-side holds).  \smallskip \noindent Then $\TTE $ is *-isomorphic to the universal unital C*-algebra generated
by a set $ \{\pr g: g\in \F \}, $ subject to the set of relations $\Rel ' \cup \SSat \cup \Rel _E$.

\Proof It clearly suffices to prove that $\TTE $ possesses the corresponding universal property.  As a first step,
observe that the elements in the range of the partial representation $\Pr $ given by \cite {PRepInTE} satisfy $\Rel '$
for obvious reasons, satisfy $\SSat $ because $\Pr $ is semi-saturated, and satisfy $\Rel _E$ by \cite {RelationsForTE}.

Since $\s a = \pr a $, for every $a $ in $E^1$, we see that the range of $\Pr $ contains $\Gen _1$.  Moreover, for every
$v\in E^0$, one may pick $a \in E^1$ such that $d(a )=v$, because $E$ has no sinks, whence $$ \p v = \s a ^*\s a = \pri
a \pr a , $$ and so we see that the range of $\Pr $ also contains $\Gen _0$, hence generating $\TTE $ as a C*-algebra.

Let us now assume we are given a unital C*-algebra $B$ which contains a set $ \{\hat \pr g: g\in \F \}, $ satisfying
relations $\Rel ' \cup \SSat \cup \Rel _E$.  As usual, for each $g$ in $\F $, we will let $\hat \e g = \hat \pr g\hat
\pri g$.

For each $a$ in $E^1$, define, $\hat \s a = \hat \pr a $, and for every $v\in E^0$, choose an edge $b \in E^1$, with
$d(b )=v$, and put $\hat \p v = \hat \ei b $.  Notice that such an edge exists because $v$ is not a sink.  Moreover, if
$a $ is another edge with $d(a )=v$, then $\hat \ei b = \hat \ei a $ by \cite {RelationsForTE/i}, a relation which is
part of $\Rel _E$, and is therefore satisfied by the $\hat \pr g$.

We now claim that the sets $$ \widehat {\Gen }_0 = \{\hat \p v : v\in E^0\} \and \widehat {\Gen }_1 = \{\hat \s a : a
\in E^1\}, $$ satisfy conditions \cite {DefineGraphAlg/i\&ii}.  In fact, since $\hat \Pr $ is clearly a partial
representation, we have that each $\hat \s a $ is a partial isometry by \cite {PrepFormedByPisos}, and hence the $\hat
\p v$ are projections.

If $v_1,v_2\in E^0$ are two distinct vertices, choose $a _1,a _2\in E^1$, with $d(a _i) = v_i$.  Then $$ \hat \p {v_1}
\hat \p {v_2} = \hat \ei {a _1} \hat \ei {a _2} = 0, $$ by \cite {RelationsForTE/ii}, so the $\p v$ are pairwise
orthogonal, as required.

In order to check \cite {CKOrthogonal}, let $a ,b \in E^1$.  Then $$ \hat \s a ^* \hat \s b = \hat \pri a \hat \pr b =
\hat \pri a \hat \pr a \hat \pri a \hat \pr b \hat \pri b \hat \pr b = \hat \pri a \hat \e a \hat \e b \hat \pr b .  $$
If $a \neq b $, the above vanishes by \cite {RelationsForTE/iv}.  Otherwise, we have $$ \hat \s a ^* \hat \s b = \hat
\pri a \hat \pr a = \hat \ei a = \hat \p {d(a )}.  $$

Let us now prove \cite {CKLessThan}, so pick $a \in E^1$.  Since $E$ has no sinks, there exists some $b \in E^1$, such
that $ d(b ) = r(a ).  $ Then $$ \hat \s a \hat \s a ^* = \hat \pr a \hat \pri a = \hat \e a \_\leq {RelationsForTE/iii}
\hat \ei b = \hat \p {d(b )} = \hat \p {r(a )}.  $$

Therefore, by the universal property of $\TE $, we conclude that there exists a *-homomor\-phism $ \varphi :\TE \to B, $
such that $$ \varphi (\p v)=\hat \p v \and \varphi (\s a )=\hat \s a , \equationmark OkOnGenerators $$ for all $v\in
E^0$, and all $a \in E^1$.  By mapping identity to identity, we may extend $\varphi $ to a unital map from $\TTE $ to
$B$.

The proof will then be concluded once we prove that $\varphi (\pr g) = \hat \pr g$, for all $g$ in $\F $, which we do by
induction on $|g|$.  The case $|g|=0$, follows since $\varphi $ preserve identities, while the case $|g|=1$, follows
from \cite {OkOnGenerators}.

If $|g|\geq 2$, write $g = g_1g_2$, with $ |g| = |g_1|+|g_2|, $ and $|g_1|,|g_2|<|g|$.  Then, since both $\Pr $ and
$\hat \Pr $ are semi-saturated, we have by induction that $$ \varphi (\pr g) = \varphi (\pr {g_1g_2}) = \varphi (\pr
{g_1}\pr {g_2}) = \varphi (\pr {g_1})\varphi (\pr {g_2}) = \hat \pr {g_1}\hat \pr {g_2} = \hat \pr {g_1g_2} = \hat \pr
{g}.  $$ This concludes the proof.  \endProof

According to \cite {DefineCstarParRel}, our last result may then be interpreted as saying that $$ \TTE \simeq
\CstarparRel \F {\SSat \cup \Rel _E}, $$ so we may apply \cite {CstarRelAsCP} in order to describe $\TTE $ as a partial
crossed product.

\state Corollary \label TTEAsCrossProd Given a graph $E$ without sinks, denote by $\OE $ the spectrum of the set of
relations $\SSat \cup \Rel _E$, and consider the partial action of\/ $\F $ on $\OE $ obtained by restricting the partial
Bernoulli action.  Then there is a *-isomor\-phism $$ \varphi : \TTE \to C(\OE )\rt \Th \F , $$ such that $$ \varphi (\s
a ) = 1_a \delta _a \for a \in E^1.  $$

\bigskip In what follows we will denote the already mentioned restriction of the partial Bernoulli action to $\OE $ by
$$ \Th _E = \big (\{\D g^E \}_{g\in G}, \{\th g^E \}_{g\in G}\big ), \equationmark BernouliTTE $$

By direct inspection, using \cite {DefineOmegaR}, it is easy to see that $$ \omega _0 = \{1\}, \equationmark
TrivalElement $$ is an element of $\OE $.  Since $\omega _0$ does not lie in any $\D g^E$ (see \cite
{DescriptBernouliDoms}), for $g\neq 1$, we have that $\{\omega _0 \}$ is an invariant subset of $\OE $, whence its
complement is an open invariant set.

This moreover implies that $C_0\big (\OE \setmenos \{\omega _0\}\big ) $ is an invariant ideal of $C(\OE )$, whence
$C_0\big (\OE \setmenos \{\omega _0\}\big ) \rt \tau \F $ is an ideal in $C(\OE )\rt \tau \F $ by \cite {ExactSeqCross}.

\state Theorem \label OnlyTEAsCrossProd In the context of \cite {TTEAsCrossProd}, the restriction of $\varphi $ to $\TE
$ is an isomorphism from $\TE $ onto $C_0\big (\OE \setmenos \{\omega _0\}\big )\rt \tau \F $.

\Proof As already seen, we have that $\omega _0\notin \D a ^E$, for every edge $a $, whence $\D a ^E\subseteq \OE
\setmenos \{\omega _0\}$.  Consequently $$ 1_a \delta _a \in C_0\big (\OE \setmenos \{\omega _0\}\big )\rt \tau \F , $$
showing that $\varphi (\TE )$ is contained in $C_0\big (\OE \setmenos \{\omega _0\}\big )\rt \tau \F $.  By \cite
{ExactSeqCross} we have that $$ {C(\OE )\rt \tau \F \over C_0\big (\OE \setmenos \{\omega _0\}\big )\rt \tau \F }\simeq
C(\{\omega _0\})\rt ?\F , $$ which is a one-dimensional algebra given that the intersection of the $\D g^E$ with
$\{\omega _0\}$ is empty, except for when $g=1$.  Thus, the co-dimension of $C_0\big (\OE \setmenos \{\omega _0\}\big
)\rt \tau \F $ in $C(\OE )\rt \tau \F $ is the same as the co-dimension of $\TE $ in $\TTE $, both being equal to 1.
Since $\varphi $ is an isomorphism, one then sees that $\varphi (\TE )$ must coincide with in $C_0\big (\OE \setmenos
\{\omega _0\}\big )\rt \tau \F $.  \endProof

\nrem In 1980, Enomoto and Watatani \ref {EnomotoWatatani/1980} realized that certain properties of the recently
introduced Cuntz-Krieger C*-algebras \ref {CuntzKrieger/1980} could be described by notions from graph theory.
Seventeen years later, Kumjian, Pask, Raeburn and Renault \ref {KumjianPaskRaeburnRenault/1997} began studying
C*-algebras for certain infinite graphs, in the wake of which many authors established an intense area of research going
by the name of ``graph C*-algebras''.  Definition \cite {DefineGraphAlg} first appeared in \ref
{FowlerLacaRaeburn/2000}.  It was inspired by ideas from \ref {KumjianPaskRaeburnRenault/1997} and \ref {ExelLaca/1999}.
The reader will find a long list of references for the theory of graph C*-algebras in \ref {Raeburn/2005}.

\chapter Path spaces

The definition of the spectrum of a set of relations given in \cite {DefineOmegaR} is a one-size-fits-all description,
having many applications to concrete examples, such as the various C*-algebras associated to semigroups of isometries
and the Toeplitz C*-algebra of a graph studied above.  However experience shows that, once we have focused on a specific
example, the description of $\OR $ given in its formal definition may often be greatly simplified, leading to a much
more intuitive and geometrically meaningful picture of this space.

The purpose of this chapter is thus to develop a careful analysis of the spectrum $\OE $ of the relations leading up to
$\TTE $, as described in \cite {TTEAsCrossProd}.

\medskip \fix As above, we will continue working with a fixed graph $E$ without sinks.  Since $\OE $ is a subset of
$\{0,1\}^\F = \Part (\F )$, each $\omega \in \OE $ may be viewed as a subset of $\F $.  A useful mental picture of an
element $\omega $ in $\OE $ is thus to imagine the vertices of the Cayley graph of $\F $ painted black whenever they
correspond to an element in $\omega $.

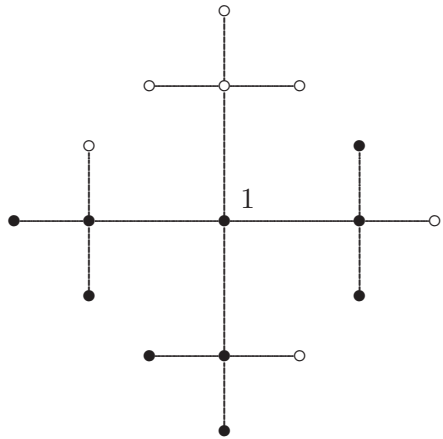
\null \hfill \beginpicture \setcoordinatesystem units <0.0020truecm, 0.0020truecm> point at 0 0 \newcount \x \newcount
\y \newcount \z \newcount \w \newcount \u \newcount \v \newcount \h

\def \shorten #1#2{\u = #1 \v = #2 \ifnum \u = \v \h = 0 \else \ifnum \v > \u \h = 1 \else \h = -1 \fi \fi \multiply \h
by 42 \advance \u by \h \advance \v by -\h }

\def \ponto #1 #2 #3 #4 #5;{\put {$#1$} at #4 #5 \shorten {#2}{#4}\x =\u \z =\v \shorten {#3}{#5}\y =\u \w =\v \plot
{\number \x } {\number \y } {\number \z } {\number \w } /}

\put {$\bullet $} at 0 0 \put {$1$} <8pt,8pt> at 0 0

\ponto {\bullet } 0 0 900 0; \ponto {\circ } 900 000 1400 000; \ponto {\bullet } 900 000 900 500; \ponto {\bullet } 900
000 900 -500;

\ponto {\circ } 000 000 000 900; \ponto {\circ } 000 900 000 1400; \ponto {\circ } 000 900 500 900; \ponto {\circ } 000
900 -500 900;

\ponto {\bullet } 000 000 000 -900; \ponto {\bullet } 000 -900 000 -1400; \ponto {\bullet } 000 -900 -500 -900; \ponto
{\circ } 000 -900 500 -900;

\ponto {\bullet } 000 000 -900 000; \ponto {\bullet } -900 000 -1400 000; \ponto {\circ } -900 000 -900 500; \ponto
{\bullet } -900 000 -900 -500;

\endpicture \hfill \null \bigskip

Observing that $\OE $ is a subspace of the space $\OuG $ defined in \cite {DefineOuG}, every $\omega \in \OE $ is a
subset of $\F $ containing the unit group element, so the vertex of the Cayley graph corresponding to 1 will always be
painted black.  Moreover, the presence of relations $\SSat $ imply that every $\omega $ is convex, according to \cite
{SemiSatVsConvex}, and one may easily show that the abstract notion of convexity given in \cite {DefineConvex/b}
coincides with the corresponding geometric notion.

Our alternative description of the elements $\omega $ belonging to $\OE $ will be aided by the following concept:

\definition Let $\omega \in \OuG $, and let $g\in \omega $.  The \subj {local configuration} of $\omega $ at $g$ is the
set $$ \locstd = \{h\in \F : |h|=1,\ gh\in \omega \}.  $$

Since the vertices in the Cayley graph adjacent to a given $g$ are precisely those of the form $gh$, with $|h|=1$, the
local configuration of $\omega $ at $g$ determines which of the adjacent vertices correspond to elements belonging to
$\omega $.  Also notice that the condition $|h|=1$ implies that $h$ is either an edge or the inverse of an edge.
Therefore $\locstd \subseteq E^1\cup (E^1)\inv $, where the inverse is evidently taken with respect to the group
structure of $\F $.

\state Proposition \label LocalConditions A given $\omega \in \OuG $ lies in $\OE $ if and only if $\omega $ is convex,
and for all $g$ in $\omega $, the local configuration of $\omega $ at $g$ is of one of the following mutually exclusive
types: \iaitem \aitem there is an edge $a \in E^1$, such that $$ \locstd = \{a \} \cup \{b \inv : b\in E^1,\ d(b ) = r(a
)\}, $$ \aitem there is a vertex $v\in E^0$, such that $$ \locstd = \phantom {\{a \} \cup } \{b \inv : b\in E^1,\ d(b )
= v\}, $$ \aitem $\locstd $ is the empty set.

\Proof Let us first prove the only if part, so let us suppose that $\omega \in \OE $.  Since the relations involved in
the definition of $\OE $ (see \cite {TTEAsCrossProd}) include $\SSat $, we have that $\OE \subseteq \ORSat $, so $\omega
$ is in $\ORSat $ and then \cite {SemiSatVsConvex} implies that $\omega $ is convex.

Let us now suppose that we are given $g$ in $G$, and there exists an edge $a \in \locstd $.  If so, we claim that the
local configuration of $\omega $ at $g$ is of type (a).  In fact, suppose by contradiction that some other edge $b $
lies in $\locstd $.  Then, taking relation \cite {RelationsForTE/iv} into account, one sees that the function $$ f = \ee
a \ee b $$ (see \cite {DefineOmegaR}) lies in ${\cal F}_{\Rel _E}$, whence $$ 0 = f(g\inv \omega ) = \bool {a \in g\inv
\omega } \bool {b \in g\inv \omega } = \bool {ga \in \omega } \bool {gb \in \omega }, $$ and since $ga \in \omega $, it
follows that $gb \notin \omega $, whence $b\notin \locstd $.  This proves that $a $ is the only edge (as opposed to the
inverse of an edge) lying in $\locstd $.

Let us now prove that $b \inv \in \locstd $, whenever the edge $b $ is such that $d(b )=r(a )$.  Under these conditions
relation \cite {RelationsForTE/iii} is in $\Rel _E$, so the function $$ f = \ee a - \ee a \ee {b \inv } = \ee a (1 - \ee
{b \inv }) $$ lies in ${\cal F}_{\Rel _E}$, whence $$ 0 = f(g\inv \omega ) = \bool {a \in g\inv \omega }\big (1-\bool {b
\inv \in g\inv \omega }\big ) \$= \bool {ga \in \omega }\big (1-\bool {gb \inv \in \omega }\big ).  $$ Since $ga \in
\omega $, it follows that $gb \inv \in \omega $, so $b\inv \in \locstd $, as desired.

To conclude the proof of the claim that the local configuration of $\omega $ at $g$ is of type (a), it now suffices to
prove that $b \inv \notin \locstd $, whenever $d(b )\neq r(a )$.  Using that $r(a )$ is not a sink, choose $c \in E^1$
such that $d(c ) = r(a )$, so $c \inv \in \locstd $ by the previous paragraph.  We then have that $d(c ) \neq d(b )$, so
relation \cite {RelationsForTE/ii} applies and letting $$ f = \ee {b \inv } \ee {c \inv }, $$ we have $$ \vrule height
10pt depth 6pt width 0pt 0= f(g\inv \omega ) = \bool {b \inv \in g\inv \omega } \bool {c \inv \in g\inv \omega } = \bool
{gb \inv \in \omega } \bool {gc \inv \in \omega }.  $$

Since $gc \inv \in \omega $, as seen above, we have that $gb \inv \notin \omega $, and consequently $b \inv \notin
\locstd $.

This proves the claim under the assumption that $\locstd $ contains an edge, so let us suppose the contrary.  If
$\locstd $ is empty, the proof is over, so we are left with the case in which $b \inv $ lies in $\locstd $ for some edge
$b $.  Denoting by $v=d(b )$, we have for an arbitrary edge $c $ that $$ \bool {gc \inv \in \omega } = \bool {gc \inv
\in \omega } \bool {gb \inv \in \omega } = \ee {c \inv }(g\inv \omega ) \ee {b \inv }(g\inv \omega ) = \cdots $$
Considering the appropriate function in ${\cal F}_{\Rel _E}$ related to \cite {RelationsForTE/i--ii}, according to
whether or not $d(c) = d(b)$, the above equals $$ \cdots = \bool {d(c )=d(b )}\, \ee {b \inv }(g\inv \omega ) = \bool
{d(c )=v}, $$ which says that $c \in \locstd $ if and only if $d(c )=v$.  This proves that $\locstd $ is of type (b),
and hence the proof of the only if part is concluded.

The proof of the converse essentially consists in reversing the arguments above and is left to the reader.  \endProof

As already observed in \cite {TrivalElement}, the element $ \omega _0 = \{1\} $ is a member of $\OE $.  By convexity,
this is the only element displaying a local configuration of type \cite {LocalConditions/c}.

By inspection of the possible local configuration types above one may easily prove the following:

\state Proposition \label LocalConditionsInterpreted Let $\omega $ be in $\OE $, and let $g$ be in $\omega $.  \izitem
\zitem If\/ $\locstd $ contains an edge $a $, then $$ b \inv \in \locstd \iff d(b ) = r(a ) \for b \in E^1.  $$ \zitem
If\/ $\locstd $ contains $c \inv $, for some edge $c $, then $$ b \inv \in \locstd \iff d(b ) = d(c ) \for b \in E^1.
$$

We now begin to work towards giving a description of the elements of $\OE $ based on local configurations.

\state Lemma \label UmPuxaOOutro Let $\omega $ be in $\OE $, and let $g$ be in $\omega $.  Given a finite path $\nu $,
with $m:=|\nu |\geq 1$, and $\nu _m\inv \in \locstd $, then $g\nu \inv \in \omega $.

\Proof We prove by (backwards) induction that $\mu \nu _m\inv \nu _{m-1}\inv \ldots \nu _k\inv \in \omega $, for every
$k$.  By the definition of local configurations it is immediate that $g\nu _m\inv \in \omega $.  Now supposing that
$k_0\leq m$, and $$ \mu \nu _m\inv \nu _{m-1}\inv \ldots \nu _k\inv \in \omega \for k\geq k_0, $$ let $h=\mu \nu _m\inv
\nu _{m-1}\inv \ldots \nu _{k_0}\inv $.  Then $\nu _{k_0}\in \loc \omega h$, and $d(\nu _{k_0-1}) = r(\nu _{k_0})$, so
$\nu _{k_0-1}\inv $ is in $\loc \omega h$ by \cite {LocalConditionsInterpreted/i}, meaning that $$ \omega \ni h\nu
_{k_0-1}\inv = \mu \nu _m\inv \nu _{m-1}\inv \ldots \nu _{k_0}\inv \nu _{k_0-1}\inv , $$ as desired.  \endProof

In our next main result we will use the above local description of the elements in $\OE $ to give a parametrization of
$\OE $ by the set of all (finite and infinite) paths.  In order to do this we must first fine tune our interpretation of
finite paths as elements of $\F $: if $\alpha $ is a finite path of positive length, then we may write $\alpha = \alpha
_1\alpha _2\ldots \alpha _n$, and there is really only one way to see $\alpha $ as an element in $\F $, namely as the
product of the $\alpha _i$.  However, when $\alpha $ has length zero then it consists of a single vertex and we will
adopt the perhaps not so obvious convention that $\alpha $ represents the unit element of $\F $.  This is not to say
that we are identifying all paths of length zero with one another, the convention only applying when paths are seen as
elements of $\F $.

\medskip The correspondence, about to defined, of elements of $\F $ with paths in $E$, will send each path $\alpha $ to
the subset $\omega _\alpha $ introduced below:

\definition \label DefineOAlp Given $\alpha \in \Path $, let $\omega _\alpha $ be the subset of $\F $ defined by $$
\omega _\alpha = \{\mu \nu \inv : \mu ,\nu \in \FinPath ,\ d(\mu )=d(\nu ),\ \mu \hbox { is a prefix of } \alpha \}.  $$

Taking $\mu $ to be the range of $\alpha $, considered as a path of length zero, as well as a prefix of $\alpha $, one
sees that $\omega _\alpha $ contains every element of the form $\nu \inv $, where $\nu $ is a path with $d(\nu ) =
r(\alpha )$.  Likewise, if $\mu $ is any prefix of $\alpha $, then, taking $\nu $ to be the path of length zero
consisting of the vertex $d(\mu )$, we see that $\mu $ is in $\omega _\alpha $.

Let us now discuss the reduced form of elements in $\omega _\alpha $.

\state Proposition \label ReducedFormOmegaAlpha For any $\alpha \in \Path $, and for every $g\in \omega _\alpha $, there
are $\mu ,\nu \in \FinPath $ such that \izitem \zitem $g=\mu \nu \inv $, \zitem $\ell (g) = |\mu | + |\nu |$, \zitem
$d(\mu )=d(\nu )$, and \zitem $\mu $ is a prefix of $\alpha $.

\Proof Given $g=\mu \nu \inv \in \omega _\alpha $, as in \cite {DefineOAlp}, there is nothing to do in case $|\mu |$ or
$|\nu |$ vanish.  So we may assume that $|\mu |$ and $|\nu |$ are both at least equal to $1$, and we may then write
write $\mu =\mu _1\mu _2\ldots \mu _n$, and $\nu =\nu _1\nu _2\ldots \nu _m$, with the $\mu _i$ and the $\nu _j$ in
$E^1$, whence $$ g = \mu _1\mu _2\ldots \mu _{n-1}\mu _n\nu _m\inv \nu _{m-1}\inv \ldots \nu _2\inv \nu _1\inv .  $$ If
$\mu _n\neq \nu _m$, then $g$ is in reduced form, so $\ell (g) = |\mu | + |\nu |$, as needed.

Otherwise let us use induction on $|\mu | + |\nu |$.  So, supposing that $\mu _n=\nu _m$, we may clearly cancel out the
term ``$\mu _n\nu _m\inv $'' above, and we are then left with $g=\mu '{\nu '}\inv $, where $\mu '=\mu _1\mu _2\ldots \mu
_{n-1}$, and $\nu '=\nu _1\nu _2\ldots \nu _{m-1}$.  In this case, notice that $$ d(\mu ') = d(\mu _{n-1}) = r(\mu _n) =
r(\nu _m) = d(\nu _{m-1}) = d(\nu '), $$ so $\mu '$ and $\nu '$ satisfy the conditions required of $\mu $ and $\nu $ in
the definition of $\omega _\alpha $, and the proof follows by induction.

In fact we still need to take into account the case in which $m$ or $n$ coincide with 1, since e.g.~when $m=1$, there is
no $\nu _{m-1}$.  Under this situation we may still cancel out the term ``$\mu _n\nu _m\inv $'', so that $g=\mu '{\nu
'}\inv $, where $\mu '$ is as above and $\nu '$ is the path of length zero consisting of the vertex $d(\mu _{n-1})$.

If, instead, one had $n=1$, we could take $\mu '=d(\nu _{m-1})$ and $\nu '$, as above, observing that $$ \mu ' = d(\nu
_{m-1}) = r(\nu _m) = r(\mu _1) = r(\mu ) = r(\alpha ), $$ so $\mu '$ is still a prefix of $\alpha $ and we could
proceed as above. Of course this still leaves out the case in which $n=m=1$, but then $g=1$ and it is enough to take
$\mu =\nu =r(\alpha )$.  \endProof

As we will soon find out, the elements of $\F $ of the form described in \cite {ReducedFormOmegaAlpha} have a special
importance to us.  We shall thus introduce a new concept to highlight these elements.

\definition \label DefineStdForm Given $g$ in $\F $ of the form $g=\mu \nu \inv $, where $\mu ,\nu \in \FinPath $, we
will say that $g$ is in \subj {standard form} if conditions \cite {ReducedFormOmegaAlpha/ii--iii} are satisfied.

The following is the main technical result leading up to our concrete description of $\OE $.

\state Theorem \label OneOneConfigPath For each finite or infinite path $\alpha $, one has that $\omega _\alpha $
belongs to $\OE $.  In addition, the correspondence $ \alpha \mapsto \omega _\alpha $ is a one-to-one mapping from
$\Path $ onto $\OE \setmenos \{\omega _0\}$, where $\omega _0$ is defined in \cite {TrivalElement}.

\Proof By checking convexity and analyzing local configurations it is easy to use \cite {LocalConditions} in order to
prove that $\omega _\alpha $ is in $\OE \setmenos \{\omega _0\}$, for every $\alpha $ in $\Path $.  It is also evident
that the mapping referred to in the statement is injective.

In order to prove surjectivity, let $\omega \in \OE \setmenos \{\omega _0\}$, and define $$ \omega _+ = \omega \cap \F
_+, $$ where $\F _+$ is the positive cone\fn {$\F _+$ is defined as the sub-semigroup of $\F $ generated by $E^1$.  We
assume that $\F _+$ also contains the unit of $\F $.}  of $\F $.  We claim that any $$ \alpha =\alpha _1\alpha _2\ldots
\alpha _n\in \omega _+, $$ where each $\alpha _i$ is an edge, and $n\geq 1$, is necessarily a path.  In fact, for every
$i=1,\ldots ,n-1$, notice that, by convexity, $g:=\alpha _1\alpha _2\ldots \alpha _i\in \omega _+$, and $$ \alpha _i\inv
,\alpha _{i+1}\in \locstd .  $$ The local configuration of $\omega $ at $g$ is therefore of type \cite
{LocalConditions/a}, since it contains an edge, whence $d(\alpha _i) = r(\alpha _{i+1})$, proving that $\alpha $ is
indeed a path.

We next claim that, given any two finite paths $\alpha ,\beta \in \omega _+$, with $1\leq |\alpha |\leq |\beta |$, then
necessarily $\alpha $ is a prefix of $\beta $.  Otherwise let $i$ be the smallest integer such that $\alpha _i\neq \beta
_i$.  Again by convexity we have that $$ g := \alpha _1\alpha _2\ldots \alpha _{i-1}=\beta _1\beta _2\ldots \beta
_{i-1}\in \omega _+, $$ (in case $i=0$, then $g$ is to be interpreted as 1), and $\alpha _i,\beta _i\in \locstd $.

Including two different edges, $\locstd $ fails to be of any of the types described in \cite {LocalConditions}, thus
bringing about a contradiction.  This proves that $\alpha $ is a prefix of $\beta $, as desired.

We will now build a path $\alpha $ which will later be proven to satisfy $\omega =\omega _\alpha $.  \izitem \zitem
Assuming that the set of finite paths of positive length in $\omega _+$ is finite and nonempty, we let $\alpha $ be the
longest such path.  \zitem Should there exist arbitrarily long paths in $\omega _+$, we let $\alpha $ be the infinite
path whose prefixes are the finite paths in $\omega $.  \zitem In case $\omega $ contains no path of positive length,
then the local configuration of $\omega $ at $1$ is necessarily of type \cite {LocalConditions/b} (it cannot be of type
\cite {LocalConditions/c} since $\omega $ is convex and $\omega \neq \omega _0$).  In this case we put $\alpha = v$,
where $v$ is the vertex referred to in \cite {LocalConditions/b}.

\bigskip As already indicated, we will prove that $\omega _\alpha =\omega $.  By construction it is clear that $\omega
_\alpha \cap \F _+=\omega \cap \F _+$.  Therefore, for any $g$ in $\F _+$, we have that $g\in \omega $, if and only if
$g$ is a prefix of $\alpha $.

In order to prove that $\omega _\alpha \subseteq \omega $, pick any $g\in \omega _\alpha $, and write $g=\mu \nu \inv $,
with $\mu $ a prefix of $\alpha $ and $\nu $ a path with $d(\mu )=d(\nu )$.

If $|\nu |=0$, then $g$ is a prefix of $\alpha $, so $g$ is in $\omega $.  We therefore suppose that $$ \nu =\nu _1\nu
_2\ldots \nu _m, $$ where $m\geq 1$.

\medskip \noindent {\tensc case 1.}\enspace Assuming that $n:=|\mu |\geq 1$, we have that $\mu $ is in $\omega $.
Moreover $\mu _n\inv \in \loc \omega \mu $, and $d(\mu _n) = d(\mu ) = d(\nu ) = d(\nu _m)$.  Therefore $\nu _m\inv \in
\locstd $ by \cite {LocalConditionsInterpreted/ii}.  Employing \cite {UmPuxaOOutro} we then have that $g=\mu \nu \inv
\in \omega $.

\medskip \noindent {\tensc case 2.}\enspace Assuming that $|\mu |=0$, we have that $\mu =v$, for some vertex $v$.  The
fact that $\mu $ is a prefix of $\alpha $ in this case means that $v=r(\alpha )$, therefore $$ d(\nu _m) =d(\nu )=d(\mu
)=v=r(\alpha ).  $$ Temporarily assuming that $|\alpha |\geq 1$, we then have that $\alpha _1\in \loc \omega 1$, and
$d(\nu _m) = r(\alpha _1)$, so $\nu _m\inv \in \loc \omega 1$, by \cite {LocalConditionsInterpreted/i}.  Another
application of \cite {UmPuxaOOutro} then gives $$ \omega \ni 1\nu \inv = \mu \nu \inv = g.  $$

Still under case (2), but now assuming that $|\alpha |=0$, write $\alpha =v$, for some vertex $v$.  This implies that we
are under situation (iii) above, whence the local configuration of $\omega $ at $1$ is of type \cite
{LocalConditions/b}.  We then have by definition that $\loc \omega 1 = \{b \inv : b\in E^1,\ d(b ) = v\}$.  Observe that
$$ d(\nu _m)= d(\nu ) =d(\mu ) = d(\alpha ) = v, $$ so $\nu _m\inv \in \loc \omega 1$, and again by \cite {UmPuxaOOutro}
we conclude that $ \omega \ni 1\nu \inv =\mu \nu \inv = g.  $

This concludes the proof that $\omega _\alpha \subseteq \omega $, so we are left with the task of proving the reverse
inclusion.  Thus, let $g\in \omega $, and write $$ g=x_1x_2\ldots x_n, $$ in reduced form. We claim that there is no $i$
for which $x_i\in (E^1)\inv $ and $x_{i+1}\in E^1$.  Otherwise, write $x_i=a \inv $, and $x_{i+1}=b $, with $a $ and $b
$ in $E^1$, and let $g=x_1x_2\ldots x_i$.  Then, by convexity, $$ g,\ ga ,\ gb \in \omega , $$ whence $a ,b \in \loc
\omega g$.  This results in a local configuration with two distinct edges, contradicting \cite {LocalConditions}.
Therefore we see that the factors of the reduced form of $g$ lying in $E^1$ must be to the left of those lying in
$(E^1)\inv $.  In other words $$ g = \mu \nu \inv , $$ where $\mu $ and $\nu $ are in $\F _+$.  One may now easily
employ \cite {LocalConditionsInterpreted} to prove that $\mu $ and $\nu $ are paths, and also that $d(\mu )=d(\nu )$.
By convexity we see that $\mu $ is in $\omega $, hence $\mu $ is a prefix of $\alpha $.  Consequently $g\in \omega
_\alpha $, thus verifying that $\omega \subseteq \omega _\alpha $, and hence concluding the proof.  \endProof

In order to obtain a model for $\OE $ in terms of paths which accounts for the exceptional element $\omega _0$, we make
the following:

\definition \iaitem \aitem We fix any element in the universe not belonging to $\Path $, denote it by $\varnothing $,
and call it the \subj {empty path}\fn {The empty path should not be confused with any of the paths of length zero we
have so far been working with.}.  \aitem No path in $\Path $ is considered to be a prefix of $\varnothing $.  \aitem The
\subj {full path space} of $E$ is the set $ \FullPath = \Path \cup \{\varnothing \}.  $ \aitem Given any $\omega \in \OE
\setmenos \{\omega _0\}$, the \subj {stem} of $\omega $, denoted $\sigma (\omega )$, is the unique element $\alpha $ in
$\Path $, such that $\omega _\alpha =\omega $, according to \cite {OneOneConfigPath}.  \aitem The stem of $\omega _0$ is
defined to be the empty path $\varnothing $.

The stem may then be seen as a function $$ \sigma : \OE \to \FullPath , $$ which is bijective thanks to \cite
{OneOneConfigPath}.

The restriction of the partial Bernoulli action to $\OE $, which gives rise to $\TTE $ according to \cite
{TTEAsCrossProd}, may therefore be transferred via the stem function over to $\FullPath $.  It is our next goal to give
a concrete description of this partial action.  We begin with a technical result.

\state Lemma \label PathCodingOfDg Let $g$ be any element of\/ $\F \setmenos \{1\}$.  Then: \izitem \zitem If $\D g^E$
is nonempty, then $g$ admits a standard form $g = \mu \nu \inv $ (see Definition \cite {DefineStdForm}), \zitem If
$g=\mu \nu \inv $ is in standard form then $$ \sigma (\D g^E) = X_\mu := \{\alpha \in \Path : \mu \hbox { is a prefix of
} \alpha \}.  $$

\Proof Assuming that $\D g^E$ is nonempty, let $\omega \in \D g^E$, and let $\alpha =\sigma (\omega )$.  So $$ g \_\in
{DescriptBernouliRelDoms} \omega =\omega _\alpha , $$ and point (i) follows from \cite {ReducedFormOmegaAlpha}.

In order to prove (ii), let $\omega \in \D g^E$, and let $\alpha =\sigma (\omega )$.  Since $g\in \omega =\omega _\alpha
$, we conclude from \cite {ReducedFormOmegaAlpha} that there is a prefix of $\alpha $, say $\mu '$, and a finite path
$\nu '$ with $d(\mu ')=d(\nu ')$, such that $$ g = \mu \nu \inv = \mu '{\nu '}\inv , $$ and $\ell (g) = |\mu '| + |\nu
'|$. By uniqueness of reduced forms, we conclude that $\mu =\mu '$, and $\nu =\nu '$.  In particular $\mu $ is a prefix
of $\alpha $, so $$ \sigma (\omega ) = \alpha \in X_\mu .  $$

This proves that $\sigma (\D g^E) \subseteq X_\mu $.  In order to prove the reverse inclusion, let $\alpha \in X_\mu $.
Then $\mu $ is a prefix of $\alpha $, so it is clear from the definition of $\omega _\alpha $ that $g\in \omega _\alpha
$, which in turn implies by \cite {DescriptBernouliRelDoms} that $\omega _\alpha \in \D g^E$.  Therefore $$ \alpha =
\sigma (\omega _\alpha ) \in \sigma ( \D g^E).  \omitDoubleDollar \endProof

Having understood the mirror images of the $\D g^E$ through the stem function, we will now describe the corresponding
partial action on the path space.

\state Proposition \label ThetaGSameAsTauMuNu Given $g=\mu \nu \inv $ in standard form, consider the mapping $$ \tau
_{\mu ,\nu }: \nu \gamma \in X_\nu \mapsto \mu \gamma \in X_\mu , $$ Then the diagram

\vskip -0.8cm \null \hfill \beginpicture \setcoordinatesystem units <0.0025truecm, -0.0025truecm> point at 0 0 \put
{$\Di g^E$} at 000 000 \put {$\D g^E$} at 1000 000 \arrow <0.11cm> [0.3,1.2] from 240 0 to 800 000 \put {$\th g^E$} at
530 -150

\put {$X_\nu $} at 000 700 \put {$X_\mu $} at 1000 700 \arrow <0.11cm> [0.3,1.2] from 240 700 to 800 700 \put {$\tau
_{\mu ,\nu }$} at 530 580

\arrow <0.11cm> [0.3,1.2] from 0 200 to 0 500 \put {$\sigma $} at -140 350 \arrow <0.11cm> [0.3,1.2] from 1000 200 to
1000 500 \put {$\sigma $} at 1140 350 \endpicture \hfill \null

\medskip \noindent commutes.

\Proof Given any $\omega $ in $\Di g^E$, let $$ \alpha =\sigma (\omega ) \and \beta =\sigma (g\omega ).  $$ Then $
\alpha \in X_\nu , $ and $ \beta \in X_\mu , $ so we may write $$ \alpha = \nu \gamma \and \beta = \mu \delta , $$ for
suitable paths $\gamma $ and $\delta $.

Given any finite prefix $\gamma '$ of $\gamma $, notice that $\nu \gamma '$ is a prefix of $\alpha =\nu \gamma $, whence
$\nu \gamma '\in \omega $, and $$ g \omega \ni g\nu \gamma ' = \mu \gamma ', $$ so we conclude that $\mu \gamma '$ is a
prefix of $\beta $.  Since $\gamma '$ may be taken to be equal to $\gamma $ when $\gamma $ is finite, or arbitrarily
large if $\gamma $ is infinite, we conclude that $\mu \gamma $ is a prefix of $\beta =\mu \delta $, whence $\gamma $ is
a prefix of $\delta $.

By repeating the above reasoning relative to $g\inv $, we similarly conclude that $\delta $ must be a prefix of $\gamma
$, which is to say that $\gamma =\delta $.  Therefore $$ \tau _{\mu ,\nu }(\alpha ) = \tau _{\mu ,\nu }(\nu \gamma ) =
\mu \gamma = \mu \delta = \beta , $$ concluding the proof.  \endProof

We will now describe a topology on $\FullPath $, relative to which the stem is a homeomorphism.  Before doing this, let
us observe that, since $\OE $ is a topological subspace of the product space $\{0,1\}^\F $, its topology is generated by
the inverse images of open subsets of $\{0,1\}$, under the canonical projections, $$ \pi _g: \OE \subseteq \{0,1\}^\F
\to \{0,1\}, $$ for $g\in \F $.  The topology of $\{0,1\}$, in turn, is generated by the open subsets $\{0\}$ and
$\{1\}$, so we may generate the topology of $\OE $ using the sets $$ \pi _g\inv (\{0\}) \and \pi _g\inv (\{1\}), $$ for
all $g$ in $\F $.  Under our identification of $\{0,1\}^\F $ with the set $\Part (\F )$ of all subsets of $\F $, recall
from \cite {VectorsAndSets} that the $\pi _g$ are given by $\pi _g(\omega )=\bool {g\in \omega }$.  We then have by
\cite {DescriptBernouliRelDoms} that $$ \pi _g\inv (\{1\})= \D g^E \and \pi _g\inv (\{0\})= \OE \setmenos \D g^E.  $$ So
we see that the topology of $\OE $ is generated by the $\D g^E$ and their complements.  However, by \cite
{PathCodingOfDg/i} we need only care for these when $g$ admits a standard form.

\state Proposition \label PathTopology Consider the topology on $\FullPath $ generated by all of the $X_\mu $, together
with their complements relative to $\FullPath $, where $\mu $ ranges in $\FinPath $.  Then the stem $$ \sigma : \OE \to
\FullPath , $$ is a homeomorphism, and consequently $\FullPath $ is a compact, Hausdorff, totally disconnected
topological space.

\Proof Follows immediately from \cite {PathCodingOfDg/ii} and the above comments about the topology of $\OE $. \endProof

\nostate \label SummaryTau Summarizing, let us give a detailed description of the restricted partial Bernoulli action on
$\OE $, once the latter is identified with the full path space of $E$.  \iaitem \aitem The acting group is the free
group $\F $ on the set $E^1$ of edges in our graph, \aitem The space $\FullPath $ consists of all finite and infinite
paths on $E$, plus a so called empty path $\varnothing $.  Thus $$ \FullPath = \Path \cup \{\varnothing \} = \FinPath
\cup \InfPath \cup \{\varnothing \}.  $$ \aitem For each $g$ in $\F $ admitting a standard form $g=\mu \nu \inv $, we
set $$ \FullPath _g=X_\mu = \{\alpha \in \Path : \mu \hbox { is a prefix of } \alpha \}.  $$ \aitem For each $g$ in $F$
not admitting a standard form, we let $\FullPath _g$ be the empty set.  \aitem For each $g$ in $\F $ admitting a
standard form $g=\mu \nu \inv $, we let $\tau _g$ be the partial homeomorphism of $\FullPath $ given by $$ \def \quad {}
\matrix {\tau _g = \tau _{\mu ,\nu }\ : \ \nu \gamma \ & \in & X_\nu & \longmapsto \ \mu \gamma \in &\ X_\mu . \cr &&
\vertequal && \vertequal \cr && \pilar {14pt} \FullPath _{g\inv } && \pilar {14pt} \FullPath _g } $$ \aitem For each $g$
in $\F $ not admitting a standard form, we let $\tau _g$ be the empty map from the empty set $\FullPath _{g\inv }$ to
the empty set $\FullPath _g$.

\bigskip \noindent We thus obtain a topological partial action $$ \tau = \big (\{\FullPath _g\}_{g\in \F },\ \{\tau
_g\}_{g\in \F }\big ), $$ of $\F $ on $\FullPath $, which is evidently topologically equivalent to the restriction of
the partial Bernoulli action to $\OE $.

\bigskip Since equivalent partial actions clearly give rise to isomorphic crossed products, we have the following
immediate consequence of \cite {TTEAsCrossProd}:

\state Theorem \label BestModelForTTe Given a graph $E$ with no sinks, consider the partial action $\tau $ described
above.  Then there exists a *-isomor\-phism $$ \varphi :\TTE \to C(\FullPath )\rt \tau \F , $$ such that $$ \varphi (\s
a ) = 1_{X_a }\delta _a , $$ for all $a \in E^1$.

Since the stem of $\omega _0$ coincides with $\varnothing $, we have that $\sigma $ establishes a covariant
homeomorphism from $\OE \setmenos \{\omega _0\}$ onto $\FullPath \setmenos \{\varnothing \} = \Path $, so we have, as in
\cite {OnlyTEAsCrossProd} that:

\state Theorem \label BestModelForTeOnPath In the context of \cite {BestModelForTTe}, the restriction of $\varphi $ to
$\TE $ is an isomorphism from $\TE $ onto $C_0(\Path )\rt \tau \F $.

Let us take a few moments to discuss the question of compactness of the various spaces appearing above.  The starting
point is of course the fact the $\OE $ is compact by \cite {ORIsCompactInvar}.  Being homeomorphic to $\OE $ by \cite
{PathTopology}, one also has that $\FullPath $ is compact.

\state Proposition \label Compactos Let $E$ be a graph without sinks.  Then the following are equivalent: \izitem \zitem
$E^0$ is finite, \zitem $\Path $ is compact, \zitem $\omega _0$ is an isolated point in $\OE $, \zitem $\varnothing $ is
an isolated point in $\FullPath $.

\Proof (i) $\Rightarrow $ (iii): For each vertex $v\in E^0$, choose an edge $a _v\in E^1$ such that $d(a _v)=v$, and
consider the compact-open subset of $\OE $ given by $$ U_v = \Di {a _v}^E = \{\omega \in \OE : a _v\inv \in \omega \}.
$$ By \cite {LocalConditionsInterpreted/ii} applied to $g=1$, we see that the definition of $U_v$ does not depend on the
choice of $a _v$.  Using \cite {LocalConditions} we than conclude that $$ \OE \setmenos \{\omega _0\} = \bigcup _{v\in
E^0}U_v.  $$

Assuming that $E^0$ is finite, we see that the above union of sets is closed, whence its complement, namely $\{\omega
_0\}$, is open, proving $\omega _0$ to be an isolated point.

\medskip \noindent (iii) $\Leftrightarrow $ (iv): Follows immediately from the fact that the stem is a homeomorphism
from $\OE $ to $\FullPath $, sending $\omega _0$ to $\varnothing $.

\medskip \noindent (iv) $\Rightarrow $ (ii): If $\varnothing $ is isolated, then $$ \Path = \FullPath \setmenos
\{\varnothing \} $$ is closed in the compact space $\FullPath $, hence $\Path $ is compact.

\medskip \noindent (ii) $\Rightarrow $ (i): Assuming (ii), the open cover $\{X_v\}_{v\in E^0}$ of $\Path $ admits a
finite subcover, say $$ \Path = X_{v_1}\cup X_{v_2}\cup \ldots \cup X_{v_n}.  $$ Given any vertex $v$, viewed as a path
of length zero, hence an element of $\Path $, there exists some $k$ such that $v\in X_{v_k}$, so $v_k$ is a prefix of
$v$, which is to say that $v_k=v$.  Thus $E^0 = \{v_1,v_2,\ldots ,v_n\}$ is a finite set.  \endProof

\definition We shall denote by $\Bun ^E$ the semi-direct product bundle for the partial action $\tau $ given in \cite
{SummaryTau}.

As a consequence of \cite {BestModelForTTe} and \cite {CPIsCrossSectAlg}, we have that $\TTE $ is isomorphic to the full
cross sectional C*-algebra of $\Bun ^E$.

Recall from \cite {PRepInTE} that there exists a semi-saturated partial representation $\Pr $ of $\F $ in $\TTE $, such
that $ \pr a = \s a , $ for all $a \in E^1$.  Composing $\Pr $ with the isomorphism given in \cite {BestModelForTTe} we
obtain a semi-saturated partial representation $$ v:\F \to C(\FullPath )\rt \tau \F $$ such that $$ v_a = 1_{X_a }\delta
_a \for a \in E^1.  \equationmark DescrPartRep $$

\state Proposition The Fell bundle $\Bun ^v$ associated to the above partial representation $v$, in line with \cite
{BunFromPrep}, is isomorphic to $\Bun ^E$.

\Proof Writing $\Bun ^E = \{B^E_g\}_{g\in \F }$, we will identify each $B^E_g$ with the corresponding grading subspace
$$ C_0( \FullPath _g)\delta _g \subseteq C(\FullPath )\rt \tau \F .  $$ With this identification, we will in fact prove
that $\Bun ^v$ is \"{equal to} $\Bun ^E$.

By \cite {DescrPartRep} we clearly have that $$ v_g \in B^E_g \for g\in E^1, $$ and, by taking adjoints, this implies
that the same holds for all $g$ in $(E^1)\inv $.

Given a general element $g$ in $\F $, write $g = x_1\cdots x_n$ in reduced form, that is, each $x_i$ lies in $E^1\cup
(E^1)\inv $, and $x_i\inv \neq x_{i+1}$.  Then, using the fact that $v$ is semi-saturated, we deduce that $$ v_g =
v_{x_1}\cdots v_{x_n} \in B^E_{x_1}\cdots B^E_{x_n} \subseteq B^E_{x_1\cdots x_n} = B^E_g.  \equationmark vgInBg $$

If we now suppose that $g=h_1\cdots h_n$, no longer necessarily in reduced form, we have that $$ v_{h_1}\cdots v_{h_n}
\in B^E_{h_1}\cdots B^E_ {h_n} \subseteq B^E_{h_1\cdots h_n} = B^E_g.  $$ By definition of $\Bun ^v$ (see \cite
{BunFromPrep}) we then see that $B^v_g\subseteq B^E_g$, for every $g$ in $\F $, and the proof will be concluded once we
prove that in fact $B^v_g$ coincides with $B^E_g$, which is to say that $B^E_g$ is the closed linear span of the set of
elements of the form $$ v_{h_1}\cdots v_{h_n}, $$ with $g=h_1\cdots h_n$.

By definition $\TE $ is generated, as a C*-algebra, by the $\s a $, for $a \in E^1$, whence $\TTE $ is generated by
$\{\s a :a \in E^1\}\cup \{1\}$.  Consequently $C(\FullPath )\rt \tau \F $ is generated by $\{1_{X_a }\delta _a :a \in
E^1\}\cup \{1\}$, and evidently also by the range of $v$.  Thus, given $y$ in any $B_g^E$, we may write $$ y = \lim
_{n\to \infty } y_n, $$ where each $y_n$ is a linear combination of elements of the form $v_{h_1}\cdots v_{h_n}$, with
$h_i$ in $\F $.

Let $P_g$ be the composition $$ P_g : C^*\big (\Bun ^E\big ) \>{\ \Lambda \ } \CstarRed \big (\Bun ^E\big ) \>{E_g}
B_g^E, $$ where $E_g$ is the Fourier coefficient operator given by \cite {Fourier}.  We then have that $$ y = P_g(y) =
\lim _{n\to \infty } P_g(y_n).  $$

By \cite {vgInBg} it is easy to see that $$ P_g(v_{h_1}\cdots v_{h_n}) = \left \{\matrix {v_{h_1}\cdots v_{h_n}, & \hbox
{ if } h_1\cdots h_n =g, \cr \pilar {16pt} 0, & \hbox { otherwise.}\hfill }\right .  $$ So, upon replacing each $y_n$ by
$P_g(y_n)$, we may assume that each $y_n$ is a linear combination of elements of the form $v_{h_1}\cdots v_{h_n}$, with
$h_1\cdots h_n=g$.  This implies that $y_n$ is in $B^v_g$, and hence also $y$ is in $B^v_g$.  This shows that
$B^v_g=B^E_g$, concluding the proof.  \endProof

A generalization of the above result could be proved for any C*-algebra of the form considered in \cite {CstarRelAsCP},
although we shall not pursue this here.

\medskip Let us now discuss the question of amenability for $\Bun ^E$.

\state Theorem \label AmenaPath Given a graph $E$ without sinks, one has that \izitem \zitem $\Bun ^E$ satisfies the
approximation property and hence is amenable, \zitem $C(\FullPath )\rt \tau \F $ is naturally isomorphic to $C(\FullPath
)\redrt \tau \F $, \zitem $C(\Path )\rt \tau \F $ is naturally isomorphic to $C(\Path )\redrt \tau \F $, \zitem $\TTE $
and $\TE $ are nuclear C*-algebras.

\Proof We have already seen that the partial representation $\Pr $ of \cite {PRepInTE} is semi-saturated.  By \cite
{CKOrthogonal} we have that $\Pr $ is also orthogonal, in the sense of \cite {ApproxBundle/i}.  Since $v$ is the
composition of $\Pr $ with an isomorphism, it is clear that $v$ is also semi-saturated and orthogonal.  It therefore
follows from \cite {ApproxBundle/i} that $\Bun ^v$, and hence also $\Bun ^E$, satisfy the approximation property.

We then conclude from \cite {Main} that $\Bun ^E$ is amenable, which is to say that the regular representation $$
\Lambda :C^*(\Bun ^E)\to \Cr {\Bun ^E} $$ is an isomorphism.  So, (ii) follows from the characterization of the full
crossed product as a full cross sectional C*-algebra (see \cite {CPIsCrossSectAlg}), and the definition of the reduced
crossed product as the reduced cross sectional C*-algebra (see \cite {DefineRedCP}).

As already noticed, $\Path $ is an invariant open subset of $\FullPath $.  Thus, if $\Abun $ denotes the semi-direct
product bundle for the restriction of $\tau $ to $\Path $, we have that $\Abun $ is naturally isomorphic to a
\subFellBun dle of $\Bun ^E$.  Since the unit fiber algebra of $\Abun $, namely $C_0(\Path )$, is an ideal in
$C(\FullPath )$, the unit fiber algebra of $\Bun $, we deduce from \cite {SubBunHasAP} that $\Abun $ also satisfies the
approximation property.  Thus (iii) follows as above.

Finally, (iv) follows from \cite {AproxImplyNuc}, in view of the fact that commutative C*-algebras are nuclear.
\endProof

\nrem The partial crossed product description of $\TE $ given in \cite {BestModelForTeOnPath} follows the general method
adopted in \ref {ExelLaca/1999} to give a similar description of Cuntz-Krieger algebras for infinite matrices.

\chapter Graph C*-algebras

Having obtained a description of $\TE $ in terms of a concrete partial dynamical system in \cite {BestModelForTeOnPath},
we will now proceed to giving a similar model for the graph C*-algebra $C^*(E)$.  There are two methods for doing this,
the most obvious and equally effective one being to repeat the above procedure, applying \cite {CstarRelAsCP} to the set
of relations defining $C^*(E)$, and then reinterpreting the spectrum of these relations as a path space.

For a change we will instead see $C^*(E)$ as a quotient of $\TE $, and then we will apply the results of chapter \cite
{FunctPASec}.

\medskip

\fix As before, let us fix a graph $E=(E^0,E^1,r,d)$, assumed to have no sinks.

Recall that $C^*(E)$ is defined in much the same way as $\TE $, the only difference being that relations \cite {CKSum}
are required to hold in $C^*(E)$ but not in $\TE $.  These are the relations $$ \p v = \soma {r(a )=v}\s a \s a ^*, $$
for each regular vertex $v\in E^0$ (see Definition \cite {DefineSinkSourceRegular/iii}).

Since we are working with graphs without sinks, we may phrase these relations in terms of edges only, as follows: given
any vertex $v$ as above, choose an edge $b$ such that $d(b) = v$, in which case $\p v = \s b ^*\s b $, by \cite
{CKOrthogonal}.  Therefore \cite {CKSum} may be rewritten as $$ \s b ^*\s b = \soma {r(a )=d(b )}\s a \s a ^*,
\equationmark NewCKSum $$ for every $b$ in $E^1$ such that $d(b )$ is a regular vertex.  These relations actually carry
a greater similarity to the original relations studied by Cuntz and Krieger \ref {CuntzKrieger/1980}.

Since relations \cite {DefineGraphAlg/\CKOrthogonalLoc --\CKLessThanLoc } are already satisfied in $\TE $, we clearly
obtain the following:

\state Proposition \label CStarQuotientOfTE Let $L$ be the closed two-sided ideal in $\TE $ generated by the elements $$
\s b ^*\s b - \soma {r(a )=d(b )}\s a \s a ^*, $$ for every $b$ in $E^1$ such that $d(b )$ is a regular vertex.  Then
there exists a *-isomorphism $$ \chi : C^*(E) \to \TE /L, $$ sending each canonical generating partial isometry of
$C^*(E)$, which we henceforth denote by $\tilde \s a $, to $\s a + L$.

Observe that, as a consequence of \cite {CKOrthogonal}, and of our assumption that $E$ has no sinks, the $\tilde \s a $
are enough to generate $C^*(E)$, whence the isomorphism referred to above is uniquely determined by the fact that $$
\chi (\tilde \s a ) = \s a + L.  $$

We will now describe the elements of $ \CoefAlg \rt \tau \F $ corresponding to the above generators of the ideal $L$ via
the isomorphism given in \cite {BestModelForTeOnPath}.  As we will see, these lie inside the canonical image of $
\CoefAlg $ in the crossed product, so we will be able to characterize the quotient algebra as a partial crossed product
using \cite {QuotientOfCP}

Given any edge $b $ in $E^1$ such that $d(b )$ is regular, consider the complex valued function on $\FullPath $ defined
by $$ f_b (\alpha ) = \bool {d(b ) \hbox { is a prefix of } \alpha } - \kern -5pt \soma {r(a )=d(b )} \bool {a \hbox {
is a prefix of } \alpha } \for \alpha \in \FullPath , $$ where brackets correspond to Boolean value.  Clearly $f_b $ is
continuous.  Since $\varnothing $ has no prefixes, we see that $f_b $ vanishes on $\varnothing $, whence $f_b \in
\CoefAlg $.

Recalling from \cite {PathCodingOfDg/ii} that $X_\mu $ is the set of all finite and infinite paths admitting the given
finite path $\mu $ as a prefix, notice that $$ \bool {\mu \hbox { is a prefix of } \alpha } = 1_{X_\mu }(\alpha ).  $$
Consequently the function $f_b $ defined above may be alternatively described as $$ f_b = 1_{X_{d(b )}} - \soma {r(a
)=d(b )} 1_{X_a }.  $$

For each edge $a $ in $E^1$, notice that the isomorphism $\varphi $ of \cite {BestModelForTTe} satisfies $$ \varphi (\s
a \s a ^*) = (1_{X_a }\delta _a ) (1_{X_a }\delta _a )^* \={LotsAFormulas} 1_{X_a }\delta _1, $$ and similarly, $$
\varphi (\s b ^*\s b ) = (1_{X_b }\delta _b)^* (1_{X_b }\delta _b) = \tau _{b\inv }(1_{X_b }) \delta _1.  $$

The standard form of $b \inv $ is clearly $d(b )b \inv $, from where we see that $ \tau _{b\inv }(1_{X_b }) = 1_{X_{d(b
)}}, $ so we deduce from the above that $$ \varphi (\s b ^*\s b ) = 1_{X_{d(b )}}\delta _1.  $$ Interpreting $f_b $
within $\CoefAlg \rt \tau \F $ via the map $\iota $ introduced in \cite {DefineStdInclusion} then produces $$ \iota (f_b
) = f_b \delta _1 = 1_{X_{d(b )}}\delta _1 - \soma {r(a )=d(b )} 1_{X_a }\delta _1 \$= \varphi (\s b ^*\s b ) - \soma
{r(a )=d(b )} \varphi (\s a \s a ^*) = \varphi \Big (\s b ^*\s b - \soma {r(a )=d(b )} \s a \s a ^*\Big ), $$ which
should be compared with \cite {NewCKSum}.  As an immediate consequence of \cite {BestModelForTeOnPath} and \cite
{CStarQuotientOfTE}, we therefore conclude that:

\state Proposition \label CStarEAsQuotient Let $W$ be the subset of $\CoefAlg $ formed by the functions $$ f_b =
1_{X_{d(b )}} - \soma {r(a )=d(b )} 1_{X_a }, $$ for each $b \in E^1$ such that $d(b )$ is regular.  Also let $K$ be the
closed two-sided ideal of $\CoefAlg \rt \tau \F $ generated by $\iota (W)$.  Then the isomorphism $\varphi $ given in
\cite {BestModelForTeOnPath} sends the ideal $L$ described in \cite {CStarQuotientOfTE} onto $K$, and consequently there
exists a *-isomorphism $$ \psi :C^*(E) \mathrel {\longrightarrow } {\CoefAlg \rt \tau \F \over K}, $$ such that $$ \psi
(\tilde \s a ) = 1_{X_a }\delta _a + K, $$ where we are again denoting by $\tilde \s a $ the standard generating partial
isometries of $C^*(E)$.

\Proof It is enough to take $\psi $ to be the composition $$ C^*(E) \quad \_{{\hbox to 1.2cm{\rightarrowfill
}}}{CStarQuotientOfTE} \quad {\TE \over L} \quad \_{{\hbox to 1.2cm{\rightarrowfill }}}{BestModelForTeOnPath} \quad
{\CoefAlg \rt \tau \F \over K}.  \omitDoubleDollar \endProof

We are now precisely under the hypotheses of \cite {QuotientOfCP}, but before invoking it we will give a concrete
description of the ideal $J$ mentioned there.

\state Proposition \label DescribeSpecAllRel Let\/ $\CstarPath $ be the subset of\/ $\Path $ consisting of all paths
$\alpha $ which satisfy any one of the following conditions: \izitem \zitem $\alpha $ is infinite, \zitem $\alpha $ is
finite and $r\inv \big (d(\alpha )\big )$ is empty, \zitem $\alpha $ is finite and $r\inv \big (d(\alpha )\big )$ is
infinite.  \medskip \noindent Then $\CstarPath $ is closed and $\tau $-invariant.  Moreover, denoting by $U$ the
complement of\/ $\CstarPath $ in $\Path $, we have that $C_0(U)$ is the smallest $\tau $-invariant ideal of $\CoefAlg $
containing the set $W$ referred to in \cite {CStarEAsQuotient}.

\Proof We claim that every $\alpha $ in $U$ is an isolated point of $\Path $.  In fact, if $\alpha $ is not in
$\CstarPath $, then $\alpha $ is necessarily a finite path and $r\inv \big (d(\alpha )\big )$ is a nonempty finite set,
so $d(\alpha )$ is a regular vertex and we may write $$ r\inv \big (d(\alpha )\big ) = \{a _1,\ldots ,a _n\}.  $$

Observe that the set $$ X_\alpha \cap \bigcap _{i=1}^n \Path {\setminus } X_{\alpha a _i} $$ is open by \cite
{PathTopology}.  It consists of all paths admitting $\alpha $ as a prefix, but not admitting as a prefix any one of the
paths $\alpha a _i$.  The unique such path is evidently $\alpha $, so the above set coincides with $\{\alpha \}$, thus
proving that $\alpha $ is an isolated point.  Consequently $U$ is open, whence $\CstarPath $ is closed relative to
$\Path $.

Although this is not relevant to us at the moment, notice that any path $\alpha $ satisfying (ii) above is also an
isolated point, since $ X_\alpha = \{\alpha \}.  $

Noticing that the conditions (i--iii) above are related to the ``right end'' of $\alpha $, while $\tau $ affects its
``left end'', one may easily show that $\CstarPath $ is $\tau $-invariant.

We next claim that every $f_b $ in $W$ vanishes on $\CstarPath $.  In fact, given any $\alpha $ in $\CstarPath $, and
given any edge $b $ such that $d(b )$ is regular, we have $$ f_b (\alpha ) = \bool {d(b ) \hbox { is a prefix of }
\alpha } - \soma {r(a )=d(b )} \bool {a \hbox { is a prefix of } \alpha }.  $$

If $d(b )$ is not a prefix of $\alpha $, meaning that $r(\alpha )\neq d(b )$, then evidently $r(\alpha )\neq r(a )$, for
all edges $a $ considered in the above sum.  This implies that none of these $a $'s are prefixes of $\alpha $, given
that \cite {DefinePrefix/i} fails.  Therefore all terms making up $f_b (\alpha )$ vanish, whence $f_b (\alpha )$ itself
vanishes.

The remaining case to be treated is when $d(b )$ is a prefix of $\alpha $, that is, when $r(\alpha )=d(b )$.  Should
$\alpha $ be a path of length zero, necessarily consisting of the vertex $d(b )$, then $d(\alpha )=d(b )$, so $d(\alpha
)$ is regular by assumption, and then $\alpha $ will not satisfy any one of conditions (i--iii) above, contradicting the
fact that $\alpha $ was taken in $\CstarPath $.  This said we see that $\alpha $ must not have length zero, so we may
write $$ \alpha = \alpha _1\alpha _2\ldots $$ Given that $r(\alpha _1) = r(\alpha ) = d(b )$, we see that $\alpha _1$ is
one of the edges $a $ considered in the sum appearing in the definition of $f_b $ above, and for this edge, the value of
``$a $ is a prefix of $\alpha $'' is evidently 1.  This shows that $f_b (\alpha )=0$, showing our claim that $f_b $
vanishes on $\CstarPath $.

We will now prove that $\CstarPath $ is the biggest invariant subset of $\Path $ where the $f_b $ vanish.  For this,
suppose that $\Lambda $ is an invariant set properly containing $\CstarPath $.  Choosing any $\alpha $ in $\Lambda
\setminus \CstarPath $, we see that $\alpha $ is a finite path, so we may let $v=d(\alpha )$, noticing that $v$ is
regular.  The standard form of $ g := \alpha \inv $ is given by $g = d(\alpha )\alpha \inv $ so, by invariance of
$\Lambda $, we have $$ \Lambda \ni \tau _g(\alpha ) = d(\alpha ) = v.  $$ Choosing any edge $b $ with $d(b)=v$, we have
that $f_b $ lies in $W$, and clearly $f_b (v) = 1$, because $v$ is a prefix of itself, while $v$ can have no prefix of
positive length.  This implies that $f_b $ does not vanish on $\Lambda $, proving that indeed $\CstarPath $ is the
biggest invariant subset of $\Path $ where the $f_b $ vanish.  This also shows that $U$ is the smallest invariant subset
outside of which the $f_b $ vanish, thus concluding the proof.  \endProof

The above classification of paths in types (i), (ii) and (iii) has important consequences in what follows and in
particular the first two kinds play a special role, justifying the introduction of the following terminology:

\definition A path $\alpha $ in $\FullPath $ is said to be \subjex {maximal}{maximal path} if it satisfies \cite
{DescribeSpecAllRel/i} or \cite {DescribeSpecAllRel/ii}.

This terminology is justified because a maximal path $\alpha $ cannot be enlarged, either because $d(\alpha )$ is a
source or because $\alpha $ is already infinite.

\medskip Let us take a few moments to study the topology of $\CstarPath $.  Recall from \cite {PathTopology} that the
topology of $\FullPath $ is generated by all the $X_\mu $ plus their complements.  Therefore, the sets of the form $$
\bigcap _{i=1}^n X_{\mu _i} \cap \med \bigcap _{j=1}^m \FullPath \setmenos X_{\nu _j} $$ where $\mu _1,\ldots ,\mu _n$
and $\nu _1,\ldots ,\nu _m$ are finite paths, constitute a basis for the topology of $\FullPath $.  Consequently the
intersections of these with $\CstarPath $ form a basis of open sets for the latter.  However, maximal paths have a
neighborhood basis of a simpler nature, as we will now see.

\state Proposition \label SimpleNbd Given a graph $E$ without sinks, and given a maximal path $\alpha $ in $\CstarPath
$, the collection of sets $$ \{X_\mu \cap \CstarPath : \mu \in \FinPath ,\ \mu \hbox { is a prefix of }\alpha \} $$ is a
neighborhood basis for $\alpha $ in $\CstarPath $.

\Proof We first assume that $\alpha $ satisfies \cite {DescribeSpecAllRel/ii}, that is, $\alpha $ cannot be extended any
further due to the fact that $d(\alpha )$ is a source.  Then the only path admitting $\alpha $ as a prefix is $\alpha $
itself, meaning that $$ X_\alpha =\{\alpha \}, \equationmark XalphaSingle $$ so the result follows trivially.

Assume now that $\alpha $ satisfies \cite {DescribeSpecAllRel/i}, that is, $\alpha $ is infinite.  Given any
neighborhood $U$ of $\alpha $, there are finite paths $\mu _1,\ldots ,\mu _n$ and $\nu _1,\ldots ,\nu _m$ such that $$
\alpha \in \med \bigcap _{i=1}^n X_{\mu _i} \cap \med \bigcap _{j=1}^m \FullPath \setmenos X_{\nu _j}\subseteq U.  $$
This implies that every $\mu _i$ is a prefix of $\alpha $, while no $\nu _j$ has this property.

Let $k$ be any integer larger that the length of every $\mu _i$ and every $\nu _j$, and let $\mu $ be the finite path
formed by the first $k$ edges of $\alpha $.  It is therefore evident that each $\mu _i$ is a prefix of $\mu $, so one
has that $$ \alpha \in X_\mu \subseteq \med \bigcap _{i=1}^n X_{\mu _i}.  \equationmark BetaNeighBase $$

Observe that $X_\mu $ is disjoint from $X_{\nu _j}$, for every $j$ because, otherwise there is a path $\gamma \in X_\mu
\cap X_{\nu _j}$, and hence both $\mu $ and $\nu _j$ are prefixes of $\gamma $.  But since $|\mu |\geq |\nu _j|$, we see
that $\nu _j$ is a prefix of $\mu $, which in turn is a prefix of $\alpha $.  It follows that $\nu _j$ is a prefix of
$\alpha $, whence $\alpha \in X_{\nu _j}$, a contradiction.  So $X_\mu \subseteq \FullPath \setmenos X_{\nu _j}$ and,
building on top of \cite {BetaNeighBase}, we conclude that $$ \alpha \in X_\mu \subseteq \med \bigcap _{i=1}^n X_{\mu
_i} \cap \med \bigcap _{j=1}^m \FullPath \setmenos X_{\nu _j}\subseteq U.  $$ This concludes the proof.  \endProof

Another useful fact about the topology of $\CstarPath $ is as follows:

\state Proposition \label DensityIandII Given a graph $E$ without sinks, the subset of $\CstarPath $ formed by the
maximal paths is dense in $\CstarPath $.

\Proof Before starting the proof, observe that, by \cite {XalphaSingle}, the paths satisfying \cite
{DescribeSpecAllRel/ii} are isolated points, so they cannot be left out of any dense subset!

In order to prove the statement, it is enough to prove that any path $\alpha $ satisfying \cite {DescribeSpecAllRel/iii}
is an accumulation point of maximal paths.

Given an arbitrary neighborhood $U$ of $\alpha $, pick finite paths $\mu _1,\ldots ,\mu _n$ and $\nu _1,\ldots ,\nu _m$
such that $$ \alpha \in \med \bigcap _{i=1}^n X_{\mu _i} \cap \med \bigcap _{j=1}^m \FullPath \setmenos X_{\nu
_j}\subseteq U $$ (notice that we are unfortunately unable to use the simplification provided by \cite {SimpleNbd}
here).

Thus every $\mu _i$ is a prefix of $\alpha $, hence we may assume without loss of generality that the $\mu _i$ all
coincide with $\alpha $, so in fact $$ \alpha \in X_{\alpha } \cap \med \bigcap _{j=1}^m \FullPath \setmenos X_{\nu
_j}\subseteq U.  \equationmark BadNeighOfAlpha $$

By assumption we have that $r\inv \big (d(\alpha )\big )$ is infinite, so me may pick an edge $a $ such that $\alpha a $
is a path which is not a prefix of any $\nu _j$.

We next extend $\alpha a $ as far as possible, obtaining a path of the form $$ \beta = \alpha a \gamma , $$ which is
either infinite or cannot be extended any further, and hence is maximal.  We then claim that $$ \beta \in X_{\alpha }
\cap \med \bigcap _{j=1}^m \FullPath \setmenos X_{\nu _j}.  \equationmark ClaimWhereBetaIs $$ Since it is obvious that
$\alpha $ is a prefix of $\beta $, it suffices to check that $\beta \not \in X_{\nu _j}$, for all $j$.  Arguing by
contradiction, suppose that some $\nu _j$ is a prefix of $\beta $.  If $|\nu _j|\leq |\alpha |$, then $\nu _j$ is a
prefix of $\alpha $, contradicting \cite {BadNeighOfAlpha}.  Thus $|\nu _j|>|\alpha |$, and then necessarily $\alpha a $
is a prefix of $\nu _j$, which is again a contradiction, by the choice of $a $.  This verifies \cite {ClaimWhereBetaIs},
so in particular $\beta \in U$, concluding the proof.  \endProof

Let us now give a description of $C^*(E)$ as a partial crossed product.

\state Theorem \label CStarEAsCrossProd Let $E$ be a graph with no sinks, and let $\CstarPath $ be the $\tau $-invariant
subset of $\Path $ described in \cite {DescribeSpecAllRel}, equipped with the restricted partial action.  Then there is
a *-isomorphism $$ \rho :C^*(E) \to C_0(\CstarPath ) \rt \tau \F , $$ such that $$ \rho (\tilde \s a ) = 1_a \delta _a
\for a \in E^1, $$ where $1_a $ is the characteristic function of $X_a \cap \CstarPath $.

\Proof Letting $C_0(U)$ be the ideal referred to in \cite {DescribeSpecAllRel}, we have $$ C^*(E) \ \_\simeq
{CStarEAsQuotient} \ {\CoefAlg \rt \tau \F \over K} \ \_\simeq {QuotientOfCP} \ \left ({\CoefAlg \over C_0(U)}\right )
\rt \tau \F \ \simeq \ C_0(\CstarPath ) \rt \tau \F .  $$ That the isomorphism resulting from the composition of the
above isomorphisms does satisfy the last assertion in the statement is of easy verification and is left to the reader.
\endProof

An important special case of great importance is when, besides having no sinks, every vertex in $E^0$ is regular,
according to \cite {DefineSinkSourceRegular/iii} (so $E$ cannot have any sources either).  Then there are no paths in
$\CstarPath $ satisfying \cite {DescribeSpecAllRel/ii} or \cite {DescribeSpecAllRel/iii}, so we get the following
immediate consequence of \cite {CStarEAsCrossProd}:

\state Proposition If $E$ is a graph with no sinks and such that every vertex $v\in E^0$ is regular, then $\CstarPath
=\InfPath $, and consequently there is a *-isomorphism $$ \rho :C^*(E) \to C_0(\InfPath ) \rt \tau \F , $$ such that $$
\rho (\tilde \s a ) = 1_a \delta _a \for a \in E^1, $$ where $1_a $ is the characteristic function of $X_a \cap \InfPath
$.

Observe that if $E^0$ is finite then $\CstarPath $ is compact, because the latter is closed in $\Path $ by \cite
{DescribeSpecAllRel}, and $\Path $ is compact by \cite {Compactos}.  If, in addition, we are under the hypotheses of the
last result, then clearly $\InfPath $ is also compact.

\medskip Having described $C^*(E)$ as a partial crossed product, we now have a wide range of tools to study its
structure.  We begin with amenability.

\state Theorem \label GraphAmena Let $E$ be a graph without sinks.  Then: \iaitem \aitem the semi-direct product bundle
for the partial action describing $C^*(E)$ in \cite {CStarEAsCrossProd} satisfies the approximation property and hence
is amenable, \aitem there are natural isomorphisms $$ C^*(E) \simeq C_0(\CstarPath ) \rt \tau \F \simeq C_0(\CstarPath )
\redrt \tau \F , $$ \aitem $C^*(E)$ is nuclear.

\Proof The Fell bundle mentioned in (a) is clearly a quotient of the Fell bundle discussed in \cite {AmenaPath}, so the
conclusions in (a) follow from \cite {QuoAmena/ii} and \cite {Main}.  Point (b) follows from \cite {CStarEAsCrossProd}
and (a), while (c) is a direct consequence of (b) and \cite {AproxImplyNuc}).  \endProof

We will next study fixed points for the partial action $\tau $ of $\F $ on $\FullPath $.  For this we should recall that
by \cite {PathCodingOfDg/i}, unless a given element $g$ in $\F $ has a standard form, $\tau _g$ is the empty map and
hence $g$ cannot possibly have any fixed points.  For this reason elements not admitting a standard form are left out of
the next result.

\state Proposition \label GraFixPoint Let $g\in \F $, with $g\neq 1$, and suppose that $g$ admits a standard form $g=\mu
\nu \inv $.  \izitem \zitem Then $g$ has at most one fixed point in $\FullPath $.  \zitem If $g$ has a fixed point, then
$|\mu |\neq |\nu |$.  \zitem If $g$ has a fixed point and $|\mu |>|\nu |$, then there exists a cycle\fn {A finite path
$\gamma \in \FinPath $ is said to be a \subjex {cycle}{cycle in a graph} if $|\gamma |>1$, and $r(\gamma ) = d(\gamma
)$.  Notice that in this case we may concatenate $\gamma $ with itself as many times as we wish, forming a finite of
infinite periodic path $\gamma \gamma \gamma \ldots $.}  $\gamma $ such that $\mu =\nu \gamma $ (whence $g = \nu \gamma
\nu \inv $).  In addition, the infinite path $$ \alpha = \nu \gamma \gamma \gamma \ldots $$ is the unique fixed point
for $g$ in $\FullPath $.

\Proof Assume that $\alpha $ is a fixed point for $g$.  Then $$ \alpha \in \FullPath _{g\inv }\cap \FullPath _g = X_\nu
\cap X_\mu , $$ so we may write $$ \alpha = \nu \varepsilon = \mu \zeta , $$ for suitable (finite or infinite) paths
$\varepsilon $ and $\zeta $, and moreover $$ \alpha = \tau _g(\alpha ) = \tau _g(\nu \varepsilon ) = \mu \varepsilon .
$$

We then conclude that $\nu \varepsilon = \mu \varepsilon $, so either $\nu $ is a prefix of $\mu $, or vice versa.  If
$|\mu | = |\nu |$, one would then necessarily have $\mu =\nu $, whence $g=1$, contradicting the hypothesis.  This proves
(ii).

Speaking of (iii), let us add to the above assumptions that $|\mu |>|\nu |$, so we may write $\mu = \nu \gamma $, for
some finite path $\gamma $ such that $|\gamma |>0$.  We then have $$ \nu \varepsilon = \mu \varepsilon = \nu \gamma
\varepsilon , $$ so $\varepsilon = \gamma \varepsilon $, whence $|\varepsilon | = |\gamma | + |\varepsilon |$, from
where we deduce that $|\varepsilon |=\infty $.  Moreover notice that this implies that $$ \varepsilon =\gamma \gamma
\gamma \ldots , $$ so $\gamma $ is necessarily a cycle.  Consequently $ \alpha = \nu \varepsilon = \nu \gamma \gamma
\gamma \ldots , $ proving (iii), and hence also (i).  \endProof

\bigskip The reader must have noticed the omission of the case $|\mu |<|\nu |$, above.  However, in this situation
$g\inv $ has precisely the same fixed points as $g$, and the condition expressed in \cite {GraFixPoint/iii} evidently
holds for $g\inv $.

\bigskip Having understood fixed points, we will now study topological freeness.  We must therefore analyze the interior
of the fixed point sets $F_g$ which, as seen above, have each at most one point.  However, a singleton can only have a
nonempty interior if its unique point is an isolated point.  We must therefore discuss conditions under which a fixed
point is isolated, and this hinges on the following important graph-theoretical concept:

\vskip 2.5cm \eject

\definition Let $\gamma =\gamma _1\ldots \gamma _n$ be a cycle in $E$.  We say that $\gamma $ has an \subjex
{entry}{entry to a cycle}, if the range of some $\gamma _i$ is the range of an edge other that $\gamma _i$.

\null \hfill \beginpicture \setcoordinatesystem units <0.010truecm, 0.010truecm> point at 1000 1000 \put {$\bullet $} at
99 173 \arrow <0.15cm> [0.4,1] from 200 0 to 144 95 \plot 144 95 99 173 / \put {$\gamma _6$} at 178 115 \put {$\bullet
$} at -100 173 \arrow <0.15cm> [0.4,1] from 99 173 to -10 173 \plot -10 173 -100 173 / \put {$\gamma _5$} at -10 212
\put {$\bullet $} at -200 0 \arrow <0.15cm> [0.4,1] from -100 173 to -155 77 \plot -155 77 -200 0 / \put {$\gamma _4$}
at -189 97 \put {$\bullet $} at -99 -173 \arrow <0.15cm> [0.4,1] from -200 0 to -144 -95 \plot -144 -95 -99 -173 / \put
{$\gamma _3$} at -178 -115 \put {$\bullet $} at 100 -173 \arrow <0.15cm> [0.4,1] from -99 -173 to 10 -173 \plot 10 -173
100 -173 / \put {$\gamma _2$} at 10 -212 \put {$\bullet $} at 200 0 \arrow <0.15cm> [0.4,1] from 100 -173 to 155 -77
\plot 155 -77 200 0 / \put {$\gamma _1$} at 189 -97 \setdashes <1pt> \put {$\bullet $} at 400 0 \arrow <0.15cm> [0.4,1]
from 200 0 to 310 0 \plot 310 0 400 0 / \put {$\bullet $} at -200 346 \arrow <0.15cm> [0.4,1] from -100 173 to -155 268
\plot -155 268 -200 346 / \put {$\bullet $} at -200 0 \arrow <0.15cm> [0.4,1] from -400 0 to -290 0 \plot -290 0 -200 0
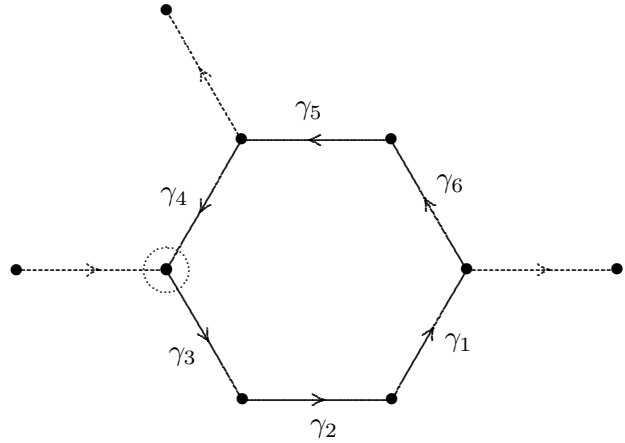
/ \put {$\bullet $} at -400 0 \setdots <1.5pt> \circulararc 360 degrees from -200 30 center at -200 0 \put {\eightrm An
entry to a cycle.} at 0 -300 \endpicture \hfill \null

\bigskip

\state Proposition \label IsolatedPoints Let $\alpha $ be an infinite path of the form $$ \alpha = \nu \gamma \gamma
\gamma \ldots , $$ so that necessarily $d(\nu ) = r(\gamma )$, and $\gamma $ is a cycle.  Then \izitem \zitem $\alpha $
(or any other infinite path) is never an isolated point in $\FullPath $, \zitem $\alpha $ is an isolated point in
$\CstarPath $ if and only if $\gamma $ has no entry.

\Proof An infinite path $\alpha =\alpha _1\alpha _2\alpha _3\ldots $ is always the limit of the sequence of finite paths
obtained by truncating $\alpha $, that is $$ \alpha = \lim _{n\to \infty }\alpha _1\ldots \alpha _n, $$ so $\alpha $ is
an accumulation point, hence not isolated.

Focusing now on (ii), suppose that $\alpha $ is an isolated point of $\CstarPath $.  So, by \cite {SimpleNbd}, there is
a large enough prefix $\mu $ of $\alpha $ such that $$ \{\alpha \} = X_\mu \cap \CstarPath .  \equationmark
AlphaIsolated $$ Therefore $\mu $ is a prefix of $\alpha $, and by enlarging $\mu $ a bit, we may clearly assume that it
is of the form $$ \mu = \nu \underbrace {\gamma \gamma \ldots \gamma }_{n}, $$ for some integer $n$.  Assuming by
contradiction that $\gamma $ has an entry, we can find a path of the form $$ \gamma _{_1}\gamma _{_2}\ldots \gamma
_{_k}\varepsilon , $$ with $\varepsilon \neq \gamma _{_k+1}$ (here $k+1$ stands for $1$, in case $k=|\gamma |$), which
we may extend as far as possible, obtaining a maximal path $$ \beta = \gamma _{_1}\gamma _{_2}\ldots \gamma
_{_k}\varepsilon a _{_1}a _{_2}a _{_3}\ldots $$ Considering the path $$ \alpha '= \mu \beta = \nu \underbrace {\gamma
\gamma \ldots \gamma }_{n}\gamma _{_1}\gamma _{_2}\ldots \gamma _{_k}\varepsilon a _{_1}a _{_2}a _{_3}\ldots $$ notice
that $\alpha '$ is clearly also maximal hence $\alpha '\in \CstarPath $.  Since $\mu $ is a prefix of $\alpha '$, we see
that $$ \alpha '\in X_\mu \cap \CstarPath , $$ contradicting \cite {AlphaIsolated}.  This proves that $\gamma $ has no
entry.

Conversely, supposing that $\gamma $ has no entry, it is easy to see that the only path in $\CstarPath $ extending $$
\mu :=\nu \gamma $$ is $\alpha $, itself.  So $$ X_\mu \cap \CstarPath = \{\alpha \}, $$ showing that $\alpha $ is
indeed isolated in $\CstarPath $.  \endProof

Putting together our findings so far we have:

\state Proposition \label GraphTopFree For every graph $E$ with no sinks, one has: \izitem \zitem The partial action
$\tau $ of\/ $\F $ on $\FullPath $, described in \cite {SummaryTau}, is topologically free.  \zitem The restriction of\/
$\tau $ to $\CstarPath $ is topologically free if and only if every cycle in $E$ has an entry.

\Proof The main ingredients of this proof have already been taken care of, so all we need is to piece them together.

Given $g$ in $\F $, one has by \cite {GraFixPoint} that the fixed point set $F_g$ is either empty or consists exactly of
one infinite path.  In the second case this infinite path is never isolated in $\FullPath $ by \cite {IsolatedPoints/i},
so the interior of $F_g$ is always empty, and hence $\tau $ is topologically free on $\FullPath $.

We next prove the only if part of (ii) via the contra-positive.  If there exists a cycle $\gamma $ with no entry, then
$$ F_\gamma = \{\gamma \gamma \gamma \ldots \}, $$ by \cite {GraFixPoint}, and moreover $\gamma \gamma \gamma \ldots $
is isolated in $\CstarPath $ by \cite {IsolatedPoints/ii}, whence $F_g$ is open and consequently $\tau $ is not
topologically free on $\CstarPath $.

Conversely, observe that any nonempty fixed point set necessarily looks like $$ F_g = \{\nu \gamma \gamma \gamma \ldots
\}, $$ by \cite {GraFixPoint/iii} (recall that if the hypothesis ``$|\mu |<|\nu |$'' in \cite {GraFixPoint/iii} fails,
then it will hold for $g\inv $, while $F_g=F_{g\inv }$), where $\gamma $ is a cycle.  If every cycle has an entry, then
$\nu \gamma \gamma \gamma \ldots $ is not isolated in $\CstarPath $ by \cite {IsolatedPoints/ii}, so the interior of
$F_g$ is empty.  This proves that $\tau $ is topologically free on $\CstarPath $.  \endProof

We may then use the above result in conjunction with \cite {topfree} to obtain the following result, sometimes referred
to as the \"{uniqueness Theorem for graph C*-algebras}:

\state Theorem Let $E$ be a graph with no sinks such that every cycle in $E$ has an entry.  Also let $B$ be a C*-algebra
which is generated by a set $$ \{\p v' : v\in E^0\} \cup \{\s a ' : a \in E^1\}, $$ where the $\p v'$ are mutually
orthogonal projections and the $\s a '$ are partial isometries satisfying the relations defining $C^*(E)$, namely \cite
{DefineGraphAlg/i-iii}.  Then $B$ is naturally isomorphic to $C^*(E)$ provided $\p v'\neq 0$, for every $v$ in $E^0$.

\Proof Consider the *-homomorphism $$ \varphi :C^*(E)\to B $$ mapping the $\p v$ to $\p v'$, and the $\s a $ to $\s a
'$, given by the universal property of $C^*(E)$.  Since $B$ is supposed to be generated by the $\p v'$ and the $\s a '$
it is clear that $\varphi $ is surjective.  So the proof will be concluded once we prove that $\varphi $ is injective.
In order to do this we argue by contradiction, supposing that $$ \Ker (\varphi )\neq \{0\}.  $$

Identifying $C^*(E)$ and $C_0(\CstarPath ) \redrt \tau \F $, by \cite {GraphAmena/b}, we may see $\Ker (\varphi )$ as an
ideal of the latter algebra.  So we have by \cite {topfree} and \cite {GraphTopFree/ii}, that $$ K:= \Ker (\varphi )\cap
C_0(\CstarPath ) \neq \{0\}.  $$

Moreover $K$ is an invariant ideal of $C_0(\CstarPath )$ by \cite {IntersInvar}.  So, writing $K = C_0(U)$, where $U$ is
an open subset of $\CstarPath $, we have that $U$ is nonempty and invariant.

Using \cite {DensityIandII} we may find a maximal path $\alpha $ in $U$, and then by \cite {SimpleNbd} we see that there
exists a finite path $\mu $ such that $$ \alpha \in X_\mu \cap \CstarPath \subseteq U.  $$ Adopting the notation \def \X
{X^\flat }$ \X _\mu := X_\mu \cap \CstarPath , $ we rewrite the above as $$ \alpha \in \X _\mu \subseteq U.  $$ The
standard form of $\mu \inv $ is $d(\mu )\mu \inv $, so $\tau _{\mu \inv }$ maps $\X _\mu $ to $\X _{d(\mu )}$.
Therefore the invariance of $U$ implies that $$ \X _{d(\mu )} = \tau _{\mu \inv }(\X _\mu ) = \tau _{\mu \inv }(\X _\mu
\cap U) \subseteq U.  $$ Consequently the characteristic function of $\X _{d(\mu )}$ lies in $C_0(U)$, and hence also in
the null space of $\varphi $.  Noticing that said function identifies with $\p {d(\mu )}$ under \cite {GraphAmena/b}, we
then conclude that $$ \p {d(\mu )}' = \varphi (\p {d(\mu )}) = 0, $$ a contradiction.  \endProof

Let us now discuss the question of simplicity of $C^*(E)$.

\definition \label DefTransitive We shall say that a graph $E$ is: \izitem \zitem \subjex {transitive}{transitive graph}
if, for any two vertices $v$ and $w$ in $E^0$, there exists a finite path $\delta $ such that $d(\delta )=v$, and
$r(\delta )=w$.  \zitem \subjex {weakly transitive}{weakly transitive graph} (also known as \subjex {co-final}{co-final
graph}) if for any path $\alpha $ in $\CstarPath $ and any vertex $w$, there exists a finite path $\delta $ such that
$d(\delta )$ is some vertex along $\alpha $, and $r(\delta )=w$.

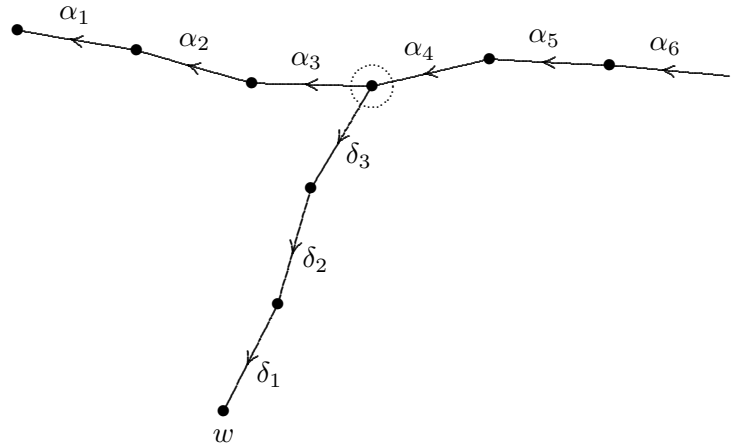
\null \hfill \beginpicture \setcoordinatesystem units <0.0080truecm, -0.0080truecm> point at 1000 1000 \setplotarea x
from -600 to 600, y from 0 to 600 \put {$\bullet $} at 301 -18 \arrow <0.15cm> [0.4,1] from 500 0 to 390 -10 \plot 390
-10 301 -18 / \put {$\alpha _{6}$} at 394 -50 \put {$\bullet $} at 101 -28 \arrow <0.15cm> [0.4,1] from 301 -18 to 191
-23 \plot 191 -23 101 -28 / \put {$\alpha _{5}$} at 193 -63 \put {$\bullet $} at -94 17 \arrow <0.15cm> [0.4,1] from 101
-28 to -6 -3 \plot -6 -3 -94 17 / \put {$\alpha _{4}$} at -15 -42 \put {$\bullet $} at -294 12 \arrow <0.15cm> [0.4,1]
from -94 17 to -204 14 \plot -204 14 -294 12 / \put {$\alpha _{3}$} at -203 -26 \put {$\bullet $} at -486 -43 \arrow
<0.15cm> [0.4,1] from -294 12 to -400 -18 \plot -400 -18 -486 -43 / \put {$\alpha _{2}$} at -389 -57 \put {$\bullet $}
at -684 -76 \arrow <0.15cm> [0.4,1] from -486 -43 to -595 -61 \plot -595 -61 -684 -76 / \put {$\alpha _{1}$} at -588
-101 \put {$\bullet $} at -196 188 \arrow <0.15cm> [0.4,1] from -94 17 to -150 111 \plot -150 111 -196 188 / \put
{$\delta _3$} at -116 132 \put {$\bullet $} at -251 381 \arrow <0.15cm> [0.4,1] from -196 188 to -226 294 \plot -226 294
-251 381 / \put {$\delta _2$} at -188 305 \put {$\bullet $} at -341 559 \arrow <0.15cm> [0.4,1] from -251 381 to -301
479 \plot -301 479 -341 559 / \put {$\delta _1$} at -265 497 \put {$w$} at -341 599 \setdots <1.5pt> \circulararc 360
degrees from -59 17 center at -94 17 \endpicture \hfill \null

\nostate \label PictorialWeakTran A pictorial way to understand weak transitivity is to imagine that you are standing at
the end (range) of the path $\alpha $, attempting to reach a given vertex $w$.  Even if there is no path to get you
there, in a weakly transitive graph you are able to back up a few edges along $\alpha $ before finding a path $\delta $
taking you to $w$.

\state Proposition \label GrapgMinimal The partial action $\tau $ of\/ $\F $ on $\CstarPath $ is minimal if and only if
$E$ is weakly transitive.

\Proof Supposing that $E$ is weakly transitive, we will show that the orbit of every $\alpha $ in $\CstarPath $ is
dense.  In view of \cite {DensityIandII}, it is enough to prove that any maximal path lies in the closure of $\Orb
(\alpha )$.

So, let $U$ be a neighborhood of a given maximal path $\beta $, and we must then prove that $U$ has a nonempty
intersection with the orbit of $\alpha $.  By \cite {SimpleNbd} there exists a finite path $\mu $ such that $$ \beta \in
X_\mu \cap \CstarPath \subseteq U.  $$ \"{En passant} we observe that $\mu $ is a prefix of $\beta $.

We next apply the hypothesis that $E$ is weakly transitive to $\alpha $ and $d(\mu )$, obtaining a path $\delta $ whose
source is some vertex $v$ along $\alpha $, and such that $r(\delta ) = d(\mu )$.

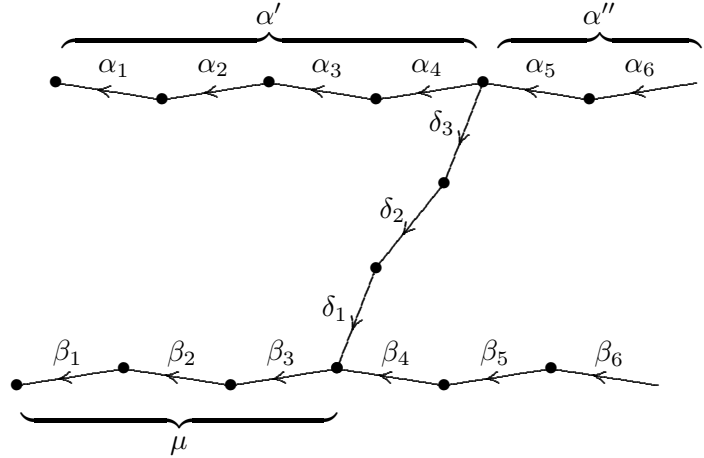
\null \hfill \beginpicture \setcoordinatesystem units <0.0080truecm, 0.0080truecm> point at 0 0 \setplotarea x from -900
to 300, y from 0 to 600 \put {$\bullet $} at -900 0 \arrow <0.15cm> [0.4,1] from -722 28 to -829 11 \plot -829 11 -900 0
/ \put {$\beta _1$} at -815 54 \put {$\bullet $} at -722 28 \arrow <0.15cm> [0.4,1] from -544 -0 to -651 17 \plot -651
17 -722 28 / \put {$\beta _2$} at -625 53 \put {$\bullet $} at -544 -0 \arrow <0.15cm> [0.4,1] from -367 28 to -473 11
\plot -473 11 -544 -0 / \put {$\beta _3$} at -460 54 \put {$\bullet $} at -367 28 \arrow <0.15cm> [0.4,1] from -189 -0
to -296 17 \plot -296 17 -367 28 / \put {$\beta _4$} at -269 53 \put {$\bullet $} at -189 -0 \arrow <0.15cm> [0.4,1]
from -11 28 to -118 11 \plot -118 11 -189 -0 / \put {$\beta _5$} at -104 54 \put {$\bullet $} at -11 28 \arrow <0.15cm>
[0.4,1] from 167 -0 to 60 17 \plot 60 17 -11 28 / \put {$\beta _6$} at 86 53 \put {$\bullet $} at -367 28 \arrow
<0.15cm> [0.4,1] from -302 196 to -341 95 \plot -341 95 -367 28 / \put {$\delta _1$} at -371 129 \put {$\bullet $} at
-302 196 \arrow <0.15cm> [0.4,1] from -189 336 to -257 252 \plot -257 252 -302 196 / \put {$\delta _2$} at -275 293 \put
{$\bullet $} at -189 336 \arrow <0.15cm> [0.4,1] from -124 504 to -163 403 \plot -163 403 -189 336 / \put {$\delta _3$}
at -193 437 \put {$\bullet $} at -124 504 \arrow <0.15cm> [0.4,1] from 53 476 to -53 493 \plot -53 493 -124 504 / \put
{$\alpha _5$} at -27 529 \put {$\bullet $} at 53 476 \arrow <0.15cm> [0.4,1] from 231 504 to 125 487 \plot 125 487 53
476 / \put {$\alpha _6$} at 138 530 \put {$\bullet $} at -302 476 \arrow <0.15cm> [0.4,1] from -124 504 to -231 487
\plot -231 487 -302 476 / \put {$\alpha _4$} at -217 530 \put {$\bullet $} at -480 504 \arrow <0.15cm> [0.4,1] from -302
476 to -409 493 \plot -409 493 -480 504 / \put {$\alpha _3$} at -383 529 \put {$\bullet $} at -658 476 \arrow <0.15cm>
[0.4,1] from -480 504 to -587 487 \plot -587 487 -658 476 / \put {$\alpha _2$} at -573 530 \put {$\bullet $} at -835 504
\arrow <0.15cm> [0.4,1] from -658 476 to -764 493 \plot -764 493 -835 504 / \put {$\alpha _1$} at -738 529 \put
{$\underbrace {\hbox to 4.2truecm{\hfill }}_{\textstyle \mu }$} at -630 -80 \put {$\overbrace {\hbox to 5.5truecm{\hfill
}}^{\textstyle \alpha '}$} at -480 600 \put {$\overbrace {\hbox to 2.7truecm{\hfill }}^{\textstyle \alpha ''}$} at 70
600 \endpicture \hfill \null

\bigskip Splitting $\alpha $ at $v$, let us write $\alpha = \alpha '\alpha ''$, with $d(\alpha ') = v = r(\alpha '')$.
In particular $\mu \delta \alpha ''$ is a well defined path.  Assuming we have \"{backed up along} $\alpha $ as few
steps as possible, as described in \cite {PictorialWeakTran}, we have that the rightmost edge of $\alpha '$ and the
rightmost edge of $\delta $ differ, so that $\delta {\alpha '}\inv $ is in standard form.  It then follows that $\mu
\delta {\alpha '}\inv $ is also in standard form and the conclusion follows from the fact that $$ \tau _{\mu \delta
{\alpha '}\inv }(\alpha ) = \mu \delta \alpha '' \in X_\mu \cap \CstarPath \cap \Orb (\alpha ) \subseteq U \cap \Orb
(\alpha ).  $$

Conversely, suppose that $E$ is not weakly transitive.  So we may pick a path $\alpha $ and a vertex $w$ violating the
condition in \cite {DefTransitive/ii}.  This means that it is impossible to replace a prefix of $\alpha $ by another, in
order to obtain a path with range $w$.  Observe that, by \cite {ThetaGSameAsTauMuNu}, the process of replacing prefixes
is nothing but applying a partial homeomorphism $\tau _g$ to a path, so we deduce from the above that $$ \Orb (\alpha
)\cap X_w=\emptyset .  $$ It follows that $\Orb (\alpha )$ is not a dense set, so its closure is an invariant subset,
refuting the minimality of $\tau $.  \endProof

As a consequence we obtain:

\state Theorem Let $E$ be a weakly transitive graph with no sinks such that every cycle in $E$ has an entry.  Then
$C^*(E)$ is a simple C*-algebra.

\Proof Still referring to the partial action $\tau $ of $\F $ on $\CstarPath $ giving rise to $C^*(E)$ in \cite
{CStarEAsCrossProd}, we have that $\tau $ is topologically free by \cite {GraphTopFree/b}, and minimal by \cite
{GrapgMinimal}.  Thus $C_0(\CstarPath )\redrt \tau \F $ is simple by \cite {TofFreeMinimal}.  The conclusion then
follows from \cite {GraphAmena/b}.  \endProof

Our next main result will be an application of \cite {ClassIdeal} to graph algebras.  In order to succeed in this
endeavor, we must therefore find sufficient conditions for $\tau $ to be topologically free on every closed invariant
subset of $\CstarPath $.

With this in mind we will now introduce another important graph theoretical concept.

\definition \label DefineRecurrent Let $\gamma $ be a cycle in $E$.  We say that $\gamma $ is: \izitem \zitem \subjex
{recurrent}{recurrent cycle} if there exists another cycle $\beta $, with the same range and source as $\gamma $, such
$$ \gamma \beta \gamma \gamma \gamma \ldots \neq \gamma \gamma \gamma \ldots .  $$ \zitem \subjex
{transitory}{transitory cycle} if it is not recurrent.

In a recurrent cycle, the existence of the above path $\beta $ means that it is possible to exit $\gamma $ and return to
it at a later time.  In a transitory cycle this is impossible.

\null \hfill \beginpicture \setcoordinatesystem units <0.0080truecm, 0.0080truecm> point at 0 0 \setplotarea x from -100
to 300, y from -100 to 300 \put {$\bullet $} at 200 200 \arrow <0.15cm> [0.4,1] from 200 0 to 200 120 \plot 200 120 200
200 / \put {$\gamma _9$} at 240 100 \put {$\bullet $} at 71 353 \arrow <0.15cm> [0.4,1] from 200 200 to 123 292 \plot
123 292 71 353 / \put {$\gamma _8$} at 167 302 \put {$\bullet $} at -126 388 \arrow <0.15cm> [0.4,1] from 71 353 to -47
374 \plot -47 374 -126 388 / \put {$\gamma _7$} at -20 410 \put {$\bullet $} at -299 288 \arrow <0.15cm> [0.4,1] from
-126 388 to -229 328 \plot -229 328 -299 288 / \put {$\gamma _6$} at -232 373 \put {$\bullet $} at -367 100 \arrow
<0.15cm> [0.4,1] from -299 288 to -340 175 \plot -340 175 -367 100 / \put {$\gamma _5$} at -371 208 \put {$\bullet $} at
-299 -88 \arrow <0.15cm> [0.4,1] from -367 100 to -326 -13 \plot -326 -13 -299 -88 / \put {$\gamma _4$} at -371 -8 \put
{$\bullet $} at -126 -188 \arrow <0.15cm> [0.4,1] from -299 -88 to -195 -148 \plot -195 -148 -126 -188 / \put {$\gamma
_3$} at -232 -173 \put {$\bullet $} at 71 -153 \arrow <0.15cm> [0.4,1] from -126 -188 to -7 -167 \plot -7 -167 71 -153 /
\put {$\gamma _2$} at -20 -210 \put {$\bullet $} at 200 0 \arrow <0.15cm> [0.4,1] from 71 -153 to 149 -61 \plot 149 -61
200 0 / \put {$\gamma _1$} at 166 -103 \put {$\bullet $} at -0 200 \arrow <0.15cm> [0.4,1] from 200 200 to 80 200 \plot
80 200 -0 200 / \put {$\nu _3$} at 100 240 \put {$\bullet $} at 0 -0 \arrow <0.15cm> [0.4,1] from -0 200 to 0 80 \plot 0
80 0 -0 / \put {$\nu _2$} at -40 100 \put {$\bullet $} at 200 0 \arrow <0.15cm> [0.4,1] from 0 -0 to 120 0 \plot 120 0
200 0 / \put {$\nu _1$} at 100 -40 \put {\eightrm A recurrent cycle $\gamma $.} at -50 -300 \endpicture \hfill \null

\bigskip The cycle $\gamma $ in the diagram above is recurrent because the path $$ \gamma \nu _1\nu _2\nu _3\gamma
_9\gamma \gamma \gamma $$ temporarily exits $\gamma $ before returning to it.

The comparison between infinite paths in \cite {DefineRecurrent/i} cannot be replaced by the comparison of $\gamma \beta
\gamma $ with a finite path of the form $\gamma \gamma \ldots \gamma $.  In fact, given a recurrent cycle $\sigma $,
consider the cycle $\gamma =\sigma \sigma $.  Then, taking $\beta =\sigma $, we have that the infinite paths $\gamma
\beta \gamma \gamma \gamma \ldots $ and $\gamma \gamma \gamma \gamma \ldots $ coincide, whereas $\gamma \beta \gamma $
does not coincide with any path of the form $\gamma \gamma \ldots \gamma $.  This is because $$ |\gamma \beta \gamma | =
5|\sigma | \and |\gamma \gamma \ldots \gamma | = 2n|\sigma |, $$ where $n$ is the number of repetitions.  Incidentally,
notice that $\sigma \sigma $ is recurrent, provided $\sigma $ is.

Observe that a recurrent cycle must necessarily have an entry.

\state Proposition \label GraphTopFreeEveryInvar Let $E$ be a graph with no sinks.  Then $\tau $ is topologically free
on every closed $\tau $-invariant subset of $\CstarPath $ if and only every cycle is recurrent.

\Proof Suppose that every cycle in $E$ is recurrent and let $C$ be a closed invariant subset of $\CstarPath $.  To show
that $\tau $ is topologically free on $C$ we need to prove that, for every $g$ in $\F \setmenos \{1\}$, the interior of
$F_g\cap C$, relative to $C$, is empty.  This is obviously true if $F_g\cap C$ itself is empty, so let us assume
otherwise.

Using \cite {GraFixPoint} we may assume that $g$ has a standard form $g=\nu \gamma \nu \inv $, where $\gamma $ is a
cycle, and that $ F_g\cap C $ consists of the single point $$ \alpha = \nu \gamma \gamma \gamma \ldots $$ (for this we
might need to replace $g$ by $g\inv $, if necessary, observing that $F_g=F_{g\inv }$).  In order to show that the
interior of $F_g\cap C$ relative to $C$ is empty, it is clearly enough to show that $\alpha $ is not isolated in $C$.

By hypothesis $\gamma $ is recurrent, so we may choose a finite path $\beta $ as in \cite {DefineRecurrent/i}.  For each
$n$, consider the path $$ \alpha _n = \tau _{\nu \gamma ^n\beta \nu \inv }(\alpha ) = \nu \gamma ^n\beta \gamma \gamma
\gamma \ldots = \nu \underbrace {\gamma \gamma \ldots \gamma }_{n}\beta \gamma \gamma \gamma \ldots $$ Since $C$ is
invariant, it is clear that $\alpha _n$ lies in $C$.  It is also easy to see that $\alpha _n\neq \alpha $, and $$ \alpha
_n\convrg n \alpha , $$ so $\alpha $ is an accumulation point of $C$, hence not isolated.

Conversely, suppose that $\tau $ is topologically free on every closed invariant subset and let $\gamma $ be a cycle in
$E$.  Let $$ \alpha =\gamma \gamma \gamma \ldots \in \CstarPath , $$ and let $C$ be the closure of the orbit of $\alpha
$.  Then it is evident that $\alpha $ is fixed by $\gamma $, so we have by \cite {GraFixPoint} that $$ F_\gamma \cap
C=\{\alpha \}.  $$ Since $\tau $ is topologically free on $C$, by assumption, we have that the interior of $F_\gamma
\cap C$ relative to $C$ is empty.  In particular $X_\gamma \cap C$, which is an open neighborhood of $\alpha $, cannot
be contained in $F_\gamma \cap C$.  Consequently $$ \emptyset \neq (X_\gamma \cap C)\setmenos (F_\gamma \cap C) =
(X_\gamma \cap C)\setmenos \{\alpha \} = \big (X_\gamma \setmenos \{\alpha \}\big )\cap C.  $$ So, $X_\gamma \setmenos
\{\alpha \}$ is an open set intersecting $C$, from where we see that $$ \big (X_\gamma \setmenos \{\alpha \}\big ) \cap
\Orb (\alpha ) \neq \emptyset .  $$

Picking some $\alpha '\neq \alpha $ in $X_\gamma \cap \Orb (\alpha )$, notice that, since the action of $\F $ consists
of prefix replacements (see \cite {ThetaGSameAsTauMuNu}), we have that $\alpha '$, just like $\alpha $, is an eventually
periodic path ending in $\gamma \gamma \gamma \ldots $.  Since $\alpha '$ is in $X_\gamma $, we see that $\gamma $ is a
prefix of $\alpha '$.  So $\alpha '$ must be a path of the form $$ \alpha ' = \gamma \beta \gamma \gamma \gamma \ldots ,
$$ proving that $\gamma $ is a recurrent cycle.  \endProof

With this we obtain the following important classification of ideals in graph C*-algebras.

\state Theorem \label ClassIdealGraph Let $E$ be a graph with no sinks, such that every cycle is recurrent.  Then there
is a one-to-one correspondence between open $\tau $-invariant subsets of $\CstarPath $ and closed two-sided ideals in
$C^*(E)$.

\Proof Observe that the partial action in point satisfies both \cite {ClassIdeal/i\&ii}, by \cite {GraphAmena/a} and
because the free group is known to be exact.

So the statement follows from \cite {GraphTopFreeEveryInvar} and \cite {ClassIdeal}, once we realize $C^*(E)$ as
$C_0(\CstarPath ) \redrt \tau \F $, by \cite {GraphAmena}.  \endProof

There is a lot more that can be said about the ideals of $C^*(E)$ under the above hypothesis.  In \ref
{BatesHongRaeburnSzymanski/2002} these are characterized in terms of \"{hereditary directed subsets} of $E$, plus some
extra ingredients.  Also, the quotient of $C^*(E)$ by any ideal is shown to be again a graph C*-algebra.

With some extra effort, these results may also be obtained via our techniques.  We leave them as exercises for the
interested reader.

\nrem Most of the results in this chapter may be generalized to Exel-Laca algebras, by employing very similar techniques
\ref {ExelLaca/1999}.

The study of the ideal structure of graph C*-algebras was initiated in \ref {HuefRaeburn/1997} for finite graphs,
followed by \ref {BatesHongRaeburnSzymanski/2002}, where the general infinite case is considered.

Graphs in which every cycle has an entry are said to satisfy condition (I).  It is a generalization of a similar
condition introduced in \ref {CuntzKrieger/1980}.  On the other hand, graphs in which every cycle is recurrent are said
to satisfy condition (K) \ref {KumjianPaskRaeburnRenault/1997}.  This in turn generalizes condition (II) of \ref
{Cuntz/1981}.

Most of the above results describing the various algebras associated to a graph $E$ as partial crossed products,
beginning with \cite {BestModelForTTe}, could also be proved directly by first \"{guessing} the appropriate partial
dynamical system and then proving that the corresponding crossed product algebra has the correct universal properties.
This strategy would evidently be logically correct and a lot shorter than the one adopted here, but we have opted to
follow the above more constructive, and hopefully more pedagogical approach, not least because it may be used in other
situations where guesswork might not be an option.

\endgroup

\newcount \bibno \def \jrn #1, #2 (#3), #4-#5;{\sl #1, \bf #2 \rm (#3), #4--#5} \def \Bibitem #1/#2 #3; #4; #5 \par
{\medbreak \global \advance \bibno by 1 \immediate \write 3 {\string \def \textbackslash #1\ntoP #2;{\number \bibno
}}\item {[\Htarget {paper.#1/#2}{\number \bibno }]} #3, ``#4'', #5.}

\def \bibitem #1/#2 #3; #4; #5 \par {\Bibitem #1/#2 #3; #4; \jrn #5; \par } \def \nobibitem #1\par {}

\tolerance 700

\vfill \eject

\Headlines {}{} \null \bigbreak \vskip 1cm

\noindent \Htarget {extra.References}{\bigrm References} \write 2 {\string \OtherIndexEntry {References} {\number
\pageno }}

\vskip 2cm \begingroup \eightpoint \nobreak \medskip \frenchspacing

\Bibitem Abadie/1999 F. Abadie; Sobre a\c c\~oes Parciais, Fibrados de Fell e Grup\'oides; PhD Thesis, University of
S\~ao Paulo, 1999

\bibitem Abadie/2003 F. Abadie; Enveloping actions and Takai duality for partial actions; J. Funct. Analysis, 197
(2003), 14-67

\Bibitem Abadie/2004 F. Abadie; On partial actions and groupoids; Proc. Amer. Math. Soc., 132 (2004), 1037-1047

\Bibitem AbadieDokuchaevExelSimon/2008 F. Abadie, M. Dokuchaev, R. Exel and J. J. Sim\'on; Morita equivalence of partial
group actions and globalization; Trans. Amer. Math. Soc., to appear

\bibitem AbadieMarti/2009 F. Abadie and L. Mart\'{\i } P\'erez; On the amenability of partial and enveloping actions;
Proc. Amer. Math. Soc., 137 (2009), no. 11, 3689-3693

\bibitem AnantharamanDelaroche/1987 C. Anantharaman-Delaroche; Syst\`emes dynamiques non commutatifs et
moyenna\-bilit\'e; Math. Ann., 279 (1987), 297-315

\Bibitem AraExelKatsura/2011 P. Ara, R. Exel and T. Katsura; Dynamical systems of type $(m,n)$ and their C*-algebras;
Ergodic Theory Dynam. Systems, 33 (2013), no. 5, 1291-1325

\bibitem ArchboldSpielberg/1993 R. J. Archbold and J. S. Spielberg; Topologically free actions and ideals in discrete
C*-dynamical systems; Proc. Edinb. Math. Soc., 37 (1993), 119-124

\Bibitem Arveson/1981 W. Arveson; An Invitation to C*-Algebras; Springer-Verlag, 1981

\bibitem Arveson/1991 W. Arveson; C*-algebras associated with sets of semigroups of isometries; Internat. J. Math., 2
(1991), no. 3, 235-255

\bibitem BagioLazzarinPaques/2010 D. Bagio, J. Lazzarin and A. Paques; Crossed products by twisted partial actions:
separability, semisimplicity, and Frobenius properties; Comm. Algebra, 38 (2010), no. 2, 496-508

\bibitem BatesHongRaeburnSzymanski/2002 T. Bates, J. H. Hong, I. Raeburn and W. Szyma\'nski; The ideal structure of the
C∗-algebras of infinite graphs; Illinois J. Math., 46 (2002), no. 4, 1159-1176

\bibitem Blackadar/1985 B. Blackadar; Shape theory for C*-algebras; Math. Scand., 56 (1985), 249-275

\Bibitem BoavaExel/2010 G. Boava and R. Exel; Partial crossed product description of the C*-algebras associated with
integral domains; Proc. Amer. Math. Soc., to appear (2010), arXiv: 1010.0967v2

\bibitem BostConnes/1995 J.-B. Bost and A. Connes; Hecke algebras, type III factors and phase transitions with
spontaneous symmetry breaking in number theory; Selecta Math. (N.S.), 1 (1995), no. 3, 411-457

\bibitem Brown/1977 L. G. Brown; Stable Isomorphism of Hereditary Subalgebras of C*-Algebras; Pacific J. Math., 71
(1977), 335-348

\bibitem BrownGreenRieffel/1977 L. G. Brown, P. Green and M. A. Rieffel; Stable isomorphism and strong Morita
equivalence of C*-algebras; Pacific J. Math., 71 (1977), 349-363

\bibitem BrownMigoShen/1994 L. G. Brown, J. A. Mingo, and N. T. Shen; Quasi-multipliers and embeddings of Hilbert
C*-bimodules; Canad. J. Math., 46 (1994), 1150-1174

\Bibitem BrownOzawa/2008 N.P. Brown and N. Ozawa; C*-algebras and finite-dimensional approximations; Graduate Studies in
Mathematics, 88, American Mathematical Society, Providence, RI, 2008

\bibitem Burdak/2004 Z. Burdak; On decomposition of pairs of commuting isometries; Ann. Polon. Math., 84 (2004), 121-135

\bibitem BussExel/2011 A. Buss and R. Exel; Twisted actions and regular Fell bundles over inverse semigroups;
Proc. London Math. Soc., 103 (2011), 235-270, arXiv: 1003.0613

\bibitem BussExel/2012 A. Buss and R. Exel; Fell bundles over inverse semigroups and twisted \'etale groupoids;
J. Operator Theory, 67 (2012), 153-205, arXiv: 0903.3388

\bibitem BussExelMeyer/2012 A. Buss, R. Exel and R. Meyer; Inverse semigroup actions as groupoid actions; Semigroup
Forum, 85 (2012), 227-243, arXiv: 1104.0811v2

\nobibitem Cartan/1904 \'E. Cartan; Sur la structure des groupes infinis de transformation; Ann. Sci. \'Ecole
Norm. Sup., 21 (1904), 153-206

\bibitem Coburn/1967 L.A. Coburn; The C*-algebra generated by an isometry; Bull. Amer. Math. Soc., 73 (1967), 722-726

\bibitem Combes/1984 F. Combes; Crossed products and Morita equivalence; Proc. London Math. Soc., 49 (1984), no. 2,
289-306

\bibitem CrispLaca/2002 J. Crisp and M. Laca; On the Toeplitz algebras of right-angled and finite-type Artin groups;
J. Aust. Math. Soc., 72 (2002), 223-245

\bibitem CrispLaca/2007 J. Crisp and M. Laca; Boundary quotients and ideals of Toeplitz C*-algebras of Artin groups;
J. Funct. Analysis, 242 (2007), 127-156

\bibitem CuntzKrieger/1980 J. Cuntz and W. Krieger; A class of C*-algebras and topological Markov chains; Invent. Math.,
56 (1980), no. 3, 251-268

\bibitem Cuntz/1981 J. Cuntz; A class of C*-algebras and topological Markov chains II: Reducible chains and the
Ext-functor for C*-algebras; Invent. Math., 63 (1981), 25-40

\bibitem CurtoMurphyWilliams/1984 R. Curto, P.  Muhly and D. Williams; Cross products of strongly Morita equivalent
C∗-algebras; Proc. Amer. Math. Soc., 90 (1984), no. 4, 528-530

\Bibitem Davidson/1996 K. R. Davidson; C*-Algebras by Example; Fields Institute Monographs, 1996

\bibitem DokuchaevExel/2005 M. Dokuchaev and R. Exel; Associativity of crossed products by partial actions, enveloping
actions and partial representations; Trans. Amer. Math. Soc., 357 (2005), 1931-1952 (electronic), arXiv: math.RA/0212056

\bibitem DokuchaevExelPiccione/1999 M. Dokuchaev, R. Exel and P. Piccione; Partial representations and partial group
algebras; J. Algebra, 226 (2000), 505-532, arXiv: math.GR/9903129

\bibitem DokuchaevExelSimon/2008 M. Dokuchaev, R. Exel and J. J. Sim\'on; Crossed products by twisted partial actions
and graded algebras; J. Algebra, 320 (2008), no. 8, 3278-3310

\bibitem DokuchaevExelSimon/2010 M. Dokuchaev, R. Exel and J. J. Sim\'on; Globalization of twisted partial actions;
Trans. Amer. Math. Soc., 362 (2010), no. 8, 4137-4160

\bibitem DokuchaevNovikov/2010 M. Dokuchaev and B. Novikov; Partial projective representations and partial actions;
J. Pure Appl. Algebra, 214 (2010), 251-268

\bibitem DokuchaevNovikov/2012 M. Dokuchaev and B. Novikov; Partial projective representations and partial actions II;
J. Pure Appl. Algebra, 216 (2012), no. 2, 438-455

\bibitem Douglas/1969 R. G. Douglas; On extending commutative semigroups of isometries; Bull. London Math. Soc., 1
(1969), 157-159

\bibitem Douglas/1972 R. G. Douglas; On the C*-algebra of a one-parameter semigroup of isometries; Acta Math., 128
(1972), 143-151

\bibitem EnomotoWatatani/1980 M. Enomoto and Y. Watatani; A graph theory for C*-algebras; Math. Japon., 25 (1980),
435-442

\Bibitem EffrosHahn/1967 E. G. Effros and F. Hahn; Locally compact transformation groups and C*-algebras; Memoirs of the
American Mathematical Society, No. 75 American Mathematical Society, Providence, R.I. 1967 92 pp.

\bibitem Elliott/1980 G. A. Elliott; Some simple C*-algebras constructed as crossed products with discrete outer
automorphism groups; Publ. Res. Inst. Math. Sci., 16 (1980), 299-311

\bibitem Erdelyi/1968 I. Erd\'elyi; Partial isometries closed under multiplication on Hilbert spaces;
J. Math. Anal. Appl., 22 (1968), 546-551

\bibitem Exel/1994a R. Exel; Circle actions on C*-algebras, partial automorphisms and a generalized Pimsner-Voiculescu
exact sequence; J. Funct. Analysis, 122 (1994), 361-401, arXiv: funct-an/9211001

\bibitem Exel/1994b R. Exel; The Bunce-Deddens algebras as crossed products by partial automorphisms;
Bull. Braz. Math. Soc. (N.S.), 25 (1994), 173-179, arXiv: funct-an/9302001

\bibitem Exel/1995 R. Exel; Approximately finite C*-algebras and partial automorphisms; Math. Scand., 77 (1995),
281-288, arXiv: funct-an/9211004

\bibitem Exel/1997a R. Exel; Twisted partial actions, a classification of regular C*-algebraic bundles; Proc. London
Math. Soc., 74 (1997), 417-443, arXiv: funct-an/9405001

\bibitem Exel/1997b R. Exel; Amenability for Fell bundles; J. Reine Angew. Math., 492 (1997), 41-73, arXiv:
funct-an/9604009

\bibitem Exel/1998 R. Exel; Partial actions of groups and actions of inverse semigroups; Proc. Amer. Math. Soc., 126
(1998), no. 12, 3481-3494

\bibitem Exel/2000 R. Exel; Partial representations and amenable Fell bundles over free groups; Pacific J. Math., 192
(2000), 39-63, arXiv: funct-an/9706001

\Bibitem Exel/2002p R. Exel; Exact Groups, Induced Ideals, and Fell Bundles; preprint, arXiv: math.\break OA/0012091

\bibitem Exel/2002 R. Exel; Exact groups and Fell bundles; Math. Ann., 323 (2002), 259-266

\bibitem Exel/2008a R. Exel; Hecke algebras for protonormal subgroups; J. Algebra, 320 (2008), no. 5, 1771-1813

\bibitem Exel/2008b R. Exel; Inverse semigroups and combinatorial C*-algebras; Bull. Braz. Math. Soc. (N.S.), 39 (2008),
191-313, arXiv: math.OA/0703182

\bibitem Exel/2011 R. Exel; Non-Hausdorff \'etale groupoids; Proc. Amer. Math. Soc., 139 (2011), no. 3, 897-907

\bibitem ExelLaca/1999 R. Exel and M. Laca; Cuntz-Krieger algebras for infinite matrices; J. reine angew. Math., 512
(1999), 119-172, arXiv: funct-an/9712008

\bibitem ExelLaca/2003 R. Exel and M. Laca; Partial dynamical systems and the KMS condition; Comm. Math. Phys., 232
(2003), no. 2, 223-277

\Bibitem ExelLacaQuigg/1997 R. Exel, M. Laca and J. Quigg; Partial dynamical systems and C*-algebras generated by
partial isometries; preprint, arXiv: funct-an/9712007

\bibitem ExelLacaQuigg/2002 R. Exel, M. Laca and J. Quigg; Partial dynamical systems and C*-algebras generated by
partial isometries; J. Operator Theory, 47 (2002), 169-186

\bibitem ExelVieira/2010 R. Exel and F. Vieira; Actions of inverse semigroups arising from partial actions of groups;
J. Math. Anal. Appl., 363 (2010), no. 1, 86-96

\Bibitem FellDoran/1969 J. M. G. Fell; An extension of Mackey's method to Banach *-algebraic bundles; \sl
Mem. Am. Math. Soc., \bf 90 \rm (1969)

\Bibitem FellDoran/1988 J. M. G. Fell and R. S. Doran; Representations of *-algebras, locally compact groups, and Banach
*-algebraic bundles; Pure and Applied Mathematics, vol 125 and 126, Academic Press, 1988

\bibitem FowlerLacaRaeburn/2000 N. Fowler, M. Laca and I. Raeburn; The C*-algebras of infinite graphs;
Proc. Amer. Math. Soc., 128 (2000), no. 8, 2319-2327

\Bibitem Greenleaf/1969 F. P. Greenleaf; Invariant means on topological groups; Mathematical Studies, vol. 16, van
Nostrand-Reinhold, 1969

\bibitem Gromov/2003 M. Gromov; Random walk in random groups; Geom. Funct. Anal., 13 (2003), 73-146

\bibitem HalmosMcLaughlin/1963 P. R. Halmos and J. E.  McLaughlin; Partial isometries; Pacific J. Math., 13 (1963),
585-596

\Bibitem HewittRoss/1970 E. Hewitt and K. A. Ross; Abstract Harmonic Analysis II; Springer-Verlag, 1970

\Bibitem WienerHopf/1931 E. Hopf and N. Wiener; \umlaut Uber eine Klasse singul\umlaut arer Integralgleichungen; \sl
S.-B. Preuß. Akad. Wiss., Phys.-Math. Kl. \bf 30/32 \rm (1931), 696--706

\bibitem HorakMuller/1989 K. Hor\'ak and V. M\umlaut uller; Functional model for commuting isometries; Czechoslovak
Math. J., 39(114) (1989), 370-379

\bibitem HuefRaeburn/1997 A. an Huef and I. Raeburn; The ideal structure of Cuntz-Krieger algebras; Ergodic Theory
Dynam. Systems, 17 (1997), 611-624

\Bibitem JensenThomsen/1991 K. Jensen and K. Thomsen; Elements of KK-Theory; Birkhauser, 1991

\bibitem KatsuraMuhlySimsTomforde/2010 T. Katsura, P. S. Muhly, A. Sims and M. Tomforde; Graph algebras, Exel-Laca
algebras, and ultragraph algebras coincide up to Morita equivalence; J. Reine Angew. Math., 640 (2010), 135-165

\bibitem KawamuraTomiyama/1990 S. Kawamura and J. Tomiyama; Properties of topological dynamical systems and
corresponding C*-algebras; Tokyo J. Math., 13 (1990), 251-257

\bibitem KumjianPaskRaeburnRenault/1997 A. Kumjian, D. Pask, I. Raeburn and J. Renault; Graphs, groupoids, and
Cuntz-Krieger algebras; J. Funct. Anal., 144 (1997), no. 2, 505-541

\bibitem Laca/2000 M. Laca; From endomorphisms to automorphisms and back: dilations and full corners; J. London
Math. Soc., 61 (2000), 893-904

\bibitem LacaRaeburn/1996 M. Laca and I. Raeburn; Semigroup Crossed Products and the Toeplitz Algebras of Nonabelian
Groups; J. Funct. Anal., 139 (1996), 415-440

\Bibitem Lance/1995 E. C. Lance; Hilbert C*-Modules: A Toolkit for Operator Algebraists; London Mathematical Society
Lecture Note Series, 1995

\Bibitem Lawson/1998 M. V. Lawson; Inverse semigroups, the theory of partial symmetries; World Scientific, 1998

\Bibitem Matsumura/2012 M. Matsumura; A characterization of amenability of group actions on C*-algebras; preprint,
arXiv: 1204.3050v1

\bibitem McClanahan/1995 K. McClanahan; K-theory for partial crossed products by discrete groups; J. Funct. Anal., 130
(1995), no. 1, 77-117

\bibitem Morita/1958 K. Morita; Duality for modules and its applications to the theory of rings with minimum condition;
Sci. Rep. Tokyo Kyoiku Daigaku Sect. A, 6 (1958), 83-142

\Bibitem Murphy/1990 G. J. Murphy; C*-algebras and operator theory; Academic Press, 1990

\bibitem Nica/1992 A. Nica; C*-algebras generated by isometries and Wiener-Hopf operators; J. Operator Theory, 27
(1992), 17-52

\nobibitem Ozawa/2006 N. Ozawa; Amenable actions and applications; International Congress of Mathematicians, vol. II,
1563--1580, Eur. Math. Soc., Zurich, 2006

\bibitem PackerRaeburn/1989 J. Packer and I. Raeburn; Twisted crossed products of C*-algebras;
Math. Proc. Camb. Phil. Soc., 106 (1989), 293-311

\bibitem PaquesSantAna/2010 A. Paques and A. Sant'Ana; When is a crossed product by a twisted partial action Azumaya?;
Comm. Algebra, 38 (2010), no. 3, 1093-1103

\bibitem Paschke/1973 W. L. Paschke; Inner product modules over B*-algebras; Trans. Amer. Math. Soc., 182 (1973),
443-468

\Bibitem Pedersen/1979 G. K. Pedersen; C*-Algebras and their automorphism groups; Acad. Press, 1979

\bibitem PhillipsRaeburn/1993 J. Phillips and I. Raeburn; Semigroups of isometries, Toeplitz algebras and twisted
crossed products; Integr. Equ. Oper. Theory, 17 (1993), 579-602

\bibitem Popovici/2004 D. Popovici; A Wold-type decomposition for commuting isometric pairs; Proc. Amer. Math. Soc., 132
(2004), 2303-2314

\bibitem Preston/1954 G. B. Preston; Inverse semi-groups; J. London Math. Soc., 29 (1954), 396-403

\bibitem Quigg/1996 J. C. Quigg; Discrete C*-coactions and C*-algebraic bundles; J. Austral. Math. Soc. Ser. A, 60
(1996), 204-221

\bibitem QuiggRaeburn/1997 J. C. Quigg and I. Raeburn; Characterisations of crossed products by partial actions;
J. Operator Theory, 37 (1997), 311-340

\Bibitem Raeburn/2005 I. Raeburn; Graph algebras; CBMS Regional Conference Series in Mathematics, 103, American
Mathematical Society, Providence, RI, 2005. vi+113 pp

\Bibitem Renault/1980 J. Renault; A groupoid approach to C*-algebras; Lecture Notes in Mathematics vol.~793, Springer,
1980

\bibitem Renault/1991 J. Renault; The ideal structure of groupoid crossed product C*-algebras; J. Operator Theory, 25
(1991), 3-36

\bibitem Rieffel/1974 Marc A. Rieffel; Morita equivalence for C*-algebras and W*-algebras; J. Pure Appl. Algebra, 5
(1974), 51-96

\Bibitem Scarparo/2014 E. P. Scarparo; \'Algebras de Toeplitz generalizadas; master's thesis, Universidade Federal de
Santa Catarina, 2014

\Bibitem Sehnem/2014 C. F. Sehnem; Uma classifica\c c\~ao de fibrados de Fell est\'aveis; master's thesis, Universidade
Federal de Santa Catarina, 2014

\bibitem Sieben/1997 N\'andor Sieben; C*-crossed products by partial actions and actions of inverse semigroups;
J. Austral. Math. Soc. Ser. A, 63 (1997), no. 1, 32-46

\bibitem Sieben/1998 N. Sieben; C*-crossed products by twisted inverse semigroup actions; J. Operator Theory, 39 (1998),
no. 2, 361-393

\bibitem Slocinski/1980 M. S\l oci\'nski; On the Wold type decomposition of a pair of commuting isometries;
Ann. Polon. Math., 37 (1980), 255-262

\Bibitem StratilaVoiculescu/1975 \c S. Str\breve atil\breve a and D.  Voiculescu; Representations of AF-algebras and of
the group $U(\infty )$; Lecture Notes in Mathematics, Vol. 486. Springer-Verlag, 1975, viii+169 pp

\bibitem Wagner/1952 V. V. Wagner; Generalised groups; Proc. USSR Acad. Sci., 84 (1952), 1119-1122

\Bibitem Watatani/1990 Y. Watatani; Index for C*-subalgebras; \sl Mem. Am. Math. Soc., \bf 424 \rm (1990), 117~p

\bibitem Zettl/1983 H. Zettl; A characterization of ternary rings of operators; Adv. Math., 48 (1983), 117-143

\bibitem Zimmer/1977 R. J. Zimmer; Hyperfinite factors and amenable ergodic actions; Invent. Math., 41 (1977), 23-31

\bibitem Zimmer/1978 R. J. Zimmer; Amenable ergodic group actions and an application to Poisson boundaries of random
walks; J. Funct. Anal., 27 (1978), 350-372

\par

\endgroup \def \subjentry #1#2#3{\Hyperref {subj}{#2}{#1,} #3 \par } \vfill \eject \advance \hsize by 0.1truecm \topskip
= 1cm \parindent 0pt

\obeylines

\Htarget {extra.Subject Index}{\bigrm Subject Index}
\write 2 {\string \OtherIndexEntry {Subject Index} {\number \pageno }}

\vskip 3truecm
\doublecolumns

\subjentry {$A$-valued inner-product}{$A$-valued inner-product}{112}
\subjentry {$G$-equivariant map}{$G$-equivariant map}{12}
\goodbreak
$G$-graded
\subjentry {\quad algebra}{$G$-graded algebra}{47}
\subjentry {\quad C*-algebra}{$G$-graded C*-algebra}{125}
\subjentry {$G$-grading}{$G$-grading}{47}
\subjentry {${\fam \bffam \tenbf K}$-algebra}{${\fam \bffam \tenbf K}$-algebra}{29}
\subjentry {*-algebra}{*-algebra}{29}
\subjentry {*-homomorphism}{*-homomorphism}{30}
\subjentry {*-partial representation}{*-partial representation}{51}
\subjentry {adjoint Hilbert bimodule}{adjoint Hilbert bimodule}{117}
\subjentry {adjointable operator}{adjointable operator}{113}
\goodbreak
algebraic
\subjentry {\quad Fell bundle}{algebraic Fell bundle}{202}
\subjentry {\quad globalization}{algebraic globalization}{31}
\subjentry {\quad partial action}{algebraic partial action}{30}
\subjentry {\quad partial dynamical system}{algebraic partial dynamical system}{30}
\subjentry {algebraically equivalent actions}{algebraically equivalent actions}{30}
\goodbreak
amenable
\subjentry {\quad action}{amenable action}{167}
\subjentry {\quad Fell bundle}{amenable Fell bundle}{159}
\subjentry {approximate identity}{approximate identity}{94}
\subjentry {approximation property}{approximation property}{162}
\goodbreak
Banach
\subjentry {\quad *-algebraic bundle}{Banach *-algebraic bundle}{124}
\subjentry {\quad algebraic bundle}{Banach algebraic bundle}{124}
\goodbreak
Bernoulli
\subjentry {\quad global action}{Bernoulli global action}{27}
\subjentry {\quad partial action}{Bernoulli partial action}{28}
\subjentry {Blaschke factor}{Blaschke factor}{263}
\subjentry {C*-algebra}{C*-algebra}{67}
\goodbreak
C*-algebraic
\subjentry {\quad bundle}{C*-algebraic bundle}{124}
\subjentry {\quad crossed product}{C*-algebraic crossed product}{70}
\subjentry {\quad globalization}{C*-algebraic globalization}{234}
\subjentry {\quad partial action}{C*-algebraic partial action}{69}
\subjentry {\quad partial dynamical system}{C*-algebraic partial dynamical system}{69}
\subjentry {C*-seminorm}{C*-seminorm}{70}
\subjentry {Cesaro net}{Cesaro net}{162}
\subjentry {closed under $\vee $, or $\vee $-closed}{closed under $\vee $, or $\vee $-closed}{89}
\subjentry {closure with respect to $\vee $}{closure with respect to $\vee $}{89}
\subjentry {co-final graph}{co-final graph}{337}
\subjentry {compatible partial isometries}{compatible partial isometries}{81}
\subjentry {completion of a pre-Fell-bundle}{completion of a pre-Fell-bundle}{206}
\goodbreak
conditional
\subjentry {\quad expectation in Fell bundles}{conditional expectation in Fell bundles}{180}
\subjentry {\quad expectation}{conditional expectation}{148}
\subjentry {conjugation}{conjugation}{29}
\subjentry {contractive linear map}{contractive linear map}{140}
\subjentry {convex subset of a group}{convex subset of a group}{109}
\subjentry {convolution product}{convolution product}{133}
\subjentry {covariant representation}{covariant representation}{56}
\subjentry {cross sectional C*-algebra}{cross sectional C*-algebra}{135}
\subjentry {crossed product}{crossed product}{43}
\subjentry {cycle in a graph}{cycle in a graph}{332}
\subjentry {dilation of a partial representation}{dilation of a partial representation}{257}
\subjentry {domain, or source map in a graph}{domain, or source map in a graph}{300}
\goodbreak
dual
\subjentry {\quad global action}{dual global action}{223}
\subjentry {\quad partial action}{dual partial action}{223}
\subjentry {edges in a graph}{edges in a graph}{300}
\subjentry {empty path}{empty path}{317}
\subjentry {entry to a cycle}{entry to a cycle}{334}
\subjentry {enveloping action}{enveloping action}{15}
\subjentry {equivalence of partial actions}{equivalence of partial actions}{12}
\subjentry {essential ideal}{essential ideal}{40}
\goodbreak
exact
\subjentry {\quad C*-algebra}{exact C*-algebra}{172}
\subjentry {\quad group}{exact group}{172}
\subjentry {extendable semigroup}{extendable semigroup}{267}
\subjentry {faithful conditional expectation}{faithful conditional expectation}{158}
\goodbreak
Fell
\subjentry {\quad bundle}{Fell bundle}{124}
\subjentry {\quad sub-bundle}{Fell sub-bundle}{172}
\subjentry {fiber}{fiber}{124}
\goodbreak
final
\subjentry {\quad projection}{final projection}{73}
\subjentry {\quad space}{final space}{73}
\subjentry {finite-$\vee $-closure}{finite-$\vee $-closure}{89}
\subjentry {finitely-$\vee $-closed}{finitely-$\vee $-closed}{89}
\subjentry {fixed point}{fixed point}{247}
\goodbreak
Fourier
\subjentry {\quad coefficient}{Fourier coefficient}{140}
\subjentry {\quad ideal}{Fourier ideal}{199}
\subjentry {free partial action}{free partial action}{247}
\goodbreak
full
\subjentry {\quad Hilbert module}{full Hilbert module}{115}
\subjentry {\quad path space}{full path space}{317}
\subjentry {\quad projection}{full projection}{229}
\subjentry {\quad subalgebra}{full subalgebra}{115}
\goodbreak
global
\subjentry {\quad action}{global action}{9}
\subjentry {\quad dynamical system}{global dynamical system}{9}
\subjentry {globalization}{globalization}{15}
\goodbreak
graded
\subjentry {\quad ideal}{graded ideal}{196}
\subjentry {\quad map}{graded map}{48}
\goodbreak
grading
\subjentry {\quad space}{grading space}{47}
\subjentry {\quad subspace in a C*-algebra}{grading subspace in a C*-algebra}{125}
\goodbreak
graph
\subjentry {\quad as a triple}{graph as a triple}{300}
\subjentry {\quad C*-algebra}{graph C*-algebra}{301}
\subjentry {\quad of a partial action}{graph of a partial action}{12}
\goodbreak
group
\subjentry {\quad bundle}{group bundle}{125}
\subjentry {\quad C*-algebra}{group C*-algebra}{135}
\subjentry {\quad representation}{group representation}{51}
\goodbreak
Hardy
\subjentry {\quad projection}{Hardy projection}{263}
\subjentry {\quad space}{Hardy space}{263}
\goodbreak
hereditary
\subjentry {\quad Fell sub-bundle}{hereditary Fell sub-bundle}{174}
\subjentry {\quad subalgebra}{hereditary subalgebra}{115}
\subjentry {Hilbert $A$-$B$-bi\-module}{Hilbert $A$-$B$-bi\-module}{115}
\subjentry {homogeneous space}{homogeneous space}{47}
\subjentry {ideal of a Fell bundle}{ideal of a Fell bundle}{174}
\subjentry {idempotent semilattice}{idempotent semilattice}{18}
\goodbreak
imprimitivity
\subjentry {\quad bimodule}{imprimitivity bimodule}{115}
\subjentry {\quad system}{imprimitivity system}{116}
\subjentry {induced ideal}{induced ideal}{196}
\subjentry {infinite path}{infinite path}{301}
\goodbreak
initial
\subjentry {\quad projection}{initial projection}{73}
\subjentry {\quad space}{initial space}{73}
\subjentry {integrated form}{integrated form}{136}
\subjentry {invariant set under a partial action}{invariant set under a partial action}{13}
\subjentry {inverse semigroup}{inverse semigroup}{18}
\subjentry {involution}{involution}{29}
\subjentry {isomorphic Fell bundles}{isomorphic Fell bundles}{169}
\subjentry {LCH topological space}{LCH topological space}{67}
\goodbreak
least
\subjentry {\quad upper bound in $P$}{least upper bound in $P$}{270}
\subjentry {\quad upper bound}{least upper bound}{270}
\goodbreak
left
\subjentry {\quad Hilbert $A$-module}{left Hilbert $A$-module}{113}
\subjentry {\quad pre-Hilbert $A$-module}{left pre-Hilbert $A$-module}{113}
\subjentry {length function}{length function}{21}
\subjentry {linking algebra}{linking algebra}{119}
\subjentry {local configuration}{local configuration}{311}
\subjentry {maximal path}{maximal path}{329}
\subjentry {minimal partial action}{minimal partial action}{252}
\goodbreak
Morita-Rieffel-equi\-va\-lent
\subjentry {\quad C*-algebras}{Morita-Rieffel-equi\-va\-lent C*-algebras}{115}
\subjentry {\quad partial actions}{Morita-Rieffel-equi\-va\-lent partial actions}{116}
\subjentry {Morita-Rieffel-globalization}{Morita-Rieffel-globalization}{246}
\subjentry {morphism between Fell bundles}{morphism between Fell bundles}{169}
\goodbreak
multiplier
\subjentry {\quad algebra}{multiplier algebra}{40}
\subjentry {\quad pair}{multiplier pair}{39}
\subjentry {NCC}{NCC}{273}
\subjentry {Nica's covariance condition}{Nica's covariance condition}{273}
\goodbreak
non-degenerate
\subjentry {\quad *-homomorphism}{non-degenerate *-homomorphism}{67}
\subjentry {\quad representation}{non-degenerate representation}{93}
\subjentry {non-degenerate{ }algebra}{non-degenerate{ }algebra}{40}
\subjentry {nuclear C*-algebra}{nuclear C*-algebra}{215}
\subjentry {orbit}{orbit}{15}
\subjentry {Ore sub-semigroup}{Ore sub-semigroup}{265}
\goodbreak
orthogonal
\subjentry {\quad Fell bundle}{orthogonal Fell bundle}{168}
\subjentry {\quad partial representation}{orthogonal partial representation}{168}
\goodbreak
partial
\subjentry {\quad action}{partial action}{9}
\subjentry {\quad automorphism}{partial automorphism}{68}
\subjentry {\quad crossed product}{partial crossed product}{43}
\subjentry {\quad dynamical system}{partial dynamical system}{9}
\subjentry {\quad group algebra}{partial group algebra}{62}
\subjentry {\quad group C*-algebra}{partial group C*-algebra}{108}
\subjentry {\quad homeomorphism}{partial homeomorphism}{23}
\subjentry {\quad isometry}{partial isometry}{72}
\subjentry {\quad representation}{partial representation}{51}
\subjentry {\quad symmetry}{partial symmetry}{19}
\goodbreak
path
\subjentry {\quad in a graph}{path in a graph}{301}
\subjentry {\quad of length zero}{path of length zero}{301}
\subjentry {pre-Fell-bundle}{pre-Fell-bundle}{203}
\subjentry {prefix of a path}{prefix of a path}{302}
\subjentry {projection}{projection}{72}
\subjentry {proper map}{proper map}{67}
\subjentry {pseudo-metric}{pseudo-metric}{109}
\subjentry {quasi-lattice}{quasi-lattice}{270}
\goodbreak
range
\subjentry {\quad map in a graph}{range map in a graph}{300}
\subjentry {\quad of a path}{range of a path}{301}
\subjentry {recurrent cycle}{recurrent cycle}{339}
\goodbreak
reduced
\subjentry {\quad cross sectional C*-algebra}{reduced cross sectional C*-algebra}{140}
\subjentry {\quad crossed product}{reduced crossed product}{142}
\subjentry {\quad form}{reduced form}{20}
\subjentry {\quad group C*-algebra}{reduced group C*-algebra}{140}
\goodbreak
regular
\subjentry {\quad representation of $C^*({\mathchoice {\hbox {\rs B}} {\hbox {\rs B}} {\hbox {\rssmall B}} {\hbox {\rssmall B}}})$}{regular representation of $C^*({\mathchoice {\hbox {\rs B}} {\hbox {\rs B}} {\hbox {\rssmall B}} {\hbox {\rssmall B}}})$}{140}
\subjentry {\quad representation of a Fell bundle}{regular representation of a Fell bundle}{139}
\subjentry {\quad semigroup of isometries}{regular semigroup of isometries}{262}
\subjentry {\quad vertex}{regular vertex}{306}
\subjentry {restricted smash product}{restricted smash product}{222}
\goodbreak
restriction
\subjentry {\quad of a global action}{restriction of a global action}{15}
\subjentry {\quad of a unitary representation}{restriction of a unitary representation}{257}
\goodbreak
right
\subjentry {\quad Hilbert $A$-module}{right Hilbert $A$-module}{113}
\subjentry {\quad pre-Hilbert $A$-module}{right pre-Hilbert $A$-module}{112}
\subjentry {saturated Fell bundle}{saturated Fell bundle}{128}
\subjentry {section of a Fell bundle}{section of a Fell bundle}{133}
\subjentry {segment in a group}{segment in a group}{109}
\subjentry {self-adjoint}{self-adjoint}{30}
\subjentry {semi-direct product bundle}{semi-direct product bundle}{126}
\goodbreak
semi-saturated
\subjentry {\quad Fell bundle}{semi-saturated Fell bundle}{128}
\subjentry {\quad partial action}{semi-saturated partial action}{22}
\subjentry {\quad partial representation}{semi-saturated partial representation}{54}
\subjentry {semigroup of isometries}{semigroup of isometries}{261}
\subjentry {semilattice}{semilattice}{18}
\subjentry {separable Fell bundle}{separable Fell bundle}{230}
\subjentry {sink vertex}{sink vertex}{306}
\subjentry {smash product}{smash product}{220}
\goodbreak
source
\subjentry {\quad of a path}{source of a path}{301}
\subjentry {\quad vertex}{source vertex}{306}
\subjentry {spectrum of a set of relations}{spectrum of a set of relations}{103}
\subjentry {stable C*-algebra}{stable C*-algebra}{226}
\subjentry {stably isomorphic}{stably isomorphic}{116}
\subjentry {standard form}{standard form}{314}
\subjentry {stem}{stem}{317}
\subjentry {strong Morita equivalence}{strong Morita equivalence}{121}
\subjentry {summation process}{summation process}{165}
\subjentry {tame set of partial isometries}{tame set of partial isometries}{75}
\goodbreak
Toeplitz
\subjentry {\quad C*-algebra of a graph}{Toeplitz C*-algebra of a graph}{301}
\subjentry {\quad operator}{Toeplitz operator}{263}
\goodbreak
topological
\subjentry {\quad globalization}{topological globalization}{24}
\subjentry {\quad grading}{topological grading}{155}
\subjentry {\quad partial action}{topological partial action}{23}
\subjentry {\quad partial dynamical system}{topological partial dynamical system}{23}
\goodbreak
topologically
\subjentry {\quad equivalent partial actions}{topologically equivalent partial actions}{24}
\subjentry {\quad free partial action}{topologically free partial action}{247}
\subjentry {total space}{total space}{124}
\subjentry {transitive graph}{transitive graph}{337}
\subjentry {transitory cycle}{transitory cycle}{339}
\subjentry {unconditionally convergent}{unconditionally convergent}{145}
\subjentry {unilateral shift}{unilateral shift}{264}
\subjentry {unit fiber algebra}{unit fiber algebra}{124}
\goodbreak
unitary
\subjentry {\quad representation}{unitary representation}{51}
\subjentry {\quad semigroup dilation}{unitary semigroup dilation}{269}
\goodbreak
universal
\subjentry {\quad partial representation}{universal partial representation}{63}
\subjentry {\quad representation of a Fell bundle}{universal representation of a Fell bundle}{135}
\subjentry {upper bound}{upper bound}{270}
\subjentry {vanish at $\infty $}{vanish at $\infty $}{67}
\subjentry {vertices in a graph}{vertices in a graph}{300}
\subjentry {weak quasi-lattice}{weak quasi-lattice}{270}
\subjentry {weakly transitive graph}{weakly transitive graph}{337}
\subjentry {Wiener-Hopf algebra}{Wiener-Hopf algebra}{292}
\subjentry {word length}{word length}{21}

\bye